\documentclass[reqno]{amsart}

\usepackage{a4,amsthm,amsfonts,amssymb,euscript}
\usepackage{latexsym, multicol, fancybox}
\usepackage{graphicx}
\usepackage{color}
\usepackage{amsmath, amsthm, amssymb, bm}
\usepackage{epstopdf}
\usepackage{caption}
\usepackage{psfrag}

\usepackage{mathrsfs}
\usepackage{enumitem}
\usepackage{hyperref}

\usepackage{slashed,ulem,cancel}

 \usepackage{tikz}

\usepackage{a4,amsmath,amssymb,amsthm,slashed}
\usepackage[top=1in, bottom=1in, left=1.35in, right=1.35in]{geometry}
\newtheorem{theorem}{Theorem}[section]
\newtheorem{lemma}[theorem]{Lemma}
\newtheorem{proposition}[theorem]{Proposition}
\newtheorem{corollary}[theorem]{Corollary}
\newtheorem{definition}[theorem]{Definition}
\newtheorem{remark}[theorem]{Remark}
\newtheorem{conjecture}[theorem]{Conjecture}
\newtheorem{assumption}[theorem]{Assumption}

\numberwithin{equation}{section}

\newcommand{\bea}{\begin{eqnarray}}
\newcommand{\eea}{\end{eqnarray}}
\def\beaa{\begin{eqnarray*}}
\def\eeaa{\end{eqnarray*}}
\def\ba{\begin{array}}
\def\ea{\end{array}}
\def\be#1{\begin{equation} \label{#1}}
\def \eeq{\end{equation}}
\newcommand{\bsub}{\begin{subequations}}
\newcommand{\esub}{\end{subequations}}

\def\lab{\label}
\newcommand{\nn}{\nonumber}
\def\les{\lesssim}
\def\c{\cdot}
\def\f12{{\frac 1 2}}
\def\dual{{\,\,^*}}
\def\hot{\widehat{\otimes}}
\def\rhod{\,\dual\rho}

\def\err{\mbox{Err}}
\def\ov{\overline}

\def\a{{\alpha}}
\def\b{{\beta}}
\def\ga{\gamma}
\def\Ga{\Gamma}
\def\de{\delta}
\def\De{\Delta}
\def\ep{\epsilon}

\def\la{\lambda}
\def\La{\Lambda}

\def\Si{\Sigma}
\def\om{\omega}
\def\Om{\Omega}
\def\Th{\Theta}

\def\th{\theta}

\def\ze{\zeta}

\def\nab{\nabla}
\def\pr{{\partial}}
\def\div{{\mbox div\,}}
\def\curl{{\mbox curl\,}}

\def\f{\frac}

\def\th{\theta}

\def\AA{{\mathcal A}}
\def\BB{{\mathcal B}}

\def\DD{{\mathcal D}}

\def\FF{{\mathcal F}}

\def\HH{{\mathcal H}}
\def\II{{\mathcal I}}

\def\LL{{\mathcal L}}
\def\Lie{{\mathcal L}}
\def\MM{{\mathcal M}}
\def\NN{{\mathcal N}}
\def\OO{{\mathcal O}}
\def\PP{{\mathcal P}}
\def\QQ{{\mathcal Q}}

\def\A{{\bf A}}
\def\B{{\bf B}}

\def\D{{\bf D}}
\def\E{{\bf E}}
\def\F{{\bf F}}

\def\g{{\bf g}}
\def\H{{\bf H}}

\def\K{{\bf K}}
\def\L{{\bf L}}
\def\M{{\bf M}}

\def\O{{\bf O}}

\def\T{{\bf T}}

\def\R{{r^2 +a^2}}

\def\tr{\mbox{tr}}
\def\trch{{\mbox tr}\, \chi}
\def\chih{{\hat \chi}}
\def\chib{{\underline \chi}}
\def\chibh{{\underline{\chih}}}
\def\trchb{{\tr \,\chib}}
\def\etab{{\underline \eta}}
\def\omb{{\underline{\om}}}
\def\bb{{\underline{\b}}}
\def\aa{{\underline{\a}}}
\def\xib{{\underline \xi}}

\def\Hb{\,\underline{H}}

\DeclareFontFamily{U}{mathx}{\hyphenchar\font45}
\DeclareFontShape{U}{mathx}{m}{n}{
      <5> <6> <7> <8> <9> <10>
      <10.95> <12> <14.4> <17.28> <20.74> <24.88>
      mathx10
      }{}
\DeclareSymbolFont{mathx}{U}{mathx}{m}{n}
\DeclareFontSubstitution{U}{mathx}{m}{n}
\DeclareMathAccent{\widecheck}{0}{mathx}{"71}

\def\omc{\widecheck{\om}}

\def\etac{\widecheck{\eta}}

\def\Gac{\widecheck{\mathbf{\Ga}}}

\def\dk{\mathfrak{d}}

\def\reg{s}

\def\Reals{\mathbb{R}}

\providecommand{\abs}[1]{\lvert#1\rvert}

\def\trap{\text{trap}}
\def\nontrap{\cancel{\trap}}

\providecommand{\bf}[1]{\textbf{#1}}

\def\Opw{\mathbf{Op}_w}

\def\xit{\xi_{\tt}}
\def\xiphi{\xi_{\tphi}}
\def\xir{\xi_r}

\def\tt{\tau}
\def\tphi{\tilde{\phi}}
\def\tmod{t_{\text{mod}}}
\def\phimod{\phi_{\text{mod}}}
\def\dbl{\delta_{\textbf{BL}}}

\def\dred{\delta_{\text{red}}}
\def\dhor{\delta_{\HH}}
\def\dec{\de_{\text{dec}}}

\def\gcheck{\widecheck{\g}}

\def\gam{\g_{a,m}}

\def\Vref{V_{\text{ref}}}

\def\qs{|q|^2}

\def\EM{{\bf EM}}
\def\EMF{{\bf EMF}}

\def\N{{\bf N}}
\def\EF{{\bf EF}}

\newcounter{mnotecount}[section]

\def\Rbf{{\bf R}}
\def\Rdot{\dot{\Rbf}}

\def\Lied{\dot{\Lie}}
\def\Ddot{\dot{\D}}
\def\Db{\dot{\D}}
\def\squared{\dot{\square}}

\def\atrch{\atr\chi}
\def\atrchb{\atr\chib}
\def\atr{\,^{(a)}\mbox{tr}}

\def\sk{\mathfrak{s}}

\def\piX{\, ^{(X)}\pi}

\def\Div{\mbox{Div}}

\def\Jk{\mathfrak{J}}

\def\Xb{\protect\underline{X}}

\def\Ab{\protect\underline{A}}
\def\Bb{\protect\underline{B}}
\def\Xb{\protect\underline{X}}
\def\Hb{\,\underline{H}}

\def\Xib{\underline{\Xi}}
\def\Xh{\widehat{X}}
\def\Xbh{\widehat{\Xb}}
\def\squared{\dot{\square}}
     
          \def\Mtrap{\,\MM_{\trap}}
\def\Mntrap{{\MM_{\nontrap}}}

\def\tauu{\underline{\tau}}
\def\tauut{\widetilde{\tauu}}

\def\NNt{\widetilde{\NN}}
\def\NNtaux{\widetilde{\NN}_{\text{aux}}}
\def\NNtmora{\NNt_{\text{Mora}}}
\def\NNtdemora{\NNt_{\text{Mora}, \de}}
\def\NNtener{\NNt_{\text{Ener}}}

\def\NNtlede{\NNt_{\text{le},\de}}
\newcommand{\NNttotalp}[1]{\NNt_{#1,\de, \text{total}}[\pmb\phi_{#1}]}

\newcommand{\NNttotalpp}[1]{\NNt'_{#1,\de, \text{total}}}

\def\MF{{\bf MF}}

\def\bsplit{\begin{split}}

\def\prtan{\pr_{\text{tan}}}

\def\bsplit{\begin{split}}

\newcommand{\Lieb}{\Lie \mkern-10mu /\,}

\def\W{\mathbf{W}}
\def\Ab{\underline{A}}

\newcommand{\psiplus}[1]{\pmb\psi_{+2}^{(#1)}}
\newcommand{\psiminus}[1]{\pmb\psi_{-2}^{(#1)}}
\newcommand{\psis}[1]{\pmb\psi_{s}^{(#1)}}
\newcommand{\psiss}[2]{\psi_{s,#1}^{(#2)}}

\newcommand{\phiplus}[1]{\pmb\phi_{+2}^{(#1)}}
\newcommand{\phiminus}[1]{\pmb\phi_{-2}^{(#1)}}
\newcommand{\phis}[1]{\pmb\phi_{s}^{(#1)}}

\newcommand{\phipluss}[2]{\phi_{+2,#1}^{(#2)}}
\newcommand{\phiminuss}[2]{\phi_{-2,#1}^{(#2)}}
\newcommand{\phiss}[2]{\phi_{s,#1}^{(#2)}}

\def\reg{\mathbf{k}}

\def\Err{{\bf{Err}}}

\def\Xcal{\mathcal{X}}

\newcommand{\EMFtotalpp}[1]{\widetilde{\EMF}'_{#1,\de, \text{total}}}

\newcommand{\EMFtotalp}[1]{\widetilde{\EMF}_{#1, \de, \text{total}}[\pmb\phi_{#1}]}

\newcommand{\EMFtotalps}[1]{{\EMF}_{#1, \de, \text{total}}[\pmb\phi_{#1}]}

\def\Jk{\mathfrak{J}}

\def\Xb{\protect\underline{X}}

\def\be{{\beta}}

\def\lap{{\triangle}}

\def\NNtlocal{\NNt_{\text{local}}}

\def\Zc{\widecheck{Z}}
\def\Hc{\widecheck{H}}
\def\Hbc{\widecheck{\Hb}}
\def\trXc{\widecheck{\tr X}}
\def\trXbc{\widecheck{\tr\Xb}}
\def\Pc{\widecheck{P}}

\def\Gac{\widecheck{\Ga}}

\def\omc{\widecheck \omega}

\newcommand{\IE}[1]{{\bf{IE}}[#1]}
\def\Ao{\protect\overline{\A}}
\def\Extra{\B_{\de}}

\def\DDs{ \, \DD \hspace{-2.4pt}\dual    \mkern-20mu /}
\def\DDd{ \, \DD \hspace{-2.4pt}    \mkern-8mu /}

\newcommand{\Rmic}{R_0}
\newcommand{\Nmic}{N_0}
\newcommand{\tmic}{\tau_{\Nmic}}
\newcommand{\Iti}{I_{\Nmic}}

\newcommand{\psish}[2]{\pmb\psi_{s}^{(#1), #2}}
\newcommand{\psissh}[3]{\psi_{s,#1}^{(#2), #3}}

\newcommand{\phish}[2]{\pmb\phi_{s}^{(#1), #2}}

\newcommand{\phiplussh}[3]{\phi_{+2,#1}^{(#2),#3}}
\newcommand{\phiminussh}[3]{\phi_{-2,#1}^{(#2), #3}}
\newcommand{\phissh}[3]{\phi_{s,#1}^{(#2), #3}}

\newcommand{\NNttotalh}[2]{\NNt_{#1,\de, \text{total}}^{(#2)}[\pmb\phi_{#1}]}
\newcommand{\EMFtotalh}[2]{\widetilde{\EMF}_{#1,\de, \text{total}}[#2\pmb\phi_{#1}]}
\newcommand{\EMFtotalhps}[2]{\EMF_{#1,\de, \text{total}}[#2\pmb\phi_{#1}]}

\newcommand{\NNttotalph}[2]{\NNt_{#1, \de, \text{total}}[#2\pmb\phi_{#1}]}

\def\dkb{ \, \mathfrak{d}     \mkern-9mu /}

\def\Errprod{\Err_{0\times1}}

\def\Errdefect{\E_{\mathrm{defect}}}

\def\good{\mathrm{Good}}

\newcommand{\Bulk}[1]{\mathbf{Bulk}_{#1,\pr_{\tt}}[\pmb\psi]}
\newcommand{\Bulkxw}[1]{\mathbf{Bulk}_{#1,(X,w)}[\pmb\psi]}

\def\prtphihat{\widehat{\pr}_{\tphi}}

\newcommand{\IEde}[1]{{\bf{IE}_{\de}}[#1]}

\def\phihp{\widehat{\pmb\phi}_{+2}}
\def\Nhp{\widehat{\N}_{+2}}

\begin{document}

\title{Energy-Morawetz estimates for Teukolsky equations in perturbations of Kerr}
\author{Siyuan Ma and J\'{e}r\'{e}mie Szeftel}

\begin{abstract}
In this paper, we prove energy and Morawetz estimates for solutions to Teukolsky equations in spacetimes with metrics that are perturbations, compatible with nonlinear applications, of Kerr metrics in the full subextremal range. The Teukolsky equations are written in tensorial form using the non-integrable formalism in \cite{GKS22}, and we follow the approach in \cite{Ma} of relying on a Teukolsky wave/transport system. The estimates are proved by extending the ideas from our earlier result \cite{MaSz24} on the corresponding problem for the scalar wave, notably the use of $r$-foliation-adapted microlocal multipliers for the wave part, and by incorporating techniques from \cite{Ma} to control the linear coupling terms between the components of the Teukolsky wave/transport system. Additionally, in order to adapt the methodology of \cite{MaSz24} to tensorial waves, we introduce a well-suited regular scalarization procedure which is of independent interest. This result, alongside our companion paper \cite{MaSz24}, is an essential step towards extending the current proof of Kerr stability in \cite{GCM1} \cite{GCM2} \cite{KS:Kerr} \cite{GKS22} \cite{Shen}, valid in the slowly rotating case, to a complete resolution of the Kerr stability conjecture, i.e., the statement that the Kerr family of spacetimes is nonlinearly stable for all subextremal angular momenta.
\end{abstract}

\maketitle

\tableofcontents


\section{Introduction}



\subsection{Kerr stability conjecture}


We begin with introducing the Einstein vacuum equations, the Kerr solutions and the Kerr stability conjecture.

The Einstein vacuum equations are given by
\bea\lab{eq:EVE}
\textbf{Ric}(\g)=0,
\eea
where $(\MM, \g)$ is a four-dimensional Lorentzian manifold, and where $\textbf{Ric}(\g)$ denotes the Ricci curvature tensor of the metric $\g$. The Kerr spacetimes \cite{Kerr63} represent a family of asymptotically flat, stationary, axially symmetric black hole solutions to the Einstein vacuum equations \eqref{eq:EVE}. The metrics of the \textit{subextremal} Kerr spacetimes, parameterized by the mass $m$ and an angular momentum per unit mass $a$ with the strict inequality $|a|<m$, take the following form in the Boyer--Lindquist  \cite{BL67} coordinates $(t,r,\th, \phi)$
\bea\lab{eq:expressionofKerrmetricinBLcoordinates:intro}
\gam=-\frac{\Delta \qs}{\Sigma^2} dt^2 + \frac{\sin^2\th\Sigma^2 }{\qs}\bigg(d\phi - \frac{2amr}{\Sigma^2} dt\bigg)^2 +\frac{\qs}{\Delta} dr^2 + \qs d\th^2,
\eea
with functions
\bea
\Delta = r^2 - 2mr +a^2, \quad \qs=r^2+a^2\cos^2\th, \quad \Sigma^2=(\R)^2 - a^2\sin^2\th \Delta.
\eea
In particular, such a subextremal Kerr spacetime contains a black hole $\{r<r_+\}$ with a nondegenerate event horizon located at $\{r=r_+\}$ where $r_+:=m+\sqrt{m^2-a^2}$ is the larger root of $\Delta=\De(r)$,  see Figure \ref{fig:penrosediagramofKerr} for the corresponding Penrose diagram.

\begin{figure}[htbp]
  \begin{center}
\begin{tikzpicture}[scale=1.1]
\tikzstyle{every node}=[font=\Small]
      \draw[dashed, color=red, thick] (0.05,3.95) -- (3.2,0.8);
  \fill[yellow!50] (-0.05,3.95)--(-3.2,0.8) arc(225:315: 4.51 and 4.51) -- (0.05,3.95);
   \draw[dashed, color=red, thick] (-0.05,3.95)--(-3.2,0.8) ;
     \node at (-2.8,3) {Black hole region};
     \node at (0.1,4.3) {$i_+$};
       \node[rotate=315]  at (2.0,2.4) {$\II_+$};
       \node[rotate=45] at (-1.6,2.1) {Event horizon $r=r_+$};
         \draw[] (0,4) circle (0.05);
         \node at (0,1) {Domain of outer communication};  
\end{tikzpicture}
\end{center}
\caption{\footnotesize{Penrose diagram of subextremal Kerr spacetimes.}}
\lab{fig:penrosediagramofKerr}
\end{figure}
 Note that the special case $a=0$ with $m>0$ corresponds to the family of Schwarzschild spacetimes, introduced by Schwarzschild \cite{Sch16} in 1916.

The \textit{Kerr stability conjecture}, one of the central open problems in general relativity, aims to prove the following statement.

\begin{conjecture}[Kerr stability conjecture]
The maximal Cauchy development of any initial data set for Einstein vacuum equations, that is sufficiently close to a subextremal Kerr initial data in a suitable sense, has a complete future null infinity and a domain of outer communication\footnote{The domain of outer communication is the causal past of future null infinity.} which is asymptotic to a nearby member of the subextremal Kerr family. 
\end{conjecture}

The most recent breakthrough towards a resolution of the Kerr stability conjecture is its proof in the slowly rotating case (that is, $|a|/m\ll 1$), established in the series of works \cite{GCM1} \cite{GCM2} \cite{KS:Kerr} \cite{GKS22} \cite{Shen}. A complete resolution of the conjecture, i.e., removing the restriction on the angular momentum parameter $a$, requires establishing energy and Morawetz estimates for both the scalar wave equation and the Teukolsky equations\footnote{The Teukolsky equations constitute two fundamental equations within the system of Einstein vacuum equations, see Section \ref{subsect:Teukolskywavetransport:intro}.} in suitable perturbations of any subextremal Kerr spacetime. Such estimates for the scalar wave equation were shown in our companion paper \cite{MaSz24}, and the goal of the present paper is to establish corresponding estimates for the Teukolsky equations.


\subsection{Teukolsky wave/transport system}
\lab{subsect:Teukolskywavetransport:intro}


Kerr spacetimes possess a distinguished pair of null vectorfields known as the principal null pair, see \eqref{def:e3e4inKerr}, that diagonalizes the curvature tensor. In perturbations of Kerr, we consider a pair of null vectorfields $(e_3, e_4)$, normalized by $\g(e_3, e_4)=-2$, which is a suitable perturbation of the principal null pair of Kerr. We then consider an orthonormal pair of spacelike vectorfields $(e_1, e_2)$ spanning the horizontal bundle $\{e_3, e_4\}^\perp$, see \eqref{def:e1e2inKerr} in Kerr, so that $(e_3, e_4, e_a)$, $a=1,2$, forms a null frame of the spacetime $(\MM, \g)$. As in \cite{GKS20} \cite{GKS22}, we associate to the null pair $(e_3, e_4)$ horizontal tensors and denote in particular by $\sk_2(\Reals)$ the set of symmetric traceless horizontal real $2$-tensors, see Section \ref{subsection:review-horiz.structures}.

Next, we denote the curvature components $\a, \aa\in\sk_2(\Reals)$ by
\beaa
\a_{ab}={\bf R}_{a4b4},\qquad \aa_{ab}={\bf R}_{a3b3},\qquad a,b=1,2,
\eeaa
where ${\bf R}_{\a\b\mu\nu}$ denotes the curvature tensor of the spacetime $(\MM, \g)$. Also, we define the complexified curvature components $A, \Ab\in \sk_2(\mathbb{C})$ as
\beaa
A=\a+i\dual\a, \qquad \Ab=\aa+i\dual\aa,
\eeaa
where $\sk_2(\mathbb{C})$ is introduced in Definition \ref{def:skC:horizontaltensors}  as the set of symmetric, traceless, anti-self-dual horizontal complex $2$-tensors. The \textit{Teukolsky equations} \cite{Teu72}, the governing equations for these curvature components, are, in the tensorial formalism introduced in \cite{GKS20} \cite{GKS22}, given by\footnote{See \eqref{eq:TeukolskyequationforAandAbintensorialforminKerr} for the form of \eqref{eq:Teu:intro} in Kerr with the normalization \eqref{eq:defintionoftensorialTeukolskyscalarspsipm2}.}
\bea
\lab{eq:Teu:intro}
\mathcal{T}_{+2, \g} A = \N_A, \qquad \mathcal{T}_{-2,\g}\Ab=\N_{\Ab},
\eea
where $\mathcal{T}_{\pm 2,\g}$ are tensorial Teukolsky wave operators in $(\MM, \g)$ and where $\N_{A}$ and $\N_{\Ab}$ are source terms\footnote{For the explicit formulas of $\N_{A}$ and $\N_{\Ab}$ in terms of the Ricci coefficients and curvature components of a perturbation of Kerr $(\MM, \g)$, see Sections 5.1.1 and 5.3.1 in \cite{GKS22}.}. 

The heart of the analysis in this paper relies on a \textit{wave/transport hierarchy} constructed from the Teukolsky equations \eqref{eq:Teu:intro}. 
Following \cite{Ma}, adapted to the tensorial formalism of \cite{GKS22},  we consider  tensors $\pmb\phi_s^{(p)}\in\sk_2(\mathbb{C})$, $s=\pm 2$\footnote{In this paper, $s$ refers to the spin weight of the tensors.}, $p=0,1,2$, with $\pmb\phi_s^{(0)}$ given by 
\bsub
\lab{def:TensorialTeuScalars:wavesystem:Kerrperturbation:intro}
\bea
\pmb\phi_{+2}^{(0)}=\frac{\ov{q}}{q}A, \qquad \pmb\phi_{-2}^{(0)}=\frac{q}{\ov{q}}\left(\frac{\De}{|q|^2}\right)^2\Ab
\eea
and with $\phis{p}$, $s=\pm 2$, $p=0,1,2$, satisfying the Teukolsky transport equations
 \bea
\lab{def:TensorialTeuScalars:wavesystem:Kerrperturbation:+2:intro}
\nab_3 \left(\frac{r\bar{q}}{q}\left(\frac{r^2}{|q|^2}\right)^{p-2}\pmb\phi_{+2}^{(p)}\right)&=&\frac{\bar{q}}{rq}\left(\frac{r^2}{|q|^2}\right)^{p-1}\pmb\phi_{+2}^{(p+1)}+\N_{T,+2}^{(p)}, \quad p=0,1,\\
\lab{def:TensorialTeuScalars:wavesystem:Kerrperturbation:-2:intro}
\nab_4\left(\frac{rq}{\bar{q}}\left(\frac{{r^2}}{|q|^2}\right)^{p-2}\pmb\phi_{-2}^{(p)}\right)&=&\frac{q}{r\bar{q}}\left(\frac{r^2}{|q|^2}\right)^{p-1}\frac{\De}{\qs}\pmb\phi_{-2}^{(p+1)}+\N_{T,-2}^{(p)}, \,\,\,\, p=0,1,
\eea
\esub
where $q:=r+ia\cos\th$ and where $\N_{T, s}^{(p)}$, $s=\pm 2$, $p=0,1$, are source terms in the transport equations. These tensors $\phis{p}$ satisfy the coupled Teukolsky wave equations
\bsub
\lab{eq:TensorialTeuSysandlinearterms:rescaleRHScontaine2:general:Kerrperturbation:intro}
\bea
\lab{eq:TensorialTeuSys:Kerrpert:intro}
\bigg(\squared_2 -\frac{4ia\cos\th}{|q|^2}\nab_{\pr_t}- \frac{4-2\de_{p0}}{\qs}\bigg){\phis{p}} = \L_{s}^{(p)}[\pmb\phi_{s}]+\N_{W,s}^{(p)}, \quad s=\pm2, \quad p=0,1,2,
\eea
where $\squared_2$ is the tensorial wave operator for tensors in $\sk_2(\mathbb{C})$, see \eqref{eq:def=squared-2}, and where the linear coupling terms $\L_{s}^{(p)}[{\pmb\phi_s}]$ have the following schematic form
\bea
\lab{eq:tensor:Lsn:onlye_2present:general:Kerrperturbation:intro}
\bsplit
{\L_{s}^{(0)}[\pmb\phi_{s}]}={}& (2sr^{-3} +O(mr^{-4}))\phis{1}+ O(mr^{-3}) \nab_{\Xcal_s}^{\leq 1}\phis{0},\\
{\L_{s}^{(1)}[\pmb\phi_{s}]}={}& (sr^{-3} +O(mr^{-4}))\phis{2}+ O(mr^{-3}) \nab_{\Xcal_s}^{\leq 1}  \phis{1}+O(mr^{-2})\nab_{\pr_{\phi}+a\pr_t}^{\leq 1}\phis{0},\\
{\L_{s}^{(2)}[\pmb\phi_{s}]}={}&O(mr^{-3})\phis{2}+O(mr^{-2})\nab_{\pr_{\phi}+a\pr_t}^{\leq 1}\phis{1}+O(m^2 r^{-2})\phis{0},
\end{split}
\eea
\esub
with $\Xcal_s$, $s=\pm 2$, being regular vectorfields that are horizontal in the case of Kerr\footnote{See \eqref{eq:formofregularhorizontalvectorfieldmathcalXs} for the form of $\Xcal_s$, $s=\pm 2$, in Kerr.}. The equations \eqref{def:TensorialTeuScalars:wavesystem:Kerrperturbation:intro} \eqref{eq:TensorialTeuSysandlinearterms:rescaleRHScontaine2:general:Kerrperturbation:intro} correspond to the tensorial \textit{Teukolsky wave/transport system} in perturbations of Kerr considered throughout this paper.


\subsection{State of the art on energy-Morawetz estimates for Teukolsky equations}
\lab{subsect:literature:Teukolsky:intro}


The analysis of the Teukolsky equations is central to understanding the dynamical evolution of Kerr spacetimes and fundamentally builds upon the framework developed for the scalar wave equation. For an in-depth review of the literature concerning scalar waves, we direct the reader to the introduction in our companion  paper \cite{MaSz24}. In this section, we review the literature pertaining to energy-Morawetz estimates for solutions to Teukolsky equations.


\subsubsection{Teukolsky equations in Kerr}
\lab{subsubsect:literature:TeukolskyinKerr:intro}


In order to derive energy-Morawetz estimates, one must first address the question of mode stability for solutions to Teukolsky equations in Kerr spacetimes. This was achieved in 1989 in the seminal work of Whiting \cite{Whiting}, who demonstrated, under  no incoming radiation assumption, that no exponentially growing mode solutions exist. It was later extended in \cite{AMPW17}, see also \cite{TdC20}, to show the absence of non-trivial mode solutions with real frequencies.

In Schwarzschild spacetimes,  energy-Morawetz estimates for the Teukolsky equations were first obtained by Dafermos-Holzegel-Rodnianski \cite{DHRT19Schw}. The proof relies on a physical space analog of the Chandrasekhar's transformation \cite{Chan75} that converts the Teukolsky equations into a Regge-Wheeler  type wave equation \cite{RW57}, to which the techniques developed for the scalar wave equation can be directly applied. Generalizations to Kerr spacetimes were achieved in the slowly rotating case by Ma \cite{Ma} and Dafermos-Holzegel-Rodnianski \cite{DHR19}, and for the full subextremal range by Shlapentokh-Rothman-Teixeira da Costa \cite{SRTdC20, SRTdC23} and Millet \cite{Millet}\footnote{While \cite{Millet} derives sharp decay estimates for solutions to Teukolsky equations that do not rely on energy-Morawetz estimates, one can easily adapt the methodology in that paper to derive such estimates, though with a loss of several derivatives, see Section \ref{sec:weakMorawetzestimatesforTeukolskyinKerrusingMillet}.}.


\subsubsection{Teukolsky equations in perturbations of Kerr with $|a|\ll m$}
\lab{subsubsect:literature:TeukolskyinpertKerr:intro}


To address the nonlinear stability of Kerr, it is important to extend the energy-Morawetz estimates for Teukolsky equations in Kerr reviewed in Section \ref{subsubsect:literature:TeukolskyinKerr:intro} to small perturbations of Kerr.  This has been achieved in the context of the recent proofs of the nonlinear stability of Schwarzschild and of Kerr spacetimes for $|a|\ll m$: see Chapter 10 of \cite{KS20} in the context of the nonlinear stability of Schwarzschild under polarized axisymmetry, Chapters 12 and 13 of \cite{DHRT21} in the context of the nonlinear stability of Schwarzschild spacetimes for a codimension-3 set of initial data, and Chapter 9 of \cite{GKS22} in the context of the nonlinear stability of slowly rotating Kerr, i.e., with $|a|\ll m$. Generalizing these results to perturbations of any subextremal Kerr spacetime remains open and is the objective of the present paper.


\subsection{First version of the main result}


Given constants $(a, m)$ with $|a|<m$ and $0<\dhor\ll m-|a|$,  let the spacetime $(\MM, \g)$, whose Penrose diagram is depicted in Figure \ref{fig:penrosediagramofM}, be such that:
\begin{itemize}
\item $\MM=\{(\tau, r, \om)\,/\,\tau\in\Reals,  r_+(1-\dhor)\leq r<+\infty,  \om \in\mathbb{S}^2\}$ is a four dimensional manifold, where $(\tau, r)$ are two coordinates on $\MM$ and $r_+=m+\sqrt{m^2-a^2}$, 

\item the boundary $\AA:=\{r=r_+(1-\dhor)\}$ of $\MM$ is spacelike,

\item the level sets of the time function $\tau$ are spacelike, transversal to the hypersurface $\AA$ and asymptotically null as $r\to +\infty$.
\end{itemize}

\begin{figure}[htbp]
  \begin{center}
\begin{tikzpicture}[scale=1.1]
\tikzstyle{every node}=[font=\Small]
  \draw[color=black] (-3.2,1.84) arc(300:324.65:10 and 7.5);
    \draw[dashed, color=red] (-0.05,3.95)--(-3.2,0.8) ;
      \draw[dashed, color=red] (0.05,3.95) -- (3.2,0.8);
     \node[rotate=32] at (-2.3,2.58) {$\AA$};
     \node at (0.1,4.3) {$i_+$};
       \node[rotate=315]  at (2.0,2.4) {$\II_+$};
       \node[rotate=45] at (-1.5,2.1) {$\HH_+$};
         \draw[] (0,4) circle (0.05);
  \draw[blue] (-2.88,1.98) arc(220:313.5:3.53 and 2.9);
         \node at (0,0.7) {$\Sigma(\tau_1)$};
    \draw[blue] (-1.77,2.595) arc(240:302.7:3 and 2.2);
      \node at (0,2.58) {$\Sigma(\tau_2)$}; 
\end{tikzpicture}
\end{center}
\caption{\footnotesize{Penrose diagram of $(\MM, \g)$. $\Si(\tau_1)$ and $\Si(\tau_2)$ are two spacelike and asymptotically null level hypersurfaces of a function $\tau$, and $\AA=\{r=r_+(1-\dhor)\}$ is spacelike.}}\lab{fig:penrosediagramofM}
\end{figure}

Our main result is the derivation of energy-Morawetz-flux (EMF) estimates for solutions to the Teukolsky wave/transport system \eqref{def:TensorialTeuScalars:wavesystem:Kerrperturbation:intro} \eqref{eq:TensorialTeuSysandlinearterms:rescaleRHScontaine2:general:Kerrperturbation:intro} in spacetimes $(\MM,\g)$, where $\g$ is a perturbation of a Kerr metric $\gam$ with $|a|<m$. We provide below a rough version of our main theorem, see Theorem \ref{thm:main} for the precise version.

\begin{theorem}[Main theorem, rough version] 
\label{thm:main:Teu:rough}
Let $(\MM, \g)$ be a perturbation\footnote{More precisely, the spacetime $(\MM,\g)$ is assumed to satisfy the assumptions on a null pair, the metric perturbation,  and a regular triplet of horizontal vectorfields made in Sections \ref{subsect:assumps:perturbednullframe}, \ref{subsubsect:assumps:perturbedmetric} and \ref{sec:regulartripletinperturbationsofKerr}, respectively.} of a Kerr spacetime with metric $\gam$ satisfying $|a|<m$ in the sense of Sections \ref{subsect:assumps:perturbednullframe}, \ref{subsubsect:assumps:perturbedmetric} and \ref{sec:regulartripletinperturbationsofKerr}. Then, we have for solutions to the Teukolsky wave/transport system \eqref{def:TensorialTeuScalars:wavesystem:Kerrperturbation:intro}-\eqref{eq:TensorialTeuSysandlinearterms:rescaleRHScontaine2:general:Kerrperturbation:intro} the following EMF estimates, for $s=\pm2$, any $1\leq\tau_1<\tau_2 <+\infty$ and any given $0<\de\leq \frac{1}{3}$,
\bea
&&\sum_{p=0,1,2}\bigg(\sup_{\tau\in [\tau_1, \tau_2]}\E^{(\reg_{s})}[\phis{p}](\tau) + {\M^{(\reg_{s})}_\de}[\phis{p}](\tau_1,\tau_2) +{\F}^{(\reg_{s})}[\phis{p}](\tau_1, \tau_2) \bigg)\nn\\
&\les&\sum_{p=0,1,2}\E^{(\reg_{s})}[\phis{p}](\tau_1)+\sum_{p=0,1,2}{\mathcal{N}^{(\reg_{s})}_\de}[\phis{p}, \N_{W,s}^{(p)}, \N_{T,s}^{(p)}](\tau_1, \tau_2).
\eea
Here, $\reg_s$ are integers measuring regularity, $\E^{(\reg_s)}[{\cdot}](\tau)$, ${\F}^{(\reg_s)}[{\cdot}](\tau_1, \tau_2)$ and $\M^{(\reg_s)}_\de[\cdot](\tau_1,\tau_2)$ are $\reg_s$-th order energy on a constant-$\tau$ hypersurface $\Sigma(\tau)$,  fluxes on both $\AA(\tau_1,\tau_2)$ and $\II_+(\tau_1,\tau_2)$ and Morawetz terms over $\MM(\tau_1,\tau_2)$, the term ${\mathcal{N}^{(\reg_{s})}_\de}[\phis{p}, \N_{W,s}^{(p)}, \N_{T,s}^{(p)}](\tau_1, \tau_2)$ corresponds to the contribution of the inhomogeneous terms $\N_{W,s}^{(p)}$ and  $\N_{T,s}^{(p)}$, and the implicit constants in $\lesssim$ are independent of  $\tau_1$ and $\tau_2$, and depend only on the black hole parameters $a$ and $m$, as well as on the constants $\dhor$ and $\de$. 
\end{theorem}

\begin{remark}
Here are some comments on the statement of Theorem \ref{thm:main:Teu:rough}: 
\begin{itemize}
\item the assumptions on  the spacetime $(\MM,\g)$ made in the theorem are consistent with the estimates in the proof of the nonlinear stability of Kerr for small angular momentum in \cite{KS:Kerr};

\item though our estimates are closed with a specific choice of the pair of regularity integers $\reg_s$, $s=\pm 2$, namely\footnote{See the last item in Remark \ref{remark:ofmainthm7.1} for an explanation of this choice.} $\reg_{+2}=11$ and $\reg_{-2}=14$, the extension to higher-order derivatives can be derived in the same manner as in our proof;

\item the methodology of deriving pointwise decay estimates in perturbations of Kerr in the full range $|a|<m$ starting from Theorem \ref{thm:main:Teu:rough} is by now standard;

\item the statement of Theorem \ref{thm:main:Teu:rough} is new even when restricted to subextremal Kerr:
\begin{itemize}
\item while the energy-Morawetz estimates in subextremal Kerr in \cite{SRTdC20, SRTdC23} are proved for the homogenous case, i.e., in the case $\N_{W,s}^{(p)}=\N_{T,s}^{(p)}=0$, $s=\pm 2$, $p=0,1,2$, Theorem \ref{thm:main:Teu:rough} restricted to subextremal Kerr holds in the general inhomogeneous case;

\item the energy-Morawetz estimates in Theorem \ref{thm:main:Teu:rough} require only standard energy bounds of the initial data for the Teukolsky wave system, which, in terms of the fall-off of the initial data, is both optimal and weaker than all corresponding results \cite{Ma} \cite{DHR19} \cite{SRTdC20, SRTdC23} in Kerr spacetimes.
\end{itemize}
\end{itemize}
\end{remark}

\begin{remark}[Relevance to the Kerr stability conjecture in the full subextremal range]
\lab{rem:thm:rough:Teu:intro}
In the proof of Kerr stability for $|a|\ll m$ in \cite{GCM1,GCM2,KS:Kerr,GKS22,Shen},  the assumption $|a|\ll m$ is only needed in \cite{GKS22} for the derivation of the main energy-Morawetz estimates for the scalar wave equation, Teukolsky equations, and Bianchi identities in perturbations of Kerr. Our main result in Theorem \ref{thm:main:Teu:rough} thus fills a crucial gap in extending the proof of Kerr stability for $|a|\ll m$ in \cite{GCM1,GCM2,KS:Kerr,GKS22,Shen} to the full subextremal range $|a|<m$.
\end{remark}


\subsection{Strategy of the proof}
\lab{subsect:proof:intro}


In this section, we provide an outline of the strategy of the proof of our main Theorem \ref{thm:main:Teu:rough}.


\subsubsection{Regular scalarization of the tensorial Teukolsky wave/transport system}


In order to establish energy-Morawetz estimates for the tensorial Teukolsky wave/transport system  \eqref{def:TensorialTeuScalars:wavesystem:Kerrperturbation:intro}  \eqref{eq:TensorialTeuSysandlinearterms:rescaleRHScontaine2:general:Kerrperturbation:intro} on $\MM(\tau_1,\tau_2)$, we will rely on microlocal energy-Morawetz multipliers which in turn requires an extension to a semi-global in time problem. Now, such an extension is significantly more cumbersome at the level of tensors than for scalars and we will thus need to scalarize the Teukolsky wave/transport system  \eqref{def:TensorialTeuScalars:wavesystem:Kerrperturbation:intro}  \eqref{eq:TensorialTeuSysandlinearterms:rescaleRHScontaine2:general:Kerrperturbation:intro}. A possible candidate is the Newman-Penrose formalism \cite{NP62, NP63errata} where the scalarization is done with respect to (w.r.t.) the horizontal frame $(e_1, e_2)$. However, the horizontal vectors $e_1, e_2$ are singular at $\th=0,\pi$, and hence, the scalars defined within the Newman-Penrose formalism are not globally regular, thus presenting a major obstacle for our analysis.

Motivated by this, we introduce in Section \ref{sec:regularscalarization} a novel approach of independent interest for the scalarization of horizontal tensors which both leads to regular scalars and is amenable to the extension to a semi-global problem and to the use of microlocal multipliers. Its fundamental difference with the Newman-Penrose formalism lies in the choice of the horizontal vectorfields used for the scalarization: we replace the singular horizontal frame $(e_1, e_2)$ with a regular triplet of horizontal vectorfields $\{\Om_i\}_{i=1,2,3}$, which spans the horizontal bundle and satisfies on $(\MM, \g)$ the following set of algebraic identities 
\bea\lab{eq:fundamentalpropertiesof1formsOmi:intro}
x^i\Om_i=0, \qquad (\Om^i)_a(\Om_i)_b=\de_{ab}, \qquad \Om_i\c\Om_j=\de_{ij}-x^ix^j, \qquad \Om_i\c\dual\Om_j=\in_{ijk}x^k,
\eea
see Definition \ref{def:regulartripletinKerrOmii=123} for the choice of such a regular triplet in Kerr and for the definition of $x^i$, $i=1,2,3$. This procedure allows us to project a tensor $\pmb\psi \in \sk_2(\mathbb{C})$ onto this regular triplet to produce  a family of regular complex-valued scalars $\psi_{ij}=\pmb\psi(\Om_i, \Om_j)$  satisfying the following set of algebraic constraints 
\bea\lab{eq:tensorizationconstraints:intro}
\psi_{ij}=\psi_{ji}, \qquad x^i\psi_{ij}=0, \qquad (\de^{ij}-x^ix^j)\psi_{ij}=0, \qquad \in_{ikl}x^l\psi_{kj}+i\psi_{ij}=0,
\eea
and, reciprocally, to reconstruct the tensor $\pmb\psi\in \sk_2(\mathbb{C})$ from the scalars $\psi_{ij}$ satisfying such a set of constraints, thus providing a one-to-one correspondence between tensors in $\sk_2(\mathbb{C})$ and families of regular scalars satisfying the set of algebraic constraints \eqref{eq:tensorizationconstraints:intro}, see Lemma \ref{lemma:backandforthbetweenhorizontaltensorsk2andscalars:complex}.

Relying on this scalarization procedure, the tensorial Teukolsky wave/transport system \eqref{def:TensorialTeuScalars:wavesystem:Kerrperturbation:intro}  \eqref{eq:TensorialTeuSysandlinearterms:rescaleRHScontaine2:general:Kerrperturbation:intro} is transformed into a scalarized Teukolsky wave/transport system\footnote{See Lemma \ref{lemma:formoffirstordertermsinscalarazationtensorialwaveeq} for the scalarization of the tensorial wave operator $\squared_2$, Lemma \ref{lem:scalarizedTeukolskywavetransportsysteminKerr:Omi} for the scalarized Teukolsky wave/transport system in Kerr and Lemma \ref{lem:scalarizedTeukolskywavetransportsysteminKerrperturbation:Omi} for the corresponding one in perturbations of Kerr.} with the scalarized Teukolsky wave equations having the following form
\bea
\lab{eq:ScalarizedTeuSys:general:Kerrperturbation:intro} 
\square_\g\phiss{ij}{p}- \widehat{S}(\phi_s^{(p)})_{ij} -(\widehat{Q}\phi_s^{(p)})_{ij} -\frac{4-2\de_{p0}}{\qs} \phiss{ij}{p} = {L_{s,ij}^{(p)}}+N_{W,s,ij}^{(p)}, 
\eea
where 
\bea
\lab{eq:scalarizationoftensors:intro}
\phiss{ij}{p}=\pmb\phi_s^{(p)}(\Om_i,\Om_j), \quad L_{s,ij}^{(p)}=\L_{s}^{(p)}[{\pmb\phi_s}](\Om_i,\Om_j), \quad N_{W,s,ij}^{(p)}=\N_{W,s}^{(p)}(\Om_i,\Om_j),
\eea
and where $\widehat{S}(\phi_s^{(p)})_{ij}$ and $(\widehat{Q}\phi_s^{(p)})_{ij}$ denote the first- and zeroth-order scalarization coupling terms, respectively.
While this system involves a significantly larger family of equations (each tensorial equation is transformed into a system of coupled scalar equations), each equation is essentially a standard scalar wave equation modulo
lower-order terms\footnote{Both the scalarization coupling terms $\widehat{S}(\phi_s^{(p)})_{ij}$ and $(\widehat{Q}\phi_s^{(p)})_{ij}$ and the linear coupling terms $L_{s,ij}^{(p)}$ are of first- or zeroth-order.}, and is hence suitable for a semi-global extension and for applying the microlocal energy-Morawetz estimates developed in our previous work \cite{MaSz24} on scalar waves. In order to derive energy-Morawetz estimates for the scalarized Teukolsky wave/transport system \eqref{eq:ScalarizedTeuSys:general:Kerrperturbation:intro}, the remaining obstacle is then to deal with the two kinds of coupling terms, i.e., the scalarization coupling terms $\widehat{S}(\phi_s^{(p)})_{ij}$ and $(\widehat{Q}\phi_s^{(p)})_{ij}$, and the scalarized linear coupling terms ${L_{s,ij}^{(p)}}$.


\subsubsection{Extension to a semi-global problem}
\lab{subsubsect:Extensiontoglobal:intro}


Building on the approach from our previous work \cite{MaSz24} for a single inhomogeneous scalar wave equation, we need to extend the local problem for the scalarized Teukolsky wave/transport system \eqref{eq:ScalarizedTeuSys:general:Kerrperturbation:intro} to a semi-global setting. This extension is necessary to use microlocal multipliers\footnote{This is due to the fact that our multipliers are microlocal in particular w.r.t. the time coordinate $\tau$.} within a bounded \( r \)-region which contains the trapping region.  To this end, we extend the scalarized Teukolsky wave equations \eqref{eq:ScalarizedTeuSys:general:Kerrperturbation:intro}  for $\phiss{ij}{p}$, which are originally defined in $\MM(\tau_1,\tau_2)$, to a semi-global in time problem on $\MM(\tmic, +\infty)$\footnote{Here, $\tmic$ is a time parameter that is smaller than but remains close to the initial time parameter $\tau_1$.} for a system of coupled scalar wave equations for $\psiss{ij}{p}$ of the type
\bea\label{eq:waveeqwidetildepsi1:intro}
\bigg({\square}_{{\g_{\chi_{\tau_1, \tau_2}}}} -\frac{4-2\delta_{p0}}{|q|^{2}}\bigg)\psi^{(p)}_{s,ij} = \widehat{F}_{s,ij}^{(p)},
\eea
see Section \ref{sect:extensiontoglobalproblem:Teu}, where:
\begin{enumerate}
\item The interpolated metric $\g_{\chi_{\tau_1, \tau_2}}$ equals $\g$ in $\MM(\tau_1+1, \tau_2-2)$, matches the Kerr metric $\gam$ in $\MM\setminus\MM(\tau_1, \tau_2-1)$, and continues to satisfy the metric perturbation assumptions of Section \ref{subsubsect:assumps:perturbedmetric}.

\item The scalars $\psiss{ij}{p}$ satisfy $\psiss{ij}{p}=\phiss{ij}{p}$ on $\MM(\tau_1+1, \tau_2-3)$.

\item The right-hand side $\widehat{F}_{s,ij}^{(p)}$ contains in particular contributions from
\begin{enumerate}
\item the scalarization coupling terms $\widehat{S}(\psi_s^{(p)})_{ij}$ and $(\widehat{Q}\psi_s^{(p)})_{ij}$ appearing in \eqref{eq:ScalarizedTeuSys:general:Kerrperturbation:intro},
\item the scalarization $L_{s,ij}^{(p)}$ of the linear coupling terms $\L_{s}^{(p)}$ of Teukolsky,
\item and the scalarization $N_{W,s,ij}^{(p)}$ of the inhomogeneous terms $\N_{W,s}^{(p)}$ of Teukolsky.
\end{enumerate}

\item The right-hand side $\widehat{F}_{s,ij}^{(p)}$ is chosen such that the scalars $\psiss{ij}{p}$ correspond on $\MM(\tau_2,+\infty)$ to the scalarization of $\pmb\psi_{s}^{(p)}\in\sk_2(\mathbb{C})$, a solution to the following tensorial wave equation in Kerr
\beaa
\squared_2\pmb\psi_{s}^{(p)} -\frac{4ia\cos\th}{|q|^2}\nab_{\pr_\tt}\pmb\psi_{s}^{(p)} -\bigg(\frac{4}{\qs}-\frac{4a^2\cos^2\th}{|q|^6}\Big(|q|^2+6mr\Big)\bigg)\pmb\psi_{s}^{(p)}=0.
\eeaa

\item The extension procedure does not preserve the algebraic constraints \eqref{eq:tensorizationconstraints:intro} in the region $\MM(\tmic,\tau_1+1)\cup \MM(\tau_2-2,\tau_2)$, so that the scalars $\psiss{ij}{p}$ are not obtained from the scalarization of a tensor in $\sk_2(\mathbb{C})$ there.
\end{enumerate}

We refer to the system \eqref{eq:waveeqwidetildepsi1:intro} as the \textit{extended scalarized Teukolsky wave system}.


\subsubsection{Main steps of the proof of Theorem \ref{thm:main:Teu:rough}}


We prove Theorem \ref{thm:main:Teu:rough} in Section \ref{subsect:proofofThm4.1} based on the following main intermediary results:
\begin{enumerate}
\item Theorem \ref{th:main:intermediary}  on high-order unweighted energy-Morawetz estimates for the combined system consisting of the scalarized Teukolsky wave/transport system for $\phiss{ij}{p}$ on $\MM(\tau_1,\tau_2)$ and the extended scalarized Teukolsky wave system \eqref{eq:waveeqwidetildepsi1:intro} for $\psiss{ij}{p}$ on $\MM(\tmic,+\infty)$, conditional on the control of zeroth-order derivatives  of $\phiss{ij}{p}$ and $\psiss{ij}{p}$,

\item Proposition \ref{prop:EnergyMorawetznearinfinitytensorialTeuk} on the high-order weighted energy-Morawetz estimates near infinity for the tensorial Teukolsky wave/transport system \eqref{def:TensorialTeuScalars:wavesystem:Kerrperturbation:intro} \eqref{eq:TensorialTeuSysandlinearterms:rescaleRHScontaine2:general:Kerrperturbation:intro} for $\phis{p}$ on $\MM(\tau_1,\tau_2)$, 

\item  and Theorem \ref{cor:weakMorawetzforTeukolskyfromMillet:bis} on a Morawetz estimate for  solutions to the inhomogeneous tensorial Teukolsky equation \eqref{eq:Teu:intro} on a fixed subextremal Kerr background.
\end{enumerate} 

More precisely, Theorem \ref{th:main:intermediary} combined with Proposition \ref{prop:EnergyMorawetznearinfinitytensorialTeuk} provides high-order weighted energy-Morawetz estimates for $\phis{p}$ in $\MM(\tau_1,\tau_2)$, conditional on the control of zeroth-order derivatives of $\phiss{ij}{p}$ and $\psiss{ij}{p}$. On the other hand,  by moving the quasilinear terms in the Teukolsky wave equation in perturbations of Kerr to the right-hand side such that it is recast into the form of an inhomogeneous Teukolsky equation in a  subextremal Kerr spacetime,   Theorem \ref{cor:weakMorawetzforTeukolskyfromMillet:bis} can be  applied to provide the control of zeroth-order derivatives with a loss of finitely many derivatives. Theorem \ref{thm:main:Teu:rough} then follows directly from combining these two types of estimates which may be regarded as high frequency and low frequency estimates, respectively. 

The proof of Proposition \ref{prop:EnergyMorawetznearinfinitytensorialTeuk}, carried out in Section \ref{sec:proofofprop:EnergyMorawetznearinfinitytensorialTeuk}, is more or less standard, while Theorem \ref{cor:weakMorawetzforTeukolskyfromMillet:bis} is a somewhat straightforward consequence of the spectral estimates in the work of Millet \cite{Millet}, see Section \ref{sec:MorawetzestimatesforTeukolskys=pm2inKerrfromMillet}\footnote{Notice however that the proofs of Proposition \ref{prop:EnergyMorawetznearinfinitytensorialTeuk} and Theorem \ref{cor:weakMorawetzforTeukolskyfromMillet:bis} rely on a refined bound for the squared $L^2$-norm of $\{(r\nab)^{\leq 1} \dk^{\leq \reg}\phis{p}\}_{p=0,1}$; see also observation \ref{observation:improvedestiforangularderi} for a further application of such a refined bound.}. The core of the proof of Theorem \ref{thm:main:Teu:rough} is thus the one of Theorem \ref{th:main:intermediary} which is carried out in Sections \ref{sect:microlocalenergyMorawetztensorialwaveequation}--\ref{sec:proofofth:main:intermediary}.

In the following Sections \ref{subsubsect:EMFscalarizedwave:intro}--\ref{subsubsect:pfofmainthm:intro}, we illustrate the key ideas in establishing Theorem \ref{th:main:intermediary} which requires to deal with the scalarization coupling terms $\widehat{S}(\psi_s^{(p)})_{ij}$ and $(\widehat{Q}\psi_s^{(p)})_{ij}$ appearing in \eqref{eq:ScalarizedTeuSys:general:Kerrperturbation:intro} and the linear coupling terms $\L_{s}^{(p)}$ of Teukolsky, and to extend energy-Morawetz estimates for the extended scalarized Teukolsky wave system \eqref{eq:waveeqwidetildepsi1:intro} to higher order derivatives.


\subsubsection{EMF estimates for inhomogeneous scalarized tensorial wave equations}
\lab{subsubsect:EMFscalarizedwave:intro}


Recall that the right-hand side $\widehat{F}_{s,ij}^{(p)}$ of the extended scalarized Teukolsky wave system \eqref{eq:waveeqwidetildepsi1:intro} contains both the scalarization coupling terms $\widehat{S}(\psi_s^{(p)})_{ij}$ and $(\widehat{Q}\psi_s^{(p)})_{ij}$ appearing in \eqref{eq:ScalarizedTeuSys:general:Kerrperturbation:intro} and the linear coupling terms $\L_{s}^{(p)}$ of Teukolsky. In Section \ref{sect:microlocalenergyMorawetztensorialwaveequation}, we first deal with the contribution coming from $\widehat{S}(\psi)_{ij}$ and $(\widehat{Q}\psi)_{ij}$, i.e., we start with the part of \eqref{eq:ScalarizedTeuSys:general:Kerrperturbation:intro} corresponding to the scalarization of the tensorial wave operator $\squared_2$. To this end, we use the global microlocal EMF estimates proven in our previous work \cite{MaSz24} for the inhomogeneous scalar wave equation as a black box, and the main difficulty is to show that the linear scalarization coupling terms $\widehat{S}(\psi)_{ij}$ and $(\widehat{Q}\psi)_{ij}$ generate lower order terms. This is possible thanks to the following observations:
\begin{enumerate}
\item \textit{Region $\MM(\tau_1+1, \tau_2-2)\cup\MM(\tau_2,+\infty)$.} In this region, we heavily rely on the fact that the scalars $\psi_{ij}$ are derived from a tensor $\pmb\psi\in\sk_2(\mathbb{C})$. The results in \cite{MaSz24} naturally divide $\MM$ into the regions $r\leq R_0$ and $r\geq R_0$ with $R_0\geq 20m$ large enough, and we thus treat these two regions separately:
\begin{enumerate}
\item \textit{Region $r\geq R_0$.} The blackbox EMF estimate from \cite{MaSz24} is purely in physical space in the region $r\geq R_0$. The proof is thus significantly easier in that region since it is away from the trapping and since we can work directly at the level of the tensorial wave equation.
\item \textit{Region $r\leq R_0$.} Given that the EMF estimate from \cite{MaSz24} is microlocal in this region, we must deal with the contribution coming from $\widehat{S}(\psi)_{ij}$ and $(\widehat{Q}\psi)_{ij}$ and show that they generate a good integration by parts structure in order to produce lower order terms. The main contribution comes from the terms 
\beaa
2M_{j}^{l\a}\pr_\a(\psi_{lk}) +2M_{k}^{l\a}\pr_\a(\psi_{jl})
\eeaa
where $M_{i\a}^j:=(\Ddot_\a\Om_i)\c\Om^j$. The main idea is then to decompose $M_{i\a}^j$ into its symmetric part $(M_{S})_{i\a}^{j}$ and its antisymmetric part $(M_{A})_{i\a}^{j}$. $(M_{A})_{i\a}^{j}$ naturally leads to a good integration by parts structure. For $(M_{S})_{i\a}^{j}$, we rely crucially on the formula
\beaa
(M_{S})_{i\a}^{j}=-\frac{1}{2}\pr_\a(x^ix^j)=-\frac{1}{2}x^i\pr_\a(x^j) -\frac{1}{2}x^j\pr_\a(x^i)
\eeaa
which ultimately generates lower order terms using the fact that $x^i\psi_{ij}=x^j\psi_{ij}=0$ in view of  the first two identities of \eqref{eq:tensorizationconstraints:intro} and the fact that the scalars $\psi_{ij}$ are derived from a tensor $\pmb\psi\in\sk_2(\mathbb{C})$ in $\MM(\tau_1+1, \tau_2-2)\cup\MM(\tau_2,+\infty)$.
\end{enumerate}
\item \textit{Region $\MM(\tmic, \tau_1+1)\cup\MM(\tau_2-2,\tau_2)$.} Recall from the extension procedure of Section \ref{subsubsect:Extensiontoglobal:intro} that the scalars $\psi_{ij}$ are not obtained from the scalarization of a tensor in $\sk_2(\mathbb{C})$ in $\MM(\tmic, \tau_1+1)\cup\MM(\tau_2-2,\tau_2)$. On this region, we thus need to estimate what we call the ``tensorization defect," see Definition \ref{def:definitionofthenotationerrforthescalarizationdefect}, measuring the discrepancy of a general set of scalars $\psi_{ij}$ to satisfy the constraints \eqref{eq:tensorizationconstraints:intro} and quantitatively characterized by a set of error terms $\err_{\textrm{TDefect}}[\psi]$. The important observations are that $\err_{\textrm{TDefect}}[\psi]$ satisfies a well-behaved system of wave equations, see Lemma \ref{lemma:waveequationsfortensordeffects}, and is supported in $\MM(\tmic,\tau_1+1)\cup \MM(\tau_2-2,\tau_2)$ where we can rely on local in time energy estimates.
\end{enumerate}


\subsubsection{EMF estimates conditional on the control of low frequencies}
\lab{subsubsect:microlocalEMF:intro}


We now consider the combined system consisting of the tensorial Teukolsky wave/transport system \eqref{def:TensorialTeuScalars:wavesystem:Kerrperturbation:intro} \eqref{eq:TensorialTeuSysandlinearterms:rescaleRHScontaine2:general:Kerrperturbation:intro} on $\MM(\tau_1,\tau_2)$ and the  system of extended scalarized Teukolsky wave equations \eqref{eq:waveeqwidetildepsi1:intro} on  $\MM(\tmic, +\infty)$.

Applying the global microlocal EMF estimates of Section \ref{subsubsect:EMFscalarizedwave:intro} to the system of extended scalarized Teukolsky wave equations \eqref{eq:waveeqwidetildepsi1:intro}, there is yet another kind of coupling terms remaining to be controlled, that is, the linear coupling terms $L_{s,ij}^{(p)}$ (equivalently, $\L_{s}^{(p)}[{\pmb\phi_s}]$) present in both the tensorial Teukolsky wave system \eqref{def:TensorialTeuScalars:wavesystem:Kerrperturbation:intro}  and the extended scalarized Teukolsky wave system \eqref{eq:waveeqwidetildepsi1:intro}. To close the EMF estimates, we rely on the following observations on the structure of the linear coupling terms $\L_{s}^{(p)}[{\pmb\phi_s}]$:
\begin{enumerate}[label=\roman*)]
\item\lab{observation:1} The definitions \eqref{def:TensorialTeuScalars:wavesystem:Kerrperturbation:intro} for $\phis{p}$, $s=\pm 2$, $p=0, 1,2$, which differ slightly from the ones in \cite{Ma}, are such that the derivatives on $\phis{p}$ appearing in each $\L_{s}^{(p)}[\pmb\phi_{s}]$ are in the direction of the vectorfields $\Xcal_s$ which coincide with horizontal vectorfields in Kerr. 

\item\lab{observation:2}  The Teukolsky wave equations \eqref{eq:TensorialTeuSys:Kerrpert:intro} with the linear coupling terms $\L_{s}^{(p)}[\pmb\phi_{s}]$ in \eqref{eq:tensor:Lsn:onlye_2present:general:Kerrperturbation:intro} exhibit a lower-triangular structure up to zeroth-order terms\footnote{That is, the equations for $\phis{p}$ are coupled only with derivatives of $\phis{p'}$ with $0\leq p'<p$.}.

\item The coefficients in front of the terms $\phis{2}$ and $\nab_{\pr_{\tphi}+a\pr_{\tt}}\phis{1}$ on the RHS of the formula for ${\L_{s}^{(2)}[\pmb\phi_{s}]}$ in \eqref{eq:tensor:Lsn:onlye_2present:general:Kerrperturbation} are real-valued functions.
\end{enumerate}
Furthermore, we exploit the following key observations:
\begin{enumerate}[label=\arabic*)]
\item Following the observation made in \cite{Ma}, the EMF estimates for $\{\phis{p}\}_{p=0,1}$, whose spacetime integrands degenerate in the trapping region, can be refined into globally nondegenerate EMF estimates, conditional on the control of $\phis{p+1}$ itself. This nondegeneracy in the Morawetz estimates for $\{\phis{p}\}_{p=0,1}$ in the trapping region is a manifestation of the well-known fact that trapping degeneracy is only present in the highest-order derivatives\footnote{Recall from the Teukolsky transport equations \eqref{def:TensorialTeuScalars:wavesystem:Kerrperturbation:intro} that $\phis{p+1}$ controls a null derivative of $\phis{p}$.}.

\item\lab{observation:improvedestiforangularderi} By rewriting the principal part of the Teukolsky wave equation for $\phis{p}$, $p=0,1$, as $\De_2\phis{p}$ plus a null derivative of $\phis{p+1}$ using the Teukolsky transport equations \eqref{eq:TensorialTeuSysandlinearterms:rescaleRHScontaine2:general:Kerrperturbation:intro}, we achieve a refined bound for the squared $L^2$-norm of $\{\nab \phis{p}\}_{p=0,1}$ (and hence also for the derivatives $\{\nab_{\Xcal_s}\phis{p}\}_{p=0,1}$ and $\{\nab_{\pr_{\phi} + a\pr_{\tau}}\phis{p}\}_{p=0,1}$\footnote{Note that $\pr_{\phi} + a\pr_{\tau}=(\pr_{\phi}+a(\sin\th)^2\pr_\tau)+a(\cos\th)^2\pr_\tau$, where $\nab_{\pr_{\phi}+a(\sin\th)^2\pr_\tau}\phis{p}$ is controlled by $r\nab\phis{p}$ in Kerr and where $\nab_{\pr_\tau}\phis{p}$ satisfies estimates with stronger $r$-weights.} appearing in the formulas of $\L_{s}^{(p)}[{\pmb\phi_s}]$ in \eqref{eq:tensor:Lsn:onlye_2present:general:Kerrperturbation:intro}) in terms of the $L^2$-norm of $\phis{p+1}$ multiplied by the nondegenerate EMF norm of $\phis{p}$. This allows us to absorb the terms involving $\nab_{\Xcal_s}\phis{p}$ and $\nab_{\pr_{\phi} + a\pr_{\tau}}\phis{p}$ by the left-hand side.

\item In the trapping region, we follow \cite{Ma}  to transform the error integral of the form
\beaa
\int \Re\Big(X\phiss{ij}{2} \ov{(\pr_{\tphi}+a\pr_{\tt})\phiss{ij}{1}}\Big)
\eeaa
into integral of products of first-order (pseudodifferential or differential) derivatives. This is realized by using the Teukolsky transport equations \eqref{def:TensorialTeuScalars:wavesystem:Kerrperturbation:+2:intro} \eqref{def:TensorialTeuScalars:wavesystem:Kerrperturbation:-2:intro} for $p=1$, integrating by parts, and noticing that the coefficient in front of the term $\nab_{\pr_{\tphi}+a\pr_{\tt}}\phis{1}$ on the RHS of the formula of ${\L_{s}^{(2)}[\pmb\phi_{s}]}$ in \eqref{eq:tensor:Lsn:onlye_2present:general:Kerrperturbation:intro} is a real-valued function.
\end{enumerate}

The above observations allow us to derive EMF estimates for the combined system 
\eqref{def:TensorialTeuScalars:wavesystem:Kerrperturbation:intro} \eqref{eq:TensorialTeuSysandlinearterms:rescaleRHScontaine2:general:Kerrperturbation:intro} \eqref{eq:waveeqwidetildepsi1:intro} conditional on the control of zeroth-order derivatives of $\phis{p}$ and $\psis{p}$, see Section \ref{sec:proofofth:main:intermediary:0order}.

\begin{remark}
Unlike \cite{Ma} \cite{DHR19} \cite{SRTdC20, SRTdC23} in Kerr, and Part II of \cite{GKS22} in perturbations of Kerr for $|a|\ll m$, we do not rely on transport estimates for the Teukolsky transport equations \eqref{def:TensorialTeuScalars:wavesystem:Kerrperturbation}. See also Remark \ref{remark:ofmainthm7.1}.
\end{remark}


\subsubsection{EMF estimates for higher order derivatives}
\lab{subsubsect:pfofmainthm:intro}


To prove high-order EMF estimates,  we now commute the system consisting of the scalarized Teukolsky wave/transport system\footnote{This scalarized system is equivalent to the tensorial Teukolsky wave/transport system \eqref{def:TensorialTeuScalars:wavesystem:Kerrperturbation:intro} \eqref{eq:TensorialTeuSysandlinearterms:rescaleRHScontaine2:general:Kerrperturbation:intro}.} on $\MM(\tau_1,\tau_2)$ and the system of extended scalarized wave equations \eqref{eq:waveeqwidetildepsi1:intro} on  $\MM(\tmic, +\infty)$, with suitably chosen derivatives. Following our previous work \cite{MaSz24}, it suffices to control the high-order derivatives $(\pr_{\tt}, \chi_0(r)\pr_{\tphi})$ from which the control of arbitrary high-order unweighted derivatives can be recovered, where $\chi_0$ is a cut-off function which equals $1$ for $r\leq 11m$ and vanishes for $r\geq 12m$. While the action of $\pr_{\tt}$ on $\phiss{ij}{p}$ preserves the identities \eqref{eq:tensorizationconstraints:intro},  the action of $\pr_{\tphi}$ on $\phiss{ij}{p}$ does not; hence, the family of scalars $\pr_{\tphi}\phiss{ij}{p}$  does not arise from the scalarization of a tensor in  $\sk_2(\mathbb{C})$. To overcome this difficulty, we introduce an alternative derivative $\widehat{\pr}_{\tphi}$\footnote{In Kerr, $\widehat{\pr}_{\tphi}$ acting on families of scalars corresponds to the horizontal Lie derivative $\Lieb_{\pr_{\tphi}}$ acting on horizontal tensors, see Lemma \ref{lemma:InKerrlinkbetweenprtauwidehatprtphiandLieb[rtauLiebprtphi} for this correspondence, and Section \ref{subsect:horizontalliederivatives} for the definition of horizontal Lie derivatives.}, given by 
\beaa
\widehat{\pr}_{\tphi}(\psi)_{ij}:=\pr_{\tphi}(\psi_{ij}) +\in_{ik3}\psi_{kj} +\in_{jk3}\psi_{ki},
\eeaa
when acting on a family of complex-valued scalars $\psi_{ij}$. It turns out that differentiation with respect to $\widehat{\pr}_{\tphi}$ preserves the identities \eqref{eq:tensorizationconstraints:intro}, see Lemma \ref{lemma:differentiatingwrtwidehatprtphipreservetheidentitiesscaloftensors}, and, hence, the family of scalars $\widehat{\pr}_{\tphi}\phiss{ij}{p}$  indeed correspond to the scalarization of a tensor in  $\sk_2(\mathbb{C})$. 

The derivatives $(\pr_{\tt}, \chi_0(r)\widehat{\pr}_{\tphi})$ are then applied as commutators for the combined scalarized system which finally leads to the high-order unweighted EMF estimates conditional on the control of lower-order derivatives. As in \cite{MaSz24}, these can be further refined to high-order unweighted EMF estimates conditional only on the control of zeroth-order derivatives, see Section \ref{sec:proofofth:main:intermediary}.


\subsection{Overview of the paper}


We review in Section \ref{sec:nonintergrableformalism} the non-integrable formalism of \cite{GKS20} \cite{GKS22}  and define the relevant geometric quantities. The scalarization of tensors and tensorial wave equations using a regular triplet, the notion of tensorization defect and the definition of the $\widehat{\pr}_{\tphi}$ derivative are presented in Section \ref{sec:regularscalarization}.
In Section \ref{sec:TeuinKerr}, we collect and relate the various forms of the Teukolsky equations in Kerr: for scalars in Newman-Penrose formalism, for horizontal tensors, and for scalars obtained within our framework using a regular triplet. Next, we introduce in Section \ref{sect:TeuinKerrperturbation} the prerequisites required  to state our main theorem: the assumptions on the spacetime metric, the null pair and the regular triplet, the tensorial and scalarized Teukolsky wave/transport systems in perturbations of Kerr, and the energy, Morawetz and flux norms. Some useful basic estimates for wave and transport equations in perturbations of Kerr are then provided in Section \ref{sect:basicestimatesforwaveequations}. 

Afterwards, Section \ref{sect:maintheoremandproof} is devoted to constructing an extension of the scalarized Teukolsky wave system  to a semi-global in time problem, stating a precise version of our main theorem, and proving it under the assumption that Theorem \ref{th:main:intermediary},  Proposition \ref{prop:EnergyMorawetznearinfinitytensorialTeuk} and Theorem \ref{cor:weakMorawetzforTeukolskyfromMillet:bis} hold. Theorem \ref{th:main:intermediary} on global-in-time energy-Morawetz estimates for high-order unweighted derivatives of solutions to the scalarized Teukolsky wave/transport system, conditional on the control of zeroth-order derivative terms, is proved in Sections \ref{sect:microlocalenergyMorawetztensorialwaveequation}-\ref{sec:proofofth:main:intermediary}. To this end, we first recall in Section \ref{sect:microlocalenergyMorawetztensorialwaveequation} the energy-Morawetz estimates for solutions to the scalar wave equation proved in \cite{MaSz24}, based on microlocal multipliers adapted to the $r$-foliation of the spacetime, and use them as a black box to derive microlocal energy-Morawetz estimates for the scalarization of tensorial wave equations. These estimates are then applied to the scalarized Teukolsky wave/transport system to complete the proof of Theorem \ref{th:main:intermediary} in Sections \ref{sec:proofofth:main:intermediary:0order} and \ref{sec:proofofth:main:intermediary}. 
Finally,  Proposition \ref{prop:EnergyMorawetznearinfinitytensorialTeuk}, on energy-Morawetz estimates near infinity for high-order weighted derivatives, and Theorem \ref{cor:weakMorawetzforTeukolskyfromMillet:bis}, on a Morawetz estimate for Teukolsky equations in subextremal Kerr, are proved in Sections \ref{sec:proofofprop:EnergyMorawetznearinfinitytensorialTeuk} and \ref{sec:MorawetzestimatesforTeukolskys=pm2inKerrfromMillet}, respectively.


\subsection{Acknowledgments}


The first author is much indebted to Fei Wang for her constant  support,  encouragement and grooviest songs. The second author is supported by the ERC grant ERC-2023 AdG 101141855 BlaHSt.


\section{Non-integrable formalism}
\lab{sec:nonintergrableformalism}


In this section, we briefly review part of the formalism for non-integrable structures introduced in \cite{GKS20} \cite{GKS22}. This will be used to write Teukolsky equations in tensorial form in Section  \ref{sec:TeukolskyinKerrintensorialform}.


\subsection{Null pairs and horizontal structures}
\lab{subsection:review-horiz.structures}


Consider a fixed null pair $e_3, e_4$, i.e., $\g(e_3, e_3)=\g(e_4, e_4)=0$,   $\g(e_3, e_4)=-2,$ and
 denote  by  $\O(\MM)$ the vectorspace  of horizontal vectorfields $X$  on $\MM$, i.e., $\g(e_3, X)= \g(e_4, X)=0$.
  Given a fixed   orientation  on $\MM$,  with corresponding  volume form  $\in$,  we define  the induced 
 volume form on   $\O(\MM)$ by,  
 \beaa
  \in(X, Y):=\frac 1 2\in(X, Y, e_3, e_4). 
  \eeaa
 A null  frame on $\MM$ consists of a choice of horizontal vectorfields  $e_1, e_2$, such that\footnote{We use greek 
 indices $\a, \b, \ga$ for $1,2,3,4$ and latin indices $a,b$ for $1,2$.}
 \beaa
 \g(e_a, e_b)=\de_{ab}\qquad  a, b=1,2.
 \eeaa  
 The commutator $[X,Y]$ of two horizontal vectorfields
may fail however to be horizontal. We say that the pair $(e_3, e_4 )$ is integrable if   $\O(\MM)$  forms an integrable distribution, i.e., $X, Y\in\O(\MM) $ implies that $[X,Y]\in\O(\MM)$. As is well-known,  the  principal null pair in Kerr fails to be integrable.
Given an arbitrary vectorfield $X$, we denote by $^{(h)}X$
its  horizontal projection, 
\beaa
{}^{(h)}X := X+ \frac 1 2 \g(X,e_3)e_4+ \frac 1 2   \g(X,e_4) e_3. 
\eeaa
A  $k$-covariant tensor-field $U$ is said to be horizontal,  $U\in \O_k(\MM)$,
if  for any $X_1,\ldots, X_k$ we have 
\beaa
U(X_1,\ldots, X_k)=U( ^{(h)} X_1,\ldots, {}^{(h)}X_k).
\eeaa

\begin{definition}\label{definition-SS-real}
We denote by $\sk_0=\sk_0(\MM, \mathbb{R})$ the set of pairs of real scalar functions on $\MM$,  by $\sk_1=\sk_1(\MM, \mathbb{R})$ the  set of real horizontal $1$-forms  on $\MM$   and by $\sk_2=\sk_2(\MM, \mathbb{R})$ the set of symmetric traceless   horizontal real $2$-tensors on $\MM$.
\end{definition}

\begin{definition}\label{definition-hodge-duals}
We define the dual of $\xi\in\sk_1$ and $U\in\sk_2$ by
\beaa
\dual \xi_{a}:=\in_{ab}\xi_b,\qquad \dual U_{ab}:=\in_{ac} U_{cb}.
\eeaa
\end{definition}

Note that given $\xi, \eta\in\sk_1$ and $U\in\sk_2$, we have 
\beaa
\dual(\dual \xi)=-\xi, \qquad \dual (\dual U)=-U,\qquad \dual\xi \c  \eta=-\xi\c\dual\eta.
\eeaa
 Also, given  $\xi, \eta\in\sk_1 $,  $U, V\in \sk_2$  we denote
\beaa
\xi\c \eta&:=&\de^{ab} \xi_a\eta_b,\qquad 
\xi\wedge\eta:=\in^{ab} \xi_a\eta_b=\xi\c\dual \eta,\qquad 
(\xi\hot \eta)_{ab}:=\xi_a \eta_b +\xi_b \eta_a-\de_{ab} \xi\c \eta,\\
(\xi\c U)_a&:=&\de^{bc} \xi_b U_{ac}, \qquad  (U\wedge V)_{ab} := \ep^{ab}U_{ac}V_{cb}.
\eeaa

For any $ X, Y\in \O(\MM)$ we define  the induced metric $g(X, Y):=\g(X, Y)$ and the null second fundamental forms
\bea
\chi(X,Y):=\g(\D_Xe_4 ,Y), \qquad \chib(X,Y):=\g(\D_Xe_3,Y).
\eea
Observe that  $\chi$ and $\chib$  are  symmetric if and only if   the horizontal structure is integrable. Indeed this  follows easily from the following formulas
 \beaa
 \chi(X,Y)-\chi(Y,X)&=&\g(\D_X e_4, Y)-\g(\D_Ye_4,X)=-\g(e_4, [X,Y]),\\
 \chib(X,Y)-\chib(Y,X)&=&\g(\D_X e_3, Y)-\g(\D_Ye_3,X)=-\g(e_3, [X,Y]).
\eeaa
  Note  that  we  can view  $\chi$ and $\chib$ as horizontal 2-covariant tensor-fields
 by extending their definition to arbitrary vectorfields  $X, Y$  by setting  $\chi(X, Y)= \chi( ^{(h)}X, ^{(h)}Y)$,  $\chib(X, Y)= \chib( ^{(h)}X, ^{(h)}Y)$.
 Given an horizontal 2-tensor $U$  we define its trace $\tr U$  and anti-trace $\atr U$
\beaa
\tr (U):=\de^{ab}U_{ab}, \qquad \atr U:=\in^{ab} U_{ab}.
\eeaa
Accordingly we  decompose $\chi, \chib$ as follows,
\beaa
\chi_{ab}&=&\chih_{ab} +\frac 1 2 \de_{ab} \trch+\frac 1 2 \in_{ab}\atrch,\\
\chib_{ab}&=&\chibh_{ab} +\frac 1 2 \de_{ab} \trchb+\frac 1 2 \in_{ab}\atrchb,
\eeaa
where $\chih$ and $\chibh$ denote respectively the symmetric traceless part of $\chi$ and $\chib$.

We define the horizontal covariant operator $\nab$ as follows. Given $X, Y\in \O(\MM)$
 \bea
 \nab_X Y&:=&^{(h)}(\D_XY)=\D_XY- \frac 1 2 \chib(X,Y)e_4 -  \frac 1 2 \chi(X,Y) e_3.
 \eea
 In particular, for  all  $X,Y, Z\in \O(\MM)$,
 \beaa
 Z g (X,Y)=g(\nab_Z X, Y)+ g(X, \nab_ZY).
 \eeaa

In the integrable case, $\nab$ coincides with the Levi-Civita connection
 of the metric induced on the integral surfaces of   $\O(\MM)$.  
 Given $X$ horizontal, $\D_4X$ and $\D_3 X$ are in general not horizontal.
 We define $\nab_4 X$ and $\nab_3 X$  to be the horizontal projections
 of the former.  More precisely,
 \beaa
 \nab_4 X&:=&^{(h)}(\D_4 X)=\D_4 X- \frac 1 2 \g(X, \D_4 e_3 ) e_4- \frac 1 2  \g(X, \D_4 e_4)  e_3 ,\\
 \nab_3 X&:=&^{(h)}(\D_3 X)=\D_3 X-   \frac 1 2 \g(X, \D_3e_3) e_3 - \frac 1 2   \g(X, \D_3 e_4 ) e_3. 
 \eeaa
The definition can be easily extended to arbitrary  $  \O_k(\MM) $ tensor-fields  $U$ 
\beaa
 \nab_4U(X_1,\ldots, X_k)&:=&e_4 (U(X_1,\ldots, X_k))- \sum_i U( X_1,\ldots, \nab_4 X_i, \ldots, X_k),\\
  \nab_3 U(X_1,\ldots, X_k)&:=&e_3 (U(X_1,\ldots, X_k)) -\sum_i U( X_1,\ldots, \nab_3 X_i, \ldots, X_k).
 \eeaa


\subsection{Ricci and curvature  coefficients}


Given a null frame $(e_1, e_2, e_3, e_4)$ we define the following connection coefficients,
 \bea
 \begin{split}
\chib_{ab}&:=\g(\D_ae_3, e_b),\qquad \,\,\,\,\,\,\,\chi_{ab}:=\g(\D_ae_4, e_b),\\
\xib_a&:=\frac 1 2 \g(\D_3 e_3 , e_a),\qquad\,\,\,\,\, \xi_a:=\frac 1 2 \g(\D_4 e_4, e_a),\\
\omb&:=\frac 1 4 \g(\D_3e_3 , e_4),\qquad\,\,\,\,\,\,\, \om:=\frac 1 4 \g(\D_4 e_4, e_3),\qquad \\
\etab_a&:=\frac 1 2\g(\D_4 e_3, e_a),\qquad \quad \eta_a:=\frac 1 2 \g(\D_3 e_4, e_a),\qquad\\
 \ze_a&:=\frac 1 2 \g(\D_{e_a}e_4,  e_3),
 \end{split}
\eea
which account for all the  connection coefficients except $\g(\D_{e_\mu} e_b, e_a)$, $\mu=1,2,3,4$, $a, b=1,2$.

We have the Ricci formulas 
\bea
\lab{eq:Ricciformula}
\D_a e_b&=&\nab_a e_b+\frac 1 2 \chi_{ab} e_3+\frac 1 2  \chib_{ab}e_4,\nn\\
\D_a e_4&=&\chi_{ab}e_b -\ze_a e_4,\nn\\
\D_a e_3&=&\chib_{ab} e_b +\ze_ae_3,\nn\\
\D_3 e_a&=&\nab_3 e_a +\eta_a e_3+\xib_a e_4,\nn\\
\D_3 e_3&=& -2\omb e_3+ 2 \xib_b e_b,\label{ricci}\\
\D_3 e_4&=&2\omb e_4+2\eta_b e_b,\nn\\
\D_4 e_a&=&\nab_4 e_a +\etab_a e_4 +\xi_a e_3,\nn\\
\D_4 e_4&=&-2 \om e_4 +2\xi_b e_b,\nn\\
\D_4 e_3&=&2 \om e_3+2\etab_b e_b.\nn
\eea 

For a given horizontal   1-form $\xi$, we  define the frame independent   operators
\bea\lab{eq:defintiondivcurlandnabhot}
\div\xi:=\de^{ab}\nab_b\xi_a,\qquad 
\curl\xi:=\in^{ab}\nab_a\xi_b,\qquad 
(\nab\hot \xi)_{ba}:=\nab_b\xi_a+\nab_a  \xi_b-\de_{ab}( \div \xi).
\eea
We also define the curvature components 
\bea
\bsplit
\a_{ab}&:={\bf R}_{a4b4},\quad \b_a:=\frac 12 {\bf R}_{a434}, \quad \rho:=\frac 1 4 {\bf R}_{3434}, \quad\rhod:=\frac 1 4 \dual{\bf R}_{3434},\\
\bb_a& :=\frac 1 2 {\bf R}_{a334}, \quad \aa_{ab}:={\bf R}_{a3b3},
\end{split}
\eea
where $\dual{\bf R}$ denotes the Hodge dual of the curvature tensor ${\bf R}$.


\subsection{Commutation formulas}


 \begin{lemma}
   \lab{lemma:comm}
   The following commutation formulas hold true:
   \begin{enumerate}
\item Given   $f \in \sk_0$, we have
       \bea\label{eq:comm-nab3-nab4-naba-f-general}
       \begin{split}
        \,[\nab_3, \nab_a] f &=-\frac 1 2 \left(\trchb \nab_a f+\atrchb \dual \nab_a f\right)+(\eta_a-\ze_a) \nab_3 f-\chibh_{ab}\nab_b f  +\xib_a \nab_4 f,\\
         \,[\nab_4, \nab_a] f &=-\frac 1 2 \left(\trch \nab_a f+\atrch \dual \nab_a f\right)+(\etab_a+\ze_a) \nab_4 f-\chih_{ab}\nab_b f  +\xi_a \nab_3 f, \\
         \, [\nab_4, \nab_3] f&= 2(\etab-\eta ) \c \nab f + 2 \om \nab_3 f -2\omb \nab_4 f. 
         \end{split}
       \eea

  \item   Given  $u\in \sk_1$, we have
    \bea\label{commutator-3-a-u-b}\label{commutator-u-in-SS1}
         \bsplit            
\,  [\nab_3,\nab_a] u_b    &=-\frac 1 2 \trchb \big( \nab_a u_b+\eta_b u_a-\de_{ab} \eta \c u \big) -\frac 1 2 \atrchb \big( \dual \nab_a u_b+\eta_b \dual u_a-\in_{ab} \eta\c u\big) \\
&+(\eta-\ze)_a \nab_3 u_b+\err_{3ab}[u],\\
  \err_{3ab}[u] &=-\dual \bb_a\dual u_b+\xib_a\nab_4 u_b-\xib_b \chi_{ac} u_c+\chi_{ab} \,\xib\c u-\chibh_{ac}\nab_c u_b-\eta_b\chibh_{ac}u_c+\chibh_{ab}\eta\c u,
   \end{split}
   \eea
   \bea\label{commutator-4-a-u-b}
   \bsplit
\,  [\nab_4,\nab_a] u_b    &=-\frac 1 2 \trch \big( \nab_a u_b+\etab_b u_a-\de_{ab} \etab \c u \big) \\
&-\frac 1 2 \atrch \big( \dual \nab_a u_b+\etab_b \dual u_a-\in_{ab} \etab\c u\big)+(\etab +\ze)_a \nab_4 u_b +\err_{4ab}[u],\\
   \err_{4ab}[u]&=\dual \b_a\dual u_b+\xi_a\nab_3 u_b-\xi_b \chib_{ac} u_c+\chib_{ab} \,\xi\c u-\chih_{ac}\nab_c u_b-\etab_b\chih_{ac}u_c+\chih_{ab}\etab\c u, 
      \end{split}
   \eea
   \bea
   \bsplit
 \, [\nab_4, \nab_3] u_a&=2 \om \nab_3 u_a -2\omb \nab_4 u_a+ 2(\etab_b-\eta_b ) \nab_b u_a +2(\etab \c u ) \eta_{a} -2 (\eta \c u )\etab_{a}\\
 &  -2 \dual \rho \dual u_a +\err_{43a}[u],\\
 \err_{43a}[u]&= 2 \big( \xib_{a}  \xi_b- \xi_{a}  \xib_b )u^b.
\end{split}
\eea

\item  Given  $u\in \sk_2$, we have 
    \bea\label{commutator-u-in-SS2}\label{commutator-3-a-u-bc}
         \bsplit            
\,  [\nab_3,\nab_a] u_{bc}    &=-\frac 1 2 \trchb\, (\nab_a u_{bc}+\eta_bu_{ac}+\eta_c u_{ab}-\de_{a b}(\eta \c u)_c-\de_{a c}(\eta \c u)_b )\\
&-\frac 1 2 \atrchb\, (\dual \nab_a u_{bc} +\eta_b\dual u_{ac}+\eta_c\dual u_{ab}- \in_{a b}(\eta \c u)_c- \in_{a c}(\eta \c u)_b )\\
&+(\eta_a-\ze_a)\nab_3 u_{bc}+\err_{3abc}[u],\\
\err_{3abc}[u]&= -2\dual \bb_a \dual u_{bc}+\xib_a \nab_4 u_{bc} -\xib_b\chi_{ad}u_{dc} -\xib_c\chi_{ad}u_{bd}+\chi_{ab}\xib_d u_{dc} \\
&+\chi_{ac}\xib_d u_{bd}-\chibh_{ad} \nab_d u_{bc} -\eta_b\chibh_{ad}u_{dc} - \eta_c\chibh_{ad}u_{bd}+\chibh_{ab}\eta_du_{dc} +\chibh_{ac}\eta_du_{bd},
   \end{split}
   \eea
   \bea\label{commutator-4-a-u-bc}
   \bsplit
\,  [\nab_4,\nab_a] u_{bc}    &=-\frac 1 2 \trch\, (\nab_a u_{bc}+\etab_bu_{ac}+\etab_c u_{ab}-\de_{a b}(\etab \c u)_c-\de_{a c}(\etab \c u)_b )\\
&-\frac 1 2 \atrch\, (\dual \nab_a u_{bc} +\etab_b\dual u_{ac}+\etab_c\dual u_{ab}- \in_{a b}(\etab \c u)_c- \in_{a c}(\etab \c u)_b )\\
&+(\etab_a+\ze_a)\nab_4 u_{bc}+\err_{4abc}[u],\\
\err_{4abc}[u]&= 2\dual \b_a \dual u_{bc}+\xi_a \nab_3 u_{bc} -\xi_b\chib_{ad}u_{dc} -\xi_c\chib_{ad}u_{bd}+\chib_{ab}\xi_d u_{dc} +\chib_{ac}\xi_d u_{bd}\\
& -\chih_{ad} \nab_d u_{bc} -\etab_b\chih_{ad}u_{dc} - \etab_c\chih_{ad}u_{bd}+\chih_{ab}\etab_du_{dc} +\chih_{ac}\etab_du_{bd}, 
     \end{split}
   \eea
   \bea\label{commutator-4-3-u-bc}
   \bsplit
   \, [\nab_4, \nab_3] u_{ab} &=2 \om \nab_3 u_{ab} -2\omb \nab_4 u_{ab} + 2(\etab_c-\eta_c ) \nab_c u_{ab} + 4 \eta \hot (\etab \c u)  \\
   &-4 \etab \hot (\eta \c u)-4 \dual \rho \dual u_{ab}+\err_{43ab}[u],\\
\err_{43ab}[u]&= 2 \big( \xib_{a}  \xi_c- \xi_{a}  \xib_c )u^c\,_{b}+2 \big( \xib_{b}  \xi_c- \xi_{b}  \xib_c )u_{a} \,^c.
   \end{split}
\eea
       \end{enumerate}
 \end{lemma}

\begin{proof}
See the proof of Lemma 2.2.8 in \cite{GKS22}.
\end{proof}


\subsection{Complex notations}
\lab{sec:complexnotationsRicciandcurvature}


Recall Definition \ref{definition-SS-real} of the set of real horizontal tensors $\sk_k=\sk_k(\MM, \mathbb{R})$ on $\MM$ for $k=0,1,2$. We now define the corresponding complexified versions.

\begin{definition} 
\lab{def:skC:horizontaltensors}
We denote by $\sk_k(\mathbb{C})$, $k=0,1,2$, the following set of  horizontal tensors on $\MM$: 
\beaa
\bsplit
& a+ i b \in \sk_0(\mathbb{C})\quad\textrm{if}\quad (a, b) \in \sk_0, \qquad F= f+ i \dual f  \in \sk_1(\mathbb{C})\quad\textrm{if}\quad f \in \sk_1, \\ 
& U=u + i \dual u \in \sk_2(\mathbb{C})\quad\textrm{if}\quad u \in \sk_2,
\end{split}
\eeaa
where $F\in\sk_1(\mathbb{C})$ and $U\in\sk_2(\mathbb{C})$ are anti-self dual, i.e., $\dual F=-iF$ and $\dual U=-iU$. 
\end{definition}

\begin{definition}
We define the following complexified curvature components 
\beaa
A:=\a+i\dual\a, \quad B:=\b+i\dual\b, \quad P:=\rho+i\dual\rho,\quad \Bb:=\bb+i\dual\bb, \quad \Ab:=\aa+i\dual\aa,
\eeaa      
with $A, \Ab\in\sk_2(\mathbb{C})$, $B, \Bb\in\sk_1(\mathbb{C})$, $P\in\sk_0(\mathbb{C})$, and the following complexified Ricci coefficients     
\beaa
&& X:=\chi+i\dual\chi, \quad \Xb:=\chib+i\dual\chib, \quad H:=\eta+i\dual \eta, \quad \Hb:=\etab+i\dual \etab,  \\ 
&& Z:=\ze+i\dual\ze, \quad \Xi:=\xi+i\dual\xi, \quad \Xib:=\xib+i\dual\xib,
\eeaa    
with $\widehat{X}, \Xbh\in\sk_2(\mathbb{C})$, $H, \Hb, Z, \Xi, \Xib\in\sk_1(\mathbb{C})$, where $\widehat{X}, \Xbh$, as well as $\tr X, \tr\Xb$ are given by
\beaa
\tr X := \trch-i\atrch, \quad \widehat{X}:=\chih+i\dual\chih, \quad \tr\Xb:=\trchb -i\atrchb, \quad \Xbh:=\chibh+i\dual\chibh.
\eeaa
\end{definition}

\begin{definition}
We define derivatives of complex quantities as follows
\begin{itemize}
\item For two scalar functions $a$ and $b$, we define
\beaa
\DD(a+ib) &:=& (\nabla+i\dual\nabla)(a+ib).
\eeaa

\item For a 1-form $f$, we define
\beaa
\ov{\DD}\c(f+i\dual f) &:=& (\nabla-i\dual\nabla)\c(f+i\dual f)
\eeaa
and  
\beaa
\DD\hot(f+i\dual f) &:=& (\nabla+i\dual\nabla)\hot(f+i\dual f).
\eeaa

\item For a symmetric traceless 2-tensor $u$, we define
\beaa
\ov{\DD}\c(u+i\dual u) &:=& (\nabla-i\dual\nabla)\c(u+i\dual u).
\eeaa
\end{itemize}
\end{definition}


\subsection{The tensorial wave operator}


In order to define the tensorial wave operator in a covariant way, we first introduce the covariant derivative $\Ddot$ acting on mixed tensors of the type $\T_k (\MM)\otimes   \O_l (\MM)$, i.e., tensors  of the form  $U_{\nu_1\ldots \nu_k,  a_1\ldots a_l}$, 
for which we define
\beaa
\Ddot_\mu U_{\nu_1\ldots \nu_k,  a_1\ldots a_l}&:=& e_\mu(U_{\nu_1\ldots \nu_k,  a_1\ldots a_l}) -U_{\D_\mu e_{\nu_1}\ldots \nu_k,  a_1\ldots a_l}-\ldots- U_{\nu_1\ldots  \D_\mu e_{\nu_k},  a_1\ldots a_l}\\
&-& U_{\nu_1\ldots \nu_k,   ^{(h)}(\D_\mu e_{a_1})\ldots a_l}-  U_{\nu_1\ldots \nu_k,   a_1 \ldots ^{(h)}(\D_\mu e_{a_l})}.
\eeaa

\begin{proposition}
\lab{Proposition:commutehorizderivatives}
For a tensor $\Psi\in \O_1 (\MM)$, we   have  the following formula
 \bea
( \Ddot _\mu\Ddot_\nu -\Ddot_\nu\Ddot _\mu)\Psi_a=\Rdot_{a b  \mu\nu}\Psi^b
 \eea
with an immediate generalization to tensors $\Psi\in \O_l (\MM)$, where, with $(\La_\a)_{\b\ga}= \g(\D_\a e_\ga, e_\b)$,
 \bea
 \lab{eq:DefineRdot}
 \bsplit
 \Rdot_{ab   \mu\nu}&:= {\bf R}_{ab    \mu\nu}+ \frac 1 2  \B_{ab   \mu\nu},\\
  \B_{ab   \mu\nu} &:=  (\La_\mu)_{3a} (\La_\nu)_{b4}+  (\La_\mu)_{4a} (\La_\nu)_{b3}- (\La_\nu)_{3a} (\La_\mu)_{b4}-  (\La_\nu)_{4a} (\La_\mu)_{b3}.
  \end{split}
 \eea 
 \end{proposition}

\begin{proof}
See Proposition 2.1.27 in \cite{GKS22}.
\end{proof}

\begin{proposition}
\lab{proposition:componentsofB}
The components of $\B$   are given   by the following formulas:
\bea
\begin{split}
\B_{ a   b  c 3}&=     -  \trchb  \big( \de_{ca}\eta_b-  \de_{cb} \eta_a\big)  -  \atrchb \big( \in_{ca}  \eta_b -  \in_{cb}  \eta_a\big) \\
&+ 2 \big(- \chibh_{ca}  \eta_b + \chibh_{cb} \eta_a-  \chi_{ca} \xib_b+  \chi_{cb} \xib_a\big),\\
\B_{ a   b  c 4}&=     -  \trch  \big( \de_{ca}\etab_b-  \de_{cb} \etab_a\big)  -  \atrch \big( \in_{ca}  \etab_b -  \in_{cb}  \etab_a\big) \\
&+ 2 \big(- \chih_{ca}  \etab_b + \chih_{cb} \etab_a-  \chib_{ca} \xi_b+  \chib_{cb} \xi_a\big),\\
\B_{ a   b  3 4} &=4\big(-\eta_a \etab_b+\etab_a\eta_b -\xib_a \xi_b+\xi_a \xib_b\big),\\
\B_{abcd} &= \left(- \frac 12  \trch \trchb-\frac 1 2 \atrch \atrchb+\chih \c \chibh\right)\in_{ab}\in_{cd}.
\end{split}
\eea
\end{proposition}

\begin{proof}
See Proposition 2.2.4 in \cite{GKS22}.
\end{proof}

Then, we  define  the wave operator for $\psi \in \sk_k(\mathbb{C})$, $k=0,1,2$, to be, see Definition 2.3.1 in \cite{GKS22},  
 \bea\label{eq:def=squared-2}
 \squared_k\psi:= \g^{\mu\nu} \Ddot_\mu\Ddot_ \nu \psi.
\eea

The following lemma provides the decomposition of $\squared_k$ in null frames. 
\begin{lemma}\label{lemma:expression-wave-operator}
The wave operator for $\psi\in {\sk_k(\mathbb{C})}$, $k=0,1,2$, is given by
\bea\label{eq:wave-squared}
\begin{split}
\squared_k \psi&=-\nab_4 \nab_3 \psi  -\frac 1 2 \trchb \nab_4\psi+\left(2\om -\frac 1 2 \trch\right) \nab_3\psi+\lap_k\psi+2\etab \c\nab \psi \\
&+ ki \left( \rhod- \eta \wedge \etab \right) \psi+(\Ga_b \c \Ga_g) \c \psi,\\
\squared_k \psi&=-\nab_3\nab_4\psi +\left(2\omb -\frac 1 2 \trchb\right)\nab_4\psi -\frac 1 2 \trch \nab_3\psi+\lap_k\psi+2\eta \c\nab \psi \\
&- ki \left( \rhod- \eta \wedge \etab \right) \psi+(\Ga_b \c \Ga_g) \c \psi,
\end{split}
\eea
where $\lap_k=\nab^a \nab_a$ denotes the horizontal Laplacian for $k$-tensors. 
\end{lemma}

\begin{proof}
See Lemma 4.7.5 in \cite{GKS22} for the first identity of \eqref{eq:wave-squared} in the case $k=2$. The proof of Lemma 4.7.5 in \cite{GKS22} immediately extends to $k=0,1,2$ and to the second identity of \eqref{eq:wave-squared}.
\end{proof}


\subsection{Horizontal Lie derivatives}
\label{subsect:horizontalliederivatives}


Recall that the Lie derivative of a $k$-covariant tensor $U$ relative to a vectorfield  $X$ is given by
\beaa
\Lie_X{U}\big(e_{\a_1}, \ldots , e_{\a_k}\big) = X\big(U_{\a_1\ldots\a_k}\big) -  U\big(\Lie_Xe_{\a_1}, \ldots,e_{\a_k}\big) - U\big(e_{\a_1}, \ldots, \Lie_Xe_{\a_k}\big),
\eeaa
where $\Lie_X Y=[X, Y]$. We define horizontal Lie derivatives as follows, see Definition 2.2.12 in \cite{GKS22}.

\begin{definition}[Horizontal Lie derivatives]\label{definition:hor-Lie-derivative}
Given   vectorfields  $X$, $Y$,  the horizontal Lie  derivative  $\Lieb_XY$ is given by
 \beaa
 \Lieb_X Y :=\Lie_X Y+ \frac 1 2 \g(\Lie_XY, e_3) e_4+  \frac 1 2 \g(\Lie_XY, e_4) e_3.
 \eeaa
 Given  a horizontal covariant k-tensor $U$,  the horizontal  Lie derivative $\Lieb_X U $ is defined   to be the projection of $\Lie_X U$ to the  horizontal space, i.e.,
  \beaa
 \Lieb_X U\big(e_{a_1}, \ldots,e_{a_k}\big) := X\big(U_{a_1\ldots a_k}\big)- U\big(\Lieb_Xe_{a_1},\ldots, e_{a_k} \big)-\ldots -U\big(e_{a_1},\ldots,  \Lieb_Xe_{a_k} \big).
  \eeaa 

Also, given a mixed tensor $U$ of the type $\T_k (\MM)\otimes   \O_l (\MM)$, we define the general horizontal derivative $\Lied_XU$ as follows 
\beaa
&& \Lied_XU\big(e_{\a_1},\ldots, e_{\a_k},  e_{a_1},\ldots, e_{a_l}\big) \\
&:=& X\big(U_{\a_1\ldots \a_k,  a_1\ldots a_l}\big) -U\big(\Lie_Xe_{\a_1},\ldots, e_{\a_k},  e_{a_1}, \ldots, e_{a_l}\big) -\ldots  - U\big(e_{\a_1},\ldots.  \Lie_Xe_{\a_k},  e_{a_1},\ldots, e_{a_l}\big)\\
&& -U\big(e_{\a_1}, \ldots, e_{\a_k},  \Lieb_Xe_{a_1},\ldots, e_{a_l}\big) - \ldots -  U\big(e_{\a_1},\ldots, e_{\a_k},   e_{a_1}, \ldots, \Lieb_X e_{a_l}\big).
\eeaa
\end{definition}


\subsection{Kerr values}
\label{subsect:Kerrspacetime}



\subsubsection{Normalized coordinates in Kerr spacetimes}
\label{subsect:normalizedcoords}


The Kerr metric in Boyer--Lindquist coordinates $(t,r,\th,\phi)$ is given by
\begin{align}
\gam={}& \g_{tt}dt^2 +\g_{rr}dr^2+(\g_{t\phi}+\g_{\phi t})dtd\phi +\g_{\phi\phi}d\phi^2 +\g_{\th\th}d\th^2,
\end{align}
where
\begin{equation}
\begin{split}
\g_{tt}={}&-\frac{\Delta-a^2\sin^2\theta}{|q|^2}, \quad \g_{t\phi}={}\g_{\phi t}=-\frac{2amr\sin\theta}{|q|^2}, \quad \g_{rr}={}\frac{|q|^2}{\Delta},\\
\g_{\phi\phi}={}&\frac{(\R)^2-a^2\sin^2\theta\Delta}{|q|^2}\sin^2\theta, \quad \g_{\th\th}={}|q|^2,
\end{split}
\end{equation}
with 
\bea
\Delta:=r^2-2mr+a^2,\qquad |q|^2:=r^2+a^2\cos^2\th.
\eea
{In particular, $\partial_{t}$ and $\partial_{\phi}$ are Killing vectorfields and $|\det(\gam)|=|q|^4\sin^2\theta$.} The larger root 
\begin{align}
r_+:=m + \sqrt{m^2 -a^2}
\end{align}
of $\Delta=\De(r)$ corresponds to the location of the event horizon. For convenience, we define 
\bea
\mu :=\frac{\Delta}{\R}.
\eea 

The nontrivial components of the inverse metric are
\begin{equation}
\begin{split}
\label{Kerrmetric:inverse:BL}
\g^{tt}={}&-\frac{(\R)^2-a^2\sin^2\theta\Delta}{|q|^2\Delta},  \quad \g^{rr}={}\frac{\Delta}{|q|^2},\\
\g^{\phi\phi}={}&\frac{\Delta-a^2\sin^2\theta}{|q|^2\Delta\sin^2\th}, \quad \g^{\th\th}={}\frac{1}{|q|^2}, \quad \g^{t\phi}={}\g^{\phi t}=-\frac{2amr}{|q|^2\Delta}.
\end{split}
\end{equation}

We define as well a tortoise coordinate $r^*$ by 
$$dr^*= \mu^{-1}dr, \qquad r^*(3m)=0.$$
Without confusion, we call $(t,r^*,\th,\phi)$ the tortoise coordinates and we denote $\pr_{r^*}$ as the coordinate derivative in this tortoise coordinate system.

It is well-known that the metric is singular on the event horizon in both the Boyer--Lindquist and the tortoise coordinates. To extend the Kerr metric beyond the future event horizon, we define the ingoing Eddington--Finkelstein coordinates $(v_+, r,\th,\phi_+)$ by 
\bea\lab{eq:definitionofingoingEFcoordiantesvplusandphiplus}
dv_+=dt+\mu^{-1}dr, \quad d\phi_+=d\phi+\frac{a}{\Delta}dr \,\,\, \text{mod } 2\pi.
\eea
The Kerr metric in this coordinate system is 
\begin{align}\lab{eq:KerrmetriciningoingEddigtonFinkelstein}
\gam={}&-\bigg(1-\frac{2mr}{|q|^2}\bigg)dv_+^2 +2dr dv_+ -\frac{4amr\sin^2\th}{|q|^2}dv_+d\phi_+ -2a\sin^2\th dr d\phi_+\nn\\
&+|q|^2 d\th^2+\frac{(\R)^2-a^2\sin^2\theta\Delta}{|q|^2}\sin^2\th d\phi_+^2.
\end{align}

In the following lemma, we introduce coordinates systems, referred to as \textit{normalized coordinates}, and used in particular in the statement of the main result of this paper. 

\begin{lemma}[Normalized coordinates]
\label{lem:specificchoice:normalizedcoord}
We fix constants $\dhor$ and $\dbl$ such that 
$$0<\dhor\ll \dbl\ll 1-\frac{|a|}{m}.$$ 
There exists a choice of smooth functions $\tmod=\tmod(r)$ and $\phimod=\phimod(r)$ such that the coordinate systems $(\tt, r, x^1_0, x^2_0)$ and $(\tt, r, x^1_p, x^2_p)$, defined respectively on $\th\neq 0, \pi$ and $\th\neq\frac{\pi}{2}$, with 
\bea
\tau=v_+-\tmod, \quad \tphi=\phi_+ -  \phimod, \quad x^1_0=\th, \quad x^2_0=\tphi, \quad x^1_p=\sin\th\cos\tphi, \quad x^2_p=\sin\th\sin\tphi,
\eea
satisfy the following properties:
\begin{enumerate}
\item defining the causal spacetime region $\MM$ and corresponding spacelike boundary $\AA$ by 
\beaa
\bsplit
{}\qquad\MM&:=\big(\{(\tt, r, x^1_0, x^2_0),\,\, \th\neq 0, \pi\}\cup\{(\tt, r, x^1_p, x^2_p),\,\, \th\neq \pi/2\}\big)\cap\{r\geq r_+(1-\dhor)\}, \\ 
{}\qquad\AA&:=\pr\MM=\big(\{(\tt, r, x^1_0, x^2_0),\,\, \th\neq 0, \pi\}\cup\{(\tt, r, x^1_p, x^2_p),\,\, \th\neq \pi/2\}\big)\cap\{r=r_+(1-\dhor)\},
\end{split}
\eeaa
$\MM$ is covered by $(\tt, r, x^1_0, x^2_0)$ and $(\tt, r, x^1_p, x^2_p)$ with the metric components and  inverse metric components being smooth on their respective coordinate patch, 

\item $(\tt, r, x^1_0, x^2_0)$ coincides with Boyer-Lindquist coordinates\footnote{In particular, we have
\beaa
\tmod'(r)=\mu^{-1}, \qquad \phimod'(r)=\frac{a}{\De}\quad\textrm{on}\quad r\in [r_+(1+2\dbl), 12m].
\eeaa}
 in $r\in [r_+(1+2\dbl),  12m]$, 

\item for $r\notin (r_+(1+\dbl), 13m)$, we choose 
\beaa
\begin{split}
\tmod'(r) &=\frac{m^2}{r^2}, \qquad \phimod'(r)=0\quad\textrm{on}\quad r\leq r_+(1+\dbl),\\
\tmod'(r)&=2\mu^{-1}-\frac{m^2}{r^2}, \qquad \phimod'(r)=\frac{2a}{\De} \quad\textrm{on}\quad r\geq 13m,
\end{split}
\eeaa

\item the level sets of $\tt$ in $\MM$ are globally spacelike, transverse to the future event horizon $\HH_+$ and the spacelike boundary $\AA$, and asymptotically null to future null infinity $\II_+$.
\end{enumerate}

Furthermore, the nontrivial inverse metric components in the coordinate system $(\tt, r, \th, \tphi)$ are
\begin{align}
\label{eq:inverse:hypercoord}
\gam^{\tt \tt}={}&\frac{a^2\sin^2\th}{|q|^2} -\frac{2(\R)}{|q|^2}\tmod'+\frac{\Delta}{|q|^2}(\tmod')^2, \,\, \gam^{rr}=\frac{\Delta}{|q|^2}, 
\,\, \gam^{\tt r}=\gam^{r \tt}=\frac{\R}{|q|^2}(1-\mu \tmod'), \nn\\
\gam^{r\tphi} ={}&\gam^{\tphi r} =\frac{a}{|q|^2}-\frac{\Delta}{|q|^2}\phimod', \quad \gam^{\tt\tphi}= \gam^{\tphi \tt}= \frac{a}{|q|^2} (1-\tmod')-\phimod'\frac{\R}{|q|^2} (1-\mu \tmod') ,\nn\\
\gam^{\th\th}={}&\frac{1}{|q|^2}, \quad \gam^{\tphi\tphi}=\frac{1}{|q|^2\sin^2\th}-\frac{2a}{|q|^2}\phimod' +\frac{\Delta}{|q|^2}(\phimod')^2,
\end{align}
and the volume form verifies, with $(x^1, x^2)$ denoting either $(x^1_0, x^2_0)$ or $(x^1_p, x^2_p)$,
\bea\lab{eq:assymptiticpropmetricKerrintaurxacoord:volumeform}
\sqrt{|\det(\gam)|}d\tau dr dx^1dx^2 = |q|^2\sqrt{\det(\mathring{\ga})}d\tau dr dx^1dx^2,
\eea
where $\mathring{\ga}$ denotes the metric on the standard unit 2-sphere.

Finally, for $r\geq 13m$, the inverse metric and metric components  satisfy the following asymptotics on their respective coordinate patch, with $(x^1, x^2)$ denoting either $(x^1_0, x^2_0)$ or $(x^1_p, x^2_p)$,
\bea\lab{eq:assymptiticpropmetricKerrintaurxacoord:1}
\bsplit
\gam^{rr}=&1+O(mr^{-1}), \qquad \gam^{r \tt}=-1+O(m^2r^{-2}),\qquad \gam^{rx^a}=O(mr^{-2}), \\
\gam^{\tt \tt}=&O(m^2r^{-2}), \qquad \gam^{\tau x^a}=O(mr^{-2}),\qquad \gam^{x^ax^b}=\frac{1}{r^2}\mathring{\ga}^{x^ax^b}+O(m^2r^{-4})
\end{split}
\eea
and
\bea\lab{eq:assymptiticpropmetricKerrintaurxacoord:2}
\bsplit
(\gam)_{rr}=&O(m^2r^{-2}), \qquad (\gam)_{r \tt}=-1+O(m^2r^{-2}),\qquad (\gam)_{rx^a}=O(m),\\
(\gam)_{\tt \tt} =& -1+O(mr^{-1}),\qquad (\gam)_{\tau x^a}=O(m),\qquad  (\gam)_{x^ax^b}=r^2\mathring{\ga}_{x^ax^b}+O(m^2).
 \end{split}
\eea
\end{lemma}

\begin{remark}\lab{rmk:phimodprimeisproportionaltoa!!}
Additionally, we may choose $\phimod$ such that 
\beaa
\phimod'(r)=a\phi_{\textrm{mod},0}'(r), \qquad \phi_{\textrm{mod},0}'(r)\geq 0\quad \forall r\in(r_+(1-\dhor), +\infty),
\eeaa
so that $\phimod'(r)$ has the same sign as $a$. From now on, we will assume that our choice of $\phimod$ satisfies this property. In view of \eqref{eq:inverse:hypercoord}, it implies that the inverse metric coefficients $\gam^{\a\b}$ in the normalized coordinates system $(\tau, r, x^1, x^2)$ are invariant under the change $(a, \tphi)\to (-a, -\tphi)$. 
\end{remark}

\begin{proof}
See the proof of Lemma 2.1 in \cite{MaSz24}.
\end{proof} 
  
Next, we consider the asymptotics of the induced metric on the level sets of $\tau$ in normalized coordinates in the region $r\geq 13m$.
 \begin{lemma}\lab{lemma:specificchoice:normalizedcoord:inducedmetricSitau}
Let $g_{a,m}$ denote the metric induced by $\gam$ on the level sets of $\tau$. Then, we have in the normalized coordinate systems $(r, x^1_0, x^2_0)$ and $(r, x^1_p, x^2_p)$, in $r\geq 13m$, 
\beaa
\bsplit
(g_{a,m})_{rr}&=O(m^2r^{-2}), \qquad (g_{a,m})_{rx^a}=O(m), \qquad (g_{a,m})_{x^ax^b}=O(r^2), \\
g_{a,m}^{rr}&=O(m^{-2}r^2), \qquad g_{a,m}^{rx^a}=O(m^{-1}), \qquad g_{a,m}^{x^ax^b}=O(r^{-2}),
\end{split} 
\eeaa
and
\beaa
\sqrt{\det(g_{a,m})}drdx^1dx^2 &= r\sqrt{2m^2 -a^2\sin^2\th +O(m^3r^{-1})}\sqrt{\det(\mathring{\ga})}drdx^1dx^2, 
\eeaa
with $(x^1, x^2)$ denoting either $(x^1_0, x^2_0)$ or $(x^1_p, x^2_p)$.
\end{lemma}
  
\begin{proof}
See the proof of Lemma 2.3 in \cite{MaSz24}.
\end{proof}
  
Next, we consider the induced metric on $\AA$.
\begin{lemma}\lab{lemma:specificchoice:normalizedcoord:inducedmetricAA}
Let $(g_{\AA})_{a,m}$ denote the metric induced by $\gam$ on $\AA$. Then,
\beaa
\sqrt{\det((g_{\AA})_{a,m})}d\tau dx^1 dx^2\simeq m\sqrt{\dhor} \sqrt{\det(\mathring{\ga})}d\tau dx^1 dx^2
\eeaa
with $(x^1, x^2)$ denoting either $(x^1_0, x^2_0)$ or $(x^1_p, x^2_p)$.
\end{lemma}

\begin{proof}
See the proof of Lemma 2.4 in \cite{MaSz24}.
\end{proof}


\subsubsection{Principal null pair in Kerr}
\label{subsect:principalnullpairinKerr}


We consider the principal null pair of Kerr which is regular across the future even horizon, i.e., in Boyer-Lindquist coordinates, 
\bea
\lab{def:e3e4inKerr}
 e_4 = \frac{r^2+a^2}{|q|^2} \pr_t +\frac{\De}{|q|^2} \pr_r +\frac{a}{|q|^2} \pr_{\phi}, \qquad 
 e_3=\frac{r^2+a^2}{\De} \pr_t -\pr_r +\frac{a}{\De} \pr_{\phi}.
\eea
Also, we consider its associated horizontal bundle $\{e_3, e_4\}^\perp$, which, for $\th\neq 0, \pi$, is spanned by 
\bea
\lab{def:e1e2inKerr}
 e_1=\frac{1}{|q|}\pr_\th,\quad e_2=\frac{a\sin\th}{|q|}\pr_t+\frac{1}{|q|\sin\th}\pr_\phi,
\eea
and we define the complex-valued scalar $q$ and the complex horizontal  $1$-forms $\Jk$ and $\Jk_\pm$ as 
\bea\lab{eq:def:Jkandq}
\bsplit
q=& r+ia\cos\th, \qquad \Jk=j+i\dual j, \qquad \Jk_\pm=j_\pm+i\dual j_\pm, \qquad j_1=0, \quad j_2=\frac{\sin\th}{|q|}, \\
(j_+)_1=&\frac{1}{|q|} \cos\th\cos\phi_+, \,\,\,\, (j_+)_2=-\frac{1}{|q|} \sin\phi_+, \,\,\,\, (j_-)_1 =\frac{1}{|q|} \cos\th\sin\phi_+, \,\,\,\, (j_-)_2=\frac{1}{|q|}  \cos\phi_+,
\end{split}
 \eea
where $\Jk$ and $\Jk_\pm$ are regular (even at the axis) as well as anti-self dual, i.e., $\Jk,\, \Jk_\pm\in\sk_1(\mathbb{C})$, and where the coordinate $\phi_+$ involved in the definition of $j_\pm$ has been introduced in \eqref{eq:definitionofingoingEFcoordiantesvplusandphiplus}. Note in particular the following identities 
\bea\lab{eq:usefulalgebraicidentitiesinvolvingscalarproductsReJkReJkpm}
\bsplit
\Re(\Jk)\c\Re(\Jk) &=\frac{(\sin\th)^2}{|q|^2}, \qquad \Re(\Jk)\c\Re(\Jk_+)=-\frac{x^2}{|q|^2},\qquad \Re(\Jk)\c\Re(\Jk_-)=\frac{x^1}{|q|^2},\\
\dual(\Re(\Jk))\c\Re(\Jk_+) &=\frac{\cos\th x^1}{|q|^2},\qquad \dual(\Re(\Jk))\c\Re(\Jk_-)=\frac{\cos\th x^2}{|q|^2}.
\end{split}
\eea

The complexified Ricci coefficients w.r.t. this principal null pair are given by  
\bea\lab{eq:KerrvaluesofcomplexifiedRicci}
\begin{split}
&\Xh=\Xbh=\Xi=\Xib=\omb=0, \qquad  \tr X=\frac{2\De \ov{q}}{|q|^4}, \qquad \tr\Xb=-\frac{2}{\ov{q}}, \qquad \om=- \frac 12\pr_r\left(\frac{\Delta}{|q|^2}\right),\\
&H=Z=\frac{a}{\ov{q}}\Jk=\frac{aq}{|q|^2}\Jk, \qquad \Hb=-\frac{a}{q}\Jk= -\frac{a\ov{q}}{|q|^2}\Jk.
\end{split}
\eea

The complexified curvature components are given by
\bea
A=B=\Bb=\Ab=0,\quad P=-\frac{2m}{q^3}.
\eea

Also, the principal null frame acts on the the normalized coordinates of Lemma \ref{lem:specificchoice:normalizedcoord} as follows 
\bea\lab{eq:actionofingoingprincipalnullframeonnormalizedcoordinates}
\bsplit
e_3(r)&=-1, \qquad\quad e_4(r)=\frac{\De}{|q|^2}, \qquad\qquad\qquad\qquad\quad\, e_1(r)=0, \quad\,\,\,\,\, e_2(r)=0,\\
e_3(\tau)&=\tmod'(r), \quad e_4(\tau)=\frac{2(r^2+a^2) - \De\tmod'(r)}{|q|^2}, \quad e_1(\tau)=0, \quad\,\,\,\, e_2(\tau)=\frac{a\sin\th}{|q|},\\
e_3(\th)&=0, \qquad\quad\,\,\,\,\, e_4(\th)=0, \qquad\qquad\qquad\qquad\qquad\,\, e_1(\th)=\frac{1}{|q|}, \quad e_2(\th)=0,\\
e_3(\tphi)&=\phimod'(r), \quad e_4(\tphi)=\frac{2a -\De\phimod'(r)}{|q|^2},\qquad\quad\,\,\,\,\, e_1(\tphi)=0, \quad\,\,\,\, e_2(\tphi)=\frac{1}{|q|\sin\th}.
\end{split}
\eea
In particular, recalling that $x^1_p=\sin\th\cos\tphi$ and $x^2_p=\sin\th\sin\tphi$, we have
\bea\lab{eq:actionofingoingprincipalnullframeonnormalizedcoordinates:bis}
\bsplit
e_3(x^1_p)&=-\phimod'(r)x^2_p, \qquad\, e_4(x^1_p)=-\frac{2a -\De\phimod'(r)}{|q|^2}x^2_p,\\ 
e_3(x^2_p)&=\phimod'(r)x^1_p, \qquad\quad e_4(x^2_p)=\frac{2a -\De\phimod'(r)}{|q|^2}x^1_p,
\end{split}
\eea
and, in view of the definition of $\Jk$ and $\Jk_\pm$,  
\bea\lab{eq:actionofingoingprincipalnullframeonnormalizedcoordinates:ter}
\DD(\tau)=a\Jk, \qquad \DD(\cos\th)=i\Jk,\qquad \DD(x^1_p)=\Jk_+, \qquad \DD(x^2_p)=\Jk_-.
\eea

Moreover the derivatives of $\Jk$ and $\Jk_\pm$ w.r.t. the principal null frame satisfy the following
\bea
\bsplit
\nab_3\Jk &=\frac{1}{\ov{q}}\Jk, \qquad \nab_4\Jk =- \frac{\De \ov{q}}{|q|^4}\Jk, \qquad \nab_3\Jk_\pm =\frac{1}{\ov{q}}\Jk_\pm, \qquad \nab_4 \Jk_\pm =- \frac{\De \ov{q}}{|q|^4}\Jk_{\pm} \mp  \frac{2a}{|q|^2}\Jk_{\mp},\\
\ov{\DD}\c\Jk &=\frac{4i(r^2+a^2)\cos\th}{|q|^4},\qquad \DD\hot\Jk=0, \\ 
\ov{\DD}\c \Jk_+ &= - \frac{4r^2 }{|q|^4}x^1_p - \frac{4ia^2\cos\th}{|q|^4}x^2_p, \qquad \ov{\DD}\c \Jk_- = - \frac{4r^2 }{|q|^4}x^2_p + \frac{4ia^2\cos\th}{|q|^4}x^1_p, \qquad \DD\hot \Jk_{\pm} =0.
\end{split}
\eea


\section{Regular scalarization of tensorial wave equations}
\lab{sec:regularscalarization}


The Teukolsky equation involves the tensorial wave operator $\squared_2$ on $\sk_2(\mathbb{C})$, see Section  \ref{sec:TeukolskyinKerrintensorialform}. In this paper, we will need to: 
\begin{itemize}
\item derive microlocal energy-Morawetz estimates for tensorial wave equations, see Sections \ref{sect:microlocalenergyMorawetztensorialwaveequation} and \ref{sec:proofofth:main:intermediary},

\item extend tensorial wave equations on $\tau\in(\tau_1, \tau_2)$ to $\tau\in(\tmic, +\infty)$, see Section  \ref{sect:extensiontoglobalproblem:Teu}, 

\item estimate the difference between $\squared_2$ expressed respectively in perturbations of Kerr and in Kerr in order to apply the energy-Morawetz estimates in Kerr of Section \ref{sec:energyMorwetzestimatesinKerr}.
\end{itemize}
While the above could likely be done directly at the level of tensors, it will be nevertheless easier to work with a scalarized version of the tensorial wave operator $\squared_2$ on $\sk_2(\mathbb{C})$. A well-known procedure is to use Newman-Penrose formalism \cite{NP62, NP63errata}, see Section  \ref{sec:Teukoslkytransportwavesystemwithcomplexscalars}, where the scalarization is performed using the horizontal vectorfields $(e_1, e_2)$. However, this comes at the expense of generating scalars which are irregular\footnote{An alternative is to use the GHP formalism, but this would again generate tensorial equations.} at the axis of symmetry, i.e., at $\th=0,\pi$. Instead, we propose here an alternative which generates regular scalars at the expense of generating more scalar wave equations.


\subsection{Scalarization using a regular triplet $\Om_i$, $i=1,2,3$}
\lab{subsect:introduceOm_i}


In order to scalarize horizontal tensors, we will rely on the following definition. 

\begin{definition}[Regular triplet]
\lab{def:definitionofregulartripletOmii=123}
Let $(\MM, \g)$ be a spacetime, $(e_3, e_4)$ be a null pair, and consider the corresponding horizontal structure $\O(\MM)$ introduced in Section  \ref{subsection:review-horiz.structures}. We say that vectorfields $\Om_i$, $i=1,2,3$, identified with elements of $\sk_1$, form a regular triplet if they are regular and satisfy the following identities 
\bea\lab{eq:fundamentalpropertiesof1formsOmi}
x^i\Om_i=0, \qquad (\Om^i)_a(\Om_i)_b=\de_{ab}, \qquad \Om_i\c\Om_j=\de_{ij}-x^ix^j, \qquad \Om_i\c\dual\Om_j=\in_{ijk}x^k,
\eea
where, by convention, we denote $\Om^i=\Om_i$, i.e., the $i$-index is lowered or raised using $\de_{ij}$ or $\de^{ij}$.
\end{definition}

\begin{remark}\lab{rmk:generalcontructionofregulartripletsingivenspacetime}
Given a spacetime $(\MM, \g)$, a null pair $(e_3, e_4)$ and the corresponding horizontal structure $\O(\MM)$, one can easily generate regular triplets by enforcing the identities \eqref{eq:fundamentalpropertiesof1formsOmi} on one given topological sphere in $\MM$, which can then be propagated to $\MM$ by defining $\Om_i$ based on well-chosen transport equations consistent with the horizontal structure of $\MM$. For a specific choice of a regular triplet in Kerr, see Section  \ref{sec:regulartripletsinKerr}.
\end{remark}

Next, we introduce the following 1-forms on $\MM$.
\begin{definition}
\lab{def:Mialphaj:Kerr}
Let $(\MM, \g)$ be a spacetime, $(e_3, e_4)$ be a null pair, and consider the corresponding horizontal structure $\O(\MM)$ introduced in Section  \ref{subsection:review-horiz.structures}. Let $\Om_i$, $1,2,3$ be a regular triplet in the sense of Definition \ref{def:definitionofregulartripletOmii=123}. We  define the following 1-forms on $\MM$ 
\bea\lab{eq:definitionofMalphaijwithoutambiguity}
M_{i\a}^j:=(\Ddot_\a\Om_i)\c\Om^j, \quad \forall \a,i,j.
\eea
Further, we define $M_{i}^{j \a}:=\g^{\a\b} M_{i\b}^j$.
\end{definition}

\begin{lemma}\lab{lemma:introductionandpropertiesoftheMalphaij}
Let $M_{i\a}^j$ be the 1-forms on $(\MM, \g)$ as defined in Definition \ref{def:Mialphaj:Kerr}.  Then we have
\bea
\label{def:Mialphaj}
\Ddot_\a\Om_i=M_{i\a}^j\Om_j.
\eea
\end{lemma}

\begin{proof}
Using \eqref{eq:fundamentalpropertiesof1formsOmi}, we have
\beaa
M_{i\a}^j(\Om_j)_a=(\Ddot_\a\Om_i)_b(\Om^j)_b(\Om_j)_a=(\Ddot_\a\Om_i)_b\de_{ab}=(\Ddot_\a\Om_i)_a,
\eeaa
and hence $\Ddot_\a\Om_i=M_{i\a}^j\Om_j$ as stated.
\end{proof}

The following lemma provides a useful property for the symmetric part of $M_{i\a}^j$.
\begin{lemma}
\label{lem:property:Omi}
The symmetric part $(M_{S})_{i\a}^{j}$ of $M_{i\a}^{j}$ w.r.t. $(i,j)$ satisfies 
\bea
\lab{formula:symmetricpartofMmatrices}
(M_{S})_{i\a}^{j}:=\frac{1}{2}(M_{i\a}^{j}+M_{j\a}^{i}), \qquad (M_{S})_{i\a}^{j}=-\frac{1}{2}\pr_\a(x^ix^j).
\eea
\end{lemma}

\begin{proof}
Using \eqref{def:Mialphaj} and the third identity in \eqref{eq:fundamentalpropertiesof1formsOmi}, we have
\beaa
2(M_{S})_{i\a}^{j} &=& M_{i\a}^{j}+M_{j\a}^{i}=(\Ddot_\a\Om_i)\c\Om_j+\Om_i\c(\Ddot_\a\Om_j)\\
&=& \pr_\a(\Om_i\c\Om_j)=\pr_\a(\de_{ij}-x^ix^j)=-\pr_\a(x^ix^j)
\eeaa
as stated in \eqref{formula:symmetricpartofMmatrices}.
This concludes the proof of Lemma \ref{lem:property:Omi}.
\end{proof}

The following lemma calculates the inner products of {horizontal 1-forms and horizontal symmetric 2-tensors} in terms of their scalarization w.r.t. the regular triplet $\{\Om_{i}\}_{i=1,2,3}$. 
\begin{lemma}
\label{lem:ucdotv:product}
{If $u$ and $v$ are} two horizontal $1$-forms, then
\bsub
\bea
u\c v={u(\Om^i)v(\Om_i)}.
\eea
{Also, if $u$ and $v$ are two symmetric horizontal  $2$-tensors}, then
\bea
u\c v={u(\Om^i, \Om^j)v(\Om_i, \Om_j)}.
\eea
\esub
\end{lemma}

\begin{proof}
Consider first the case {where} $u$ and $v$ are horizontal $1$-forms. Using \eqref{eq:fundamentalpropertiesof1formsOmi}, we {have}
\beaa
u(\Om^i) v(\Om_i) = (\Om^i)^b(\Om_i)^cu(e_b) v(e_c) =\de^{bc}u(e_b)v(e_c)=
u\c v,
\eeaa
{as stated. The identity for horizontal symmetric $2$-tensors} follows in the same manner.
\end{proof}


\subsection{From tensors to regular scalars and back}


The following lemma allows to pass from horizontal tensors in $\sk_2$ to scalars and reciprocally.
\begin{lemma}\lab{lemma:backandforthbetweenhorizontaltensorsk2andscalars}
Let $(\MM, \g)$ be a spacetime, $(e_3, e_4)$ be a null pair, and consider the corresponding horizontal structure $\O(\MM)$ introduced in Section \ref{subsection:review-horiz.structures}. Assume that $\Om_i$, $1,2,3$ is a regular triplet in the sense of Definition \ref{def:definitionofregulartripletOmii=123}. Then, the following holds:
\begin{enumerate}
\item Let $\pmb\psi\in\sk_2$ and define the scalars $\psi_{ij}:=\pmb\psi(\Om_i, \Om_j)$, $i,j=1,2,3$. Then:
\begin{itemize}
\item The real-valued scalars $\psi_{ij}$ satisfy
\bea\lab{eq:fundamentalidentitiestoderivefromscalarizationoftensor:realcase}
\psi_{ij}=\psi_{ji}, \qquad x^i\psi_{ij}=0, \qquad (\de^{ij}-x^ix^j)\psi_{ij}=0.
\eea

\item We may recover the tensor $\pmb\psi$ from the real-valued scalars $\psi_{ij}$ by the formula 
\beaa
\pmb\psi_{ab}=\psi_{ij}(\Om^i)_a(\Om^j)_b.
\eeaa
\end{itemize}

\item Reciprocally, let $\psi_{ij}$ be real-valued scalars satisfying the identities \eqref{eq:fundamentalidentitiestoderivefromscalarizationoftensor:realcase}, 
and introduce the real-valued horizontal 2-tensor $\pmb\psi$ by $\pmb\psi_{ab}:=\psi_{ij}(\Om^i)_a(\Om^j)_b$, $a,b=1,2$. Then, we have $\pmb\psi\in\sk_2$ and $\pmb\psi(\Om_i, \Om_j)=\psi_{ij}$ for all $i,j=1,2,3$.
\end{enumerate}
\end{lemma}

\begin{proof}
First, assume that $\pmb\psi\in\sk_2$ and define the real-valued scalars $\psi_{ij}:=\pmb\psi(\Om_i, \Om_j)$, $i,j=1,2,3$. Then, we have $\psi_{ij}=\psi_{ji}$, and, using \eqref{eq:fundamentalpropertiesof1formsOmi},
\beaa
x^i\psi_{ij} &=& x^i\pmb\psi(\Om_i, \Om_j)=\pmb\psi(x^i\Om_i, \Om_j)=\pmb\psi(0, \Om_j)=0,\\
(\de^{ij}-x^ix^j)\psi_{ij} &=& (\Om^i\c\Om^j)\pmb\psi(\Om_i, \Om_j)=(\Om^i)_a(\Om_i)_b(\Om^j)^a(\Om_j)_c\pmb\psi_{bc}=\de_{ab}\de^a\,\!_c\pmb\psi_{bc}=\tr\pmb\psi=0,\\
\psi_{ij}(\Om^i)_a(\Om^j)_b &=& (\Om^i)_a(\Om^j)_b\pmb\psi(\Om_i, \Om_j)=(\Om^i)_a(\Om_i)_c(\Om^j)_b(\Om_j)_d\pmb\psi_{cd}=\de_{ac}\de_{bc}\pmb\psi_{cd}=\pmb\psi_{ab},
\eeaa
as stated in \eqref{eq:fundamentalidentitiestoderivefromscalarizationoftensor:realcase}.

Reciprocally, assume that $\psi_{ij}$ are real-valued scalars satisfying the identities \eqref{eq:fundamentalidentitiestoderivefromscalarizationoftensor:realcase}, and introduce the real-valued horizontal 2-tensor $\pmb\psi$ by $\pmb\psi_{ab}:=\psi_{ij}(\Om^i)_a(\Om^j)_b$, $a,b=1,2$. Then, $\pmb\psi$ is symmetric and, using \eqref{eq:fundamentalpropertiesof1formsOmi},
\beaa
\tr\pmb\psi=\de_{ab}\de^a\,\!_c\pmb\psi_{bc}=(\Om^i)_a(\Om_i)_b(\Om^j)^a(\Om_j)_c\pmb\psi_{bc}=(\Om^i\c\Om^j)\pmb\psi(\Om_i, \Om_j)=(\de^{ij}-x^ix^j)\psi_{ij}=0,
\eeaa
so that $\pmb\psi\in\sk_2$ as stated. Also, using again Lemma \ref{lemma:fundamentalpropertiesof1formsOmi}, we have 
\beaa
\pmb\psi(\Om_i, \Om_j) &=& \pmb\psi_{ab}(\Om_i)^a(\Om_j)^b=\psi_{kl}(\Om^k)_a(\Om^l)_b(\Om_i)^a(\Om_j)^b=\psi_{kl}(\Om^k\c\Om_i)(\Om^l\c\Om_j)\\
&=& (\de_{ik}-x^ix^k)(\de_{jl}-x^jx^l)\psi_{kl}=\psi_{ij}
\eeaa
as stated, where we used several times the identities $x^i\psi_{ij}=0$ and $x^j\psi_{ij}=0$ in the last equality. This concludes the proof of Lemma \ref{lemma:backandforthbetweenhorizontaltensorsk2andscalars}.
\end{proof}

We now consider the analogous problem for horizontal tensors in $\sk_2(\mathbb{C})$.
\begin{lemma}\lab{lemma:backandforthbetweenhorizontaltensorsk2andscalars:complex}
Let $(\MM, \g)$ be a spacetime, $(e_3, e_4)$ be a null pair, and consider the corresponding horizontal structure $\O(\MM)$ introduced in Section \ref{subsection:review-horiz.structures}. Assume that $\Om_i$, $1,2,3$ is a regular triplet in the sense of Definition \ref{def:definitionofregulartripletOmii=123}. Then, the following holds:
\begin{enumerate}
\item\lab{item1:sk2Csatisfyconditions} Let $\pmb\psi\in\sk_2(\mathbb{C})$ and define the complex-valued scalars $\psi_{ij}:=\pmb\psi(\Om_i, \Om_j)$, $i,j=1,2,3$. Then:
\begin{itemize}
\item The complex-valued scalars $\psi_{ij}$ satisfy 
\bea\lab{eq:fundamentalidentitiestoderivefromscalarizationoftensor:complexcase}
\psi_{ij}=\psi_{ji}, \qquad x^i\psi_{ij}=0, \qquad (\de^{ij}-x^ix^j)\psi_{ij}=0, \qquad \in_{ikl}x^l\psi_{kj}+i\psi_{ij}=0.
\eea
\item We may recover the tensor $\pmb\psi$ from the scalars $\psi_{ij}$ by the formula 
\beaa
\pmb\psi_{ab}=\psi_{ij}(\Om^i)_a(\Om^j)_b.
\eeaa
\end{itemize}
\item Reciprocally, let $\psi_{ij}$ be complex-valued scalars satisfying the identities \eqref{eq:fundamentalidentitiestoderivefromscalarizationoftensor:complexcase}, and introduce the complex-valued horizontal 2-tensor $\pmb\psi$ by $\pmb\psi_{ab}:=\psi_{ij}(\Om^i)_a(\Om^j)_b$, $a,b=1,2$. Then, we have $\pmb\psi\in\sk_2(\mathbb{C})$ and $\pmb\psi(\Om_i, \Om_j)=\psi_{ij}$ for all $i,j=1,2,3$.
\end{enumerate}
\end{lemma}

\begin{proof}
First, assume that $\pmb\psi\in\sk_2(\mathbb{C})$ and define the complex-valued scalars $\psi_{ij}:=\pmb\psi(\Om_i, \Om_j)$, $i,j=1,2,3$. Then, as $\dual\pmb\psi=-i\pmb\psi$, we have, using \eqref{eq:fundamentalpropertiesof1formsOmi},
\beaa
0&=& \dual\pmb\psi(\Om_i, \Om_j)+i\pmb\psi(\Om_i, \Om_j)=\in_{ac}\pmb\psi_{cb}(\Om_i)_a(\Om_j)_b+i\psi_{ij}=-\pmb\psi_{cb}(\dual\Om_i)_c(\Om_j)_b+i\psi_{ij}\\
&=& -\pmb\psi_{cb}(\dual\Om_i\c\Om^k)(\Om_k)_c(\Om_j)_b+i\psi_{ij}=(\Om_i\c\dual\Om^k)\psi_{kj}+i\psi_{ij}=\in_{ikl}x^l\psi_{kj}+i\psi_{ij}
\eeaa
as stated, and the other identities follow from Lemma \ref{lemma:backandforthbetweenhorizontaltensorsk2andscalars}.

Reciprocally, assume that $\psi_{ij}$ are complex-valued scalars satisfying the identities \eqref{eq:fundamentalidentitiestoderivefromscalarizationoftensor:complexcase}, and introduce the complex-valued horizontal 2-tensor $\pmb\psi$ by $\pmb\psi_{ab}:=\psi_{ij}(\Om^i)_a(\Om^j)_b$, $a,b=1,2$. Then, by Lemma \ref{lemma:backandforthbetweenhorizontaltensorsk2andscalars}, $\Re(\pmb\psi),\, \Im(\pmb\psi)\in\sk_2$ and $\pmb\psi(\Om_i, \Om_j)=\psi_{ij}$ for all $i,j=1,2,3$. Furthermore, we have for all $a,b=1,2$, using \eqref{eq:fundamentalpropertiesof1formsOmi},
\beaa
\dual\pmb\psi_{ab}+i\pmb\psi_{ab} &=& \in_{ac}\pmb\psi_{cb}+i\psi_{ij}(\Om^i)_a(\Om^j)_b=\in_{ac}\psi_{ij}(\Om^i)_c(\Om^j)_b+i\psi_{ij}(\Om^i)_a(\Om^j)_b\\
&=& \Big(\psi_{ij}(\dual\Om^i)_a+i\psi_{ij}(\Om^i)_a\Big)(\Om^j)_b=\Big(\psi_{ij}(\dual\Om^i\c\Om_k)(\Om^k)_a+i\psi_{ij}(\Om^i)_a\Big)(\Om^j)_b\\
&=& \Big(-\psi_{ij}\in_{ikl}x^l(\Om^k)_a+i\psi_{ij}(\Om^i)_a\Big)(\Om^j)_b= \Big(\in_{ikl}x^l\psi_{kj}+i\psi_{ij}\Big)(\Om^i)_a(\Om^j)_b=0
\eeaa
so that $\pmb\psi\in\sk_2(\mathbb{C})$ as stated. This concludes the proof of Lemma \ref{lemma:backandforthbetweenhorizontaltensorsk2andscalars:complex}.
\end{proof}


\subsection{Scalarization of the tensorial wave operator $\squared_2$}
\lab{sec:scalarizationofthetensorialwaveoperatorsquared2}


The following lemma provides the scalarization of the tensorial wave operator $\squared_2$.
\begin{lemma}\lab{lemma:formoffirstordertermsinscalarazationtensorialwaveeq}
Let $(\MM, \g)$ be a spacetime, $(e_3, e_4)$ be a null pair, and consider the corresponding horizontal structure $\O(\MM)$ introduced in Section  \ref{subsection:review-horiz.structures}. Assume that $\Om_i$, $1,2,3$ is a regular triplet in the sense of Definition \ref{def:definitionofregulartripletOmii=123}. Also, let $\pmb\psi\in\sk_2(\mathbb{C})$ and let $\psi_{ij}$ be the scalars associated to it in view of Lemma \ref{lemma:backandforthbetweenhorizontaltensorsk2andscalars:complex}. Then, we have
\bea
\squared_2\pmb\psi(\Om_i, \Om_j) &=& \square_\g(\psi_{ij}) -S(\psi)_{ij} - (Q\psi)_{ij}
\eea
where 
\bsub
\label{SandV}
\begin{align}
S(\psi)_{ij} ={}& 2M_{i}^{k\a}\pr_\a(\psi_{kj}) +2M_{j}^{k\a}\pr_\a(\psi_{ik}),\\
(Q\psi)_{ij} ={}& (\Ddot^\a M_{i\a}^k)\psi_{kj}+(\Ddot^\a M_{j\a}^k)\psi_{ik} -M_{i\a}^kM_k^{l\a}\psi_{lj}-2M_{i\a}^kM_{j}^{l\a}\psi_{kl}-M_{j\a}^kM_k^{l\a}\psi_{il},
\end{align}
\esub
with the 1-forms $M_{i\a}^j$ defined by \eqref{eq:definitionofMalphaijwithoutambiguity}.
\end{lemma}

\begin{proof}
Using repeatedly Lemma \ref{lemma:introductionandpropertiesoftheMalphaij}, we have
\beaa
\square_\g(\psi_{ij}) &=& \square_\g(\pmb\psi(\Om_i,\Om_j))\\ 
&=& \squared_2\pmb\psi(\Om_i, \Om_j)+2\g^{\a\b}\Ddot_\a\pmb\psi(\Ddot_\b\Om_i, \Om_j)+2\g^{\a\b}\Ddot_\a\pmb\psi(\Om_i, \Ddot_\b\Om_j)+\pmb\psi(\squared_1\Om_i, \Om_j)\\
&&+2\g^{\a\b}\pmb\psi(\Ddot_\a\Om_i, \Ddot_\b\Om_j)+\pmb\psi(\Om_i, \squared_1\Om_j)\\
&=& \squared_2\pmb\psi(\Om_i, \Om_j)+2\g^{\a\b}M_{i\b}^k\Ddot_\a\pmb\psi(\Om_k, \Om_j)+2\g^{\a\b}M_{j\b}^k\Ddot_\a\pmb\psi(\Om_i, \Om_k)\\
&&+\pmb\psi(\Ddot^\a(M_{i\a}^k\Om_k), \Om_j)+2\g^{\a\b}M_{i\a}^kM_{j\b}^l\psi_{kl}+\pmb\psi(\Om_i, \Ddot^\a(M_{j\a}^k\Om_k))\\
&=& \squared_2\pmb\psi(\Om_i, \Om_j)+2\g^{\a\b}M_{i\b}^k\Big(\pr_\a(\psi_{kj}) - M_{k\a}^l\psi_{lj} - M_{j\a}^l\psi_{kl}\Big)\\
&&+2\g^{\a\b}M_{j\b}^k\Big(\pr_\a(\psi_{ik}) - M_{i\a}^l\psi_{lk} - M_{k\a}^l\psi_{il}\Big)+(\Ddot^\a M_{i\a}^k)\psi_{kj}+M_{i\a}^kM_k^{l\a}\psi_{lj}\\
&&+2M_{i\a}^kM_{j}^{l\a}\psi_{kl}+(\Ddot^\a M_{j\a}^k)\psi_{ik}+M_{j\a}^kM_k^{l\a}\psi_{il}\\
&=& \squared_2\pmb\psi(\Om_i, \Om_j)+2M_{i}^{k\a}\pr_\a(\psi_{kj}) +2M_{j}^{k\a}\pr_\a(\psi_{ik})+(\Ddot^\a M_{i\a}^k)\psi_{kj}+(\Ddot^\a M_{j\a}^k)\psi_{ik}\\
&&-M_{i\a}^kM_k^{l\a}\psi_{lj}-2M_{i\a}^kM_{j}^{l\a}\psi_{kl}-M_{j\a}^kM_k^{l\a}\psi_{il}
\eeaa
so that
\beaa
\squared_2\pmb\psi(\Om_i, \Om_j) &=& \square_\g(\psi_{ij}) -S(\psi)_{ij} - (Q\psi)_{ij}
\eeaa
where 
\beaa
S(\psi)_{ij} &=& 2M_{i}^{k\a}\pr_\a(\psi_{kj}) +2M_{j}^{k\a}\pr_\a(\psi_{ik}),\\
(Q\psi)_{ij} &=& (\Ddot^\a M_{i\a}^k)\psi_{kj}+(\Ddot^\a M_{j\a}^k)\psi_{ik} -M_{i\a}^kM_k^{l\a}\psi_{lj}-2M_{i\a}^kM_{j}^{l\a}\psi_{kl}-M_{j\a}^kM_k^{l\a}\psi_{il},
\eeaa
as stated. This concludes the proof of Lemma \ref{lemma:formoffirstordertermsinscalarazationtensorialwaveeq}.
\end{proof}


\subsection{Tensorization defect}
\lab{subsect:tensorizationdefect}


We now consider more general families of scalars $\psi_{ij}$ that are not necessarily derived from the scalarization\footnote{This will naturally occur when extending Teukolsky from $\tau\in (\tau_1, \tau_2)$ to $\tau\in(\tmic, +\infty)$, see Section \ref{sect:extensiontoglobalproblem:Teu}.} of a tensor in $\sk_2(\mathbb{C})$, and we aim at estimating their tensorization defect which is defined as follows. 

\begin{definition}[Tensorization defect]
\lab{def:definitionofthenotationerrforthescalarizationdefect}
For a general family of complex-valued scalars $\psi_{ij}$, define the following error term which estimates the corresponding tensorization defect 
\bea\lab{eq:definitionofthenotationerrforthescalarizationdefect}
\bsplit
\err_{\textrm{TDefect}}[\psi]:=&\Big(\err_{\textrm{TDefect},n}[\psi],\,n=1,2,3,4,5\Big),\\
(\err_{\textrm{TDefect},1}[\psi])_{ij}:=&\psi_{ij}-\psi_{ji},\qquad (\err_{\textrm{TDefect},2}[\psi])_j:=x^i\psi_{ij}, \\ 
(\err_{\textrm{TDefect},3}[\psi])_j:=&x^i\psi_{ji},\qquad \err_{\textrm{TDefect},4}[\psi]:=(\de^{ij}-x^ix^j)\psi_{ij},\\
(\err_{\textrm{TDefect},5}[\psi])_{ij}:=& \in_{ikl}x^l\psi_{kj}+i\psi_{ij}.
\end{split}
\eea
\end{definition}

\begin{remark}
In view of Lemma \ref{lemma:backandforthbetweenhorizontaltensorsk2andscalars:complex}, a family of complex-valued scalars $\psi_{ij}$ comes from the scalarization of a tensor in $\sk_2(\mathbb{C})$ if and only if $\err_{\textrm{TDefect}}[\psi]=0$.
\end{remark}

\begin{remark}
Note that 
\bea
(\err_{\textrm{TDefect},3}[\psi])_j=(\err_{\textrm{TDefect},2}[\psi])_j-x^i(\err_{\textrm{TDefect},1}[\psi])_{ij}.
\eea
\end{remark}

We also approximate general families of scalars by families that are generated by scalarization of a tensor in $\sk_2(\mathbb{C})$. 
\begin{lemma}\lab{lemma:computationerrorscalarizationdeffect}
Let $\psi_{ij}$ be a general family of complex-valued scalars and let $\Pi_2[\psi]$ be defined by
\bea\lab{eq:computationerrorscalarizationdeffect:defPi2}
(\Pi_2[\psi])_{ij}:=\frac{1}{2}\Big((\pi_2[\psi])_{ij}+i\in_{ikl}x^l(\pi_2[\psi])_{kj}\Big),
\eea
where 
\bea
\nn (\pi_2[\psi])_{ij} &:=& \frac{1}{2}(\psi_{ij}+\psi_{ji}) -\frac{1}{2}\Big(x^k(\psi_{kj}+\psi_{jk})x^i+x^k(\psi_{ki}+\psi_{ik})x^j\Big) +x^kx^l\psi_{kl}x^ix^j\\
&&-\frac{1}{2}(\de^{kl}-x^kx^l)\psi_{kl}(\de^{ij}-x^ix^j).
\eea
Then, $\Pi_2[\psi]$ is generated by the scalarization of a tensor in $\sk_2(\mathbb{C})$, i.e., we have 
\beaa
(\Pi_2[\psi])_{ij}=\pmb\Pi_2[\psi](\Om_i, \Om_j), \quad i,j=1,2,3, \quad \pmb\Pi_2[\psi]\in\sk_2(\mathbb{C}).
\eeaa
Furthermore, we have
\beaa
\psi_{ij}=(\Pi_2[\psi])_{ij} +(\widetilde{\err}_0[\psi])_{ij}
\eeaa
where 
\bea
\bsplit
(\widetilde{\err}_0[\psi])_{ij}:=& -\frac{i}{2}(\err_{\textrm{TDefect},5}[\psi])_{ij}+\frac{1}{2}\Big((\err_0[\psi])_{ij}+i\in_{ikl}x^l(\err_0[\psi])_{kj}\Big),\\
(\err_0[\psi])_{ij} :=& \frac{1}{2}(\err_{\textrm{TDefect},1}[\psi])_{ij} +\frac{1}{2}\Big((\err_{\textrm{TDefect},2}[\psi])_j+(\err_{\textrm{TDefect},3}[\psi])_j\Big)x^i\\
& +\frac{1}{2}\Big((\err_{\textrm{TDefect},2}[\psi])_i+(\err_{\textrm{TDefect},3}[\psi])_i\Big)x^j\\
& -x^k(\err_{\textrm{TDefect},2}[\psi])_kx^ix^j +\frac{1}{2}\err_{\textrm{TDefect},4}[\psi](\de^{ij}-x^ix^j).
\end{split}
\eea
\end{lemma}

\begin{proof}
A straightforward computation yields
\beaa
(\pi_2[\psi])_{ij}=(\pi_2[\psi])_{ji}, \qquad x^i(\pi_2[\psi])_{ij}=0, \qquad (\de^{ij}-x^ix^j)(\pi_2[\psi])_{ij}=0.
\eeaa
In view of Lemma \ref{lemma:backandforthbetweenhorizontaltensorsk2andscalars} applied to the families of real-valued scalars $\Re((\pi_2[\psi])_{ij})$ and $\Im((\pi_2[\psi])_{ij})$, we infer the existence of real-valued scalars $\textbf{u},\,\textbf{v}\in\sk_2$ such that 
\beaa
(\pi_2[\psi])_{ij}=\big(\textbf{u}+i\textbf{v}\big)(\Om_i, \Om_j), \quad i,j=1,2,3.
\eeaa
Then, we have
\beaa
\in_{ikl}x^l(\pi_2[\psi])_{kj} &=& (\Om_i\c\dual\Om^k)\big(\textbf{u}+i\textbf{v}\big)(\Om_k, \Om_j)=-\big(\textbf{u}+i\textbf{v}\big)\big((\dual\Om_i\c\Om^k)\Om_k, \Om_j\big)\\
&=& -\big(\textbf{u}+i\textbf{v}\big)(\dual\Om_i, \Om_j)=-\big(\textbf{u}+i\textbf{v}\big)_{ab}\in_{ac}(\Om_i)_c(\Om_j)_b\\
&=& \dual\big(\textbf{u}+i\textbf{v}\big)(\Om_i, \Om_j), \quad i,j=1,2,3,
\eeaa
and hence 
\beaa
(\Pi_2[\psi])_{ij}&=&\frac{1}{2}\Big((\pi_2[\psi])_{ij}+i\in_{ikl}x^l(\pi_2[\psi])_{kj}\Big)\\
&=& \left(\frac{1}{2}\big(\textbf{u}-\dual\textbf{v}\big)+\frac{i}{2}\dual\big(\textbf{u}-\dual\textbf{v}\big)\right)(\Om_i, \Om_j)= \pmb\Pi_2[\psi](\Om_i, \Om_j), \quad i,j=1,2,3, 
\eeaa
where 
\beaa
\pmb\Pi_2[\psi]:=\frac{1}{2}\big(\textbf{u}-\dual\textbf{v}\big)+\frac{i}{2}\dual\big(\textbf{u}-\dual\textbf{v}\big)\in\sk_2(\mathbb{C})
\eeaa
as stated.

Next, notice that we have
\beaa
(\pi_2[\psi])_{ij} &=& \psi_{ij}-\frac{1}{2}(\err_{\textrm{TDefect},1}[\psi])_{ij} -\frac{1}{2}\Big((\err_{\textrm{TDefect},2}[\psi])_j+(\err_{\textrm{TDefect},3}[\psi])_j\Big)x^i\\
&& -\frac{1}{2}\Big((\err_{\textrm{TDefect},2}[\psi])_i+(\err_{\textrm{TDefect},3}[\psi])_i\Big)x^j\\
&& +x^k(\err_{\textrm{TDefect},2}[\psi])_kx^ix^j -\frac{1}{2}\err_{\textrm{TDefect},4}[\psi](\de^{ij}-x^ix^j)
\eeaa
which we rewrite as 
\beaa
(\pi_2[\psi])_{ij} &=& \psi_{ij} - (\err_0[\psi])_{ij}
\eeaa
with the notation 
\beaa
(\err_0[\psi])_{ij} &:=& \frac{1}{2}(\err_{\textrm{TDefect},1}[\psi])_{ij} +\frac{1}{2}\Big((\err_{\textrm{TDefect},2}[\psi])_j+(\err_{\textrm{TDefect},3}[\psi])_j\Big)x^i\\
&& +\frac{1}{2}\Big((\err_{\textrm{TDefect},2}[\psi])_i+(\err_{\textrm{TDefect},3}[\psi])_i\Big)x^j\\
&& -x^k(\err_{\textrm{TDefect},2}[\psi])_kx^ix^j +\frac{1}{2}\err_{\textrm{TDefect},4}[\psi](\de^{ij}-x^ix^j).
\eeaa
Then, we infer
\beaa
(\Pi_2[\psi])_{ij} &=& \frac{1}{2}\Big((\pi_2[\psi])_{ij}+i\in_{ikl}x^l(\pi_2[\psi])_{kj}\Big)\\
&=& \frac{1}{2}\Big(\psi_{ij}+i\in_{ikl}x^l\psi_{kj}\Big) - \frac{1}{2}\Big((\err_0[\psi])_{ij}+i\in_{ikl}x^l(\err_0[\psi])_{kj}\Big)\\
&=& \psi_{ij} -(\widetilde{\err}_0[\psi])_{ij}
\eeaa
where we have introduced the notation 
\beaa
(\widetilde{\err}_0[\psi])_{ij} &=& \frac{1}{2}\Big(\psi_{ij}-i\in_{ikl}x^l\psi_{kj}\Big)+\frac{1}{2}\Big((\err_0[\psi])_{ij}+i\in_{ikl}x^l(\err_0[\psi])_{kj}\Big)\\
&=& -\frac{i}{2}(\err_{\textrm{TDefect},5}[\psi])_{ij}+\frac{1}{2}\Big((\err_0[\psi])_{ij}+i\in_{ikl}x^l(\err_0[\psi])_{kj}\Big).
\eeaa
This concludes the proof of Lemma \ref{lemma:computationerrorscalarizationdeffect}.
\end{proof}

Finally, we derive wave equations for the tensorization defect. 
\begin{lemma}\lab{lemma:waveequationsfortensordeffects}
Let $\psi_{ij}$ be a family of complex-valued scalars satisfying 
\beaa
\big(\square_\g+V\big)\psi_{ij} &=& S(\psi)_{ij} + (Q\psi)_{ij},
\eeaa
where $S$ and $Q$ are defined in \eqref{SandV} and where $V$ is a real-valued function, and let $\err_{\textrm{TDefect}}[\psi]$ be associated to the family $\psi_{ij}$ as in \eqref{eq:definitionofthenotationerrforthescalarizationdefect}. Then, $\err_{\textrm{TDefect},1}[\psi]$ satisfies 
\bea\lab{eq:waveeqpsiijminuspsiji}
(\square_\g+V)\big((\err_{\textrm{TDefect},1}[\psi])_{ij}\big) &=& S(\err_{\textrm{TDefect},1}[\psi])_{ij} + (Q\err_{\textrm{TDefect},1}[\psi])_{ij},
\eea
$\err_{\textrm{TDefect},2}[\psi]$ satisfies 
\bea\lab{eq:waveeqxipsiij}
\nn(\square_\g+V)((\err_{\textrm{TDefect},2}[\psi])_j) &=& 2M_j^{k\a}\pr_\a((\err_{\textrm{TDefect},2}[\psi])_k) +(\Ddot^\a M_{j\a}^k)(\err_{\textrm{TDefect},2}[\psi])_k\\
&& -M_{j\a}^kM_k^{l\a}(\err_{\textrm{TDefect},2}[\psi])_l,
\eea
$\err_{\textrm{TDefect},4}[\psi]$ satisfies 
\bea\lab{eq:waveeqfortraceofpsi:deltaijpsiij}
\nn &&(\square_\g+V)(\err_{\textrm{TDefect},4}[\psi])\\ 
\nn &=& -2\pr^{\a} (x^i)\pr_{\a}((\err_{\textrm{TDefect},2}[\psi])_i)  -2\pr^{\a} (x^i)\pr_{\a}((\err_{\textrm{TDefect},3}[\psi])_i)\\
&& +\Big(-2\square_\g(x^i) - x^j(\Ddot^\a M_{j\a}^i) -\pr_\a(x^k)M_k^{i\a}\Big)(\err_{\textrm{TDefect},2}[\psi])_i\nn\\
&& +\Big(- 2\square_\g(x^i) -  x^j(\Ddot^\a M_{j\a}^i)  + \pr_\a(x^k)M_k^{i\a}\Big)(\err_{\textrm{TDefect},3}[\psi])_i
\eea
and $\err_{\textrm{TDefect},5}[\psi]$ satisfies 
\bea\lab{eq:waveeqforlastdefecttensor:antiselfdual}
\nn&& (\square_\g+V)((\err_{\textrm{TDefect},5}[\psi])_{ij})\\ 
\nn&=& S(\err_{\textrm{TDefect},5}[\psi])_{ij} +(Q\err_{\textrm{TDefect},5}[\psi])_{ij} + 2\in_{ikl}x^k\g^{\a\b}\pr_\a(x^l)\pr_\b((\err_{\textrm{TDefect},2}[\psi])_j)\\ 
\nn&& +\Big(\in_{ikl}x^k\square_\g(x^l)  +\in_{ikl}\g^{\a\b}\pr_\a(x^l)\pr_\b(x^k) -\in_{knl}x^nM_{i\a}^k\pr^\a(x^l)\Big)(\err_{\textrm{TDefect},2}[\psi])_j\\
&&  -2\in_{ikl}x^kM_{j}^{m\a}\pr_\a(x^l)(\err_{\textrm{TDefect},2}[\psi])_m.
\eea
\end{lemma}

\begin{proof}
For simplicity, we assume that $V=0$ as the general case is completely analogous. 
We start with the proof of \eqref{eq:waveeqpsiijminuspsiji}. We have
\beaa
S(\psi)_{ij} - S(\psi)_{ji} &=& 2M_{i}^{k\a}\pr_\a(\psi_{kj}) +2M_{j}^{k\a}\pr_\a(\psi_{ik})-2M_{j}^{k\a}\pr_\a(\psi_{ki}) -2M_{i}^{k\a}\pr_\a(\psi_{jk})\\
&=& 2M_{i}^{k\a}\pr_\a(\psi_{kj}-\psi_{jk})  +2M_{j}^{k\a}\pr_\a(\psi_{ik}-\psi_{ki})
\eeaa
and 
\beaa
(Q\psi)_{ij} - (Q\psi)_{ji} &=& 2M_{i}^{k\a}\pr_\a(\psi_{kj}) +2M_{j}^{k\a}\pr_\a(\psi_{ik})-2M_{j}^{k\a}\pr_\a(\psi_{ki}) -2M_{i}^{k\a}\pr_\a(\psi_{jk})\\
&=& (\Ddot^\a M_{i\a}^k)\psi_{kj}+(\Ddot^\a M_{j\a}^k)\psi_{ik} -M_{i\a}^kM_k^{l\a}\psi_{lj}-2M_{i\a}^kM_{j}^{l\a}\psi_{kl}-M_{j\a}^kM_k^{l\a}\psi_{il}\\
&& -(\Ddot^\a M_{j\a}^k)\psi_{ki}-(\Ddot^\a M_{i\a}^k)\psi_{jk} +M_{j\a}^kM_k^{l\a}\psi_{li}+2M_{j\a}^kM_{i}^{l\a}\psi_{kl}+M_{i\a}^kM_k^{l\a}\psi_{jl}\\
&=& (\Ddot^\a M_{i\a}^k)(\psi_{kj}-\psi_{jk})+(\Ddot^\a M_{j\a}^k)(\psi_{ik}-\psi_{ki}) -M_{i\a}^kM_k^{l\a}(\psi_{lj}-\psi_{jl})\\
&&-2M_{i\a}^kM_{j}^{l\a}(\psi_{kl}-\psi_{lk})-M_{j\a}^kM_k^{l\a}(\psi_{il}-\psi_{li})
\eeaa
and hence
\beaa
\square_\g(\psi_{ij}-\psi_{ji}) &=& S(\psi)_{ij} +(Q\psi)_{ij} - S(\psi)_{ji} - (Q\psi)_{ji}\\
&=& 2M_{i}^{k\a}\pr_\a(\psi_{kj}-\psi_{jk})  +2M_{j}^{k\a}\pr_\a(\psi_{ik}-\psi_{ki})\\
&&+(\Ddot^\a M_{i\a}^k)(\psi_{kj}-\psi_{jk})+(\Ddot^\a M_{j\a}^k)(\psi_{ik}-\psi_{ki}) -M_{i\a}^kM_k^{l\a}(\psi_{lj}-\psi_{jl})\\
&&-2M_{i\a}^kM_{j}^{l\a}(\psi_{kl}-\psi_{lk})-M_{j\a}^kM_k^{l\a}(\psi_{il}-\psi_{li}). 
\eeaa
Since $(\err_{\textrm{TDefect},1}[\psi])_{ij}=\psi_{ij}-\psi_{ji}$, we infer
\beaa
\square_\g\big((\err_{\textrm{TDefect},1}[\psi])_{ij}\big) &=& S(\err_{\textrm{TDefect},1}[\psi])_{ij} + (Q\err_{\textrm{TDefect},1}[\psi])_{ij},
\eeaa
as stated in \eqref{eq:waveeqpsiijminuspsiji}. 

Next, we prove \eqref{eq:waveeqxipsiij}. We compute
\bea\lab{eq:waveeqxipsiij:aux1}
\nn\square_\g((\err_{\textrm{TDefect},2}[\psi])_j)  &=& \square_\g(x^i\psi_{ij})\\
\nn &=& x^i\square_\g(\psi_{ij})+2\g^{\a\b}\pr_\b(x^i)\pr_\b(\psi_{ij})+\square_\g(x^i)\psi_{ij}\\
&=& x^i(S(\psi)_{ij}+(Q\psi)_{ij})+2\g^{\a\b}\pr_\b(x^i)\pr_\b(\psi_{ij})+\square_\g(x^i)\psi_{ij},
\eea
and we first  consider the first-order terms which we rewrite as follows 
\beaa
&& \frac{1}{2}x^iS(\psi)_{ij}+\g^{\a\b}\pr_\b(x^i)\pr_\b(\psi_{ij})\\
&=& x^i\big(M_i^{k\a}\pr_\a(\psi_{kj})+M_j^{k\a}\pr_\a(\psi_{ik})\big)+\g^{\a\b}\pr_\b(x^i)\pr_\b(\psi_{ij})\\
&=& x^iM_i^{k\a}\pr_\a(\psi_{kj})+M_j^{k\a}\pr_\a(x^i\psi_{ik}) -M_j^{k\a}\pr_\a(x^i)\psi_{ik}  +\g^{\a\b}\pr_\b(x^i)\pr_\b(\psi_{ij}).
\eeaa
Now, since $x^i\Om_i=0$ and $\Om_i\c\Om_j=\de^{ij}-x^ix^j$ by assumption, we have
\bea
\lab{eq:xiMialphakformula}
x^iM_{i\a}^k &=& x^i(\Ddot_\a\Om_i)\c\Om_k = \Ddot_\a(x^i\Om_i)\c\Om_k -\pr_\a(x^i)(\Om_i\c\Om_k)\nn\\
&=& -\pr_\a(x^i)(\de^{ik}-x^ix^k) =-\pr_\a(x^k)
\eea
where we also used $x^i\pr_\a(x_i)=\pr_\a(x^ix_i)=\pr_\a(1)=0$, and hence
\beaa
&& \frac{1}{2}x^iS(\psi)_{ij}+\g^{\a\b}\pr_\b(x^i)\pr_\b(\psi_{ij})\\
&=& \g^{\a\b}x^iM_{i\b}^{k}\pr_\a(\psi_{kj})+M_j^{k\a}\pr_\a(x^i\psi_{ik}) -M_j^{k\a}\pr_\a(x^i)\psi_{ik}  +\g^{\a\b}\pr_\b(x^i)\pr_\b(\psi_{ij})\\
&=& -\g^{\a\b}\pr_\b(x^k)\pr_\a(\psi_{kj})+M_j^{k\a}\pr_\a(x^i\psi_{ik}) -M_j^{k\a}\pr_\a(x^i)\psi_{ik}  +\g^{\a\b}\pr_\b(x^i)\pr_\b(\psi_{ij})\\
&=& M_j^{k\a}\pr_\a(x^i\psi_{ik}) -M_j^{k\a}\pr_\a(x^i)\psi_{ik}.
\eeaa
Plugging the above in \eqref{eq:waveeqxipsiij:aux1}, we infer
\bea\lab{eq:waveeqxipsiij:aux2}
\square_\g((\err_{\textrm{TDefect},2}[\psi])_j)  = 2M_j^{k\a}\pr_\a(x^i\psi_{ik}) - 2M_j^{k\a}\pr_\a(x^i)\psi_{ik} +x^i(Q\psi)_{ij}+\square_\g(x^i)\psi_{ij}.
\eea

Next, we simplify the before to last term on the RHS of \eqref{eq:waveeqxipsiij:aux2} which is given by
\beaa
x^i(Q\psi)_{ij} = x^i(\Ddot^\a M_{i\a}^k)\psi_{kj}+(\Ddot^\a M_{j\a}^k)x^i\psi_{ik} - x^iM_{i\a}^kM_k^{l\a}\psi_{lj}-2x^iM_{i\a}^kM_{j}^{l\a}\psi_{kl}-M_{j\a}^kM_k^{l\a}x^i\psi_{il}.
\eeaa
Using again the above identity $x^iM_{i\a}^k  =-\pr_\a(x^k)$, we obtain 
\beaa
x^i(Q\psi)_{ij} = x^i(\Ddot^\a M_{i\a}^k)\psi_{kj}+(\Ddot^\a M_{j\a}^k)x^i\psi_{ik} +\pr_\a(x^k)M_k^{l\a}\psi_{lj} +2\pr_\a(x^k)M_{j}^{l\a}\psi_{kl}-M_{j\a}^kM_k^{l\a}x^i\psi_{il}.
\eeaa 
Plugging in \eqref{eq:waveeqxipsiij:aux2}, this yields 
\bea\lab{eq:waveeqxipsiij:aux3}
\nn\square_\g((\err_{\textrm{TDefect},2}[\psi])_j)  &=& 2M_j^{k\a}\pr_\a(x^i\psi_{ik}) +(\Ddot^\a M_{j\a}^k)x^i\psi_{ik} -M_{j\a}^kM_k^{l\a}x^i\psi_{il}\\ 
&&+x^i(\Ddot^\a M_{i\a}^k)\psi_{kj} +\pr_\a(x^k)M_k^{l\a}\psi_{lj}  +\square_\g(x^i)\psi_{ij}.
\eea

Finally, we simplify the second line of the RHS in \eqref{eq:waveeqxipsiij:aux3}. Using again  the above identity $x^iM_{i\a}^k  =-\pr_\a(x^k)$, we have
\beaa
x^i(\Ddot^\a M_{i\a}^k)\psi_{kj} &=& \Ddot^\a(x^iM_{i\a}^k)\psi_{kj} -\pr_\a(x^i)M_{i\a}^k\psi_{kj}=-\D_\a\D_\a(x^k)\psi_{kj} -\pr_\a(x^i)M_{i\a}^k\psi_{kj}\\
&=& -\square_\g(x^k)\psi_{kj} -\pr_\a(x^i)M_{i\a}^k\psi_{kj}.
\eeaa
Hence, plugging in \eqref{eq:waveeqxipsiij:aux3}, we see that the second line of the RHS in \eqref{eq:waveeqxipsiij:aux3} cancels and we deduce
\beaa
\square_\g((\err_{\textrm{TDefect},2}[\psi])_j)  &=& 2M_j^{k\a}\pr_\a(x^i\psi_{ik}) +(\Ddot^\a M_{j\a}^k)x^i\psi_{ik} -M_{j\a}^kM_k^{l\a}x^i\psi_{il}
\eeaa
or 
\beaa
\nn\square_\g((\err_{\textrm{TDefect},2}[\psi])_j) &=& 2M_j^{k\a}\pr_\a((\err_{\textrm{TDefect},2}[\psi])_k) +(\Ddot^\a M_{j\a}^k)(\err_{\textrm{TDefect},2}[\psi])_k\\
&& -M_{j\a}^kM_k^{l\a}(\err_{\textrm{TDefect},2}[\psi])_l,
\eeaa
as stated in \eqref{eq:waveeqxipsiij}.

Next, we prove \eqref{eq:waveeqfortraceofpsi:deltaijpsiij}. We compute
\bea
\lab{eq:waveactingontraceofpsi:pf}
&&\square_\g(\err_{\textrm{TDefect},4}[\psi])=\square_\g((\de^{ij}-x^ix^j)\psi_{ij}) =\square_\g(\de^{ij}\psi_{ij}) 
-\square_{\g} (x^i x^j \psi_{ij})\nn\\
&=&\de^{ij}\square_\g(\psi_{ij}) -x^ix^j\square_\g(\psi_{ij}) - 2\pr^{\a}(x^i x^j)\pr_{\a}(\psi_{ij})\nn\\ 
&& - \big(2\pr^{\a} (x^i)\pr_{\a}(x^j) +\square_\g(x^i)x^j +\square_\g(x^j)x^i\big)  \psi_{ij} \nn\\
&=& (\de^{ij}-x^i x^j)\big(S(\psi)_{ij} + (Q\psi)_{ij}\big)  -2\pr^{\a} (x^i)\pr_{\a}(x^j\psi_{ij})-2\pr^{\a} (x^j)\pr_{\a}(x^i\psi_{ij})\nn\\
&& +4\pr^{\a}(x^i)\pr_{\a}(x^j)\psi_{ij} -  2\pr^{\a}(x^i)\pr_{\a}(x^j) \psi_{ij} 
- \square_\g(x^i)x^j\psi_{ij} -\square_\g(x^j)x^i\psi_{ij} \nn\\
&=& (\de^{ij}-x^i x^j)S(\psi)_{ij}  -2\pr^{\a} (x^i)\pr_{\a}(x^j\psi_{ij})-2\pr^{\a} (x^j)\pr_{\a}(x^i\psi_{ij})\nn\\
&& +(\de^{ij}-x^i x^j)(Q\psi)_{ij}+2\pr^{\a}(x^i)\pr_{\a}(x^j)\psi_{ij}  
- \square_\g(x^i)x^j\psi_{ij} -\square_\g(x^j)x^i\psi_{ij} 
\eea
where we have substituted in the wave equation for $\psi_{ij}$ in the before to the last step.  We then compute the first term in the before to last line of  \eqref{eq:waveactingontraceofpsi:pf}, given the forms of $S(\psi)_{ij}$ in \eqref{SandV}:
\bea
\lab{eq:traceofSpsipart:pf}
 &&(\de^{ij}-x^i x^j ) S(\psi)_{ij} \nn\\
 &=&4M_{i}^{k\a}\pr_\a(\psi_{ki}) 
-2M_{i}^{k\a}\pr_\a(\psi_{ki}-\psi_{ik}) 
 - 2x^iM_{i}^{k\a}x^j\pr_\a(\psi_{kj}) -2x^jM_{j}^{k\a}x^i\pr_\a(\psi_{ik})\nn\\
 &=&4(M_{S})_{i}^{k\a}\pr_\a(\psi_{ki}) + 2\pr^{\a}(x^k)x^j\pr_\a(\psi_{kj})+ 2\pr^{\a}(x^k)x^i\pr_\a(\psi_{ik}) \nn\\
 &=&-2\pr^{\a}(x^i x^k) \pr_\a(\psi_{ki}) + 2\pr^{\a}(x^k)x^j\pr_\a(\psi_{kj})+ 2\pr^{\a}(x^k)x^i\pr_\a(\psi_{ik}) =0
\eea
where we used \eqref{eq:xiMialphakformula} and the formula \eqref{formula:symmetricpartofMmatrices} for the symmetric part $(M_{S})_{i\a}^{j}$ of $M_{i\a}^{j}$.

Next, we compute the first term in the last line of \eqref{eq:waveactingontraceofpsi:pf} using the formula \eqref{SandV}:
\beaa
&&(\de^{ij}-x^i x^j ) (Q\psi)_{ij}\nn\\
&=&(\de^{ij}-x^i x^j ) \Big( (\Ddot^\a M_{i\a}^k)\psi_{kj}+(\Ddot^\a M_{j\a}^k)\psi_{ik} -M_{i\a}^kM_k^{l\a}\psi_{lj}-2M_{i\a}^kM_{j}^{l\a}\psi_{kl}-M_{j\a}^kM_k^{l\a}\psi_{il}\Big)\nn\\
&=&2(\Ddot^\a M_{i\a}^k)\psi_{ki}
+(\Ddot^\a M_{i\a}^k)(\psi_{ik}-\psi_{ki}) - 2M_{i\a}^kM_k^{l\a}\psi_{li}
-2M_{i\a}^kM_{i}^{l\a}\psi_{kl} + 2x^i x^jM_{i\a}^kM_{j}^{l\a}\psi_{kl}\nn\\
&&-  \Big( x^i(\Ddot^\a M_{i\a}^k)x^j\psi_{kj}+x^j(\Ddot^\a M_{j\a}^k)x^i\psi_{ik} - x^i M_{i\a}^kM_k^{l\a}x^j\psi_{lj} -x^j M_{j\a}^kM_k^{l\a}x^i\psi_{il}\Big)\nn\\
&=&2(\Ddot^\a M_{i\a}^k)\psi_{ki} +(\Ddot^\a M_{i\a}^k)(\psi_{ik}-\psi_{ki})- 2M_{i\a}^kM_k^{l\a}\psi_{li}
-2M_{i\a}^kM_{i}^{l\a}\psi_{kl} + 2\pr_{\a}(x^k)\pr^{\a}(x^l)\psi_{kl}\nn\\
&&-  \Big( x^i(\Ddot^\a M_{i\a}^k)x^j\psi_{kj}+x^j(\Ddot^\a M_{j\a}^k)x^i\psi_{ik} +\pr_\a(x^k)M_k^{l\a}x^j\psi_{lj} +\pr_\a(x^k)M_k^{l\a}x^i\psi_{il}\Big),
\eeaa
where, in the last step, we have used the formula \eqref{eq:xiMialphakformula}. Substituting this together with equation \eqref{eq:traceofSpsipart:pf} into \eqref{eq:waveactingontraceofpsi:pf}, we deduce
\bea\lab{eq:waveactingontraceofpsi:v2:pf}
&&\square_\g(\err_{\textrm{TDefect},4}[\psi])\nn\\
&=&   -2\pr^{\a} (x^i)\pr_{\a}(x^j\psi_{ij})-2\pr^{\a} (x^j)\pr_{\a}(x^i\psi_{ij}) - \square_\g(x^i)x^j\psi_{ij} -\square_\g(x^j)x^i\psi_{ij}\nn\\
\nn&&-  \Big( x^i(\Ddot^\a M_{i\a}^k)x^j\psi_{kj}+x^j(\Ddot^\a M_{j\a}^k)x^i\psi_{ik} +\pr_\a(x^k)M_k^{l\a}x^j\psi_{lj} +\pr_\a(x^k)M_k^{l\a}x^i\psi_{il}\Big)\\
&&+2(\Ddot^\a M_{i\a}^k)\psi_{ki} +(\Ddot^\a M_{i\a}^k)(\psi_{ik}-\psi_{ki})- 2M_{i\a}^kM_k^{l\a}\psi_{li}
-2M_{i\a}^kM_{i}^{l\a}\psi_{kl}\nn\\
&&+4\pr^{\a}(x^i)\pr_{\a}(x^j)\psi_{ij}.
 \eea
 
Now, we compute the before to last line of the RHS of \eqref{eq:waveactingontraceofpsi:v2:pf}. For the two first terms in the before to last line on the RHS of \eqref{eq:waveactingontraceofpsi:v2:pf}, we have, using again \eqref{formula:symmetricpartofMmatrices},
\bea
\lab{estimate:thirdlastterms:tracevanish}
&& 2(\Ddot^\a M_{i\a}^k)\psi_{ki}+(\Ddot^\a M_{i\a}^k)(\psi_{ik}-\psi_{ki})\nn\\
&=& 2(\Ddot^\a (M_S)_{i\a}^k)\psi_{ki} + (\Ddot^\a M_{i\a}^k)(\psi_{ki} - \psi_{ik})+(\Ddot^\a M_{i\a}^k)(\psi_{ik}-\psi_{ki})\nn\\
 &=&- \Ddot^\a\Ddot_{\a} (x^i x^k)\psi_{ki} =-\square_{\g}(x^i x^k)\psi_{ki}\nn\\
 &=&  -\square_\g(x^i)x^j\psi_{ij}  -\square_\g(x^j)x^i\psi_{ij} - 2\pr^{\a}(x^i)\pr_{\a}(x^j)\psi_{ij}
\eea
and for the last two terms on the before to last line of the RHS of \eqref{eq:waveactingontraceofpsi:v2:pf}, we have
\bea
\lab{estimate:lasttwoterms:tracevanish}
&&- 2M_{i\a}^kM_k^{l\a}\psi_{li} -2M_{i\a}^kM_{i}^{l\a}\psi_{kl}\nn \\
&=& -2M_{i\a}^kM_k^{l\a}\psi_{li} -2M_{k\a}^iM_{k}^{l\a}\psi_{il}\nn\\
&=&-2M_{k}^{l\a}\big(2(M_S)_{i\a}^k\psi_{li} +M_{k\a}^i (\psi_{il}-\psi_{li})\big)\nn\\
&=&2M_{k}^{l\a}\pr_{\a}(x^i x^k) \psi_{li} - 2M_{k}^{l\a}M_{k\a}^i (\psi_{il}-\psi_{li})\nn\\
&=&2M_{k}^{l\a}\pr_{\a}(x^k)x^i \psi_{li} - 2\pr^{\a}(x^l)\pr_{\a} (x^i) \psi_{li}
\eea
where in the before to last step we have used \eqref{formula:symmetricpartofMmatrices},  and in the last step we have used \eqref{eq:xiMialphakformula}. Plugging in the equalities \eqref{estimate:thirdlastterms:tracevanish} and \eqref{estimate:lasttwoterms:tracevanish} into equation \eqref{eq:waveactingontraceofpsi:v2:pf}, we infer
\beaa
&&\square_\g(\err_{\textrm{TDefect},4}[\psi])\nn\\
&=&   -2\pr^{\a} (x^i)\pr_{\a}(x^j\psi_{ij})-2\pr^{\a} (x^j)\pr_{\a}(x^i\psi_{ij}) - 2\square_\g(x^i)x^j\psi_{ij} -2\square_\g(x^j)x^i\psi_{ij}\nn\\
\nn&& -x^i(\Ddot^\a M_{i\a}^k)x^j\psi_{kj}-x^j(\Ddot^\a M_{j\a}^k)x^i\psi_{ik} +\pr_\a(x^k)M_k^{l\a}x^j\psi_{lj} -\pr_\a(x^k)M_k^{l\a}x^i\psi_{il}
 \eeaa
or
\beaa
\square_\g(\err_{\textrm{TDefect},4}[\psi]) &=&  -2\pr^{\a} (x^i)\pr_{\a}((\err_{\textrm{TDefect},2}[\psi])_i) -2\pr^{\a} (x^i)\pr_{\a}((\err_{\textrm{TDefect},3}[\psi])_i)\\
&& +\Big(-2\square_\g(x^i) - x^j(\Ddot^\a M_{j\a}^i) -\pr_\a(x^k)M_k^{i\a}\Big)(\err_{\textrm{TDefect},2}[\psi])_i\nn\\
&& +\Big(- 2\square_\g(x^i) -  x^j(\Ddot^\a M_{j\a}^i)  + \pr_\a(x^k)M_k^{i\a}\Big)(\err_{\textrm{TDefect},3}[\psi])_i 
 \eeaa
as stated in \eqref{eq:waveeqfortraceofpsi:deltaijpsiij}. 

Finally, we prove \eqref{eq:waveeqforlastdefecttensor:antiselfdual}. To this end, we  introduce the auxiliary family of scalars $\widetilde{\psi}_{ij}$ given by 
\beaa
\widetilde{\psi}_{ij}:=\in_{ikl}x^l\psi_{kj}.
\eeaa
We compute 
\bea\lab{eq:waveeqinijkxlpsikj:aux0}
\nn\square_\g(\widetilde{\psi}_{ij}) &=& \nn\square_\g(\in_{ikl}x^l\psi_{kj})\\
\nn&=& \in_{ikl}x^l\square_\g(\psi_{kj})+2\in_{ikl}\g^{\a\b}\pr_\a(x^l)\pr_\b(\psi_{kj})+\in_{ikl}\square_\g(x^l)\psi_{kj}\\
\nn&=& \in_{ikl}x^l\big(S(\psi)_{kj}+(Q\psi)_{kj}\big)+2\in_{ikl}\g^{\a\b}\pr_\a(x^l)\pr_\b(\psi_{kj})\\
&&+\in_{ikl}\square_\g(x^l)\psi_{kj}
\eea
and we first  consider the first-order terms which we rewrite as follows 
\beaa
&&\frac{1}{2}\in_{ikl}x^lS(\psi)_{kj}+\in_{ikl}\g^{\a\b}\pr_\a(x^l)\pr_\b(\psi_{kj})\\
&=& \in_{ikl}x^l\Big(M_{k}^{n\a}\pr_\a(\psi_{nj}) +M_{j}^{n\a}\pr_\a(\psi_{kn})\Big)+\in_{ikl}\g^{\a\b}\pr_\a(x^l)\pr_\b(\psi_{kj})\\
&=& \in_{ikl}x^lM_{k}^{n\a}\pr_\a(\psi_{nj}) +M_{j}^{n\a}\pr_\a(\in_{ikl}x^l\psi_{kn}) -M_{j}^{n\a}\in_{ikl}\pr_\a(x^l)\psi_{kn}+\in_{ikl}\g^{\a\b}\pr_\a(x^l)\pr_\b(\psi_{kj})\\
&=& \in_{ikl}x^lM_{k}^{n\a}\pr_\a(\psi_{nj}) +M_{j}^{n\a}\pr_\a(\widetilde{\psi}_{in}) -M_{j}^{n\a}\in_{ikl}\pr_\a(x^l)\psi_{kn}+\in_{ikl}\g^{\a\b}\pr_\a(x^l)\pr_\b(\psi_{kj}). 
\eeaa
Now, since 
\bea\lab{eq:waveeqinijkxlpsikj:aux0:00}
\nn\in_{ikl}x^lM_{k\a}^{n} &=& \in_{ikl}x^l(\Ddot_\a\Om_k)\c\Om_n=\Ddot_\a(\in_{ikl}x^l\Om_k)\c\Om_n- \in_{ikl}\pr_\a(x^l)\Om_k\c\Om_n\\
&=&\Ddot_\a(\in_{ikl}x^l\Om_k)\c\Om_n- \in_{ikl}\pr_\a(x^l)(\de^{kn}-x^kx^n),
\eea
we infer
\beaa
&&\frac{1}{2}\in_{ikl}x^lS(\psi)_{kj}+\in_{ikl}\g^{\a\b}\pr_\a(x^l)\pr_\b(\psi_{kj})\\
&=& \g^{\a\b}\Big(\Ddot_\a(\in_{ikl}x^l\Om_k)\c\Om_n- \in_{ikl}\pr_\a(x^l)(\de^{kn}-x^kx^n)\Big)\pr_\b(\psi_{nj})\\
&& +M_{j}^{n\a}\pr_\a(\widetilde{\psi}_{in}) -M_{j}^{n\a}\in_{ikl}\pr_\a(x^l)\psi_{kn}+\in_{ikl}\g^{\a\b}\pr_\a(x^l)\pr_\b(\psi_{kj})\\
&=& \g^{\a\b}\Ddot_\a(\in_{ikl}x^l\Om_k)\c\Om_n\pr_\b(\psi_{nj}) +M_{j}^{n\a}\pr_\a(\widetilde{\psi}_{in}) -M_{j}^{n\a}\in_{ikl}\pr_\a(x^l)\psi_{kn}\\
&& + \g^{\a\b}\in_{ikl}\pr_\a(x^l)x^kx^n\pr_\b(\psi_{nj}) 
\eeaa
and hence
\bea\lab{eq:waveeqinijkxlpsikj:aux1}
\nn&&\frac{1}{2}\in_{ikl}x^lS(\psi)_{kj}+\in_{ikl}\g^{\a\b}\pr_\a(x^l)\pr_\b(\psi_{kj})\\
\nn&=& \g^{\a\b}\Ddot_\a(\in_{ikl}x^l\Om_k)\c\Om_n\pr_\b(\psi_{nj}) +M_{j}^{n\a}\pr_\a(\widetilde{\psi}_{in}) + \in_{ikl}x^k\g^{\a\b}\pr_\a(x^l)\pr_\b(x^n\psi_{nj})\\
&& -M_{j}^{n\a}\in_{ikl}\pr_\a(x^l)\psi_{kn} -   \in_{ikl}x^k\g^{\a\b}\pr_\a(x^l)\pr_\b(x^n)\psi_{nj}. 
\eea

Next, recalling that $\in_{ikl}x^l=\Om_i\c\dual\Om_k$, we have 
\beaa
\Ddot_\a(\in_{ikl}x^l\Om_k)\c\Om_n &=& \Ddot_\a((\Om_i\c\dual\Om_k)\Om_k)\c\Om_n =-\Ddot_\a((\Om_k\c\dual\Om_i)\Om_k)\c\Om_n.
\eeaa
Since 
\beaa
(\Om_k\c\dual\Om_i)(\Om_k)_a=(\Om_k)_a(\Om_k)_b(\c\dual\Om_i)_b=\de_{ab}(\c\dual\Om_i)_b=(\dual\Om_i)_a,
\eeaa
where we used $(\Om_k)_a(\Om_k)_b=\de_{ab}$, we infer $(\Om_k\c\dual\Om_i)\Om_k=\dual\Om_i$ and hence
\bea\lab{eq:waveeqinijkxlpsikj:aux1:00}
\nn\Ddot_\a(\in_{ikl}x^l\Om_k)\c\Om_n &=&-\Ddot_\a(\dual\Om_i)\c\Om_n=\Ddot_\a\Om_i\c\dual\Om_n=(\Ddot_\a\Om_i\c\Om_k)(\Om_k\c\dual\Om_n)\\
&=& M_{i\a}^k\in_{knl}x^l.
\eea
Plugging in \eqref{eq:waveeqinijkxlpsikj:aux1}, we deduce
\beaa
\nn&&\frac{1}{2}\in_{ikl}x^lS(\psi)_{kj}+\in_{ikl}\g^{\a\b}\pr_\a(x^l)\pr_\b(\psi_{kj})\\
\nn&=& M_{i}^{k\a}\in_{knl}x^l\pr_\a(\psi_{nj}) +M_{j}^{n\a}\pr_\a(\widetilde{\psi}_{in}) + \in_{ikl}x^k\g^{\a\b}\pr_\a(x^l)\pr_\b(x^n\psi_{nj})\\
&& -M_{j}^{n\a}\in_{ikl}\pr_\a(x^l)\psi_{kn} -   \in_{ikl}x^k\g^{\a\b}\pr_\a(x^l)\pr_\b(x^n)\psi_{nj}\\
 &=& M_{i}^{k\a}\pr_\a(\widetilde{\psi}_{kj}) +M_{j}^{k\a}\pr_\a(\widetilde{\psi}_{ik}) + \in_{ikl}x^k\g^{\a\b}\pr_\a(x^l)\pr_\b(x^n\psi_{nj})\\
&& -M_{j}^{n\a}\in_{ikl}\pr_\a(x^l)\psi_{kn} -   \in_{ikl}x^k\g^{\a\b}\pr_\a(x^l)\pr_\b(x^n)\psi_{nj}  -M_{i}^{k\a}\in_{knl}\pr_\a(x^l)\psi_{nj}\\
 &=& \frac{1}{2}S(\widetilde{\psi})_{ij} + \in_{ikl}x^k\g^{\a\b}\pr_\a(x^l)\pr_\b(x^n\psi_{nj})\\
&& -M_{j}^{n\a}\in_{ikl}\pr_\a(x^l)\psi_{kn} -   \in_{ikl}x^k\g^{\a\b}\pr_\a(x^l)\pr_\b(x^n)\psi_{nj}  -M_{i}^{k\a}\in_{knl}\pr_\a(x^l)\psi_{nj}
\eeaa
which together with \eqref{eq:waveeqinijkxlpsikj:aux0} implies
\bea\lab{eq:waveeqinijkxlpsikj:aux2}
\nn&&\square_\g(\widetilde{\psi}_{ij})\\ 
\nn&=& S(\widetilde{\psi})_{ij}  + 2\in_{ikl}x^k\g^{\a\b}\pr_\a(x^l)\pr_\b(x^n\psi_{nj}) -2M_{j}^{n\a}\in_{ikl}\pr_\a(x^l)\psi_{kn}\\
\nn&& -   2\in_{ikl}x^k\g^{\a\b}\pr_\a(x^l)\pr_\b(x^n)\psi_{nj}  -2M_{i}^{k\a}\in_{knl}\pr_\a(x^l)\psi_{nj} +\in_{ikl}x^l(Q\psi)_{kj}\\
&&+\in_{ikl}\square_\g(x^l)\psi_{kj},
\eea

Next, we consider the before to last term of the RHS in \eqref{eq:waveeqinijkxlpsikj:aux2}. We have
\beaa
\in_{ikl}x^l(Q\psi)_{kj} &=& \in_{ikl}x^l(\Ddot^\a M_{k\a}^n)\psi_{nj}+\in_{ikl}x^l(\Ddot^\a M_{j\a}^n)\psi_{kn} -\in_{ikl}x^lM_{k\a}^nM_n^{m\a}\psi_{mj}\\
&&-2\in_{ikl}x^lM_{k\a}^nM_{j}^{m\a}\psi_{nm}-\in_{ikl}x^lM_{j\a}^nM_n^{m\a}\psi_{km}\\
&=&  \Ddot^\a(\in_{ikl}x^lM_{k\a}^n)\psi_{nj} - \in_{ikl}\pr_\a(x^l)M_{k}^{n\a}\psi_{nj}
+(\Ddot^\a M_{j\a}^n)\in_{ikl}x^l\psi_{kn} \\
&&-\in_{ikl}x^lM_{k\a}^nM_n^{m\a}\psi_{mj}-2\in_{ikl}x^lM_{k\a}^nM_{j}^{m\a}\psi_{nm}-M_{j\a}^nM_n^{m\a}\in_{ikl}x^l\psi_{km}\\
&=&  \Ddot^\a(\in_{ikl}x^lM_{k\a}^n)\psi_{nj} - \in_{ikl}\pr_\a(x^l)M_{k}^{n\a}\psi_{nj}
+(\Ddot^\a M_{j\a}^k)\widetilde{\psi}_{ik}\\
&&-\in_{ikl}x^lM_{k\a}^nM_n^{m\a}\psi_{mj}-2\in_{ikl}x^lM_{k\a}^nM_{j}^{m\a}\psi_{nm}-M_{j\a}^kM_k^{l\a}\widetilde{\psi}_{il}.
\eeaa
Since we have, in view of \eqref{eq:waveeqinijkxlpsikj:aux0:00} and \eqref{eq:waveeqinijkxlpsikj:aux1:00}  
\bea\lab{eq:waveeqinijkxlpsikj:aux2:00}
\in_{ikl}x^lM_{k\a}^{n} &=& M_{i\a}^k\in_{knl}x^l - \in_{ikl}\pr_\a(x^l)(\de^{kn}-x^kx^n),
\eea
we infer
\beaa
&&\in_{ikl}x^l(Q\psi)_{kj}\\ 
&=&  \Ddot^\a\Big(M_{i\a}^k\in_{knl}x^l - \in_{ikl}\pr_\a(x^l)(\de^{kn}-x^kx^n)\Big)\psi_{nj} - \in_{ikl}\pr_\a(x^l)M_{k}^{n\a}\psi_{nj}\\
&&+(\Ddot^\a M_{j\a}^k)\widetilde{\psi}_{ik} -\Big(M_{i\a}^k\in_{knl}x^l - \in_{ikl}\pr_\a(x^l)(\de^{kn}-x^kx^n)\Big)M_n^{m\a}\psi_{mj}\\
&&-2\Big(M_{i\a}^k\in_{knl}x^l - \in_{ikl}\pr_\a(x^l)(\de^{kn}-x^kx^n)\Big)M_{j}^{m\a}\psi_{nm} -M_{j\a}^kM_k^{l\a}\widetilde{\psi}_{il}
\eeaa
i.e.,
\beaa
&&\in_{ikl}x^l(Q\psi)_{kj}\\ 
&=&  (\Ddot^\a M_{i\a}^k)\widetilde{\psi}_{kj}+M_{i\a}^k\in_{knl}\pr_\a(x^l)\psi_{nj} - \square_\g(x^l)(\in_{inl} -\in_{ikl}x^kx^n)\psi_{nj}\\
&& +\in_{ikl}\g^{\a\b}\pr_\a(x^l)\pr_\b(x^kx^n)\psi_{nj} - \in_{ikl}\pr_\a(x^l)M_{k}^{n\a}\psi_{nj}+(\Ddot^\a M_{j\a}^k)\widetilde{\psi}_{ik}\\
&& -M_{i\a}^k\in_{knl}x^lM_n^{m\a}\psi_{mj} +\pr_\a(x^l)(\in_{inl}-\in_{ikl}x^kx^n)M_n^{m\a}\psi_{mj}\\
&& -2M_{i\a}^k\in_{knl}x^lM_{j}^{m\a}\psi_{nm} +2\pr_\a(x^l)(\in_{inl} -\in_{ikl}x^kx^n)M_{j}^{m\a}\psi_{nm} -M_{j\a}^kM_k^{l\a}\widetilde{\psi}_{il}.
\eeaa
Hence, using again \eqref{eq:waveeqinijkxlpsikj:aux2:00}, we obtain 
\beaa
&&\in_{ikl}x^l(Q\psi)_{kj}\\ 
&=&  (\Ddot^\a M_{i\a}^k)\widetilde{\psi}_{kj}+M_{i\a}^k\in_{knl}\pr_\a(x^l)\psi_{nj} - \square_\g(x^l)(\in_{inl} -\in_{ikl}x^kx^n)\psi_{nj}\\
&& +\in_{ikl}\g^{\a\b}\pr_\a(x^l)\pr_\b(x^kx^n)\psi_{nj} - \in_{ikl}\pr_\a(x^l)M_{k}^{n\a}\psi_{nj}+(\Ddot^\a M_{j\a}^k)\widetilde{\psi}_{ik}\\
&& -M_{i\a}^k\Big(M_{k}^{n\a}\in_{nml}x^l - \in_{knl}\pr^\a(x^l)(\de^{nm}-x^nx^m)\Big)\psi_{mj}\\
&& +\pr_\a(x^l)(\in_{inl}-\in_{ikl}x^kx^n)M_n^{m\a}\psi_{mj} -2M_{i\a}^kM_{j}^{m\a}\widetilde{\psi}_{km}\\ 
&&+2\pr_\a(x^l)(\in_{inl} -\in_{ikl}x^kx^n)M_{j}^{m\a}\psi_{nm} -M_{j\a}^kM_k^{l\a}\widetilde{\psi}_{il},
\eeaa
or, using \eqref{eq:xiMialphakformula} as well, 
\beaa
&&\in_{ikl}x^l(Q\psi)_{kj}\\ 
&=&  (\Ddot^\a M_{i\a}^k)\widetilde{\psi}_{kj}+(\Ddot^\a M_{j\a}^k)\widetilde{\psi}_{ik}  -M_{i\a}^kM_{k}^{n\a}\widetilde{\psi}_{nj} -2M_{i\a}^kM_{j}^{m\a}\widetilde{\psi}_{km} -M_{j\a}^kM_k^{l\a}\widetilde{\psi}_{il} \\
&& +\square_\g(x^l)\in_{ikl}x^kx^n\psi_{nj}  +\in_{ikl}\g^{\a\b}\pr_\a(x^l)\pr_\b(x^k)x^n\psi_{nj}-M_{i\a}^k\in_{knl}\pr^\a(x^l)x^nx^m\psi_{mj}\\
&& -2\pr_\a(x^l)\in_{ikl}x^kM_{j}^{m\a}x^n\psi_{nm} +M_{i\a}^k\in_{knl}\pr^\a(x^l)\psi_{nj} - \square_\g(x^l)\in_{inl}\psi_{nj}\\
&&+\in_{ikl}\g^{\a\b}\pr_\a(x^l)x^k\pr_\b(x^n)\psi_{nj} - \in_{ikl}\pr_\a(x^l)M_{k}^{n\a}\psi_{nj}  +M_{i\a}^k\in_{knl}\pr^\a(x^l)\psi_{nj}\\
&& +\pr_\a(x^l)\in_{inl}M_n^{m\a}\psi_{mj} +\in_{ikl}x^k\g^{\a\b}\pr_\a(x^l)\pr_\b(x^m)\psi_{mj} +2\pr_\a(x^l)\in_{inl}M_{j}^{m\a}\psi_{nm},
\eeaa
which we rewrite as
\beaa
&&\in_{ikl}x^l(Q\psi)_{kj}\\ 
&=&  (Q\widetilde{\psi})_{ij}  +\square_\g(x^l)\in_{ikl}x^kx^n\psi_{nj}  +\in_{ikl}\g^{\a\b}\pr_\a(x^l)\pr_\b(x^k)x^n\psi_{nj}-M_{i\a}^k\in_{knl}\pr^\a(x^l)x^nx^m\psi_{mj}\\
&& -2\pr_\a(x^l)\in_{ikl}x^kM_{j}^{m\a}x^n\psi_{nm} +M_{i\a}^k\in_{knl}\pr^\a(x^l)\psi_{nj} - \square_\g(x^l)\in_{inl}\psi_{nj}\\
&&+\in_{ikl}\g^{\a\b}\pr_\a(x^l)x^k\pr_\b(x^n)\psi_{nj} - \in_{ikl}\pr_\a(x^l)M_{k}^{n\a}\psi_{nj}  +M_{i\a}^k\in_{knl}\pr^\a(x^l)\psi_{nj}\\
&& +\pr_\a(x^l)\in_{inl}M_n^{m\a}\psi_{mj} +\in_{ikl}x^k\g^{\a\b}\pr_\a(x^l)\pr_\b(x^m)\psi_{mj} +2\pr_\a(x^l)\in_{inl}M_{j}^{m\a}\psi_{nm}.
\eeaa

Plugging the last identity in the RHS of \eqref{eq:waveeqinijkxlpsikj:aux2}, we obtain 
\beaa
\nn&&\square_\g(\widetilde{\psi}_{ij})\\ 
\nn&=& S(\widetilde{\psi})_{ij} +(Q\widetilde{\psi})_{ij} + 2\in_{ikl}x^k\g^{\a\b}\pr_\a(x^l)\pr_\b(x^n\psi_{nj})  +\square_\g(x^l)\in_{ikl}x^kx^n\psi_{nj}\\
&&  +\in_{ikl}\g^{\a\b}\pr_\a(x^l)\pr_\b(x^k)x^n\psi_{nj}-M_{i\a}^k\in_{knl}\pr^\a(x^l)x^nx^m\psi_{mj}  -2\pr_\a(x^l)\in_{ikl}x^kM_{j}^{m\a}x^n\psi_{nm} \\
&&+F_{ij}
\eeaa
where we have introduced the family of scalars $F_{ij}$ given by
\beaa
F_{ij} &:=& -2M_{j}^{n\a}\in_{ikl}\pr_\a(x^l)\psi_{kn} -   2\in_{ikl}x^k\g^{\a\b}\pr_\a(x^l)\pr_\b(x^n)\psi_{nj}  -2M_{i}^{k\a}\in_{knl}\pr_\a(x^l)\psi_{nj}\\
&& +\in_{ikl}\square_\g(x^l)\psi_{kj} +M_{i\a}^k\in_{knl}\pr^\a(x^l)\psi_{nj} - \square_\g(x^l)\in_{inl}\psi_{nj}\\
&&+\in_{ikl}\g^{\a\b}\pr_\a(x^l)x^k\pr_\b(x^n)\psi_{nj} - \in_{ikl}\pr_\a(x^l)M_{k}^{n\a}\psi_{nj}  +M_{i\a}^k\in_{knl}\pr^\a(x^l)\psi_{nj}\\
&& +\pr_\a(x^l)\in_{inl}M_n^{m\a}\psi_{mj} +\in_{ikl}x^k\g^{\a\b}\pr_\a(x^l)\pr_\b(x^m)\psi_{mj} +2\pr_\a(x^l)\in_{inl}M_{j}^{m\a}\psi_{nm}.
\eeaa
Noticing that all terms in $F_{ij}$ cancel, we obtain $F_{ij} =0$ and hence
\beaa
\nn&&\square_\g(\widetilde{\psi}_{ij})\\ 
\nn&=& S(\widetilde{\psi})_{ij} +(Q\widetilde{\psi})_{ij} + 2\in_{ikl}x^k\g^{\a\b}\pr_\a(x^l)\pr_\b(x^n\psi_{nj})  +\square_\g(x^l)\in_{ikl}x^kx^n\psi_{nj}\\
&&  +\in_{ikl}\g^{\a\b}\pr_\a(x^l)\pr_\b(x^k)x^n\psi_{nj}-M_{i\a}^k\in_{knl}\pr^\a(x^l)x^nx^m\psi_{mj}  -2\pr_\a(x^l)\in_{ikl}x^kM_{j}^{m\a}x^n\psi_{nm}
\eeaa
or
\beaa
\nn&&\square_\g(\widetilde{\psi}_{ij})\\ 
\nn&=& S(\widetilde{\psi})_{ij} +(Q\widetilde{\psi})_{ij} + 2\in_{ikl}x^k\g^{\a\b}\pr_\a(x^l)\pr_\b((\err_{\textrm{TDefect},2}[\psi])_j)\\
&&  +\square_\g(x^l)\in_{ikl}x^k(\err_{\textrm{TDefect},2}[\psi])_j  +\in_{ikl}\g^{\a\b}\pr_\a(x^l)\pr_\b(x^k)(\err_{\textrm{TDefect},2}[\psi])_j\\
&& -M_{i\a}^k\in_{knl}\pr^\a(x^l)x^n(\err_{\textrm{TDefect},2}[\psi])_j  -2\pr_\a(x^l)\in_{ikl}x^kM_{j}^{m\a}(\err_{\textrm{TDefect},2}[\psi])_m.
\eeaa
Since $(\err_{\textrm{TDefect},5}[\psi])_{ij}= \in_{ikl}x^l\psi_{kj}+i\psi_{ij}=\widetilde{\psi}_{ij}+i\psi_{ij}$, we deduce 
\beaa
\nn&&\square_\g((\err_{\textrm{TDefect},5}[\psi])_{ij})\\ 
\nn&=& S(\err_{\textrm{TDefect},5}[\psi])_{ij} +(Q\err_{\textrm{TDefect},5}[\psi])_{ij} + 2\in_{ikl}x^k\g^{\a\b}\pr_\a(x^l)\pr_\b((\err_{\textrm{TDefect},2}[\psi])_j)\\ 
&& +\Big(\in_{ikl}x^k\square_\g(x^l)  +\in_{ikl}\g^{\a\b}\pr_\a(x^l)\pr_\b(x^k) -\in_{knl}x^nM_{i\a}^k\pr^\a(x^l)\Big)(\err_{\textrm{TDefect},2}[\psi])_j\\
&&  -2\in_{ikl}x^kM_{j}^{m\a}\pr_\a(x^l)(\err_{\textrm{TDefect},2}[\psi])_m.
\eeaa
as stated in \eqref{eq:waveeqforlastdefecttensor:antiselfdual}. This concludes the proof of Lemma \ref{lemma:waveequationsfortensordeffects}. 
\end{proof}


\subsection{Differentiation with respect to $\pr_\tau$ and $\widehat{\pr}_{\tphi}$ preserving identities   \eqref{eq:fundamentalidentitiestoderivefromscalarizationoftensor:complexcase}}


We start by noticing that differentiation w.r.t. $\pr_\tau$ preserves the identities  \eqref{eq:fundamentalidentitiestoderivefromscalarizationoftensor:complexcase}. 
\begin{lemma}\lab{lemma:differentiatingwrtprtaupreservetheidentitiesscaloftensors}
Let $\psi_{ij}$ be a family of complex-valued scalars satisfying the identities \eqref{eq:fundamentalidentitiestoderivefromscalarizationoftensor:complexcase}. Then, $\pr_\tau(\psi_{ij})$ satisfies the identities \eqref{eq:fundamentalidentitiestoderivefromscalarizationoftensor:complexcase} as well.
\end{lemma}

\begin{proof}
This follows immediately from the fact that $\pr_\tau(x^i)=0$, $i=1,2,3$. 
\end{proof}

While $\pr_{\tphi}(x^3)=0$, we have $\pr_{\tphi}(x^1)=-x^2$ and $\pr_{\tphi}(x^2)=x^1$. Hence, differentiation w.r.t.  $\pr_{\tphi}$ does not preserve the identities \eqref{eq:fundamentalidentitiestoderivefromscalarizationoftensor:complexcase} and we will instead use the following modification. 

\begin{definition}\lab{def:widehatprtphi}
Let $\widehat{\pr}_{\tphi}$ denote the first-order operator acting on families of complex-valued scalars $\psi_{ij}$ as follows 
\beaa
\widehat{\pr}_{\tphi}(\psi)_{ij}:=\pr_{\tphi}(\psi_{ij}) +\in_{ik3}\psi_{kj} +\in_{jk3}\psi_{ki}.
\eeaa
\end{definition}

\begin{remark}
In Kerr, if $\psi_{ij}=\pmb\psi(\Om_i, \Om_j)$ with $\pmb\psi\in\sk_2$, then we have $\widehat{\pr}_{\tphi}(\psi)_{ij}=\Lieb_{\pr_{\tphi}}\pmb\psi(\Om_i, \Om_j)$, see Lemma \ref{lemma:InKerrlinkbetweenprtauwidehatprtphiandLieb[rtauLiebprtphi}, where  the horizontal Lie derivative $\Lieb$ has been introduced in Definition \ref{definition:hor-Lie-derivative}. This motivates Definition \ref{def:widehatprtphi}.
\end{remark}

The following lemma proves that differentiation w.r.t. $\widehat{\pr}_{\tphi}$ preserves the identities \eqref{eq:fundamentalidentitiestoderivefromscalarizationoftensor:complexcase}. 
\begin{lemma}\lab{lemma:differentiatingwrtwidehatprtphipreservetheidentitiesscaloftensors}
Let $\psi_{ij}$ be a family of complex-valued scalars satisfying the identities \eqref{eq:fundamentalidentitiestoderivefromscalarizationoftensor:complexcase} and let $\widehat{\pr}_{\tphi}$ be as in Definition \ref{def:widehatprtphi}. Then, $\widehat{\pr}_{\tphi}(\psi)_{ij}$ satisfies the identities \eqref{eq:fundamentalidentitiestoderivefromscalarizationoftensor:complexcase} as well.
\end{lemma}

\begin{proof}
Since $\psi_{ij}$ satisfies the identities \eqref{eq:fundamentalidentitiestoderivefromscalarizationoftensor:complexcase}, we immediately have
\beaa
\widehat{\pr}_{\tphi}(\psi)_{ji} = \widehat{\pr}_{\tphi}(\psi)_{ij}.
\eeaa
Also, we have
\beaa
x^i\widehat{\pr}_{\tphi}(\psi)_{ij} = -\pr_{\tphi}(x^i)\psi_{ij}+x^i\in_{ik3}\psi_{kj}=x^2\psi_{1j}-x^1\psi_{2j}+x^1\psi_{2j}-x^2\psi_{1j}=0
\eeaa
and 
\beaa
(\de^{ij}-x^ix^j)\widehat{\pr}_{\tphi}(\psi)_{ij}=\de^{ij}\widehat{\pr}_{\tphi}(\psi)_{ij}=2\de^{ij}\in_{ik3}\psi_{kj}=\in_{ik3}\psi_{ki}=0
\eeaa
where we used the antisymmetry of $\in_{ik3}$ and the symmetry of $\psi_{ki}$ w.r.t. $(i,k)$. 

It remains to check the last identity of \eqref{eq:fundamentalidentitiestoderivefromscalarizationoftensor:complexcase}. We compute
\beaa
&&\in_{ikl}x^l\widehat{\pr}_{\tphi}(\psi)_{kj}+i\widehat{\pr}_{\tphi}(\psi)_{ij}\\ 
&=& \in_{ikl}x^l\Big(\pr_{\tphi}(\psi_{kj}) +\in_{kn3}\psi_{nj} +\in_{jn3}\psi_{nk}\Big)+i\Big(\pr_{\tphi}(\psi_{ij}) +\in_{il3}\psi_{lj} +\in_{jl3}\psi_{li}\Big)\\
&=& -\in_{ikl}\pr_{\tphi}(x^l)\psi_{kj}+ \in_{ikl}x^l\in_{kn3}\psi_{nj} +\in_{ikl}x^l\in_{jn3}\psi_{nk}+i\Big(\in_{il3}\psi_{lj} +\in_{jl3}\psi_{li}\Big).
\eeaa
Since 
\beaa
\in_{ikl}x^l\in_{jn3}\psi_{nk} = \in_{jn3}\big(\in_{ikl}x^l\psi_{kn}\big)=-i\in_{jn3}\psi_{in}=-i\in_{jl3}\psi_{li},
\eeaa
we infer
\beaa
\in_{ikl}x^l\widehat{\pr}_{\tphi}(\psi)_{kj}+i\widehat{\pr}_{\tphi}(\psi)_{ij} &=& -\in_{ikl}\pr_{\tphi}(x^l)\psi_{kj}+ \in_{ikl}x^l\in_{kn3}\psi_{nj} +i\in_{il3}\psi_{lj}\\
&=&  -\in_{ikl}\pr_{\tphi}(x^l)\psi_{kj}+ \in_{ikl}x^l\in_{kn3}\psi_{nj} +\in_{il3}(-\in_{lkn}x^n\psi_{kj})\\
&=&  \big(-\in_{ikl}\pr_{\tphi}(x^l)+ \in_{inl}x^l\in_{nk3} -\in_{il3}\in_{lkn}x^n\big)\psi_{kj}\\
&=& A_{ik}\psi_{kj}
\eeaa
where $A_{ik}$ is given by 
\beaa
A_{ik}:=-\in_{ikl}\pr_{\tphi}(x^l) - \in_{inl}x^l\in_{kn3} +\in_{in3}\in_{knl}x^l.
\eeaa
Now, since $A_{ik}$ is antisymmetric w.r.t. $(i,k)$ and since 
\beaa
A_{12} &=& -\in_{12l}\pr_{\tphi}(x^l) - \in_{1nl}x^l\in_{2n3} +\in_{1n3}\in_{2nl}x^l\\
&=& -\in_{123}\pr_{\tphi}(x^3) - \in_{11l}x^l\in_{213} +\in_{123}\in_{22l}x^l=0,\\
A_{13} &=& -\in_{13l}\pr_{\tphi}(x^l) - \in_{1nl}x^l\in_{3n3} +\in_{1n3}\in_{3nl}x^l\\
&=& -\in_{132}\pr_{\tphi}(x^2)  +\in_{123}\in_{32l}x^l=x^1+\in_{321}x^1=x^1-x^1=0,\\
A_{23} &=& -\in_{23l}\pr_{\tphi}(x^l) - \in_{2nl}x^l\in_{3n3} +\in_{2n3}\in_{3nl}x^l\\
&=&  -\in_{231}\pr_{\tphi}(x^1) +\in_{213}\in_{31l}x^l=x^2 - \in_{312}x^2=x^2-x^2=0,
\eeaa
we infer that $A_{ik}=0$ for all $i,k=1,2,3$ and hence 
\beaa
\in_{ikl}x^l\widehat{\pr}_{\tphi}(\psi)_{kj}+i\widehat{\pr}_{\tphi}(\psi)_{ij} = 0.
\eeaa
Thus, $\widehat{\pr}_{\tphi}(\psi)_{ij}$ satisfies the identities  \eqref{eq:fundamentalidentitiestoderivefromscalarizationoftensor:complexcase} as stated. This concludes the proof of Lemma \ref{lemma:differentiatingwrtwidehatprtphipreservetheidentitiesscaloftensors}.
\end{proof}
 
\begin{remark}
\lab{rem:prtauandprtphihatpreservesk2C}
In view of Lemmas \ref{lemma:backandforthbetweenhorizontaltensorsk2andscalars:complex},  \ref{lemma:differentiatingwrtprtaupreservetheidentitiesscaloftensors} and   \ref{lemma:differentiatingwrtwidehatprtphipreservetheidentitiesscaloftensors}, we immediately infer the fact that 
if $\psi_{ij}=\pmb\psi(\Om_i, \Om_j)$ for $\pmb\psi\in \sk_2(\mathbb{C})$, then for any $k,l\in\mathbb{N}$, there exists $\pmb\psi_{(k,l)}\in\sk_2(\mathbb{C})$ such that $\pr_\tau^k\widehat{\pr}_{\tphi}^l(\psi)_{ij}=\pmb\psi_{(k,l)}(\Om_i,\Om_j)$.
\end{remark}


\subsection{Regular triplet in Kerr}
\lab{sec:regulartripletsinKerr}


\begin{definition}\lab{def:regulartripletinKerrOmii=123}
Let 
\bea
x^1:=\cos\phi\sin\th, \qquad x^2:=\sin\phi\sin\th, \qquad x^3:=\cos\th.
\eea
Then, we define the following horizontal vectorfields $\Om^i$ in Kerr by 
\bea
\Om^i :=|q|\dual\nab(x^i), \quad i=1,2,3.
\eea
\end{definition}

\begin{lemma}\lab{lemma:fundamentalpropertiesof1formsOmi}
The horizontal vectorfields $\Om_i$ in Kerr introduced in Definition \ref{def:regulartripletinKerrOmii=123} satisfy \eqref{eq:fundamentalpropertiesof1formsOmi}. In particular, they form a regular triplet in Kerr in the sense of Definition \ref{def:definitionofregulartripletOmii=123}.
\end{lemma}

\begin{proof}
We start with the first identity in \eqref{eq:fundamentalpropertiesof1formsOmi}. Since $(x^1)^2+(x^2)^2+(x^3)^2=1$, we have
\beaa
x^i\Om_i=\sum_{i=1}^3|q|x^i\dual\nab(x^i)=\frac{1}{2}|q|\dual\nab\left(\sum_{i=1}^3(x^i)^2\right)=\frac{1}{2}|q|\dual\nab(1)=0
\eeaa
as stated.

Next, we consider the second identity in \eqref{eq:fundamentalpropertiesof1formsOmi}. We have
\bsub
\label{eq:derionxi}
\begin{align}
|q|e_1(x^1)=& \cos\th\cos\phi, & |q|e_1(x^2)=& \cos\th\sin\phi, & |q|e_1(x^3)=& -\sin\th,\\
|q|e_2(x^1)=& -\sin\phi, & |q|e_2(x^2)=& \cos\phi,& |q|e_2(x^3)=& 0,
\end{align} 
\esub
which yields 
\beaa
(\Om^j)_1(\Om_j)_1 &=& \sum_{j=1}^3|q|^2\dual\nab_1(x^j)\dual\nab_1(x^j)=\sum_{j=1}^3(|q|e_2(x^j))^2=1,\\
(\Om^j)_1(\Om_j)_2 &=& \sum_{j=1}^3|q|^2 \dual\nab_1(x^j)\dual\nab_2(x^j)=-\sum_{j=1}^3|q|^2e_1(x^j)e_2(x^j)=0,\\
(\Om^j)_2(\Om_j)_2 &=& \sum_{j=1}^3|q|^2 \dual\nab_2(x^j)\dual\nab_2(x^j)=\sum_{j=1}^3(|q|e_1(x^j))^2=1,
\eeaa
and hence $(\Om^i)_a(\Om_i)_b=\de_{ab}$ as stated.

Next, we consider the third identity in \eqref{eq:fundamentalpropertiesof1formsOmi}. We have
\beaa
\Om_i\c\Om_j=|q|^2\dual\nab(x^i)\c\dual\nab(x^j)=|q|^2\nab(x^i)\c\nab(x^j).
\eeaa
Together with \eqref{eq:derionxi}, we infer
\beaa
\bsplit
\Om_1\c\Om_1&=1-(\sin\th)^2(\cos\phi)^2=1-(x^1)^2, \qquad \Om_1\c\Om_2=-(\sin\th)^2\cos\phi\sin\phi=-x^1x^2,\\
\Om_1\c\Om_3&=-\cos\th\sin\th\cos\phi=-x^1x^3, \qquad \Om_2\c\Om_2=1-(\sin\th)^2(\sin\phi)^2=1-(x^2)^2, \\
\Om_2\c\Om_3&=-\cos\th\sin\th\sin\phi=-x^2x^3, \qquad \Om_3\c\Om_3=1-(\cos\th)^2=1-(x^3)^2,
\end{split}
\eeaa
so that $\Om_i\c\Om_j=\de_{ij}-x^ix^j$ as stated. 

Finally, we consider the fourth identity in \eqref{eq:fundamentalpropertiesof1formsOmi}. We have
\beaa
\Om_i\c\dual\Om_j=|q|^2\dual\nab(x^i)\c\dual\dual\nab(x^j)=-|q|^2\dual\nab(x^i)\c\nab(x^j)=|q|^2\nab(x^i)\c\dual\nab(x^j).
\eeaa
Together with \eqref{eq:derionxi}, we infer
\beaa
\bsplit
\Om_i\c\dual\Om_j&=-\Om_j\c\dual\Om_i, \qquad  \Om_1\c\dual\Om_2=\cos\th=x^3\\
\Om_1\c\dual\Om_3&=-\sin\th\sin\phi=-x^2,\qquad \Om_2\c\dual\Om_3=\sin\th\cos\phi=x^1
\end{split}
\eeaa
so that $\Om_i\c\dual\Om_j=\in_{ijk}x^k$ as stated. This concludes the proof of Lemma \ref{lemma:fundamentalpropertiesof1formsOmi}.
\end{proof}

We also derive the following properties of the 1-forms $M_{i\a}^j$ in Kerr.
\begin{lemma}\lab{lemma:computationoftheMialphajinKerr}
Let $\Om_i$, $i=1,2,3$, be the regular triplet in Kerr of Definition \ref{def:regulartripletinKerrOmii=123}, and let $M_{i\a}^j$ be the corresponding 1-forms in Kerr given by \eqref{eq:definitionofMalphaijwithoutambiguity}. Then, we have, for $i,j=1,2,3$,  
\beaa
M_{i3}^j = \frac{a\cos\th}{|q|^2}\in_{ijk}x^k,\quad M_{i4}^j = \frac{a\cos\th\De}{|q|^4}\in_{ijk}x^k, \quad M_{i\a}^j(\pr_\tau)^\a=-\frac{2amr\cos\th}{|q|^4}\in_{ijk}x^k,
\eeaa
and the following asymptotic holds, for $r$ large and $i,j=1,2,3$,
\beaa
M_{ia}^j=O(r^{-1}), \quad a=1,2.
\eeaa
Also, denoting by $\pr^{\textrm{BL}}_r$ the corresponding Boyer-Lindquist coordinate vectorfield, we have
\beaa
M_{i\a}^j(\pr^{\textrm{BL}}_r)^\a=0,
\eeaa
which implies in particular for $r\in [r_+(1+2\dbl),  12m]$
\beaa
M_{i\a}^j(\pr_r)^\a=0,\qquad \g^{r\a}M_{i\a}^j=0, \quad i,j=1,2,3.
\eeaa
\end{lemma}

\begin{proof}
The asymptotic for $M_{ia}^j$, $a=1,2$ is immediate. Next, we focus on the computation of $M_{i\a}^j(\pr_\tau)^\a$. In view of Lemma 4.3.6 and (C.5.3) in \cite{GKS22}, we have for a scalar function $f \in \sk_0$ the following commutator in Kerr, where $\pr_\tau$ is Killing,
\beaa
[\nab_{\pr_\tau}, |q|\nab_a]f=\frac{2amr\cos\th}{|q|^4}|q|\dual\nab_a(f)
\eeaa
and hence, taking the dual, we obtain 
\beaa
[\nab_{\pr_\tau}, |q|\dual\nab]f=-\frac{2amr\cos\th}{|q|^4}|q|\nab(f)
\eeaa
which we apply to $f=x^i$, $i=1,2,3$, to infer
\beaa
\nab_{\pr_\tau}\Om_i &=& \nab_{\pr_\tau}|q|\dual\nab(x^i)=|q|\dual\nab(\pr_\tau(x^i))+[\nab_{\pr_\tau}, |q|\dual\nab]x^i\\
&=& -\frac{2amr\cos\th}{|q|^4}|q|\nab(x^i)=\frac{2amr\cos\th}{|q|^4}\dual\Om_i.
\eeaa
This yields, in view of \eqref{eq:definitionofMalphaijwithoutambiguity},
\beaa
M_{i\a}^j(\pr_\tau)^\a=(\nab_{\pr_\tau}\Om_i\c\Om^j)=\frac{2amr\cos\th}{|q|^4}(\dual\Om_i\c\Om^j)=-\frac{2amr\cos\th}{|q|^4}\in_{ijk}x^k
\eeaa
as stated. 

Next, we focus on the computation of $M_{i3}^j$ and $M_{i4}^j$. According to the commutation formulas \eqref{eq:comm-nab3-nab4-naba-f-general}, we have for $f \in \sk_0$
       \beaa
       \begin{split}
        \,[\nab_3, \nab_a] f &=-\frac 1 2 \left(\trchb \nab_a f+\atrchb \dual \nab_a f\right),\\
         \,[\nab_4, \nab_a] f &=-\frac 1 2 \left(\trch \nab_a f+\atrch \dual \nab_a f\right)+(\etab_a+\eta_a) \nab_4 f,
         \end{split}
       \eeaa
where we have used the fact that we have in Kerr $\chih=\chibh=\xi=\xib=0$, as well as $\ze=\eta$ in the regular frame of Kerr. This yields 
       \beaa
       \begin{split}
        \,[\nab_3, |q|\nab_a] f =&-\frac 1 2 \left(\left(\trchb -\frac{e_3(|q|^2)}{|q|^2}\right)|q|\nab_a f+\atrchb |q|\dual \nab_a f\right),\\
         \,[|q|^2\nab_4, |q|\nab_a] f =&-\frac{|q|^2}{2}\left(\left(\trch -\frac{e_4(|q|^2)}{|q|^2}\right)|q|\nab_a f+\atrch |q|\dual \nab_a f\right)\\
         &+|q|^3\left(\etab_a+\eta_a -\frac{\nab_a(|q|^2)}{|q|^2}\right) \nab_4 f.
         \end{split}
       \eeaa
In view of the following consequence of \eqref{eq:KerrvaluesofcomplexifiedRicci}
\beaa
\bsplit
\trch&=\frac{2\De r}{|q|^4},\quad \atrch=\frac{2a \De \cos\th}{|q|^4}, \quad  \trchb=-\frac{2r}{|q|^2}, \quad \atrchb=\frac{2a\cos\th}{|q|^2},\\
\eta_1&= -\frac{a^2\sin\th \cos\th}{|q|^3}, \qquad  \eta_2=\frac{a r \sin\th }{|q|^3}, \qquad \etab_1= -\frac{a^2\sin\th \cos\th }{|q|^3}, \qquad \etab_2 =-\frac{ar\sin\th }{|q|^3},
\end{split}
\eeaa
and the fact that, in view of \eqref{eq:actionofingoingprincipalnullframeonnormalizedcoordinates},  
\beaa
e_3(|q|^2)=-2r, \qquad e_4(|q|^2)=\frac{2r\De}{|q|^2}, \qquad \nab_1(|q|^2)=-\frac{2a^2\sin\th\cos\th}{|q|}, \qquad \nab_2(|q|^2)=0,
\eeaa
we infer 
       \beaa
       \,[\nab_3, |q|\nab_a] f = - \frac{a\cos\th}{|q|^2} |q|\dual \nab_a f,\qquad   \,[|q|^2\nab_4, |q|\nab_a] f = - \frac{a \De \cos\th}{|q|^2} |q|\dual \nab_a f,
       \eeaa
and hence, taking the dual, 
       \beaa
       \,[\nab_3, |q|\dual\nab_a] f = \frac{a\cos\th}{|q|^2} |q|\nab_a f,\qquad   \,[|q|^2\nab_4, |q|\dual\nab_a] f = \frac{a \De \cos\th}{|q|^2} |q|\nab_a f.
       \eeaa
Applying the above commutators to $f=x^i$, $i=1,2,3$, and using Definition \ref{def:regulartripletinKerrOmii=123}, we infer
\beaa
(\nab_3\Om_i)_a &=& \nab_3|q|\dual\nab_a(x^i)=|q|\dual\nab_a(e_3(x^i))+\,[\nab_3, |q|\dual\nab_a]x^i\\
&=& \frac{a\cos\th}{|q|^2} |q|\nab_a(x^i)=-\frac{a\cos\th}{|q|^2}(\dual\Om_i)_a,
\eeaa
and
\beaa
(|q|^2\nab_4\Om_i)_a &=& |q|^2\nab_4|q|\dual\nab_a(x^i)=|q|\dual\nab_a(|q|^2e_4(x^i))+\,[|q|^2\nab_4, |q|\dual\nab_a]x^i\\
&=& \frac{a \De \cos\th}{|q|^2} |q|\nab_a(x^i)=-\frac{a \De \cos\th}{|q|^2}(\dual\Om_i)_a,
\eeaa
where we have used the fact that $\nab(e_3(x^i))=0$ and $\nab(|q|^2e_4(x^i))=0$ for $i=1,2,3$ in view of the definition of $x^i$ and \eqref{eq:actionofingoingprincipalnullframeonnormalizedcoordinates}. In view of \eqref{eq:definitionofMalphaijwithoutambiguity}, we deduce 
\beaa
M_{i3}^j &=& (\nab_3\Om_i)\c\Om^j = -\frac{a\cos\th}{|q|^2}(\dual\Om_i\c\Om^j)= \frac{a\cos\th}{|q|^2}\in_{ijk}x^k,\\
M_{i4}^j &=& (\nab_4\Om_i)\c\Om^j = -\frac{a\cos\th\De}{|q|^4}(\dual\Om_i\c\Om^j)= \frac{a\cos\th\De}{|q|^4}\in_{ijk}x^k,
\eeaa
as stated. 

Finally, the Boyer-Lindquist coordinate vectorfield $\pr^{\textrm{BL}}_r$ is given by 
\beaa
\pr^{\textrm{BL}}_r=\frac{1}{2}\left(\frac{|q|^2}{\De}e_4 - e_3\right),
\eeaa
and hence 
\beaa
M_{i\a}^j(\pr^{\textrm{BL}}_r)^\a &=&  \frac{1}{2}M_{i\a}^j\left(\frac{|q|^2}{\De}e_4 - e_3\right)^\a\\
&=& \frac{|q|^2}{2\De}\left(M_{i4}^j -\frac{\De}{|q|^2}M_{i3}^j\right)=0
\eeaa
as stated. In particular, since the normalized coordinates coincide with the Boyer-Lindquist coordinates for $r\in [r_+(1+2\dbl),  12m]$, we have 
\beaa
\pr_r=\pr^{\textrm{BL}}_r, \qquad \g^{r\a}=0\quad\textrm{for}\quad x^\a\neq r, \qquad \forall r\in [r_+(1+2\dbl),  12m], 
\eeaa
and hence, for $r\in [r_+(1+2\dbl),  12m]$,  
\beaa
M_{i\a}^j(\pr_r)^\a=M_{i\a}^j(\pr^{\textrm{BL}}_r)^\a=0, \qquad \g^{r\a}M_{i\a}^j=\g^{rr}M_{ir}^j=\g^{rr}M_{i\a}^j(\pr_r)^\a=0
\eeaa
as stated. This concludes the proof of Lemma \ref{lemma:computationoftheMialphajinKerr}. 
\end{proof}

\begin{lemma}\lab{lemma:InKerrlinkbetweenprtauwidehatprtphiandLieb[rtauLiebprtphi}
In Kerr, if $\psi_{ij}=\pmb\psi(\Om_i, \Om_j)$ with $\pmb\psi\in\sk_2$, then 
\beaa
\Lieb_{\pr_\tau}\Om_i=0, \qquad \Lieb_{\pr_{\tphi}}\Om_i = -\in_{ij3}\Om_j,
\eeaa
and
\beaa
\pr_\tau(\psi_{ij})=\Lieb_{\pr_\tau}\pmb\psi(\Om_i, \Om_j), \qquad \widehat{\pr}_{\tphi}(\psi)_{ij}=\Lieb_{\pr_{\tphi}}\pmb\psi(\Om_i, \Om_j),
\eeaa
where the horizontal Lie derivative $\Lieb$ has been introduced in Definition \ref{definition:hor-Lie-derivative}, and where $\widehat{\pr}_{\tphi}$ has been introduced in Definition \ref{def:widehatprtphi}.
\end{lemma}

\begin{proof}
In view of Definition \ref{def:regulartripletinKerrOmii=123}, we have
\beaa
\Lieb_{\pr_\tau}\Om_i &=& \Lieb_{\pr_\tau}|q|\dual\nab(x^i)=\dual [\Lieb_{\pr_\tau},|q|\nab](x^i)=0,\\
\Lieb_{\pr_{\tphi}}\Om_i &=& \Lieb_{\pr_{\tphi}}|q|\dual\nab(x^i)=|q|\dual\nab(\pr_{\tphi}(x^i))+\dual [\Lieb_{\pr_{\tphi}},|q|\nab](x^i)\\
&=& -\in_{ij3}|q|\dual\nab(\pr_{\tphi}(x^j))= -\in_{ij3}\Om_j,
\eeaa
as stated, where we used the fact that $[\Lieb_{\pr_\tau},|q|\nab]=0$, $[\Lieb_{\pr_{\tphi}},|q|\nab]=0$, $\pr_\tau(x^i)=0$ and $\pr_{\tphi}(x^i)=-\in_{ij3}x^j$. We infer
\beaa
\pr_\tau(\psi_{ij}) = \pr_\tau\big(\pmb\psi(\Om_i, \Om_j)\big)=\Lieb_{\pr_\tau}\pmb\psi(\Om_i, \Om_j)+\pmb\psi(\Lieb_{\pr_\tau}\Om_i, \Om_j)+\pmb\psi(\Om_i, \Lieb_{\pr_\tau}\Om_j)=\Lieb_{\pr_\tau}\pmb\psi(\Om_i, \Om_j)
\eeaa
and
\beaa
\pr_{\tphi}(\psi_{ij}) &=& \pr_{\tphi}\big(\pmb\psi(\Om_i, \Om_j)\big)=\Lieb_{\pr_{\tphi}}\pmb\psi(\Om_i, \Om_j)+\pmb\psi(\Lieb_{\pr_{\tphi}}\Om_i, \Om_j)+\pmb\psi(\Om_i, \Lieb_{\pr_{\tphi}}\Om_j)\\
&=& \Lieb_{\pr_{\tphi}}\pmb\psi(\Om_i, \Om_j)-\in_{ik3}\psi_{kj}-\in_{jk3}\psi_{ik},
\eeaa
and hence
\beaa
\pr_\tau(\psi_{ij})=\Lieb_{\pr_\tau}\pmb\psi(\Om_i, \Om_j), \qquad \widehat{\pr}_{\tphi}(\psi)_{ij}=\Lieb_{\pr_{\tphi}}\pmb\psi(\Om_i, \Om_j),
\eeaa
as stated. This concludes the proof of Lemma \ref{lemma:InKerrlinkbetweenprtauwidehatprtphiandLieb[rtauLiebprtphi}.
\end{proof}


\section{Teukolsky equations in Kerr}
\lab{sec:TeuinKerr}


In this section, we recall the form of Teukolsky equations in Kerr. We first recall the classical formulation in the Newman-Penrose (NP) formalism in terms of complex-valued scalars. We then show how to go from the NP formalism to tensorial equations and deduce the form of Teukolsky equations in tensorial form. Finally, relying on Section  \ref{sec:regularscalarization}, we provide Teukolsky equations in a regular scalarized form.


\subsection{Teukolsky equations in Kerr using Newman-Penrose formalism}
\lab{sec:Teukoslkytransportwavesystemwithcomplexscalars}


We recall in this section the Teukolsky equations in NP formalism. Let $(e_3, e_4, e_1, e_2)$ be defined as in \eqref{def:e3e4inKerr} and \eqref{def:e1e2inKerr}. Denoting by $\W$ the linearized Weyl curvature tensor, we introduce the complex-valued scalars $\phi_{\pm 2,\text{NP}}$ as follows 	
\bea\lab{def:phiplusminus2:complexscalars:Kerr}
\bsplit
\phi_{+2,\text{NP}} &:= \frac{1}{2}|q|^2\ov{q}^2\W(e_4, e_1+ie_2, e_4, e_1+ie_2),\\  
\phi_{-2,\text{NP}} &:= \frac{1}{2}|q|^{-2}\ov{q}^2\W(e_3, e_1-ie_2, e_3, e_1-ie_2).
\end{split}
 \eea
 In a Kerr spacetime, the scalars $\{\phi_{s,\text{NP}}\}_{s=\pm 2}$ defined in \eqref{def:phiplusminus2:complexscalars:Kerr} solve the following Teukolsky equation in the Boyer--Lindquist coordinates, see \cite{Teu72} (see also equation (22) in \cite{Ma}):
\bea\label{eq:TME}
|q|^2\square_{s}\phi_{s,\text{NP}}  -  s\phi_{s,\text{NP}} 
+2s[(r-m)e_3-2r\partial_t]\phi_{s,\text{NP}}=0,
\eea
where the operator $\square_s$, $s=\pm 2$, is a spin-weighted wave operator defined by
\bea\lab{eq:spinweightedwaveoperator}
|q|^2\square_{s} := |q|^2 \square + 2is\bigg(\frac{\cos\th}{\sin^2\th}\pr_{\phi} -a\cos\th\pr_{t}\bigg)-s^2\cot^2\th.
\eea


\subsubsection{Teukolsky wave/transport system in \cite{Ma}}


In this section, we recall the wave/transport system derived in \cite{Ma}. Based on $\phi_{\pm 2,\text{NP}}$ introduced in \eqref{def:phiplusminus2:complexscalars:Kerr}, we define 
\bsub\lab{def:ComplexTeuScalars:wavesystem:Kerr:Ma17}
 \bea
\bsplit
\dot{\phi}_{+2,\text{NP}}^{(0)}:=\frac{1}{r^4}\phi_{+2,\text{NP}},\qquad \dot{\phi}_{+2,\text{NP}}^{(p+1)}:= \big({r e_3 r}\big)\dot{\phi}_{+2,\text{NP}}^{(p)}, \quad p=0,1,
\end{split}
\eea
and\footnote{Note that $\dot{\phi}_{-2,\text{NP}}^{(1)}$ in \eqref{def:ComplexTeuScalars:wavesystem:Kerr:Ma17:caseminus2} has the opposite sign convention compared to the corresponding quantity in \cite[Equation (22)]{Ma}. Hence, the equations \eqref{def:ComplexTeuScalars:wavesystem:Kerr:LspasinMa17:minus2case} for $\dot{L}^{(p)}_{-2,\text{NP}}[\dot{\phi}_{-2,\text{NP}}]$ have some sign discrepancies with the corresponding ones in  \cite[Equations (25)]{Ma}.}
\bea\lab{def:ComplexTeuScalars:wavesystem:Kerr:Ma17:caseminus2}
\bsplit
\dot{\phi}_{-2,\text{NP}}^{(0)}:=\frac{\De^2}{r^4}\phi_{-2,\text{NP}},\qquad \dot{\phi}_{-2,\text{NP}}^{(p+1)}:= \left({r\frac{|q|^2}{\Delta}e_4 r}\right)\dot{\phi}_{-2,\text{NP}}^{(p)}, \quad p=0,1.
\end{split}
\eea
\esub
Then, one deduces from \eqref{eq:TME} the following spin-weighted wave equations for $\dot{\phi}_{s,\text{NP}}^{(p)}$, see \cite{Ma}, 
\bea\lab{eq:TeukolskywavetransportsystemasinMa17}
\left(\square_{s} - \frac{2\big(2r^2 -4mr+4a^2 -\de_{p0}(r^2 -6mr+6a^2)\big)}{|q|^2r^2}\right)\dot{\phi}_{s,\text{NP}}^{(p)}=\dot{L}^{(p)}_{s,NP}[\dot{\phi}_{s,\text{NP}}], 
\eea
where $\dot{L}^{(p)}_{s,NP}[\dot{\phi}_{s,\text{NP}}]$ are given, for $s=\pm 2$ and $p=0,1,2$, by 
\bsub\lab{def:ComplexTeuScalars:wavesystem:Kerr:LspasinMa17}
\bea
\bsplit
|q|^2\dot{L}^{(0)}_{+2,\text{NP}}[\dot{\phi}_{+2,\text{NP}}] :=& \frac{4(r^2-3mr+2a^2)}{r^3}\dot{\phi}_{+2,\text{NP}}^{(1)}-\frac{8(a^2\pr_t+a\pr_{\phi})}{r}\dot{\phi}_{+2,\text{NP}}^{(0)},\\
|q|^2\dot{L}^{(1)}_{+2,\text{NP}}[\dot{\phi}_{+2,\text{NP}}] :=& \frac{2(r^2-3mr+2a^2)}{r^3}\dot{\phi}_{+2,\text{NP}}^{(2)} +\frac{6mr-12a^2}{r}\dot{\phi}_{+2,\text{NP}}^{(0)}\\
& -\frac{4(a^2\pr_t+a\pr_{\phi})}{r}\dot{\phi}_{+2,\text{NP}}^{(1)} -6(a^2\pr_t+a\pr_\phi)\dot{\phi}_{+2,\text{NP}}^{(0)},\\
|q|^2\dot{L}^{(2)}_{+2,\text{NP}}[\dot{\phi}_{+2,\text{NP}}] :=& -8(a^2\pr_t+a\pr_\phi)\dot{\phi}_{+2,\text{NP}}^{(1)} -12a^2\dot{\phi}_{+2,\text{NP}}^{(0)},
\end{split}
\eea
and
\bea\lab{def:ComplexTeuScalars:wavesystem:Kerr:LspasinMa17:minus2case}
\bsplit
|q|^2\dot{L}^{(0)}_{-2,\text{NP}}[\dot{\phi}_{-2,\text{NP}}] :=& -\frac{4(r^2-3mr+2a^2)}{r^3}\dot{\phi}_{-2,\text{NP}}^{(1)}+\frac{8(a^2\pr_t+a\pr_{\phi})}{r}\dot{\phi}_{-2,\text{NP}}^{(0)},\\
|q|^2\dot{L}^{(1)}_{-2,\text{NP}}[\dot{\phi}_{-2,\text{NP}}] :=& -\frac{2(r^2-3mr+2a^2)}{r^3}\dot{\phi}_{-2,\text{NP}}^{(2)} -\frac{6mr-12a^2}{r}\dot{\phi}_{-2,\text{NP}}^{(0)}\\
& +\frac{4(a^2\pr_t+a\pr_{\phi})}{r}\dot{\phi}_{-2,\text{NP}}^{(1)} -6(a^2\pr_t+a\pr_\phi)\dot{\phi}_{-2,\text{NP}}^{(0)},\\
|q|^2\dot{L}^{(2)}_{-2,\text{NP}}[\dot{\phi}_{-2,\text{NP}}] :=& -8(a^2\pr_t+a\pr_\phi)\dot{\phi}_{-2,\text{NP}}^{(1)} -12a^2\dot{\phi}_{-2,\text{NP}}^{(0)}.
\end{split}
\eea
\esub


\subsubsection{Teukolsky wave/transport system with a different normalization}


The form of the Teukolsky  wave/transport system recalled in \eqref{def:ComplexTeuScalars:wavesystem:Kerr:Ma17} \eqref{eq:TeukolskywavetransportsystemasinMa17} is not suitable for the purposes of the paper, and we would like to trade the derivative $\pr_{\phi}+a\pr_t$ appearing in the term $(a^2\pr_t+a\pr_\phi)\dot{\phi}_{s,\text{NP}}^{(p)}$ on the RHS of the definition of $\dot{L}^{(p)}_{s,\text{NP}}[\dot{\phi}_{s,\text{NP}}]$, see \eqref{def:ComplexTeuScalars:wavesystem:Kerr:LspasinMa17}, with a horizontal derivative, see Remark \ref{rmk:explainpointofintroducingnewNPcomplexscalarphisp} below. To this end, we introduce the new complex-valued scalars $\phi_{s,\text{NP}}^{(p)}$ defined as 
\bsub\lab{def:ComplexTeuScalars:wavesystem:Kerr}
 \bea
\phi_{+2,\text{NP}}^{(0)}:= {\frac{1}{|q|^4}}\phi_{+2,\text{NP}},\qquad \phi_{+2,\text{NP}}^{(p+1)}:=\left({\frac{r^2}{|q|^2}}\right)^{1-p} \big({r e_3 r}\big)\left({\frac{r^2}{|q|^2}}\right)^{p-2} \phi_{+2,\text{NP}}^{(p)}, \quad p=0,1,
\eea
and
\bea
\phi_{-2,\text{NP}}^{(0)}:=\frac{\De^2}{|q|^4}\phi_{-2,\text{NP}},\quad\,\,\, \phi_{-2,\text{NP}}^{(p+1)}:=\left(\frac{{r^2}}{|q|^2}\right)^{1-p}\left({r\frac{|q|^2}{\Delta}e_4 r}\right)\left(\frac{{r^2}}{|q|^2}\right)^{p-2}\phi_{-2,\text{NP}}^{(p)}, \,\,\,\, p=0,1.
\eea
\esub

Comparing the definitions \eqref{def:ComplexTeuScalars:wavesystem:Kerr:Ma17} and \eqref{def:ComplexTeuScalars:wavesystem:Kerr}, one immediately checks that 
the complex-valued scalars $\phi_{s,\text{NP}}^{(p)}$ are related to $\dot{\phi}_{s,\text{NP}}^{(p)}$ by
\bea\lab{eq:relationbetweencomplexscalarsinMa17andnewnormalization}
\phi_{s,\text{NP}}^{(p)}=\left(\frac{r^2}{|q|^2}\right)^{2-p}\dot{\phi}_{s,\text{NP}}^{(p)}, \quad s=\pm 2, \quad p=0,1,2.
\eea
The following lemma provides the structure of the system of wave equations  for $\phi_{s,\text{NP}}^{(p)}$. 
\begin{lemma}\lab{lemma:structureofwavesystemforphisp:newNPform}
Let $\phi_{s,\text{NP}}^{(p)}$ be the complex-valued scalars  defined as in \eqref{def:ComplexTeuScalars:wavesystem:Kerr}. Then, $\phi_{s,\text{NP}}^{(p)}$ satisfy 
the following spin-weighted wave equations 
\bea\lab{eq:Teukolskywavetransportsystem:NPnewform}
\left(\square_{s}  - \frac{4-2\de_{p0}}{|q|^2}\right)\phi_{s,\text{NP}}^{(p)}=L^{(p)}_{s,NP}[\phi_{s,\text{NP}}], \quad s=\pm 2, \quad p=0,1,2,
\eea
where $L^{(p)}_{s,\text{NP}}[\phi_{s,\text{NP}}]$ are given by  
\bea\lab{eq:Teukolskywavetransportsystem-RHSLpsofpsips:NPnewform}
\bsplit
L^{(0)}_{s,\text{NP}}[\phi_{s,\text{NP}}] :=& \frac{2s}{r^3}(1+O(mr^{-1}))\phi_{s,\text{NP}}^{(1)}+O(mr^{-3})(\mathcal{X}_s)^{\leq 1}\phi_{s,\text{NP}}^{(0)},\\
L^{(1)}_{s,\text{NP}}[\phi_{s,\text{NP}}] :=& \frac{s}{r^3}(1+O(mr^{-1}))\phi_{s,\text{NP}}^{(2)}  +O(mr^{-3})(\mathcal{X}_s)^{\leq 1}\phi_{s,\text{NP}}^{(1)}\\
&+O(mr^{-2})(\pr_{\phi}+a\pr_t)^{\leq 1}\phi_{s,\text{NP}}^{(0)},\\
L^{(2)}_{s,\text{NP}}[\phi_{s,\text{NP}}] :=& O(mr^{-3})\phi_{s,\text{NP}}^{(2)}+O(mr^{-2})(\pr_\phi+a\pr_t)\phi_{s,\text{NP}}^{(1)} +O(m^2r^{-2})\phi_{s,\text{NP}}^{(0)},
\end{split}
\eea
with all the coefficients on  the RHS of \eqref{eq:Teukolskywavetransportsystem-RHSLpsofpsips:NPnewform} being real functions\footnote{In fact, all we need in later discussions is that the coefficients in front of the terms $\phi_{s,\text{NP}}^{(2)}$ and $(\pr_\phi+a\pr_t)\phi_{s,\text{NP}}^{(1)}$ on the RHS of expression of $L^{(2)}_{s,\text{NP}}[\phi_{s,\text{NP}}]$ are real functions.} independent of coordinates $t$ and 
$\phi$, and where $\mathcal{X}_s$ are regular horizontal vectorfields of the form\footnote{More precisely, $\mathcal{X}_s = s(\pr_\phi+a(\sin\th)^2\pr_t)-\frac{2a\cos\th}{r}\sin\th\pr_\th$, but this explicit expression will not be needed.}
\bea\lab{eq:formofregularhorizontalvectorfieldmathcalXs}
\mathcal{X}_s = s(\pr_\phi+a(\sin\th)^2\pr_t)+O(ar^{-1})\sin\th\pr_\th.
\eea
\end{lemma}

\begin{remark}\lab{rmk:explainpointofintroducingnewNPcomplexscalarphisp}
As explained above, the point of introducing the new complex-valued scalars $\phi_{s,{\text{NP}}}^{(p)}$ is to trade the derivative $\pr_{\phi}+a\pr_t$ appearing in the term $(a^2\pr_t+a\pr_\phi)\dot{\phi}_{s,\text{NP}}^{(p)}$ on the RHS of the definition of $\dot{L}^{(p)}_{s,\text{NP}}[\dot{\phi}_{s,\text{NP}}]$, see \eqref{def:ComplexTeuScalars:wavesystem:Kerr:LspasinMa17}, with the horizontal derivative $\mathcal{X}_s$ appearing on the RHS of the definition of $L^{(p)}_{s,\text{NP}}[\dot{\phi}_{s,\text{NP}}]$, see \eqref{eq:Teukolskywavetransportsystem-RHSLpsofpsips:NPnewform}.
\end{remark}

\begin{remark}
The equations \eqref{eq:Teukolskywavetransportsystem:NPnewform}-\eqref{eq:Teukolskywavetransportsystem-RHSLpsofpsips:NPnewform} and \eqref{def:ComplexTeuScalars:wavesystem:Kerr} are the Teukolsky wave and transport equations in Kerr in NP formalism, respectively.
\end{remark}

\begin{proof}
Let $\psi$ be a scalar function and $h=h(r, \cos\th)$, then, recalling \eqref{eq:spinweightedwaveoperator}, we have
\beaa
\square_s(h\psi) &=& h\square_s\psi+[\square, h]\psi = h\square_s\psi+2\g^{\a\b}e_\a(h)e_\b(\psi)+\square(h)\psi.
\eeaa
Since we have in Kerr
\beaa
e_4(r)=\frac{\De}{|q|^2}, \qquad e_3(r)=-1, \qquad \nab(r)=0, \qquad e_4(\th)=e_3(\th)=0, \qquad \nab(\cos\th)=-\dual\Re(\Jk),  
\eeaa
we infer
\beaa
2\gam^{\a\b}e_\a(h)e_\b(\psi) &=& -e_3(h)e_4(\psi) -e_4(h)e_3(\psi)+2\nab(h)\c\nab(\psi)\\
&=& \pr_r(h)\left(e_4 - \frac{\De}{|q|^2}e_3\right)\psi -2\pr_{\cos\th}(h)\dual\Re(\Jk)\c\nab\psi
\eeaa
and hence
\beaa
\square_s(h\psi)  &=& h\square_s\psi+\pr_r(h)\left(e_4 - \frac{\De}{|q|^2}e_3\right)\psi -2\pr_{\cos\th}(h)\dual\Re(\Jk)\c\nab\psi+\square(h)\psi\\
&=& h\left(\square_s+\frac{2s(2-p)}{r|q|^2}(a^2\pr_t+a\pr_\phi)\right)\psi + hA_s^{(p)}[h]\psi -2\pr_{\cos\th}(h)\dual\Re(\Jk)\c\nab\psi+\square(h)\psi,
\eeaa
where $A_s^{(p)}[h]$ denotes the first-order operator given by 
\beaa
A_s^{(p)}[h]:= -\frac{2s(2-p)}{r|q|^2}(a^2\pr_t+a\pr_\phi) +\pr_r(h)\left(e_4 - \frac{\De}{|q|^2}e_3\right).
\eeaa

Next, we evaluate $A_s[h]$. We start with the case $s=+2$ for which we have
\beaa
&& A_{+2}^{(p)}[h]+\frac{2\De}{|q|^2}\frac{\pr_r(h)}{h}e_3\\ 
&=& -\frac{4(2-p)}{r|q|^2}(a^2\pr_t+a\pr_\phi) +\frac{\pr_r(h)}{h}\left(e_4 + \frac{\De}{|q|^2}e_3\right)\\
&=& -\frac{4(2-p)}{r|q|^2}(a^2\pr_t+a\pr_\phi) +2\frac{\pr_r(h)}{h}\left(\frac{r^2+a^2}{|q|^2}\pr_t+\frac{a}{|q|^2}\pr_\phi\right)\\
&=& -\frac{2a}{r|q|^2}\left(2(2-p) - r\frac{\pr_r(h)}{h}\right)\big(\pr_\phi+a(\sin\th)^2\pr_t\big)\\
&&+\frac{2}{r|q|^2}\left(-2(2-p)a^2+r(r^2+a^2)\frac{\pr_r(h)}{h} + \left(2(2-p) - r\frac{\pr_r(h)}{h}\right)a^2(\sin\th)^2\right)\pr_t\\
&=& -\frac{2a}{r|q|^2}\left(2(2-p) - r\frac{\pr_r(h)}{h}\right)\big(\pr_\phi+a(\sin\th)^2\pr_t\big)+2\left(\frac{\pr_r(h)}{h} -2(2-p)\frac{a^2(\cos\th)^2}{r|q|^2}\right)\pr_t\\
&=&  -\frac{2a}{r|q|^2}\left(2(2-p) - r\frac{\pr_r(h)}{h}\right)\big(\pr_\phi+a(\sin\th)^2\pr_t\big)+2\pr_r\log\left(h\left(\frac{|q|^2}{r^2}\right)^{2-p}\right)\pr_t. 
\eeaa
Also, we have for the case $s=-2$
\beaa
A_{-2}^{(p)}[h] -2\frac{\pr_r(h)}{h}e_4 &=& \frac{4(2-p)}{r|q|^2}(a^2\pr_t+a\pr_\phi) -\frac{\pr_r(h)}{h}\left(e_4 + \frac{\De}{|q|^2}e_3\right)\\
&=& -\left(A_{+2}^{(p)}[h]+\frac{2\De}{|q|^2}\frac{\pr_r(h)}{h}e_3\right)\\
&=&  \frac{2a}{r|q|^2}\left(2(2-p) - r\frac{\pr_r(h)}{h}\right)\big(\pr_\phi+a(\sin\th)^2\pr_t\big)\\
&& -2\pr_r\log\left(h\left(\frac{|q|^2}{r^2}\right)^{2-p}\right)\pr_t. 
\eeaa
We thus choose from now on
\bea\lab{eq:defandpropertiesofhpinprooffromNPinMa17tonewNP}
h_p:=\left(\frac{r^2}{|q|^2}\right)^{2-p}, \quad h_p=1+O(a^2r^{-2}), \quad \pr_r(h_p)=O(a^2r^{-3}),\quad \pr_{\cos\th}(h_p)=O(a^2r^{-2}),
\eea
which yields in view of the above
\beaa
A_{+2}^{(p)}[h_p] &=& O(a^2r^{-3})e_3+ O(ar^{-3})\big(\pr_\phi+a(\sin\th)^2\pr_t\big),\\
A_{-2}^{(p)}[h_p] &=& O(a^2r^{-3})e_4+ O(ar^{-3})\big(\pr_\phi+a(\sin\th)^2\pr_t\big), 
\eeaa
and hence
\beaa
\square_{+2}(h_p\psi) &=& \big(1+O(a^2r^{-2})\big)\left(\square_{+2}+\frac{4(2-p)}{r|q|^2}(a^2\pr_t+a\pr_\phi)\right)\psi + h_pA_{+2}^{(p)}[h_p]\psi \\
&& +O(a^2r^{-2})\dual\Re(\Jk)\c\nab\psi+O(a^2r^{-4})\psi\\
&=& \big(1+O(a^2r^{-2})\big)\left(\square_{+2}+\frac{4(2-p)}{r|q|^2}(a^2\pr_t+a\pr_\phi)\right)\psi + O(a^2r^{-3})e_3(\psi)\\
&&+ O(ar^{-3})\big(\pr_\phi+a(\sin\th)^2\pr_t\big) +O(a^2r^{-2})\dual\Re(\Jk)\c\nab\psi+O(a^2r^{-4})\psi,\\
\square_{-2}(h_p\psi) &=& \big(1+O(a^2r^{-2})\big)\left(\square_{-2}-\frac{4(2-p)}{r|q|^2}(a^2\pr_t+a\pr_\phi)\right)\psi + h_pA_{-2}^{(p)}[h_p]\psi \\
&& +O(a^2r^{-2})\dual\Re(\Jk)\c\nab\psi+O(a^2r^{-4})\psi\\
&=& \big(1+O(a^2r^{-2})\big)\left(\square_{-2}-\frac{4(2-p)}{r|q|^2}(a^2\pr_t+a\pr_\phi)\right)\psi + O(a^2r^{-3})e_4(\psi)\\
&&+ O(ar^{-3})\big(\pr_\phi+a(\sin\th)^2\pr_t\big) +O(a^2r^{-2})\dual\Re(\Jk)\c\nab\psi+O(a^2r^{-4})\psi,
\eeaa
or 
\bea\lab{eq:eqforsquaresofhppsiinprooffromNPinMa17tonewNP}
\bsplit
\square_{+2}(h_p\psi) =& \big(1+O(a^2r^{-2})\big)\left(\square_{+2}+\frac{4(2-p)}{r|q|^2}(a^2\pr_t+a\pr_\phi)\right)\psi + O(a^2r^{-3})e_3(\psi)\\
&+ O(ar^{-3})\mathcal{X}_{+2}\psi+O(a^2r^{-4})\psi,\\
\square_{-2}(h_p\psi) =& \big(1+O(a^2r^{-2})\big)\left(\square_{-2}-\frac{4(2-p)}{r|q|^2}(a^2\pr_t+a\pr_\phi)\right)\psi + O(a^2r^{-3})e_4(\psi)\\
&+ O(ar^{-3})\mathcal{X}_{-2}\psi+O(a^2r^{-4})\psi,
\end{split}
\eea
where $\mathcal{X}_s$ are regular horizontal vectorfields of the form
\beaa
\mathcal{X}_s = s\big(\pr_\phi+a(\sin\th)^2\pr_t\big)+O(ar^{-1})\sin\th\pr_\th.
\eeaa
 
We now rewrite the system of spin-weighted wave equations \eqref{eq:TeukolskywavetransportsystemasinMa17} as follows 
\beaa
\left(\square_{s} +\frac{2s(2-p)}{r|q|^2}(a^2\pr_t+a\pr_\phi) - \frac{4-2\de_{p0}}{|q|^2}\right)\dot{\phi}_{s,\text{NP}}^{(p)}=\ddot{L}^{(p)}_s[\dot{\phi}_{s,\text{NP}}], 
\eeaa
where, in view of \eqref{def:ComplexTeuScalars:wavesystem:Kerr:LspasinMa17}, $\ddot{L}^{(p)}_{s,\text{NP}}[\dot{\phi}_{s,\text{NP}}]$ are given by 
\beaa
\bsplit
\ddot{L}^{(0)}_{s,\text{NP}}[\dot{\phi}_{s,\text{NP}}] :=& \frac{2s}{r^3}(1+O(mr^{-1}))\dot{\phi}_{s,\text{NP}}^{(1)}+O(mr^{-3})\dot{\phi}_{+2,\text{NP}}^{(0)},\\
\ddot{L}^{(1)}_{s,\text{NP}}[\dot{\phi}_{s,\text{NP}}] :=& \frac{s}{r^3}(1+O(mr^{-1}))\dot{\phi}_{s,\text{NP}}^{(2)}  +O(mr^{-2})(\pr_{\phi}+a\pr_t)^{\leq 1}\dot{\phi}_{+2,\text{NP}}^{(0)} +O(mr^{-3})\dot{\phi}_{+2,\text{NP}}^{(1)},\\
\ddot{L}^{(2)}_{s,\text{NP}}[\dot{\phi}_{s,\text{NP}}] :=& O(mr^{-2})(\pr_\phi+a\pr_t)\dot{\phi}_{+2,\text{NP}}^{(1)} +O(m^2r^{-2})\dot{\phi}_{s,\text{NP}}^{(0)}+O(mr^{-3})\dot{\phi}_{s,\text{NP}}^{(2)}.
\end{split}
\eeaa
Plugging in \eqref{eq:eqforsquaresofhppsiinprooffromNPinMa17tonewNP}, and using the fact that, in view of \eqref{eq:relationbetweencomplexscalarsinMa17andnewnormalization} 
and \eqref{eq:defandpropertiesofhpinprooffromNPinMa17tonewNP}, we have
\beaa
\phi_{s,\text{NP}}^{(p)}=h_p\dot{\phi}_{s,\text{NP}}^{(p)},
\eeaa
we infer
\beaa
\left(\square_{s}  - \frac{4-2\de_{p0}}{|q|^2}\right)\phi_{s,\text{NP}}^{(p)}=L^{(p)}_{s,\text{NP}}[\phi_{s,\text{NP}}], 
\eeaa
where, in view of \eqref{def:ComplexTeuScalars:wavesystem:Kerr:LspasinMa17}, $L^{(p)}_{s,\text{NP}}[\phi_{s,\text{NP}}]$ are given by 
\beaa
\bsplit
L^{(0)}_{s,\text{NP}}[\phi_{s,\text{NP}}] :=& \frac{2s}{r^3}(1+O(mr^{-1}))\phi_{s,\text{NP}}^{(1)}+O(mr^{-3})(\mathcal{X}_s)^{\leq 1}\phi_{+2,\text{NP}}^{(0)},\\
L^{(1)}_{s,\text{NP}}[\phi_{s,\text{NP}}] :=& \frac{s}{r^3}(1+O(mr^{-1}))\phi_{s,\text{NP}}^{(2)}  +O(mr^{-3})(\mathcal{X}_s)^{\leq 1}\phi_{+2,\text{NP}}^{(1)}+O(mr^{-2})(\pr_{\phi}+a\pr_t)^{\leq 1}\phi_{+2,\text{NP}}^{(0)},\\
L^{(2)}_{s,\text{NP}}[\phi_{s,\text{NP}}] :=& O(mr^{-3})\phi_{s,\text{NP}}^{(2)}+O(mr^{-2})(\pr_\phi+a\pr_t)\phi_{+2,\text{NP}}^{(1)} +O(m^2r^{-2})\phi_{s,\text{NP}}^{(0)}.
\end{split}
\eeaa
Note that we also used the fact that, in view of \eqref{def:ComplexTeuScalars:wavesystem:Kerr}, we have, for $p=1,2$, 
\beaa
e_3(\phi_{+2,\text{NP}}^{(p)}) &=& O(r^{-2})\phi_{+2,\text{NP}}^{(p+1)}+O(r^{-1})\phi_{+2,\text{NP}}^{(p)},\\
e_4(\phi_{-2,\text{NP}}^{(p)}) &=& O(r^{-2})\phi_{-2,\text{NP}}^{(p+1)}+O(r^{-1})\phi_{-2,\text{NP}}^{(p)}.
\eeaa
This concludes the proof of Lemma \ref{lemma:structureofwavesystemforphisp:newNPform}.
\end{proof}


\subsection{Link between equations in tensorial form and in Newman-Penrose formalism}


In order to deduce Teukolsky equations in tensorial form in Section  \ref{sec:TeukolskyinKerrintensorialform} from the ones in NP formalism given by Lemma \ref{lemma:structureofwavesystemforphisp:newNPform}, we explain in this section how to go from equations in NP formalism to tensorial equations. To this end, with $e_1, e_2$ given as in \eqref{def:e1e2inKerr}, we consider ${\pmb\psi_{\pm 2}}\in\sk_2(\mathbb{C})$ and associate the complex-valued scalar $\psi_{\pm 2,\text{NP}}$ as follows 
\bea\lab{def:psiplusminus2:complexscalars:Kerr}
\psi_{+2,\text{NP}} := {\pmb\psi_{+2}}(e_1, e_1), \qquad \psi_{-2,\text{NP}} := \ov{{\pmb\psi_{-2}}(e_1, e_1)}.
\eea
Also, we have the following relations for $\pmb\psi\in\sk_2(\mathbb{C})$
\bea\lab{eq:relationspsi12psi22topsi11}
\pmb\psi(e_1, e_2)=\pmb\psi(e_2, e_1)=-i\pmb\psi(e_1, e_1), \qquad \pmb\psi(e_2, e_2)=-\pmb\psi(e_1, e_1).
\eea
Finally, we have the following formulas for connection coefficients in Kerr, see (14.13) in \cite{KS:Kerr}:
\bea
 \lab{eq:Kerr.La_abc}
 \bsplit
 \gam(\D_1e_1, e_2) &=0, \qquad \qquad \qquad  \gam(\D_2e_1, e_2) =\frac{r^2+a^2}{|q|^3} \cot \th, \,\, \\
 \gam(\D_3e_1, e_2) &=  \frac{a\cos \th}{|q|^2}, \qquad  \gam(\D_4e_1, e_2)= \frac{a\De\cos \th}{|q|^4}.
  \end{split}
 \eea
 
We start with the following lemma. 
\begin{lemma}\lab{lemma:scalarizationoncomplexscalarsofsquared2}
For $\pmb\psi\in\sk_2(\mathbb{C})$, we have
\beaa
\square_{\gam}(\pmb\psi(e_1, e_1)) &=&  (\squared_2\pmb\psi)(e_1, e_1) - \frac{4i\cos\th}{\sin^2\th|q|^2}\nab_{\pr_{\phi}}\pmb\psi(e_1, e_1)\\
&&+\left(\frac{4a^2\cos^2\th\De }{|q|^6} -\frac{4(\R)^2}{|q|^6} \cot^2\th\right)\pmb\psi(e_1, e_1).
\eeaa
\end{lemma}

\begin{proof}
For $\pmb\psi\in\sk_2(\mathbb{C})$, we have
\beaa
&&\square_{\gam}(\pmb\psi(e_1, e_1))\\ 
&=& \gam^{\a\b}e_\a(e_\b(\pmb\psi(e_1, e_1)) -\gam^{\a\b}\D_{\D_\a e_\b}(\pmb\psi(e_1, e_1))\\
&=& \gam^{\a\b}e_\a\big(\Ddot_\b\pmb\psi(e_1, e_1)+2\pmb\psi(\Ddot_\b e_1, e_1)\big) -\gam^{\a\b}\big(\Ddot_{\Ddot_\a e_\b}\pmb\psi(e_1, e_1)+2\pmb\psi(\Ddot_{\Ddot_\a e_\b}e_1, e_1)\big)\\
&=& \gam^{\a\b}\big(\Ddot_\a\Ddot_\b\pmb\psi(e_1, e_1)+\Ddot_{\Ddot_\a e_\b}\pmb\psi(e_1, e_1)+2\Ddot_\b\pmb\psi(\Ddot_\a e_1, e_1)+2\Ddot_\a\pmb\psi(\Ddot_\b e_1, e_1)+2\pmb\psi(\Ddot_\a\Ddot_\b e_1, e_1)\\
&&+2\pmb\psi(\Ddot_{\Ddot_\a e_\b}e_1, e_1)+2\pmb\psi(\Ddot_\b e_1, \Ddot_\a e_1)\big)  -\gam^{\a\b}\big(\Ddot_{\Ddot_\a e_\b}\pmb\psi(e_1, e_1)+2\pmb\psi(\Ddot_{\Ddot_\a e_\b}e_1, e_1)\big)\\
&=& \gam^{\a\b}\Ddot_\a\Ddot_\b\pmb\psi(e_1, e_1)+4\gam^{\a\b}\Ddot_\a\pmb\psi(\Ddot_\b e_1, e_1) +2\gam^{\a\b}\pmb\psi(\Ddot_\a\Ddot_\b e_1,e _1)+2\pmb\psi(\Ddot^\a e_1, \Ddot_\a e_1)\\
&=& \squared_2\pmb\psi(e_1, e_1)+4\gam^{\a\b}{\Ddot_{\a}}\pmb\psi({\Ddot_\b} e_1, e_1)+2\pmb\psi({\Ddot^\a\Ddot_{\a}}e_1, e_1)+2\pmb\psi({\Ddot^\a}e_1, {\Ddot_\a}e_1).
\eeaa

Next, we compute
\beaa
\gam^{\a\b}{\Ddot_\a}\pmb\psi({\Ddot_\b} e_1, e_1) &=& \gam^{\a\b}\gam({\Ddot_\b} e_1, e_2){\Ddot_\a}\pmb\psi(e_2, e_1)\\
&=& -i\gam^{\a\b}\gam({\Ddot_\b} e_1, e_2){\Ddot_\a}\pmb\psi(e_1, e_1)\\
&=&  \frac{i}{2}\gam(\nab_4e_1, e_2)\nab_3\pmb\psi(e_1, e_1) +\frac{i}{2}\gam(\nab_3e_1, e_2)\nab_4\pmb\psi(e_1, e_1)\\
&& - i\gam(\nab_2e_1, e_2)\nab_2\pmb\psi(e_1, e_1)\\
&=&  i\frac{a\De\cos \th}{2|q|^4}\nab_3\pmb\psi(e_1, e_1) +i\frac{a\cos \th}{2|q|^2}\nab_4\pmb\psi(e_1, e_1) - i\frac{r^2+a^2}{|q|^3} \cot \th\nab_2\pmb\psi(e_1, e_1).
\eeaa
Since we have
\beaa
e_4+\frac{\De}{|q|^2}e_3 &=& \frac{2(r^2+a^2)}{|q|^2}\pr_t+\frac{2a}{|q|^2}\pr_\phi,\qquad e_2=\frac{a\sin\th}{|q|}\pr_t+\frac{1}{|q|\sin\th}\pr_\phi,
\eeaa
we infer
\beaa
\gam^{\a\b}{\Ddot_\a}\pmb\psi({\Ddot_\b} e_1, e_1) &=&  \frac{ai(\R)\cos \th}{|q|^4}\nab_{\pr_t}\pmb\psi(e_1, e_1) +\frac{ai\cos \th}{|q|^2}\frac{a}{|q|^2}\nab_{\pr_\phi}\pmb\psi(e_1, e_1)\\
&& - i\frac{r^2+a^2}{|q|^3}\frac{1}{|q|\sin\th}\cot \th\nab_{\pr_{\phi}}\pmb\psi(e_1, e_1) - i\frac{r^2+a^2}{|q|^3}\frac{a\sin\th}{|q|}\cot \th\nab_{\pr_t}\pmb\psi(e_1, e_1)\\
&=& - i\frac{\R -a^2\sin^2\th}{|q|^4}\frac{\cos\th}{\sin^2\th}\nab_{\pr_{\phi}}\pmb\psi(e_1, e_1)\\
&=& - i\frac{1}{|q|^2}\frac{\cos\th}{\sin^2\th}\nab_{\pr_{\phi}}\pmb\psi(e_1, e_1)  
\eeaa
and hence
\beaa
&&\square_{\gam}(\pmb\psi(e_1, e_1)) \nn\\
&=& (\squared_2\pmb\psi)(e_1, e_1)+4\gam^{\a\b}{\Ddot_\a}\pmb\psi({\Ddot_\b} e_1, e_1)+2\pmb\psi({\Ddot^\a\Ddot_\a} e_1, e_1)+2\pmb\psi({\Ddot^\a e_1}, {\Ddot_\a e_1})\\
&=& (\squared_2\pmb\psi)(e_1, e_1) - i\frac{4}{|q|^2}\frac{\cos\th}{\sin^2\th}\nab_{\pr_{\phi}}\pmb\psi(e_1, e_1)+2\pmb\psi({\Ddot^\a\Ddot_\a} e_1, e_1)+2\pmb\psi({\Ddot^\a e_1}, {\Ddot_\a e_1})\\
&=& (\squared_2\pmb\psi)(e_1, e_1) - i\frac{4}{|q|^2}\frac{\cos\th}{\sin^2\th}\nab_{\pr_{\phi}}\pmb\psi(e_1, e_1)+2\gam({\Ddot^\a\Ddot_\a} e_1, e_b)\pmb\psi(e_b, e_1)+2\pmb\psi({\Ddot^\a e_1}, {\Ddot_\a e_1})\\
&=& (\squared_2\pmb\psi)(e_1, e_1) - i\frac{4}{|q|^2}\frac{\cos\th}{\sin^2\th}\nab_{\pr_{\phi}}\pmb\psi(e_1, e_1)+2\gam({\Ddot^\a\Ddot_\a} e_1, e_1)\pmb\psi(e_1, e_1)\\
&& +2\gam({\Ddot^\a\Ddot_\a} e_1, e_2)\pmb\psi(e_2, e_1)+2\pmb\psi({\Ddot^\a e_1}, {\Ddot_\a e_1})\\
&=& (\squared_2\pmb\psi)(e_1, e_1) - i\frac{4}{|q|^2}\frac{\cos\th}{\sin^2\th}\nab_{\pr_{\phi}}\pmb\psi(e_1, e_1)+2\gam({\Ddot^\a\Ddot_\a} e_1, e_1)\pmb\psi(e_1, e_1)\\
&& -2i\gam({\Ddot^\a\Ddot_\a} e_1, e_2)\pmb\psi(e_1, e_1)+2\pmb\psi({\Ddot^\a e_1}, {\Ddot_\a e_1}).
\eeaa

Also, we have
\beaa
\gam({\Ddot^\a\Ddot_\a} e_1, e_b) &=& \gam^{\a\b}\gam({\Ddot_\a}{\Ddot_\b} e_1, e_b)\\
&=& \gam^{\a\b}\gam(\Ddot_\a\Ddot_\b e_1, e_b)\\
&=& \gam^{\a\b}e_\a(\gam(\Ddot_\b e_1, e_b)) -\gam^{\a\b}\gam(\Ddot_\b e_1, \Ddot_\a e_b)
-\gam(\Ddot_{\Ddot^\a e_{\a}} e_1, e_b)\\
&=& \gam^{\a\b}e_\a(\gam({\Ddot_\b} e_1, e_b)) -\gam^{\a\b}\gam({\Ddot_\b} e_1, {\Ddot_\a} e_b)
-\gam(\nab_{\D^\a e_{\a}} e_1, e_b).
\eeaa
We compute the before to last term respectively for $b=1$ and $b=2$ and obtain
\beaa
 -\gam^{\a\b}\gam({\Ddot_\b} e_1, {\Ddot_\a} e_1) &=& \gam(\nab_3 e_1, \nab_4 e_1)  - \gam(\nab_1e_1, \nab_1e_1) - \gam(\nab_2e_1, \nab_2e_1)\\
 &=&  \gam(\nab_4e_1, \nab_3e_1) -\gam(\nab_2e_1, \nab_2e_1) \\
 &=& \gam(\nab_4e_1, e_2)\gam(\nab_3e_1,e_2) -(\gam(\nab_2e_1,e_2))^2\\
 &=& \frac{a^2\cos^2\th\De}{|q|^6} -\frac{(\R)^2}{|q|^6} \cot^2\th ,\\
 -\gam^{\a\b}\gam({\Ddot_\b} e_1, {\Ddot_\a} e_2) 
 &=& \frac{1}{2}\gam(\nab_3e_1, \nab_4 e_2) + \frac{1}{2}\gam(\nab_4e_1, \nab_3e_2) -\gam(\nab_2e_1, \nab_2e_2) \\
 &=&0.
\eeaa
Also, the last term is given by
\beaa
-\gam(\nab_{\D^\a e_{\a}} e_1, e_1)&=&0,\\
-\gam(\nab_{\D^\a e_{\a}} e_1, e_2)&=& -\gam\bigg(\nab_{-\frac{r}{|q|^2}e_4 +\left(\frac{1}{2}\pr_r\left(\frac{\De}{|q|^2}\right) +\frac{\De r}{|q|^4}\right)e_3  +\big(-\eta_1-\etab_1+\gam(\nab_2e_2, e_1)\big)e_1}e_1, e_2\bigg)\\
&=& \frac{r}{|q|^2}\frac{a\De\cos\th}{|q|^4}-\left(\frac{1}{2}\pr_r\left(\frac{\De}{|q|^2}\right) +\frac{\De r}{|q|^4}\right)\frac{a\cos\th}{|q|^2}\\
&=& -\frac{1}{2}\pr_r\left(\frac{\Delta}{|q|^2}\right)\frac{a\cos\th}{|q|^2}
\eeaa
where, in addition to \eqref{eq:Kerr.La_abc}, we used
\beaa
\D^\a e_{\a} &=& -\frac{1}{2}\D_3e_4 -\frac{1}{2}\D_4e_3+\D_1e_1+\D_2e_2\\
&=& -\frac{1}{2}\big(2\omb e_4+2\eta_be_b) -\frac{1}{2}\big(2\om e_3+2\etab_be_b)
+\nab_1e_1+\nab_2e_2+\frac{1}{2}\trch e_3+\frac{1}{2}\trchb e_4\\
&=& \left(-\omb+\frac{1}{2}\trchb\right)e_4 +\left(-\om+\frac{1}{2}\trch\right)e_3  \nn\\
&&+\big(-\eta_1-\etab_1+\gam(\nab_2e_2, e_1)\big)e_1 -(\eta_2+\etab_2)e_2\\
&=&  -\frac{r}{|q|^2}e_4 +\left(\frac{1}{2}\pr_r\left(\frac{\De}{|q|^2}\right) +\frac{\De r}{|q|^4}\right)e_3  +\big(-\eta_1-\etab_1+\gam(\nab_2e_2, e_1)\big)e_1
\eeaa
which itself relies on the following explicit values in Kerr{, see \eqref{eq:KerrvaluesofcomplexifiedRicci},}
\beaa
\bsplit
\trch=&\Re\left(\frac{2\De\ov{q}}{|q|^4}\right)=\frac{2\De r}{|q|^4}, \qquad \trchb=\Re\left(-\frac{2q}{|q|^2}\right)=-\frac{2r}{|q|^2}, \qquad\om = -\frac{1}{2}\pr_r\left(\frac{\De}{|q|^2}\right), \qquad \omb=0,\\
\eta_2=&\Re\left(\frac{aq}{|q|^2}\Jk_2\right)=\Re\left(\frac{aq}{|q|^2}\frac{\sin\th}{|q|}\right)=\frac{ar\sin\th}{|q|^3},\\
\etab_2=& -\Re\left(\frac{a\ov{q}}{|q|^2}\Jk_2\right)=-\Re\left(\frac{a\ov{q}}{|q|^2}\frac{\sin\th}{|q|}\right)=-\frac{ar\sin\th}{|q|^3}.
\end{split}
\eeaa
Plugging in the above, this yields
\beaa
\gam({\Ddot^\a\Ddot_\a} e_1, e_1) &=& \gam^{\a\b}e_\a(\gam({\Ddot_\b} e_1, e_1)) -\gam^{\a\b}\gam({\Ddot_\b} e_1, {\Ddot_\a} e_1)-\gam(\nab_{\Ddot^\a e_{\a}} e_1, e_1)\\
&=& -\gam^{\a\b}\gam({\Ddot_\b} e_1, {\Ddot_\a} e_1)\\
&=& \frac{a^2\cos^2\th\De}{|q|^6} -\frac{(\R)^2}{|q|^6} \cot^2\th
\eeaa
and 
\beaa
&&\gam({\Ddot^\a\Ddot_\a} e_1, e_2) \nn\\
&=&\gam^{\a\b}e_\a(\gam({\Ddot_\b} e_1, e_2)) -\gam^{\a\b}\gam({\Ddot_\b} e_1, {\Ddot_\a} e_2)-\gam(\nab_{\D^\a e_a} e_1, e_2)\\
&=& \gam^{\a\b}e_\a(\gam({\Ddot_\b} e_1, e_2))-\gam(\nab_{\D^\a e_a} e_1, e_2)\\
&=& -\frac{1}{2}e_4(\gam(\nab_3e_1, e_2)) -\frac{1}{2}e_3(\gam(\nab_4e_1, e_2))+e_2(\gam(\nab_2e_1, e_2))-\frac{1}{2}\pr_r\bigg(\frac{\Delta}{|q|^2}\bigg)\frac{a\cos\th}{|q|^2}\\
&=& -\frac{1}{2}e_4\left(\frac{a\cos \th}{|q|^2}\right) -\frac{1}{2}e_3\left(\frac{a\De\cos \th}{|q|^4}\right)+e_2\left(\frac{r^2+a^2}{|q|^3} \cot \th\right)-\frac{1}{2}\pr_r\bigg(\frac{\Delta}{|q|^2}\bigg)\frac{a\cos\th}{|q|^2}\\
&=& -\frac{1}{2}\frac{\De}{|q|^2}\pr_r\left(\frac{a\cos \th}{|q|^2}\right) +\frac{1}{2}\pr_r\left(\frac{a\De\cos \th}{|q|^4}\right)-\frac{1}{2}\pr_r\bigg(\frac{\Delta}{|q|^2}\bigg)\frac{a\cos\th}{|q|^2}\\
&=&0.
\eeaa
We also have
\beaa
\pmb\psi({\Ddot^\a e_1}, {\Ddot_\a e_1})&=&\gam^{\a\b}\pmb\psi(\Ddot_{\b}e_1, {\Ddot_\a e_1})\nn\\
&=&\gam^{\a\b}\gam(\Ddot_{\b}e_1,e_2)\gam({\Ddot_\a e_1},e_2)\pmb\psi(e_2, e_2)\nn\\
&=& \bigg(\frac{a^2\cos^2\th\De}{|q|^6}-\frac{(\R)^2}{|q|^6} \cot^2\th\bigg)\pmb\psi(e_1, e_1).
\eeaa
Hence, we deduce 
\beaa
\square_{\gam}(\pmb\psi(e_1, e_1)) &=& (\squared_2\pmb\psi)(e_1, e_1) - i\frac{4}{|q|^2}\frac{\cos\th}{\sin^2\th}\nab_{\pr_{\phi}}\pmb\psi(e_1, e_1)+2\gam({\Ddot^\a\Ddot_\a} e_1, e_1)\pmb\psi(e_1, e_1)\\
&&  -2i\gam({\Ddot^\a\Ddot_\a} e_1, e_2)\pmb\psi(e_1, e_1) 
+2\pmb\psi({\Ddot^\a e_1}, {\Ddot_\a e_1})\\
&=&  (\squared_2\pmb\psi)(e_1, e_1) - i\frac{4}{|q|^2}\frac{\cos\th}{\sin^2\th}\nab_{\pr_{\phi}}\pmb\psi(e_1, e_1)\\
&&+\left(\frac{4a^2\cos^2\th\De }{|q|^6} -\frac{4(\R)^2}{|q|^6} \cot^2\th\right)\pmb\psi(e_1, e_1)
\eeaa
as stated. This concludes the proof of Lemma \ref{lemma:scalarizationoncomplexscalarsofsquared2}.
\end{proof}

\begin{lemma}\lab{lemma:scalarizationoncomplexscalarsofnabprtannabprphi}
For $\pmb\psi\in\sk_2(\mathbb{C})$, we have
\beaa
\pr_t(\pmb\psi(e_1, e_1)) &=& \nab_{\pr_t}\pmb\psi(e_1, e_1)  +\frac{4iamr\cos\th}{|q|^4}\pmb\psi(e_1, e_1),
\eeaa
\beaa
\pr_\phi(\pmb\psi(e_1, e_1)) &=& \nab_{\pr_\phi}\pmb\psi(e_1, e_1)   -\frac{2i((\R)^2-a^2\sin^2\th\De)\cos\th}{|q|^4}\pmb\psi(e_1, e_1)
\eeaa
and
\beaa
e_3(\pmb\psi(e_1,e_1)) &=& \nab_3\pmb\psi(e_1,e_1)-\frac{2ia\cos\th}{|q|^2}\pmb\psi(e_1,e_1), \\ 
e_4(\pmb\psi(e_1,e_1)) &=& \nab_4\pmb\psi(e_1,e_1)-\frac{2ia\De\cos\th}{|q|^4}\pmb\psi(e_1,e_1).
\eeaa
\end{lemma}

\begin{proof}
For $\pmb\psi\in\sk_2(\mathbb{C})$, and $X$ a vectorfield spanned by $(e_2, e_3, e_4)$, we have
\beaa
X(\pmb\psi(e_1, e_1)) &=& \nab_X(\pmb\psi(e_1, e_1))\\
&=& \nab_X\pmb\psi(e_1, e_1)+2\pmb\psi(\nab_Xe_1, e_1)\\
&=& \nab_X\pmb\psi(e_1, e_1)+2\gam(\nab_Xe_1, e_2)\pmb\psi(e_2, e_1)\\
&=& \nab_X\pmb\psi(e_1, e_1)-2i\gam(\nab_Xe_1, e_2)\pmb\psi(e_1, e_1)\\
&=& \nab_X\pmb\psi(e_1, e_1) -2iX^4\gam(\nab_4e_1, e_2)\pmb\psi(e_1, e_1) -2iX^3\gam(\nab_3e_1, e_2)\pmb\psi(e_1, e_1)\\
&& -2iX^2\gam(\nab_2e_1, e_2)\pmb\psi(e_1, e_1)
\eeaa
and hence
\beaa
X(\pmb\psi(e_1, e_1)) &=& \nab_X\pmb\psi(e_1, e_1) -2iX^4\frac{a\De\cos \th}{|q|^4}\pmb\psi(e_1, e_1) -2iX^3\frac{a\cos \th}{|q|^2}\pmb\psi(e_1, e_1)\\
&& -2iX^2\frac{r^2+a^2}{|q|^3} \cot \th\pmb\psi(e_1, e_1).
\eeaa
Since 
\beaa
2\pr_t &=& e_4+\frac{\De}{|q|^2}e_3 -\frac{2a\sin\th}{|q|}e_2,\\
2\pr_\phi &=& \frac{2(\R)\sin\th}{|q|}e_2 -a(\sin\th)^2{e_4} -\frac{a(\sin\th)^2\De}{|q|^2}{e_3},
\eeaa
we infer
\beaa
\pr_t(\pmb\psi(e_1, e_1)) &=& \nab_{\pr_t}\pmb\psi(e_1, e_1) -i\frac{a\De\cos \th}{|q|^4}\pmb\psi(e_1, e_1) -i\frac{\De}{|q|^2}\frac{a\cos \th}{|q|^2}\pmb\psi(e_1, e_1)\\
&& +i\frac{2a\sin\th}{|q|}\frac{r^2+a^2}{|q|^3} \cot \th\pmb\psi(e_1, e_1)\\
&=& \nab_{\pr_t}\pmb\psi(e_1, e_1) -\frac{2ai\De\cos \th}{|q|^4}\pmb\psi(e_1, e_1)  +\frac{2ai\cos\th(\R)}{|q|^4}\pmb\psi(e_1, e_1)\\
&=& \nab_{\pr_t}\pmb\psi(e_1, e_1)  +\frac{2ai\cos\th(\R-\De)}{|q|^4}\pmb\psi(e_1, e_1)\\
&=& \nab_{\pr_t}\pmb\psi(e_1, e_1)  +\frac{4amri\cos\th}{|q|^4}\pmb\psi(e_1, e_1),
\eeaa
\beaa
\pr_\phi(\pmb\psi(e_1, e_1)) &=& \nab_{\pr_\phi}\pmb\psi(e_1, e_1) +ia(\sin\th)^2\frac{a\De\cos \th}{|q|^4}\pmb\psi(e_1, e_1) +i\frac{a(\sin\th)^2\De}{|q|^2}\frac{a\cos \th}{|q|^2}\pmb\psi(e_1, e_1)\\
&& -i\frac{2(\R)\sin\th}{|q|}\frac{r^2+a^2}{|q|^3} \cot \th\pmb\psi(e_1, e_1)\\
&=& \nab_{\pr_\phi}\pmb\psi(e_1, e_1) +\frac{2ia^2(\sin\th)^2\De\cos \th}{|q|^4}\pmb\psi(e_1, e_1)  -\frac{2i(\R)^2\cos\th}{|q|^4}\pmb\psi(e_1, e_1)\\
&=& \nab_{\pr_\phi}\pmb\psi(e_1, e_1)   -\frac{2i((\R)^2-a^2\sin^2\th\De)\cos\th}{|q|^4}\pmb\psi(e_1, e_1)
\eeaa
and
\beaa
e_3(\pmb\psi(e_1,e_1)) &=& \nab_3\pmb\psi(e_1,e_1)-\frac{2ia\cos\th}{|q|^2}\pmb\psi(e_1,e_1), \\ 
e_4(\pmb\psi(e_1,e_1)) &=& \nab_4\pmb\psi(e_1,e_1)-\frac{2ia\De\cos\th}{|q|^4}\pmb\psi(e_1,e_1),
\eeaa
as stated. This concludes the proof of Lemma \ref{lemma:scalarizationoncomplexscalarsofnabprtannabprphi}.
\end{proof}

\begin{lemma}
\lab{lem:scalarizationofspinweightedwaveoperator:Kerr}
For $\pmb\psi\in\sk_2(\mathbb{C})$, we have
\beaa
&& \square_{\gam}(\pmb\psi(e_1, e_1)) + \frac{4i}{|q|^2}\left(\frac{\cos\th}{\sin^2\th}\pr_{\phi} -a\cos\th\pr_{t}\right)(\pmb\psi(e_1, e_1)) -\frac{4}{|q|^2}(\cot\th)^2\pmb\psi(e_1, e_1)\\ 
&=& \left(\squared_2\pmb\psi -\frac{4ia\cos\th}{|q|^2}\nab_{\pr_t}\pmb\psi\right)(e_1, e_1) +\frac{4a^2\cos^2\th}{|q|^6}\Big(|q|^2+6mr\Big)\pmb\psi(e_1, e_1).
\eeaa
\end{lemma}

\begin{proof}
For $\pmb\psi\in\sk_2(\mathbb{C})$, recall from Lemma \ref{lemma:scalarizationoncomplexscalarsofsquared2} that we have
\beaa
\square_{\gam}(\pmb\psi(e_1, e_1)) &=&  (\squared_2\pmb\psi)(e_1, e_1) - i\frac{4}{|q|^2}\frac{\cos\th}{\sin^2\th}\nab_{\pr_{\phi}}\pmb\psi(e_1, e_1)\\
&&+\left(\frac{4a^2\cos^2\th\De }{|q|^6} -\frac{4(\R)^2}{|q|^6} \cot^2\th\right)\pmb\psi(e_1, e_1).
\eeaa
Also, recall from Lemma \ref{lemma:scalarizationoncomplexscalarsofnabprtannabprphi} that we have
\beaa
\pr_t(\pmb\psi(e_1, e_1)) &=& \nab_{\pr_t}\pmb\psi(e_1, e_1)  +\frac{4amri\cos\th}{|q|^4}\pmb\psi(e_1, e_1)
\eeaa
and 
\beaa
\pr_\phi(\pmb\psi(e_1, e_1)) &=& \nab_{\pr_\phi}\pmb\psi(e_1, e_1)   -\frac{2i((\R)^2-a^2\sin^2\th\De)\cos\th}{|q|^4}\pmb\psi(e_1, e_1),
\eeaa
hence, we infer
\beaa
&&\square_{\gam}(\pmb\psi(e_1, e_1)) + \frac{4i}{|q|^2}\left(\frac{\cos\th}{\sin^2\th}\pr_{\phi} -a\cos\th\pr_{t}\right)(\pmb\psi(e_1, e_1))\\
&=&  (\squared_2\pmb\psi)(e_1, e_1) - i\frac{4}{|q|^2}\frac{\cos\th}{\sin^2\th}\nab_{\pr_{\phi}}\pmb\psi(e_1, e_1)+\left(\frac{4a^2\cos^2\th\De }{|q|^6} -\frac{4(\R)^2}{|q|^6} \cot^2\th\right)\pmb\psi(e_1, e_1)
\\
&& + \frac{4i}{|q|^2}\frac{\cos\th}{\sin^2\th}\left(\nab_{\pr_\phi}\pmb\psi(e_1, e_1)   -\frac{2i((\R)^2-a^2\sin^2\th\De)\cos\th}{|q|^4}\pmb\psi(e_1, e_1)\right)\\
&& -\frac{4ia\cos\th}{|q|^2}\left( \nab_{\pr_t}\pmb\psi(e_1, e_1)  +\frac{4amri\cos\th}{|q|^4}\pmb\psi(e_1, e_1)\right)\\
&=& (\squared_2\pmb\psi)(e_1, e_1) -\frac{4ia\cos\th}{|q|^2}\nab_{\pr_t}\pmb\psi(e_1, e_1)\\
&&+\frac{2}{|q|^6}\Big(2a^2\cos^2\th\De+8a^2mr\cos^2\th\\
&&\qquad\quad+\big(-{2}(\R)^2+4((\R)^2-a^2\sin\th^2\De)\big)\cot^2\th\Big)\pmb\psi(e_1, e_1)
\eeaa
and hence
\beaa
\square_{\gam}(\pmb\psi(e_1, e_1)) + \frac{4i}{|q|^2}\left(\frac{\cos\th}{\sin^2\th}\pr_{\phi} -a\cos\th\pr_{t}\right)(\pmb\psi(e_1, e_1)) = \left(\squared_2\pmb\psi -\frac{4ia\cos\th}{|q|^2}\nab_{\pr_t}\pmb\psi+V\pmb\psi\right)(e_1, e_1)
\eeaa
where 
\beaa
V&:=&\frac{2}{|q|^6}\Big(2a^2\cos^2\th\De+8a^2mr\cos^2\th+\Big(-{2}(\R)^2+4((\R)^2-a^2\sin\th^2\De)\Big)\cot^2\th\Big)\\
&=& \frac{1}{|q|^6}\Big(4(\R)^2\cot^2\th -4a^2\cos^2\th\De +16a^2mr\cos^2\th\Big)\\
&=& \frac{4}{|q|^2}\cot^2\th+\frac{1}{|q|^6}\Big(4((\R)^2-|q|^4)\cot^2\th -4a^2\cos^2\th\De +16a^2mr\cos^2\th\Big)\\
&=& \frac{4}{|q|^2}\cot^2\th+\frac{1}{|q|^6}\Big(4(\R+|q|^2)a^2\cos^2\th -4a^2\cos^2\th\De +16a^2mr\cos^2\th\Big)\\ 
&=&  \frac{4}{|q|^2}\cot^2\th +\frac{4a^2\cos^2\th}{|q|^6}\Big(|q|^2+6mr\Big).
\eeaa
This yields 
\beaa
&& \square_{\gam}(\pmb\psi(e_1, e_1)) + \frac{4i}{|q|^2}\left(\frac{\cos\th}{\sin^2\th}\pr_{\phi} -a\cos\th\pr_{t}\right)(\pmb\psi(e_1, e_1)) -\frac{4}{|q|^2}(\cot\th)^2\pmb\psi(e_1, e_1)\\ 
&=& \left(\squared_2\pmb\psi -\frac{4ia\cos\th}{|q|^2}\nab_{\pr_t}\pmb\psi\right)(e_1, e_1) +\frac{4a^2\cos^2\th}{|q|^6}\Big(|q|^2+6mr\Big)\pmb\psi(e_1, e_1)
\eeaa
as stated, which concludes the proof of Lemma \ref{lem:scalarizationofspinweightedwaveoperator:Kerr}.
\end{proof}

The following corollary will allow us to deduce Teukolsky wave/transport systems in tensorial form in Section  \ref{sec:TeukolskyinKerrintensorialform} from the ones in NP formalism given by Lemma \ref{lemma:structureofwavesystemforphisp:newNPform}.
\begin{corollary}\lab{cor:howtogofromequationsinNPformtotensorialequations}
Let ${\pmb\psi_{\pm 2}}\in\sk_2(\mathbb{C})$, and let $\psi_{\pm 2,\text{NP}}$ be the complex-valued scalars introduced in \eqref{def:psiplusminus2:complexscalars:Kerr}, i.e., 
\beaa
\psi_{+2,\text{NP}} := {\pmb\psi_{+2}}(e_1, e_1), \qquad \psi_{-2,\text{NP}} := \ov{{\pmb\psi_{-2}}(e_1, e_1)}.
\eeaa
Then, we have 
\beaa
\square_{+2}(\psi_{+2,\text{NP}}) &=& \left(\squared_2{\pmb\psi_{+2}} -\frac{4ia\cos\th}{|q|^2}\nab_{\pr_t}{\pmb\psi_{+2}}\right)(e_1, e_1) +\frac{4a^2\cos^2\th}{|q|^6}\Big(|q|^2+6mr\Big)\psi_{+2,\text{NP}},\\
\ov{\square_{-2}(\psi_{-2,\text{NP}})} &=& \left(\squared_2{\pmb\psi_{-2}} -\frac{4ia\cos\th}{|q|^2}\nab_{\pr_t}{\pmb\psi_{-2}}\right)(e_1, e_1) +\frac{4a^2\cos^2\th}{|q|^6}\Big(|q|^2+6mr\Big)\ov{\psi_{-2,\text{NP}}},
\eeaa
and 
\beaa
e_3(\psi_{+2,\text{NP}}) = \frac{q}{\ov{q}}\nab_3\left(\frac{\ov{q}}{q}{\pmb\psi_{+2}}\right)(e_1, e_1),\qquad e_4(\ov{\psi_{-2,\text{NP}}}) = \frac{\ov{q}}{q}\nab_4\left(\frac{q}{\ov{q}}{\pmb\psi_{-2}}\right)(e_1, e_1).
\eeaa
Also, if $\pmb\psi(e_1, e_1)=0$ for all $\th\in (0,\pi)$, then $\pmb\psi$ vanishes identically. 
\end{corollary}

\begin{proof}
The statements concerning wave equations and transport equations follow immediately from Lemmas \ref{lem:scalarizationofspinweightedwaveoperator:Kerr} and \ref{lemma:scalarizationoncomplexscalarsofnabprtannabprphi} respectively. Also, in view of the identities \eqref{eq:relationspsi12psi22topsi11}, if $\pmb\psi(e_1, e_1)=0$ for all $\th\in (0,\pi)$, then $\pmb\psi$ vanishes identically for all $\th\in (0,\pi)$ and hence everywhere by continuity.
\end{proof}


\subsection{Teukolsky equation in Kerr in tensorial form}
\lab{sec:TeukolskyinKerrintensorialform}


Denoting by $\W$ the linearized Weyl curvature tensor, we introduce $\pmb\phi_{\pm 2}\in\sk_2(\mathbb{C})$, see Definition \ref{def:skC:horizontaltensors}, as follows
\bea\lab{eq:defintionoftensorialTeukolskyscalarspsipm2}
\bsplit
{\pmb\phi_{+2}} :=&|q|^2\ov{q}^2(\a+i\dual\a),  \qquad\,\,\,\, \a_{bc}:=\W(e_4, e_b, e_4, e_c),\\
{\pmb\phi_{-2}} :=&|q|^{-2}{q^2}(\aa+i\dual\aa),  \qquad \aa_{bc}:=\W(e_3, e_b, e_3, e_c),
\end{split}
\eea
where $\a$ and  $\aa$ are horizontal symmetric traceless 2-tensors, and where the principal null pair $(e_3, e_4)$ is normalized as in \eqref{def:e3e4inKerr}. In particular, notice from the definition of $\phi_{\pm 2,\text{NP}}$ in \eqref{def:phiplusminus2:complexscalars:Kerr} that
\bea\lab{eq:linkteukolskyscalarsandtensors}
\phi_{+2,\text{NP}}=\pmb\phi_{+2}(e_1, e_1), \qquad \phi_{-2,\text{NP}}=\ov{\pmb\phi_{-2}(e_1, e_1)}.
\eea

The following lemma provides Teukolsky equations for $\pmb\phi_{\pm 2}$ in tensorial form.
\begin{lemma}\lab{lemma:Teukolskyintensorialform}
Let $\pmb\phi_{\pm 2}\in\sk_2(\mathbb{C})$ be given by 
\eqref{eq:defintionoftensorialTeukolskyscalarspsipm2}. Then, Teukolsky equations in tensorial form, for $s=\pm 2$, are given by 
\bea\lab{eq:TeukolskyequationforAandAbintensorialforminKerr}
\nn\left(\squared_2 -\frac{4ia\cos\th}{|q|^2}\nab_{\pr_t} -  \frac{s}{|q|^2}\right)\pmb\phi_s +\frac{2s}{|q|^2}(r-m)\nab_3\pmb\phi_s  -\frac{4sr}{|q|^2}\nab_{\pr_t}\pmb\phi_s\\
+\frac{4a\cos\th}{|q|^6}\Big(a\cos\th\big(|q|^2+6mr\big) - is\big((r-m)|q|^2+4mr^2\big)\Big)\pmb\phi_s &=& 0.
\eea
\end{lemma}

\begin{proof}
We start with the case $s=+2$. We have, see \eqref{eq:TME},
\beaa
\square_{+2}\phi_{+2,\text{NP}}  -  \frac{2}{|q|^2}\phi_{+2,\text{NP}} 
+\frac{4}{|q|^2}[(r-m)e_3-2r\partial_t]\phi_{+2,\text{NP}}=0.
\eeaa
In view of \eqref{eq:linkteukolskyscalarsandtensors}, we may apply Corollary \ref{cor:howtogofromequationsinNPformtotensorialequations} and 
Lemma \ref{lemma:scalarizationoncomplexscalarsofnabprtannabprphi}. This yields
\beaa
0 &=& \left(\squared_2\pmb\phi_{+2} -\frac{4ia\cos\th}{|q|^2}\nab_{\pr_t}\pmb\phi_{+2}\right)(e_1, e_1) +\frac{4a^2\cos^2\th}{|q|^6}\Big(|q|^2+6mr\Big)\pmb\phi_{+2}(e_1, e_1)\\
&& -  \frac{2}{|q|^2}\pmb\phi_{+2}(e_1, e_1)  +\frac{4}{|q|^2}(r-m)\left(\nab_3\pmb\phi_{+2}(e_1,e_1)-\frac{2ia\cos\th}{|q|^2}\pmb\phi_{+2}(e_1,e_1)\right)\\
&& -\frac{8r}{|q|^2}\left(\nab_{\pr_t}\pmb\phi_{+2}(e_1, e_1)  +\frac{4iamr\cos\th}{|q|^4}\pmb\phi_{+2}(e_1, e_1)\right)\\
&=& \Bigg[\left(\squared_2 -\frac{4ia\cos\th}{|q|^2}\nab_{\pr_t} -  \frac{2}{|q|^2}\right)\pmb\phi_{+2} +\frac{4}{|q|^2}(r-m)\nab_3\pmb\phi_{+2}  -\frac{8r}{|q|^2}\nab_{\pr_t}\pmb\phi_{+2}\\
&& +\left(\frac{4a^2\cos^2\th}{|q|^6}\Big(|q|^2+6mr\Big) - \frac{8ia(r-m)\cos\th}{|q|^4}  -\frac{32iamr^2\cos\th}{|q|^6}\right)\pmb\phi_{+2}\Bigg](e_1, e_1)
\eeaa
and hence
\beaa
\Bigg[\left(\squared_2 -\frac{4ia\cos\th}{|q|^2}\nab_{\pr_t} -  \frac{2}{|q|^2}\right)\pmb\phi_{+2} +\frac{4}{|q|^2}(r-m)\nab_3\pmb\phi_{+2}  -\frac{8r}{|q|^2}\nab_{\pr_t}\pmb\phi_{+2}\\
+\frac{4a\cos\th}{|q|^6}\Big(a\cos\th\big(|q|^2+6mr\big) - 2i\big((r-m)|q|^2+4mr^2\big)\Big)\pmb\phi_{+2}\Bigg](e_1, e_1) &=& 0.
\eeaa
Applying again Corollary \ref{cor:howtogofromequationsinNPformtotensorialequations}, we infer
\beaa
\left(\squared_2 -\frac{4ia\cos\th}{|q|^2}\nab_{\pr_t} -  \frac{2}{|q|^2}\right)\pmb\phi_{+2} +\frac{4}{|q|^2}(r-m)\nab_3\pmb\phi_{+2}  -\frac{8r}{|q|^2}\nab_{\pr_t}\pmb\phi_{+2}\\
+\frac{4a\cos\th}{|q|^6}\Big(a\cos\th\big(|q|^2+6mr\big) - 2i\big((r-m)|q|^2+4mr^2\big)\Big)\pmb\phi_{+2} &=& 0
\eeaa
as stated.

Next, we consider the case $s=-2$. We have, see \eqref{eq:TME},
\beaa
\square_{-2}\phi_{-2,\text{NP}}  +\frac{2}{|q|^2}\phi_{-2,\text{NP}} 
-\frac{4}{|q|^2}[(r-m)e_3-2r\partial_t]\phi_{-2,\text{NP}}=0,
\eeaa
or
\beaa
\ov{\square_{-2}\phi_{-2,\text{NP}}}  +\frac{2}{|q|^2}\ov{\phi_{-2,\text{NP}}} 
-\frac{4}{|q|^2}[(r-m)e_3-2r\partial_t]\ov{\phi_{-2,\text{NP}}}=0.
\eeaa
In view of \eqref{eq:linkteukolskyscalarsandtensors}, we may apply Corollary \ref{cor:howtogofromequationsinNPformtotensorialequations} and 
Lemma \ref{lemma:scalarizationoncomplexscalarsofnabprtannabprphi}. This yields
\beaa
0 &=& \left(\squared_2\pmb\phi_{-2} -\frac{4ia\cos\th}{|q|^2}\nab_{\pr_t}\pmb\phi_{-2}\right)(e_1, e_1) +\frac{4a^2\cos^2\th}{|q|^6}\Big(|q|^2+6mr\Big)\pmb\phi_{-2}(e_1, e_1)\\
&& +  \frac{2}{|q|^2}\pmb\phi_{-2}(e_1, e_1)  -\frac{4}{|q|^2}(r-m)\left(\nab_3\pmb\phi_{-2}(e_1,e_1)-\frac{2ia\cos\th}{|q|^2}\pmb\phi_{-2}(e_1,e_1)\right)\\
&& +\frac{8r}{|q|^2}\left(\nab_{\pr_t}\pmb\phi_{-2}(e_1, e_1)  +\frac{4iamr\cos\th}{|q|^4}\pmb\phi_{-2}(e_1, e_1)\right)\\
&=& \Bigg[\left(\squared_2 -\frac{4ia\cos\th}{|q|^2}\nab_{\pr_t} +  \frac{2}{|q|^2}\right)\pmb\phi_{-2} -\frac{4}{|q|^2}(r-m)\nab_3\pmb\phi_{-2}  +\frac{8r}{|q|^2}\nab_{\pr_t}\pmb\phi_{-2}\\
&& +\left(\frac{4a^2\cos^2\th}{|q|^6}\Big(|q|^2+6mr\Big) + \frac{8ia(r-m)\cos\th}{|q|^4}  +\frac{32iamr^2\cos\th}{|q|^6}\right)\pmb\phi_{-2}\Bigg](e_1, e_1)
\eeaa
and hence
\beaa
\Bigg[\left(\squared_2 -\frac{4ia\cos\th}{|q|^2}\nab_{\pr_t} +  \frac{2}{|q|^2}\right)\pmb\phi_{-2} -\frac{4}{|q|^2}(r-m)\nab_3\pmb\phi_{+2}  +\frac{8r}{|q|^2}\nab_{\pr_t}\pmb\phi_{-2}\\
+\frac{4a\cos\th}{|q|^6}\Big(a\cos\th\big(|q|^2+6mr\big) + 2i\big((r-m)|q|^2+4mr^2\big)\Big)\pmb\phi_{-2}\Bigg](e_1, e_1) &=& 0.
\eeaa
Applying again Corollary \ref{cor:howtogofromequationsinNPformtotensorialequations}, we infer
\beaa
\left(\squared_2 -\frac{4ia\cos\th}{|q|^2}\nab_{\pr_t} +  \frac{2}{|q|^2}\right)\pmb\phi_{-2} -\frac{4}{|q|^2}(r-m)\nab_3\pmb\phi_{+2}  +\frac{8r}{|q|^2}\nab_{\pr_t}\pmb\phi_{-2}\\
+\frac{4a\cos\th}{|q|^6}\Big(a\cos\th\big(|q|^2+6mr\big) + 2i\big((r-m)|q|^2+4mr^2\big)\Big)\pmb\phi_{-2} &=& 0
\eeaa
as stated. This concludes the proof of Lemma \ref{lemma:Teukolskyintensorialform}.
\end{proof}

Next, we introduce $\pmb\phi_s^{(p)}\in\sk_2(\mathbb{C})$ defined as 
\bsub\lab{def:TensorialTeuScalars:wavesystem:Kerr}
 \bea
\pmb\phi_{+2}^{(0)}:= \frac{1}{|q|^4}\pmb\phi_{+2},\qquad \pmb\phi_{+2}^{(p+1)}:=\left(\frac{r^2}{|q|^2}\right)^{1-p}\left({\frac{q}{\bar{q}}}r \nab_3 r{\frac{\bar{q}}{q}}\right)\left(\frac{r^2}{|q|^2}\right)^{p-2}\pmb\phi_{+2}^{(p)}, \quad p=0,1,
\eea
and
\bea
\pmb\phi_{-2}^{(0)}:=\frac{\De^2}{|q|^4}\pmb\phi_{-2},\quad\,\,\, \pmb\phi_{-2}^{(p+1)}:=\left(\frac{{r^2}}{|q|^2}\right)^{1-p}\left(\frac{\bar{q}}{q}r\frac{\qs}{\De}\nab_4 r\frac{q}{\bar{q}}\right)\left(\frac{{r^2}}{|q|^2}\right)^{p-2}\pmb\phi_{-2}^{(p)}, \,\,\,\, p=0,1.
\eea
\esub
Note, in view of \eqref{def:ComplexTeuScalars:wavesystem:Kerr}, \eqref{eq:linkteukolskyscalarsandtensors} and Corollary \ref{cor:howtogofromequationsinNPformtotensorialequations}, that we have
\bea\lab{eq:linkteukolskyscalarsandtensors:p=012case}
\phi^{(p)}_{+2,\text{NP}}=\pmb\phi^{(p)}_{+2}(e_1, e_1), \qquad \phi^{(p)}_{-2,\text{NP}}=\ov{\pmb\phi^{(p)}_{-2}(e_1, e_1)}, \quad p=0,1,2.
\eea

The following lemma provides the tensorial wave equations satisfied by $\pmb\phi^{(p)}_s$, $s=\pm 2$, $p=0,1,2$.  
\begin{lemma}\lab{lemma:derivationoftensorialwaveequationforpmbphispfromNPcase}
Let $\pmb\phi^{(p)}_s$, $s=\pm 2$, $p=0,1,2$, in $\sk_2(\mathbb{C})$ be given by \eqref{def:TensorialTeuScalars:wavesystem:Kerr}. Then, they satisfy 
\bsub
\lab{eq:TensorialTeuSysandlinearterms:rescaleRHScontaine2:general:Kerr}
\bea
\lab{eq:TensorialTeuSys:rescaleRHScontaine2:general:Kerr}
\bigg(\squared_2 -\frac{4ia\cos\th}{|q|^2}\nab_{\pr_{\tt}}- \frac{4-2\de_{p0}}{\qs}\bigg){\phis{p}} &=& \L_{s}^{(p)}[\pmb\phi_{s}], \quad s=\pm2, \quad p=0,1,2,
\eea
where the linear coupling terms $\L_{s}^{(p)}[{\pmb\phi_s}]$ have the following schematic forms
\bea
\lab{eq:tensor:Lsn:onlye_2present:general:Kerr}
\bsplit
{\L_{s}^{(0)}[\pmb\phi_{s}]}={}& (2sr^{-3} +O(mr^{-4}))\phis{1}+ O(mr^{-3}) \nab_{\Xcal_s}^{\leq 1}\phis{0},\\
{\L_{s}^{(1)}[\pmb\phi_{s}]}={}& (sr^{-3} +O(mr^{-4}))\phis{2}+ O(mr^{-3}) \nab_{\Xcal_s}^{\leq 1}  \phis{1}+O(mr^{-2})\nab_{\pr_{\tphi}+a\pr_{\tt}}^{\leq 1}\phis{0},\\
{\L_{s}^{(2)}[\pmb\phi_{s}]}={}&O(mr^{-3})\phis{2}+O(mr^{-2})\nab_{\pr_{\tphi}+a\pr_{\tt}}^{\leq 1}\phis{1}+O(m^2 r^{-2})\phis{0},
\end{split}
\eea
\esub
with $\Xcal_s$, $s=\pm 2$, being the regular horizontal vectorfields introduced in \eqref{eq:formofregularhorizontalvectorfieldmathcalXs}, with all the coefficients in  \eqref{eq:tensor:Lsn:onlye_2present:general:Kerr} being independent of coordinates $\tau$ and 
$\tphi$, and with the coefficients in front of the terms $\phis{2}$ and $\nab_{\pr_{\tphi}+a\pr_{\tt}}\phis{1}$ on the RHS of equation of ${\L_{s}^{(2)}[\pmb\phi_{s}]}$ in \eqref{eq:tensor:Lsn:onlye_2present:general:Kerr} being real functions. 
\end{lemma}

\begin{remark}
The equations \eqref{eq:TensorialTeuSysandlinearterms:rescaleRHScontaine2:general:Kerr} and \eqref{def:TensorialTeuScalars:wavesystem:Kerr} are the Teukolsky wave and transport equations in Kerr in the tensorial form, respectively.
\end{remark}

\begin{proof}
We associate to $\pmb\psi\in\sk_2(\mathbb{C})$ the complex-valued scalar $\psi_{\pm 2,\text{NP}}$ as in \eqref{def:psiplusminus2:complexscalars:Kerr}, i.e., 
\beaa
\psi_{+2,\text{NP}} := \pmb\psi(e_1, e_1), \qquad \psi_{-2,\text{NP}} := \ov{\pmb\psi(e_1, e_1)}.
\eeaa
Then, in view of Corollary \ref{cor:howtogofromequationsinNPformtotensorialequations}, we have
\beaa
\square_{+2}(\psi_{+2,\text{NP}}) &=& \left(\squared_2\pmb\psi -\frac{4ia\cos\th}{|q|^2}\nab_{\pr_t}\pmb\psi\right)(e_1, e_1) +O(a^2r^{-4})\psi_{+2,\text{NP}},\\
\ov{\square_{-2}(\psi_{-2,\text{NP}})} &=& \left(\squared_2\pmb\psi -\frac{4ia\cos\th}{|q|^2}\nab_{\pr_t}\pmb\psi\right)(e_1, e_1) +O(a^2r^{-4})\ov{\psi_{-2,\text{NP}}}.
\eeaa
Also, in view of Lemma \ref{lemma:scalarizationoncomplexscalarsofnabprtannabprphi}, we have
\beaa
\pr_t(\pmb\psi(e_1, e_1)) &=& \nab_{\pr_t}\pmb\psi(e_1, e_1)  +O(amr^{-3})\pmb\psi(e_1, e_1),\\
\pr_\phi(\pmb\psi(e_1, e_1)) &=& \nab_{\pr_\phi}\pmb\psi(e_1, e_1)+O(1)\pmb\psi(e_1, e_1),
\eeaa
and, since $\g(\D_1e_1, e_2)=0$, we have
\beaa
\sin\th\pr_\th(\pmb\psi(e_1, e_1)) = \nab_{\sin\th\pr_\th}\pmb\psi(e_1, e_1),
\eeaa
and hence, recalling the form of the regular horizontal vectorfields $\mathcal{X}_s$, $s=\pm 2$, see \eqref{eq:formofregularhorizontalvectorfieldmathcalXs}, we infer
\beaa
\mathcal{X}_s(\pmb\psi(e_1, e_1)) &=& \nab_{\mathcal{X}_s}\pmb\psi(e_1, e_1)+O(1)\pmb\psi(e_1, e_1).
\eeaa
The system of tensorial wave equations \eqref{eq:TensorialTeuSysandlinearterms:rescaleRHScontaine2:general:Kerr} for $\pmb\phi^{(p)}_s$ then follows immediately from the above identities, \eqref{eq:linkteukolskyscalarsandtensors:p=012case}, and the system of spin-weighted wave equations \eqref{eq:Teukolskywavetransportsystem:NPnewform} \eqref{eq:Teukolskywavetransportsystem-RHSLpsofpsips:NPnewform} for $\phi_{s,\text{NP}}^{(p)}$. Finally, the fact that the coefficients in front of the terms $\phis{2}$ and $\nab_{\pr_{\tphi}+a\pr_{\tt}}\phis{1}$ on the RHS of equation of ${\L_{s}^{(2)}[\pmb\phi_{s}]}$ in \eqref{eq:tensor:Lsn:onlye_2present:general:Kerr} are real functions follows from the property that the coefficients in front of the terms $\phi_{s,\text{NP}}^{(2)}$ and of $(\pr_\phi+a\pr_t)\phi_{s,\text{NP}}^{(1)}$ on the RHS of the expression of $L^{(2)}_{s,\text{NP}}[\phi_{s,\text{NP}}]$ in \eqref{eq:Teukolskywavetransportsystem-RHSLpsofpsips:NPnewform} are real functions. This concludes the proof of Lemma \ref{lemma:derivationoftensorialwaveequationforpmbphispfromNPcase}.
\end{proof}


\subsection{Scalarization of Teukolsky  equations in Kerr using  regular triplets}


The two forms of the Teukolsky transport/wave systems introduced respectively in Sections  \ref{sec:TeukolskyinKerrintensorialform} and \ref{sec:Teukoslkytransportwavesystemwithcomplexscalars} will not be suitable to prove our main result for the following reasons:
\begin{itemize}
\item Part of this paper will rely on microlocal methods adapted to the $r$-foliation of the spacetime $\MM$ which will in turn force us to extend the Teukolsky  transport/wave systems from finite time intervals to global in time solutions. This procedure would be complicated to perform for the tensorial formulation of Section   \ref{sec:TeukolskyinKerrintensorialform} and turns out to be significantly simpler for scalar equations.  
\item The extension to perturbations of Kerr of the formulation of Section  \ref{sec:Teukoslkytransportwavesystemwithcomplexscalars}, relying on the scalarization \eqref{def:psiplusminus2:complexscalars:Kerr} of the tensorial equations based on the frame $(e_1, e_2)$, seems a priori problematic due to the fact that $(e_1, e_2)$ is singular at $\th=0,\pi$. 
\end{itemize}
The above observations suggest to apply an alternative way to scalarize the tensorial Teukolsky transport/wave systems of Section  \ref{sec:TeukolskyinKerrintensorialform} which, unlike the one in Section  \ref{sec:Teukoslkytransportwavesystemwithcomplexscalars}, is regular everywhere on $\MM$. To this end, instead of using the irregular basis $(e_1, e_2)$, we will scalarize the Teukolsky equations  using a family of regular triplets $(\Om_i)_{i=1,2,3}$ as introduced in Definition \ref{def:definitionofregulartripletOmii=123}.

Applying Lemma \ref{lemma:formoffirstordertermsinscalarazationtensorialwaveeq} to the tensorial Teukolsky wave/transport systems \eqref{eq:TensorialTeuSysandlinearterms:rescaleRHScontaine2:general:Kerr}  \eqref{def:TensorialTeuScalars:wavesystem:Kerr}, we deduce the scalarized Teukolsky wave/transport systems based on regular triplets.

\begin{lemma}[Scalarized Teukolsky wave/transport systems in Kerr using regular triplets]
\lab{lem:scalarizedTeukolskywavetransportsysteminKerr:Omi}
Let a regular triplet $(\Om_i)_{i=1,2,3}$ be in the sense of Definition \ref{def:definitionofregulartripletOmii=123}. Then, the Teukolsky wave system \eqref{eq:TensorialTeuSysandlinearterms:rescaleRHScontaine2:general:Kerr} scalarized using $(\Om_i)_{i=1,2,3}$ takes the following form 
\bea
\lab{eq:ScalarizedTeuSys:general:Kerr} 
\widehat\square_{\gam}(\phi_s^{(p)})_{ij} -\frac{4-2\de_{p0}}{\qs} \phiss{ij}{p} = {L_{s,ij}^{(p)}}, \quad s=\pm2, \quad p=0,1,2,
\eea
where we have defined
\bea
\widehat\square_{\gam}(\phi_s^{(p)})_{ij}  :=\square_{\gam}\phiss{ij}{p}- \widehat{S}(\phi_s^{(p)})_{ij} -(\widehat{Q}\phi_s^{(p)})_{ij}
\eea
with
\bsub
\begin{align}
\widehat{S}(\phi_s^{(p)})_{ij} ={}&S(\phi_s^{(p)})_{ij} +\frac{4ia\cos\th}{|q|^2} \pr_t\phiss{ij}{p},\\
(\widehat{Q}\phi_s^{(p)})_{ij} ={}& (Q\phi_s^{(p)})_{ij} 
-\frac{4ia\cos\th}{|q|^2}\big(M_{it}^l \phiss{lj}{p}+M_{jt}^l \phiss{il}{p}\big),
\end{align}
\esub
where the linear coupling terms ${L_{s,ij}^{(p)}}=(\L_{s}^{(p)}[\pmb\phi_s])_{ij}$ are given by
\bsub
\lab{eq:linearterms:ScalarizedTeuSys:general:Kerr}
\bea
{L_{s,ij}^{(0)}}&=& (2sr^{-3} +O(mr^{-4}))\phiss{ij}{1} + O(mr^{-3}){\Xcal_s}\phiss{ij}{0}+\sum_{k,l=1,2,3}O(mr^{-3})\phiss{kl}{0},\\
{L_{s,ij}^{(1)}}&=&(sr^{-3} +O(mr^{-4}))\phiss{ij}{2} + O(mr^{-3}){\Xcal_s}\phiss{ij}{1}+O(mr^{-2})(\pr_\phi +a\pr_t)\phiss{ij}{0}\nn\\
  &&+\sum_{k,l=1,2,3}\Big(O(mr^{-3})\phiss{kl}{1}+O(mr^{-2})\phiss{kl}{0}\Big),\\
{L_{s,ij}^{(2)}}&=&O(mr^{-3}) \phiss{ij}{2}+O(mr^{-2})(\pr_\phi +a\pr_t)\phiss{ij}{1}\nn\\
&&+\sum_{k,l=1,2,3}O(mr^{-2})\phiss{kl}{1}+O(m^2r^{-2})\phiss{ij}{0}
\eea
\esub
with ${\Xcal_s}$ defined as in \eqref{eq:formofregularhorizontalvectorfieldmathcalXs},  {and with the coefficients in front of the terms $\phiss{ij}{2}$ and $(\pr_\phi+a\pr_t)\phiss{ij}{1}$ on the RHS of equation of $L_{s,ij}^{(2)}$ in \eqref{eq:linearterms:ScalarizedTeuSys:general:Kerr} being real functions,} where $S(\phi_s^{(p)})_{ij}$ and $(Q\phi_s^{(p)})_{ij}$ are given as in \eqref{SandV}, and where 
the complex-valued scalars $\phiss{ij}{p}$, $s=\pm 2$, $p=0, 1,2$, $i,j=1,2,3$, satisfy the scalarized Teukolsky transport system of equations
\bsub
 \lab{eq:ScalarizedQuantitiesinTeuSystem:Kerr}
\bea
\bsplit
\phipluss{ij}{0}:=&{\frac{1}{|q|^4}\pmb\phi_{+2}}(\Om_i,\Om_j),\\
\phipluss{ij}{p+1}:=&\bigg({\frac{r^2}{|q|^2}}\bigg)^{1-p}\left({\frac{q}{\bar{q}}}r e_3r{\frac{\bar{q}}{q}}\right)\bigg(\bigg({\frac{r^2}{|q|^2}}\bigg)^{p-2}\phipluss{ij}{p} \bigg)\\
&- {|q|^2}\Big(M_{i3}^k \phipluss{kj}{p} + M_{j3}^k \phipluss{ik}{p}\Big)
\end{split}
\eea
and 
\bea
\bsplit
\phiminuss{ij}{0}:=&\frac{\De^2}{|q|^4}{\pmb\phi_{-2}}(\Om_i,\Om_j),\\
\phiminuss{ij}{p+1}:=&\left(\frac{r^2}{|q|^2}\right)^{1-p}\left({\frac{\bar{q}}{q}}r \frac{\qs}{\De}e_4 r {\frac{q}{\bar{q}}}\right) \bigg(\bigg(\frac{r^2}{|q|^2}\bigg)^{p-2}\phiminuss{ij}{p}\bigg)\\
&-\frac{|q|^4}{\Delta} \Big(M_{i4}^k \phiminuss{kj}{p} + M_{j4}^k \phiminuss{ik}{p}\Big).
\end{split}
\eea
\esub
\end{lemma}

\begin{remark}
As outlined at the beginning of this section, the advantages of this formulation of the Teukolsky wave/transport systems \eqref{eq:ScalarizedTeuSys:general:Kerr} and \eqref{eq:ScalarizedQuantitiesinTeuSystem:Kerr}  are twofold:
\begin{itemize}
\item the systems are for scalars (as opposed to horizontal tensors), which makes it easily amenable to extensions from local in time to global in time problems and using microlocal calculus;
\item the formulation with the scalarization {using a regular triplet $\Om_i$, $i=1,2,3$,} is regular everywhere unlike the formulation using the NP formalism for complex scalars. 
\end{itemize}
\end{remark}

\begin{proof}
We project the Teukolsky wave system \eqref{eq:TensorialTeuSysandlinearterms:rescaleRHScontaine2:general:Kerr} {using a regular triplet $(\Om_i)_{i=1,2,3}$}, and in view of Lemma \ref{lemma:formoffirstordertermsinscalarazationtensorialwaveeq}, we infer the LHS equals
\begin{align*}
&\bigg(\bigg(\squared_2-\frac{4ia\cos\th}{|q|^2} \nab_{\pr_t}- \frac{4-2\de_{p0}}{\qs}\bigg)\phis{p}\bigg)(\Om_i,\Om_j)\nn\\
=&\big(\squared_2\phis{p}\big)(\Om_i,\Om_j)-\frac{4ia\cos\th}{|q|^2}\big( \nab_{\pr_t}\phis{p}\big)(\Om_i,\Om_j)- \frac{4-2\de_{p0}}{\qs}\phis{p}(\Om_i,\Om_j)\nn\\
=&\big(\squared_2\phis{p}\big)_{ij}-\frac{4ia\cos\th}{|q|^2}\pr_t\phiss{ij}{p}+\frac{4ia\cos\th}{|q|^2}\phis{p}(\nab_{\pr_t}\Om_i, \Om_j)
+\frac{4ia\cos\th}{|q|^2}\phis{p}(\Om_i, \nab_{\pr_t}\Om_j)\\
&- \frac{4-2\de_{p0}}{\qs}\phiss{ij}{p}\nn\\
=&\square_{\gam}(\phiss{ij}{p})- S(\phi_s^{(p)})_{ij} -(Q\phi_s^{(p)})_{ij}  -\frac{4ia\cos\th}{|q|^2}\pr_t\phiss{ij}{p} 
+\frac{4ia\cos\th}{|q|^2}\big(M_{it}^l \phiss{lj}{p}+M_{jt}^l \phiss{il}{p}\big)\\
&- \frac{4-2\de_{p0}}{\qs}\phiss{ij}{p}\nn\\
=&\widehat\square_{\gam}(\phi_s^{(p)})_{ij} -\frac{4-2\de_{p0}}{\qs} \phiss{ij}{p},
\end{align*}
which hence proves equation \eqref{eq:ScalarizedTeuSys:general:Kerr}. The formulas \eqref{eq:linearterms:ScalarizedTeuSys:general:Kerr} are manifest by expanding out ${L_{s,ij}^{(p)}}=(\L_{s}^{(p)}[\pmb\phi_s])_{ij}$ using \eqref{eq:tensor:Lsn:onlye_2present:general:Kerr} {and Lemma \ref{lemma:introductionandpropertiesoftheMalphaij}, as well as the fact that 
\beaa
M_{i\a}^j(\Xcal_s)^\a=O(1), \qquad M_{i\a}^j(\pr_\phi+a\pr_t)^\a=O(1),
\eeaa
in view of Lemma \ref{lemma:computationoftheMialphajinKerr}. Finally}, the formulas \eqref{eq:ScalarizedQuantitiesinTeuSystem:Kerr} follow from projecting both sides of the Teukolsky transport system 
\eqref{def:TensorialTeuScalars:wavesystem:Kerr} onto the basis $(\Om_i)_{i=1,2,3}$ and using {Lemma \ref{lemma:introductionandpropertiesoftheMalphaij}}.
\end{proof}


\section{Teukolsky equations in perturbations of Kerr}
\lab{sect:TeuinKerrperturbation}


The goal of this section is to introduce the Teukolsky wave/transport system in perturbations of Kerr which will be studied in this paper. To this end, we first provide our main assumptions on our perturbed Kerr spacetime $(\MM, \g)$. The Teukolsky wave/transport system  in $(\MM, \g)$ is then introduced in Section  \ref{sec:TeukolskyWavesysteminperturbationsofKerr}. Finally, we conclude the section by a discussion of future null infinity and energy-Morawetz norms.


\subsection{Choices of constants}
\lab{sec:smallnesconstants}


The following constants  are  involved in  the statement and in the proof of our main result:
\begin{itemize}
\item The constants $m>0$ and $a$, with $|a|<m$, are the mass and the angular momentum per unit mass of the Kerr solution relative to which the perturbation of the metric $\g$ is measured. 

\item The size of the metric perturbation  is measured by $\ep\geq 0$. 

\item The constant $\dhor$ is tied to the boundary of $\MM$ given by $\pr\MM=\AA=\{r=r_+(1-\dhor)\}$. 

\item The constant $\dred$ measures the width of the redshift region.

\item The constant $\dbl$ appears in the construction of normalized coordinates, see Lemma \ref{lem:specificchoice:normalizedcoord}.

\item The constant $\dec$ is tied to decay estimates in $(r, \tau)$ of the perturbed metric coefficients, see Section \ref{subsubsect:assumps:perturbedmetric}. 

\item The constant $\de$ is tied to $r$-weights in the Morawetz norm {$\M_\de[\psi]$}, see \eqref{def:variousMorawetzIntegrals},

\item The large integer $\Nmic$ is tied to the choice of a contant\footnote{The contant $\Rmic$ is used to define the region $\MM_{r\leq\Rmic}$ on which we will derive microlocal energy-Morawetz estimates in Section  \ref{sect:microlocalenergyMorawetztensorialwaveequation}.} $\Rmic\in [\Nmic m, (\Nmic+1)m]$, see Remark \ref{rmk:choiceofconstantRbymeanvalue}.
\end{itemize}

These  constants are chosen such that \bea\lab{eq:constraintsonthemainsmallconstantsepanddelta}
0<\ep\ll\dhor\ll\dred\ll \dbl \ll 1-\frac{|a|}{m}, \qquad \ep\ll \de\leq \frac{1}{3}, \qquad \ep\ll \dec, \qquad \ep\ll \Nmic^{-1}.
\eea

From now on, in the rest of the paper, $\lesssim$ means bounded by a positive constant multiple, with this positive constant depending only on universal constants (such as constants arising from Sobolev embeddings, elliptic estimates,...) as well as the constants 
$$m,\,\, a, \,\, \dhor,\,\, \dred,\,\, \dbl, \,\, \dec,\,\,\de,\,\,\Nmic, $$
\textit{but not on} $\ep$. {Also, note that the} constants $\dhor, \dred$ and $\dbl$ can be {chosen} to be only dependent on $m$ and $a$.

Throughout this paper, ``LHS" and ``RHS" are {abbreviations} for ``left-hand side" and ``right-hand side", respectively,  ``w.r.t." is {an abbreviation} for ``with respect to", ``EMF" is {an abbreviation} for ``energy-Morawetz-flux",  and $\Re(\cdot)$ and $\Im(\cdot)$ mean taking the real part and the imaginary part, respectively.


\subsection{Subregions and hypersurfaces of $\MM$}
\lab{sect:subregionsandhypersurfaces}


Let $(\MM, \g)$ be a four dimensional Lorentzian manifold covered by coordinate systems $(\tt, r, x^1_0, x^2_0)$ and {$(\tt, r, x^1_p, x^2_p)$}, defined respectively on $\th\neq 0, \pi$ and $\th\neq\frac{\pi}{2}$, with 
\beaa
\tau\in\mathbb{R}, \quad r_+(1-\dhor)\leq r<+\infty, \quad x^1_0=\th, \quad x^2_0=\tphi, \quad {x^1_p=\sin\th\cos\tphi, \quad x^2_p=\sin\th\sin\tphi}.
\eeaa

We define a few subregions and hypersurfaces of $\MM$. 

 \begin{definition}
 Define the following subregions and hypersurfaces of $\MM$:
 \bsub
 \begin{align}
 \MM(\tt_1,\tt_2):={}&\MM\cap\{\tt_1\leq \tt\leq \tt_2\}, \quad \forall\tau_1<\tau_2,\\
\MM_{r_1,r_2}:={}&\MM\cap \{r_1\leq r\leq r_2\}, \quad \forall r_+(1-\dhor)\leq r_1<r_2,\\
 \Sigma(\tt_1):={}&\MM\cap\{\tt=\tt_1\}, \quad \forall \tt_1\in \Reals,\\
  \Sigma_{r_1,r_2}(\tt_1):={}& \Sigma(\tt_1)\cap \{r_1\leq r\leq r_2\}, \quad \forall \tt_1\in \Reals, \forall r_+(1-\dhor)\leq r_1<r_2,\\
 H_{r_1}:={}&\MM\cap \{r=r_1\},\quad  \forall r_1\geq r_+(1-\dhor), \\
 \AA:={}&\MM\cap\{r=r_+(1-\dhor)\},\\
 \MM_{red}:={}&\MM\cap\{r\leq r_+(1+\dred)\},\\
 \Mtrap:={}&\MM_{r_+(1+2\dbl), 10m},\\ 
 \Mntrap:={}&\MM\setminus\Mtrap.
 \end{align}
 \esub
 \end{definition}

 
 \subsection{Assumptions on the null pair $(e_3, e_4)$ and consequences}
 \label{subsect:assumps:perturbednullframe}
 

We assume that $(e_3, e_4)$ is a null pair defined on $\MM$ and we consider its corresponding horizontal structure, Ricci coefficients and curvature components as introduced in Section  \ref{sec:nonintergrableformalism}. We also assume the existence on $\MM$ of complex horizontal 1-forms $\Jk,\,\Jk_\pm\in\sk_1(\mathbb{C})$, and, using the coordinates $(r, \th)$ on $\MM$, we define the complex-valued scalar $q$ by 
\beaa
q:=r+ia\cos\th.
\eeaa 
Finally, we assume that $\Jk$ and $\Jk_\pm$ satisfy the same algebraic identities as in \eqref{eq:usefulalgebraicidentitiesinvolvingscalarproductsReJkReJkpm}, i.e.,
\bea\lab{eq:usefulalgebraicidentitiesinvolvingscalarproductsReJkReJkpm:Kerrpert}
\bsplit
\Re(\Jk)\c\Re(\Jk) &=\frac{(\sin\th)^2}{|q|^2}, \qquad \Re(\Jk)\c\Re(\Jk_+)=-\frac{x^2}{|q|^2},\qquad \Re(\Jk)\c\Re(\Jk_-)=\frac{x^1}{|q|^2},\\
\dual(\Re(\Jk))\c\Re(\Jk_+) &=\frac{\cos\th x^1}{|q|^2},\qquad \dual(\Re(\Jk))\c\Re(\Jk_-)=\frac{\cos\th x^2}{|q|^2}.
\end{split}
\eea

\begin{remark}\lab{rmk:howtoenforceusefulalgebraicidentitiesinvolvingscalarproductsReJkReJkpm}
Given a spacetime $(\MM, \g)$, a null pair $(e_3, e_4)$ and the corresponding horizontal structure $\O(\MM)$, one can easily generate complex horizontal 1-forms $\Jk,\,\Jk_\pm\in\sk_1(\mathbb{C})$ verifying \eqref{eq:usefulalgebraicidentitiesinvolvingscalarproductsReJkReJkpm:Kerrpert} by enforcing these identities 
 on one given topological sphere in $\MM$, which can then be propagated to $\MM$ by defining $\Jk$, $\Jk_\pm$ based on well-chosen transport equations consistent with the horizontal structure of $\MM$. 
\end{remark}

 
 \subsubsection{Definition of linearized quantities}
 

Recall that the constants $m>0$ and $a$, with $|a|<m$, are the mass and the angular momentum per unit mass of the Kerr solution relative to which the perturbation of the metric $\g$ is measured. In view of the Kerr values in Section  \ref{subsect:principalnullpairinKerr}, we introduce the following linearized quantities.

\begin{definition}
\lab{def:renormalizationofallnonsmallquantitiesinPGstructurebyKerrvalue}
We  define  the following renormalizations.
\begin{enumerate}
\item Linearization of the complex-valued Ricci and curvature coefficients:
\beaa
\bsplit
\trXc &:= \tr X-\frac{2\ov{q}\De}{|q|^4}, \qquad\trXbc := \tr\Xb+\frac{2}{\ov{q}},\\ 
\Pc &:= P+\frac{2m}{q^3},\qquad\qquad\, \omc  := \om  + \frac{1}{2}\pr_r\left(\frac{\De}{|q|^2} \right),\\
 \Hc &:= H-\frac{aq}{|q|^2}\Jk, \qquad\quad \Hbc:=\Hb+\frac{a\ov{q}}{|q|^2}\Jk,\qquad\quad \Zc := Z-\frac{aq}{|q|^2}\Jk.
 \end{split}
\eeaa

\item Linearization of derivatives of the scalar functions $r$, $\cos\th$, $q$, $\tau$, $x^1_p$ and $x^2_p$:
\beaa
\bsplit
\widecheck{e_3(r)} :=& e_3(r)+1, \qquad\qquad\qquad\,\,\,\,\, \widecheck{e_4(r)} := e_4(r)-\frac{\Delta}{|q|^2},\\
\widecheck{e_3(\tau)}:=& e_3(\tau)-\tmod'(r), \qquad\qquad\, \widecheck{e_4(\tau)}:=e_4(\tau)-\frac{2(r^2+a^2) - \De\tmod'(r)}{|q|^2},\\
\widecheck{e_3(x^1_p)}:=& e_3(x^1_p)+\phimod'(r)x^2_p, \qquad \widecheck{e_4(x^1_p)}:=  e_4(x^1_p)+\frac{2a -\De\phimod'(r)}{|q|^2}x^2_p, \\ \widecheck{e_3(x^2_p)}:=& e_3(x^2_p)-\phimod'(r)x^1_p, \qquad \widecheck{e_4(x^2_p)}:=e_4(x^2_p)-\frac{2a -\De\phimod'(r)}{|q|^2}x^1_p,\\
\widecheck{\DD q} :=& \DD q+a\Jk, \qquad\qquad\, \widecheck{\DD \ov{q}} :=\DD \ov{q}-a\Jk,\qquad \,\,\,\widecheck{\DD(\cos\th)} := \DD(\cos\th) -i\Jk,\\
\widecheck{\DD(\tau)}:=& \DD(\tau)-a\Jk, \qquad  \widecheck{\DD(x^1_p)} := \DD(x^1_p) - \Jk_+, \qquad \widecheck{\DD(x^1_p)} := \DD(x^2_p) - \Jk_-.
\end{split}
\eeaa

\item Linearization of derivatives of the complex 1-form $\Jk$:
\beaa
\bsplit
\widecheck{\nab_3\Jk}:=&\nab_3\Jk -\frac{1}{\ov{q}}\Jk, \qquad\qquad  \widecheck{\nab_4\Jk}:=\nab_4\Jk +\frac{\De \ov{q}}{|q|^4}\Jk,\qquad \widecheck{\ov{\DD}\c\Jk}:= \ov{\DD}\c\Jk-\frac{4i(r^2+a^2)\cos\th}{|q|^4},\\
\widecheck{\nab_3\Jk_\pm}:=&\nab_3\Jk_\pm - \frac{1}{\ov{q}}\Jk_\pm,\qquad \widecheck{\nab_4 \Jk_\pm}:=\nab_4 \Jk_\pm + \frac{\De \ov{q}}{|q|^4}\Jk_{\pm} \pm  \frac{2a}{|q|^2}\Jk_{\mp},\\
\widecheck{\ov{\DD}\c \Jk_+}:=&\ov{\DD}\c \Jk_+ + \frac{4r^2 }{|q|^4}x^1_p + \frac{4ia^2\cos\th}{|q|^4}x^2_p,\qquad \widecheck{\ov{\DD}\c \Jk_-}:=\ov{\DD}\c \Jk_- + \frac{4r^2 }{|q|^4}x^2_p - \frac{4ia^2\cos\th}{|q|^4}x^1_p.
\end{split}
\eeaa
 \end{enumerate}
\end{definition}

 
 \subsubsection{Notations $\Ga_b$ and $\Ga_g$ for error terms and assumptions}
 

\begin{definition}
\lab{definition.Ga_gGa_b}
The set of all linearized quantities is of the form $\Ga_g\cup \Ga_b$ with  $\Ga_g,  \Ga_b$
 defined as follows.
 \begin{enumerate}
\item 
 The set $\Ga_g$ is given by $\Ga_g=\Ga_{g,1}\cup \Ga_{g, 2}\cup\Ga_{g,3}$   with
 \bea
 \bsplit
 \Ga_{g,1} &= \Big\{\Xi, \quad \omc, \quad\trXc,\quad  \Xh,\quad \Zc,\quad \Hbc, \quad \trXbc , \quad r\Pc, \quad  rB, \quad  rA\Big\},\\
 \Ga_{g,2} &= \Big\{\widecheck{e_4(r)}, \quad r^{-1}\nab(r), \quad \widecheck{e_4(\tau)}, \quad r^{-1}\widecheck{\DD(\tau)}, \quad r^{-1}\widecheck{e_3(\tau)},  \quad e_4(\cos\th), \quad \widecheck{e_4(x^1_p)}, \quad \widecheck{e_4(x^2_p)}\Big\},\\
  \Ga_{g,3} &= \Big\{r\widecheck{\nab_4\Jk}, \quad r\widecheck{\nab_4\Jk_\pm}\Big\}.
 \end{split}
 \eea
 
 \item The set $\Ga_b$ is given by $\Ga_b=\Ga_{b,1}\cup \Ga_{b, 2}\cup \Ga_{b,3}$   with
 \bea
 \bsplit
 \Ga_{b,1}&= \Big\{\Hc, \quad \Xbh, \quad \omb, \quad \Xib,\quad  r\Bb, \quad \Ab\Big\},\\
  \Ga_{b, 2}&= \Big\{r^{-1}\widecheck{e_3(r)}, \quad  \widecheck{\DD(\cos\th)}, \quad e_3(\cos\th),  \quad  \widecheck{\DD(x^1_p)}, \quad {\widecheck{e_3(x^1_p)}}, \quad  \widecheck{\DD(x^2_p)}, \quad {\widecheck{e_3(x^2_p)}}\Big\}, \\
   \Ga_{b,3}&=\bigg\{ r\,\widecheck{\ov{\DD}\c\Jk}, \quad r\,\DD\hot\Jk, \quad r\,\widecheck{\nab_3\Jk}, \quad r\,\widecheck{\ov{\DD}\c\Jk_\pm}, \quad r\,\DD\hot\Jk_\pm, \quad r\,\widecheck{\nab_3\Jk_\pm}\bigg\}. 
   \end{split}
 \eea
\end{enumerate}
\end{definition}

We assume that the null pair $(e_3, e_4)$ and its associated horizontal structure satisfy on $\MM$ the following bounds, stated in terms of $\Ga_g$ and $\Ga_b$ introduced in Definition \ref{definition.Ga_gGa_b},
\bea\lab{eq:decaypropertiesofGabGag}
|\dk^{\leq {15}}\Ga_g|\les \min\left\{\frac{\ep}{r^2\tau^{\frac{1+\dec}{2}}}, \, \frac{\ep}{r\tau^{1+\dec}}\right\}, \qquad |\dk^{\leq {15}}\Ga_b|\les \frac{\ep}{r\tau^{1+\dec}},\qquad |\dk^{\leq {15}}\xi|\les \frac{\ep}{r^3},
\eea
where $\dec>0$, and where the weighted derivatives $\dk$ are defined by
\bea\lab{eq:defweightedderivative}
\dk:=\{\nab_3,\, r\nab_4,\, r\nabla\}.
\eea
For convenience, we also define the unweighted derivatives $\pr$ as follows
\bea\lab{eq:defunweightedderivative}
\pr:=\{\nab_3,\, \nab_4,\, \nabla\}.
\eea

\begin{remark}
The decay assumptions \eqref{eq:decaypropertiesofGabGag} are consistent with the decay estimates derived in the proof of the nonlinear stability of Kerr for small angular momentum in \cite{KS:Kerr}.
\end{remark}

\begin{remark}
In view of \eqref{eq:decaypropertiesofGabGag}, we note that $\Ga_g$ satisfies the assumptions of $\Ga_b$ and that $r^{-1}\Ga_b$ satisfies the assumptions of $\Ga_g$. Thus, in the rest of the paper, we will systematically replace $\Ga_g+\Ga_b$ by $\Ga_b$ and $r^{-1}\Ga_b+\Ga_g$ by $\Ga_g$.  
\end{remark}

\begin{remark}
In contrast to the scalar wave work \cite{MaSz24} which requires control of $2$ derivatives of $\Ga_g$ and $\Ga_b$, the present work for the Teukolsky equations requires control of $15$ derivatives of $\Ga_g$ and $\Ga_b$.
\end{remark}

\begin{lemma}\lab{lemma:relationsbetweennullframeandcoordinatesframe}
For $r\geq 13m$, we have
\bea\lab{eq:relationsbetweennullframeandcoordinatesframe:1}
\begin{split}
e_4 &= \left(1+O(mr^{-1})\right)\pr_r + O(m^2r^{-2})\pr_\tau+O(\ep r^{-2})\pr_{x^a},\\
e_3 &= \big(-1+O(\ep)\big)\pr_r + \big(2+O(mr^{-1})\big)\pr_\tau+\big(O(mr^{-2})+O(\ep r^{-1})\big)\pr_{x^a},\\
e_a &= O(\ep r^{-1})\pr_r + O(mr^{-1})\pr_\tau+r^{-1}\big((1+O(mr^{-1}))\mathring{\ga}_a\!\,^b+O(\ep)\big)\pr_{x^b},
\end{split}
\eea   
as well as
\bea\lab{eq:relationsbetweennullframeandcoordinatesframe:2}
\begin{split}
\pr_r &= \left(1+O(mr^{-1})\right)e_4 + O(m^2r^{-2})e_3+\big(O(m^3r^{-3})+O(\ep r^{-1})\big)e_a,\\
\pr_\tau &= \frac{1}{2}\big(1+O(mr^{-1}) +O(\ep)\big)e_4+\frac{1}{2}\big(1+O(mr^{-1})\big)e_3+\big(O(mr^{-1})+O(\ep)\big)e_a,\\
\pr_{x^a} &= O(m)e_4+O(m)e_3+r\big((1+O(mr^{-1}))(\mathring{\ga}_a\!\,^b)^{-1}+O(\ep)\big)e_b,
\end{split}
\eea
where $\mathring{\ga}_a\!\,^b$ is the invertible matrix denoting the decomposition of $(\pr_\th, \sin\th^{-1}\pr_{\tphi})$ on $(\pr_{x^1}, \pr_{x^2})$.
\end{lemma}

\begin{proof}
We decompose the null frame $(e_3, e_4, e_a)$ into the normalized coordinates vectorfields   
\beaa
e_\a=e_\a(x^\b)\pr_\b=[e_\a(x^\b)]_K\pr_\b+\widecheck{e_\a(x^\b)}\pr_\b, \quad \a=3,4,a, \quad a=1,2,
\eeaa
where $[e_\a(x^\b)]_K$ denotes the corresponding Kerr values given by \eqref{eq:actionofingoingprincipalnullframeonnormalizedcoordinates} \eqref{eq:actionofingoingprincipalnullframeonnormalizedcoordinates:bis} \eqref{eq:actionofingoingprincipalnullframeonnormalizedcoordinates:ter}. Together with Definition \ref{definition.Ga_gGa_b}, we infer for $r\geq 13m$, 
\bea\lab{eq:relationsbetweennullframeandcoordinatesframe2:moreprecise:00}
\bsplit
e_4 =& \left(1+O(mr^{-1})+\Ga_g\right)\pr_r + \big(O(m^2r^{-2})+\Ga_g\big)\pr_\tau+\Ga_g\pr_{x^a},\\
e_3 =& \left(-1+r\Ga_b\right)\pr_r + \big(2+O(mr^{-1})+r\Ga_g\big)\pr_\tau+\big(O(mr^{-2})+\Ga_b\big)\pr_{x^a},\\
e_a =& r\Ga_g\pr_r + \big(O(mr^{-1})+r\Ga_g\big)\pr_\tau+\big(r^{-1}(1+O(mr^{-1}))\mathring{\ga}_a\!\,^b+\Ga_b\big)\pr_{x^b},
\end{split}
\eea
where $\mathring{\ga}_a\!\,^b$ is the invertible matrix denoting the decomposition of $(\pr_\th, \sin\th^{-1}\pr_{\tphi})$ on $(\pr_{x^1}, \pr_{x^2})$, which implies
\bea\lab{eq:relationsbetweennullframeandcoordinatesframe2:moreprecise}
\bsplit
\pr_r =& \left(1+O(mr^{-1})+\Ga_g\right)e_4 + \left(O(m^2r^{-2})+\Ga_g\right)e_3+\big(O(m^3r^{-3})+r\Ga_g\big)e_a,\\
\pr_\tau =& \frac{1}{2}\big(1+O(mr^{-1})+r\Ga_b\big)e_4+\frac{1}{2}\big(1+O(mr^{-1})+r\Ga_g\big)e_3+\big(O(mr^{-1})+r\Ga_b\big)e_a,\\
\pr_{x^a} =& \big(O(m)+r^2\Ga_g\big)e_4+\big(O(m)+r^2\Ga_g\big)e_3+r\big((1+O(mr^{-1}))(\mathring{\ga}_a\!\,^b)^{-1}+r\Ga_b\big)e_b.
\end{split}
\eea
Together with \eqref{eq:decaypropertiesofGabGag}, this yields for $r\geq 13m$, 
\beaa
e_4 &=& \left(1+O(mr^{-1})\right)\pr_r + O(m^2r^{-2})\pr_\tau+O(\ep r^{-2})\pr_{x^a},\\
e_3 &=& \big(-1+O(\ep)\big)\pr_r + \big(2+O(mr^{-1})\big)\pr_\tau+\big(O(mr^{-2})+O(\ep r^{-1})\big)\pr_{x^a},\\
e_a &=& O(\ep r^{-1})\pr_r + O(mr^{-1})\pr_\tau+r^{-1}\big((1+O(mr^{-1}))\mathring{\ga}_a\!\,^b+O(\ep)\big)\pr_{x^b},
\eeaa   
as well as
\beaa
\pr_r &=& \left(1+O(mr^{-1})\right)e_4 + O(m^2r^{-2})e_3+\big(O(m^3r^{-3})+O(\ep r^{-1})\big)e_a,\\
\pr_\tau &=& \frac{1}{2}\big(1+O(mr^{-1}) +O(\ep)\big)e_4+\frac{1}{2}\big(1+O(mr^{-1})\big)e_3+\big(O(mr^{-1})+O(\ep)\big)e_a,\\
\pr_{x^a} &=& O(m)e_4+O(m)e_3+r\big((1+O(mr^{-1}))(\mathring{\ga}_a\!\,^b)^{-1}+O(\ep)\big)e_b,
\eeaa   
as stated.
\end{proof}

\begin{lemma}\lab{lemma:formofregularvectorfieldmathcalXs:Kerr:copy:horizontal}
For $s=\pm 2$, let the vectorfields $\mathcal{X}_s$ be given by
\bea
\lab{eq:formofregularhorizontalvectorfieldmathcalXs:Kerrperturbation}
\mathcal{X}_s = s\big(\pr_{\tphi}+a(\sin\th)^2\pr_\tau\big)+O(ar^{-1})\sin\th\pr_\th
\eea
where the coefficient $O(ar^{-1})$ is the regular function of $(r,\cos\th)$ appearing in \eqref{eq:formofregularhorizontalvectorfieldmathcalXs}. Then, we have
\bea
\lab{eq:formofregularvectorfieldmathcalXs:Kerr:copy:horizontal}
\mathcal{X}_s = \widetilde{\mathcal{X}}_s + r^2\Ga_g\dk, \qquad \widetilde{\mathcal{X}}_s\in\OO(\MM), \qquad |\widetilde{\mathcal{X}}_s|\les r,
\eea
i.e., $\mathcal{X}_s$ coincides with a horizontal vectorfield $\widetilde{\mathcal{X}}_s$ up to an error term. Also, we have
\bea
\lab{eq:formofprphi+aprtt:horizontal}
\pr_{\tphi} + a\pr_{\tt} =\mathcal{Y} +r^2\Ga_g \dk + O(|a|)\pr_{\tt}, \qquad {\mathcal{Y}}\in\OO(\MM), \qquad |{\mathcal{Y}}|\les r, 
\eea
i.e., $\pr_{\tphi} + a\pr_{\tt} $ coincides with a horizontal vectorfield $\mathcal{Y}$ up to error terms.
\end{lemma}

\begin{proof}
We consider the horizontal vectorfields $\widetilde{\mathcal{X}}_s$, $s=\pm 2$, defined by 
\beaa
\widetilde{\mathcal{X}}_s := s|q|^2\Re(\Jk)^be_b + O(ar^{-1})|q|^2\dual(\Re(\Jk))^be_b, \qquad \widetilde{\mathcal{X}}_s\in\OO(\MM), \qquad |\widetilde{\mathcal{X}}_s|\les r,
\eeaa
where the coefficient $O(ar^{-1})$ is the regular function of $(r,\cos\th)$ appearing in \eqref{eq:formofregularhorizontalvectorfieldmathcalXs}. Then, in view of Definition \ref{definition.Ga_gGa_b} and the identities \eqref{eq:usefulalgebraicidentitiesinvolvingscalarproductsReJkReJkpm:Kerrpert}, we have
\beaa
\widetilde{\mathcal{X}}_s &=& \widetilde{\mathcal{X}}_s(r)\pr_r+\widetilde{\mathcal{X}}_s(\tau)\pr_\tau+\widetilde{\mathcal{X}}_s(x^b_p)\pr_{x^b_p}\\
&=&  as|q|^2\Re(\Jk)\c\Re(\Jk)\pr_\tau+\Big(s|q|^2\Re(\Jk) + O(ar^{-1})|q|^2\dual(\Re(\Jk))\Big)\c\Re(\Jk_+)\pr_{x^1_p}\\
&&+\Big(s|q|^2\Re(\Jk) + O(ar^{-1})|q|^2\dual(\Re(\Jk))\Big)\c\Re(\Jk_-)\pr_{x^2_p}+r^2\Ga_g\dk\\
&=& s\Big(a(\sin\th)^2\pr_\tau -x^2_p\pr_{x^1_p}+x^1_p\pr_{x^2_p}\Big)+O(ar^{-1})\cos\th\big(x^1_p\pr_{x^1_p}+x^2_p\pr_{x^2_p}\big)+r^2\Ga_g\dk\\
&=& s\big(\pr_{\tphi}+a(\sin\th)^2\pr_\tau\big)+O(ar^{-1})\sin\th\pr_\th +r^2\Ga_g\dk = \mathcal{X}_s+r^2\Ga_g\dk
\eeaa
which proves the desired formula \eqref{eq:formofregularvectorfieldmathcalXs:Kerr:copy:horizontal}.

Next, we consider the horizontal vectorfields $\mathcal{Y}$ defined by
\beaa
\mathcal{Y}:=|q|^2\Re(\Jk)^be_b,\qquad {\mathcal{Y}}\in\OO(\MM), \qquad |{\mathcal{Y}}|\les r, 
\eeaa
and obtain, in view of Definition \ref{definition.Ga_gGa_b} and the identities \eqref{eq:usefulalgebraicidentitiesinvolvingscalarproductsReJkReJkpm:Kerrpert}, 
\beaa
\mathcal{Y} &=& \mathcal{Y}(r)\pr_r+\mathcal{Y}(\tau)\pr_\tau+\mathcal{Y}(x^b_p)\pr_{x^b_p}\\
&=& a|q|^2\Re(\Jk)\c\Re(\Jk)\pr_\tau + |q|^2\Re(\Jk)\c\Re(\Jk_+)\pr_{x^1} + |q|^2\Re(\Jk)\c\Re(\Jk_-)\pr_{x^2} +r^2\Ga_g\dk\\
&=& a(\sin\th)^2\pr_\tau -x^2\pr_{x^1} + x^1\pr_{x^2} +r^2\Ga_g\dk = \pr_{\tphi}+a(\sin\th)^2\pr_\tau +r^2\Ga_g\dk.
\eeaa
We deduce
\beaa
\pr_{\tphi} + a\pr_{\tau} = \big(\pr_{\tphi}+a(\sin\th)^2\pr_\tau\big)+ a (\cos\th)^2 \pr_{\tt} = \mathcal{Y} +r^2\Ga_g \dk  + a (\cos\th)^2 \pr_{\tt}
\eeaa
which proves the desired formula \eqref{eq:formofprphi+aprtt:horizontal}. This concludes the proof of Lemma \ref{lemma:formofregularvectorfieldmathcalXs:Kerr:copy:horizontal}.
\end{proof}


\subsubsection{Commutation formulas revisited}


We start by revisiting the commutators of Lemma \ref{lemma:comm}.
 \begin{corollary}
   \lab{cor:corollaryofLemmacomm}
   The following commutation formulas hold true:
   \begin{enumerate}
\item Given   $f \in \sk_0$, we have
       \bea\label{eq:comm-nab3-nab4-naba-f-general:cor}
       \begin{split}
        \,[\nab_3, \nab_a] f &=\frac{1}{r}\nab_a f +\big(O(mr^{-2})+\etac\big)\nab_3 f +\big(O(mr^{-3})+\Ga_g\big)\dk f,\\
         \,[\nab_4, \nab_a] f &=-\frac{1}{r}\nab_a f +\xi_a \nab_3 f +\big(O(mr^{-3})+r^{-1}\Ga_g\big)\dk f, \\
         \, [\nab_4, \nab_3] f&= O(mr^{-2})\nab_3 f +\big(O(mr^{-3})+\Ga_g\big)\dk f,
         \end{split}
       \eea
       and in particular 
          \bea\label{eq:comm-nab3-nab4-naba-f-general:cor:bis}
       \begin{split}
       \,[\nab_3, r\nab_a] f &=\big(O(mr^{-1})+r\etac\big)\nab_3 f +\big(O(mr^{-2})+r\Ga_g\big)\dk f,\\
        \,[\nab_4, r\nab_a] f &=r\xi_a \nab_3 f +\big(O(mr^{-2})+\Ga_g\big)\dk f.
       \end{split}
       \eea

  \item   Given  $u\in \sk_1$ or $u\in\sk_2$, we have
    \bea\label{commutator-3-a-u-b:cor}
         \bsplit            
\,  [\nab_3,\nab_a] u    &=\frac{1}{r}\nab_au +\big(O(mr^{-2})+\etac\big)\nab_3u +\big(O(mr^{-3})+\Ga_g\big)\dk^{\leq 1}u,\\
\,  [\nab_4,\nab_a] u    &=-\frac{1}{r}\nab_au +\xi_a \nab_3u +\big(O(mr^{-3})+r^{-1}\Ga_g\big)\dk^{\leq 1}u,\\
 \, [\nab_4, \nab_3] u &=O(mr^{-2})\nab_3u +\big(O(mr^{-3})+\Ga_g\big)\dk^{\leq 1}u,
\end{split}
\eea
and in particular
   \bea\label{commutator-3-a-u-b:cor:bis}
         \bsplit            
\,  [\nab_3,r\nab_a] u    &=\big(O(mr^{-1})+r\etac\big)\nab_3u +\big(O(mr^{-2})+r\Ga_g\big)\dk^{\leq 1}u,\\
\,  [\nab_4,r\nab_a] u    &=r\xi_a \nab_3u +\big(O(mr^{-2})+\Ga_g\big)\dk^{\leq 1}u.
\end{split}
\eea
       \end{enumerate}
Finally, we also have, for $u\in\sk_k$, $k=1,2$, 
\bea\label{commutator-nab-a-nab-b:cor}
[\nab_a, \nab_b]u = \big(O(m r^{-2})+\Ga_g\big)\nab_3u+\big(O(m r^{-2})+\Ga_g\big)\nab_4u+\big(O(r^{-2})+r^{-1}\Ga_g\big)u.
\eea 
 \end{corollary}

\begin{proof}
The proof of \eqref{eq:comm-nab3-nab4-naba-f-general:cor}--\eqref{commutator-3-a-u-b:cor:bis} follows immediately from Lemma \ref{lemma:comm} and the definition of $(\Ga_b, \Ga_g)$, while \eqref{commutator-nab-a-nab-b:cor} is a non-sharp consequence of Proposition 2.1.43 in \cite{GKS22}.
\end{proof}

\begin{corollary}\lab{cor:commutatorweightedderivativesrnabandnab4rwithnabprtau}
For $u\in\sk_k$, $k=1,2$, we have
\beaa
\,[\nab_{\pr_\tau}, r\nab_b]u &=& \big(O(mr^{-1})+r\dk^{\leq 1}\Ga_b\big)\nab_3u+\big(O(mr^{-2})+r\dk^{\leq 1}\Ga_g\big)\dk^{\leq 1}u,\\
\,[\nab_{\pr_\tau}, \nab_4r]u &=& \big(O(mr^{-1})+r\dk^{\leq 1}\Ga_g\big)\nab_3u+\big(O(mr^{-2})+r\dk^{\leq 1}\Ga_g\big)\dk^{\leq 1}u.
\eeaa
\end{corollary}

\begin{proof}
Recalling from \eqref{eq:relationsbetweennullframeandcoordinatesframe2:moreprecise} that 
\beaa
\pr_\tau = \frac{1}{2}\big(1+O(mr^{-1})+r\Ga_b\big)e_4+\frac{1}{2}\big(1+O(mr^{-1})+r\Ga_g\big)e_3+\big(O(mr^{-1})+r\Ga_b\big)e_a,
\eeaa
we have, for $u\in\sk_k$, $k=1,2$, 
\beaa
\bsplit
[\nab_{\pr_\tau}, r\nab_b]u =& \frac{1}{2}\big(1+O(mr^{-1})+r\Ga_b\big)[\nab_4, r\nab_b]u+\frac{1}{2}\big(1+O(mr^{-1})+r\Ga_g\big)[\nab_3, r\nab_b]u\\
&+\big(O(mr^{-1})+r\Ga_b\big)[\nab_a, r\nab_b]u+\big(O(mr^{-1})+r\dk\Ga_g\big)\nab_3u+\big(O(mr^{-2})+\dk\Ga_b\big)\dk u.
\end{split}
\eeaa
Together with \eqref{commutator-3-a-u-b:cor:bis} and \eqref{commutator-nab-a-nab-b:cor}, we infer
\beaa
[\nab_{\pr_\tau}, r\nab_b]u &=& \big(O(mr^{-1})+r\dk^{\leq 1}\Ga_b\big)\nab_3u+\big(O(mr^{-2})+r\dk^{\leq 1}\Ga_g\big)\dk^{\leq 1}u
\eeaa
as stated.

Similarly, we have, for $u\in\sk_k$, $k=1,2$,
\beaa
\bsplit
[\nab_{\pr_\tau}, \nab_4r]u =& \frac{1}{2}\big(1+O(mr^{-1})+r\Ga_b\big)[\nab_4, \nab_4r]u+\frac{1}{2}\big(1+O(mr^{-1})+r\Ga_g\big)[\nab_3, \nab_4r]u\\
&+\big(O(mr^{-1})+r\Ga_b\big)[\nab_a, \nab_4r]u+\big(O(mr^{-1})+r\dk\Ga_g\big)\nab_3u+\big(O(mr^{-2})+\dk\Ga_b\big)\dk u\\
=&  \frac{1}{2}(e_4(r)+e_3(r))\nab_4 u +\frac{1}{2}\big(1+O(mr^{-1})+r\Ga_g\big)[\nab_3, \nab_4](ru)\\
&+\big(O(mr^{-1})+r\Ga_b\big)[\nab_a, \nab_4](ru)+\big(O(mr^{-1})+r\dk\Ga_g\big)\nab_3u\\
&+\big(O(mr^{-2})+\dk^{\leq 1}\Ga_b\big)\dk u
\end{split}
\eeaa
and hence
\beaa
\bsplit
[\nab_{\pr_\tau}, \nab_4r]u =& \frac{1}{2}\big(1+O(mr^{-1})+r\Ga_g\big)[\nab_3, \nab_4](ru)+\big(O(mr^{-1})+r\Ga_b\big)[\nab_a, \nab_4](ru)\\
&+\big(O(mr^{-1})+r\dk\Ga_g\big)\nab_3u+\big(O(mr^{-2})+\dk^{\leq 1}\Ga_b\big)\dk u
\end{split}
\eeaa
which together with \eqref{commutator-3-a-u-b:cor} implies 
\beaa
[\nab_{\pr_\tau}, \nab_4r]u =\big(O(mr^{-1})+r\dk^{\leq 1}\Ga_g\big)\nab_3u+\big(O(mr^{-2})+r\dk^{\leq 1}\Ga_g\big)\dk^{\leq 1}u
\eeaa
as stated.
\end{proof}

 
 \subsection{Assumptions on the perturbed metric and consequences}
\label{subsect:assumps:perturbedmetric}
 

In this section, we recall from \cite[Section 2.4]{MaSz24} the assumptions for the metric perturbations relative to a subextremal Kerr and further estimates for the metric under these metric perturbation assumptions. All the statements and estimates in this section are from \cite[Section 2.4]{MaSz24}.


\subsubsection{Assumptions on the inverse metric perturbation}
 \label{subsubsect:assumps:perturbedmetric}


We introduce our assumptions on the perturbed inverse metric.

\begin{assumption}[Inverse metric assumptions]
\label{intro:assump:metric}
Let a subextremal Kerr metric $\gam$ be given and define, in the normalized coordinates $(\tt,r,x^1_0, x^2_0)$ and {$(\tt,r,x^1_p, x^2_p)$}, the inverse metric difference 
\bea
\gcheck^{\a\b}:=\g^{\a\b}-\gam^{\a\b}.
\eea
Then, with $(\Ga_b, \Ga_g)$ verifying \eqref{eq:decaypropertiesofGabGag}, we assume that $\gcheck^{\a\b}$ satisfies the following\footnote{These estimates hold in fact on each coordinate patch, i.e., in the coordinates $(\tt,r,x^1_0, x^2_0)$ for $\th\in[\frac{\pi}{4}, \frac{3\pi}{4}]$, and in the coordinates {$(\tt,r,x^1_p, x^2_p)$} for $\th\in [0,\pi]\setminus(\frac{\pi}{3}, \frac{2\pi}{3})$.}:
\bea\lab{eq:controloflinearizedinversemetriccoefficients}
\widecheck{\g}^{rr}=r\Ga_b, \quad \widecheck{\g}^{r\tau}=r\Ga_g,\quad \widecheck{\g}^{\tau\tau}=\Ga_g, \quad \widecheck{\g}^{ra}=\Ga_b, \quad \widecheck{\g}^{\tau a}=\Ga_g,\quad \widecheck{\g}^{ab}=r^{-1}\Ga_g.
\eea 
\end{assumption}

\begin{remark}
The decay assumptions \eqref{eq:controloflinearizedinversemetriccoefficients} on $\gcheck^{\a\b}$ are consistent with the decay estimates derived in the proof of the nonlinear stability of Kerr for small angular momentum in \cite{KS:Kerr}.
\end{remark}

The following immediate non-sharp consequence of \eqref{eq:assymptiticpropmetricKerrintaurxacoord:1}, \eqref{eq:controloflinearizedinversemetriccoefficients} and  \eqref{eq:decaypropertiesofGabGag} will be useful
\bea\lab{eq:consequenceasymptoticKerrandassumptionsinverselinearizedmetric}
\begin{split}
\g^{rr}&=O(1), \qquad \g^{r\tau}=O(1), \qquad \g^{ra}=O(r^{-1}),\\
\g^{\tau\tau}&=O(m^2r^{-2}), \qquad \g^{\tau a}=O(mr^{-2}), \qquad \g^{ab}=O(r^{-2}).
\end{split}
\eea


\subsubsection{Control of the metric perturbation}


The following lemma provides the control of the perturbed metric coefficients which follows from the assumption \eqref{eq:controloflinearizedinversemetriccoefficients} on the perturbed inverse metric coefficients. 
\begin{lemma}\lab{lemma:controlofmetriccoefficients:bis}
Assume that $\widecheck{\g}^{\a\b}$ verifies \eqref{eq:controloflinearizedinversemetriccoefficients}. Then, $\widecheck{\g}_{\a\b}:=\g_{\a\b}-(\gam)_{\a\b}$ verifies  
\bea\lab{eq:controloflinearizedmetriccoefficients}
\begin{split}
\widecheck{\g}_{rr} &=\Ga_g, \qquad\quad \widecheck{\g}_{r\tau}=r\Ga_g, \qquad\quad \widecheck{\g}_{\tau\tau}=r\Ga_b, \\
\widecheck{\g}_{\tau a} &=r^2\Ga_b, \qquad\, \widecheck{\g}_{ra}=r^2\Ga_g, \qquad\,\,\,\, \widecheck{\g}_{ab}=r^3\Ga_g.
\end{split}
\eea
Also, we have
\bea\lab{eq:controloflinearizedmetriccoefficients:det}
\widecheck{\det(\g)}=\det(\g_{a,m})r^2\Ga_g, \qquad \widecheck{\det(\g)}:=\det(\g)-\det(\g_{a,m}).
\eea
\end{lemma}

The following immediate non-sharp consequence of \eqref{eq:assymptiticpropmetricKerrintaurxacoord:2}, \eqref{eq:controloflinearizedmetriccoefficients} and  \eqref{eq:decaypropertiesofGabGag} will be useful
\bea\lab{eq:consequenceasymptoticKerrandassumptionsinverselinearizedmetric:bis}
\bsplit
\g_{rr}&=O(m^2r^{-2}), \qquad \g_{r \tt}=O(1),\qquad \g_{ra}=O(m),\\
\g_{\tt \tt}&= O(1),\qquad\qquad\,\, \g_{\tau a}=O(r),\qquad \, \g_{ab}=O(r^2).
 \end{split}
\eea


\subsubsection{Control of the induced metric on $\Si(\tau)$ and $\AA$}


The following lemma provides the control of the induced metric on $\Si(\tau)$. 
\begin{lemma}\lab{lemma:controllinearizedmetric:inducedmetricSitau}
Let $g$ denote the metric induced by $\g$ on the level sets of $\tau$. Assume that $\widecheck{\g}^{\a\b}$ verifies \eqref{eq:controloflinearizedinversemetriccoefficients}. Then, $\widecheck{g}_{ij}:=g_{ij}-(g_{a,m})_{ij}$ and $\widecheck{g}^{ij}:=g^{ij}-g_{a,m}^{ij}$ verify
\beaa
\bsplit
\widecheck{g}_{rr} &=\Ga_g,  \qquad \widecheck{g}_{ra}=r^2\Ga_g, \qquad \widecheck{g}_{ab}=r^3\Ga_g,\\
\widecheck{g}^{rr} &=r^4\Ga_g,  \qquad \widecheck{g}^{ra}=r^2\Ga_g, \qquad \widecheck{g}^{ab}=\Ga_g.
\end{split}
\eeaa
Also, we have $\widecheck{\det(g)}=r^4\Ga_g$, with $\widecheck{\det(g)}:=\det(g)-\det(g_{a,m})$. 
\end{lemma}

The following lemma provides the control of the determinant of the induced metric on $\AA$. 
\begin{lemma}\lab{lemma:controllinearizedmetric:inducedmetricAA}
Let $g_\AA$ denote the metric induced by $\g$ on the spacelike hypersurface $\AA$.  Assume that $\widecheck{\g}^{\a\b}$ verifies \eqref{eq:controloflinearizedinversemetriccoefficients}. Then, $\widecheck{\det(g_\AA)}=O(\ep\tau^{-1-\de})$, with $\widecheck{\det(g_\AA)}:=\det(g)-\det(g_\AA)$. 
\end{lemma}


\subsubsection{Further consequences of the metric assumptions}


In this section, we draw further consequences of the assumption \eqref{eq:controloflinearizedinversemetriccoefficients} on the perturbed inverse metric coefficients. 

\begin{lemma}\lab{lemma:computationofthederiveativeofsrqtg}
Let the 1-form $N_{det}$ be defined by  
\beaa
(N_{det})_\mu:=\frac{1}{\sqrt{|\g|}}\pr_\mu\sqrt{|\g|} - \frac{1}{\sqrt{|\g_{a,m}|}}\pr_\mu\sqrt{|\g_{a,m}|}. 
\eeaa
Then, we have 
\beaa
(N_{det})_r=\dk^{\leq 1}\Ga_g, \qquad (N_{det})_\tau=r\dk^{\leq 1}\Ga_g, \qquad (N_{det})_{x^a}=r\dk^{\leq 1}\Ga_g,
\eeaa 
and 
\beaa
(N_{det})^r=r\dk^{\leq 1}\Ga_g, \qquad (N_{det})^\tau=\dk^{\leq 1}\Ga_g, \qquad (N_{det})^{x^a}=r^{-1}\dk^{\leq 1}\Ga_g.
\eeaa
\end{lemma}

We have the following corollary of Lemma \ref{lemma:computationofthederiveativeofsrqtg}. 
\begin{corollary}\lab{cor:controloflinearizeddivergencecoordvectorfields}
We have
\beaa
\widecheck{\textbf{\textrm{Div}}(\pr_r)}=\dk^{\leq 1}\Ga_g, \qquad \widecheck{\textbf{\textrm{Div}}(\pr_\tau)}=r\dk^{\leq 1}\Ga_g, \qquad \widecheck{\textbf{\textrm{Div}}(\pr_{x^a})}=r\dk^{\leq 1}\Ga_g, \quad a=1,2.
\eeaa
\end{corollary}

Next, we provide the control of deformation tensors involved in energy-Morawetz estimates{, where the deformation tensor of a vectorfield $X$ is given by
\bea
\label{def:deformationtensor:lastsect}
{}^{(X)}\pi_{\a\b} :=  \D_{\a}X_{\b} + \D_{\b}X_{\a}=\LL_X\g_{\a\b}.
\eea}

\begin{lemma}\lab{lemma:controlofdeformationtensorsforenergyMorawetz}
The deformation {tensors of $\pr_\tau$ and $\pr_{\tphi}$ satisfy}
\beaa
\begin{split}
{}^{(\pr_\tau)}\pi_{rr}{, {}^{(\pr_{\tphi})}\pi_{rr}} &=\dk^{\leq 1}\Ga_g, \qquad {}^{(\pr_\tau)}\pi_{r\tau}{, {}^{(\pr_{\tphi})}\pi_{r\tau}}=r\dk^{\leq 1}\Ga_g, \qquad {}^{(\pr_\tau)}\pi_{\tau\tau}{, {}^{(\pr_{\tphi})}\pi_{\tau\tau}}=r\dk^{\leq 1}\Ga_b, \\
{}^{(\pr_\tau)}\pi_{\tau a}{, {}^{(\pr_{\tphi})}\pi_{\tau a}} &=r^2\dk^{\leq 1}\Ga_b, \qquad {}^{(\pr_\tau)}\pi_{ra}{, {}^{(\pr_{\tphi})}\pi_{ra}}=r^2\dk^{\leq 1}\Ga_g, \qquad {}^{(\pr_\tau)}\pi_{ab}{, {}^{(\pr_{\tphi})}\pi_{ab}}=r^3\dk^{\leq 1}\Ga_g,
\end{split}
\eeaa
and 
\beaa
\begin{split}
{}^{(\pr_\tau)}\pi^{rr}{, {}^{(\pr_{\tphi})}\pi^{rr}} &=r\dk^{\leq 1}\Ga_b, \qquad {}^{(\pr_\tau)}\pi^{r\tau}{, {}^{(\pr_{\tphi})}\pi^{r\tau}}=r\dk^{\leq 1}\Ga_g, \qquad {}^{(\pr_\tau)}\pi^{\tau\tau}{, {}^{(\pr_{\tphi})}\pi^{\tau\tau}}=\dk^{\leq 1}\Ga_g, \\
{}^{(\pr_\tau)}\pi^{\tau a}{, {}^{(\pr_{\tphi})}\pi^{\tau a}} &=\dk^{\leq 1}\Ga_g, \qquad {}^{(\pr_\tau)}\pi^{ra}{, {}^{(\pr_{\tphi})}\pi^{ra}}=\dk^{\leq 1}\Ga_b, \qquad {}^{(\pr_\tau)}\pi^{ab}{, {}^{(\pr_{\tphi})}\pi^{ab}}=r^{-1}\dk^{\leq 1}\Ga_g.
\end{split}
\eeaa

Also, the perturbed deformation tensor of $\pr_r$ satisfies 
\beaa
\begin{split}
\widecheck{{}^{(\pr_r)}\pi}_{rr} &=r^{-1}\dk^{\leq 1}\Ga_g, \qquad \widecheck{{}^{(\pr_r)}\pi}_{r\tau}=\dk^{\leq 1}\Ga_g, \qquad \widecheck{{}^{(\pr_r)}\pi}_{\tau\tau}=\dk^{\leq 1}\Ga_b, \\
\widecheck{{}^{(\pr_r)}\pi}_{\tau a} &=r\dk^{\leq 1}\Ga_b, \qquad \widecheck{{}^{(\pr_r)}\pi}_{ra}=r\dk^{\leq 1}\Ga_g, \qquad \widecheck{{}^{(\pr_r)}\pi}_{ab}=r^2\dk^{\leq 1}\Ga_g,
\end{split}
\eeaa
and
\beaa
\begin{split}
\widecheck{{}^{(\pr_r)}\pi}^{rr} &=\dk^{\leq 1}\Ga_b, \qquad \widecheck{{}^{(\pr_r)}\pi}^{r\tau}=\dk^{\leq 1}\Ga_g, \qquad \widecheck{{}^{(\pr_r)}\pi}^{\tau\tau}=r^{-1}\dk^{\leq 1}\Ga_g, \\
\widecheck{{}^{(\pr_r)}\pi}^{\tau a} &=r^{-1}\dk^{\leq 1}\Ga_g, \qquad \widecheck{{}^{(\pr_r)}\pi}^{ra}=r^{-1}\dk^{\leq 1}\Ga_b, \qquad \widecheck{{}^{(\pr_r)}\pi}^{ab}=r^{-2}\dk^{\leq 1}\Ga_g.
\end{split}
\eeaa
\end{lemma}


\subsection{Assumptions on the regular triplet $\Om_i$, $i=1,2,3$}
\lab{sec:regulartripletinperturbationsofKerr}


Recall from the discussion at the beginning of Section  \ref{subsect:introduceOm_i} that it will be convenient to scalarize the Teukolsky equations, i.e., to transform a system of tensorial wave equations on $\sk_2(\mathbb{C})$ to a coupled system of scalar wave equations. To this end, we assume on $\MM$ the existence of a regular triplet $\Om_i$, $i=1,2,3$, in the sense of Definition \ref{def:definitionofregulartripletOmii=123}, i.e., $\Om_i\in\sk_1$ and  satisfy the following identities 
\beaa
x^i\Om_i=0, \qquad (\Om^i)_a(\Om_i)_b=\de_{ab}, \qquad \Om_i\c\Om_j=\de_{ij}-x^ix^j,\qquad \Om_i\c\dual\Om_j=\in_{ijk}x^k,
\eeaa
see Remark \ref{rmk:generalcontructionofregulartripletsingivenspacetime} on how to generate regular triplets on a given spacetime $\MM$. 

Next, recall from Definition \ref{def:Mialphaj:Kerr} that we associate to the regular triplet above the following 1-forms on $\MM$ 
\beaa
M_{i\a}^j:=(\Ddot_\a\Om_i)\c\Om^j, \quad \forall \a,i,j.
\eeaa
Our assumptions on the regular triplet $\Om_i$, $i=1,2,3$ are the following
\bsub\lab{eq:assumptionsonregulartripletinperturbationsofKerr}
\bea
&&\lab{eq:assumptionsonregulartripletinperturbationsofKerr:0}
\widecheck{M_{i4}^j}=\Ga_g, \qquad \widecheck{M_{i3}^j}=\Ga_b,  \qquad \widecheck{M_{ia}^j}=\Ga_b, \quad \forall\, i,j,a,\\
&&\lab{eq:assumptionforLiebprtphiOmiinKerrperturbation}
\Lieb_{\pr_{\tphi}}\Om_i +\in_{ij3}\Om^j=r\Ga_b, \quad\textrm{for}\quad i=1,2,3.
\eea
\esub 

\begin{remark}
The decay assumptions \eqref{eq:assumptionsonregulartripletinperturbationsofKerr} are consistent with the decay estimates derived in the proof of the nonlinear stability of Kerr for small angular momentum in \cite{KS:Kerr}.
\end{remark}

In view of Lemma \ref{lemma:computationoftheMialphajinKerr} for $(M_{i\a}^j)_K$ and the assumption \eqref{eq:assumptionsonregulartripletinperturbationsofKerr:0}
 for $\widecheck{M_{i\a}^j}$, we have the following estimates for $M_{i\a}^j$.
\begin{lemma}
\lab{lem:estimatesforMialphaj:Kerrpert}
Under the assumption \eqref{eq:assumptionsonregulartripletinperturbationsofKerr:0}  for $\widecheck{M_{i\a}^j}$, we have
 \begin{equation}
 \lab{estimates:Mialphaj:Kerrperturbations}
 \begin{split}
& M_{i4}^j={O(r^{-2})}, \quad
 M_{ia}^j = {O(r^{-1})}, \quad  M_{i3}^j=O(r^{-2})+\Ga_b,\\
 &M_{i\a}^j(\pr_{\tau})^{\a}=O(r^{-3}) + \Ga_b, \quad M_{i\a}^j( \pr_{r}^{\text{BL}})^{\a}=\Ga_b
 \end{split}
 \end{equation}
 and 
 \bea
 \g^{r\a}M_{i\a}^{j}=\Ga_b, \quad r\in [r_+(1+2\dbl), 12m].
 \eea
\end{lemma}

In view of \eqref{eq:assumptionforLiebprtphiOmiinKerrperturbation} and the definition of $\widehat{\pr}_{\tphi}$, we have the following lemma.

\begin{lemma}\lab{lemma:differencebetweenwidehatprtphiandLiebprtrphiisanerrorterm}
Let  $\psi_{ij}=\pmb\psi(\Om_i, \Om_j)$ with $\pmb\psi\in\sk_2$. Then, under the assumption \eqref{eq:assumptionforLiebprtphiOmiinKerrperturbation}, we have 
\beaa
\widehat{\pr}_{\tphi}(\psi)_{ij}=\Lieb_{\pr_{\tphi}}\pmb\psi(\Om_i, \Om_j)+r\Ga_b\c\pmb\psi, \quad \forall\, i,j,  
\eeaa
where $\widehat{\pr}_{\tphi}$ has been introduced in Definition \ref{def:widehatprtphi}.
\end{lemma}

\begin{proof}
We have
\beaa
\pr_{\tphi}(\psi_{ij}) = \pr_{\tphi}(\pmb\psi(\Om_i, \Om_j))= \Lieb_{\pr_{\tphi}}\pmb\psi(\Om_i, \Om_j)+\pmb\psi(\Lieb_{\pr_{\tphi}}\Om_i, \Om_j)+\pmb\psi(\Om_i, \Lieb_{\pr_{\tphi}}\Om_j)
\eeaa
which together with  \eqref{eq:assumptionforLiebprtphiOmiinKerrperturbation} implies 
\beaa
\pr_{\tphi}(\psi_{ij}) =  \Lieb_{\pr_{\tphi}}\pmb\psi(\Om_i, \Om_j) -\in_{ik3}\psi_{kj} -\in_{jk3}\psi_{ik} +r\Ga_b\c\pmb\psi,
\eeaa
and hence, in view of Definition \ref{def:widehatprtphi}, 
\beaa
\widehat{\pr}_{\tphi}(\psi)_{ij}=\Lieb_{\pr_{\tphi}}\pmb\psi(\Om_i, \Om_j) +r\Ga_b\c\pmb\psi
\eeaa
as stated.
\end{proof}


\subsection{Teukolsky wave/transport systems in perturbations of Kerr}
\lab{sec:TeukolskyWavesysteminperturbationsofKerr}


We are now ready to provide the form of the Teukolsky wave/transport systems in perturbations of Kerr, i.e., in $(\MM, \g)$ where $\g$ is a perturbation of $\gam$ with $|a|<m$. We will start with the original tensorial form as derived in \cite{GKS22}, and we will then provide the corresponding scalarized form using the regular triplet of Section  \ref{sec:regulartripletinperturbationsofKerr}.


\subsubsection{Tensorial Teukolsky wave/transport systems in perturbations of Kerr}
\lab{sec:TeukolskyWavesysteminperturbationsofKerr:tensorialform}


We consider $\pmb\phi_s^{(p)}\in\sk_2(\mathbb{C})$, $s=\pm 2$, $p=0,1,2$, with $\pmb\phi_s^{(0)}$ given by 
\bea\lab{eq:defintionphipm2p=0:Kerrperturbation}
\pmb\phi_{+2}^{(0)}:=\frac{\ov{q}}{q}A, \qquad \pmb\phi_{-2}^{(0)}:=\frac{q}{\ov{q}}\left(\frac{\De}{|q|^2}\right)^2\Ab.
\eea
Then, the tensorial Teukolsky wave equations in perturbations of Kerr are given by 
\bsub
\lab{eq:TensorialTeuSysandlinearterms:rescaleRHScontaine2:general:Kerrperturbation}
\bea
\lab{eq:TensorialTeuSys:rescaleRHScontaine2:general:Kerrperturbation}
\bigg(\squared_2 -\frac{4ia\cos\th}{|q|^2}\nab_{\pr_{\tt}}- \frac{4-2\de_{p0}}{\qs}\bigg){\phis{p}} = \L_{s}^{(p)}[\pmb\phi_{s}]+\N_{W,s}^{(p)}, \quad s=\pm2, \quad p=0,1,2,
\eea
where the linear coupling terms $\L_{s}^{(p)}[{\pmb\phi_s}]$ have the following schematic forms
\bea
\lab{eq:tensor:Lsn:onlye_2present:general:Kerrperturbation}
\bsplit
{\L_{s}^{(0)}[\pmb\phi_{s}]}={}& (2sr^{-3} +O(mr^{-4}))\phis{1}+ O(mr^{-3}) \nab_{\Xcal_s}^{\leq 1}\phis{0},\\
{\L_{s}^{(1)}[\pmb\phi_{s}]}={}& (sr^{-3} +O(mr^{-4}))\phis{2}+ O(mr^{-3}) \nab_{\Xcal_s}^{\leq 1}  \phis{1}+O(mr^{-2})\nab_{\pr_{\tphi}+a\pr_{\tt}}^{\leq 1}\phis{0},\\
{\L_{s}^{(2)}[\pmb\phi_{s}]}={}&O(mr^{-3})\phis{2}+O(mr^{-2})\nab_{\pr_{\tphi}+a\pr_{\tt}}^{\leq 1}\phis{1}+O(m^2 r^{-2})\phis{0},
\end{split}
\eea
\esub
with $\Xcal_s$, $s=\pm 2$, being the regular vectorfields introduced in \eqref{eq:formofregularhorizontalvectorfieldmathcalXs:Kerrperturbation}, with all the coefficients in  \eqref{eq:tensor:Lsn:onlye_2present:general:Kerrperturbation} being independent of coordinates $\tau$ and 
$\tphi$, and with the coefficients in front of the terms $\phis{2}$ and $\nab_{\pr_{\tphi}+a\pr_{\tt}}\phis{1}$ on the RHS of equation of ${\L_{s}^{(2)}[\pmb\phi_{s}]}$ in \eqref{eq:tensor:Lsn:onlye_2present:general:Kerrperturbation} being real functions. 

Moreover, the tensorial Teukolsky transport equations in perturbations of Kerr are given by
\bsub\lab{def:TensorialTeuScalars:wavesystem:Kerrperturbation}
 \bea
\lab{def:TensorialTeuScalars:wavesystem:Kerrperturbation:+2}
\nab_3 \left(\frac{r\bar{q}}{q}\left(\frac{r^2}{|q|^2}\right)^{p-2}\pmb\phi_{+2}^{(p)}\right)=\frac{\bar{q}}{rq}\left(\frac{r^2}{|q|^2}\right)^{p-1}\pmb\phi_{+2}^{(p+1)}+\N_{T,+2}^{(p)}, \quad p=0,1,
\eea
and
\bea
\lab{def:TensorialTeuScalars:wavesystem:Kerrperturbation:-2}
\nab_4\left(\frac{rq}{\bar{q}}\left(\frac{{r^2}}{|q|^2}\right)^{p-2}\pmb\phi_{-2}^{(p)}\right)=\frac{q}{r\bar{q}}\left(\frac{r^2}{|q|^2}\right)^{p-1}\frac{\De}{\qs}\pmb\phi_{-2}^{(p+1)}+\N_{T,-2}^{(p)}, \,\,\,\, p=0,1.
\eea
\esub
The equations \eqref{eq:TensorialTeuSysandlinearterms:rescaleRHScontaine2:general:Kerrperturbation}  and \eqref{def:TensorialTeuScalars:wavesystem:Kerrperturbation} then correspond to the tensorial Teukolsky wave/transport systems in perturbations of Kerr.

\begin{remark}
For the proof of our main result, the only relevant property of the vectorfields $\Xcal_s$ and $\pr_{\tphi}+a\pr_{\tt}$ appearing in  \eqref{eq:tensor:Lsn:onlye_2present:general:Kerrperturbation} is that they satisfy the structure exhibited in \eqref{eq:formofregularvectorfieldmathcalXs:Kerr:copy:horizontal} \eqref{eq:formofprphi+aprtt:horizontal}.
\end{remark}

\begin{remark}\lab{rmk:introductionofpmbphip=0s=minus2nodegtocontrolAbnearr=rplus}
The terms $\N_{W,s}^{(p)}$ and $\N_{T,s}^{(p)}$, respectively on the RHS of \eqref{eq:TensorialTeuSysandlinearterms:rescaleRHScontaine2:general:Kerrperturbation} and \eqref{def:TensorialTeuScalars:wavesystem:Kerrperturbation}, correspond to nonlinear terms generated when deriving the tensorial Teukolsky wave/transport systems in perturbations of Kerr. In particular, we have $\N_{W,s}^{(p)}=\N_{T,s}^{(p)}=0$ in Kerr in which case \eqref{eq:TensorialTeuSysandlinearterms:rescaleRHScontaine2:general:Kerrperturbation}  \eqref{def:TensorialTeuScalars:wavesystem:Kerrperturbation} coincide with the tensorial Teukolsky wave/transport systems in Kerr, see \eqref{eq:TensorialTeuSysandlinearterms:rescaleRHScontaine2:general:Kerr} and  \eqref{def:TensorialTeuScalars:wavesystem:Kerr}. 
\end{remark}

\begin{remark}\lab{rmk:pmbphiminus2p=0isdegenerateatr=rplus}
Notice in view of \eqref{eq:defintionphipm2p=0:Kerrperturbation} that $\pmb\phi_{-2}^{(0)}$ degenerates at $r=r_+$ and does thus not allow to recover estimates for $\Ab$ near $r=r_+$. To remedy this problem, we will rely on the following form of Teukolsky equation for $\Ab$ in the redshift region $r\leq r_+(1+2\dred)$:
\bea\lab{eq:waveequationpmbphip=0sminus2nodeginredshiftregion}
\nn\squared_2\nab_4^p\Ab &=& (2-p)\pr_r\left(\frac{\De}{|q|^2}\right)\nab_3\nab_4^p\Ab
+ O(1)\big(\nab_{4}\nab_4^{\leq p}\Ab,\nab\nab_4^{\leq p}\Ab, \nab_4^{\leq p}\Ab\big)\\
&&+\N_{\nab_4^p\Ab}, \qquad p=0,1,2,
\eea
where the case $p=0$ is a non-sharp consequence\footnote{Note that, for $r\leq r_+(1+2\dred)$,  we have $\N_{\Ab}=\dk^{\leq 1}(\Ga_g\c\Ga_b)$ in view of Lemma 5.3.3 in \cite{GKS22}.}  of Lemma 5.3.3 in \cite{GKS22}, and where the cases $p=1,2$ follow from the case $p=0$, commutation with\footnote{Note that, for $r\leq r_+(1+2\dred)$,  we have $\N_{\nab_4^p\Ab}=\nab_4^p\N_{\Ab}+\dk^{\leq p+1}(\Ga_g\c\Ab)$, $p=1,2$, in view of Lemma \ref{lemma:commutator-nab4-square-redshift}.} $\nab_4^p$ and Lemma \ref{lemma:commutator-nab4-square-redshift}. The sign $\pr_r(\frac{\De}{|q|^2})>0$ in the region $r\leq r_+(1+2\dred)$ will allow us to control $\nab_4^p\Ab$, $p=0,1,2$,  near $r=r_+$ using redshift estimates, see Corollary \ref{cor:redshift:Ab:highregularity}.
\end{remark}


\subsubsection{Scalarized Teukolsky wave/transport systems in perturbations of Kerr}
\lab{sec:TeukolskyWavesysteminperturbationsofKerr:scalarizedform}


Applying Lemma \ref{lemma:formoffirstordertermsinscalarazationtensorialwaveeq} to the tensorial Teukolsky wave/transport systems \eqref{eq:TensorialTeuSysandlinearterms:rescaleRHScontaine2:general:Kerrperturbation}  \eqref{def:TensorialTeuScalars:wavesystem:Kerrperturbation} by using the regular triplet $\Om_i$, $i=1,2,3$ of Section  \ref{sec:regulartripletinperturbationsofKerr}, we obtain  the following scalarized Teukolsky wave/transport systems in perturbations of Kerr.  

\begin{lemma}[Scalarized Teukolsky wave/transport systems in perturbations of Kerr]
\lab{lem:scalarizedTeukolskywavetransportsysteminKerrperturbation:Omi}
Let $\Om_i$, $i=1,2,3$ be the regular triplet  of Section  \ref{sec:regulartripletinperturbationsofKerr}, and define the complex-valued scalars $\phiss{ij}{p}$ by 
\bea
\phiss{ij}{p}:=\pmb\phi_s^{(p)}(\Om_i,\Om_j), \quad s=\pm 2,\quad i,j=1,2,3, \quad p=0,1,2, 
\eea
where $\pmb\phi_s^{(p)}\in\sk_2(\mathbb{C})$, $s=\pm 2$, $p=0,1,2$, satisfy the tensorial Teukolsky wave/transport systems \eqref{eq:TensorialTeuSysandlinearterms:rescaleRHScontaine2:general:Kerrperturbation}  \eqref{def:TensorialTeuScalars:wavesystem:Kerrperturbation}. Then, the tensorial Teukolsky wave systems  \eqref{eq:TensorialTeuSysandlinearterms:rescaleRHScontaine2:general:Kerrperturbation} scalarized using $(\Om_i)_{i=1,2,3}$ take the following form 
\bea
\lab{eq:ScalarizedTeuSys:general:Kerrperturbation} 
\widehat\square_\g(\phi_s^{(p)})_{ij} -\frac{4-2\de_{p0}}{\qs} \phiss{ij}{p} = {L_{s,ij}^{(p)}}+N_{W,s,ij}^{(p)}, \quad s=\pm2, \quad i,j=1,2,3, \quad p=0,1,2,
\eea
where we have defined 
\bea
\lab{eq:defofwidehatsquaregoperator}
\widehat\square_\g(\phi_s^{(p)})_{ij}&:=&\square_\g\phiss{ij}{p}- \widehat{S}(\phi_s^{(p)})_{ij} -(\widehat{Q}\phi_s^{(p)})_{ij}
\eea
with
\bsub\lab{eq:definitionwidehatSandwidehatQperturbationsofKerr}
\begin{align}
\widehat{S}(\phi_s^{(p)})_{ij} ={}&S(\phi_s^{(p)})_{ij} +\frac{4ia\cos\th}{|q|^2} \pr_\tau\phiss{ij}{p},\\
(\widehat{Q}\phi_s^{(p)})_{ij} ={}&(Q\phi_s^{(p)})_{ij} 
-\frac{4ia\cos\th}{|q|^2}\big(M_{i\tau}^l \phiss{lj}{p}+M_{j\tau}^l \phiss{il}{p}\big),
\end{align}
\esub
where the linear coupling terms ${L_{s,ij}^{(p)}}:=(\L_{s}^{(p)}[\pmb\phi_s])_{ij}$ are given by
\bsub
\lab{eq:linearterms:ScalarizedTeuSys:general:Kerrperturbation}
\bea
{L_{s,ij}^{(0)}}&=& (2sr^{-3} +O(mr^{-4}))\phiss{ij}{1} + O(mr^{-3}){\Xcal_s}\phiss{ij}{0}\nn\\
&&+\sum_{k,l=1,2,3}O(mr^{-3})\phiss{kl}{0},\\
{L_{s,ij}^{(1)}}&=&(sr^{-3} +O(mr^{-4}))\phiss{ij}{2} + O(mr^{-3}){\Xcal_s}\phiss{ij}{1}+O(mr^{-2})(\pr_{\tphi} +a\pr_\tau)\phiss{ij}{0}\nn\\
  &&+\sum_{k,l=1,2,3}\Big(O(mr^{-3})\phiss{kl}{1}+O(mr^{-2})\phiss{kl}{0}\Big),\\
{L_{s,ij}^{(2)}}&=&O(mr^{-3}) \phiss{ij}{2}+O(mr^{-2})(\pr_{\tphi} +a\pr_\tau)\phiss{ij}{1} +O(m^2r^{-2})\phiss{ij}{0}\nn\\
&&+\sum_{k,l=1,2,3}O(mr^{-2})\phiss{kl}{1}
\eea
\esub
with $\Xcal_s$, $s=\pm 2$, being the regular vectorfields introduced in \eqref{eq:formofregularhorizontalvectorfieldmathcalXs:Kerrperturbation}, with the coefficients in front of the terms $\phiss{ij}{2}$ and $(\pr_{\tphi}+a\pr_{\tau})\phiss{ij}{1}$ on the RHS of equation of $L_{s,ij}^{(2)}$ in \eqref{eq:linearterms:ScalarizedTeuSys:general:Kerrperturbation} being real functions, and with all the coefficients in the first line\footnote{The coefficients on the second line of the three equations in \eqref{eq:linearterms:ScalarizedTeuSys:general:Kerrperturbation} involve $M_{i\a}^j$ so that they are independent of $\tau$ but depend on $\tphi$.} of the three equations in  \eqref{eq:linearterms:ScalarizedTeuSys:general:Kerrperturbation} being independent of coordinates $\tau$ and 
$\tphi$, where $S(\phi_s^{(p)})_{ij}$ and $(Q\phi_s^{(p)})_{ij}$ are given as in \eqref{SandV}, and where the complex-valued scalars $N_{W,s,ij}^{(p)}$ are given by 
\bea
N_{W,s,ij}^{(p)}:=\N_{W,s}^{(p)}(\Om_i, \Om_j),\quad s=\pm 2, \quad i,j=1,2,3, \quad p=0,1,2.
\eea

Moreover, the tensorial Teukolsky transport equations \eqref{def:TensorialTeuScalars:wavesystem:Kerrperturbation} scalarized using $(\Om_i)_{i=1,2,3}$ take the following form 
\bsub
 \lab{eq:ScalarizedQuantitiesinTeuSystem:Kerrperturbation}
\bea
\bsplit
& e_3\bigg(\frac{r\bar{q}}{q}\bigg({\frac{r^2}{|q|^2}}\bigg)^{p-2}\phipluss{ij}{p} \bigg) - \frac{r\bar{q}}{q}\bigg({\frac{r^2}{|q|^2}}\bigg)^{p-2}\Big(M_{i3}^k \phipluss{kj}{p} + M_{j3}^k \phipluss{ik}{p}\Big)\\
=&\frac{\bar{q}}{rq}\bigg({\frac{r^2}{|q|^2}}\bigg)^{p-1}\phipluss{ij}{p+1}+N_{T,+2,ij}^{(p)}
\end{split}
\eea
and 
\bea
\bsplit
& e_4\bigg(\frac{rq}{\bar{q}}\bigg(\frac{r^2}{|q|^2}\bigg)^{p-2}\phiminuss{ij}{p}\bigg) - \frac{rq}{\bar{q}}\bigg(\frac{r^2}{|q|^2}\bigg)^{p-2}\Big(M_{i4}^k \phiminuss{kj}{p} + M_{j4}^k \phiminuss{ik}{p}\Big)\\
=& \frac{q}{r\bar{q}}\left(\frac{r^2}{|q|^2}\right)^{p-1}\frac{\De}{\qs}\phiminuss{ij}{p+1} +N_{T,-2,ij}^{(p)},
\end{split}
\eea
\esub
where the complex-valued scalars $N_{T,s,ij}^{(p)}$ are given by 
\bea
N_{T,s,ij}^{(p)}:=\N_{T,s}^{(p)}(\Om_i, \Om_j),\quad s=\pm 2, \quad i,j=1,2,3, \quad p=0,1.
\eea
The equations \eqref{eq:ScalarizedTeuSys:general:Kerrperturbation} \eqref{eq:ScalarizedQuantitiesinTeuSystem:Kerrperturbation} then correspond to the scalarized Teukolsky wave/transport systems in perturbations of Kerr.
\end{lemma}

\begin{proof}
The proof follows along the same line as the one of Lemma \ref{lem:scalarizedTeukolskywavetransportsysteminKerr:Omi}.
\end{proof}


\subsection{Future null infinity of the perturbed {spacetime}}
\label{subsect:nullinf}


All the statements and estimates in this section are from \cite[Section 2.5]{MaSz24}.  We start by constructing an auxiliary ingoing optical function $\tauu$ in a subregion of $(\MM, \g)$. 
\begin{lemma}\lab{lemma:constructionoftheingoingopticalfunctiontauu}
There exists an ingoing optical function $\tauu$ defined in $\MM\cap\{r\geq |\tau|+ 10m\}$ by 
\bea
\tauu:=\tauu_0 +\tauut, \qquad \tauu_0:=\tau+2r+4m\log\left(\frac{r}{2m}\right),
\eea
where $\tauut$ satisfies 
\bea
|\dk^{\leq 2}\tauut|\les r^{-1}+\ep\quad\textrm{in}\quad\MM\cap\{r\geq |\tau|+ 10m\}.
\eea
\end{lemma}

Making use of the ingoing optical function $\tauu$, we may now define $\II_+$. 
\begin{definition}[Definition of $\II_+$]
\lab{def:howmathcalIplusisdefinedinMM} 
Consider the coordinates $(\tauu, \tau, x^1, x^2)$ covering the spacetime region $\MM\cap\{r\geq |\tau|+ 10m\}$, where $\tauu$ is the ingoing optical function constructed in Lemma \ref{lemma:constructionoftheingoingopticalfunctiontauu}. Then, the future null infinity of $(\MM, \g)$ is defined as
\bea
\II_+:=\MM\cap\{\tauu=+\infty\}.
\eea
\end{definition}

The following lemma provides the control of the induced geometry on $\II_+$ in the perturbed spacetime $(\MM, \g)$. 
\begin{lemma}\lab{lemma:controllinearizedmetric:inducedmetricII+}
Let $\II_+$ be given by Definition \ref{def:howmathcalIplusisdefinedinMM}. Consider the coordinates system $(\tau, x^1, x^2)$ covering $\II_+$, and denote by $(\pr_\tau^{\II_+}, \pr_{x^1}^{\II_+}, \pr_{x^2}^{\II_+})$ the corresponding coordinate vectorfields. Then, 
\begin{enumerate}
\item the coordinate vectorfields $\pr_{x^a}^{\II_+}$, $a=1,2,$ satisfy 
\bea\label{expression:prxaIIplus:nullinf}
\pr_{x^a}^{\II_+}=\pr_{x^a}+O(\ep)\pr_r,\,\,\, a=1,2,
\eea

\item the spheres $S^{\II_+}(\tau_1):=\II_+\cap\{\tau=\tau_1\}$ foliating $\II_+$ are round,

\item $\pr_\tau^{\II_+}$ is ingoing null and there exists a scalar function $b^r$ such that
\bea\label{expression:prtauIIplus:nullinf}
\pr_\tau^{\II_+}=\pr_\tau -\frac{1}{2}(1+b^r)\pr_r+O(\ep)\nab, \qquad |\dk^{\leq 1}b^r|\les \ep,
\eea

\item $\pr_r$ is an outgoing null vectorfield on $\II_+$ and satisfies 
\bea\label{expression:prrIIplus:nullinf}
\g(\pr_\tau^{\II_+}, \pr_r)=-1, \qquad \g(r^{-1}\pr_{x^a}^{\II_+}, \pr_r)=0.
\eea
\end{enumerate}
\end{lemma}


\subsection{Energy, Morawetz and flux norms}
\label{subsect:norms}


We introduce in this section the energy, Morawetz and flux norms needed to state our main result. First, given any $(\tau, r)$, and for any scalar function $F$ on the spheres $S(\tau, r)$ of constant $\tau$ and $r$, we introduce the following notation 
\beaa
\int_{\mathbb{S}^2}F(\tau, r, \om)d\mathring{\ga} := \int F(\tau, r, x^1, x^2)\sqrt{\det(\mathring{\ga})}dx^1dx^2,
\eeaa 
as well as the corresponding notation for the spheres $S^{\II_+}(\tau)$ of constant $\tau$ on $\II_+$. We start by introducing energy, Morawetz and flux norms for horizontal tensors in $\sk_k(\mathbb{C})$, $k=1,2$.


\subsubsection{Energy, Morawetz and flux norms for $\pmb\psi\in\sk_k(\mathbb{C})$, $k=1,2$}
\lab{sec:energyMorawetzFluxnormsforhorizontaltensorsinsk2C}


For $\pmb\psi\in\sk_k(\mathbb{C})$, $k=1,2$, and $\tau_1<\tau_2$, we define flux norms\footnote{For $\F_{\II_+}[\psi](\tau_1,\tau_2)$, recall that $\II_+=\MM\cap\{\tauu=+\infty\}$ where the ingoing optical function $\tauu$ has been constructed in Lemma \ref{lemma:constructionoftheingoingopticalfunctiontauu}, and recall that the notations $\pr_\tau^{\II_+}$ and $\pr_{x^a}^{\II_+}$ on $\II_+$ have been introduced in Lemma \ref{lemma:controllinearizedmetric:inducedmetricII+}.}
\bsub
\label{def:variousMorawetzIntegrals}
\bea
\bsplit
\F_{\AA}[\pmb\psi](\tau_1,\tau_2):=& \int_{\tau_1}^{\tau_2}\int_{\mathbb{S}^2}\big(|\mu| |\nab_3\pmb\psi|^2 +|\nab_4\pmb\psi|^2+|\nabla\pmb\psi|^2+|\pmb\psi|^2\big)(\tau, r=r_+(1-\dhor), \om){d\mathring{\ga}d\tau},\\
\F_{\II_+}[\pmb\psi](\tau_1,\tau_2):=& \int_{\tau_1}^{\tau_2}\int_{\mathbb{S}^2}\Bigg(|\nab_{\pr_\tau^{\II_+}}\pmb\psi|^2\\
&\qquad\qquad\quad +r^{-2}\left(|\nabla_{\pr_{x^1}^{\II_+}}\pmb\psi|^2+|\nabla_{\pr_{x^2}^{\II_+}}\pmb\psi|^2+|\pmb\psi|^2\right)\Bigg)(\tauu=+\infty, \tau, \om)r^2d\mathring{\ga}d\tau,\\
\F[\pmb\psi](\tau_1,\tau_2):=& \F_{\II_+}[\pmb\psi](\tau_1,\tau_2)+\F_{\AA}[\pmb\psi](\tau_1,\tau_2),
\end{split}
\eea
the energy norm
\bea
\E[\pmb\psi](\tau):= \int_{r_+(1-\dhor)}^{+\infty}\int_{\mathbb{S}^2}\Big(|\nab_4\pmb\psi|^2+|\nab\pmb\psi|^2+r^{-2}|\nab_3\pmb\psi|^2+r^{-2}|\pmb\psi|^2\Big)r^2{d\mathring{\ga}dr},
\eea
and the Morawetz norms
\bea
\bsplit
\M_\de[\pmb\psi](\tau_1,\tau_2):=&\int_{\MM_{\nontrap}(\tau_1,\tau_2)}\left(\frac{\abs{\nab_{\partial_{\tt}}\pmb\psi}^2}{r^{1+\de}} +\frac{\abs{\nabla\pmb\psi}^2}{r}\right)
+\int_{\MM(\tau_1,\tau_2)}\left(\frac{\abs{\nab_{\partial_r}\pmb\psi}^2}{r^{1+\de}} +\frac{\abs{\pmb\psi}^2}{r^3}\right),\\
\M[\pmb\psi](\tau_1,\tau_2):=&\int_{\MM_{\nontrap}(\tau_1,\tau_2)}\left(\frac{\abs{\nab_{\partial_{\tt}}\pmb\psi}^2}{r^2} +\frac{\abs{\nabla\pmb\psi}^2}{r}\right)
+\int_{\MM(\tau_1,\tau_2)}\left(\frac{\abs{\nab_{\partial_r}\pmb\psi}^2}{r^2} +\frac{\abs{\pmb\psi}^2}{r^3}\right),
\end{split}
\eea
\esub
for any given $0\leq\de\leq 1$. We also define, for $\H\in\sk_k(\mathbb{C})$, $k=1,2$,
\bea\lab{eq:defmathcalNpsif}
&&\mathcal{N}_\de[\pmb\psi, \H](\tau_1, \tau_2)\nn\\ 
&:=& \!\!\int_{\MM(\tau_1, \tau_2)}r^{1+\de}|\H|^2\nn\\
&+&\!\!\!\!\min\left[\left(\int_{\Mtrap(\tau_1, \tau_2)}|\H|^2\right)^{\frac{1}{2}} \left(\int_{\Mtrap(\tau_1, \tau_2)}|\dk{^{\leq 1}}\pmb\psi|^2\right)^{\frac{1}{2}}, 
\int_{\Mtrap(\tau_1, \tau_2)}\tau^{1+\de}|\H|^2\right]\!,
\eea
and  
\bea\lab{eq:defintionwidehatmathcalNfpsinormRHS}
&&\widehat{\mathcal{N}}[\pmb\psi, \H](\tau_1, \tau_2) \nn\\
\nn&:=&\sup_{\tau_1< \tau'<\tau''< \tau_2}\bigg|\int_{\Mntrap(\tau', \tau'')}\nab_{\pr_\tau}\pmb\psi\c\ov{\H}\bigg|+\int_{\Mntrap(\tau_1, \tau_2)}r^{-1}|\dk^{\leq 1}\pmb\psi||\H|+\int_{\MM(\tau_1, \tau_2)}|\H|^2\nn\\
&+&\!\!\!\!\min\left[\left(\int_{\Mtrap(\tau_1, \tau_2)}|\H|^2\right)^{\frac{1}{2}} \left(\int_{\Mtrap(\tau_1, \tau_2)}|\dk{^{\leq 1}}\pmb\psi|^2\right)^{\frac{1}{2}}, 
\int_{\Mtrap(\tau_1, \tau_2)}\tau^{1+\de}|\H|^2\right].
\eea
Also, we define $\widehat{\mathcal{N}}_{r\geq r_0}[\pmb\psi, \H](\tau_1, \tau_2)$, for any $r_0\geq r_+(1-\dhor)$, similarly to the formula \eqref{eq:defintionwidehatmathcalNfpsinormRHS} of  $\widehat{\mathcal{N}}[\pmb\psi, \H](\tau_1, \tau_2)$ but with all integrals further restricted to $r\geq r_0$, i.e.,
\bea\lab{eq:defintionwidehatmathcalNfpsinormRHS:rgeqr0}
&&\widehat{\mathcal{N}}_{r\geq r_0}[\pmb\psi, \H](\tau_1, \tau_2) \nn\\
\nn&:=&\sup_{\tau_1< \tau'<\tau''< \tau_2}\bigg|\int_{\MM_{\nontrap, r\geq r_0}(\tau', \tau'')}\nab_{\pr_\tau}\pmb\psi\c\ov{\H}\bigg|+\int_{\MM_{\nontrap, r\geq r_0}(\tau_1, \tau_2)}r^{-1}|\dk^{\leq 1}\pmb\psi||\H|\\
&+&\!\!\!\!\min\left[\left(\int_{\MM_{\trap, r\geq r_0}(\tau_1, \tau_2)}|\H|^2\right)^{\frac{1}{2}} \left(\int_{\MM_{\trap, r\geq r_0}(\tau_1, \tau_2)}|\dk{^{\leq 1}}\pmb\psi|^2\right)^{\frac{1}{2}}, 
\int_{\MM_{\trap, r\geq r_0}(\tau_1, \tau_2)}\tau^{1+\de}|\H|^2\right]\nn\\
&+&\int_{\MM_{r\geq r_0}(\tau_1, \tau_2)}|\H|^2.
\eea

\begin{remark}\lab{rmk:controlofwidehatNfpsibyNfpsi}
In view of the above definitions, we immediately deduce the following bound
\beaa
\widehat{\mathcal{N}}[\pmb\psi, \H](\tau_1, \tau_2)\les \Big(\M_\de[\pmb\psi](\tau_1,\tau_2)\Big)^{\frac{1}{2}}\Big(\mathcal{N}_\de[\pmb\psi, \H](\tau_1, \tau_2)\Big)^{\frac{1}{2}}+\mathcal{N}_\de[\pmb\psi, \H](\tau_1, \tau_2).
\eeaa
Also, note that $\M[\pmb\psi](\tau_1,\tau_2)=\M_1[\pmb\psi](\tau_1,\tau_2)$. 
\end{remark}

Next, for any nonnegative integer $\reg$, let
\bea\lab{eq:definitionofhigherordernormsFregEregMdeltaregMregNNdeltaregwidehatNreg}
\bsplit
\F^{(\reg)}[\pmb\psi](\tau_1,\tau_2)&:=\F[\dk^{\leq\reg}\pmb\psi](\tau_1,\tau_2),\qquad\qquad\qquad\quad \E^{(\reg)}[\pmb\psi](\tau):=\E[\dk^{\leq\reg}\pmb\psi](\tau), \\ 
\M^{(\reg)}_\de[\pmb\psi](\tau_1,\tau_2)&:=\M_\de[\dk^{\leq\reg}\pmb\psi](\tau_1,\tau_2),\qquad\qquad\, \M^{(\reg)}[\pmb\psi](\tau_1,\tau_2):=\M[\dk^{\leq\reg}\pmb\psi](\tau_1,\tau_2),\\
\mathcal{N}^{(\reg)}_\de[\pmb\psi, \H](\tau_1, \tau_2)&:=\mathcal{N}_\de[\dk^{\leq\reg}\pmb\psi, \dk^{\leq \reg}\H](\tau_1, \tau_2),\quad\widehat{\mathcal{N}}^{(\reg)}[\pmb\psi, \H](\tau_1, \tau_2):=\widehat{\mathcal{N}}[\dk^{\leq\reg}\pmb\psi, \dk^{\leq \reg}\H](\tau_1, \tau_2).
\end{split}
\eea

Finally, we define for any nonnegative integer $\reg$ the following combined norms 
\bea\lab{eq:definitionofhigherordernormsFregEregMdeltaregMregNNdeltaregwidehatNreg:combinednorms}
\bsplit
\EMF^{(\reg)}_\de[\pmb\psi](\tau_1,\tau_2) := \sup_{\tt\in [\tau_1, \tau_2]} \E^{(\reg)}[\pmb\psi](\tt) + \M^{(\reg)}_\de[\pmb\psi](\tau_1,\tau_2)+\F^{(\reg)}[\pmb\psi](\tau_1,\tau_2),\\
\EMF^{(\reg)}[\pmb\psi](\tau_1,\tau_2) := \sup_{\tt\in [\tau_1, \tau_2]} \E^{(\reg)}[\pmb\psi](\tt) + \M^{(\reg)}[\pmb\psi](\tau_1,\tau_2)+\F^{(\reg)}[\pmb\psi](\tau_1,\tau_2),
\end{split}
\eea
with $\EM^{(\reg)}_\de[\pmb\psi](\tau_1,\tau_2)$, $\EM^{(\reg)}[\pmb\psi](\tau_1,\tau_2)$, $\MF^{(\reg)}[\pmb\psi](\tau_1,\tau_2)$ and $\EF^{(\reg)}[\pmb\psi](\tau_1,\tau_2)$ being defined in a similar way.


\subsubsection{Energy, Morawetz and flux norms for scalars}
\lab{sec:energyMorawetzFluxnormsforscalars}


For any scalars $\psi$ and $H$, we define the norms $\F_\AA[\psi]$, $\F_{\II_+}[\psi]$, $\F[\psi]$, $\E[\psi]$, $\M_\de[\psi]$, $\M[\psi]$, $\mathcal{N}_\de[\psi, H]$, $\widehat{\mathcal{N}}[\psi, H]$, as well as the corresponding higher order derivative norms and combined quantities as in Section  \ref{sec:energyMorawetzFluxnormsforhorizontaltensorsinsk2C} by replacing $\pmb\psi, \,\H\in\sk_k(\mathbb{C})$, $k=1,2$, with scalars $\psi$, $H$ in the formulas.


\subsubsection{Energy, Morawetz and flux norms for scalarized tensors using regular triplets}


We start with the following lemma.
\begin{lemma}\lab{lemma:equivalentofmodulussquarepmbpsiandsumijpisijalsowithdkreg}
Let $\pmb\psi\in\sk_2(\mathbb{C})$ and let $\psi_{ij}$ be the corresponding scalars given by $\psi_{ij}:=\pmb\psi(\Om_i, \Om_j)$ for $i,j=1,2,3$ where $\Om_i$, $i=1,2,3$ is the regular triplet introduced in Section  \ref{sec:regulartripletinperturbationsofKerr}. Then, for any integer $\reg\leq {15}$, we have, for $\a=a, 3, 4$, $a=1,2$, 
\beaa
|\dk^{\leq \reg}\pmb\psi|^2\simeq \sum_{i,j=1}^3|\dk^{\leq\reg}(\psi_{ij})|^2, \qquad |\Ddot_\a\dk^{\leq\reg}\pmb\psi|^2\simeq \sum_{i,j=1}^3|e_\a(\dk^{\leq\reg}(\psi_{ij}))|^2+O(r^{-2})\sum_{i,j=1}^3|\dk^{\leq\reg}(\psi_{ij})|^2.
\eeaa
\end{lemma}

\begin{proof}
In view of Lemma \ref{lem:ucdotv:product}, we have
\beaa
|\pmb\psi|^2=\sum_{i,j}|\psi_{ij}|^2, \qquad |\Ddot_\a\pmb\psi|^2=\sum_{i,j}|(\Ddot_\a\pmb\psi)(\Om_i, \Om_j)|^2, \quad \a=a, 3, 4, \,\, a=1,2.
\eeaa
Also, we have 
\beaa
(\Ddot_\a\pmb\psi)(\Om_i, \Om_j)=e_\a(\psi_{ij}) - M_{i\a}^k\psi_{kj}-M_{j\a}^k\psi_{ik},
\eeaa
and since $M_{i\a}^j=O(r^{-1})$ as a consequence of Lemma \ref{lemma:computationoftheMialphajinKerr} and \eqref{eq:assumptionsonregulartripletinperturbationsofKerr:0}, we infer
\beaa
(\Ddot_\a\pmb\psi)(\Om_i, \Om_j)=e_\a(\psi_{ij}) +O(r^{-1})\psi_{kl}.
\eeaa
Hence, in view of the above, we have 
\beaa
|\pmb\psi|^2=\sum_{i,j}|\psi_{ij}|^2, \qquad |\Ddot_\a\pmb\psi|^2\simeq\sum_{i,j}|e_\a(\psi_{ij})|^2+O(r^{-2})\sum_{i,j}|\psi_{ij}|^2, \quad \a=a, 3, 4, \,\, a=1,2.
\eeaa
The general case for $\reg\leq 15$ follows by iteration on $\reg$, using again the fact that $\dk^{\leq 15}M_{i\a}^j=O(r^{-1})$ again as a consequence of Lemma \ref{lemma:computationoftheMialphajinKerr} and \eqref{eq:assumptionsonregulartripletinperturbationsofKerr:0}. This concludes the proof of the lemma.
\end{proof}

Then, recalling the norms introduced in Sections \ref{sec:energyMorawetzFluxnormsforhorizontaltensorsinsk2C} and \ref{sec:energyMorawetzFluxnormsforscalars}, we infer  immediately from Lemma \ref{lemma:equivalentofmodulussquarepmbpsiandsumijpisijalsowithdkreg} the following equivalence relations between the norms for tensors and the norms for the scalars obtained by scalarization using regular triplets.
\begin{corollary}
[Equivalence of the norms for tensors and for scalars using regular triplets]
\lab{coro:equivalenceofnormsfortensorsandscalars}
Let $\pmb\psi, \,\H\in\sk_2(\mathbb{C})$, and let $\psi_{ij}, H_{ij}$ be the corresponding scalars given by $\psi_{ij}:=\pmb\psi(\Om_i, \Om_j)$ and $H_{ij}:=\H(\Om_i, \Om_j)$, for $i,j=1,2,3$, respectively, where $\Om_i$, $i=1,2,3$ is the regular triplet introduced in Section  \ref{sec:regulartripletinperturbationsofKerr}. 
Then we have the equivalence relations, for any $\reg\leq 14$, $\de\in[0,1]$ and $\tau_1<\tau_2$,
\bea\lab{eq:equivalenceofnormsfortensorsandscalars}
\E^{(\reg)}[\pmb\psi](\tau)\simeq \sum_{i,j=1}^3\E^{(\reg)}[\psi_{ij}](\tau), \qquad \mathcal{N}^{(\reg)}_\de[\pmb\psi, \H](\tau_1, \tau_2)\simeq \sum_{i,j=1}^3\mathcal{N}^{(\reg)}_\de[\psi_{ij}, H_{ij}](\tau_1, \tau_2),
\eea
and similarly for all the other norms appearing in Section  \ref{sec:energyMorawetzFluxnormsforhorizontaltensorsinsk2C}.
\end{corollary}

\begin{remark}\lab{rmk:abusenotationbetweentensorandscalarizedversioninEMFnorms}
More generally, we will also consider family of scalars $\psi_{ij}$, $i,j=1,2,3,$ that are not generated\footnote{Recall from Lemma \ref{lemma:backandforthbetweenhorizontaltensorsk2andscalars:complex} that complex-valued scalars $\psi_{ij}$, $i,j=1,2,3,$ come from the scalarization of a tensor in $\sk_2(\mathbb{C})$ if and only if they satisfy the identities stated in the second item of Lemma \ref{lemma:backandforthbetweenhorizontaltensorsk2andscalars:complex}.} by the scalarization of a tensor in $\sk_2(\mathbb{C})$. In that case\footnote{This will be the case in $\MM(\tau_2-2, \tau_2)$ due to the semi global extension procedure of Proposition  \ref{prop:extensionprocedureoftheTeukolskywaveequations} which does not preserve the identities stated in the second item of Lemma \ref{lemma:backandforthbetweenhorizontaltensorsk2andscalars:complex}.}, by a slight abuse of notation, we will still denote the norms appearing on the RHS of the identities in  \eqref{eq:equivalenceofnormsfortensorsandscalars} by the ones appearing on the LHS (even though the corresponding tensors $\pmb\psi$, $\H$, do not exist).
\end{remark}


\section{Basic estimates for wave and transport equations in perturbations of Kerr}
\lab{sect:basicestimatesforwaveequations}


In this section, we prove basic estimates for (scalar and tensorial) wave equations, as well as transport equations,  in perturbations of Kerr in the range $|a|<m$.


\subsection{Standard calculation for generalized currents}


Consider  variational  wave equations  for  tensors  $\pmb\psi\in \sk_k(\mathbb{C})$, $k=0,1,2$, of the form 
\bea\lab{eq:Gen.RW-general}
\squared_k \pmb\psi-V\pmb\psi=\pmb{N},
\eea
where $V$ is a real potential. The variational wave equation \eqref{eq:Gen.RW-general} has Lagrangian
\beaa
 \LL[\pmb\psi]:= \g^{\mu\nu}\Re\left(\Db_\mu \pmb\psi\c\ov{\Db_\nu  \pmb\psi }\right)+ V |\pmb\psi |^2,
 \eeaa
  where  the dot product   here denotes full contraction with respect to the  horizontal indices.
The  corresponding   energy-momentum tensor associated to \eqref{eq:Gen.RW-general} is given by 
 \bea\label{eq:definition-QQ-mu-nu}
\nn  \QQ_{\mu\nu}[\pmb\psi] &:=& \Re\Big(\Db_\mu  \pmb\psi \c \ov{\Db _\nu  \pmb\psi }\Big)
          -\frac 12 \g_{\mu\nu} \left( \Re\left(\Db_\la  \pmb\psi\c\ov{\Db^\la  \pmb\psi }\right)+ V |\pmb\psi |^2\right)\\
&=& \Re\Big(\Db_\mu  \pmb\psi \c \ov{\Db _\nu  \pmb\psi }\Big)  -\frac 12 \g_{\mu\nu}  \LL[\pmb\psi]. 
 \eea
Also, recall from \eqref{def:deformationtensor:lastsect} that the deformation tensor of a vectorfield $X$ is defined by
\beaa
{}^{(X)}\pi_{\a\b} =  \D_{\a}X_{\b} + \D_{\b}X_{\a}.
\eeaa

\begin{lemma}\label{lemma-divergence-QQ}
Given a solution $\pmb\psi\in \sk_k(\mathbb{C})$, $k=0,1,2$, to equation  \eqref{eq:Gen.RW-general} we have
 \bea\lab{eq:divergence-QQ:complexcase}
 \D^\nu\QQ_{\mu\nu}[\pmb\psi]
  &=& \Re\Big(  \left(\squared_k \pmb\psi- V\pmb\psi\right)\c \ov{\Db_\mu  \pmb\psi }\Big)+ 
   \frac{k}{2}\in^{ab}\Rdot_{ ab\nu\mu}\Im\Big(\pmb\psi\c\ov{\Ddot^{\nu}\pmb\psi}\Big) -\frac 1 2 \D_\mu V |\pmb\psi|^2.
 \eea
\end{lemma}

\begin{proof}
We have from Lemma 4.7.1 in \cite{GKS22} that for a solution $\phi\in \sk_k$ to equation \eqref{eq:Gen.RW-general}, 
\beaa
 \D^\nu\QQ_{\mu\nu}[\phi] &=& \Db_\mu\phi \c  \left(\squared_k\phi- V\phi\right)+ E_k[\phi] -\frac 1 2 \D_\mu V |\phi|^2,
 \eeaa
 where $E_k[\phi]$, $k=0,1,2$, is given by 
 \beaa
 E_0[\phi]=0, \qquad E_1[\phi]=\Db^\nu\phi ^a\Rdot_{ a   b   \nu\mu}\phi^b, \qquad E_2[\phi]=\Db^\nu\phi ^{ac}\Rdot_{ab\nu\mu}\phi^b\,\!_{c}.
 \eeaa
 Using the antisymmetry in $(a,b)$ of $\Rdot_{ab\nu\mu}$, we infer
 \bea\lab{eq:divergence-QQ:realcase}
 \D^\nu\QQ_{\mu\nu}[\phi] &=& \Db_\mu\phi \c  \left(\squared_k\phi- V\phi\right)+ \frac{k}{2}\in^{ab}\Rdot_{ ab\nu\mu}\Big(\dual\phi\c\Ddot^{\nu}\phi\Big) -\frac 1 2 \D_\mu V |\phi|^2.
 \eea
 
Next, we consider $\pmb\psi \in \sk_k(\mathbb{C})$, which, according to Definition \ref{def:skC:horizontaltensors}, is given by $\pmb\psi=\psi+i\dual\psi$ with $\psi\in\sk_k$. Notice that $\QQ_{\mu\nu}[\pmb\psi]=\QQ_{\mu\nu}[\psi]+\QQ_{\mu\nu}[\dual\psi]$ so that the proof follows from applying \eqref{eq:divergence-QQ:realcase} respectively with $\phi=\psi$ and $\phi=\dual\psi$ and summing the two resulting identities.
\end{proof}

We collect here some general calculations for generalized currents associated to equation \eqref{eq:Gen.RW-general}.

\begin{proposition}\lab{prop-app:stadard-comp-Psi}
 Let   $\pmb\psi\in \mathfrak{s}_k(\mathbb{C})$, $k=0,1,2$, and  let $X$ be  a real-valued vectorfield and $w$ a real scalar function.  Define the $1$-form $\PP_\mu[\pmb\psi](X, w)$ by
 \bea
\lab{definitionofcurrentPPmuXw:generaltensor}
\PP_\mu[\pmb\psi](X, w):=\QQ_{\mu\nu}[\pmb\psi] X^\nu +\frac 1 2  w \Re\Big(\pmb\psi \c \ov{\Db_\mu \pmb\psi }\Big)-\frac 1 4|\pmb\psi|^2   \pr_\mu w,
  \eea
and define the 1-form ${}^{(X)}A_\nu$  by
\bea\lab{eq:thespacetime1formXA}
{}^{(X)}A_\nu &:=& X^\mu \in^{ab}\Rdot_{ ab   \nu\mu}.
\eea
Then, we have
 \bea\lab{eq:DivofPPmu:tensor:RW:prop}
  \D^\mu  \PP_\mu[\pmb\psi](X, w)&=& \frac 1 2 \QQ[\pmb\psi]  \c\piX - \frac 1 2 X( V ) |\pmb\psi|^2 +\frac{k}{2}{}^{(X)}A_\nu\Im\Big(\pmb\psi\c\ov{\Ddot^{\nu}\pmb\psi}\Big)\nn\\
  &&+\frac 12  w \LL[\pmb\psi] -\frac 1 4|\pmb\psi|^2   \square_\g  w   +  \Re\bigg(\ov{\bigg(\nab_X\pmb\psi +\frac 1 2   w \pmb\psi\bigg)}\c \left(\squared_k \pmb\psi- V\pmb\psi\right)\bigg).
 \eea
\end{proposition}

\begin{proof}
This is an extension to $\sk_k(\mathbb{C})$ of part of Proposition 4.7.2 in \cite{GKS22}. We have
\beaa
\Ddot^\mu\left[\frac 1 2  w \Re\Big(\pmb\psi \c \ov{\Db_\mu \pmb\psi }\Big)-\frac 1 4|\pmb\psi|^2\pr_\mu w \right] &=& \frac{1}{2}w\Re\left(\pmb\psi\c\ov{\squared_k\pmb\psi}\right) +\frac{1}{2}w\Re\left(\Ddot^\mu\pmb\psi\c\ov{\Ddot_\mu\pmb\psi}\right) -\frac{1}{4}|\pmb\psi|^2\square_\g w\\
&=& \frac 12  w \LL[\pmb\psi] -\frac 1 4|\pmb\psi|^2   \square_\g  w   +  \frac 1 2   w\Re\left(\ov{\pmb\psi}\c\left(\squared_k \pmb\psi- V\pmb\psi\right)\right)
\eeaa
which together with \eqref{eq:divergence-QQ:complexcase} and the definition of $^{(X)}A_\nu$ in \eqref{eq:thespacetime1formXA} concludes the proof of the proposition.
\end{proof}

  
\subsubsection{Term of the type $\pmb\psi\c\Ddot_\nu\pmb\psi$ in \eqref{eq:DivofPPmu:tensor:RW:prop}}
\label{sec:somekerrvalues}
  

In this section, we control the term of the type $\pmb\psi\c\Ddot_\nu\pmb\psi$ in \eqref{eq:DivofPPmu:tensor:RW:prop}. 

\begin{lemma}\lab{lemma:computationofthecomponentsofthetensorAinKerr}
Let $(\MM, \g)$ satisfy the assumptions of Section \ref{subsect:assumps:perturbednullframe}. Let ${}^{(X)}A_{\nu}$ be the spacetime 1-form given by \eqref{eq:thespacetime1formXA}. Then, 
we have
\bsub
\begin{align*}
{}^{(X)}A_4 ={}& -4\rhod X^3 -4(\etab\wedge\eta)X^3 +  \trch\big({}^{(h)}X\wedge\etab)  - \atrch\big(\etab\c {}^{(h)}X) {+\Ga_b\c\Ga_g X^3+r^{-1}\Ga_g {}^{(h)}X},\\
{}^{(X)}A_3 ={}& 4\rhod X^4 +4(\etab\wedge\eta)X^4 +  \trchb\big({}^{(h)}X\wedge\eta)  - \atrchb\big(\eta\c {}^{(h)}X\big) {+\Ga_b\c\Ga_g X^4+r^{-1}\Ga_b {}^{(h)}X},\\
{}^{(X)}A_e ={}& \Big( -  \trchb\dual\eta_e   + \atrchb\eta_e\Big)X^3 +\Big( -  \trch\dual\etab_e   + \atrch\etab_e\Big)X^4\nn\\
{}& -\frac{1}{2}\Big(4\rho  + \trch\trchb+\atrch\atrchb\Big)\dual ({}^{(h)}X)_e {+r^{-1}\Ga_b X^3 +r^{-1}\Ga_g X^4+\Ga_b\c\Ga_g {}^{(h)}X }.
\end{align*}
\esub
\end{lemma}

\begin{proof}
We rewrite ${}^{(X)}A_\mu$ as 
\beaa
{}^{(X)}A_\mu &=& \in^{bc}\Rdot_{bc \mu 3}X^3 +\in^{bc}\Rdot_{bc \mu 4}X^4 +\in^{bc}\Rdot_{bc \mu d}X^d.
\eeaa
Next, we compute the various components of ${}^{(X)}A_\mu$. {We have}, using the horizontal tensor ${}^{(h)}X$ defined by $({}^{(h)}X)_b=X_b$, the definition  \eqref{eq:DefineRdot} of $\Rdot$,  and Proposition \ref{proposition:componentsofB},
\bsub
\begin{align*}
{}^{(X)}A_4 ={}& \in^{bc}\Rdot_{bc43}X^3  +\in^{bc}\Rdot_{bc4d}X^d\nn\\
={}& \in^{bc}\Big(-2\in_{bc}\dual\rho -2(\etab_b\eta_c-\eta_b\etab_c) {+\Ga_b\c\Ga_g}\Big)X^3\nn\\
{}& +  \frac{1}{2}\in^{bc}\Big(\trch(\de_{db}\etab_c-\de_{dc}\etab_b) + \atrch(\in_{db}\etab_c - \in_{dc}\etab_b) {+r^{-1}\Ga_g}\Big)X^d\nn\\
={}& -4\rhod X^3 -4(\etab\wedge\eta)X^3 {+\Ga_b\c\Ga_g X^3+r^{-1}\Ga_g {}^{(h)}X}\nn\\
{}&+  \frac{1}{2}\in^{bc}\Big(\trch(\etab_c X_b - \etab_b X_c) + \atrch(-\etab_c\dual ({}^{(h)}X)_b + \dual ({}^{(h)}X)_c \etab_b)\Big)\nn\\
={}& -4\rhod X^3 -4(\etab\wedge\eta)X^3 +  \trch\big({}^{(h)}X\wedge\etab)  - \atrch\big(\etab\c {}^{(h)}X) {+\Ga_b\c\Ga_g X^3+r^{-1}\Ga_g {}^{(h)}X},\\
{}^{(X)}A_3 ={}& \in^{bc}\Rdot_{bc34}X^4 +\in^{bc}\Rdot_{bc3d}X^d\nn\\
={}& 4\rhod X^4 +4(\etab\wedge\eta)X^4 +  \trchb\big({}^{(h)}X\wedge\eta)  - \atrchb\big(\eta\c {}^{(h)}X\big) {+\Ga_b\c\Ga_g X^4+r^{-1}\Ga_b {}^{(h)}X},\\
{}^{(X)}A_e ={}& \in^{bc}\Rdot_{bce3}X^3 +\in^{bc}\Rdot_{bce4}X^4 +\in^{bc}\Rdot_{bced}X^d\nn\\
={}& \frac{1}{2}\in^{bc}\Big( -  \trchb  \big( \de_{eb}\eta_c -  \de_{ec} \eta_b\big)  -  \atrchb \big( \in_{eb}  \eta_c -  \in_{ec}  \eta_b\big) {+r^{-1}\Ga_b}\Big)X^3\nn\\
{}&+ \frac{1}{2}\in^{bc}\Big( -  \trch  \big( \de_{eb}\etab_c -  \de_{ec} \etab_b\big)  -  \atrch \big( \in_{eb}  \etab_c -  \in_{ec}  \etab_b\big) {+r^{-1}\Ga_g}\Big)X^4\nn\\
{}& +\frac{1}{2}\in^{bc}\left(-2\in_{bc}\in_{ed}\rho  -\frac{1}{2}\left(\trch\trchb+\atrch\atrchb\right)\in_{bc}\in_{ed} {+\Ga_b\c\Ga_g}\right)X^d\nn\\
={}& \Big( -  \trchb\dual\eta_e   + \atrchb\eta_e\Big)X^3 +\Big( -  \trch\dual\etab_e   + \atrch\etab_e\Big)X^4\nn\\
{}& -\frac{1}{2}\Big(4\rho  + \trch\trchb+\atrch\atrchb\Big)\dual ({}^{(h)}X)_e {+r^{-1}\Ga_b X^3 +r^{-1}\Ga_g X^4+\Ga_b\c\Ga_g {}^{(h)}X }
\end{align*}
\esub
as stated. This concludes the proof of Lemma \ref{lemma:computationofthecomponentsofthetensorAinKerr}.
\end{proof}

We infer the following corollaries.

\begin{corollary}\lab{cor:computationofthecomponentsofthetensorAforexactenergyconservationKerr}
Let $(\MM, \g)$ satisfy the assumptions of Section \ref{subsect:assumps:perturbednullframe}. Then, we have
\bea
{}^{({\pr_\tau})}A_\mu=-\D_\mu\left(\Im\left(\frac{2m}{q^2}\right)\right) {+r^{-1}\Ga_b}.
\eea
\end{corollary}

\begin{proof}
{In view of \eqref{eq:relationsbetweennullframeandcoordinatesframe:2}, we have
\beaa
\pr_\tau = O(1)e_4+O(1)e_3+\big(O(mr^{-1})+O(\ep)\big)e_a.
\eeaa
Hence, applying Lemma \ref{lemma:computationofthecomponentsofthetensorAinKerr} with $(\pr_\tau)^3=O(1)$, $(\pr_\tau)^4=O(1)$ and $(\pr_\tau)^b=O(1)$, we infer
\beaa
{}^{(\pr_\tau)}A_\mu=\big[{}^{(\pr_\tau)}A_\mu\big]_K+r^{-1}\Ga_b
\eeaa
where $[{}^{(\pr_\tau)}A_\mu]_K$ denotes the corresponding Kerr value. The proof then  follows from the value of $[{}^{(\pr_\tau)}A_\mu]_K$ in Corollary 7.7.3 in \cite{GKS22}.}
\end{proof}

\begin{corollary}\lab{cor:asymtpticbehavioroftheRdottermintensorialenergyidentity:larger}
Let $(\MM, \g)$ satisfy the assumptions of Section \ref{subsect:assumps:perturbednullframe}. For $X$ such that 
\beaa
X^4=O(1), \qquad X^3=O(1), \qquad X^b=O(r^{-1}), \,\, b=1,2,
\eeaa
we have 
\beaa
{}^{(X)}A_\nu\Im\Big(\pmb\psi\c\ov{\Ddot^{\nu}\pmb\psi}\Big) = O(r^{-3})\Im\Big(\pmb\psi\c\ov{\nab_4\pmb\psi }\Big)+ O(r^{-3})\Im\Big(\pmb\psi\c\ov{\nab_3\pmb\psi }\Big) + \Big(O(r^{-3})+r^{-1}\Ga_b\Big)\Im\Big(\pmb\psi\c\ov{\nab\pmb\psi }\Big).  
\eeaa
\end{corollary}

\begin{proof}
 Using Lemma \ref{lemma:computationofthecomponentsofthetensorAinKerr}, the assumptions on the components of $X$ and the non-sharp asymptotic behavior 
 \beaa
\trch,\, \trchb=O(r^{-1}), \qquad  \eta=O(r^{-2})+\Ga_b,\qquad \atrch,\, \atrchb,\, \etab=O(r^{-2}), \qquad \rho,\, \dual\rho=O(r^{-3}), 
 \eeaa
 we have the following non-sharp asymptotic behavior for the components of ${}^{(X)}A_\nu$
 \beaa
 {}^{(X)}A_4=O(r^{-3}), \qquad {}^{(X)}A_3=O(r^{-3}), \qquad {}^{(X)}A_a=O(r^{-3})+r^{-1}\Ga_b, 
 \eeaa
 which together with the identity
 \beaa
 {}^{(X)}A_\nu\Im\Big(\pmb\psi\c\ov{\Ddot^{\nu}\pmb\psi}\Big) &=& -\frac{1}{2}{}^{(X)}A_3\Im\Big(\pmb\psi\c\ov{\nab_4\pmb\psi }\Big) -\frac{1}{2}{}^{(X)}A_4\Im\Big(\pmb\psi\c\ov{\nab_3\pmb\psi }\Big)+{}^{(X)}A_a\Im\Big(\pmb\psi\c\ov{\nab_a\pmb\psi }\Big)
 \eeaa
 concludes the proof of the corollary.
 \end{proof}


\subsubsection{An energy identity for a tensorial wave equation}


We consider solutions to the following tensorial wave equation for $\pmb\psi\in \sk_k(\mathbb{C})$, $k=0,1,2$,
\bea
\lab{eq:tensorialwaveRW:withprttderivative}
 \squared_k \pmb\psi-\frac{4ia\cos\th}{\qs}\nab_{\pr_{\tt}}\pmb\psi-V\pmb\psi=\pmb{N}
 \eea
and, as in Section  7.3 of \cite{GKS22}, we make use of Corollary \ref{cor:computationofthecomponentsofthetensorAforexactenergyconservationKerr} to derive an energy identity.

\begin{lemma}\lab{cor:modifiedcurrentsforprtandprvphi}
Let $(\MM, \g)$ satisfy the assumptions of Section \ref{subsect:assumps:perturbednullframe}. Let $\pmb\psi\in\sk_k(\mathbb{C})$, $k=0,1,2$, be a solution to the tensorial wave equation \eqref{eq:tensorialwaveRW:withprttderivative}, with the real potential $V$ satisfying $\pr_{\tt}V=0$.
Define 
\bea
\widetilde{w} := \Im\left(\frac{m}{q^2}\right)= -\frac{2amr\cos\th}{|q|^4},
\eea
and define the following modified current associated to the vectorfield $\pr_{\tt}$:
\bea
{}^{(\pr_{\tt})}\widetilde{\PP}_\mu[\pmb\psi] :=  {\PP}_\mu[\pmb\psi](\pr_{\tt}, 0) + k\widetilde{w}\Im\left( \pmb\psi\c\ov{\Ddot_\mu\pmb\psi}\right)+(\pr_{\tt})_{\mu}\frac{2a\cos\th}{\qs}k\widetilde{w}|\pmb\psi|^2.
\eea
Then, we have
\bea
\D^\mu{}^{(\pr_{\tt})}\widetilde{\PP}_\mu [\pmb\psi]&=&  \Re\bigg(\ov{\big(\nab_{\pr_{\tt}}\pmb\psi-ik\widetilde{w}\pmb\psi\big)}\c \left(\squared_k\pmb\psi-\frac{4ia\cos\th}{\qs}\nab_{\pr_{\tt}}\pmb\psi-V\pmb\psi\right)\bigg)\nn\\
&&+\frac 1 2 \QQ[\pmb\psi]  \c {^{(\pr_{\tt})}}\pi+\Div (\pr_{\tt})\frac{2a\cos\th}{\qs}k\widetilde{w}|\pmb\psi|^2+r^{-1}\Ga_b\Im\Big(\pmb\psi\c\ov{\Ddot^{\nu}\pmb\psi}\Big).
\eea
\end{lemma}

\begin{remark}
In the case $\g=\gam$ and $\pmb{N}=0$, this induces a conservation of energy.
\end{remark}

\begin{proof}
In view of \eqref{eq:DivofPPmu:tensor:RW:prop}, \eqref{eq:thespacetime1formXA}, and the fact that $\pr_{\tt}V=0$, we have
\beaa
 \D^\mu  \PP_\mu[\pmb\psi](\pr_{\tt}, 0)&=& \frac 1 2 \QQ[\pmb\psi]  \c{}^{(\pr_{\tt})}\pi +\frac{k}{2}{}^{(\pr_{\tt})}A_\nu\Im\Big(\pmb\psi\c\ov{\Ddot^{\nu}\pmb\psi}\Big) +  \Re\bigg(\ov{\nab_{\pr_{\tt}}\pmb\psi}\c \left(\squared_k \pmb\psi- V\pmb\psi\right)\bigg)\\
 &=& \frac 1 2 \QQ[\pmb\psi]  \c{}^{(\pr_{\tt})}\pi -k\D_\nu(\widetilde{w})\Im\Big(\pmb\psi\c\ov{\Ddot^{\nu}\pmb\psi}\Big) \\
 &&+  \Re\bigg(\ov{\nab_{\pr_{\tt}}\pmb\psi}\c \left(\squared_k \pmb\psi -\frac{4ia\cos\th}{\qs}\nab_{\pr_{\tt}}\pmb\psi - V\pmb\psi\right)\bigg)+r^{-1}\Ga_b\Im\Big(\pmb\psi\c\ov{\Ddot^{\nu}\pmb\psi}\Big)
 \eeaa
where we also used Corollary \ref{cor:computationofthecomponentsofthetensorAforexactenergyconservationKerr} and the definition of $\tilde{w}$.

On the other hand, noticing that $\pr_\tau(r)=\pr_\tau(\th)=\pr_\tau(q)=\pr_\tau(\widetilde{w})=0$, we have
\beaa
\bsplit
&\Ddot^\mu\left(k\widetilde{w}\Im\left( \pmb\psi\c\ov{\Ddot_\mu\pmb\psi}\right)+(\pr_{\tt})_{\mu}\frac{2a\cos\th}{\qs}k\widetilde{w}|\pmb\psi|^2\right)\\
=& k\widetilde{w}\Im\left( \pmb\psi\c\ov{\squared_k\pmb\psi}\right)+k\D^\mu(\widetilde{w})\Im\left( \pmb\psi\c\ov{\Ddot_\mu\pmb\psi}\right)+k\Div (\pr_{\tt})\frac{2a\cos\th}{\qs}\widetilde{w}|\pmb\psi|^2 +\frac{4a\cos\th}{\qs}k\widetilde{w}\Re(\pmb\psi\c\ov{\nab_{\pr_{\tt}}\pmb\psi})\\
=& k\widetilde{w}\Im\left( \pmb\psi\c\ov{\left(\squared_k \pmb\psi-\frac{4ia\cos\th}{\qs}\nab_{\pr_{\tt}}\pmb\psi-V\pmb\psi\right)}\right)+k\D^\mu(\widetilde{w})\Im\left( \pmb\psi\c\ov{\Ddot_\mu\pmb\psi}\right)+k\Div (\pr_{\tt})\frac{2a\cos\th}{\qs}\widetilde{w}|\pmb\psi|^2. 
\end{split}
\eeaa
Adding the two above identities yields 
\beaa
\D^\mu{}^{(\pr_{\tt})}\widetilde{\PP}_\mu[\pmb\psi] &=&  \D^\mu{\PP}_\mu[\pmb\psi](\pr_{\tt}, 0) + \D^\mu\left(k\widetilde{w}\Im\left( \pmb\psi\c\ov{\Ddot_\mu\pmb\psi}\right)+(\pr_{\tt})_{\mu}\frac{2a\cos\th}{\qs}k\widetilde{w}|\pmb\psi|^2\right)\\
&=& \Re\bigg(\ov{\big(\nab_{\pr_{\tt}}\pmb\psi-ik\widetilde{w}\pmb\psi\big)}\c \left(\squared_k\pmb\psi-\frac{4ia\cos\th}{\qs}\Ddot_{\pr_{\tt}}\pmb\psi-V\pmb\psi\right)\bigg)\\
&&+\frac 1 2 \QQ[\pmb\psi]  \c{}^{(\pr_{\tt})}\pi +k\Div (\pr_{\tt})\frac{2a\cos\th}{\qs}\widetilde{w}|\pmb\psi|^2 +r^{-1}\Ga_b\Im\Big(\pmb\psi\c\ov{\Ddot^{\nu}\pmb\psi}\Big)
\eeaa
as stated. This concludes the proof of Lemma \ref{cor:modifiedcurrentsforprtandprvphi}.
\end{proof}

We have the following corollary of Lemma \ref{cor:modifiedcurrentsforprtandprvphi}.
\begin{corollary}\lab{cor:energyestimatetensoriallevel:scalarizedversion}
Let $(\MM, \g)$ satisfy the assumptions of Sections \ref{subsect:assumps:perturbednullframe} and \ref{sec:regulartripletinperturbationsofKerr}. Let $\pmb\psi\in\sk_2(\mathbb{C})$, and let $\psi_{ij}$ be the corresponding scalars given by $\psi_{ij}:=\pmb\psi(\Om_i, \Om_j)$, for $i,j=1,2,3$, where $\Om_i$, $i=1,2,3$ is the regular triplet introduced in Section  \ref{sec:regulartripletinperturbationsofKerr}. Under the assumptions of Lemma \ref{cor:modifiedcurrentsforprtandprvphi}, and given ${}^{(\pr_{\tt})}\widetilde{\PP}_\mu[\pmb\psi]$ and $\widetilde{w}$ defined as in Lemma \ref{cor:modifiedcurrentsforprtandprvphi}, we have
\bea
&&\D^\mu{}^{(\pr_{\tt})}\widetilde{\PP}_\mu [\pmb\psi]\nn\\
&=& \Re\left(\Big(\square_{\g}(\psi^{ij})+\widehat{S}(\psi)^{ij}+(\widehat{Q}\psi)^{ij}-V\psi^{ij}\Big)\ov{\Big(\pr_{\tt}(\psi_{ij}) - M_{i\tt}^k\psi_{kj} - M_{j\tt}^k\psi_{ik} -i2\widetilde{w}\psi_{ij}\Big)}\right)\nn\\
&&+\frac 1 2 \QQ[\pmb\psi]  \c {^{(\pr_{\tt})}}\pi+\Div (\pr_{\tt})\frac{4a\cos\th}{\qs}\widetilde{w}|\pmb\psi|^2+r^{-1}\Ga_b\Im\Big(\pmb\psi\c\ov{\Ddot^{\nu}\pmb\psi}\Big).
\eea
\end{corollary}

\begin{proof}
We have, in view of Lemmas \ref{lem:ucdotv:product} and \ref{lemma:formoffirstordertermsinscalarazationtensorialwaveeq}, 
\beaa
&&\Re\bigg( \left(\squared_2\pmb\psi - V\pmb\psi -\frac{4ia\cos\th}{|q|^2}\nab_{\pr_\tt}\pmb\psi\right)\c\ov{\big(\nab_{\pr_\tt}\pmb\psi - i2\widetilde{w}\pmb\psi\big)}\bigg)\\
 &=&\Re\bigg( \left(\squared_2\pmb\psi - V\pmb\psi -\frac{4ia\cos\th}{|q|^2}\nab_{\pr_\tt}\pmb\psi\right)(\Om_i, \Om_j)\ov{\big(\nab_{\pr_\tt}\pmb\psi - i2\widetilde{w}\pmb\psi\big)(\Om^i, \Om^j)}\bigg)\\
&=& \Re\left(\Big(\square_{\g}(\psi_{ij}) -S(\psi)_{ij}-(Q\psi)_{ij} -\frac{4ia\cos\th}{|q|^2}(\nab_{\pr_\tt}\pmb\psi)_{ij} -V\psi_{ij}\Big)\ov{\big((\nab_{\pr_\tt}\pmb\psi)_{ij} - i2\widetilde{w}\psi_{ij}\big)}\right).
\eeaa
Next, using 
\beaa
(\nab_{\pr_\tt}\pmb\psi)_{ij} &=& \pr_{\tt}(\psi_{ij}) - M_{i\tt}^k\psi_{kj} - M_{j\tt}^k\psi_{ik},
\eeaa
we infer, in view of the definition \eqref{eq:definitionwidehatSandwidehatQperturbationsofKerr} for $\widehat{S}$ and $\widehat{Q}$,
\beaa
&&\Re\bigg( \left(\squared_2\pmb\psi - V\pmb\psi -\frac{4ia\cos\th}{|q|^2}\nab_{\pr_\tt}\pmb\psi\right)\c\ov{\big(\nab_{\pr_\tt}\pmb\psi - i2\widetilde{w}\pmb\psi\big)}\bigg)\\
&=& \Re\left(\Big(\square_{\g}(\psi^{ij})-\widehat{S}(\psi)^{ij}-(\widehat{Q}\psi)^{ij}-V\psi^{ij}\Big)\ov{\Big(\pr_{\tt}(\psi_{ij}) - M_{i\tt}^k\psi_{kj} - M_{j\tt}^k\psi_{ik} -i2\widetilde{w}\psi_{ij}\Big)}\right),
\eeaa
which together with Lemma \ref{cor:modifiedcurrentsforprtandprvphi} yields the desired identity.
\end{proof}


\subsection{Control of error terms}
\lab{sec:controlforerrortermsinNRGMorawetz}


The following two lemmas, taken respectively from \cite[Lemma 3.3]{MaSz24} and \cite[Lemma 3.5]{MaSz24}, will allow us to control the error terms arising in the derivation of energy-Morawetz estimates in $\MM(\tau_1,\tau_2)$.  

\begin{lemma}\lab{lemma:basiclemmaforcontrolNLterms:ter}
Let $h\in r^{-1}\dk^{\leq 1}\Ga_b$ be a scalar function and let 
$M^{\a\b}$ be symmetric and satisfy
\beaa
&& M^{rr}\in r\dk^{\leq 1}\Ga_b, \qquad M^{r\tau}\in r\dk^{\leq 1}\Ga_g, \qquad M^{\tau\tau}\in \dk^{\leq 1}\Ga_g,\\
&& M^{rx^a}\in \dk^{\leq 1}\Ga_b, \qquad M^{\tau x^a}\in \dk^{\leq 1}\Ga_g,  \qquad M^{x^ax^b}\in r^{-1}\dk^{\leq 1}\Ga_g,
\eeaa
where $a,b=1,2$. Then, the following estimate holds
\beaa
\int_{\MM(\tau_1, \tau_2)}\Big(\big|M^{\a\b}\pr_\a\psi\pr_\b\psi\big|+h|\psi|^2\Big) &\les& \ep\EM[\psi](\tau_1, \tau_2).
\eeaa
\end{lemma}

\begin{lemma}\lab{lemma:basiclemmaforcontrolNLterms:bis}
Let $M^{\a\b}$ be symmetric and satisfy
\beaa
&& M^{rr}\in r\dk^{\leq 1}\Ga_b, \qquad M^{r\tau}\in r\dk^{\leq 1}\Ga_g, \qquad M^{\tau\tau}\in \dk^{\leq 1}\Ga_g,\\
&& M^{rx^a}\in \dk^{\leq 1}\Ga_b, \qquad M^{\tau x^a}\in \dk^{\leq 1}\Ga_g,  \qquad M^{x^ax^b}\in r^{-1}\dk^{\leq 1}\Ga_g,
\eeaa
where $a,b=1,2$. Then, the following estimate holds
\beaa
\int_{\MM(\tau_1, \tau_2)}\big|M^{\a\b}\pr_\a\pr_\b\psi\big|^2
+\bigg|\int_{\MM(\tau_1, \tau_2)}M^{\a\b}\pr_\a\pr_\b\psi \pr_\tau(\pr^{\leq 1}\psi)\bigg|\\
+\int_{\MM(\tau_1, \tau_2)}r^{-1}\big|M^{\a\b}\pr_\a\pr_\b\psi\big| \big|\dk^{\leq 1}\pr^{\leq 1}\psi\big| &\les& \ep\EM[\pr^{\leq 1}\psi](\tau_1, \tau_2). 
\eeaa

Also, let $N$ be a spacetime vectorfield such that we have 
\beaa
N^r\in r\dk^{\leq 2}\Ga_g, \qquad N^\tau\in \dk^{\leq 2}\Ga_g, \qquad N^{x^a}\in \dk^{\leq 2}\Ga_g.
\eeaa
Then, the following holds
\beaa
\int_{\MM(\tau_1, \tau_2)}\big|N^\a\pr_\a\psi\big|^2+ \int_{\MM(\tau_1, \tau_2)}\big|N^\a\pr_\a\psi\big| \Big|\big(\pr_\tau, \pr_r, r^{-1}\pr_{x^a}, r^{-1}\big)\pr^{\leq 1}\psi\Big| \les \ep\EM[\pr^{\leq 1}\psi](\tau_1, \tau_2). 
\eeaa
\end{lemma}

\begin{remark}
In practice, concerning the quantities estimated in Lemma \ref{lemma:basiclemmaforcontrolNLterms:bis}:
\begin{itemize}
\item $\big(\pr_\tau, \pr_r, r^{-1}\pr_{x^a}, r^{-1}\big)\pr^{\leq 1}\psi$ will be due to  energy-Morawetz multipliers, 

\item  $M^{\a\b}\pr_\a\pr_\b\psi$ and $N^\a\pr_\a\psi$ will come from the RHS of the wave equation, in particular after commutation with {various vectorfields such as} $\pr_\tau$.
\end{itemize}
\end{remark}


\subsection{Horizontal Hodge operators}
\lab{sect:HorizontalHodgeopes}


We introduce the following horizontal Hodge operators.
\begin{definition}
\lab{def:HorizontalHodgeopes}   
We define the following horizontal Hodge type operators
\begin{itemize}
\item  $\DDd_2 $ takes $\sk_2$ into $\sk_1$:
\beaa
(\DDd_2 \xi)_a := \nab^b \xi_{ab}.
\eeaa
\item Recalling \eqref{eq:defintiondivcurlandnabhot}, $\DDs_2 $ takes  $\sk_1$ into $\sk_2$:
\beaa
\DDs_2 \xi := -\frac 1 2 \nab\hot \xi.
\eeaa
\end{itemize}   
\end{definition}

Next, define the higher order weighted horizontal Hodge operators $\dkb^j$, $j\in \mathbb{N}$, as follows
\begin{equation}\lab{eq:definitiondkbpowerj}
\begin{split}
\dkb^j:=&(r\DDs_2\,\, r\DDd_2)^{\frac{j}{2}}, \quad \text{if $j$ is even,}\\
\dkb^{j}:=&r\DDd_2(r\DDs_2\,\, r\DDd_2)^{\frac{j-1}{2}}, \quad \text{if $j$ is odd,}
\end{split}
\end{equation}
and define a set of high-order weighted covariant derivatives, for any $\reg\in\mathbb{N}$,
\bea
\lab{def:highorderweightedderivatives:containhorizontal}
\widecheck{\dk}^{\reg}:= {\left\{\dkb^{\reg_1}(\nab_{\pr_{\tt}})^{\reg_2}(\nab_{4}r)^{\reg_3} ,  \,\,\, \reg_1+\reg_2+\reg_3=\reg\right\}}.
\eea

\begin{proposition}\lab{prop:equivalencerelationonweightednorms}
Let $(\MM, \g)$ satisfy the assumptions of Sections \ref{subsect:assumps:perturbednullframe} and  \ref{subsubsect:assumps:perturbedmetric}. Then, the following estimate holds for any $\pmb\psi\in\sk_2(\mathbb{C})$, any $\reg\leq 15$, and any sphere $S(\tt, r)\subset\MM$, with $r\geq R$ and $R\gg 20m$ sufficiently large,
\bea
\lab{eq:equivalencerelationonweightednorms}
\int_{S(\tt,r)}  |\dk^{\leq \reg}\pmb\psi|^2 \les \int_{S(\tt,r)}  |\widecheck{\dk}^{\leq \reg}\pmb\psi|^2.
\eea
\end{proposition}

\begin{proof}
First, note from the assumptions of Sections \ref{subsect:assumps:perturbednullframe} and  \ref{subsubsect:assumps:perturbedmetric}, and in particular of Lemma \ref{lemma:relationsbetweennullframeandcoordinatesframe}, that we have, for $r\geq R$ and $1\leq\reg\leq 15$,
\beaa
\int_{S(\tt, r)}|(r\nab)^\reg\pmb\psi|^2&\les& \int_{S(\tt, r)}|(\nab_{\pr_{x^1}}, \nab_{\pr_{x^2}})^\reg\pmb\psi|^2 +R^{-1}\int_{S(\tt, r)}|(r\nab)^\reg\pmb\psi|^2\\
&&+\int_{S(\tt, r)}|(\nab_4, \nab_{\pr_\tau})(r\nab)^{\reg-1}\pmb\psi|^2+\int_{S(\tt, r)}|\dk^{\leq\reg-1}\pmb\psi|^2.
\eeaa
Also, denoting by $\DDd_2^{\,\,S}$, $\DDs_2^{\,\,\,\,S}$ and $(\dkb^S)^j$ the Hodge operators for tangential tensors on $S(\tt, r)$ that are the analog of the corresponding horizontal Hodge operators of Definition \ref{def:HorizontalHodgeopes} and \eqref{eq:definitiondkbpowerj}, and using the fact that the spheres $S(\tt, r)$ are close to round spheres for $r\geq R$ with $R$ large enough in view of \eqref{eq:assymptiticpropmetricKerrintaurxacoord:2} and \eqref{eq:controloflinearizedmetriccoefficients}, we have, see for instance Lemma 5.1.27 in \cite{KS:Kerr}, 
\beaa
\int_{S(\tt, r)}|(\nab_{\pr_{x^1}}, \nab_{\pr_{x^2}})^\reg\pmb\psi|^2 &\les& \int_{S(\tt, r)}|(\dkb^S)^{\leq\reg}\pmb\psi|^2
\eeaa
and hence, using again the assumptions of Sections \ref{subsect:assumps:perturbednullframe} and  \ref{subsubsect:assumps:perturbedmetric}, and in particular of Lemma \ref{lemma:relationsbetweennullframeandcoordinatesframe}, we obtain, for $r\geq R$ and $1\leq\reg\leq 15$, 
\beaa
\int_{S(\tt, r)}|(\nab_{\pr_{x^1}}, \nab_{\pr_{x^2}})^\reg\pmb\psi|^2 &\les& \int_{S(\tt, r)}|\dkb^{\reg}\pmb\psi|^2+R^{-1}\int_{S(\tt, r)}|(r\nab)^\reg\pmb\psi|^2\\
&&+\int_{S(\tt, r)}|(\nab_4, \nab_{\pr_\tau})(r\nab)^{\reg-1}\pmb\psi|^2+\int_{S(\tt, r)}|\dk^{\leq\reg-1}\pmb\psi|^2.
\eeaa
Plugging in the above, this implies, for $R$ large enough
\beaa
\int_{S(\tt, r)}|(r\nab)^\reg\pmb\psi|^2&\les& \int_{S(\tt, r)}|\dkb^{\reg}\pmb\psi|^2 +\int_{S(\tt, r)}|(\nab_4, \nab_{\pr_\tau})(r\nab)^{\reg-1}\pmb\psi|^2+\int_{S(\tt, r)}|\dk^{\leq\reg-1}\pmb\psi|^2.
\eeaa
Together with the commutator formula \eqref{commutator-3-a-u-b:cor:bis}, we deduce 
\beaa
\int_{S(\tt, r)}|(r\nab)^\reg\pmb\psi|^2&\les& \int_{S(\tt, r)}|\dkb^{\reg}\pmb\psi|^2 +\int_{S(\tt, r)}|(r\nab)^{\reg-1}(\nab_4, \nab_{\pr_\tau})\pmb\psi|^2+\int_{S(\tt, r)}|\dk^{\leq\reg-1}\pmb\psi|^2.
\eeaa
Then, arguing by iteration on $\reg$, we immediately infer, for all $\reg\leq 15$, 
\beaa
\int_{S(\tt,r)}  |(r\nab)^{\reg}\pmb\psi|^2 &\les& \int_{S(\tt,r)}  |\widecheck{\dk}^{\leq \reg}\pmb\psi|^2,
\eeaa
and then
\beaa
\int_{S(\tt,r)}  |\dk^{\leq\reg}\pmb\psi|^2 &\les& \int_{S(\tt,r)}  |\widecheck{\dk}^{\leq \reg}\pmb\psi|^2
\eeaa
as stated in \eqref{eq:equivalencerelationonweightednorms}. This concludes the proof of Proposition \ref{prop:equivalencerelationonweightednorms}.
\end{proof}


\subsection{Commutators with the D'Alembertian}
\lab{sect:commutatorwithDalembertian}


The following two lemmas provide the structure of commutators between first-order derivatives and the scalar wave operator, respectively for unweighted and weighted derivatives.

\begin{lemma}
\lab{lem:commutatorwithwave:firstorderderis:0}
Let $(\MM, \g)$ satisfy the assumptions of Section \ref{subsubsect:assumps:perturbedmetric}. Then, the  commutator between $\square_{\g}$ and $\pr_{\tau}$ satisfies 
\bea
\label{esti:commutatorBoxgandT:general:0}
\, [ \pr_{\tau}, \square_{\g}]\psi = \pr_{\tau}(\gcheck^{\a\b})\pr_{\a}\pr_{\b}\psi+
\dk^{\leq 2}\Ga_g\c\dk\psi,
\eea
and the commutator between $\square_\g$ and $(\pr_r,r^{-1}\pr_{x^a})$ satisfies 
\bea
\lab{eq:localwavecommutators:withfirstordergoodderis:0}
{[(\pr_r,r^{-1}\pr_{x^a}), \square_{\g}]\psi} = O(r^{-2}) (\pr_{\tt}, \pr_{r}, \pr_{x^a})^{\leq 1}\pr\psi 
+r^{-1}\dk^{\leq 1}\Ga_b\dk\pr\psi+r^{-1}\dk^{\leq 2}\Ga_g\dk \psi.
\eea
\end{lemma}

\begin{proof}
The first estimate \eqref{esti:commutatorBoxgandT:general:0} is proven in \cite[Lemma 3.7]{MaSz24}. For the second estimate \eqref{eq:localwavecommutators:withfirstordergoodderis:0}, we recall, from the very end of the proof of Lemma 3.7 in \cite{MaSz24},
\beaa
[(\pr_r,r^{-1}\pr_{x^a}), \square_{\g}] &=& [(\pr_r,r^{-1}\pr_{x^a}), \square_{\gam}] + r^{-1}\dk(\gcheck^{\mu\nu})\pr_{\mu}\pr_{\nu}+\dk^{\leq 2}\Ga_g\c\dk\\
&=& O(r^{-2}) (\pr_{\tt}, \pr_{r}, \pr_{x^a})^{\leq 1}\pr\psi +r^{-1}\dk^{\leq 1}\Ga_b\dk\pr\psi +\dk^{\leq 2}\Ga_g\c\dk\psi
\eeaa
where in the last step we used the fact that
$$
[(\pr_r,r^{-1}\pr_{x^a}), \square_{\gam}] = O(r^{-2}) (\pr_{\tt}, \pr_{r}, \pr_{x^a})^{\leq 1}\pr\psi.
$$
This proves \eqref{eq:localwavecommutators:withfirstordergoodderis:0}  and hence concludes the proof of Lemma \ref{lem:commutatorwithwave:firstorderderis:0}.
\end{proof}

\begin{lemma}
\lab{lem:commutatorwithwave:firstorderderis}
Let $(\MM, \g)$ satisfy the assumptions of Section \ref{subsubsect:assumps:perturbedmetric}. Then, the commutator between $\square_{\g}$ and $\pr_{\tau}$ satisfies
\bea
\label{esti:commutatorBoxgandT:general}
\, [ \pr_{\tau}, \square_{\g}]\psi = \dk^{\leq 1}(\dk^{\leq 1}\Ga_g\c\dk\psi),
\eea
and the commutator between $\square_\g$ and $(r\pr_r,\pr_{x^a})$ satisfies 
\bea
\lab{eq:localwavecommutators:withfirstordergoodderis}
{[(r\pr_r, \pr_{x^a}), \square_{\g}]\psi} = \big(-\square_{\g}\psi, 0\big)+O(r^{-2})\dk^{\leq 1}\dk\psi + \dk^{\leq 1}(\dk^{\leq 1}\Ga_g\c\dk\psi).
\eea
\end{lemma}

\begin{proof}
Using \eqref{eq:controloflinearizedinversemetriccoefficients}--\eqref{eq:consequenceasymptoticKerrandassumptionsinverselinearizedmetric} 
 and Lemma \ref{lemma:computationofthederiveativeofsrqtg}, we have
 \bea\lab{eq:comparisionbetweenscalarwaveoperatorinKerrandinKerrpert}
\nn\square_\g\psi &=& \frac{1}{\sqrt{|\g|}} \pr_{\mu}\Big(\sqrt{|\g|}\Big) \g^{\mu\nu}\pr_{\nu}\psi
+\pr_{\mu} (\g^{\mu\nu}\pr_{\nu}\psi)\\
\nn&=& \left((N_{det})_{\mu}\g^{\mu\nu}+\frac{1}{\sqrt{|\gam|}} \pr_{\mu}\Big(\sqrt{|\gam|}\Big) \gcheck^{\mu\nu}\right)\pr_{\nu}\psi
+ \pr_{\mu}(\gcheck^{\mu\nu}\pr_{\nu}\psi)+\square_{\gam}\psi\\
&=& \square_{\gam}\psi+\dk^{\leq 1}(\Ga_g\c\dk\psi) 
\eea
and hence
\beaa
\,[(\pr_\tau, r\pr_r, \pr_{x^a}), \square_{\g}]\psi = [(\pr_\tau, r\pr_r, \pr_{x^a}), \square_{\gam}]\psi + \dk^{\leq 1}(\dk^{\leq 1}\Ga_g\c\dk\psi).
\eeaa

We deduce, since $[\pr_\tau, \square_{\gam}]=0$, 
\beaa
[ \pr_{\tau}, \square_{\g}]\psi&=&\dk^{\leq 1}(\dk^{\leq 1}\Ga_g\c\dk\psi),
\eeaa
as stated in \eqref{esti:commutatorBoxgandT:general}, as well as 
\beaa
[(r\pr_r, \pr_{x^a}), \square_{\g}] &=& [(r\pr_r, \pr_{x^a}), \square_{\gam}] + \dk^{\leq 1}(\dk^{\leq 1}\Ga_g\c\dk\psi)\\
&=& \big(-\square_{\g}\psi, 0\big)+O(r^{-2})\dk^{\leq 1}\dk\psi + \dk^{\leq 1}(\dk^{\leq 1}\Ga_g\c\dk\psi)
\eeaa
as stated in \eqref{eq:localwavecommutators:withfirstordergoodderis}, where we used in particular the fact that
\beaa
\square_{\gam}\psi=-2\pr_r\pr_\tau\psi-\frac{2}{r}\pr_\tau\psi + O(r^{-2})\dk^{\leq 1}\dk\psi
\eeaa
which follows from \eqref{eq:assymptiticpropmetricKerrintaurxacoord:volumeform} \eqref{eq:assymptiticpropmetricKerrintaurxacoord:1}. This concludes the proof of Lemma \ref{lem:commutatorwithwave:firstorderderis}.
\end{proof}

Next, we provide the structure of the commutators between the horizontal Hodge operators introduced in Definition \ref{def:HorizontalHodgeopes} and the tensorial wave operators $\squared_{k}$, $k=1,2$.

\begin{lemma}\lab{lemma:commutationofhodgeellipticorder1withsqaured2fdiluhs} 
Let $(\MM, \g)$ satisfy the assumptions of Sections \ref{subsect:assumps:perturbednullframe} and  \ref{subsubsect:assumps:perturbedmetric}. Then, the following commutation formula holds true for $\psi \in \sk_2$
\bea
\lab{commutator:rDDd2withrescaledwave}
&&r\DDd_2(r^2\squared_2\psi) - r^2\squared_1(r\DDd_2\psi) \nn\\
&=& 3r\DDd_2\psi +O(m)\nab_3\dk^{\leq 1}\psi+O(mr^{-1})\dk^{\leq 2}\psi
+ \dk^{\leq 2} (r^2 \Ga_g \c \psi)+r\Ga_b\c r^2\squared_2 \psi,
\eea
and the following commutation formula holds true for $\psi \in \sk_1$
\bea
\lab{commutator:rDDs2withrescaledwave}
&&r\DDs_2 \, (r^2\squared_1\psi) - r^2\squared_2 (r\DDs_2\,\psi) \nn\\
&=&- 3r\DDs_2\,\psi+O(m)\nab_3\dk^{\leq 1}\psi+O(mr^{-1})\dk^{\leq 2}\psi
+\dk^{\leq 2}(r^2 \Ga_g\c\psi)+r\Ga_b\c r^2\squared_1 \psi .
\eea
\end{lemma}

\begin{proof}
In view of Lemma \ref{lemma:expression-wave-operator}, we have for $k=1,2$, 
\bea\lab{eq:simplifiedversionsquaredusefulinlargerregion}
\nn\squared_k \psi &=& -\nab_3\nab_4\psi +\frac{1}{r}\nab_4\psi -\frac{1}{r}\nab_3\psi+\lap_k\psi +O(mr^{-2})\nab_3\psi+O(mr^{-3})\dk^{\leq 1}\psi\\
&&+\dk^{\leq 1}(\Ga_g\c\psi),
\eea
which yields
\beaa
&& r\DDd_2\squared_2\psi - \squared_1r\DDd_2\psi\\
&=& -[r\DDd_2, \nab_3]\nab_4\psi - \nab_3[r\DDd_2, \nab_4]\psi +\frac{1}{r}[r\DDd_2, \nab_4]\psi -\frac{1}{r}[r\DDd_2, \nab_3]\psi +r\big(\DDd_2 \lap_2\psi - \lap_1\DDd_2\big)\psi\\
&&+O(mr^{-2})\nab_3\dk^{\leq 1}\psi+O(mr^{-3})\dk^{\leq 2}\psi+\dk^{\leq 2}(\Ga_g\c\psi).
\eeaa
Now, according to the proof of Lemma 4.7.13 in \cite{GKS22}, we have
\beaa
\DDd_2 \lap_2\psi - \lap_1 \DDd_2\psi &=&  \frac{3}{r^2}\DDd_2\psi +O(mr^{-3})\nab_3\psi+O(mr^{-4})\dk^{\leq 2}\psi
\eeaa
and hence
\beaa
&& r\DDd_2\squared_2\psi - \squared_1r\DDd_2\psi\\
&=& \frac{3}{r^2}r\DDd_2\psi -[r\DDd_2, \nab_3]\nab_4\psi - \nab_3[r\DDd_2, \nab_4]\psi +\frac{1}{r}[r\DDd_2, \nab_4]\psi -\frac{1}{r}[r\DDd_2, \nab_3]\psi \\
&&+O(mr^{-2})\nab_3\dk^{\leq 1}\psi+O(mr^{-3})\dk^{\leq 2}\psi+\dk^{\leq 2}(\Ga_g\c\psi).
\eeaa

Next, using the commutator formulas \eqref{commutator-3-a-u-b:cor:bis}, we have
\beaa
&& -[r\DDd_2, \nab_3]\nab_4\psi - \nab_3[r\DDd_2, \nab_4]\psi +\frac{1}{r}[r\DDd_2, \nab_4]\psi -\frac{1}{r}[r\DDd_2, \nab_3]\psi \\
&=& r\etac\left(\nab_3\nab_4\psi +\frac{1}{r}\nab_3\psi\right) +\nab_3(\xi\c\nab_3\psi)  +O(mr^{-2})\nab_3\dk^{\leq 1}\psi+O(mr^{-3})\dk^{\leq 2}\psi+\dk^{\leq 2}(\Ga_g\c\psi)
\eeaa
which together with \eqref{eq:simplifiedversionsquaredusefulinlargerregion} implies 
\beaa
&& -[r\DDd_2, \nab_3]\nab_4\psi - \nab_3[r\DDd_2, \nab_4]\psi +\frac{1}{r}[r\DDd_2, \nab_4]\psi -\frac{1}{r}[r\DDd_2, \nab_3]\psi \\
&=& -r\etac\squared_k \psi +\nab_3(\xi\c\nab_3\psi)  +O(mr^{-2})\nab_3\dk^{\leq 1}\psi+O(mr^{-3})\dk^{\leq 2}\psi+\dk^{\leq 2}(\Ga_g\c\psi)
\eeaa
and hence, plugging in the above, 
\beaa
r\DDd_2\squared_2\psi - \squared_1r\DDd_2\psi &=& \frac{3}{r^2}r\DDd_2\psi +O(mr^{-2})\nab_3\dk^{\leq 1}\psi+O(mr^{-3})\dk^{\leq 2}\psi\\
&&+ \dk^{\leq 2} (\Ga_g \c \psi)+r\Ga_b\c \squared_2 \psi + \nab_3 (r\xi \c \nab_3 \psi ).
\eeaa
Since we have $\DDd_2(r)=r\Ga_g$ from  Definition \ref{definition.Ga_gGa_b} and $\dk^{\leq 15}\xi=O(\ep r^{-3})$ from \eqref{eq:decaypropertiesofGabGag}, it then follows
\beaa
&&r\DDd_2(r^2\squared_2\psi) - r^2\squared_1(r\DDd_2\psi) \nn\\
&=& 3r\DDd_2\psi +O(m)\nab_3\dk^{\leq 1}\psi+O(mr^{-1})\dk^{\leq 2}\psi
+ \dk^{\leq 2} (r^2 \Ga_g \c \psi)+r\Ga_b\c r^2\squared_2 \psi,
\eeaa
as stated. The proof of the second identity is similar and left to the reader. 
\end{proof}

\begin{lemma}\label{lemma:commutator-nab3-nab4-square} 
Let $(\MM, \g)$ satisfy the assumptions of Sections \ref{subsect:assumps:perturbednullframe} and  \ref{subsubsect:assumps:perturbedmetric}. Then, the following commutation formula holds true for $\psi\in\sk_k$, $k=1,2$:
\bea\label{eq:comm-rnab4-squared2}
[\nab_4 r, r^2\squared_k]\psi &=& -2r\nab_{\pr_{\tt}}\nab_4(r\psi) +O(m)\nab_3\dk^{\leq 1}\psi+O(mr^{-1})\dk^{\leq 2}\psi\nn\\
&&+O(mr^{-1})r^2\squared_k\psi +\dk^{\leq 2}(r^2\Ga_g\c\psi),
\eea
and 
the following commutation formula holds true for $\psi\in\sk_k$, $k=1,2$:
\bea\label{eq:comm-nabT-squared2}
[\nab_{\pr_{\tt}}, r^2\squared_k]\psi = O(mr^{-2})\dk^{\leq 1}\psi +\dk^{\leq 1}(\dk^{\leq 1}(r^2\Ga_g)\c\dk\psi).
\eea
\end{lemma}

\begin{proof}
First, notice that 
\beaa
r^{-1}\nab_3\nab_4(r\psi) &=& r^{-1}\nab_3\big(r\nab_4\psi+e_4(r)\psi\big)\\
&=& \nab_3\nab_4\psi+\frac{e_3(r)}{r}\nab_4\psi+r^{-1}\nab_3\Big(\big(1+O(mr^{-1})+\Ga_g\big)\psi\Big)\\
&=& \nab_3\nab_4\psi -\frac{1}{r}\nab_4\psi+\frac{1}{r}\nab_3\psi+O(mr^{-2})\nab_3\psi+O(mr^{-3})\dk^{\leq 1}\psi +\dk^{\leq 1}(\Ga_g\c\psi)
\eeaa
which together with \eqref{eq:simplifiedversionsquaredusefulinlargerregion} implies
\bea\lab{eq:simplifiedversionsquaredusefulinlargerregion:bis}
\squared_k \psi = -r^{-1}\nab_3\nab_4(r\psi)+\lap_k\psi +O(mr^{-2})\nab_3\psi+O(mr^{-3})\dk^{\leq 1}\psi+\dk^{\leq 1}(\Ga_g\c\psi).
\eea

Next, relying on \eqref{eq:simplifiedversionsquaredusefulinlargerregion:bis}, we have
\beaa
\nab_4(r\squared_k\psi) &=& -\nab_4\nab_3\nab_4(r\psi)+\nab_4(r\Delta_k\psi)+O(mr^{-2})\nab_3\dk^{\leq 1}\psi+O(mr^{-3})\dk^{\leq 2}\psi+\dk^{\leq 2}(\Ga_g\c\psi)
\eeaa
and relying on  \eqref{eq:simplifiedversionsquaredusefulinlargerregion} we have
\beaa
\squared_k\nab_4(r\psi) &=& -\nab_3\nab_4\nab_4(r\psi) +\frac{1}{r}\nab_4\nab_4(r\psi) -\frac{1}{r}\nab_3\nab_4(r\psi) +\lap_k\nab_4(r\psi) \\
&&+O(mr^{-2})\nab_3\dk^{\leq 1}\psi+O(mr^{-3})\dk^{\leq 2}\psi+\dk^{\leq 2}(\Ga_g\c\psi).
\eeaa
We deduce 
\beaa
[\nab_4 r, \squared_k]\psi &=&  -[\nab_4,\nab_3]\nab_4(r\psi) -\frac{1}{r}\nab_4\nab_4(r\psi) +\frac{1}{r}\nab_3\nab_4(r\psi) +[\nab_4r,\Delta_k]\psi\\
&&+O(mr^{-2})\nab_3\dk^{\leq 1}\psi+O(mr^{-3})\dk^{\leq 2}\psi+\dk^{\leq 2}(\Ga_g\c\psi).
\eeaa

Next, using \eqref{eq:simplifiedversionsquaredusefulinlargerregion:bis} and the commutator identities \eqref{commutator-3-a-u-b:cor}, we have
\beaa
&&  -[\nab_4,\nab_3]\nab_4(r\psi) +\frac{1}{r}\nab_3\nab_4(r\psi) +[\nab_4r,\Delta_k]\psi\\
&=& -\squared_k\psi -\Delta_k\psi +O(mr^{-2})\nab_3\dk^{\leq 1}\psi+O(mr^{-3})\dk^{\leq 2}\psi+\dk^{\leq 2}(\Ga_g\c\psi)
\eeaa
so that 
\beaa
[\nab_4 r, \squared_k]\psi = -\frac{1}{r}\nab_4\nab_4(r\psi) -\squared_k\psi -\lap_k\psi +O(mr^{-2})\nab_3\dk^{\leq 1}\psi+O(mr^{-3})\dk^{\leq 2}\psi+\dk^{\leq 2}(\Ga_g\c\psi).
\eeaa
Substituting \eqref{eq:simplifiedversionsquaredusefulinlargerregion:bis} into this formula to rewrite $\lap_k\psi$, we deduce
\beaa
[\nab_4 r, \squared_k]\psi = -\frac{1}{r}(\nab_4+\nab_3)\nab_4(r\psi) -2\squared_k\psi +O(mr^{-2})\nab_3\dk^{\leq 1}\psi+O(mr^{-3})\dk^{\leq 2}\psi+\dk^{\leq 2}(\Ga_g\c\psi).
\eeaa
which, in view of $\nab_4(r)=\frac{\De}{\qs} + \Ga_g$ and 
\beaa
\nab_4+\nab_3 =2\nab_{\pr_{\tt}}+ O(mr^{-1})\nab_3 +O(mr^{-2})\dk  +r\Ga_g\dk,
\eeaa
which follows from \eqref{eq:relationsbetweennullframeandcoordinatesframe2:moreprecise:00}, yields
\beaa
[\nab_4 r, r^2\squared_k]\psi &=& -2r\nab_{\pr_{\tt}}\nab_4(r\psi) +O(m)\nab_3\dk^{\leq 1}\psi+O(mr^{-1})\dk^{\leq 2}\psi\nn\\
&&+O(mr^{-1})r^2\squared_k\psi +\dk^{\leq 2}(r^2\Ga_g\c\psi)
\eeaa
as stated. This proves the first identity \eqref{eq:comm-rnab4-squared2}.

For the second identity, we make use of the approximate Killing vectorfield $T$ introduced in \cite[Definition 4.3.1]{GKS22} which is given by
\beaa
T:=\frac{1}{2}\left(e_4+\frac{\De}{|q|^2}e_3-2a\Re(\Jk)^be_b\right)
\eeaa
and which satisfies, in view of\footnote{The case $k=2$ is in \cite[Proposition 4.3.3]{GKS22} while the case $k=1$ follows in a similar manner.} \cite[Proposition 4.3.3]{GKS22}, 
\beaa
[\nab_T, \squared_k]\psi = O(mr^{-4})\dk^{\leq 1}\psi + \dk^{\leq 1}(\Ga_g\c\dk\psi).
\eeaa
Since $T=\pr_\tau$ in Kerr, and using \eqref{eq:relationsbetweennullframeandcoordinatesframe2:moreprecise}, we have $\pr_\tau = T+r\Ga_g\dk$, and hence, we infer 
\beaa
[\nab_{\pr_{\tt}}, \squared_k]\psi &=& [r\Ga_g\dk, \squared_k]\psi+O(mr^{-4})\dk^{\leq 1}\psi + \dk^{\leq 1}(\Ga_g\c\dk\psi)\\
&=& O(mr^{-4})\dk^{\leq 1}\psi +\dk^{\leq 1}(\dk^{\leq 1}\Ga_g\c\dk\psi).
\eeaa
As $\pr_\tau(r)=0$, we deduce
\beaa
[\nab_{\pr_{\tt}}, r^2\squared_k]\psi = O(mr^{-2})\dk^{\leq 1}\psi +\dk^{\leq 1}(\dk^{\leq 1}(r^2\Ga_g)\c\dk\psi)
\eeaa
as stated. This concludes the proof of Lemma \ref{lemma:commutator-nab3-nab4-square}.
\end{proof}

\begin{lemma}\label{lemma:commutator-nab4-square-redshift} 
Let $(\MM, \g)$ satisfy the assumptions of Section \ref{subsect:assumps:perturbednullframe}. Then, the following commutation formula holds true for $\psi\in\sk_2$ in the redshift region $r\leq r_+(1+2\dred )$:
\bea\label{eq:comm-nab4-squared2-redshiftregion}
\begin{split}
[\nab_4, \squared_2]\psi =& \pr_r\left(\frac{\De}{|q|^2}\right)\nab_3\nab_4\psi+O(|r-r_+|)\big(\squared_2\psi, \nab_3\nab_4^{\leq 1}\psi\big)\\
&+O(1)(\nab_4\nab_4^{\leq 1}\psi, \nab\nab_4^{\leq 1}\psi, \nab_4^{\leq 1}\psi)+\dk^{\leq 2}(\Ga_g\c\psi).
\end{split}
\eea
\end{lemma}

\begin{proof}
In view of Lemma \ref{lemma:expression-wave-operator}, we have in the redshift region $r\leq r_+(1+2\dred )$
\beaa
\squared_2\psi =-\nab_3\nab_4\psi +\lap_2\psi +O(r-r_+)\nab_3\psi +O(1)\big(\nab_4\psi, \nab \psi, \psi\big)+\Ga_g\c\dk^{\leq 1}\psi
\eeaa
which implies
\beaa
[\nab_4, \squared_2]\psi &=& -[\nab_4, \nab_3]\nab_4\psi +[\nab_4, \lap_2]\psi +O(r-r_+)[\nab_3, \nab_4]\psi+O(r-r_+)\nab_3\psi \\
&&+O(1)[\nab_4, \nab]\psi +O(1)\big(\nab_4\psi, \nab \psi, \psi\big)+\dk^{\leq 2}(\Ga_g\c\psi).
\eeaa
Next, we have in view of the commutator identities \eqref{commutator-4-a-u-bc} \eqref{commutator-4-3-u-bc}, in the redshift region $r\leq r_+(1+2\dred )$, 
\beaa
[\nab_4,\nab_a]\psi &=& -\frac 1 2 \trch\nab_a\psi -\frac 1 2 \atrch\dual\nab_a\psi +O(1)\nab_4\psi+O(1)\psi+\Ga_g\c\dk^{\leq 1}\psi
\eeaa
and 
\beaa
[\nab_4, \nab_3]\psi &=& 2\om\nab_3\psi+O(1)\nab\psi+O(1)\psi+\Ga_g\c\dk^{\leq 1}\psi
\eeaa
which yield
\beaa
[\nab_4, \squared_2]\psi &=& -2\om\nab_3\nab_4\psi -\trch\Delta_2\psi - \atrch\in^{ab}\nab_a\nab_b\psi \\
&&+O(r-r_+)\nab_3\nab_4^{\leq 1}\psi+O(1)\big(\nab_4\nab_4^{\leq 1}\psi, \nab\nab_4^{\leq 1}\psi, \nab_4^{\leq 1}\psi\big)+\dk^{\leq 2}(\Ga_g\c\psi)\\
&=& \pr_r\left(\frac{\De}{|q|^2}\right)\nab_3\nab_4\psi +O(r-r_+)\Delta_2\psi +O(r-r_+)\in^{ab}\nab_a\nab_b\psi \\
&&+O(r-r_+)\nab_3\nab_4^{\leq 1}\psi+O(1)\big(\nab_4\nab_4^{\leq 1}\psi, \nab\nab_4^{\leq 1}\psi, \nab_4^{\leq 1}\psi\big)+\dk^{\leq 2}(\Ga_g\c\psi).
\eeaa
 Plugging the above formula for $\squared_2\psi$ and using the commutator identity \eqref{commutator-nab-a-nab-b:cor} then implies 
\beaa
\begin{split}
[\nab_4, \squared_2]\psi =& \pr_r\left(\frac{\De}{|q|^2}\right)\nab_3\nab_4\psi+O(r-r_+)\big(\squared_2\psi, \nab_3\nab_4^{\leq 1}\psi\big)\\
&+O(1)(\nab_4\nab_4^{\leq 1}\psi, \nab\nab_4^{\leq 1}\psi, \nab_4^{\leq 1}\psi)+\dk^{\leq 2}(\Ga_g\c\psi)
\end{split}
\eeaa
 as stated.
\end{proof}


\subsection{Local energy estimate}


We have the following basic local (in time) energy estimates for systems of wave equations.

\begin{lemma}[Local energy estimate]
\lab{lemma:localenergyestimate}
Let $\g$ satisfy the assumptions of Section \ref{subsubsect:assumps:perturbedmetric}. Let $(\psi)_{ij}$, $i,j=1,2,3$,  satisfy the following coupled system of scalar wave equations
\bea
\lab{eq:eqsforlocalenergyestimatelemma:general}
\square_{\g}\psi_{ij}={N_{ij}, \qquad N_{ij}:=D_1r^{-1}\pr_{\tt} \psi_{ij}+O(r^{-2})\dk^{\leq 1}\psi_{kl} +F_{ij},}
\eea
with the constant $D_1\geq 0$.
For any $\tau_0\in\mathbb{R}$, $0\leq\reg\leq 14$, $\de\in (0,1]$ and $q>0$, we have the following future directed local energy estimates
\bsub
\bea
\label{eq:localenergyestimate:future}
\EMF_\de^{(\reg)}[\psi](\tau_0, \tau_0+q) &\les_q & \E^{(\reg)}[\psi](\tau_0)  +\NNtlede^{(\reg)}[\psi, F](\tau_0, \tau_0+q),\\
\EMF_\de^{(\reg)}[\psi](\tau_0, \tau_0+q) &\les_q & \E^{(\reg)}[\psi](\tau_0) +\sum_{i,j}\int_{\MM(\tau_0, \tau_0+q)}r^{1+\de}|\dk^{\leq \reg}F_{ij}|^2,\label{eq:localenergyestimate:future:bis}
\eea
\esub
where for any $\tau'<\tau''$, 
\begin{align}\lab{def:NNtleinlocalenergyestimate}
\NNtlede^{(\reg)}[\psi, F](\tau', \tau'')&:=\sum_{i,j}\bigg(\int_{\MM_{r\geq 10m}(\tau', \tau'')}r^{-1}|\dk^{\leq \reg}F_{ij}| |\dk^{\leq \reg+1}\psi_{ij}| 
+\int_{\MM(\tau', \tau'')}|\dk^{\leq \reg}F_{ij}|^2\bigg)  \nn\\
+& \sup_{\tau'''\in[\tau', \tau'']}\bigg|\sum_{i,j}\int_{\MM_{r\geq 10m}(\tau', \tau''')} \Re\Big(\ov{\dk^{\leq \reg}F_{ij}}\big(1+O(r^{-\de})\big)\pr_{\tt} \dk^{\leq \reg}\psi_{ij}\Big)\bigg| .
\end{align}
\end{lemma}

\begin{proof}
First, we commute \eqref{eq:eqsforlocalenergyestimatelemma:general} by $(\pr_\tau, r\pr_r, \pr_{{x^a}})^{\leq\reg}$ and obtain in view of Lemma \ref{lem:commutatorwithwave:firstorderderis}, for $\reg\leq 14$, 
\beaa
\square_{\g}\dk^{\leq\reg}\psi_{ij}= D_1 \big(r^{-1}\pr_{\tt} \dk^{\leq\reg}\psi_{ij} + O(r^{-1})\pr_{\tt} \dk^{\leq\reg-1}\psi_{ij}
\big)+O(r^{-2})\dk^{\leq \reg+1}\psi_{kl} +\dk^{\leq\reg}F_{ij}.
\eeaa
We apply Proposition \ref{prop-app:stadard-comp-Psi} to $\dk^{\leq\reg}\psi_{ij}\in\sk_0(\mathbb{C})$ with $V=0$, $\N=O(r^{-2})\dk^{\leq \reg+1}\psi_{kl} +\dk^{\leq\reg}F_{ij}$, and we choose $w=0$, and a vector field $X$ that is globally uniformly timelike in $\MM$ and equals\footnote{Note, {in view of \eqref{eq:controloflinearizedmetriccoefficients}}, that $\g(\pr_\tau, \pr_\tau)=(\gam)_{\tau\tau}+r\Ga_b=-(1-\frac{2mr}{|q|^2})+O(\ep)\les -1$ in $r\geq 3m$.} $\pr_{\tau}$ for $r\geq 3m$. By integrating over $\MM(\tau_0, \tau)$, for $\tau\in[\tau_0, \tau_0+q]$, we infer
\bea
\lab{eq:localenergyestimate:zeroorder:general}
\nn&&\EF^{(\reg)}[\psi](\tau_0, \tau) 
+D_1\int_{\MM(\tau_0,\tau)} r^{-1}|\pr_{\tt} \dk^{\leq \reg}\psi|^2\\ &\les& \E^{(\reg)}[\psi](\tau_0) +\bigg| \int_{\MM(\tau_0, \tau)} \Re\Big(X(\dk^{\leq\reg}\psi)\ov{\big(O(r^{-2})\dk^{\leq \reg+1}\psi_{kl} +\dk^{\leq\reg}F_{ij}\big)}\Big) \bigg|\nn\\
&&+\frac{1}{2} \bigg|\int_{\MM(\tau_0, \tau)}{}^{(X)} \pi \cdot \QQ[\dk^{\leq\reg}\psi]\bigg|.
\eea
Next, we estimate the last integral in \eqref{eq:localenergyestimate:zeroorder:general}. Since we have in view of Lemma \ref{lemma:controlofdeformationtensorsforenergyMorawetz}, 
\beaa
&& \big({}^{(\pr_{\tau})} \pi\big)^{rr}\in r\dk^{\leq 1}\Ga_b, \qquad \big({}^{(\pr_{\tau})} \pi\big)^{r\tau}\in r\dk^{\leq 1}\Ga_g, \qquad \big({}^{(\pr_{\tau})} \pi\big)^{\tau\tau}\in \dk^{\leq 1}\Ga_g,\\  
&& \big({}^{(\pr_{\tau})} \pi\big)^{rx^a}\in \dk^{\leq 1}\Ga_b, \qquad \big({}^{(\pr_{\tau})} \pi\big)^{\tau x^a}\in \dk^{\leq 1}\Ga_g,  \qquad \big({}^{(\pr_{\tau})} \pi\big)^{x^ax^b}\in r^{-1}\dk^{\leq 1}\Ga_g,
\eeaa
we deduce, as $X=\pr_\tau$ for $r\geq 3m$, 
\beaa
\bigg|\int_{\MM(\tau_0, \tau)}{}^{(X)} \pi \cdot \QQ[\dk^{\leq\reg}\psi]\bigg| &\les& \int_{\MM_{r_+(1-\dhor), 3m}(\tau_0, \tau)}|\pr^{\leq\reg+1}\psi|^2+\int_{\MM_{3m,+\infty}(\tau_0, \tau)}|\dk^{\leq 1}\Ga_g||\dk^{\leq\reg+1}\psi|^2\\
&\les& \int_{\MM(\tau_0, \tau)}r^{-2}|\dk^{\leq\reg+1}\psi|^2\\
&\les& \int_{\tau_0}^{\tau} \E^{(\reg)}[\psi](\tau') d\tau'
\eeaa
which thus yields
\bea
\lab{eq:localenergyfluxestimate:intheproof}
&&\EF^{(\reg)}[\psi](\tau_0, \tau)
+D_1\int_{\MM(\tau_0,\tau)} r^{-1}|\pr_{\tt} \dk^{\leq \reg}\psi|^2\nn\\
& \les &\E^{(\reg)}[\psi](\tau_0) +\int_{\tau_0}^{\tau} \E^{(\reg)}[\psi](\tau') d\tau'+ \bigg|\int_{\MM(\tau_0, \tau)}\Re\big({X(\dk^{\leq\reg}\psi)\ov{\dk^{\leq\reg}F}}\big)\bigg|\nn\\
& \les &\E^{(\reg)}[\psi](\tau_0) +\int_{\tau_0}^{\tau} \E^{(\reg)}[\psi](\tau') d\tau'+\NNtlede^{(\reg)}[\psi, F](\tau_0, \tau_0+q).
\eea

Next, we control the Morawetz part which can also be estimated in a standard way. We apply Proposition \ref{prop-app:stadard-comp-Psi} with the choice $(X=X_{\de}, w=w_{\de})$ given by \eqref{def:Xdeandwde:improvedMorawetz:generaltensorialwave:Kerrpert}. This yields the  estimate \eqref{eq:Moraestimate:case1D1=0:0order:proof:step1}, and hence, for $\de\in(0,1]$, 
\beaa
\M_{\de}^{(\reg)}[\psi](\tau_0, \tau)& \les &{\EF^{(\reg)}}[\psi](\tau_0,\tau) +D_1\int_{\MM(\tau_0,\tau)} r^{-1}|\pr_{\tt} \dk^{\leq \reg}\psi|^2\\
&& + \bigg|\int_{\MM(\tau_0, \tau)}\Re\Big(\big(X_{\de}(\dk^{\leq\reg}\psi) +w_{\de}\dk^{\leq\reg}\psi\big)\ov{O(r^{-2})\dk^{\leq \reg+1}\psi_{kl}}\Big)\bigg|
 \nn\\
&&+\int_{\tau_0}^{\tau} {\E^{(\reg)}}[\psi](\tau')d\tau'
+ \bigg|\int_{\MM(\tau_0, \tau)}\Re\Big(\big(X_{\de}(\dk^{\leq\reg}\psi) +w_{\de}\dk^{\leq\reg}\psi\big)\ov{\dk^{\leq\reg}F}\Big)\bigg|
\nn\\
& \les & {\EF^{(\reg)}}[\psi](\tau_0,\tau) +D_1\int_{\MM(\tau_0,\tau)} r^{-1}|\pr_{\tt} \dk^{\leq \reg}\psi|^2 +\int_{\tau_0}^{\tau}{\E^{(\reg)}}[\psi]( \tau')d\tau'\\
&&+\NNtlede^{(\reg)}[\psi, F](\tau_0, \tau),
\eeaa
where we used the fact that $X_\de=-\pr_\tau+O(r^{-\de})\pr_\tau+O(r^{-1})\dk$ and $w_\de=O(r^{-1})$ in view of \eqref{def:Xdeandwde:improvedMorawetz:generaltensorialwave:Kerrpert}. Combining this with the local energy-flux estimate \eqref{eq:localenergyfluxestimate:intheproof},  and applying Gr\"onwall's inequality, we obtain \eqref{eq:localenergyestimate:future}. Finally, the estimate
 \eqref{eq:localenergyestimate:future:bis} follows from \eqref{eq:localenergyestimate:future} in view of the following straightforward estimate for $\NNtlede^{(\reg)}[\psi, F](\tau_1,\tau_2)$
 \beaa
 \NNtlede^{(\reg)}[\psi, F](\tau_1,\tau_2) \les \Big(\EMF_{\de}^{(\reg)}[\psi](\tau_1,\tau_2)\Big)^{\frac{1}{2}}
 \bigg( \sum_{i,j}\int_{\MM(\tau_1,\tau_2)}r^{1+\de}|\dk^{\leq \reg}F_{ij}|^2\bigg)^{\frac{1}{2}} +\int_{\MM(\tau_1, \tau_2)}|\dk^{\leq \reg}F|^2.
 \eeaa
This concludes the proof of Lemma \ref{lemma:localenergyestimate}.
\end{proof}

The following lemma is the analog of Lemma \ref{lemma:localenergyestimate} upon replacing weighted derivatives $\dk$ with unweighted derivatives $\pr$.
\begin{lemma}[Local energy estimate with unweighted derivatives]
\lab{lemma:localenergyestimate:unweigthedderivatives}
Let $\g$ satisfy the assumptions of Section \ref{subsubsect:assumps:perturbedmetric}. Let $(\psi)_{ij}$, $i,j=1,2,3$,  satisfy the coupled system of scalar wave equations \eqref{eq:eqsforlocalenergyestimatelemma:general} with the constant $D_1\geq 0$.
For any $\tau_0\in\mathbb{R}$, $0\leq\reg\leq 14$, $\de\in (0,1]$ and $q>0$, we have the following future directed local energy estimates
\bsub
\bea
\label{eq:localenergyestimate:future:unweigthedderivatives}
\EMF_\de[\pr^{\leq\reg}\psi](\tau_0, \tau_0+q) &\les_q & \E[\pr^{\leq\reg}\psi](\tau_0)  +\NNtlede[\pr^{\leq\reg}\psi, \pr^{\leq\reg}F](\tau_0, \tau_0+q),\\
\EMF_\de[\pr^{\leq\reg}\psi](\tau_0, \tau_0+q) &\les_q & \E[\pr^{\leq\reg}\psi](\tau_0) +\sum_{i,j}\int_{\MM(\tau_0, \tau_0+q)}r^{1+\de}|\pr^{\leq\reg}F_{ij}|^2,\label{eq:localenergyestimate:future:unweigthedderivatives:bis}
\eea
\esub
and the following past directed local energy estimates 
\bsub
\begin{align}
\label{eq:localenergyestimate:past:unweigthedderivatives}
\EMF_\de[\pr^{\leq\reg}\psi](\tau_0-q, \tau_0) \les_q & \,\E[\pr^{\leq\reg}\psi](\tau_0) +\F[\pr^{\leq\reg}\psi](\tau_0-q, \tau_0)\nn\\
&+ \NNtlede[\pr^{\leq\reg}\psi, \pr^{\leq\reg}F](\tau_0-q, \tau_0),\\
\EMF_\de[\pr^{\leq\reg}\psi](\tau_0-q, \tau_0) \les_q & \,\E[\pr^{\leq\reg}\psi](\tau_0) +\F[\pr^{\leq\reg}\psi](\tau_0-q, \tau_0)\nn\\
&+ \sum_{i,j}\int_{\MM(\tau_0-q, \tau_0)}r^{1+\de}|\pr^{\leq \reg}F_{ij}|^2,\label{eq:localenergyestimate:past:bis:unweigthedderivatives}
\end{align}
\esub
where $\NNtlede$ has been introduced in \eqref{def:NNtleinlocalenergyestimate}.
\end{lemma}

\begin{proof}
We start with the proof of \eqref{eq:localenergyestimate:future:unweigthedderivatives} \eqref{eq:localenergyestimate:future:unweigthedderivatives:bis}. First, we commute \eqref{eq:eqsforlocalenergyestimatelemma:general} by $(\pr_\tau, \pr_r, r^{-1}\pr_{{x^a}})^{\leq\reg}$ and obtain in view of Lemma \ref{lem:commutatorwithwave:firstorderderis:0}, for $\reg\leq 14$, 
\beaa
\square_{\g}\pr^{\leq\reg}\psi_{ij}= D_1r^{-1}\pr_{\tt}\pr^{\leq\reg}\psi_{ij} +O(r^{-2})\dk^{\leq 1}\pr^{\leq \reg}\psi_{kl} +\pr^{\leq\reg}F_{ij}+\pr^{\leq\reg}(\gcheck^{\a\b})\pr_{\a}\pr_{\b}\pr^{\leq\reg -1}\psi.
\eeaa
Arguing as in the proof of Lemma \ref{lemma:localenergyestimate}, we then infer for $\tau\in[\tau_0, \tau_0+q]$
\bea\lab{eq:intemediaryestimateforproof:eq:localenergyestimate:future:unweigthedderivatives}
&&\EMF_\de[\pr^{\leq\reg}\psi](\tau_0, \tau) +D_1\int_{\MM(\tau_0,\tau)} r^{-1}|\pr_{\tt} \pr^{\leq \reg}\psi|^2\nn\\ 
& \les & \E[\pr^{\leq\reg}\psi](\tau_0) +\int_{\tau_0}^{\tau} \E[\pr^{\leq\reg}\psi](\tau') d\tau' +\NNtlede[\pr^{\leq\reg}\psi, \pr^{\leq\reg}F](\tau_0, \tau_0+q)\nn\\
&&+\left|\int_{\MM(\tau_0, \tau)}\pr^{\leq\reg}(\gcheck^{\a\b})\pr_{\a}\pr_{\b}\pr^{\leq\reg -1}\psi\ov{(\pr, r^{-1})\pr^{\leq\reg}\psi}\right|.
\eea

Next, we estimate the last term on the RHS of \eqref{eq:intemediaryestimateforproof:eq:localenergyestimate:future:unweigthedderivatives}. To this end, we notice in view of \eqref{eq:controloflinearizedinversemetriccoefficients} that we may apply Lemma \ref{lemma:basiclemmaforcontrolNLterms:bis} which yields
\beaa
&& \left|\int_{\MM(\tau_0, \tau)}\pr^{\leq\reg}(\gcheck^{\a\b})\pr_{\a}\pr_{\b}\pr^{\leq\reg -1}\psi\ov{(\pr, r^{-1})\pr^{\leq\reg}\psi}\right|\\
&\les& \left(\int_{\MM(\tau_0, \tau)}|\pr^{\leq\reg}(\gcheck^{\a\b})\pr_{\a}\pr_{\b}\pr^{\leq\reg -1}\psi|^2\right)^{\frac{1}{2}}\left(\int_{\tau_0}^{\tau} \E[\pr^{\leq\reg}\psi](\tau') d\tau'\right)^{\frac{1}{2}}\\
&&+\left|\int_{\MM(\tau_0, \tau)}\pr^{\leq\reg}(\gcheck^{\a\b})\pr_{\a}\pr_{\b}\pr^{\leq\reg -1}\psi\ov{\pr_\tau(\pr^{\leq\reg}\psi)}\right|\\
&\les& \int_{\tau_0}^{\tau} \E[\pr^{\leq\reg}\psi](\tau') d\tau'+\ep\EM[\pr^{\leq\reg}\psi](\tau_0,\tau).
\eeaa
Plugging this estimate in \eqref{eq:intemediaryestimateforproof:eq:localenergyestimate:future:unweigthedderivatives}, we infer for $\ep>0$ small enough 
\beaa
&&\EMF_\de[\pr^{\leq\reg}\psi](\tau_0, \tau) +D_1\int_{\MM(\tau_0,\tau)} r^{-1}|\pr_{\tt} \pr^{\leq \reg}\psi|^2\nn\\ 
& \les & \E[\pr^{\leq\reg}\psi](\tau_0) +\int_{\tau_0}^{\tau} \E[\pr^{\leq\reg}\psi](\tau') d\tau' +\NNtlede[\pr^{\leq\reg}\psi, \pr^{\leq\reg}F](\tau_0, \tau_0+q).
\eeaa
Finally, applying Gr\"onwall's inequality, we obtain \eqref{eq:localenergyestimate:future:unweigthedderivatives}, and relying on 
 \beaa
 \NNtlede[\pr^{\leq\reg}\psi, \pr^{\leq\reg}F](\tau_1,\tau_2)&\les& \Big(\EMF_{\de}[\pr^{\leq\reg}\psi](\tau_1,\tau_2)\Big)^{\frac{1}{2}}
 \bigg( \sum_{i,j}\int_{\MM(\tau_1,\tau_2)}r^{1+\de}|\pr^{\leq \reg}F_{ij}|^2\bigg)^{\frac{1}{2}}\\
 &&+\int_{\MM(\tau_1, \tau_2)}|\pr^{\leq \reg}F|^2,
 \eeaa
we deduce \eqref{eq:localenergyestimate:future:unweigthedderivatives:bis}. The estimates \eqref{eq:localenergyestimate:past:unweigthedderivatives} \eqref{eq:localenergyestimate:past:bis:unweigthedderivatives} follow in the same manner. This concludes the proof of Lemma \ref{lemma:localenergyestimate:unweigthedderivatives}.
\end{proof}


\subsection{Redshift estimates}


This section is devoted to proving redshift estimates, which are useful in removing the degeneracy of the energy in a neighborhood of the event horizon, for wave equations in perturbations of Kerr.

Let us recall from \cite[Lemma 3.12]{MaSz24} the following general redshift estimates for scalar waves.
\begin{lemma}[General redshift estimates for scalar wave]
\lab{lemma:redshiftestimatesscalarwave:general}
Let $\g$ satisfy the assumptions of Section \ref{subsubsect:assumps:perturbedmetric}. Let a scalar function $\psi$ satisfy a wave equation which,  in the redshift region $r\leq r_+(1+2\dred)$, can be written in the  form,
\bea
\lab{eq:Redshift-gen.scalarwaveequation}
\square_\g\psi={-}\left(C_++O\left(\bigg|{\frac{r}{r_+}}-1\bigg|\right)\right)\pr_r\psi +O(1)(\pr_\tau\psi, \pr_{x^a}\psi, \psi)+F{,}
\eea
where $C_+$   is a {function} satisfying
 \bea
 \lab{eq:Redshift-gen.scalarwaveequation-1}
 C_+\ge   0, \qquad |\pr^{\leq \reg} C_+|\les 1.
 \eea
 Then, for any $1\leq\tau_1<\tau_2 <+\infty$, we have, for any $0\leq \reg\leq {14}$,  
\bea
\lab{eq:redshift:general:scalarwave:highorderregularity}
\nn\EMF^{(\reg)}_{r\leq r_+(1+\dred)}[\psi](\tau_1, \tau_2)&\les & \E^{(\reg)}[\psi](\tau_1)+\dred^{-1}\M_{r_+(1+\dred), r_+(1+2\dred)}^{(\reg)}[\psi](\tau_1, \tau_2)\\
&&+\int_{\MM_{r\leq r_+(1+2\dred)}(\tau_1, \tau_2)}|\pr^{\leq \reg}F|^2.
\eea
\end{lemma}

We generalize the above redshift estimates for a scalar wave equation to a system of coupled scalar wave equations in the following way.
\begin{lemma}[Redshift estimates for a system of coupled scalar wave equations]
\lab{lemma:redshiftestimatesscalarwaveeqs:general}
Let $\g$ satisfy the assumptions of Section \ref{subsubsect:assumps:perturbedmetric}. Let $0\leq \reg\leq {14}$. Let scalar functions $\psi_{ij}$, $1\leq i,j\leq 3$, satisfy a system of coupled wave equations which, in the redshift region $r\leq r_+(1+2\dred)$, can be written in the form,
\bea
\lab{eq:Redshift-gen.scalarwaveeqs}
\square_\g\psi_{ij}= -C_{+}\pr_r\psi_{ij}
+\sum_{k,l=1}^3O\left(\bigg|{\frac{r}{r_+}}-1\bigg|\right)\pr_r\psi _{kl}+\sum_{k,l=1}^3O(1)(\pr_\tau\psi_{kl},\pr_{x^a}\psi_{kl},\psi_{kl})+F_{ij},
\eea
where $C_{+}$ is a function satisfying
\bea
 \lab{eq:Redshift-gen.scalarwaveeqs-1}
C_{+}\ge   0, \qquad |\pr^{\leq \reg}C_{+}|\les 1.
\eea
Then, for any $1\leq\tau_1<\tau_2 <+\infty$, we have
\bea
\lab{eq:redshift:general:scalarwaveeqs:highorderregularity}
&&\sum_{i,j=1}^3\EMF^{(\reg)}_{r\leq r_+(1+\dred)}[\psi_{ij}](\tau_1, \tau_2)\\
&\les& \sum_{i,j=1}^3\bigg(\E^{(\reg)}[\psi_{ij}](\tau_1)+\dred^{-1}\M_{r_+(1+\dred), r_+(1+2\dred)}^{(\reg)}[\psi_{ij}](\tau_1, \tau_2)\nn\\
&&+\int_{\MM_{r\leq r_+(1+2\dred)}(\tau_1, \tau_2)}|\pr^{\leq \reg}F_{ij}|^2\bigg).
\eea
\end{lemma}

\begin{proof}
The proof is a direct adaptation from the one of Lemma  \ref{lemma:redshiftestimatesscalarwave:general}.
\end{proof}

Next, we state a redshift estimate  near the event horizon for a class of general tensorial wave equations in perturbations of Kerr.

\begin{lemma}[Redshift estimates for tensorial wave equations]
\lab{lemma:redshiftestimates:general:tensorial:kerrpert}
Let $0\leq \reg\leq {14}$. Let $(\MM, \g)$ satisfy the assumptions of Sections \ref{subsect:assumps:perturbednullframe},  \ref{subsubsect:assumps:perturbedmetric} and \ref{sec:regulartripletinperturbationsofKerr}. Let $\pmb\psi\in\sk_2(\mathbb{C})$ be, in the redshift region $r\leq r_+(1+2\dred)$, a solution to the tensorial wave equation 
\bea
\lab{eq:tensorialwave:general:redshiftsection}
\squared_2\pmb\psi=\left(C_{+}+O\left(\bigg|{\frac{r}{r_+}}-1\bigg|\right)\right)\nab_3\pmb\psi+O(1)\nab_{\pr_\tau}\pmb\psi+O(1)\nab_{\pr_{x^a}}\pmb\psi+O(1)\pmb\psi+\bf{F},
\eea
with $C_+\geq 0$ and $|\pr^{\leq \reg}C_{+}|\les 1$.
Then for any $1\leq \tau_1<\tau_2<+\infty$, we have
\bea
\nn\EMF^{(\reg)}_{r\leq r_+(1+\dred)}[\pmb\psi](\tau_1, \tau_2)&\les& \E^{(\reg)}[\pmb\psi](\tau_1)+\dred^{-1}\M_{r_+(1+\dred), r_+(1+2\dred)}^{(\reg)}[\pmb\psi](\tau_1, \tau_2)\\
&&+\int_{\MM_{r\leq r_+(1+2\dred)}(\tau_1, \tau_2)}|\pr^{\leq \reg}\bf{F}|^2.
\eea
\end{lemma}

\begin{proof}
In view of Lemma \ref{lemma:formoffirstordertermsinscalarazationtensorialwaveeq}, we have
\bea
\squared_2\pmb\psi(\Om_i, \Om_j) &=& \square_\g(\psi_{ij}) -S(\psi)_{ij} - (Q\psi)_{ij}
\eea
with
\bsub
\lab{eq:expansionofSandQterms:neareventhorizon}
\bea
S(\psi)_{ij} &=& 2M_{i}^{k\a}\pr_\a(\psi_{kj}) +2M_{j}^{k\a}\pr_\a(\psi_{ik})\nn\\
&=&(O(|r-r_+|)+O(\ep))e_{3}\psi + O(1)(e_{4}\psi, e_{a}\psi, \psi) \pr\psi, \quad \text{in} \,\,\MM_{r\leq 3m},\\
(Q\psi)_{ij} &=& (\Ddot^\a M_{i\a}^k)\psi_{kj}+(\Ddot^\a M_{j\a}^k)\psi_{ik} -M_{i\a}^kM_k^{l\a}\psi_{lj}-2M_{i\a}^kM_{j}^{l\a}\psi_{kl}-M_{j\a}^kM_k^{l\a}\psi_{il}\nn\\
&=&O(1)\psi_{kl} , \quad \text{in} \,\,\MM_{r\leq 3m},
\eea
\esub
where we have used Lemma \ref{lemma:computationoftheMialphajinKerr} and in particular $(M_{i4}^k)_{K}=O(|r-r_+|)$. Hence, in view of the tensorial wave equation \eqref{eq:tensorialwave:general:redshiftsection} satisfied by $\pmb\psi$ and the fact  that in $\MM_{r\leq 3m}$
\bea
\lab{eq:expansionofealpha:neareventhorizon}
e_{3}=-(1+O(\ep))\pr_{r} + O(1)\sum_{\a\neq r}\pr_{\a}, \quad e_{\a}= O(1)\sum_{\b\neq r}\pr_{\b} +(O(|r-r_+|)+O(\ep))\pr_{r}, \,\, \forall \a\neq 3, 
\eea
the scalars $\psi_{ij}$ then satisfy,  in the redshift region $r\leq r_+(1+2\dred)$, 
\beaa
\square_{\g}\psi_{ij}=-C_+\pr_{r}\psi_{ij}+\sum_{k,l}\Big(O(|r-r_+|)\pr_r\psi_{kl} + O(1)(\pr_{\tt}\psi_{kl}, \pr_{x^a}\psi_{kl}, \psi_{kl})\Big)  +F_{ij}+\sum_{k,l}O(\ep)\pr_r\psi_{kl}.
\eeaa
This system of wave equations for the scalars $\psi_{ij}$ can be put into the form of \eqref{eq:Redshift-gen.scalarwaveeqs} and, applying Lemma \ref{lemma:redshiftestimatesscalarwaveeqs:general}, we have,  for any $1\leq\tau_1<\tau_2 <+\infty$ and $\reg\leq {14}$,  and for $\ep$ suitably small,
\beaa
\nn\EMF^{(\reg)}_{r\leq r_+(1+\dred)}[\pmb\psi](\tau_1, \tau_2)&\les& \E^{(\reg)}[\pmb\psi](\tau_1)+\dred^{-1}\M_{r_+(1+\dred), r_+(1+2\dred)}^{(\reg)}[\pmb\psi](\tau_1, \tau_2)\\
&&+\int_{\MM_{r\leq r_+(1+2\dred)}(\tau_1, \tau_2)}|\pr^{\leq \reg}\bf{F}|^2
\eeaa
as desired. This proves Lemma \ref{lemma:redshiftestimates:general:tensorial:kerrpert}.
\end{proof}

As a corollary, we show a redshift estimate for $\Ab$.
\begin{corollary}
\lab{cor:redshift:Ab:highregularity}
Let $(\MM, \g)$ satisfy the assumptions of Sections \ref{subsect:assumps:perturbednullframe},  \ref{subsubsect:assumps:perturbedmetric} and \ref{sec:regulartripletinperturbationsofKerr}.  We have the following redshift estimate for any $0\leq \reg\leq {14}$ and  $\tau_1<\tau_2$
\bea
\lab{eq:redshift:Ab:highorderregularity}
\nn\sum_{p=0}^2\EMF^{(\reg)}_{r\leq r_+(1+\dred)}[\nab_4^p\Ab](\tau_1, \tau_2) &\les& \sum_{p=0}^2\E_{r\leq r_+(1+2\dred )}^{(\reg)}[\nab_4^p\Ab](\tau_1)\\
\nn&&+\dred^{-5}\sum_{p=0}^2\M_{r_+(1+\dred), r_+(1+2\dred)}^{(\reg)}[\pmb\phi_{-2}^{(p)}](\tau_1, \tau_2)\\
&&+\dred^{-5}\sum_{p=0}^1\int_{\MM_{r\leq r_+(1+2\dred)}(\tau_1, \tau_2)}|\dk^{\leq \reg+1}\N_{T,-2}^{(p)}|^2\nn\\
&&+\sum_{p=0}^2\int_{\MM_{r\leq r_+(1+2\dred)}(\tau_1, \tau_2)}|\dk^{\leq \reg}\N_{\nab_4^p\Ab}|^2.
\eea
\end{corollary}

\begin{proof}
Recall from \eqref{eq:waveequationpmbphip=0sminus2nodeginredshiftregion} that $\nab_4^p\Ab$, $p=0,1,2$, satisfies, in the redshift region $r\leq r_+(1+2\dred )$, 
\beaa
\squared_2\nab_4^p\Ab &=& (2-p)\pr_r\left(\frac{\De}{|q|^2}\right)\nab_3\nab_4^p\Ab
+ O(1)\big(\nab_{4}\nab_4^{\leq p}\Ab,\nab\nab_4^{\leq p}\Ab, \nab_4^{\leq p}\Ab\big)+\N_{\nab_4^p\Ab}\\
&=&  \left((2-p)\pr_r\left(\frac{\De}{|q|^2}\right)+O\left(\bigg|{\frac{r}{r_+}}-1\bigg|\right)\right)\nab_3\nab_4^p\Ab +O(1)\big(\nab_{\pr_\tau}\nab_4^p\Ab, \nab_{\pr_{x^a}}\nab_4^p\Ab, \nab_4^p\Ab\big)\\
&&+\widetilde{\N}_{\nab_4^p\Ab}
\eeaa
where
\beaa
\widetilde{\N}_{\Ab}=O(\ep)\nab_3\Ab+\N_{\Ab}, \qquad \widetilde{\N}_{\nab_4^p\Ab}=O(\ep)\nab_3\nab_4^p\Ab+O(1)\dk^{\leq 1}\nab_4^{\leq p-1}\Ab+\N_{\nab_4^p\Ab}, \quad p=1,2,
\eeaa
which is of the form \eqref{eq:tensorialwave:general:redshiftsection} since $\pr_r(\frac{\De}{|q|^2})>0$ near $r=r_+$. We may thus apply Lemma \ref{lemma:redshiftestimates:general:tensorial:kerrpert} which implies
\beaa
\EMF^{(\reg)}_{r\leq r_+(1+\dred)}[\nab_4^p\Ab](\tau_1, \tau_2) &\les& \E_{r\leq r_+(1+2\dred )}^{(\reg)}[\nab_4^p\Ab](\tau_1)\\
&&+\dred^{-1}\M_{r_+(1+\dred), r_+(1+2\dred)}^{(\reg)}[\nab_4^p\Ab](\tau_1, \tau_2)\\
&&+\int_{\MM_{r\leq r_+(1+2\dred)}(\tau_1, \tau_2)}|\dk^{\leq \reg}\widetilde{\N}_{\nab_4^p\Ab}|^2.
\eeaa
In view of the definition of $\widetilde{\N}_{\nab_4^p\Ab}$, we infer
\beaa
\EMF^{(\reg)}_{r\leq r_+(1+\dred)}[\Ab](\tau_1, \tau_2) &\les& \E_{r\leq r_+(1+2\dred )}^{(\reg)}[\Ab](\tau_1)+\dred^{-1}\M_{r_+(1+\dred), r_+(1+2\dred)}^{(\reg)}[\Ab](\tau_1, \tau_2)\\
&&+\int_{\MM_{r\leq r_+(1+2\dred)}(\tau_1, \tau_2)}|\dk^{\leq \reg}\N_{\Ab}|^2
\eeaa
and
\beaa
\EMF^{(\reg)}_{r\leq r_+(1+\dred)}[\nab_4^p\Ab](\tau_1, \tau_2) &\les& \E_{r\leq r_+(1+2\dred )}^{(\reg)}[\nab_4^p\Ab](\tau_1)+\dred^{-1}\M_{r_+(1+\dred), r_+(1+2\dred)}^{(\reg)}[\nab_4^p\Ab](\tau_1, \tau_2)\\
&&+\M_{r\leq r_+(1+2\dred)}^{(\reg)}[\nab_4^{\leq p-1}\Ab](\tau_1, \tau_2)\\
&&+\int_{\MM_{r\leq r_+(1+2\dred)}(\tau_1, \tau_2)}|\dk^{\leq \reg}\N_{\nab_4^p\Ab}|^2, \quad p=0,1,
\eeaa
which yields
\beaa
\sum_{p=0}^2\EMF^{(\reg)}_{r\leq r_+(1+\dred)}[\nab_4^p\Ab](\tau_1, \tau_2) &\les& \sum_{p=0}^2\E_{r\leq r_+(1+2\dred )}^{(\reg)}[\nab_4^p\Ab](\tau_1)\\
&&+\dred^{-1}\sum_{p=0}^2\M_{r_+(1+\dred), r_+(1+2\dred)}^{(\reg)}[\nab_4^p\Ab](\tau_1, \tau_2)\\
&&+\sum_{p=0}^2\int_{\MM_{r\leq r_+(1+2\dred)}(\tau_1, \tau_2)}|\dk^{\leq \reg}\N_{\nab_4^p\Ab}|^2.
\eeaa
Using the definition of $\pmb\phi_{-2}^{(0)}$ in  \eqref{eq:defintionphipm2p=0:Kerrperturbation} and the transport equations \eqref{def:TensorialTeuScalars:wavesystem:Kerrperturbation}, we infer
\beaa
\sum_{p=0}^2\EMF^{(\reg)}_{r\leq r_+(1+\dred)}[\nab_4^p\Ab](\tau_1, \tau_2) &\les& \sum_{p=0}^2\E_{r\leq r_+(1+2\dred )}^{(\reg)}[\nab_4^p\Ab](\tau_1)\\
&&+\dred^{-5}\sum_{p=0}^2\M_{r_+(1+\dred), r_+(1+2\dred)}^{(\reg)}[\pmb\phi_{-2}^{(p)}](\tau_1, \tau_2)\\
&&+\dred^{-5}\sum_{p=0}^1\int_{\MM_{r\leq r_+(1+2\dred)}(\tau_1, \tau_2)}|\dk^{\leq \reg+1}\N_{T,-2}^{(p)}|^2\\
&&+\sum_{p=0}^2\int_{\MM_{r\leq r_+(1+2\dred)}(\tau_1, \tau_2)}|\dk^{\leq \reg}\N_{\nab_4^p\Ab}|^2
\eeaa
as stated.
\end{proof}


\section{Statement and proof of the main theorem}
\lab{sect:maintheoremandproof}


 In this section, we state a precise version of our main theorem on energy-Morawetz estimates for Teukolsky equations in perturbations of Kerr and provide its proof.


\subsection{Statement of the main theorem}


We now provide a precise version of our main theorem on the derivation of energy-Morawetz estimates for solutions to Teukolsky equations on $(\MM, \g)$, where $\g$ is a perturbation of a Kerr metric $\gam$ with $|a|<m$. 

\begin{theorem}[Energy-Morawetz for Teukolsky equations, precise version] 
\label{thm:main}
Let $(\MM, \g)$ satisfy the assumptions of Sections \ref{subsect:assumps:perturbednullframe}, \ref{subsubsect:assumps:perturbedmetric} and \ref{sec:regulartripletinperturbationsofKerr}. Then, there exists a suitably small constant $\ep'>0$ such that for any $\ep\leq \ep'$, we have for solutions $\pmb\phi_s^{(p)}$, $s=\pm 2$, $p=0,1,2$, to the tensorial Teukolsky wave/transport systems \eqref{eq:TensorialTeuSysandlinearterms:rescaleRHScontaine2:general:Kerrperturbation}  \eqref{def:TensorialTeuScalars:wavesystem:Kerrperturbation}  on spacetimes $(\MM,\g)$ perturbations of Kerr the following energy-Morawetz-flux estimates, for any $1\leq\tau_1<\tau_2 <+\infty$ and any $0<\de\leq \frac{1}{3}$, 
\bea
\label{MainEnerMora:psi:plus2case}
\sum_{p=0}^2\EMF^{({11})}_{\de}[\pmb\phi_{+2}^{(p)}](\tau_1, \tau_2)
&\les&  \sum_{p=0}^2\E^{({11})}[\pmb\phi_{+2}^{(p)}](\tau_1)
+\sum_{p=0}^2\NN^{({11})}_\de[\pmb\phi_{+2}^{(p)}, \N_{W,+2}^{(p)}](\tt_1, \tt_2)\nn\\ 
&&
+\sum_{p=0}^1\int_{\MM(\tt_1, \tt_2)}r^{-1+\de}|\dk^{\leq 12}\N_{T,+2}^{(p)}|^2
\eea
and, assuming also that $\Ab$ satisfies \eqref{eq:waveequationpmbphip=0sminus2nodeginredshiftregion},
\bea
\label{MainEnerMora:psi:minus2case}
\nn&&\sum_{p=0}^2\EMF^{({14})}_\de[\pmb\phi_{-2}^{(p)}](\tau_1, \tau_2) + \sum_{p=0}^2\EMF^{({14})}_{r\leq r_+(1+\dred)}[\nab_4^p\Ab](\tau_1, \tau_2)\\
\nn&\les& \sum_{p=0}^2\E^{({14})}[\pmb\phi_{-2}^{(p)}](\tau_1)+\sum_{p=0}^2\E^{({14})}_{r\leq r_+(1+\dred)}[\nab_4^p\Ab](\tau_1) +\int_{\Si(\tau_1)}|\dk^{\leq 12}\N_{T,-2}^{(0)}|^2\\
\nn&&+\sum_{p=0}^1\int_{\MM(\tt_1, \tt_2)}r^{-1+\de}|\dk^{\leq {15}}\N_{T,-2}^{(p)}|^2+\sum_{p=0}^2\NN_\de^{({14})}[\pmb\phi_{-2}^{(p)}, \N_{W,-2}^{(p)}](\tt_1, \tt_2)\\ 
&&+\sum_{p=0}^2\int_{\MM_{r\leq r_+(1+2\dred)}(\tau_1, \tau_2)}|\dk^{\leq {14}}\N_{\nab_4^p\Ab}|^2,
\eea
where the norms $\EMF^{(\reg)}_\de[\c](\tau_1,\tau_2)$, $\E^{(\reg)}[\c](\tau_1)$ and $\mathcal{N}^{(\reg)}_\de[\c, \c](\tau_1, \tau_2)$ have been introduced in Section \ref{subsect:norms}, where $\N_{W,s}^{(p)}$, $\N_{T,s}^{(p)}$ and $\N_{\nab_4^p\Ab}$ are introduced in equations \eqref{eq:TensorialTeuSysandlinearterms:rescaleRHScontaine2:general:Kerrperturbation},   \eqref{def:TensorialTeuScalars:wavesystem:Kerrperturbation} and \eqref{eq:waveequationpmbphip=0sminus2nodeginredshiftregion}, where $\dk$ is given in \eqref{eq:defweightedderivative}, and where the implicit constant in $\lesssim$ only depends on $a$, $m$, $\de$, $\Nmic$ and $\dec$ (with $\dec$ appearing in \eqref{eq:decaypropertiesofGabGag}).
\end{theorem}

{\begin{remark}
\lab{remark:ofmainthm7.1}
Here are some comments on the statement and proof of Theorem \ref{thm:main}:
\begin{itemize}
\item The energy-Morawetz estimates proven in Theorem \ref{thm:main} require only the standard energy bound of the initial data for the Teukolsky wave system. In terms of the required $r$-decay, our assumptions on the fall-off of the initial data are both optimal and weaker than all corresponding results \cite{Ma} \cite{DHR19} \cite{SRTdC20, SRTdC23} \cite{Millet} in Kerr spacetimes.  

\item Unlike the proofs in \cite{Ma} \cite{DHR19} \cite{SRTdC20, SRTdC23} in Kerr, and Part II of \cite{GKS22} in perturbations of Kerr for $|a|\ll m$, we do not rely on transport estimates for the Teukolsky transport equations \eqref{def:TensorialTeuScalars:wavesystem:Kerrperturbation}. Instead, we use these equations only as identities in order to rewrite the Teukolsky wave equations \eqref{eq:TensorialTeuSysandlinearterms:rescaleRHScontaine2:general:Kerrperturbation}. 

\item The EMF norms for $\{\phis{p}\}_{s=\pm 2, p=0,1}$ on the LHS of  \eqref{MainEnerMora:psi:plus2case} \eqref{MainEnerMora:psi:minus2case} can, as shown in the proof of Theorem \ref{thm:main}, be replaced by EMF norms that are non-degenerate in $\Mtrap$.

\item The regularity requirements $\reg=11$ in \eqref{MainEnerMora:psi:plus2case} and $\reg=14$ in \eqref{MainEnerMora:psi:minus2case} are due to the weak Morawetz estimates  for inhomogeneous Teukolsky equations on a subextremal Kerr background stated in Theorem \ref{cor:weakMorawetzforTeukolskyfromMillet:bis}, whose proof relies on Millet's result \cite{Millet}. Since the rest of the proof of Theorem \ref{thm:main} holds for all $1\leq\reg\leq 14$, any regularity improvement of the statement of Theorem \ref{cor:weakMorawetzforTeukolskyfromMillet:bis} would lead to a corresponding improvement of our regularity requirements in Theorem \ref{thm:main}.
\end{itemize}
\end{remark}}


\subsection{Main intermediary results}
\lab{sec:statementofthemainintermediarystepsinproofofmainth}



\subsubsection{Extension to a semi-global problem}
\lab{sect:extensiontoglobalproblem:Teu}


We start with the following definition.
\begin{definition}\lab{def:definitionofthetimetauR}
Let $\Nmic$ be the large integer introduced in  Section \ref{sec:smallnesconstants}. We define $\tmic<\tau_1$ as the smallest real number such that
\bea
\Si(\tmic)\cap\{r\leq (\Nmic+1)m\}\subset D^-\big(\Si(\tau_1)\cap\{r\leq 2\Nmic m\}\big),
\eea
i.e., $\Si(\tmic)\cap\{r\leq (\Nmic+1)m\}$ is included in the past domain of dependence of $\Si(\tau_1)\cap\{r\leq 2\Nmic m\}$. We also define the interval $\Iti:=(\tmic, +\infty)$.
\end{definition}

\begin{remark}
In view of the choice of the coordinate $\tau$ in Lemma \ref{lem:specificchoice:normalizedcoord}, particularly in the region $r\geq 13m$, and the assumptions for $\g$ in Section \ref{subsect:assumps:perturbedmetric}, the choice of $\tmic$ in Definition \ref{def:definitionofthetimetauR} satisfies  
\bea
-\frac{m}{\Nmic}\les \tmic -\tau_1\les -\frac{m}{\Nmic}.
\eea
\end{remark}

Recall from Lemma \ref{lem:scalarizedTeukolskywavetransportsysteminKerrperturbation:Omi} that the tensorial Teukolsky wave/transport systems \eqref{eq:TensorialTeuSysandlinearterms:rescaleRHScontaine2:general:Kerrperturbation}  \eqref{def:TensorialTeuScalars:wavesystem:Kerrperturbation} are equivalent to the scalarized Teukolsky wave/transport systems  \eqref{eq:ScalarizedTeuSys:general:Kerrperturbation} \eqref{eq:ScalarizedQuantitiesinTeuSystem:Kerrperturbation}. The proof of Theorem \ref{thm:main} will follow from global energy-Morawetz estimates for an extension to $\tau\in\Iti$, with $\Iti$ as in Definition \ref{def:definitionofthetimetauR}, of the scalarized system of wave equations \eqref{eq:ScalarizedTeuSys:general:Kerrperturbation}. The goal of this section is to construct this extended solution.

\begin{proposition}\lab{prop:extensionprocedureoftheTeukolskywaveequations}
Let $(\MM, \g)$ satisfy the assumptions of Sections \ref{subsect:assumps:perturbednullframe}, \ref{subsubsect:assumps:perturbedmetric} and \ref{sec:regulartripletinperturbationsofKerr}. Assume that $\phi^{(p)}_{s,ij}$ satisfies the scalarized system of wave equations \eqref{eq:ScalarizedTeuSys:general:Kerrperturbation} for $\tau\in(\tau_1, \tau_2)$ with $\tau_2$ satisfying\footnote{Since the proof of Theorem \ref{thm:main} in the case $\tau_1<\tau_2\leq\tau_1+10$ follows immediately from local energy type arguments, see Step 0 in Section \ref{subsect:proofofThm4.1}, we focus here on the case where $\tt_2$ satisfies  \eqref{eq:tau2geqtau1plus10conditionisthemaincase}.} 
\bea\lab{eq:tau2geqtau1plus10conditionisthemaincase}
\tau_2\geq\tau_1+10.
\eea
Also, let $\tmic<\tau_1$ and $\Iti$ be as in Definition \ref{def:definitionofthetimetauR}, let $\chi_{\tau_1, \tau_2}=\chi_{\tau_1, \tau_2}(\tau)$ be a smooth cut-off function satisfying 
\bea\lab{eq:propertieschitoextendmetricg}
\chi_{\tau_1, \tau_2}(\tau)=0\,\,\,\textrm{on}\,\,\, \mathbb{R}\setminus(\tau_1, \tau_2-1), \,\,\,\,\chi_{\tau_1, \tau_2}(\tau)=1\,\,\,\textrm{on}\,\,\, [\tau_1+1, \tau_2-2], \,\,\,\, \|\chi_{\tau_1, \tau_2}\|_{W^{{15},+\infty}(\mathbb{R})}\les 1,
\eea
let $\chi^{(1)}_{\tau_1, \tau_2}=\chi^{(1)}_{\tau_1, \tau_2}(\tau)$ be a smooth cut-off function satisfying 
\bea\lab{eq:propertieschi:thisonetodealwithRHSofTeukolsky}
\chi^{(1)}_{\tau_1, \tau_2}(\tau)=0\,\,\,\textrm{on}\,\,\, \mathbb{R}\setminus(\tau_1, \tau_2-2), \,\,\,\,\chi^{(1)}_{\tau_1, \tau_2}(\tau)=1\,\,\,\textrm{on}\,\,\, [\tau_1+1, \tau_2-3], \,\,\,\, \|\chi^{(1)}_{\tau_1, \tau_2}\|_{W^{{15},+\infty}(\mathbb{R})}\les 1,
\eea
and define the extended metric
\bea\lab{eq:extendedmetricgchitau1tau2}
\g_{\chi_{\tau_1, \tau_2}}^{\a\b}=\chi_{\tau_1, \tau_2} \g^{\a\b}+(1-\chi_{\tau_1, \tau_2})\g_{a,m}^{\a\b}.
\eea
Then, there exists $\psi^{(p)}_{s,ij}$ satisfying the following system of scalar wave equations defined by
\bea\lab{eq:waveeqwidetildepsi1}
\big({\square}_{\g_{\chi_{\tau_1, \tau_2}}}- (4-2\de_{p0})|q|^{-2}\big)\psi^{(p)}_{s,ij} &=F_{total,s,ij}^{(p)} \quad\textrm{on}\quad\MM(\Iti),
\eea
with\footnote{Below, $\widehat{S}_K$ and $\widehat{Q}_K$ denote $\widehat{S}$ and $\widehat{Q}$ of \eqref{eq:definitionwidehatSandwidehatQperturbationsofKerr} computed using the regular triplet in Kerr of Section  \ref{sec:regulartripletsinKerr}.}
\bea\lab{def:tildef}
\bsplit
F_{total,s,ij}^{(p)}:=&\widehat{F}^{(p)}_{s,ij} +\widetilde{F}^{(p)}_{s,ij}+\underline{F}^{(p)}_{s,ij}+\breve{F}^{(p)}_{s,ij},\\
\widehat{F}^{(p)}_{s,ij}:=&\chi_{\tau_1, \tau_2}\big( \widehat{S}(\pmb\psi^{(p)}_s)_{ij} +(\widehat{Q}\pmb\psi^{(p)}_s)_{ij}\big) +(1-\chi_{\tau_1, \tau_2})\big(\widehat{S}_K(\pmb\psi^{(p)}_s)_{ij} +(\widehat{Q}_K\pmb\psi^{(p)}_s)_{ij} +f_p\psi^{(p)}_{s,ij}\big),\\
f_p:=& f_p(r, \cos\th)=\frac{2\de_{p0}}{|q|^2}- \frac{4a^2\cos^2\th(|q|^2+6mr)}{|q|^6},\\
\widetilde{F}^{(p)}_{s,ij}:=&\chi^{(1)}_{\tau_1, \tau_2}F^{(p)}_{s,ij}, \qquad F^{(p)}_{s,ij}=L^{(p)}_{s,ij}+N^{(p)}_{W,s,ij},  \quad i,j=1,2,3,\,\,\, p=0,1,2, \,\,\, s=\pm 2,
\end{split}
\eea
${L_{s,ij}^{(p)}}=(\L_{s}^{(p)}[\pmb\phi_s])_{ij}$ given by \eqref{eq:linearterms:ScalarizedTeuSys:general:Kerrperturbation}, 
and $\underline{F}^{(p)}_{s,ij}$ and $\breve{F}^{(p)}_{s,ij}$ given respectively by \eqref{def:tildef0} and \eqref{def:breveFpsij},
such that the following properties hold:
\begin{itemize}
\item $\g_{\chi_{\tau_1, \tau_2}}$ satisfies the assumptions of Section \ref{subsubsect:assumps:perturbedmetric} and coincides with Kerr in $\MM\setminus(\tau_1, \tau_2-1)$, 

\item $\psi^{(p)}_{s,ij}$ satisfies the following identities 
\bea\lab{eq:causlityrelationsforwidetildepsi1}
\bsplit
\psi^{(p)}_{s,ij}&=\phi^{(p)}_{s,ij}\quad\textrm{on}\quad \MM(\tau_1+1, \tau_2-3), \\
\dk^{\leq\reg}(\psi^{(p)}_{s,ij}) &= 0\quad\textrm{on}\quad \Si(\tmic)\cap 
\{r\leq(\Nmic+1)m\}\,\,\,\,\forall\,\reg\in\mathbb{N},
\end{split}
\eea
so that $\psi^{(p)}_{s,ij}$ can be smoothly extended\footnote{This will allow us to derive microlocal energy-Morawetz in Section  \ref{sect:microlocalenergyMorawetztensorialwaveequation} on $\MM_{r\leq\Rmic}$ where  $\Rmic\in [\Nmic m, (\Nmic+1)m]$ is introduced in Remark \ref{rmk:choiceofconstantRbymeanvalue}.} by $0$ in $\MM(-\infty, \tmic)\cap\{r\leq (\Nmic+1)m\}$.

\item $\psi^{(p)}_{s,ij}$ satisfies the following local energy estimate on $\tau\in[\tmic, \tau_1+3]$, for all $\reg\leq {14}$, 
\bea\lab{eq:localenergyforpsiontau1minus1tau1plus1}
&&\sum_{p=0}^2\sum_{p=0}^2\EMF_\de[\pr^{\leq \reg}\psi^{(p)}_{s}](\tmic, \tau_1+3)\nn\\
&\les& \sum_{i,j=1}^3\sum_{p=0}^2\bigg(\E[\pr^{\leq \reg}\phi^{(p)}_{s,ij}](\tau_1)+\int_{\MM(\tau_1, \tau_1+3)}r^{1+\de} |\pr^{\leq \reg}N_{W,s,ij}^{(p)}|^2\bigg).
\eea

\item $\psi^{(p)}_{s,ij}$ satisfies the following local energy estimates on $\tau\in[\tau_2-3, \tau_2]$, for all $\reg\leq {14}$, 
\bea\lab{eq:localenergyestimateforpsisijponSigmatau2minus3tau2:partialderivatives}
\nn&&\sum_{i,j=1}^3\sum_{p=0}^2\EMF_\de[\pr^{\leq \reg}\psi^{(p)}_{s,ij}](\tau_2-3, \tau_2)\\ 
&\les& \sum_{i,j=1}^3\sum_{p=0}^2\bigg(\E[\pr^{\leq \reg}\phi^{(p)}_{s,ij}](\tau_2-3)+\int_{\MM(\tau_2-3, \tau_2-2)}r^{1+\de} |\pr^{\leq \reg}N_{W,s,ij}^{(p)}|^2\bigg).
\eea

\item The tensorization defect $\err_{\textrm{TDefect}}[\psi^{(p)}_{s}]$ corresponding to the family of complex-valued scalars $\psi_{ij}$ as introduced in  Definition \ref{def:definitionofthenotationerrforthescalarizationdefect} satisfies
\bea\lab{eq:propertyofscalarizationdefect:prop1}
\err_{\textrm{TDefect}}[\psi^{(p)}_{s}]=0\quad\textrm{on}\quad\MM(\tau_1+1, +\infty)\setminus(\tau_2-2, \tau_2),
\eea
and, for all $\reg\leq {14}$,
\bea\lab{eq:propertyofscalarizationdefect:prop2}
\nn&& \sum_{p=0}^2\EMF_\de[\pr^{\leq \reg}\err_{\textrm{TDefect}}[\psi^{(p)}_{s}]](\tau_2-2, \tau_2)\\
&\les& \ep^2\sum_{i,j=1}^3\sum_{p=0}^2\bigg(\E[\pr^{\leq \reg}\phi^{(p)}_{s,ij}](\tau_2-3)+\int_{\MM(\tau_2-3, \tau_2-2)}r^{1+\de} |\pr^{\leq \reg}N_{W,s,ij}^{(p)}|^2\bigg).
\eea
\end{itemize}
\end{proposition}

\begin{proof}
We proceed along the following steps. 

\noindent{\bf Step 1.} Since $\g_{\chi_{\tau_1, \tau_2}}$ is defined by \eqref{eq:extendedmetricgchitau1tau2}, we have
\beaa
\widecheck{\g}_{\chi_{\tau_1, \tau_2}}^{\a\b}=\chi_{\tau_1, \tau_2} \g^{\a\b}+(1-\chi_{\tau_1, \tau_2})\gam^{\a\b} -\gam^{\a\b}=\chi_{\tau_1, \tau_2}\left(\g^{\a\b}-\gam^{\a\b}\right)=\chi_{\tau_1, \tau_2}\widecheck{\g}^{\a\b}
\eeaa
and hence
\beaa
|\dk^{\leq {15}}\widecheck{\g}_{\chi_{\tau_1, \tau_2}}^{\a\b}|\les |\dk^{\leq {15}}\chi_{\tau_1, \tau_2}||\dk^{\leq {15}}\widecheck{\g}^{\a\b}|\les \|\chi_{\tau_1, \tau_2}\|_{W^{{15},+\infty}}|\dk^{\leq {15}}\widecheck{\g}^{\a\b}| \les |\dk^{\leq {15}}\widecheck{\g}^{\a\b}|,
\eeaa
where we used the fact that $|\dk^{\leq {15}}\chi_{\tau_1, \tau_2}|\les |(\pr_\tau, r\pr_r, \pr_{x^a})^{\leq {15}}\chi_{\tau_1, \tau_2}|=|\pr_\tau^{\leq {15}}\chi_{\tau_1, \tau_2}|$ since $\chi_{\tau_1, \tau_2}=\chi_{\tau_1, \tau_2}(\tau)$ and in view of the definition of the weighted derivatives $\dk$. Since $\g$ satisfies the assumptions of Section \ref{subsubsect:assumps:perturbedmetric}, and in view of the properties \eqref{eq:propertieschitoextendmetricg} of $\chi_{\tau_1, \tau_2}$, we deduce that $\g_{\chi_{\tau_1, \tau_2}}$ also satisfies the assumptions of Section \ref{subsubsect:assumps:perturbedmetric}, and in addition coincides with $\gam$ in $\mathbb{R}\setminus(\tau_1, \tau_2)$. 

\noindent{\bf Step 2.}  Next, we introduce the solutions $\widetilde{\phi}^{(p)}_{s,ij}$ to the following auxiliary system of scalar wave equations, for $p=0,1,2$, $s=\pm 2$, and $i,j=1,2,3$,
\bea\lab{eq:waveequationdefiningfirstextensionwidetildepsi}
\bsplit
\big({\square}_{\g_{\chi_{\tau_1, \tau_2}}}-(4-2\de_{p0})|q|^{-2}\big)\widetilde{\phi}^{(p)}_{s,ij} =& \chi_{\tau_1, \tau_2}\big( \widehat{S}(\widetilde{\pmb\phi}^{(p)}_s)_{ij} +(\widehat{Q}\widetilde{\pmb\phi}^{(p)}_s)_{ij}\big)\\
& +(1-\chi_{\tau_1, \tau_2})\big( \widehat{S}_K(\widetilde{\pmb\phi}^{(p)}_s)_{ij} +(\widehat{Q}_K\widetilde{\pmb\phi}^{(p)}_s)_{ij}+f_p\widetilde{\phi}^{(p)}_{s,ij}\big)\\
&+\widetilde{F}^{(p)}_{s,ij}\,\,\,\,\textrm{on}\,\,\,\MM(\tau_1, \tau_2), \\
 \widetilde{\phi}^{(p)}_{s,ij} =& \phi^{(p)}_{s,ij}, \quad N_{\Sigma(\tau_1+1)}\widetilde{\phi}^{(p)}_{s,ij}=N_{\Sigma(\tau_1+1)}\phi^{(p)}_{s,ij}\,\,\,\,\textrm{on}\,\,\,\Si(\tau_1+1),\\
\widetilde{\phi}^{(p)}_{s,ij} =& \phi^{(p)}_{s,ij}\quad\textrm{on}\quad(\AA_+\cup\II_+)\cap\{\tau_1\leq\tau\leq\tau_1+1\}.
\end{split}
\eea 
Then, we have in view of the local energy estimate \eqref{eq:localenergyestimate:past:bis:unweigthedderivatives}  with $\reg\leq {14}$ for \eqref{eq:waveequationdefiningfirstextensionwidetildepsi}
\beaa
&&\sum_{i,j=1}^3\sum_{p=0}^2\EMF_\de[\pr^{\leq \reg}\widetilde{\phi}^{(p)}_{s,ij}](\tau_1, \tau_1+1)\nn\\
&\les& \sum_{i,j=1}^3\sum_{p=0}^2\bigg(\EF[\pr^{\leq \reg}\phi^{(p)}_{s,ij}](\tau_1, \tau_1+1)+\int_{\MM(\tau_1, \tau_1+1)}r^{1+\de} |\pr^{\leq \reg}N_{W,s,ij}^{(p)}|^2\bigg).
\eeaa
Together with the local energy estimate \eqref{eq:localenergyestimate:future:unweigthedderivatives:bis} with $\reg\leq {14}$ for $\phi^{(p)}_{s,ij}$, we infer
\bea\lab{eq:localenergyestimateforphisijponSigmatau1tau1plus1}
\nn&&\sum_{i,j=1}^3\sum_{p=0}^2\EMF_\de[\pr^{\leq \reg}\phi^{(p)}_{s,ij}](\tau_1, \tau_1+1)\\ 
&\les& \sum_{i,j=1}^3\sum_{p=0}^2\bigg(\E[\pr^{\leq \reg}\phi^{(p)}_{s,ij}](\tau_1)+\int_{\MM(\tau_1, \tau_1+1)}r^{1+\de} |\pr^{\leq \reg}N_{W,s,ij}^{(p)}|^2\bigg)
\eea
and
\bea\lab{eq:localenergyestimateforwidetildepsionSigma}
\nn&&\sum_{i,j=1}^3\sum_{p=0}^2\EMF_\de[\pr^{\leq \reg}\widetilde{\phi}^{(p)}_{s,ij}](\tau_1, \tau_1+1)\\ 
&\les& \sum_{i,j=1}^3\sum_{p=0}^2\bigg(\E[\pr^{\leq \reg}\phi^{(p)}_{s,ij}](\tau_1)+\int_{\MM(\tau_1, \tau_1+1)}r^{1+\de} |\pr^{\leq \reg}N_{W,s,ij}^{(p)}|^2\bigg).
\eea

Also, let $\chi_{\Nmic}(r)$ be a smooth cut-off function such that $\chi_{\Nmic}(r)=1$ for $r\leq 2\Nmic m$ and $\chi_{\Nmic}(r)=0$ for $r\geq 4\Nmic m$. Then, we introduce the solution $(\phi_{aux})^{(p)}_{s,ij}$ to the following auxiliary system of wave equations
\bea\lab{eq:waveequationdefiningpsiaux}
\bsplit
\big({\square}_{\g_{\chi_{\tau_1, \tau_2}}}-(4-2\de_{p0})|q|^{-2}\big)(\phi_{aux})^{(p)}_{s,ij} =&\widehat{S}_K((\pmb\phi_{aux})^{(p)}_s)_{ij} +(\widehat{Q}_K(\pmb\phi_{aux})^{(p)}_s)_{ij}\\
& +f_p(\phi_{aux})^{(p)}_{s,ij}\,\,\,\textrm{on}\,\,\,\MM(\tmic, \tau_1),\\
(\phi_{aux})^{(p)}_{ij}=&\chi_{\Nmic}(r)\widetilde{\phi}^{(p)}_{s,ij}\quad\textrm{on}\quad\Si(\tau_1),\\
N_{\Sigma(\tau_1)}(\phi_{aux})^{(p)}_{s,ij}=&\chi_{\Nmic}(r)N_{\Sigma(\tau_1)}(\widetilde{\phi}^{(p)}_{s,ij})\quad\textrm{on}\quad\Si(\tau_1),\\
(\phi_{aux})^{(p)}_{s,ij}=& (\phi_\AA)^{(p)}_{s,ij}\quad\textrm{on}\quad\AA_+\cap\{\tmic\leq\tau\leq\tau_1\},\\ 
(\phi_{aux})^{(p)}_{s,ij}=&0\quad\textrm{on}\quad\II_+\cap\{\tmic\leq\tau\leq\tau_1\},
\end{split}
\eea 
where $(\phi_\AA)^{(p)}_{s,ij}$ is a smooth extension of $\phi^{(p)}_{s,ij}$ from $\AA\cap\{\tau\geq \tau_1\}$ to $\AA\cap\{\tmic\leq\tau\leq\tau_1\}$ satisfying, for $\reg\leq {14}$,
\bea\lab{eq:propertyextensionpsiAAofpsitotauleqtau1}
\F_\AA[\pr^{\leq\reg}(\phi_\AA)^{(p)}_{s,ij}](\tmic, \tau_1) \les \F_\AA[\pr^{\leq\reg}\phi^{(p)}_{s,ij}](\tau_1, \tau_1+1).
\eea 
The local energy estimate \eqref{eq:localenergyestimate:past:bis:unweigthedderivatives} for  \eqref{eq:waveequationdefiningpsiaux} yields, for $\reg\leq {14}$, 
\beaa
&&\sum_{i,j=1}^3\sum_{p=0}^2\EMF_\de[\pr^{\leq \reg}(\phi_{aux})^{(p)}_{s,ij}](\tmic, \tau_1)\\ 
&\les& \sum_{i,j=1}^3\sum_{p=0}^2\Big(\E[\pr^{\leq \reg}(\phi_{aux})^{(p)}_{s,ij}](\tau_1)+\F[\pr^{\leq \reg}(\phi_{aux})^{(p)}_{s,ij}](\tmic, \tau_1)\Big)\\
&\les& \sum_{i,j=1}^3\sum_{p=0}^2\Big(\E[\pr^{\leq \reg}\widetilde{\phi}^{(p)}_{s,ij}](\tau_1)+\F_\AA[\pr^{\leq \reg}(\phi_\AA)^{(p)}_{s,ij}](\tmic, \tau_1)\Big)
\eeaa
which together with  \eqref{eq:localenergyestimateforwidetildepsionSigma} and \eqref{eq:propertyextensionpsiAAofpsitotauleqtau1} implies
\beaa
&&\sum_{i,j=1}^3\sum_{p=0}^2\EMF_\de[\pr^{\leq \reg}(\phi_{aux})^{(p)}_{s,ij}](\tmic, \tau_1)\\
&\les& \sum_{i,j=1}^3\sum_{p=0}^2\bigg(\EF[\pr^{\leq \reg}\phi^{(p)}_{s,ij}](\tau_1, \tau_1+1)+\int_{\MM(\tau_1, \tau_1+1)}r^{1+\de} |\pr^{\leq \reg}N_{W,s,ij}^{(p)}|^2\bigg).
\eeaa
Using \eqref{eq:localenergyestimateforphisijponSigmatau1tau1plus1}, we deduce, for $\reg\leq {14}$,
\bea\lab{eq:localenergyestimateforpsiauxonSigma}
\nn&&\sum_{i,j=1}^3\sum_{p=0}^2\EMF_\de[\pr^{\leq \reg}(\phi_{aux})^{(p)}_{s,ij}](\tmic, \tau_1)\\ 
&\les& \sum_{i,j=1}^3\sum_{p=0}^2\bigg(\E[\pr^{\leq \reg}\phi^{(p)}_{s,ij}](\tau_1)+\int_{\MM(\tau_1, \tau_1+1)}r^{1+\de} |\pr^{\leq \reg}N_{W,s,ij}^{(p)}|^2\bigg).
\eea

\noindent{\bf Step 3.} Next, we define
\bea
\label{def:tildef0}
\underline{F}^{(p)}_{s,ij} = \left\{\begin{array}{l}
\Big({\square}_{\g_{\chi_{\tau_1, \tau_2}}}-(4-2\de_{p0})|q|^{-2} -\big( \widehat{S}_K+\widehat{Q}_K +f_p\big)\Big)(\chi_{\tau_1}(\phi_{aux})^{(p)}_{s,ij})\,\,\,\textrm{on}\,\,\,\MM(\tmic, \tau_1),\\[2mm]
0\quad\textrm{on}\quad \MM\setminus\MM(\tmic, \tau_1),
\end{array}\right.
\eea
where the smooth cut-off $\chi_{\tau_1}=\chi_{\tau_1}(\tau)$ is such that $\chi_{\tau_1}=1$ for $\tau\geq \tau_1$ and 
$\chi_{\tau_1}=0$ for $\tau\leq \tmic$. In particular, \eqref{eq:waveequationdefiningpsiaux} and \eqref{def:tildef0} imply that, for all $\tau\in\mathbb{R}$, 
\beaa
\underline{F}^{(p)}_{s,ij} &=&\big({\square}_{\g_{\chi_{\tau_1, \tau_2}}}-(4-2\de_{p0})|q|^{-2} -\big( \widehat{S}_K +\widehat{Q}_K+f_p\big)\big)(\chi_{\tau_1}(\phi_{aux})^{(p)}_{s,ij}) \\
&=& 2\g_{\chi_{\tau_1, \tau_2}}^{\a\b}\pr_\a(\chi_{\tau_1})\pr_\b((\phi_{aux})^{(p)}_{s,ij})+\square_{\g_{\chi_{\tau_1, \tau_2}}}(\chi_{\tau_1})(\phi_{aux})^{(p)}_{s,ij} -[\widehat{S}_K, \chi_{\tau_1}](\phi_{aux})^{(p)}_{s,ij}.
\eeaa
Now, since $\g_{\chi_{\tau_1, \tau_2}}$ satisfies the assumptions of Section \ref{subsubsect:assumps:perturbedmetric} in view of Step 2, we easily infer the following non-sharp consequence of \eqref{eq:assymptiticpropmetricKerrintaurxacoord:1}, \eqref{eq:assymptiticpropmetricKerrintaurxacoord:volumeform},  \eqref{eq:controloflinearizedinversemetriccoefficients} and Lemma \ref{lemma:computationofthederiveativeofsrqtg} 
\bea\lab{eq:structureoftildef0}
\underline{F}^{(p)}_{s,ij} = O(r^{-1})\Big(\chi_{\tau_1}''(\tau), \chi_{\tau_1}'(\tau)\Big)\dk^{\leq 1}(\phi_{aux})^{(p)}_{s,ij}.
\eea
Also, notice from \eqref{eq:waveequationdefiningpsiaux} and finite speed of propagation that $(\phi_{aux})^{(p)}_{s,ij}$ vanishes in the past domain of dependence of $\Si(\tau_1)\cap\{r\geq 4\Nmic m\}$ which clearly includes $\{r\geq 5\Nmic m\}\cap\{\tmic\leq\tau\leq\tau_1\}$, and hence 
\bea\lab{eq:structureoftildef0:support}
\textrm{Supp}\left(\underline{F}^{(p)}_{s,ij}\right)\subset\{\tmic\leq\tau\leq\tau_1\}\cap\{r\leq 5\Nmic m\}.
\eea
By using the control of the energy of $(\phi_{aux})^{(p)}_{s,ij}$ provided by  \eqref{eq:localenergyestimateforpsiauxonSigma} together with \eqref{eq:structureoftildef0} \eqref{eq:structureoftildef0:support}, we obtain
\bea\lab{eq:localenergyestimateforunderlineFsijponSigma}
\nn&& \sum_{i,j=1}^3\sum_{p=0}^2\int_{\MM(\tmic, \tau_1)}r^{1+\de}|\pr^{\leq\reg}\underline{F}_{s,ij}^{(p)}|^2\\ 
&\les& \sum_{i,j=1}^3\sum_{p=0}^2\bigg(\E[\pr^{\leq \reg}\phi^{(p)}_{s,ij}](\tau_1)+\int_{\MM(\tau_1, \tau_1+1)}r^{1+\de}{|\pr^{\leq \reg}N_{W,s,ij}^{(p)}|^2}\bigg).
\eea

\noindent{\bf Step 4.} Next, we introduce the solutions $\breve{\phi}^{(p)}_{s,ij}$ to the following auxiliary system of scalar wave equations, for $p=0,1,2$, $s=\pm 2$, and $i,j=1,2,3$,
\bea\lab{eq:waveequationdefiningfirstextensionbrevepsi}
\bsplit
\big({\square}_{\gam}-(4-2\de_{p0})|q|^{-2}\big)\breve{\phi}^{(p)}_{s,ij} =&  \widehat{S}_K(\breve{\phi}^{(p)}_s)_{ij} +(\widehat{Q}_K\breve{\phi}^{(p)}_s)_{ij}+f_p\breve{\phi}^{(p)}_{s,ij},\,\,\,\textrm{on}\,\,\,\MM(\tau_2-1, +\infty), \\
 \breve{\phi}^{(p)}_{s,ij} =& (\Pi_2[\widetilde{\phi}_{s}^{(p)}])_{ij}\quad\textrm{on}\quad\Si(\tau_2-1), \\ 
 N_{\Sigma(\tau_2-1)}\breve{\phi}^{(p)}_{s,ij}=&N_{\Sigma(\tau_2-1)}(\Pi_2[\widetilde{\phi}^{(p)}_{s}])_{ij}\quad\textrm{on}\quad\Si(\tau_2-1),
\end{split}
\eea
where $\widetilde{\phi}_{s,ij}^{(p)}$ is the solution of \eqref{eq:waveequationdefiningfirstextensionwidetildepsi}, and where $\Pi_2$ has been  introduced in \eqref{eq:computationerrorscalarizationdeffect:defPi2}. Then, since $\err_{\textrm{TDefect}}[\Pi_2[\breve{\phi}^{(p)}_{s}]]=0$ in view of Lemma \ref{lemma:computationerrorscalarizationdeffect}, we infer from \eqref{eq:waveequationdefiningfirstextensionbrevepsi} that 
\beaa
\err_{\textrm{TDefect}}[\breve{\phi}^{(p)}_{s}] = 0, \qquad  N_{\Sigma(\tau_2-1)}\err_{\textrm{TDefect}}[\breve{\phi}^{(p)}_{s}]= 0\quad\textrm{on}\quad\Si(\tau_2-1),
\eeaa
which together with uniqueness for the system of wave equations for $\err_{\textrm{TDefect}}[\breve{\phi}^{(p)}_{s}]$ of Lemma \ref{lemma:waveequationsfortensordeffects} implies 
\bea\lab{eq:thetensordeffectofbrevephivanishesidentically}
\err_{\textrm{TDefect}}[\breve{\phi}^{(p)}_{s}] = 0\quad\textrm{on}\quad\MM(\tau_2-1, +\infty).
\eea

Also, we define $\breve{F}^{(p)}_{s,ij}$ as follows 
\bea
\label{def:breveFpsij}
\breve{F}^{(p)}_{s,ij} = \left\{\begin{array}{l}
\Big({\square}_{\gam}-(4-2\de_{p0})|q|^{-2} \\
\qquad\qquad\quad -\big( \widehat{S}_K+\widehat{Q}_K +f_p\big)\Big)(\chi_{\tau_2}\widetilde{\phi}^{(p)}_{s,ij}+(1-\chi_{\tau_2})\breve{\phi}^{(p)}_{s,ij})\,\,\,\textrm{on}\,\,\,\MM(\tau_2-1, \tau_2),\\[2mm]
0\quad\textrm{on}\quad \MM\setminus\MM(\tau_2-1, \tau_2),
\end{array}\right.
\eea
where the smooth cut-off $\chi_{\tau_2}=\chi_{\tau_2}(\tau)$ is such that $\chi_{\tau_2}=1$ for $\tau\leq \tau_2-2/3$ and 
$\chi_{\tau_2}=0$ for $\tau\geq \tau_2-1/3$. In particular, \eqref{eq:waveequationdefiningfirstextensionwidetildepsi}, \eqref{eq:waveequationdefiningfirstextensionbrevepsi} and \eqref{def:breveFpsij} imply that, for all $\tau\in\mathbb{R}$, 
\beaa
\bsplit
\breve{F}^{(p)}_{s,ij} =& 2\gam^{\a\b}\pr_\a(\chi_{\tau_2})\pr_\b((\widetilde{\phi}-\breve{\phi})^{(p)}_{s,ij})+\square_{\gam}(\chi_{\tau_2})(\widetilde{\phi}-\breve{\phi})^{(p)}_{s,ij} - [\widehat{S}_K, \chi_{\tau_2}](\widetilde{\phi}-\breve{\phi})^{(p)}_{s,ij}
\\
=&  2\gam^{\tau\tau}\chi_{\tau_2}'(\tau)\pr_\tau((\widetilde{\phi}-\breve{\phi})^{(p)}_{s,ij})+2\gam^{\tau r}\chi_{\tau_2}'(\tau)\pr_r((\widetilde{\phi}-\breve{\phi})^{(p)}_{s,ij})+2\gam^{\tau x^a}\chi_{\tau_2}'(\tau)\pr_{x^a}((\widetilde{\phi}-\breve{\phi})^{(p)}_{s,ij})\\
& +\left(\gam^{\tau\tau}\chi_{\tau_2}''(\tau)+\frac{1}{\sqrt{|\gam|}}\pr_\a(\sqrt{|\gam|}\gam^{\a\tau})\chi_{\tau_2}'(\tau)\right)(\widetilde{\phi}-\breve{\phi})^{(p)}_{s,ij} - [\widehat{S}_K, \chi_{\tau_2}](\widetilde{\phi}-\breve{\phi})^{(p)}_{s,ij}.
\end{split}
\eeaa
We infer from  \eqref{eq:assymptiticpropmetricKerrintaurxacoord:1}, \eqref{eq:assymptiticpropmetricKerrintaurxacoord:volumeform}  and Lemma  \ref{lemma:computationoftheMialphajinKerr} that 
\bea\lab{eq:structureofbreveFpsij}
\breve{F}^{(p)}_{s,ij} =  -2\chi_{\tau_2}'(\tau)r^{-1}\pr_r(r(\widetilde{\phi}-\breve{\phi})^{(p)}_{s,ij}) +\sum_{k,l}O(r^{-2})\Big(\chi_{\tau_2}''(\tau), \chi_{\tau_2}'(\tau)\Big)\dk^{\leq 1}(\widetilde{\phi}-\breve{\phi})^{(p)}_{s,kl}.
\eea

\noindent{\bf Step 5.} Next, we introduce the solution $\psi^{(p)}_{s,ij}$ to the following scalar wave equation 
\bea\lab{eq:definitionofthescalarizedsystemwaveequatationpsisijp}
\bsplit
\big({\square}_{\g_{\chi_{\tau_1, \tau_2}}}-(4-2\de_{p0})|q|^{-2}\big)\psi^{(p)}_{s,ij} &= \widehat{F}^{(p)}_{s,ij}+\widetilde{F}^{(p)}_{s,ij}+\underline{F}^{(p)}_{s,ij} +\breve{F}^{(p)}_{s,ij}\quad\textrm{on}\quad\MM(\Iti),\\
\psi^{(p)}_{s,ij} &= \widetilde{\phi}^{(p)}_{s,ij}, \quad N_{\Sigma(\tau_1)}\psi^{(p)}_{s,ij}=N_{\Sigma(\tau_1)}\widetilde{\phi}^{(p)}_{s,ij}\quad\textrm{on}\quad\Si(\tau_1),\\
\psi^{(p)}_{s,ij} &= \chi_{\tau_1}(\phi_\AA)^{(p)}_{s,ij}\quad\textrm{on}\quad\AA\cap\{\tau\leq\tau_1\},\\
\psi^{(p)}_{s,ij} &= \chi_{\tau_1}(\phi_\II)^{(p)}_{s,ij}\quad\textrm{on}\quad\II_+\cap\{\tau\leq\tau_1\},
\end{split}
\eea
where we recall that $(\phi_\AA)^{(p)}_{s,ij}$ is a smooth extension of $\phi^{(p)}_{s,ij}$  from $\AA\cap\{\tau\geq \tau_1\}$ to $\AA\cap\{\tmic\leq\tau\leq\tau_1\}$ satisfying \eqref{eq:propertyextensionpsiAAofpsitotauleqtau1}, and where $(\phi_\II)^{(p)}_{s,ij}$ is a smooth extension from $\II\cap\{\tau\geq \tau_1\}$ to $\II\cap\{\tmic\leq\tau\leq\tau_1\}$ satisfying, for $\reg\leq 14$,
\bea\lab{eq:propertyextensionpsiIIofpsitotauleqtau1}
\F_\II[\pr^{\leq \reg}(\phi_\II)^{(p)}_{s,ij}](\tmic, \tau_1) \les \F_\II[\pr^{\leq \reg}\phi^{(p)}_{s,ij}](\tau_1, \tau_1+1).
\eea
In particular, note by causality, using in particular Definition \ref{def:definitionofthetimetauR} of $\tmic$, that we have
\bsub
\lab{eq:psispij:intermsoftildephibrevephiandphiaux}
\bea
\psi^{(p)}_{s,ij} &=& \chi_{\tau_2}\widetilde{\phi}^{(p)}_{s,ij}+(1-\chi_{\tau_2})\breve{\phi}^{(p)}_{s,ij}\quad\textrm{on}\quad\MM(\tau_1, +\infty), \\ 
\psi^{(p)}_{s,ij} &=& \chi_{\tau_1}(\phi_{aux})^{(p)}_{s,ij}\quad\textrm{on}\quad\MM(\tmic, \tau_1)\cap\{r\leq (\Nmic+1)m\}.
\eea
\esub
On the other hand, we have 
\beaa
\widetilde{\phi}^{(p)}_{s,ij}=\phi^{(p)}_{s,ij}\quad\textrm{on}\quad \MM(\tau_1+1, \tau_2-3)
\eeaa
by causality in view of \eqref{eq:waveequationdefiningfirstextensionwidetildepsi}, and we thus deduce
\beaa
\bsplit
\psi^{(p)}_{s,ij}&=\phi^{(p)}_{s,ij}\quad\textrm{on}\quad \MM(\tau_1+1, \tau_2-3), \\\dk^{\leq\reg}(\psi^{(p)}_{s,ij}) &= 0\quad\textrm{on}\quad \Si(\tmic)\cap 
\{r\leq(\Nmic+1)m\}\,\,\,\,\forall\,\reg\in\mathbb{N}.
\end{split}
\eeaa

\noindent{\bf Step 6.} We now derive local energy estimates for $\psi^{(p)}_{s,ij}$ on  $\tau\in[\tmic, \tau_1+3]$ using the system of scalar wave equations in \eqref{eq:definitionofthescalarizedsystemwaveequatationpsisijp}. Using \eqref{eq:localenergyestimate:future:unweigthedderivatives:bis} for the wave system consisting the wave equations for $\psiss{ij}{p}$ in \eqref{eq:definitionofthescalarizedsystemwaveequatationpsisijp} and the wave equations of $\phiss{ij}{p}$ in $\MM(\tau_1,\tau_1+3)$, using \eqref{eq:localenergyestimate:past:bis:unweigthedderivatives} for the wave equations for $\psiss{ij}{p}$ in \eqref{eq:definitionofthescalarizedsystemwaveequatationpsisijp} in $\MM(\tmic,\tau_1)$,  and using the initialization of $\psi^{(p)}_{s,ij}$ on $\Si(\tau_1)$, see \eqref{eq:definitionofthescalarizedsystemwaveequatationpsisijp}, we infer, for all $\reg\leq {14}$, 
\beaa
&&\sum_{i,j=1}^3\sum_{p=0}^2\EMF_\de[\pr^{\leq \reg}\psi^{(p)}_{s,ij}](\tmic, \tau_1+3)+\sum_{i,j=1}^3\sum_{p=0}^2\EMF_\de[\pr^{\leq \reg}\phi^{(p)}_{s,ij}](\tt_1, \tau_1+3)\\ 
&\les& \sum_{i,j=1}^3\sum_{p=0}^2\Big(\E[\pr^{\leq \reg}\widetilde{\phi}^{(p)}_{s,ij}](\tau_1)
+\E[\pr^{\leq \reg}{\phi}^{(p)}_{s,ij}](\tau_1)+\F_\AA[\pr^{\leq \reg}(\phi_\AA)^{(p)}_{s,ij}](\tmic, \tau_1)\\
&&\qquad\qquad\quad+\F_\II[\pr^{\leq \reg}(\phi_\II)^{(p)}_{s,ij}](\tmic, \tau_1)\Big)\\
&&+\sum_{i,j=1}^3\sum_{p=0}^2\int_{\MM(\tau_1, \tau_1+3)}r^{1+\de}|\pr^{\leq\reg}N_{W,s,ij}^{(p)}|^2+\sum_{i,j=1}^3\sum_{p=0}^2\int_{\MM(\tmic, \tau_1)}r^{1+\de}|\pr^{\leq\reg}\underline{F}_{s,ij}^{(p)}|^2,
\eeaa
which together with  \eqref{eq:localenergyestimateforwidetildepsionSigma},  \eqref{eq:propertyextensionpsiAAofpsitotauleqtau1}, \eqref{eq:propertyextensionpsiIIofpsitotauleqtau1} and \eqref{eq:localenergyestimateforunderlineFsijponSigma}  yields, for all $\reg\leq {14}$, 
\beaa
&&\sum_{p=0}^2\sum_{p=0}^2\EMF_\de[\pr^{\leq \reg}\psi^{(p)}_{s}](\tmic, \tau_1+3)\nn\\
&\les& \sum_{i,j=1}^3\sum_{p=0}^2\bigg(\E[\pr^{\leq \reg}\phi^{(p)}_{s,ij}](\tau_1)+\int_{\MM(\tau_1, \tau_1+3)}r^{1+\de} |\pr^{\leq \reg}N_{W,s,ij}^{(p)}|^2\bigg).
\eeaa

\noindent{\bf Step 7.} Next, we derive local energy estimates for $\psi^{(p)}_{s,ij}$ for $\tau\in[\tau_2-3, \tau_2]$. To this end, we first derive local energy estimates for $\widetilde{\phi}^{(p)}_{s,ij}$. Applying the local energy estimate \eqref{eq:localenergyestimate:future:unweigthedderivatives:bis} with $\reg\leq {14}$ to \eqref{eq:waveequationdefiningfirstextensionwidetildepsi}, we have
\beaa
\nn&&\sum_{i,j=1}^3\sum_{p=0}^2\EMF_\de[\pr^{\leq \reg}\widetilde{\phi}^{(p)}_{s,ij}](\tau_2-3, \tau_2)\\ 
&\les& \sum_{i,j=1}^3\sum_{p=0}^2\bigg(\E[\pr^{\leq \reg}\widetilde{\phi}^{(p)}_{s,ij}](\tau_2-3)+\E[\pr^{\leq \reg}\phi^{(p)}_{s,ij}](\tau_2-3)+\int_{\MM(\tau_2-3, \tau_2-2)}r^{1+\de} |\pr^{\leq \reg}N_{W,s,ij}^{(p)}|^2\bigg)
\eeaa
which, together with the fact that $\widetilde{\phi}^{(p)}_{s,ij}=\phi^{(p)}_{s,ij}$ on $\MM(\tau_1+1, \tau_2-3)$  in view of Step 5, yields, for $\reg\leq {14}$,
\bea\lab{eq:localenergyestimateforwidetildephisijponSigmatau2minus3tau2}
\nn&&\sum_{i,j=1}^3\sum_{p=0}^2\EMF_\de[\pr^{\leq \reg}\widetilde{\phi}^{(p)}_{s,ij}](\tau_2-3, \tau_2)\\ 
&\les& \sum_{i,j=1}^3\sum_{p=0}^2\bigg(\E[\pr^{\leq \reg}\phi^{(p)}_{s,ij}](\tau_2-3)+\int_{\MM(\tau_2-3, \tau_2-2)}r^{1+\de} |\pr^{\leq \reg}N_{W,s,ij}^{(p)}|^2\bigg).
\eea

Next, applying the local energy estimate \eqref{eq:localenergyestimate:future:unweigthedderivatives:bis} with $\reg\leq {14}$ to $\breve{\phi}^{(p)}_{s,ij}$ solution of \eqref{eq:waveequationdefiningfirstextensionbrevepsi}, we have
\beaa
\sum_{i,j=1}^3\sum_{p=0}^2\EMF_\de[\pr^{\leq \reg}\breve{\phi}^{(p)}_{s,ij}](\tau_2-1, \tau_2) &\les& \sum_{i,j=1}^3\sum_{p=0}^2\E[\pr^{\leq \reg}(\Pi_2[\widetilde{\phi}_{s}^{(p)}])](\tau_2-1)\\
&\les& \sum_{i,j=1}^3\sum_{p=0}^2\E[\pr^{\leq \reg}\widetilde{\phi}_{s}^{(p)}](\tau_2-1)
\eeaa
which together with \eqref{eq:localenergyestimateforwidetildephisijponSigmatau2minus3tau2} yields, for $\reg\leq {14}$,
\bea\lab{eq:localenergyestimateforbrevephisijponSigmatau2minus1tau2}
\nn&&\sum_{i,j=1}^3\sum_{p=0}^2\EMF_\de[\pr^{\leq \reg}\breve{\phi}^{(p)}_{s,ij}](\tau_2-1, \tau_2)\\ 
&\les& \sum_{i,j=1}^3\sum_{p=0}^2\bigg(\E[\pr^{\leq \reg}\phi^{(p)}_{s,ij}](\tau_2-3)+\int_{\MM(\tau_2-3, \tau_2-2)}r^{1+\de} |\pr^{\leq \reg}N_{W,s,ij}^{(p)}|^2\bigg).
\eea
Since we have, in view of Step 5,
\beaa
\psi^{(p)}_{s,ij} &=& \chi_{\tau_2}\widetilde{\phi}^{(p)}_{s,ij}+(1-\chi_{\tau_2})\breve{\phi}^{(p)}_{s,ij}\quad\textrm{on}\quad\MM(\tau_1, +\infty), 
\eeaa
we immediately infer from \eqref{eq:localenergyestimateforwidetildephisijponSigmatau2minus3tau2} and 
\eqref{eq:localenergyestimateforbrevephisijponSigmatau2minus1tau2}, for $\reg\leq {14}$,
\beaa
\nn&&\sum_{i,j=1}^3\sum_{p=0}^2\EMF_\de[\pr^{\leq \reg}\psi^{(p)}_{s,ij}](\tau_2-3, \tau_2)\\ 
&\les& \sum_{i,j=1}^3\sum_{p=0}^2\bigg(\E[\pr^{\leq \reg}\phi^{(p)}_{s,ij}](\tau_2-3)+\int_{\MM(\tau_2-3, \tau_2-2)}r^{1+\de} |\pr^{\leq \reg}N_{W,s,ij}^{(p)}|^2\bigg),
\eeaa
which is the desired estimate \eqref{eq:localenergyestimateforpsisijponSigmatau2minus3tau2:partialderivatives}. 

\noindent{\bf Step 8.} Next, we estimate the tensorization defect $\err_{\textrm{TDefect}}[\psi^{(p)}_{s}]$ corresponding to the family of complex-valued scalars $\psi_{ij}$ as introduced in  Definition \ref{def:definitionofthenotationerrforthescalarizationdefect}. First, recall from Step 5 that 
\beaa
\psi^{(p)}_{s,ij}=\phi^{(p)}_{s,ij}\quad\textrm{on}\quad \MM(\tau_1+1, \tau_2-3), \qquad \psi^{(p)}_{s,ij} =\breve{\phi}^{(p)}_{s,ij}\quad\textrm{on}\quad\MM(\tau_2, +\infty),
\eeaa
which together with \eqref{eq:thetensordeffectofbrevephivanishesidentically}, and the fact that $\phi^{(p)}_{s,ij}=\pmb\phi^{(p)}_{s}(\Om_i, \Om_j)$ by definition, yields
\beaa
\err_{\textrm{TDefect}}[\psi^{(p)}_{s}]=0\quad\textrm{on}\quad\MM(\tau_1+1, +\infty)\setminus(\tau_2-3, \tau_2).
\eeaa

Then, we consider the range $\tau\in[\tau_2-3, \tau_2-2)$ and, to this end, we introduce the following auxiliary coupled system of tensorial wave equations, for $s=\pm 2$, $p=0,1,2$,
\beaa
\bsplit
&\bigg(\squared_2 -\frac{4ia\cos\th}{|q|^2}\nab_{\pr_{\tt}}- \frac{4-2\de_{p0}}{\qs}\bigg){(\pmb\phi_{aux,1})_s^{(p)}} = \chi^{(1)}_{\tau_1, \tau_2}\big(\L_{s}^{(p)}[(\pmb\phi_{aux,1})_s]+\N_{W,s}^{(p)}\big),  \,\,\textrm{on}\,\,\MM(\tau_2-3, \tau_2),\\
&(\pmb\phi_{aux,1})^{(p)}_s =\pmb\phi^{(p)}_{s}, \quad \nab_{N_{\Sigma(\tau_2-3)}}(\pmb\phi_{aux,1})^{(p)}_s =\nab_{N_{\Sigma(\tau_2-3)}}\pmb\phi^{(p)}_{s} \,\,\textrm{on}\,\,\Sigma(\tau_2-3).
\end{split}
\eeaa
Defining $(\phi_{aux,1})_{s,ij}^{(p)}:=(\pmb\phi_{aux,1})_s^{(p)}(\Om_i, \Om_j)$, and arguing as in Section \ref{sec:TeukolskyWavesysteminperturbationsofKerr:scalarizedform}, one easily checks that $(\phi_{aux,1})_{s,ij}^{(p)}$ satisfies the same system of wave equations as $\psi_{s,ij}^{(p)}$ on $\MM(\tau_2-3, \tau_2-2)$ and the same initial data on $\Si(\tau_2-3)$. By uniqueness for the wave equation, we infer 
\beaa
\psi_{s,ij}^{(p)}=(\phi_{aux,1})_{s,ij}^{(p)}=(\pmb\phi_{aux,1})_s^{(p)}(\Om_i, \Om_j)\quad\textrm{on}\quad\MM(\tau_2-3, \tau_2-2)
\eeaa
which implies that $\err_{\textrm{TDefect}}[\psi^{(p)}_{s}]=0$ on $\MM(\tau_2-3, \tau_2-2)$. In view of the above, we deduce 
\beaa
\err_{\textrm{TDefect}}[\psi^{(p)}_{s}]=0\quad\textrm{on}\quad\MM(\tau_1+1, +\infty)\setminus(\tau_2-2, \tau_2).
\eeaa

Next, we derive energy estimates for $\err_{\textrm{TDefect}}[\psi^{(p)}_{s}]$ on $\MM(\tau_2-2, \tau_2)$. To this end, we first derive a system of wave equations for $\err_{\textrm{TDefect}}[\widetilde{\phi}^{(p)}_{s}]$. In view of \eqref{eq:waveequationdefiningfirstextensionwidetildepsi}, we have on $\MM(\tau_2-2, \tau_2)$
\beaa
\big({\square}_{\g_{\chi_{\tau_1, \tau_2}}}-(4-2\de_{p0})|q|^{-2}\big)\widetilde{\phi}^{(p)}_{s,ij} &=& \chi_{\tau_1, \tau_2}\big( \widehat{S}(\widetilde{\pmb\phi}^{(p)}_s)_{ij} +(\widehat{Q}\widetilde{\pmb\phi}^{(p)}_s)_{ij}\big)\\
&& +(1-\chi_{\tau_1, \tau_2})\big( \widehat{S}_K(\widetilde{\pmb\phi}^{(p)}_s)_{ij} +(\widehat{Q}_K\widetilde{\pmb\phi}^{(p)}_s)_{ij}+f_p\widetilde{\phi}^{(p)}_{s,ij}\big)
\eeaa 
and hence, together with \eqref{eq:assumptionsonregulartripletinperturbationsofKerr:0}, we infer 
\bea\lab{eq:waveequationdefiningfirstextensionwidetildepsi:tau2minus2totau2}
\bsplit
\Big({\square}_{\g_{\chi_{\tau_1, \tau_2}}}+V\Big)\widetilde{\phi}^{(p)}_{s,ij} =&  \widehat{S}(\widetilde{\pmb\phi}^{(p)}_s)_{ij} +(\widehat{Q}\widetilde{\pmb\phi}^{(p)}_s)_{ij} +\Ga_g\dk^{\leq 1}\widetilde{\phi}^{(p)}_{s},\quad\textrm{on}\quad\MM(\tau_2-2, \tau_2),\\
V:=& -(4-2\de_{p0})|q|^{-2}-(1-\chi_{\tau_1, \tau_2})f_p.
\end{split}
\eea 
Also, in view of Lemma \ref{lemma:computationoftheMialphajinKerr} and \eqref{eq:assumptionsonregulartripletinperturbationsofKerr:0}, and \eqref{eq:definitionofthenotationerrforthescalarizationdefect}, we have
\beaa
&& M_{i\tau}^l\widetilde{\phi}^{(p)}_{s,lj}+M_{j\tau}^l\widetilde{\phi}^{(p)}_{s,il}\\ 
&=& -\frac{2amr\cos\th}{|q|^4}\Big(\in_{ilk}x^k\widetilde{\phi}^{(p)}_{s,lj}+\in_{jlk}x^k\widetilde{\phi}^{(p)}_{s,il}\Big)+\Ga_g\widetilde{\phi}^{(p)}_{s}\\
&=& -\frac{2amr\cos\th}{|q|^4}\Big(-2i\widetilde{\phi}^{(p)}_{s,ij}+(\err_{\textrm{TDefect},5}[\psi])_{ij}+(\err_{\textrm{TDefect},5}[\psi])_{ji}\\
&&+i(\err_{\textrm{TDefect},1}[\psi])_{ij}+\in_{jlk}x^k(\err_{\textrm{TDefect},1}[\psi])_{il}\Big)+\Ga_g\widetilde{\phi}^{(p)}_{s}.
\eeaa
Together with \eqref{eq:waveequationdefiningfirstextensionwidetildepsi:tau2minus2totau2} and Lemma \ref{lemma:waveequationsfortensordeffects}, we infer, using also Lemma \ref{lemma:computationoftheMialphajinKerr}, \eqref{eq:assumptionsonregulartripletinperturbationsofKerr:0}, and the properties of $\g_{\chi_{\tau_1, \tau_2}}$,
\beaa
\square_{\g_{\chi_{\tau_1, \tau_2}}}\err_{\textrm{TDefect}}[\widetilde{\phi}^{(p)}_{s}] =
O(r^{-2})\dk^{\leq 1}\err_{\textrm{TDefect}}[\widetilde{\phi}^{(p)}_{s}]+\Ga_g\dk^{\leq 1}\widetilde{\phi}^{(p)}_{s},\quad\textrm{on}\quad\MM(\tau_2-2, \tau_2).
\eeaa
Since the initial data for $\err_{\textrm{TDefect}}[\widetilde{\phi}^{(p)}_{s}]$ is trivial at $\tau=\tau_2-2$, using the local energy estimate \eqref{eq:localenergyestimate:future:unweigthedderivatives:bis}, we infer, for $\reg\leq 14$,
\beaa
\EMF_\de[\pr^{\leq \reg}\err_{\textrm{TDefect}}[\widetilde{\phi}^{(p)}_{s}]](\tau_2-2, \tau_2)\les \ep^2\EM_\de[\pr^{\leq \reg}\widetilde{\phi}^{(p)}_{s}](\tau_2-2, \tau_2).
\eeaa
Together with \eqref{eq:localenergyestimateforwidetildephisijponSigmatau2minus3tau2}, we deduce, for $\reg\leq 14$,
\bea\lab{eq:controlscalarizationdefectforwidetildephiontau2minus2tau2}
\nn&&\sum_{p=0}^2\EMF_\de[\pr^{\leq \reg}\err_{\textrm{TDefect}}[\widetilde{\phi}^{(p)}_{s}]](\tau_2-2, \tau_2)\\ &\les& \ep^2\sum_{i,j=1}^3\sum_{p=0}^2\bigg(\E[\pr^{\leq \reg}\phi^{(p)}_{s,ij}](\tau_2-3)+\int_{\MM(\tau_2-3, \tau_2-2)}r^{1+\de} |\pr^{\leq \reg}N_{W,s,ij}^{(p)}|^2\bigg).
\eea
Now, recalling from Step 5 that 
\beaa
\psi^{(p)}_{s,ij} = \chi_{\tau_2}\widetilde{\phi}^{(p)}_{s,ij}+(1-\chi_{\tau_2})\breve{\phi}^{(p)}_{s,ij}\quad\textrm{on}\quad\MM(\tau_1, +\infty), 
\eeaa
we infer
\beaa
\err_{\textrm{TDefect}}[\psi^{(p)}_{s}] = \chi_{\tau_2}\err_{\textrm{TDefect}}[\widetilde{\phi}^{(p)}_{s}]+(1-\chi_{\tau_2})\err_{\textrm{TDefect}}[\breve{\phi}^{(p)}_{s}]\quad\textrm{on}\quad\MM(\tau_1, +\infty), 
\eeaa
which together wit \eqref{eq:thetensordeffectofbrevephivanishesidentically} yields
\beaa
\err_{\textrm{TDefect}}[\psi^{(p)}_{s}] = \chi_{\tau_2}\err_{\textrm{TDefect}}[\widetilde{\phi}^{(p)}_{s}]\quad\textrm{on}\quad\MM(\tau_1, +\infty).
\eeaa
In view of the above control for $\err_{\textrm{TDefect}}[\widetilde{\phi}^{(p)}_{s}]$, this yields, for $\reg\leq 14$,
\beaa
\nn&& \sum_{p=0}^2\EMF_\de[\pr^{\leq \reg}\err_{\textrm{TDefect}}[\psi^{(p)}_{s}]](\tau_2-2, \tau_2)\\
&\les& \ep^2\sum_{i,j=1}^3\sum_{p=0}^2\bigg(\E[\pr^{\leq \reg}\phi^{(p)}_{s,ij}](\tau_2-3)+\int_{\MM(\tau_2-3, \tau_2-2)}r^{1+\de} |\pr^{\leq \reg}N_{W,s,ij}^{(p)}|^2\bigg).
\eeaa

In addition, it will also be useful to control $\widetilde{\phi}^{(p)}_{s,ij}- \breve{\phi}^{(p)}_{s,ij}$ on $\MM(\tau_2-1, \tau_2)$. In view of \eqref{eq:waveequationdefiningfirstextensionwidetildepsi} and \eqref{eq:waveequationdefiningfirstextensionbrevepsi}, $\widetilde{\phi}^{(p)}_{s,ij}- \breve{\phi}^{(p)}_{s,ij}$ satisfies
\beaa
\bsplit
\big({\square}_{\gam}-(4-2\de_{p0})|q|^{-2}\big)(\widetilde{\phi}-\breve{\phi})^{(p)}_{s,ij} =&  \big(\widehat{S}_K +\widehat{Q}_K+f_p\big)(\widetilde{\phi}-\breve{\phi})^{(p)}_{s,ij},\,\,\,\textrm{on}\,\,\,\MM(\tau_2-1, +\infty), \\
(\widetilde{\phi}-\breve{\phi})^{(p)}_{s,ij} =& (\widetilde{\phi}^{(p)}_{s}-\Pi_2[\widetilde{\phi}_{s}^{(p)}])_{ij}\quad\textrm{on}\quad\Si(\tau_2-1), \\ 
 N_{\Sigma(\tau_2-1)}(\widetilde{\phi}-\breve{\phi})^{(p)}_{s,ij}=&N_{\Sigma(\tau_2-1)}(\widetilde{\phi}^{(p)}_{s}-\Pi_2[\widetilde{\phi}^{(p)}_{s}])_{ij}\quad\textrm{on}\quad\Si(\tau_2-1).
\end{split}
\eeaa
Thus, applying the local energy estimate \eqref{eq:localenergyestimate:future:unweigthedderivatives:bis} with $\reg\leq {14}$, we infer
\beaa
\nn&&\sum_{i,j=1}^3\sum_{p=0}^2\EMF_\de[\pr^{\leq \reg}(\widetilde{\phi}-\breve{\phi})^{(p)}_{s,ij}](\tau_2-1, \tau_2) \les \sum_{i,j=1}^3\sum_{p=0}^2\E[\pr^{\leq \reg}(\widetilde{\phi}^{(p)}_{s}-\Pi_2[\widetilde{\phi}^{(p)}_{s}])_{ij}](\tau_2-1).
\eeaa
 Together with Lemma \ref{lemma:computationerrorscalarizationdeffect}, we deduce, for 
  $\reg\leq {14}$, 
\beaa
\sum_{i,j=1}^3\sum_{p=0}^2\EMF_\de[\pr^{\leq \reg}(\widetilde{\phi}-\breve{\phi})^{(p)}_{s,ij}](\tau_2-1, \tau_2) \les \sum_{p=0}^2\E[\pr^{\leq \reg}\err_{\textrm{TDefect}}[\widetilde{\phi}^{(p)}_{s}]](\tau_2-1).
\eeaa
Plugging \eqref{eq:controlscalarizationdefectforwidetildephiontau2minus2tau2}, we infer, for $\reg\leq {14}$, 
\bea\lab{eq:controlofwidetildephiminusbreevephiontau2minus1tau2}
\nn&&\sum_{i,j=1}^3\sum_{p=0}^2\EMF_\de[\pr^{\leq \reg}((\widetilde{\phi}-\breve{\phi})^{(p)}_{s,ij})](\tau_2-1, \tau_2)\\
&\les& \ep^2\sum_{i,j=1}^3\sum_{p=0}^2\bigg(\E[\pr^{\leq \reg}\phi^{(p)}_{s,ij}](\tau_2-3)+\int_{\MM(\tau_2-3, \tau_2-2)}r^{1+\de} |\pr^{\leq \reg}N_{W,s,ij}^{(p)}|^2\bigg).
\eea

\noindent{\bf Step 9.} Finally, we have obtained the following:
\begin{itemize}
\item in view of \eqref{eq:extendedmetricgchitau1tau2} and Step 1, 
$\g_{\chi_{\tau_1, \tau_2}}$ satisfies the assumptions of Section \ref{subsubsect:assumps:perturbedmetric} and coincides with Kerr in $\MM\setminus(\tau_1, \tau_2)$, 

\item in view of Step 5 and the fact that $F_{total,s,ij}^{(p)}=\widehat{F}^{(p)}_{s,ij}+\widetilde{F}^{(p)}_{s,ij}+\underline{F}^{(p)}_{s,ij}+\breve{F}^{(p)}_{s,ij}$, $\psi^{(p)}_{s,ij}$ satisfies \eqref{eq:waveeqwidetildepsi1} \eqref{eq:causlityrelationsforwidetildepsi1},

\item in view of Step 6, $\psi^{(p)}_{s,ij}$ satisfies \eqref{eq:localenergyforpsiontau1minus1tau1plus1},

\item in view of Step 7, $\psi^{(p)}_{s,ij}$ satisfies \eqref{eq:localenergyestimateforpsisijponSigmatau2minus3tau2:partialderivatives},

\item and in view of Step 8, $\psi^{(p)}_{s,ij}$ satisfies \eqref{eq:propertyofscalarizationdefect:prop1} and \eqref{eq:propertyofscalarizationdefect:prop2}.
\end{itemize}
This concludes the proof of Proposition \ref{prop:extensionprocedureoftheTeukolskywaveequations}.
\end{proof}


\subsubsection{Global energy-Morawetz estimates for unweighted derivatives of solutions to \eqref{eq:waveeqwidetildepsi1}}

 
In order to prove our main Theorem \ref{thm:main}, i.e., the derivation of energy-Morawetz estimates for $\tau$ in $(\tau_1, \tau_2)$, we first state in this section global energy-Morawetz estimates for \eqref{eq:waveeqwidetildepsi1}, i.e., energy-Morawetz estimates for $\tau$ in $\mathbb{R}$, that hold for unweighted derivatives $\pr$ introduced in \eqref{eq:defunweightedderivative}.

\begin{theorem}\lab{th:main:intermediary}
Let $(\MM, \g)$ satisfy the assumptions of Sections \ref{subsect:assumps:perturbednullframe}, \ref{subsubsect:assumps:perturbedmetric} and \ref{sec:regulartripletinperturbationsofKerr}, and assume that $\tt_1$ and $\tt_2$ satisfy \eqref{eq:tau2geqtau1plus10conditionisthemaincase}. Let $\pmb\phi_s^{(p)}$, $s=\pm 2$, $p=0,1,2$, be a solution to the tensorial Teukolsky wave/transport systems \eqref{eq:TensorialTeuSysandlinearterms:rescaleRHScontaine2:general:Kerrperturbation}  \eqref{def:TensorialTeuScalars:wavesystem:Kerrperturbation}  in perturbations of Kerr, 
and let $\psis{p}$, $s=\pm 2$, $p=0,1,2$, be a solution to \eqref{eq:waveeqwidetildepsi1} satisfying \eqref{eq:causlityrelationsforwidetildepsi1}. Then, we have, for any $\reg\leq 14$ and $0<\de\leq \frac{1}{3}$,
\bea
\lab{eq:mainhighorderunweightedEMF:maintheorem:pm2}
&&\sum_{p=0}^2\EMF_{\de}[\pr^{\leq \reg}\phis{p}](\tau_1,\tau_2)
+ \sum_{p=0}^2{\EMF}_{\de}[\pr^{\leq \reg}\psis{p}](\Iti)\nn\\
&\les&\sum_{p=0}^2\E[\pr^{\leq \reg}\phis{p}](\tau_1)+\sum_{p=0}^2\NN_\de[\pr^{\leq \reg}\phis{p}, \pr^{\leq \reg}\N_{W,s}^{(p)}](\tt_1, \tt_2)+\sum_{p=0}^1\int_{\MM(\tau_1,\tau_2)}r^{-3+\de}|\pr^{\leq \reg+1}\N_{T,s}^{(p)}|^2\nn\\
&&
+\sum_{p=0}^1\int_{\MM_{r\geq 10m}(\tau_1,\tau_2)} \Big(r^{-1+\de}|\pr^{\leq \reg}\N^{(p)}_{W,s}|+r^{-2+\de}|\pr^{\leq \reg+1}\N^{(p)}_{T,s}|\Big)|\pr^{\leq \reg}\phis{p}|\nn\\
&& +\sum_{p=0}^2\Big(\EMF[\psis{p}](\Iti)
+\EMF[\phis{p}](\tau_1, \tau_2)\Big).
\eea
\end{theorem}

\begin{remark}
In view of the proof of Theorem \ref{th:main:intermediary}, the terms in the last line of the RHS of \eqref{eq:mainhighorderunweightedEMF:maintheorem:pm2} can in fact be replaced by 
$$\sum\limits_{p=0}^2(\A[\psis{p}](\Iti)+\A[\phis{p}](\tau_1, \tau_2)),$$ 
with $\A[\cdot]$ given as in \eqref{def:AandAonorms:tau1tau2:phisandpsis}, so that Theorem \ref{th:main:intermediary} may be upgraded to high-order unweighed energy-Morawetz estimates conditional only on the control of zeroth-order derivatives of $\psis{p}$ and $\phis{p}$.
\end{remark}

The proof of Theorem \ref{th:main:intermediary} relies in particular on the microlocal energy-Morawetz estimates in $\Mtrap$ derived in Section \ref{sect:microlocalenergyMorawetztensorialwaveequation} and is thus postponed to Section \ref{sec:proofofth:main:intermediary}.


\subsubsection{Energy-Morawetz estimates near infinity in perturbations of Kerr}


The following proposition provides energy-Morawetz estimates for Teukolsky equations in perturbations of Kerr for $r\geq R$ with $R$ large enough.

\begin{proposition}\lab{prop:EnergyMorawetznearinfinitytensorialTeuk}
Let $(\MM, \g)$ satisfy the assumptions of Sections \ref{subsect:assumps:perturbednullframe} and \ref{subsubsect:assumps:perturbedmetric}.  We have for solutions $\pmb\phi_s^{(p)}$, $s=\pm 2$, $p=0,1,2$, to the tensorial Teukolsky wave/transport systems \eqref{eq:TensorialTeuSysandlinearterms:rescaleRHScontaine2:general:Kerrperturbation}  \eqref{def:TensorialTeuScalars:wavesystem:Kerrperturbation}  in perturbations of Kerr the following energy-Morawetz-flux estimates, for any $1\leq\tau_1<\tau_2 <+\infty$ and any $0<\de\leq \frac{1}{3}$, for $\reg\leq 14$, and for $R\geq 20 m$ large enough, 
\bea
\lab{eq:EMnearinfinity:highorderweightedderivatives:Teu:pm2}
\nn&&\sum_{p=0}^2\EMF^{(\reg)}_{\de, \geq R}[\pmb\phi_{s}^{(p)}](\tau_1, \tau_2) +\sum_{p=0}^1\int_{\MM_{r\geq R}(\tau_1,\tau_2)}r^{-3+\de} \big|(r\nab)^{\leq 1} \dk^{\leq \reg}\phis{p}\big|^2\\
\nn&\les_R& \sum_{p=0}^2\M^{(\reg)}_{R/2, R}[\pmb\phi_{s}^{(p)}](\tau_1, \tau_2)+\sum_{p=0}^2\E^{(\reg)}_{r\geq R/2}[\pmb\phi_{s}^{(p)}](\tau_1)+\sum_{p=0}^2\int_{\MM_{r\geq R/2}(\tt_1, \tt_2)}r^{1+\de}|\dk^{\leq \reg}\N_{W,s}^{(p)}|^2\\
&&+\sum_{p=0}^1\int_{\MM_{r\geq R/2}(\tt_1, \tt_2)}r^{-1+\de}|\dk^{\leq \reg+1} \N^{(p)}_{T,s}|^2.
\eea
\end{proposition}

The proof of Proposition \ref{prop:EnergyMorawetznearinfinitytensorialTeuk} is postponed to Section  \ref{sec:proofofprop:EnergyMorawetznearinfinitytensorialTeuk}.


\subsubsection{Energy-Morawetz estimates for tensorial wave equations in subextremal Kerr}
\lab{sec:energyMorwetzestimatesinKerr}


In order to control the lower order terms appearing last on the RHS of \eqref{eq:mainhighorderunweightedEMF:maintheorem:pm2} and  \eqref{eq:mainhighorderunweightedEMF:maintheorem:pm2}, we will rely on the two energy-Morawetz estimates in Kerr stated below. First, we consider solutions $\pmb\phi_s\in\sk_2(\mathbb{C})$, $s=\pm 2$, to the following tensorial wave equations on a subextremal Kerr background 
\bea\lab{eq:TeukolskyequationforAandAbintensorialforminKerr:inhomogenouscase}
\nn\left(\squared_2 -\frac{4ia\cos\th}{|q|^2}\nab_{\pr_t} -  \frac{s}{|q|^2}\right)\pmb\phi_s +\frac{2s}{|q|^2}(r-m)\nab_3\pmb\phi_s  -\frac{4sr}{|q|^2}\nab_{\pr_t}\pmb\phi_s\\
+\frac{4a\cos\th}{|q|^6}\Big(a\cos\th\big(|q|^2+6mr\big) - is\big((r-m)|q|^2+4mr^2\big)\Big)\pmb\phi_s &=& \N_s,
\eea
with $\N_s\in\sk_2(\mathbb{C})$, $s=\pm 2$, where \eqref{eq:TeukolskyequationforAandAbintensorialforminKerr:inhomogenouscase} is the inhomogeneous version of the tensorial Teukolsky equations \eqref{eq:TeukolskyequationforAandAbintensorialforminKerr}  in Kerr.

\begin{theorem}[Weak-Morawetz for Teukolsky in subextremal Kerr]
\lab{cor:weakMorawetzforTeukolskyfromMillet:bis}
Let $\tau_0\geq 1$ and assume that $\pmb\phi_s\in\sk_2(\mathbb{C})$, $s=\pm 2$, satisfies the inhomogeneous Teukosky equation 
\eqref{eq:TeukolskyequationforAandAbintensorialforminKerr:inhomogenouscase}  on a subextremal Kerr background.  Then, for any $\de\in (0,\frac{1}{3}]$, we have
\bea\lab{eq:weakMorawetzforTeukolskyfromMillet:bis:minus2}
\nn\int_{\MM(\tau\geq\tau_0)}r^{-3+\de}|\dk^{\leq 3}\pmb\phi_{-2}|^2 &\les& \E^{(13)}[\pmb\phi_{-2}](\tau_0)+\E^{(11)}[r^{\frac{1+\de}{2}}\nab_{\pr_r}(r\pmb\phi_{-2})](\tau_0)\\
&&+\int_{\MM(\tau\geq\tau_0)}r^{1+\de}|\dk^{\leq 13}\N_{-2}|^2,
\eea
and, for any $\de\in(0,\frac{1}{3}]$, any $R_0\geq 10m$ and any $\tau_1>\tau_0+1$, we have
\bea\lab{eq:weakMorawetzforTeukolskyfromMillet:plus2:nolossinr}
&&\int_{\MM(\tau_0,\tau_1-1)}r^{-11+\frac{\de}{2}}|\dk^{\leq 3}\pmb\phi_{+2}|^2\nn\\
 &\les&R_0^{1+\frac{\de}{2}}\E_{r\leq 2R_0}^{(9)}[r^{-4}\pmb\phi_{+2}](\tau_0) +\int_{\MM(\tau_0,\tau_1)}r^{-7+\de}|\dk^{\leq 10}\N_{+2}|^2 \nn\\
&&+R_0^{-\frac{\de}{2}}\bigg({\EMF}^{(11)}_{\de}[r^{-4}\pmb\phi_{+2}](\tau_0,\tau_1)+\int_{\MM(\tau_0,\tau_1)} r^{-3+\de}|\dk^{\leq 11} \nab_{3}(r^{-3}\pmb\phi_{+2})|^2\bigg).
\eea
\end{theorem}

The proof of Theorem \ref{cor:weakMorawetzforTeukolskyfromMillet:bis} is postponed to Section \ref{sec:MorawetzestimatesforTeukolskys=pm2inKerrfromMillet}. Next, we also state an energy-Morawetz estimates for $\breve{\pmb\phi}\in\sk_2(\mathbb{C})$ satisfying the following tensorial wave equation in a subextremal Kerr background
\bea\lab{eq:basictensorialwaveequationinKerr:forbrevephi}
\squared_2\breve{\pmb\phi} -\frac{4ia\cos\th}{|q|^2}\nab_{\pr_t}\breve{\pmb\phi} -\bigg(\frac{4}{\qs}-\frac{4a^2\cos^2\th}{|q|^6}\Big(|q|^2+6mr\Big)\bigg)\breve{\pmb\phi}=0.
\eea

\begin{theorem}[Energy-Morawetz for \eqref{eq:basictensorialwaveequationinKerr:forbrevephi} in subextremal Kerr]
\lab{prop:weakMorawetzfortensorialwaveeqfromDRSR}
Let $\breve{\pmb\phi}\in\sk_2(\mathbb{C})$ be a solution in a subextremal Kerr background to the tensorial wave equation \eqref{eq:basictensorialwaveequationinKerr:forbrevephi} in $\MM(\tau\geq \tau_0)$, 
where $\tau_0\geq 1$ is a constant. Then, we have 
\bea
\lab{eq:EMFestimates:tensorialwaveinsubextremalKerr:weakMorathm}
\EMF[\breve{\pmb\phi}](\tau_0,+\infty) \les \E[\breve{\pmb\phi}](\tau_0).
\eea
\end{theorem}

\begin{proof}
By subtracting $\frac{4}{\qs}{\breve{\pmb\phi}}(e_1, e_1)$ from both sides of the {identity} in Lemma \ref{lem:scalarizationofspinweightedwaveoperator:Kerr}, we deduce
\beaa
&& \square_{\gam}({\breve{\pmb\phi}}(e_1, e_1)) + \frac{4i}{|q|^2}\left(\frac{\cos\th}{\sin^2\th}\pr_{\phi} -a\cos\th\pr_{t}\right)({\breve{\pmb\phi}}(e_1, e_1)) -\frac{4}{|q|^2\sin^2\th}{\breve{\pmb\phi}}(e_1, e_1)\\ 
&=& \left(\squared_2{\breve{\pmb\phi}} -\frac{4ia\cos\th}{|q|^2}\nab_{\pr_t}{\breve{\pmb\phi}}-\bigg(\frac{4}{\qs}-\frac{4a^2\cos^2\th}{|q|^6}\Big(|q|^2+6mr\Big)\bigg){\breve{\pmb\phi}}\right)(e_1, e_1) .
\eeaa
The RHS of the above equation vanishes {since $\breve{\pmb\phi}$ satisfies  \eqref{eq:basictensorialwaveequationinKerr:forbrevephi}}, hence, by denoting
\bea
\lab{eq:psi2definitionNPfromtensor}
{\breve{\phi}_{2,\text{NP}}}:={\breve{\pmb\phi}}(e_1, e_1),
\eea
we infer that ${\breve{\phi}_{2,\text{NP}}}$ satisfies the following wave equation, in Boyer-Lindquist coordinates,
\bea
\lab{eq:spinweightedwaveequation:NPform:Kerr:proof}
 \bigg(\qs \square_{\gam}+ 4i \left(\frac{\cos\th}{\sin^2\th}\pr_{\phi} -a\cos\th\pr_{t}\right)-\frac{4}{\sin^2\th}\bigg){\breve{\phi}_{2,\text{NP}}}=0.
\eea
This is an analog\footnote{As for the scalar wave operator, this wave operator is easily seen to be separable, and by separation {of} variables in the coordinates $t$ and $\phi$, this wave equation is equivalent to a sum of a radial ODE and an angular ODE, where the radial ODE is the same as the one for the scalar wave equation, while the angular ODE is a spin-$2$ Teukolsky angular equation as opposed to the spin-$0$ angular equation of the scalar wave equation.}  of the spin-$0$ scalar wave equation $\square_{\gam}\psi=0$ and, to show energy-Morawetz estimates for this wave equation, we follow the argument in \cite{DRSR} for scalar wave equation. In \cite{DRSR}, the argument crucially relies on, after separation of variables in $t$ and $\phi$ with frequency parameters $\om$ and $M$ respectively, the following bounds of the eigenvalues $\{\la_{ML}\}_{M\in \mathbb{Z}, L\geq |M|}$, indexed by the parameter $L$, of the operator $-\pr_{\th\th} - (\sin\th)^{-2} M^2 -a^2\sin^2\th \om^2$, the analog in frequency space of the Carter operator $-(\De_{\mathbb{S}^2} + a^2\sin^2\th \pr_{tt})$:
\bea
\lab{eq:eigenvalueofCarterforscalarwaveinKerr}
\la_{ML}\geq \max\{2|a M\om|, |M| (|M|+1)\},
\eea
where $\De_{\mathbb{S}^2}$ is the spherical Laplacian on unit sphere.
For equation \eqref{eq:spinweightedwaveequation:NPform:Kerr:proof}, we can also do separation in variables in $t$ and $\phi$ with frequencies $\om$ and $M$, and consider the eigenvalues $\{\nu_{ML}\}_{M\in \mathbb{Z}, L\geq \max\{|M|,2\}}$ of the operator 
$$
-\pr_{\th\th} - \frac{M^2}{\sin^2\th}  -a^2\sin^2\th \om^2+\frac{4\cos\th}{\sin^2\th}M + 4a\cos\th \om +\frac{4}{\sin^2\th},
$$
the analog in frequency space of the Teukolsky angular operator in Boyer-Lindquist coordinates
$$
-(\De_{\mathbb{S}^2} + a^2\sin^2\th \pr_{tt})+2si \left(\frac{\cos\th}{\sin^2\th}\pr_{\phi} -a\cos\th\pr_{t}\right)-\frac{s^2}{\sin^2\th}
$$
with $s=+2$.
As shown in \cite{PTII73}\footnote{The above Teukolsky angular operator equals the sum of the Teukolsky angular operator with spin weight $s=+2$ in \cite{PTII73} and a potential $-s^2$ which equals $4$.}, these eigenvalues satisfy 
\bea
\nu_{ML}\geq \max\{2|a M\om|, |M| (|M|+1), 4\},
\eea
where, in particular, the bound \eqref{eq:eigenvalueofCarterforscalarwaveinKerr} holds in this case. Then, after these separation in variables, we obtain the exactly same radial equation as in \cite{DRSR} for scalar wave, and by taking exactly the same multiplier and arguing in the same manner as in \cite{DRSR}, we infer an analog of the energy-Morawetz-flux estimate proven in \cite{DRSR} which, together with the definition \eqref{eq:psi2definitionNPfromtensor} of ${\breve{\phi}_{2,\text{NP}}}$ and the assumption ${\breve{\pmb\phi}}\in\sk_2(\mathbb{C})$, yields
\beaa
\EMF[{\breve{\pmb\phi}}](\tau_0,+\infty) \les \E[{\breve{\pmb\phi}}](\tau_0)
\eeaa
as desired. This conclude the proof of Theorem \ref{prop:weakMorawetzfortensorialwaveeqfromDRSR}.
\end{proof}


\subsection{Proof of Theorem \ref{thm:main}}
\lab{subsect:proofofThm4.1}


Let $(\MM, \g)$ satisfy the assumptions of Sections \ref{subsect:assumps:perturbednullframe}, \ref{subsubsect:assumps:perturbedmetric} and \ref{sec:regulartripletinperturbationsofKerr}, let $\tau_1$ and $\tau_2$ be such that $1\leq\tau_1<\tau_2 <+\infty$, let $\phis{p}$, $s=\pm 2$, $p=0,1,2$, be solutions to the tensorial Teukolsky wave/transport systems \eqref{eq:TensorialTeuSysandlinearterms:rescaleRHScontaine2:general:Kerrperturbation}  \eqref{def:TensorialTeuScalars:wavesystem:Kerrperturbation}  in perturbations of Kerr, and  assume also that $\Ab$ satisfies \eqref{eq:waveequationpmbphip=0sminus2nodeginredshiftregion}. We proceed in the following steps.

\noindent{\bf Step 0.} We consider first the case $\tau_1<\tau_2\leq \tau_1+10$. As the scalarized Teukolsky wave system \eqref{eq:ScalarizedTeuSys:general:Kerrperturbation} satisfied by $\phi^{(p)}_{s,ij}=\pmb\phi_s^{(p)}(\Om_i, \Om_j)$ is of the form \eqref{eq:eqsforlocalenergyestimatelemma:general} with $\psi_{ij}=(\phi^{(p)}_{s,ij})_{p=0,1,2}$, $D_1=0$ and $F_{ij}=(N^{(p)}_{W,s,ij})_{p=0,1,2}$, we infer from \eqref{eq:localenergyestimate:future:bis}, with $\tau_0=\tau_1$, $0<q=\tau_2-\tau_1\leq 10$, $0\leq\reg\leq 14$ and $0<\de\leq \frac{1}{3}$, for any $\tau_1<\tau_2\leq \tau_1+10$,
\bea\lab{eq:intermediarylocalenergyestimatetoprovemainThoerem:particularcasetau2leqtau1plus10}
\sum_{p=0}^2\EMF_\de^{(\reg)}[\pmb\phi^{(p)}_s](\tau_1, \tau_2) &\les & \sum_{p=0}^2\left(\E^{(\reg)}[\pmb\phi^{(p)}_s](\tau_1) +\int_{\MM(\tau_1, \tau_2)}r^{1+\de}|\dk^{\leq \reg}\N^{(p)}_{W,s}|^2\right).
\eea
This estimate with $\reg=11$ yields the desired estimate \eqref{MainEnerMora:psi:plus2case} for $s=+2$.

In the case $s=-2$, we deduce from \eqref{eq:intermediarylocalenergyestimatetoprovemainThoerem:particularcasetau2leqtau1plus10} and the redshift estimate \eqref{eq:redshift:Ab:highorderregularity} with $\reg=14$ that
\beaa
\nn&&\sum_{p=0}^2\EMF^{({14})}_\de[\pmb\phi_{-2}^{(p)}](\tau_1, \tau_2) + \sum_{p=0}^2\EMF^{({14})}_{r\leq r_+(1+\dred)}[\nab_4^p\Ab](\tau_1, \tau_2)\\
\nn&\les& \sum_{p=0}^2\E^{({14})}[\pmb\phi_{-2}^{(p)}](\tau_1)+\sum_{p=0}^2\E^{({14})}_{r\leq r_+(1+\dred)}[\nab_4^p\Ab](\tau_1) +\sum_{p=0}^1\int_{\MM_{r\leq r_+(1+2\dred)}(\tau_1, \tau_2)}|\dk^{\leq 15}\N_{T,-2}^{(p)}|^2\\
\nn&&+\sum_{p=0}^2\int_{\MM(\tau_1, \tau_2)}r^{1+\de}|\dk^{\leq 14}\N^{(p)}_{W,-2}|^2+\sum_{p=0}^2\int_{\MM_{r\leq r_+(1+2\dred)}(\tau_1, \tau_2)}|\dk^{\leq {14}}\N_{\nab_4^p\Ab}|^2,
\eeaa
which immediately yields the desired estimate \eqref{MainEnerMora:psi:minus2case}. 
This concludes the proof of Theorem \ref{thm:main} in the case $\tau_1<\tau_2\leq \tau_1+10$.

\noindent{\bf Step 1.} We now assume that \eqref{eq:tau2geqtau1plus10conditionisthemaincase} holds, i.e., we consider the case $\tau_2\geq\tau_1+10$, and we begin by upgrading the estimate \eqref{eq:mainhighorderunweightedEMF:maintheorem:pm2} to energy-Morawetz estimates controlling in addition weighted derivatives of $\pmb\phi_s^{(p)}$. To this end, we combine the energy-Morawetz estimates \eqref{eq:mainhighorderunweightedEMF:maintheorem:pm2} and \eqref{eq:EMnearinfinity:highorderweightedderivatives:Teu:pm2} to obtain, for $s=\pm 2$, and for any $\reg\leq 14$ and $0<\de\leq \frac{1}{3}$,
\bea
\lab{eq:mainhighorderunweightedEMF:maintheorem:pm2:consequenceaddiingcontrolweightedderivativestau1tau2}
&&\sum_{p=0}^2\EMF_{\de}^{(\reg)}[\phis{p}](\tau_1,\tau_2) + \sum_{p=0}^2{\EMF}_{\de}[\pr^{\leq \reg}\psis{p}](\Iti)\nn\\
&&+\sum_{p=0}^1\int_{\Mntrap(\tau_1,\tau_2)}r^{-3+\de} \big|(r\nab)^{\leq 1} \dk^{\leq \reg}\phis{p}\big|^2
\nn\\
&\les&\sum_{p=0}^2\E^{(\reg)}[\phis{p}](\tau_1)+\sum_{p=0}^2\NN_\de^{(\reg)}[\phis{p}, \N_{W,s}^{(p)}](\tt_1, \tt_2)+\sum_{p=0}^1\int_{\MM(\tau_1,\tau_2)}r^{-1+\de}|\dk^{\leq \reg+1}\N_{T,s}^{(p)}|^2\nn\\
&&
+\sum_{p=0}^1\int_{\MM_{r\geq 10m}(\tau_1,\tau_2)} \Big(r^{-1+\de}|\pr^{\leq \reg}\N^{(p)}_{W,s}|+r^{-2+\de}|\pr^{\leq \reg+1}\N^{(p)}_{T,s}|\Big)|\pr^{\leq \reg}\phis{p}|\nn\\
&& +\sum_{p=0}^2\Big(\EMF[\psis{p}](\Iti)
+\EMF[\phis{p}](\tau_1, \tau_2)\Big).
\eea
Also, we have, for any $\de\in(0,\frac{1}{3}]$ and $\reg\leq 14$,
\beaa
&&\sum_{p=0}^1\int_{\MM_{r\geq 10m}(\tau_1,\tau_2)} \Big(r^{-1+\de}|\pr^{\leq \reg}\N^{(p)}_{W,s}|+r^{-2+\de}|\pr^{\leq \reg+1}\N^{(p)}_{T,s}|\Big)|\pr^{\leq \reg}\phis{p}|\nn\\
&\les& \sum_{p=0}^1\bigg(\int_{\Mntrap(\tau_1,\tau_2)}r^{-3+\de} \big|\dk^{\leq \reg}\phis{p}\big|^2\bigg)^{\frac{1}{2}}\nn\\
&&\times
\bigg(\int_{\MM(\tt_1, \tt_2)}\Big(r^{-1+\de}|\dk^{\leq {\reg}+1}\N_{T,s}^{(p)}|^2+r^{1+\de}|\dk^{\leq {\reg}}\N_{W,s}^{(p)}|^2\Big)
\bigg)^{\frac{1}{2}}.
\eeaa
Together with \eqref{eq:mainhighorderunweightedEMF:maintheorem:pm2:consequenceaddiingcontrolweightedderivativestau1tau2}, we infer, for $s=\pm 2$, and for any $\reg\leq 14$ and $0<\de\leq \frac{1}{3}$,
\bea\lab{eq:allderivativesrecoveredandonlylowerordertermonRHS:plus2caseandminus2caseatthesametime}
&&\sum_{p=0}^2\EMF_{\de}^{(\reg)}[\phis{p}](\tau_1,\tau_2) + \sum_{p=0}^2{\EMF}_{\de}[\pr^{\leq \reg}\psis{p}](\Iti)\nn\\
&\les&\sum_{p=0}^2\E^{(\reg)}[\phis{p}](\tau_1)+\sum_{p=0}^2\NN_\de^{(\reg)}[\phis{p}, \N_{W,s}^{(p)}](\tt_1, \tt_2)+\sum_{p=0}^1\int_{\MM(\tau_1,\tau_2)}r^{-1+\de}|\dk^{\leq \reg+1}\N_{T,s}^{(p)}|^2\nn\\
&& +\sum_{p=0}^2\Big(\EMF[\psis{p}](\Iti)
+\EMF[\phis{p}](\tau_1, \tau_2)\Big).
\eea

\noindent{\bf Step 2.} In this step, we remove the lower order terms which appear on the last line of the RHS of \eqref{eq:allderivativesrecoveredandonlylowerordertermonRHS:plus2caseandminus2caseatthesametime} by relying on the energy-Morawetz estimates in Kerr of Section  \ref{sec:energyMorwetzestimatesinKerr}. To this end, we first deduce from the scalarized wave equation for $\phi_{s,ij}^{(0)}$ on $(\MM, \g)$ a corresponding scalarized wave equation in Kerr.

Recall from Lemma \ref{lem:scalarizedTeukolskywavetransportsysteminKerrperturbation:Omi} that $\phi_{s,ij}^{(0)}$ satisfies 
\bea\lab{eq:scalarizedwaveequationforphis0:Kerrpert}
\widehat\square_\g(\phis{0})_{ij} -\frac{2}{\qs} \phiss{ij}{0} = {L_{s,ij}^{(0)}}+N_{W,s,ij}^{(0)}, \quad s=\pm2, \quad i,j=1,2,3, 
\eea
with
\beaa
\widehat\square_\g(\phis{0})_{ij} &=&\square_\g\phiss{ij}{0}- \widehat{S}(\phis{0})_{ij} -(\widehat{Q}\phis{0})_{ij},\\
\widehat{S}(\phis{0})_{ij} &=& S(\phis{0})_{ij} +\frac{4ia\cos\th}{|q|^2} \pr_\tau\phiss{ij}{0},\\
(\widehat{Q}\phis{0})_{ij} &=& (Q\phis{0})_{ij} 
-\frac{4ia\cos\th}{|q|^2}\big(M_{i\tau}^l \phiss{lj}{0}+M_{j\tau}^l \phiss{il}{0}\big),\\
S(\phis{0})_{ij} &=& 2M_{i}^{k\a}\pr_\a(\phi^{(0)}_{s,kj}) +2M_{j}^{k\a}\pr_\a(\phi^{(0)}_{s,ik}),\\
(Q\phis{0})_{ij} &=& (\Ddot^\a M_{i\a}^k)\phi^{(0)}_{s,kj}+(\Ddot^\a M_{j\a}^k)\phi^{(0)}_{s,ik} -M_{i\a}^kM_k^{l\a}\phi^{(0)}_{s,lj}-2M_{i\a}^kM_{j}^{l\a}\phi^{(0)}_{s,kl}-M_{j\a}^kM_k^{l\a}\phi^{(0)}_{s,il}.
\eeaa
Now, in view of \eqref{eq:assumptionsonregulartripletinperturbationsofKerr:0} and Lemma \ref{lemma:computationoftheMialphajinKerr}, we have
\beaa
\widehat{S}(\phis{0})_{ij} = \widehat{S}_K(\phis{0})_{ij}+\Ga_g\dk(\phi_s^{(0)}), \qquad (\widehat{Q}\phis{0})_{ij} = (\widehat{Q}_K\phis{0})_{ij}+\dk^{\leq 1}(\Ga_g)\phi_s^{(0)},
\eeaa
where $\widehat{S}_K$ and $\widehat{Q}_K$ correspond to $\widehat{S}$ and $\widehat{Q}$ in Kerr. Also, we have in view of \eqref{eq:comparisionbetweenscalarwaveoperatorinKerrandinKerrpert}
\beaa
\square_\g(\phi^{(0)}_{s,ij}) &=&  \square_{\gam}(\phi^{(0)}_{s,ij}) +\dk^{\leq 2}(\Ga_g\phi_s^{(0)}).
\eeaa
The above implies 
\beaa
\widehat\square_\g(\phis{0})_{ij} = \widehat{\square}_{\gam}(\phis{0})_{ij} +\dk^{\leq 2}(\Ga_g\phi_s^{(0)})
\eeaa
which together with \eqref{eq:scalarizedwaveequationforphis0:Kerrpert} yields
\bea\lab{eq:scalarizedwaveequationforphis0:threwquasilineartermRHS}
\widehat{\square}_{\gam}(\phis{0})_{ij} -\frac{2}{\qs} \phiss{ij}{0} = {L_{s,ij}^{(0)}}+N_{W,s,ij}^{(0)}+\dk^{\leq 2}(\Ga_g\phi_s^{(0)}), \quad s=\pm2, \quad i,j=1,2,3.
\eea

In order to get the analog of the Teukolsky wave equation for $\phiss{ij}{0}$ in Kerr, we still need to replace ${L_{s,ij}^{(0)}}$, appearing on the RHS \eqref{eq:scalarizedwaveequationforphis0:threwquasilineartermRHS}, by its Kerr value.  Recalling from Lemma \ref{lem:scalarizedTeukolskywavetransportsysteminKerrperturbation:Omi} that 
\beaa
{L_{s,ij}^{(0)}}&=& (2sr^{-3} +O(mr^{-4}))\phiss{ij}{1} + O(mr^{-3}){\Xcal_s}\phiss{ij}{0}+\sum_{k,l=1,2,3}O(mr^{-3})\phiss{kl}{0},
\eeaa
it suffices to replace $\phiss{ij}{1}$ by the corresponding expression in Kerr in terms of $\phiss{ij}{0}$ and first-order derivatives of $\phiss{ij}{0}$. Now, recalling from \eqref{eq:ScalarizedQuantitiesinTeuSystem:Kerrperturbation} that 
\beaa
\bsplit
& e_3\bigg(\frac{r\bar{q}}{q}\bigg({\frac{r^2}{|q|^2}}\bigg)^{-2}\phipluss{ij}{0} \bigg) - \frac{r\bar{q}}{q}\bigg({\frac{r^2}{|q|^2}}\bigg)^{-2}\Big(M_{i3}^k \phipluss{kj}{0} + M_{j3}^k \phipluss{ik}{0}\Big)\\
=&\frac{\bar{q}}{rq}\bigg({\frac{r^2}{|q|^2}}\bigg)^{-1}\phipluss{ij}{1}+N_{T,+2,ij}^{(0)}
\end{split}
\eeaa
and 
\beaa
\bsplit
& e_4\bigg(\frac{rq}{\bar{q}}\bigg(\frac{r^2}{|q|^2}\bigg)^{-2}\phiminuss{ij}{0}\bigg) - \frac{rq}{\bar{q}}\bigg(\frac{r^2}{|q|^2}\bigg)^{-2}\Big(M_{i4}^k \phiminuss{kj}{0} + M_{j4}^k \phiminuss{ik}{0}\Big)\\
=& \frac{q}{r\bar{q}}\left(\frac{r^2}{|q|^2}\right)^{-1}\frac{\De}{\qs}\phiminuss{ij}{1} +N_{T,-2,ij}^{(0)},
\end{split}
\eeaa
and relying on \eqref{eq:assumptionsonregulartripletinperturbationsofKerr:0} and Lemma \ref{lemma:computationoftheMialphajinKerr}, as well as the fact that 
\beaa
e_4=(e_4)_K+\Ga_g\dk, \qquad e_3=(e_3)_K+r\Ga_g\dk,
\eeaa
which follows immediately from \eqref{eq:relationsbetweennullframeandcoordinatesframe2:moreprecise:00}, we infer
\beaa
\bsplit
& (e_3)_K\bigg(\frac{r\bar{q}}{q}\bigg({\frac{r^2}{|q|^2}}\bigg)^{-2}\phipluss{ij}{0} \bigg) - \frac{r\bar{q}}{q}\bigg({\frac{r^2}{|q|^2}}\bigg)^{-2}\Big((M_K)_{i3}^k \phipluss{kj}{0} + (M_K)_{j3}^k \phipluss{ik}{0}\Big)\\
=&\frac{\bar{q}}{rq}\bigg({\frac{r^2}{|q|^2}}\bigg)^{-1}\phipluss{ij}{1}+N_{T,+2,ij}^{(0)}+r^2\Ga_g\dk^{\leq 1}(\phi_{+2}^{(0)})
\end{split}
\eeaa
and 
\beaa
\bsplit
& (e_4)_K\bigg(\frac{rq}{\bar{q}}\bigg(\frac{r^2}{|q|^2}\bigg)^{-2}\phiminuss{ij}{0}\bigg) - \frac{rq}{\bar{q}}\bigg(\frac{r^2}{|q|^2}\bigg)^{-2}\Big((M_K)_{i4}^k \phiminuss{kj}{0} + (M_K)_{j4}^k \phiminuss{ik}{0}\Big)\\
=& \frac{q}{r\bar{q}}\left(\frac{r^2}{|q|^2}\right)^{-1}\frac{\De}{\qs}\phiminuss{ij}{1} +N_{T,-2,ij}^{(0)} +r\Ga_g\dk^{\leq 1}(\phi_{-2}^{(0)}),
\end{split}
\eeaa
and hence, 
\beaa
L_{+2,ij}^{(0)} &=& (L_K)_{+2,ij}^{(0)} +O(r^{-2})N_{T,+2,ij}^{(0)}+\Ga_g\dk^{\leq 1}(\phi_{+2}^{(0)}),  \quad i,j=1,2,3,\\
L_{-2,ij}^{(0)} &=& (L_K)_{-2,ij}^{(0)} +O(\De^{-1})N_{T,-2,ij}^{(0)}+r\De^{-1}\Ga_g\dk^{\leq 1}(\phi_{-2}^{(0)}),  \quad i,j=1,2,3,
\eeaa
which together with \eqref{eq:scalarizedwaveequationforphis0:threwquasilineartermRHS} implies, for $i,j,=1,2,3$,
\bea\lab{eq:scalarizedwaveequationforphis0:threwquasilineartermRHS:bis:plus2case}
\widehat{\square}_{\gam}(\pmb\phi^{(0)}_{+2})_{ij} -\frac{2}{\qs} \phi^{(0)}_{+2,ij} = (L_K)_{+2,ij}^{(0)}+N_{W,+2,ij}^{(0)}+O(r^{-2})N_{T,+2,ij}^{(0)}+\dk^{\leq 2}(\Ga_g\phi_{+2}^{(0)}),
\eea
and, for $r\geq r_+(1+2\dred )$, 
\bea\lab{eq:scalarizedwaveequationforphis0:threwquasilineartermRHS:bis:minus2case}
\widehat{\square}_{\gam}(\pmb\phi^{(0)}_{-2})_{ij} -\frac{2}{\qs}\phi^{(0)}_{-2,ij} = (L_K)_{-2,ij}^{(0)}+N_{W,-2,ij}^{(0)}+O(\De^{-1})N_{T,-2,ij}^{(0)}+\frac{r^2}{\De}\dk^{\leq 2}(\Ga_g\phi_{-2}^{(0)}).
\eea

Next, using the regular triplet $(\Om_K)_i$, $i=1,2,3$, in Kerr introduced in Definition \ref{def:regulartripletinKerrOmii=123}, we define the horizontal tensors $\pmb\phi_s$, $s=\pm 2$, by 
\bea\lab{eq:definitionofphiplusminus2inKerrfromphi0plusminus2Kerrpert}
(\pmb\phi_{+2})_{ab}=|q|^4\phi_{+2,ij}^{(0)}(\Om_K^i)_a(\Om_K^j)_b, \qquad (\pmb\phi_{-2})_{ab}=\frac{|q|^4}{\De^2}\phi_{-2,ij}^{(0)}(\Om_K^i)_a(\Om_K^j)_b,
\eea
where the definition of $\pmb\phi_{-2}$ will only be used in $r\geq r_+(1+2\dred )$. Then, since $\phi_{s,ij}^{(0)}$ corresponds by construction to the scalarization of a tensor in $\sk_2(\mathbb{C})$ w.r.t. the regular triplet of Section  \ref{sec:regulartripletinperturbationsofKerr}, it thus satisfies the identities in the first item of Lemma \ref{lemma:backandforthbetweenhorizontaltensorsk2andscalars:complex}, and hence, in view of  the second item of Lemma \ref{lemma:backandforthbetweenhorizontaltensorsk2andscalars:complex}, this implies that $\pmb\phi_s$ as defined in \eqref{eq:definitionofphiplusminus2inKerrfromphi0plusminus2Kerrpert} satisfy $\pmb\phi_s\in\sk_2(\mathbb{C})$ in Kerr. Also, using Lemma \ref{lem:scalarizedTeukolskywavetransportsysteminKerr:Omi}, and the definition \eqref{eq:definitionofphiplusminus2inKerrfromphi0plusminus2Kerrpert}, we infer from \eqref{eq:scalarizedwaveequationforphis0:threwquasilineartermRHS:bis:plus2case} and \eqref{eq:scalarizedwaveequationforphis0:threwquasilineartermRHS:bis:minus2case}
 the following inhomogeneous analog of the Teukolsky equations in Kerr \eqref{eq:TeukolskyequationforAandAbintensorialforminKerr}, for $s=\pm 2$, 
\bea\lab{eq:TeukolskyequationforAandAbintensorialforminKerr:inhomogeneous} 
\nn\left(\squared_{2,K} -\frac{4ia\cos\th}{|q|^2}\nab_{\pr_t} -  \frac{s}{|q|^2}\right)\pmb\phi_s +\frac{2s}{|q|^2}(r-m)\nab_{(e_K)_3}\pmb\phi_s  -\frac{4sr}{|q|^2}\nab_{\pr_t}\pmb\phi_s\\
+\frac{4a\cos\th}{|q|^6}\Big(a\cos\th\big(|q|^2+6mr\big) - is\big((r-m)|q|^2+4mr^2\big)\Big)\pmb\phi_s &=& \N_s,  
\eea
where the inhomogeneous RHS $\N_s$, $s=\pm 2$, are given by  
\bsub\lab{eq:structureofinhomogenoustermsNplusminus2inTeukolskyforphisinKerr}
\bea
\N_{+2} &=& |q|^4\Big[\N_{W,+2}^{(0)}+O(r^{-2})\N_{T,+2}^{(0)}+\dk^{\leq 2}(\Ga_g\c\pmb\phi_{+2}^{(0)})\Big],\\
\N_{-2} &=& \frac{|q|^4}{\De^2}\Big[\N_{W,-2}^{(0)}+O(\De^{-1})\N_{T,-2}^{(0)}+\frac{r^2}{\De}\dk^{\leq 2}(\Ga_g\c\pmb\phi_{-2}^{(0)})\Big], 
\eea
\esub
with the definition for $\N_{-2}$ in \eqref{eq:structureofinhomogenoustermsNplusminus2inTeukolskyforphisinKerr} being used only for $r\geq r_+(1+2\dred)$. For, $r\leq r_+(1+2\dred)$, we rely instead on $\Ab$ which satisfies in view of \eqref{eq:waveequationpmbphip=0sminus2nodeginredshiftregion} in that region in $(\MM, \g)$, 
\beaa
\squared_2\Ab = 2\pr_r\left(\frac{\De}{|q|^2}\right)\nab_3\Ab
+ O(1)\big(\nab_{4}\Ab,\nab\Ab, \Ab\big)+\N_{\Ab}.
\eeaa
Then: 
\begin{itemize}
\item we scalarize this equation for $\Ab$ using the regular triplet of Section  \ref{sec:regulartripletinperturbationsofKerr} and Lemma \ref{lemma:formoffirstordertermsinscalarazationtensorialwaveeq},
\item similarly to \eqref{eq:definitionofphiplusminus2inKerrfromphi0plusminus2Kerrpert}, we define $\pmb\phi_{-2}\in\sk_2(\mathbb{C})$ in Kerr for $r\leq r_+(1+2\dred)$ by\footnote{The definitions of $\pmb\phi_{-2}$ in \eqref{eq:definitionofphiplusminus2inKerrfromphi0plusminus2Kerrpert} and \eqref{eq:definitionofphiplusminus2inKerrfromphi0plusminus2Kerrpert:bis} agree since $
\phi_{-2,ij}^{(0)}=q(\ov{q})^{-1}\De^2|q|^{-4}\Ab(\Om_i, \Om_j)$ in view of \eqref{eq:defintionphipm2p=0:Kerrperturbation}.}
\bea\lab{eq:definitionofphiplusminus2inKerrfromphi0plusminus2Kerrpert:bis}
(\pmb\phi_{-2})_{ab}=|q|^{-2}q^2\Ab(\Om_i, \Om_j)(\Om_K^i)_a(\Om_K^j)_b
\eea
which is motivated by the identity \eqref{eq:defintionoftensorialTeukolskyscalarspsipm2} in Kerr. 
\end{itemize}
With this definition, the inhomogeneous term $\N_{-2}$ appearing on the RHS of \eqref{eq:TeukolskyequationforAandAbintensorialforminKerr:inhomogeneous}, and given in $r\geq r_+(1+2\dred )$ by \eqref{eq:structureofinhomogenoustermsNplusminus2inTeukolskyforphisinKerr}, satisfies 
in $r\leq r_+(1+2\dred )$
\bea\lab{eq:structureofinhomogenoustermsNplusminus2inTeukolskyforphisinKerr:HHm2}
\N_{-2}=|q|^{-2}q^2\Big[\N_{\Ab}+\dk^{\leq 2}(\Ga_g\c\Ab)\Big].
\eea

Next, we introduce $\pmb\varphi_s\in\sk_2(\mathbb{C})$, $s=\pm 2$, in Kerr which satisfy the same initial condition as $\pmb\phi_s$ on $\tau=\tau_1$ and are the solutions of the following modification of \eqref{eq:TeukolskyequationforAandAbintensorialforminKerr:inhomogeneous}
\bea\lab{eq:TeukolskyequationforAandAbintensorialforminKerr:inhomogeneous:tilde} 
\nn\left(\squared_{2,K} -\frac{4ia\cos\th}{|q|^2}\nab_{\pr_t} -  \frac{s}{|q|^2}\right)\pmb\varphi_s +\frac{2s}{|q|^2}(r-m)\nab_{(e_K)_3}\pmb\varphi_s  -\frac{4sr}{|q|^2}\nab_{\pr_t}\pmb\varphi_s\\
+\frac{4a\cos\th}{|q|^6}\Big(a\cos\th\big(|q|^2+6mr\big) - is\big((r-m)|q|^2+4mr^2\big)\Big)\pmb\varphi_s &=& \widetilde{\N}_s,  
\eea
with 
\bea\lab{eq:definitiionofwidetildeNsextendingNsinKerr}
\widetilde{\N}_s=\widetilde\chi_{\tau_2}(\tau)\N_s, \qquad s=\pm 2,
\eea
where $\widetilde\chi_{\tau_2}$ is a smooth cut-off function satisfying $\widetilde\chi_{\tau_2}(\tau)=1$ for $\tau\leq \tau_2-1$ and $\widetilde\chi_{\tau_2}=0$ for $\tau\geq \tau_2$ so that we have by causality
\bea\lab{eq:causalityrelationbetweenpmbphisandpmbpsisinKerr}
\pmb\varphi_s=\pmb\phi_s, \quad s=\pm 2, \quad \textrm{for}\quad \tau\leq\tau_2-1.
\eea
Now, in view of \eqref{eq:TeukolskyequationforAandAbintensorialforminKerr:inhomogeneous:tilde}, $\pmb\varphi_s$ satisfies \eqref{eq:TeukolskyequationforAandAbintensorialforminKerr:inhomogenouscase}, and we may thus apply the weak Morawetz estimates in Kerr of Theorem \ref{cor:weakMorawetzforTeukolskyfromMillet:bis}, with $\tau_0=\tau_1$,  which implies, for any $\de\in(0,\frac{1}{3}]$,
\bea\lab{eq:weakMorawetzestimatesforphiminus2:consequenceThoerem}
\nn\int_{\MM(\tau_1,+\infty)}r^{-3+\de}|\dk^{\leq 3}\pmb\varphi_{-2}|^2 &\les& \E^{({13})}[\pmb\varphi_{-2}](\tau_1)+\E^{({11})}[r\nab_{\pr_r}(r\pmb\varphi_{-2})](\tau_1)\\
&&+\int_{\MM(\tau_1,+\infty)}r^{1+\de}|\dk^{\leq {13}}\widetilde{\N}_{-2}|^2,
\eea
and with $(\tau_0,\tau_1)=(\tau_1,\tau_2-1)$, which implies, for any $\de\in(0,\frac{1}{3}]$ and any $R_0\geq 10m$, 
\begin{align}\lab{eq:weakMorawetzestimatesforphiplus2:consequenceThoerem}
&\int_{\MM(\tau_1,\tau_2-2)}r^{-11+\frac{\de}{2}}|\dk^{\leq 3}\pmb\varphi_{+2}|^2\nn\\
 \les{}& R_0^{1+\frac{\de}{2}}\E^{({9})}[r^{-4}\pmb\varphi_{+2}](\tau_1)+\int_{\MM(\tau_1,\tau_2-1)}r^{-7+\de}|\dk^{\leq {10}}\widetilde{\N}_{+2}|^2 \nn\\
&+R_0^{-\frac{\de}{2}}\bigg({\EMF}^{(11)}_{\de}[r^{-4}\pmb\varphi_{+2}](\tau_1,\tau_2-1)+\int_{\MM(\tau_1,\tau_2-1)} r^{-3+\de}|\dk^{\leq 11} \nab_{3,\K}(r^{-3}\pmb\varphi_{+2})|^2\bigg).
\end{align}

We first derive a  weak Morawetz estimate for $\phiminus{0}$. In view of  \eqref{eq:definitionofphiplusminus2inKerrfromphi0plusminus2Kerrpert}, \eqref{eq:definitionofphiplusminus2inKerrfromphi0plusminus2Kerrpert:bis} and \eqref{eq:causalityrelationbetweenpmbphisandpmbpsisinKerr}, we have
\beaa
\int_{\MM(\tau_1, \tau_2-1)}r^{-3+\de}|\dk^{\leq 3}\pmb\phi_{-2}^{(0)}|^2 \les   \int_{\MM(\tau_1, \tau_2-1)}r^{-3+\de}|\dk^{\leq 3}\pmb\phi_{-2}|^2=\int_{\MM(\tau_1, \tau_2-1)}r^{-3+\de}|\dk^{\leq 3}\pmb\varphi_{-2}|^2
\eeaa
and\footnote{Recall, from our convention introduced in Section \ref{sec:smallnesconstants}, that we do not need to track the dependence of $\les$ on $\dred $ in this estimate.}
\beaa
\E^{({13})}[\pmb\varphi_{-2}](\tau_1)+\E^{({11})}[r\nab_{\pr_r}(r\pmb\varphi_{-2})](\tau_1) &=& \E^{({13})}[\pmb\phi_{-2}](\tau_1)+\E^{({11})}[r\nab_{\pr_r}(r\pmb\phi_{-2})](\tau_1)\\
&\les& \E^{({13})}[\pmb\phi_{-2}^{(0)}](\tau_1)+\E^{({11})}[r\nab_{\pr_r}(r\pmb\phi_{-2}^{(0)})](\tau_1)\\
&&+\E^{({13})}_{r\leq r_+(1+\dred )}[\Ab](\tau_1).
\eeaa
Also, we have in view of \eqref{eq:definitiionofwidetildeNsextendingNsinKerr} and the properties of the cut-off function $\widetilde\chi_{\tau_2}$ 
\beaa
\int_{\MM(\tau_1,+\infty)}r^{1+\de}|\dk^{\leq {13}}\widetilde{\N}_{-2}|^2\les \int_{\MM(\tau_1,\tau_2)}r^{1+\de}|\dk^{\leq {13}}\N_{-2}|^2.
\eeaa 
Plugging the above estimates in \eqref{eq:weakMorawetzestimatesforphiminus2:consequenceThoerem}, we infer
\beaa
\nn\int_{\MM(\tau_1, \tau_2-1)}r^{-3+\de}|\dk^{\leq 3}\pmb\phi_{-2}^{(0)}|^2 &\les& \E^{({13})}[\pmb\phi_{-2}^{(0)}](\tau_1)+\E^{({11})}[r\nab_{\pr_r}(r\pmb\phi_{-2}^{(0)})](\tau_1)\\
&&+\E^{({13})}_{r\leq r_+(1+\dred )}[\Ab](\tau_1)+\int_{\MM(\tau_1,\tau_2)}r^{1+\de}|\dk^{\leq {13}}\N_{-2}|^2.
\eeaa
Finally, together with \eqref{eq:structureofinhomogenoustermsNplusminus2inTeukolskyforphisinKerr} and \eqref{eq:structureofinhomogenoustermsNplusminus2inTeukolskyforphisinKerr:HHm2}, the first formula in \eqref{eq:relationsbetweennullframeandcoordinatesframe:1} which states
\beaa
e_4 = \left(1+O(mr^{-1})\right)\pr_r + O(m^2r^{-2})\pr_\tau+O(\ep r^{-2})\pr_{x^a},
\eeaa
and the fact that, in view of \eqref{def:TensorialTeuScalars:wavesystem:Kerrperturbation},  
\beaa
\E^{({11})}[r\nab_{\pr_r}(r\pmb\phi_{-2}^{(0)})](\tau_1) &\les& \E^{({11})}[\pmb\phi_{-2}^{(1)}](\tau_1)+\E^{(12)}[\pmb\phi_{-2}^{(0)}](\tau_1)+\int_{\Si(\tau_1)}|\dk^{\leq {12}}\N_{T,-2}^{(0)}|^2,
\eeaa
we deduce the following weak Morawetz estimate for $\phiminus{0}$  in perturbations of Kerr\footnote{Again, recall, from our convention introduced in Section  \ref{sec:smallnesconstants}, that we do not need to track the dependence of $\les$ on $\dred $ in the two first lines on the RHS of \eqref{eq:weakMorawetzestimatesforphiminus2:consequenceThoerem:1}. Additionally, recall that $\ep\ll\dred $ in view of \eqref{eq:constraintsonthemainsmallconstantsepanddelta} which will allow us, below, to absorb the nonlinear terms on the last line of \eqref{eq:weakMorawetzestimatesforphiminus2:consequenceThoerem:1} regardless of the dependence on $\dred $ so that we do not track this dependence in the last line either.} 
\bea\lab{eq:weakMorawetzestimatesforphiminus2:consequenceThoerem:1}
\nn&&\int_{\MM(\tau_1, \tau_2-1)}r^{-3+\de}|\dk^{\leq 3}\pmb\phi_{-2}^{(0)}|^2\\ 
\nn&\les& \E^{({13})}[\pmb\phi_{-2}^{(0)}](\tau_1)+ \E^{({11})}[\pmb\phi_{-2}^{(1)}](\tau_1)+\int_{\Si(\tau_1)}|\dk^{\leq {12}}\N_{T,-2}^{(0)}|^2 +\E^{({13})}_{r\leq r_+(1+\dred )}[\Ab](\tau_1)\\
\nn&+&\int_{\MM(\tau_1,\tau_2)}r^{1+\de}|\dk^{\leq {13}}\N_{W,-2}^{(0)}|^2+\int_{\MM(\tau_1,\tau_2)}r^{-3+\de}|\dk^{\leq {13}}\N_{T,-2}^{(0)}|^2+\int_{\MM_{r\leq r_+(1+\dred )}(\tau_1,\tau_2)}|\dk^{\leq {13}}\N_{\Ab}|^2\\
&&+\ep^2\sup_{\tau\in[\tau_1, \tau_2]}\E^{({14})}[\pmb\phi_{-2}^{(0)}](\tau)+\ep^2\sup_{\tau\in[\tau_1, \tau_2]}\E^{({14})}_{r\leq r_+(1+\dred )}[\Ab](\tau).
\eea

Next, we derive a  weak Morawetz estimate for $\phiplus{0}$. In view of \eqref{eq:definitionofphiplusminus2inKerrfromphi0plusminus2Kerrpert}, \eqref{eq:definitiionofwidetildeNsextendingNsinKerr}, \eqref{eq:causalityrelationbetweenpmbphisandpmbpsisinKerr} and the properties of the cut-off function $\widetilde\chi_{\tau_2}$, the estimate \eqref{eq:weakMorawetzestimatesforphiplus2:consequenceThoerem} implies
\beaa
&&\int_{\MM(\tau_1,\tau_2-2)}r^{-3+\frac{\de}{2}}|\dk^{\leq 3}\phiplus{0}|^2\nn\\
 &\les& R_0^{1+\frac{\de}{2}}\E^{({9})}[\phiplus{0}](\tau_1)+\int_{\MM(\tau_1,\tau_2-1)}r^{-7+\de}|\dk^{\leq {10}}{\N}_{+2}|^2 \nn\\
&&+R_0^{-\frac{\de}{2}}\bigg({\EMF}^{(11)}_{\de}[\phiplus{0}](\tau_1,\tau_2-1)+\int_{\MM(\tau_1,\tau_2-1)} r^{-3+\de}|\dk^{\leq 11} \nab_{3,\K}(r\phiplus{0})|^2\bigg).
\eeaa
Together with the second formula in \eqref{eq:relationsbetweennullframeandcoordinatesframe2:moreprecise:00} which implies
\beaa
e_{3,\K} =e_3 +r\Ga_b \pr_r + r\Ga_g \pr_{\tau} + \Ga_b \pr_{x^a}=e_3 + r\Ga_g \dk,
\eeaa
and the fact that, in view of \eqref{def:TensorialTeuScalars:wavesystem:Kerrperturbation},  
\beaa
\int_{\MM(\tau_1,\tau_2-1)} r^{-3+\de}|\dk^{\leq 11}\nab_3(r\phiplus{0})|^2&\les& 
\sum_{p=0}^1\int_{\MM(\tau_1,\tau_2-1)} r^{-5+\de}|\dk^{\leq 11} \phiplus{p}|^2\nn\\
&&+\int_{\MM(\tau_1,\tau_2-1)} r^{-3+\de}|\dk^{\leq 11} \N_{T, +2}^{(0)}|^2,
\eeaa
we deduce the following weak Morawetz estimates for $\phiplus{0}$  in perturbations of Kerr
\bea\lab{eq:weakMorawetzestimatesforphiplus2:consequenceThoerem:1}
&&\int_{\MM(\tau_1,\tau_2-2)}r^{-3+\frac{\de}{2}}|\dk^{\leq 3}\phiplus{0}|^2\nn\\
 &\les& R_0^{1+\frac{\de}{2}}\E^{({9})}[\phiplus{0}](\tau_1)+\int_{\MM(\tau_1,\tau_2-1)}\Big(r^{-7+\de}|\dk^{\leq {10}}{\N}_{+2}|^2+r^{-3+\de}|\dk^{\leq 11} \N_{T, +2}^{(0)}|^2\Big) \nn\\
&&+R_0^{-\frac{\de}{2}}\bigg({\EMF}^{(11)}_{\de}[\phiplus{0}](\tau_1,\tau_2-1)+\int_{\MM(\tau_1,\tau_2-1)} r^{-5+\de}|\dk^{\leq 11} \phiplus{1}|^2\bigg)\nn\\
&\les& R_0^{1+\frac{\de}{2}}\E^{({9})}[\phiplus{0}](\tau_1)+\int_{\MM(\tau_1,\tau_2-1)}\Big(r^{1+\de}|\dk^{\leq {10}}{\N}_{W,+2}^{(0)}|^2+r^{-3+\de}|\dk^{\leq 11} \N_{T, +2}^{(0)}|^2\Big) \nn\\
&&+\Big(R_0^{-\frac{\de}{2}}+\ep^2\Big){\EMF}^{(11)}_{\de}[\phiplus{0}](\tau_1,\tau_2-1)+R_0^{-\frac{\de}{2}}\int_{\MM(\tau_1,\tau_2-1)} r^{-5+\de}|\dk^{\leq 11} \phiplus{1}|^2,
\eea
where we have used \eqref{eq:structureofinhomogenoustermsNplusminus2inTeukolskyforphisinKerr} in the last step.

Next, notice that  we have, for any $R<+\infty$,
\beaa
\sum_{p=0}^2\EMF_{r\leq R}[\pmb\phi_{+2}^{(p)}](\tau_1, \tau_2-2)\les_{R}\int_{\MM(\tau_1, \tau_2-2)}r^{-3}|\dk^{\leq 3}\pmb\phi_{+2}^{(0)}|^2
+\sum_{p=0}^2\int_{\MM_{r\leq R+m}(\tau_1, \tau_2-2)}|\N_{W,+2}^{(p)}|^2,\\
\sum_{p=0}^2\EMF_{r\leq R}[\pmb\phi_{-2}^{(p)}](\tau_1, \tau_2-1)\les_{R} \int_{\MM(\tau_1, \tau_2-1)}r^{-3}|\dk^{\leq 3}\pmb\phi_{-2}^{(0)}|^2+\sum_{p=0}^2\int_{\MM_{r\leq R+m}(\tau_1, \tau_2-1)}|\N_{W,-2}^{(p)}|^2,
\eeaa
where
\begin{itemize} 
\item the control of the Morawetz norm is immediate, 
\item the control of the flux norm on $\AA$ follows from redshift estimates and the control of Morawetz, 
\item and the control of the energy norm follows from the mean value argument which yields the control of the energy for $\tau=\tau_n$ with $\tau_n\in[n,n+1)$ for any $[n,n+1)\subset(\tau_1, \tau_2)$ and local energy estimates. 
\end{itemize}
In view of the above and \eqref{eq:weakMorawetzestimatesforphiminus2:consequenceThoerem:1} \eqref{eq:weakMorawetzestimatesforphiplus2:consequenceThoerem:1}, we may apply Proposition \ref{prop:EnergyMorawetznearinfinitytensorialTeuk} with $\reg=0$ to deduce
\beaa
\nn&&\sum_{p=0}^2\EMF[\pmb\phi_{+2}^{(p)}](\tau_1, \tau_2-2)\\ 
\nn&\les& R_0^{1+\de}\sum_{p=0}^2\E^{({9})}[\pmb\phi_{+2}^{(p)}](\tau_1)
+\sum_{p=0}^2\int_{\MM(\tau_1,\tau_2-1)}r^{1+\de}|\dk^{\leq {10}}{\N}_{W,+2}^{(p)}|^2\\
\nn&&+\sum_{p=0}^1\int_{\MM(\tau_1,\tau_2-1)}r^{-1+\de}|\dk^{\leq 11} \N_{T, +2}^{(p)}|^2
+\Big(R_0^{-\frac{\de}{2}}+\ep^2\Big)\sum_{p=0}^1{\EMF}^{(11)}_{\de}[\phiplus{p}](\tau_1,\tau_2-1)
\eeaa
and 
\beaa
\nn&&\sum_{p=0}^2\EMF[\pmb\phi_{-2}^{(p)}](\tau_1, \tau_2-1)\\ 
\nn&\les& \sum_{p=0}^2\E^{({13})}[\pmb\phi_{-2}^{(p)}](\tau_1)+\int_{\Si(\tau_1)}|\dk^{\leq {12}}\N_{T,-2}^{(0)}|^2 +\E^{({13})}_{r\leq r_+(1+\dred )}[\Ab](\tau_1)\\
\nn&+&\sum_{p=0}^2\int_{\MM(\tau_1,\tau_2)}r^{1+\de}|\dk^{\leq {13}}\N_{W,-2}^{(p)}|^2+\sum_{p=0}^1\int_{\MM(\tau_1,\tau_2)}r^{-1+\de}|\dk^{\leq {13}}\N_{T,-2}^{(p)}|^2\\
&+&\int_{\MM_{r\leq r_+(1+\dred )}(\tau_1,\tau_2)}|\dk^{\leq {13}}\N_{\Ab}|^2+\ep^2\sup_{\tau\in[\tau_1, \tau_2]}\E^{({14})}[\pmb\phi_{-2}^{(0)}](\tau)+\ep^2\sup_{\tau\in[\tau_1, \tau_2]}\E^{({14})}_{r\leq r_+(1+\dred )}[\Ab](\tau).
\eeaa
Together with \eqref{eq:causlityrelationsforwidetildepsi1}, this implies for the scalars $\psi^{(p)}_{s,ij}$
\beaa
\nn&&\sum_{p=0}^2\EMF[\psiplus{p}](\tau_1+1, \tau_2-3)+
\sum_{p=0}^2\EMF[\pmb\phi_{+2}^{(p)}](\tau_1, \tau_2-2)\\ 
\nn&\les& R_0^{1+\de}\sum_{p=0}^2\E^{({9})}[\pmb\phi_{+2}^{(p)}](\tau_1)
+\sum_{p=0}^2\int_{\MM(\tau_1,\tau_2-1)}r^{1+\de}|\dk^{\leq {10}}{\N}_{W,+2}^{(p)}|^2\\
\nn&&+\sum_{p=0}^1\int_{\MM(\tau_1,\tau_2-1)}r^{-1+\de}|\dk^{\leq 11} \N_{T, +2}^{(p)}|^2
+\Big(R_0^{-\frac{\de}{2}}+\ep^2\Big)\sum_{p=0}^1{\EMF}^{(11)}_{\de}[\phiplus{p}](\tau_1,\tau_2-1)
\eeaa
and
\beaa
\nn&&\sum_{p=0}^2\EMF[\psiminus{p}](\tau_1+1, \tau_2-3)
+\sum_{p=0}^2\EMF[\pmb\phi_{-2}^{(p)}](\tau_1, \tau_2-1)\\ 
\nn&\les& \sum_{p=0}^2\E^{({13})}[\pmb\phi_{-2}^{(p)}](\tau_1)+\int_{\Si(\tau_1)}|\dk^{\leq {12}}\N_{T,-2}^{(0)}|^2 +\E^{({13})}_{r\leq r_+(1+\dred )}[\Ab](\tau_1)\\
\nn&&+\sum_{p=0}^2\int_{\MM(\tau_1,\tau_2)}r^{1+\de}|\dk^{\leq {13}}\N_{W,-2}^{(p)}|^2+\sum_{p=0}^1\int_{\MM(\tau_1,\tau_2)}r^{-1+\de}|\dk^{\leq {13}}\N_{T,-2}^{(0)}|^2\\
&+&\int_{\MM_{r\leq r_+(1+\dred )}(\tau_1,\tau_2)}|\dk^{\leq {13}}\N_{\Ab}|^2+\ep^2\sup_{\tau\in[\tau_1, \tau_2]}\E^{({14})}[\pmb\phi_{-2}^{(0)}](\tau)+\ep^2\sup_{\tau\in[\tau_1, \tau_2]}\E^{({14})}_{r\leq r_+(1+\dred )}[\Ab](\tau).
\eeaa
Together with the local energy estimates \eqref{eq:localenergyestimateforpsisijponSigmatau2minus3tau2:partialderivatives} with $\reg=0$, the one in  \eqref{eq:localenergyestimate:future:bis} applied with $\reg=0$, $\tau_0=\tau_2-2$,  $q=2$ and $\psi=(\phiplus{p})_{p=0,1,2}$, the one in \eqref{eq:localenergyestimate:future:bis} applied with $\reg=0$, $\tau_0=\tau_2-1$, $q=1$ and $\psi=(\phiminus{p})_{p=0,1,2}$, the one in  \eqref{eq:localenergyforpsiontau1minus1tau1plus1} with $\reg=0$, and applying Theorem \ref{prop:weakMorawetzfortensorialwaveeqfromDRSR} with $\tau_0=\tau_2$, we infer 
\bea\lab{eq:weakMorawetzestimatesforphiplus2:consequenceThoerem:2}
\nn&&\sum_{p=0}^2\EMF[\psiplus{p}](\Iti)
+\sum_{p=0}^2\EMF[\pmb\phi_{+2}^{(p)}](\tau_1, \tau_2)\\ 
\nn&\les& R_0^{1+\de}\sum_{p=0}^2\E^{({9})}[\pmb\phi_{+2}^{(p)}](\tau_1)
+\sum_{p=0}^1\int_{\MM(\tau_1,\tau_2-1)}r^{-1+\de}|\dk^{\leq 11}\N_{T, +2}^{(p)}|^2\\
&&
+\sum_{p=0}^2\int_{\MM(\tau_1,\tau_2)}r^{1+\de}|\dk^{\leq {10}}{\N}_{W,+2}^{(p)}|^2+\Big(R_0^{-\frac{\de}{2}}+\ep^2\Big)\sum_{p=0}^1{\EMF}^{(11)}_{\de}[\phiplus{p}](\tau_1,\tau_2-1)
\eea
and  
\bea\lab{eq:weakMorawetzestimatesforphiminus2:consequenceThoerem:2}
\nn&&\sum_{p=0}^2\EMF[\psiminus{p}](\Iti)
+\sum_{p=0}^2\EMF[\pmb\phi_{-2}^{(p)}](\tau_1, \tau_2)\\ 
\nn&\les& \sum_{p=0}^2\E^{({13})}[\pmb\phi_{-2}^{(p)}](\tau_1)+\int_{\Si(\tau_1)}|\dk^{\leq 12}\N_{T,-2}^{(0)}|^2 +\E^{({13})}_{r\leq r_+(1+\dred )}[\Ab](\tau_1)\\
\nn&&+\sum_{p=0}^2\int_{\MM(\tau_1,\tau_2)}r^{1+\de}|\dk^{\leq {13}}\N_{W,-2}^{(p)}|^2+\sum_{p=0}^1\int_{\MM(\tau_1,\tau_2)}r^{-1+\de}|\dk^{\leq {13}}\N_{T,-2}^{(p)}|^2\\
\nn&&+\int_{\MM_{r\leq r_+(1+\dred )}(\tau_1,\tau_2)}|\dk^{\leq {13}}\N_{\Ab}|^2\\
&&+\ep^2\sup_{\tau\in[\tau_1, \tau_2]}\E^{({14})}[\pmb\phi_{-2}^{(0)}](\tau)+\ep^2\sup_{\tau\in[\tau_1, \tau_2]}\E^{({14})}_{r\leq r_+(1+\dred )}[\Ab](\tau).
\eea

Finally, using \eqref{eq:weakMorawetzestimatesforphiplus2:consequenceThoerem:2} and \eqref{eq:weakMorawetzestimatesforphiminus2:consequenceThoerem:2} to control the lower order terms which appear on the last line of the RHS of \eqref{eq:allderivativesrecoveredandonlylowerordertermonRHS:plus2caseandminus2caseatthesametime}, we infer
\bea\lab{eq:allderivativesrecoveredevenlowerordertermonRHS:plus2case}
\nn&& \sum_{p=0}^2\EMF^{({11})}_{\de}[\pmb\phi_{+2}^{(p)}](\tau_1, \tau_2)\\ 
\nn&\les& R_0^{1+\de}\sum_{p=0}^2\E^{({11})}[\pmb\phi_{+2}^{(p)}](\tau_1)+\sum_{p=0}^1\int_{\MM(\tt_1, \tt_2)}r^{-1+\de}|\dk^{\leq {12}}\N_{T,+2}^{(p)}|^2\\
&&+\sum_{p=0}^2\NN^{({11})}_\de[\pmb\phi_{+2}^{(p)}, \N_{W,+2}^{(p)}](\tt_1, \tt_2)
+\Big(R_0^{-\frac{\de}{2}}+\ep^2\Big)\sum_{p=0}^1{\EMF}^{(11)}_{\de}[\phiplus{p}](\tau_1,\tau_2-1)
\eea
and
\bea\lab{eq:allderivativesrecoveredevenlowerordertermonRHS:minus2case}
\nn \sum_{p=0}^2\EMF^{({14})}_{\de}[\pmb\phi_{-2}^{(p)}](\tau_1, \tau_2)&\les& \sum_{p=0}^2\E^{({14})}[\pmb\phi_{-2}^{(p)}](\tau_1)+\E^{({13})}_{r\leq r_+(1+\dred )}[\Ab](\tau_1)+\int_{\Si(\tau_1)}|\dk^{\leq {12}}\N_{T,-2}^{(0)}|^2\\
\nn&&+\sum_{p=0}^2\NN_\de^{({14})}[\pmb\phi_{-2}^{(p)}, \N_{W,-2}^{(p)}](\tt_1, \tt_2)+\int_{\MM_{r\leq r_+(1+\dred )}(\tau_1,\tau_2)}|\dk^{\leq {13}}\N_{\Ab}|^2\\
\nn&&+\sum_{p=0}^1\int_{\MM(\tt_1, \tt_2)}r^{-1+\de}|\dk^{\leq {15}}\N_{T,-2}^{(p)}|^2\\
&&+\ep^2\sup_{\tau\in[\tau_1, \tau_2]}\E^{({14})}[\pmb\phi_{-2}^{(0)}](\tau)
+\ep^2\sup_{\tau\in[\tau_1, \tau_2]}\E^{({14})}_{r\leq r_+(1+\dred )}[\Ab](\tau),
\eea
where we used in the derivation of \eqref{eq:allderivativesrecoveredevenlowerordertermonRHS:plus2case} \eqref{eq:allderivativesrecoveredevenlowerordertermonRHS:minus2case} the fact that, for $s=\pm 2$ and $\reg\leq 14$,
\beaa
\sum_{p=0}^2\int_{\MM(\tau_1,\tau_2)}r^{1+\de}|\dk^{\leq \reg}{\N}_{W,s}^{(p)}|^2
\les\sum_{p=0}^2\NN^{(\reg)}_\de[\pmb\phi_{s}^{(p)}, \N_{W,s}^{(p)}](\tt_1, \tt_2)
\eeaa
in view of the definition of $\mathcal{N}^{(\reg)}_\de[\c, \c](\tau_1, \tau_2)$ in \eqref{eq:defmathcalNpsif}. Also, combining \eqref{eq:allderivativesrecoveredevenlowerordertermonRHS:minus2case} with the redshift estimates of Corollary \ref{cor:redshift:Ab:highregularity} implies 
\bea\lab{eq:allderivativesrecoveredevenlowerordertermonRHS:minus2case:withredshift}
\nn &&\sum_{p=0}^2\EMF^{({14})}_{\de}[\pmb\phi_{-2}^{(p)}](\tau_1, \tau_2)+\sum_{p=0}^2\EMF^{({14})}_{r\leq r_+(1+\dred)}[\nab_4^p\Ab](\tau_1, \tau_2)\\
\nn&\les& \sum_{p=0}^2\E^{({14})}[\pmb\phi_{-2}^{(p)}](\tau_1)+\sum_{p=0}^2\E_{r\leq r_+(1+2\dred )}^{({14})}[\nab_4^p\Ab](\tau_1)+\int_{\Si(\tau_1)}|\dk^{\leq 12}\N_{T,-2}^{(0)}|^2\\
\nn&&+\sum_{p=0}^2\NN_\de^{({14})}[\pmb\phi_{-2}^{(p)}, \N_{W,-2}^{(p)}](\tt_1, \tt_2)+\sum_{p=0}^1\int_{\MM(\tt_1, \tt_2)}r^{-1+\de}|\dk^{\leq {15}}\N_{T,-2}^{(p)}|^2\\
\nn&&+\sum_{p=0}^2\int_{\MM_{r\leq r_+(1+2\dred)}(\tau_1, \tau_2)}|\dk^{\leq {14}}\N_{\nab_4^p\Ab}|^2\\
&&+\ep^2\sup_{\tau\in[\tau_1, \tau_2]}\E^{({14})}[\pmb\phi_{-2}^{(0)}](\tau)
+\ep^2\sup_{\tau\in[\tau_1, \tau_2]}\E^{({14})}_{r\leq r_+(1+\dred )}[\Ab](\tau).
\eea
Now, for $R_0$ large enough and $\ep$ small enough, we may absorb the last terms with $\ep^2$ coefficients on the RHS of \eqref{eq:allderivativesrecoveredevenlowerordertermonRHS:minus2case:withredshift} and the last term with $R_0^{-\frac{\de}{2}}+\ep^2$ coefficient on the RHS of \eqref{eq:allderivativesrecoveredevenlowerordertermonRHS:plus2case} which yields 
\beaa
\nn&&\sum_{p=0}^2\EMF^{({11})}_{\de}[\pmb\phi_{+2}^{(p)}](\tau_1, \tau_2)\\ 
\nn&\les&  \sum_{p=0}^2\E^{({11})}[\pmb\phi_{+2}^{(p)}](\tau_1)
+\sum_{p=0}^1\int_{\MM(\tt_1, \tt_2)}r^{-1+\de}|\dk^{\leq {12}}\N_{T,+2}^{(p)}|^2+\sum_{p=0}^2\NN^{({11})}_\de[\pmb\phi_{+2}^{(p)}, \N_{W,+2}^{(p)}](\tt_1, \tt_2)
\eeaa
and
\beaa
\nn &&\sum_{p=0}^2\EMF^{({14})}_{\de}[\pmb\phi_{-2}^{(p)}](\tau_1, \tau_2)+\sum_{p=0}^2\EMF^{({14})}_{r\leq r_+(1+\dred)}[\nab_4^p\Ab](\tau_1, \tau_2)\\
\nn&\les& \sum_{p=0}^2\E^{({14})}[\pmb\phi_{-2}^{(p)}](\tau_1)+\sum_{p=0}^2\E_{r\leq r_+(1+2\dred )}^{({14})}[\nab_4^p\Ab](\tau_1)+\int_{\Si(\tau_1)}|\dk^{\leq 12}\N_{T,-2}^{(0)}|^2\\
\nn&&+\sum_{p=0}^2\NN_\de^{({14})}[\pmb\phi_{-2}^{(p)}, \N_{W,-2}^{(p)}](\tt_1, \tt_2)+\sum_{p=0}^1\int_{\MM(\tt_1, \tt_2)}r^{-1+\de}|\dk^{\leq {15}}\N_{T,-2}^{(p)}|^2\\
&&+\sum_{p=0}^2\int_{\MM_{r\leq r_+(1+2\dred)}(\tau_1, \tau_2)}|\dk^{\leq {14}}\N_{\nab_4^p\Ab}|^2,
\eeaa
as stated. This concludes the proof of Theorem \ref{thm:main}.


\subsection{Structure of the rest of the paper}


In Section \ref{sect:microlocalenergyMorawetztensorialwaveequation}, we prove global energy-Morawetz estimates for a coupled system of scalar wave equations. Next, relying on the results of Section \ref{sect:microlocalenergyMorawetztensorialwaveequation}, we prove Theorem \ref{th:main:intermediary}  in Section \ref{sec:proofofth:main:intermediary}. Then, we prove Proposition \ref{prop:EnergyMorawetznearinfinitytensorialTeuk} in Section \ref{sec:proofofprop:EnergyMorawetznearinfinitytensorialTeuk}. Finally, Theorem \ref{cor:weakMorawetzforTeukolskyfromMillet:bis} is proved in Section \ref{sec:MorawetzestimatesforTeukolskys=pm2inKerrfromMillet}.


\section{Global energy-Morawetz estimates for a system of scalar wave equations}
\lab{sect:microlocalenergyMorawetztensorialwaveequation}


In this section, we consider a system of scalar wave equations for complex-valued scalars $\psi_{ij}$
\bea
\lab{eq:ScalarizedWaveeq:general:Kerrpert} 
\big({\square}_{\g}-D_0|q|^{-2}\big)\psi_{ij} =\widehat{F}_{ij}&:=&\chi_{\tau_1, \tau_2}\big(\widehat{S}(\psi)_{ij} +(\widehat{Q}\psi)_{ij})+(1-\chi_{\tau_1, \tau_2})\big( \widehat{S}_K(\psi)_{ij} +(\widehat{Q}_K\psi)_{ij}\big)\nn\\
&&
+(1-\chi_{\tau_1, \tau_2})f_{D_0}\psi_{ij}+F_{ij} ,\quad i,j=1,2,3, \quad\textrm{on}\quad\MM,
\eea
where $(\MM, \g)$ satisfy the assumptions of Sections \ref{subsect:assumps:perturbednullframe}, \ref{subsubsect:assumps:perturbedmetric} and \ref{sec:regulartripletinperturbationsofKerr} and $\g=\gam$ for $\tau\in\Reals\setminus(\tau_1, \tau_2)$, where the cut-off $\chi_{\tau_1, \tau_2}$ satisfies \eqref{eq:propertieschitoextendmetricg}, where  $D_0>0$ is a constant, where
\bea
f_{D_0}:=\frac{4-D_0}{|q|^2}- \frac{4a^2\cos^2\th(|q|^2+6mr)}{|q|^6},
\eea
where
\bsub
\label{hatSandV:generalwave:Kerrpert}
\bea
\widehat{S}(\psi)_{ij} &=& {S}(\psi)_{ij}+\frac{4ia\cos\th}{|q|^2} \pr_{\tt}(\psi_{ij})\nn\\
&=&2M_{i}^{k\a}\pr_\a(\psi_{kj}) +2M_{j}^{k\a}\pr_\a(\psi_{ik})+\frac{4ia\cos\th}{|q|^2} \pr_{\tt}(\psi_{ij}),\\
(\widehat{Q}\psi)_{ij} &=&({Q}\psi)_{ij} -\frac{4ia\cos\th}{|q|^2}\big(M_{i\tau}^l \psi_{jl}+M_{j\tau}^l \psi_{il}\big)\nn\\
&=& (\Ddot^\a M_{i\a}^k)\psi_{kj}+(\Ddot^\a M_{j\a}^k)\psi_{ik} -M_{i\a}^kM_k^{l\a}\psi_{lj}-2M_{i\a}^kM_{j}^{l\a}\psi_{kl}\nn\\
&&-M_{j\a}^kM_k^{l\a}\psi_{il}-\frac{4ia\cos\th}{|q|^2}\big(M_{i\tau}^l \psi_{jl}+M_{j\tau}^l \psi_{il}\big)
\eea
\esub
with the 1-forms $M_{i\a}^j$ defined by \eqref{eq:definitionofMalphaijwithoutambiguity}.

\begin{remark}
Note that the solution $\psi^{(p)}_{s,ij}$ of \eqref{eq:waveeqwidetildepsi1} satisfies \eqref{eq:ScalarizedWaveeq:general:Kerrpert} with the choices $D_0=4-2\de_{p0}$ and $F_{ij}=\widetilde{F}^{(p)}_{s,ij}+\underline{F}^{(p)}_{s,ij}$, where $\widetilde{F}^{(p)}_{s,ij}$ and $\underline{F}^{(p)}_{s,ij}$ are defined respectively in \eqref{def:tildef} and \eqref{def:tildef0}. This justifies the introduction of the model problem \eqref{eq:ScalarizedWaveeq:general:Kerrpert}.
\end{remark}

The aim of this section is to prove the microlocal energy-Morawetz estimates stated in  Theorem \ref{th:mainenergymorawetzmicrolocal} (see Section \ref{sec:statemenmainresultsection8globalEMFscalarizedwave}) for the above coupled system of inhomogenous scalar wave equations \eqref{eq:ScalarizedWaveeq:general:Kerrpert} under additional assumptions on $\psi_{ij}$ and $F_{ij}$. We start by introducing the microlocal calculus that will be used in Section \ref{sec:definitionmicrolocalenergyMorawetznorms} to define microlocal energy-Morawetz norms.


\subsection{Microlocal calculus on $\MM$}
\lab{sect:microlocalcalculus}


In this section, we introduce the necessary microlocal calculus on the manifold $\MM$. The material in Sections \ref{sec:mixedsymbolsonRn}--\ref{subsubsect:WeylquantizationofmixsymbolsonMM} is taken from \cite[Section 5]{MaSz24}.


\subsubsection{Mixed symbols on $\mathbb{R}^n$ and their Weyl quantization}
\lab{sec:mixedsymbolsonRn}


In view of our latter applications to microlocal energy-Morawetz estimates, we decompose $x=(x', x^n)\in\mathbb{R}^{n-1}\times\mathbb{R}$ and consider mixed operators which are PDO in $x'$ and differential in $x^n$. We first define $x^n$-tangential symbols on $\mathbb{R}^n$.
\begin{definition}[$x^n$-tangential symbols on $\mathbb{R}^n$]
 \label{def:symbols:Rn:rtangent}
For $m\in \mathbb{R}$, let $S^m_{tan}(\mathbb{R}^n)$ denote the set of functions $a$ which are $C^{\infty}(\mathbb{R}^n\times\mathbb{R}^{n-1})$ such that for all multi-indices $\alpha$, $\beta$,
 \beaa
\forall x=(x', x^n)\in\mathbb{R}^n, \,\, \forall \xi\in\mathbb{R}^{n-1},\quad \abs{\partial_x^{\alpha} \partial_\xi^{\beta}a(x,\xi)}\leq C_{\alpha,\beta} \langle \xi\rangle^{m-\abs{\beta}},
\eeaa
with $C_{\alpha,\beta}<+\infty$. An element $a\in S^m_{tan}(\mathbb{R}^n)$ is called an $x^n$-tangential symbol of order $m$. 
\end{definition}

Next, we introduce a class of mixed symbols on $\mathbb{R}^n$.
\begin{definition}[Mixed symbols on $\mathbb{R}^n$]
\label{PDO:Rn:Shom}
For $m\in\mathbb{R}$ and $N\in\mathbb{N}$, we define the class $\widetilde{S}^{m,N}(\mathbb{R}^n)$ of symbols as $a\in C^{\infty}(\mathbb{R}^n\times\mathbb{R}^n)$ such that for all $(x, \xi)\in\mathbb{R}^n\times\mathbb{R}^n$, we have, for $\xi=(\xi', \xi_n)$,  
\beaa
a(x,\xi)=\sum_{j=0}^N v_{m-j}(x,\xi')(\xi_n)^j,  \qquad v_{m-j}\in S^{m-j}_{tan}(\mathbb{R}^n).
\eeaa
An element $a\in\widetilde{S}^{m,N}(\mathbb{R}^n)$ is called a mixed symbol of order $(m, N)$. 
\end{definition}

In this paper, we will always rely on the Weyl quantization which we recall  below.
\begin{definition}[Weyl quantization of mixed symbols on $\mathbb{R}^n$]
\label{def:PDO:Rn:Weylquan}
Let $m\in\mathbb{R}$, $N\in\mathbb{N}$, and $a\in\widetilde{S}^{m,N}(\mathbb{R}^n)$. 
Then, the Weyl quantization of $a$ is given by
\beaa
\Opw(a)u(x) :=(2\pi)^{-n} \int_{\Reals^n}  \int_{\Reals^n} e^{i(x-y)\cdot\xi} a\bigg(\frac{x+y}{2},\xi\bigg)  u(y)d\xi d y.
\eeaa
\end{definition}

\begin{remark}\lab{rmk:WeylquantizationofmixedsymbolsonRnlocalinxn}
The Weyl quantization of mixed symbols is pseudo-differential w.r.t. $x'$ but differential w.r.t. $x^n$ so that it can be applied to functions that are defined on $\mathbb{R}^{n-1}\times I$ where $I$ is an open set of $\mathbb{R}$, see Lemma 5.13 and Remark 5.14 in \cite{MaSz24}.
\end{remark}


\subsubsection{Coordinates systems on $H_r$ and $\MM$}
\lab{sec:isochorecoordinatesonHr}


We introduce local coordinates on $\mathbb{S}^2$ for which the corresponding density is the one of the Lebesgue measure. This is done in the following lemma. 
\begin{lemma}\lab{lemma:isochorecoordinates}
Let the coordinates  $(x^1_0, x^2_0)$ and $(x^1_p, x^2_p)$ be defined by  
\beaa
x_0^1=\cos\th, \qquad x_0^2=\tphi, \qquad x^1_p=\sin\th\cos\tphi, \qquad x^2_p=\arcsin\left(\frac{\sin\th\sin\tphi}{\sqrt{1-(\sin\th)^2(\cos\tphi)^2}}\right),
\eeaa
with the corresponding coordinates patches 
\beaa
\mathbb{S}^2=\mathring{U}_0\cup\mathring{U}_p, \quad \mathring{U}_0:=\left\{(x_0^1, x_0^2)\,\,/\,\,\frac{\pi}{4}<\th<\frac{3\pi}{4}\right\}, \quad \mathring{U}_p:=\left\{(x^1_p, x^2_p)\,\,/\,\,\th\in [0,\pi]\setminus\left[\frac{\pi}{3}, \frac{2\pi}{3}\right]\right\}.
\eeaa
Then, the measure of the unit round 2-sphere in these coordinates has the density of the Lebesgue measure, i.e., 
\beaa
d\mathring{\ga}=dx^1_0 dx^2_0\quad\textrm{on}\quad U_0, \qquad d\mathring{\ga}=dx^1_p dx^2_p\quad\textrm{on}\quad U_p.
\eeaa
The coordinates systems $(x^1_0, x^2_0)$ and $(x^1_p, x^2_p)$ are called isochore coordinates.
The notation $(x^1, x^2)$ will be used to denote either $(x^1_0, x^2_0)$ or $(x^1_p, x^2_p)$.
\end{lemma}

\begin{proof}
See Lemma 5.19 in \cite{MaSz24}.
\end{proof}

We consider on $H_r$ the coordinates systems $(\tau, x^1_0, x^2_0)$ and $(\tau, x^1_p, x^2_p)$ with $(x^1_0, x^2_0)$ and $(x^1_p, x^2_p)$ constructed in Lemma \ref{lemma:isochorecoordinates}. This induces on $\MM$ coordinates $(\tau, r, x^1, x^2)$ with $(x^1, x^2)$ denoting  either $(x^1_0, x^2_0)$ or $(x^1_p, x^2_p)$ with
\beaa
x=(x', r), \qquad x':=(x^0, x^1, x^2), \qquad x^0:=\tau, \qquad x^3:=r,
\eeaa
and the corresponding coordinate patches
\beaa
\MM=U_0\cup U_p, \qquad U_q=\mathbb{R}_\tau\times \mathring{U}_q\times [r_+(1-\dhor), +\infty), \quad q=0,p,
\eeaa
with $\mathring{U}_q$, $q=0,p$ the coordinate patches on $\mathbb{S}^2$ introduced in Lemma \ref{lemma:isochorecoordinates}. Also, we denote by $\varphi_q: U_q\to \varphi_q(U_q)\subset\mathbb{R}^4$, $q=0,p$, the corresponding coordinates charts.

Next, we denote by $(\chi_q)_{q=0,p}$, a partition of unity subordinated to the covering by the coordinates patches $U_q$, $q=0,p$, i.e., $\chi_q$ are smooth scalar functions on $\MM$ satisfying  
\bea\lab{eq:partictionofuniityonHr}
\chi_0+\chi_p=1\,\,\,\textrm{on}\,\,\,\MM, \qquad \textrm{supp}(\chi_q)\subset U_q,\,\, q=0,p,  \qquad \pr_r\chi_q=\pr_\tau\chi_q=0,\,\, q=0,p.
\eea
Moreover, we also introduce smooth scalar functions $\widetilde{\chi}_q$, $q=0,p$ on $\MM$ satisfying 
\bea\lab{eq:partictionofuniityonHr:bis}
\widetilde{\chi}_q=1\,\,\,\textrm{on}\,\,\,\textrm{supp}(\chi_q),\,\, q=0,p, \quad \textrm{supp}(\widetilde{\chi_q})\subset U_q,\,\, q=0,p, \quad \pr_r\widetilde{\chi}_q=\pr_\tau\widetilde{\chi}_q=0,\,\, q=0,p.
\eea

To define symbols on $T^\star H_r$, we will need to introduce a norm of a co-vector $\xi'$ on the cotangent bundle. To this end, we introduce the following Riemannian metric $h_r$ on $H_r$
\bea
h_r=(d\tau)^2+\mathring{\ga}.
\eea 
Then, we define the length of a co-vector $\xi'$ by 
\bea
\abs{\xi'} := \sqrt{h_r^{ij} \xi_i'\xi_j'}=\sqrt{(\xi_0')^2+\mathring{\ga}^{ab}\xi_a'\xi_b'},
\eea
where latin indices $i,j$ take values $0,1,2$ and $a,b$ take values $1,2$.

In addition we have the following lemma concerning the properties of $\sqrt{|\det(\g)|}$ in the coordinates systems $(\tau, r, x^1, x^2)$.

\begin{lemma}\lab{lemma:spacetimevolumeformusingisochorecoordinates}
There exists a well-defined scalar function $f_0$ on $\MM$ such that $\g$ satisfies in the coordinates systems $(\tau, r, x^1, x^2)$
\bea\lab{eq:spacetimevolumeformusingisochorecoordinates}
f_0:=\sqrt{|\det(\g)|}, \qquad f_0=|q|^2(1+r^2\Ga_g).
\eea
This implies
\bea
\lab{def:dVref}
\sqrt{|\det(\g)|}d\tau dr dx^1dx^2=f_0 d\Vref, \qquad d\Vref:=d\tau dr dx^1dx^2,
\eea
with $d\Vref$ denoting the Lebesgue measure in the coordinates system $(\tau, r, x^1, x^2)$.
\end{lemma}

\begin{proof}
See Lemma 5.22 in \cite{MaSz24}.
\end{proof}


\subsubsection{Classes of mixed symbols on $\MM$}


We first define $r$-tangential symbols on $\MM$.
\begin{definition}[$r$-tangential symbols on $\MM$]
 \label{def:symbols:mflds:rtangent}
For $m\in \mathbb{R}$, let $S^m_{tan}(\MM)$ denote the set of functions $a$ which are $C^\infty$ in $r$ with values in $C^{\infty}(T^{\star}H_r)$ such that in $x=(x',r)=(\tau, x^1, x^2, r)$ coordinates of $\MM$,  for all multi-indices $\alpha$, $\beta$, and for all $q=0,p$,
 \beaa
\forall x=(x',r)\in\varphi_q(U_q), \,\, \forall \xi'\in T^{\star}H_r,\quad \abs{\partial_x^{\alpha} \partial_{\xi'}^\beta a(\varphi_q^{-1}(x),\xi_j'(dx^j_q)_{|_{\varphi_q^{-1}(x)}})}\leq C_{\alpha,\beta} \langle \xi'\rangle^{m-\abs{\beta}},
\eeaa
with $C_{\alpha,\beta}<+\infty$ and $\langle \xi'\rangle :=\sqrt{1+\abs{\xi'}^2}$. An element $a\in S^m_{tan}(\MM)$ is called an $r$-tangential symbol of order $m$. We also denote $S^{-\infty}_{tan}(\MM) :=\cap_{m\in \Reals} S^m_{tan}(\MM)$.
 \end{definition}

Next, we introduce a class of mixed symbols on $\MM$.
\begin{definition}[Mixed symbols on $\MM$]
\label{PDO:MM:Shom}
For $m\in\mathbb{R}$ and $N\in\mathbb{N}$, we define the class $\widetilde{S}^{m,N}(\MM)$ of symbols as $a\in C^{\infty}(T^\star\MM)$ such that for all $q=0,p$, for all $x\in\varphi_q(U_q)$ and for all $\xi\in\mathbb{R}^4$, 
\beaa
a(\varphi_q^{-1}(x),\xi_\a dx^\a)=\sum_{j=0}^N v_{m-j}(\varphi_q^{-1}(x),\xi_i'(dx^i_q)_{|_{\varphi_q^{-1}(x)}})(\xi_3)^j,  \qquad v_{m-j}\in S^{m-j}_{tan}(\MM),
\eeaa 
where $\xi=(\xi', \xi_3)$. An element $a\in\widetilde{S}^{m,N}(\MM)$ is called a mixed symbol of order $(m,N)$. We also denote $\widetilde{S}^{-\infty,N}(\MM) :=\cap_{m\in \Reals}\widetilde{S}^{m,N}(\MM)$.
\end{definition}

\begin{remark}
Notice that $S^m_{tan}(\MM)=\widetilde{S}^{m,0}(\MM)$. 
\end{remark}


\subsubsection{Weyl quantization of mixed symbols on $\MM$}
\lab{subsubsect:WeylquantizationofmixsymbolsonMM}


Recall the coordinates charts $(\varphi_q)_{q=0,p}$ on $\MM$, as well as the partition of unity $(\chi_q)_{q=0,p}$, introduced in Section \ref{sec:isochorecoordinatesonHr}. For $m\in\mathbb{R}$ and $N\in\mathbb{N}$, given $a\in\widetilde{S}^{m,N}(\MM)$, we introduce the following notation, for all $q=0,p$, $x\in\varphi_q(U_q)$ and $\xi\in\mathbb{R}^4$,
\bea\lab{eq:defsymbolonRnaqchiq}
a_{q,\chi_q}(x, \xi):=\chi_q(\varphi_q^{-1}(x))a\big(\varphi_q^{-1}(x), \xi_\a(dx_q^\a)_{|_{\varphi_q^{-1}(x)}}\big), \qquad a_{q,\chi_q}\in \widetilde{S}^{m,N}(\mathbb{R}^4),
\eea
where the class of mixed symbols $\widetilde{S}^{m,N}(\mathbb{R}^n)$ has been introduced in Definition \ref{PDO:Rn:Shom}.

\begin{definition}[Weyl quantization of mixed symbols on $\MM$]
\lab{def:weylquantizationforrhomsymbolsMM}
Let $m\in\mathbb{R}$, $N\in\mathbb{N}$ and $a\in\widetilde{S}^{m,N}(\MM)$. We associate to $a$ the operator $\Opw(a)$ in the Weyl quantization as follows 
\bea\lab{eq:defintionWeylquantizationmixedsymbolonMM}
\Opw(a)\psi := \sum_{q=0,p}\widetilde{\chi}_q\varphi_q^{\#}\Opw(a_{q,\chi_q})[(\varphi_q^{-1})^\#(\widetilde{\chi}_q\psi)],
\eea
where $(\widetilde{\chi}_q)_{q=0,p}$ is given by \eqref{eq:partictionofuniityonHr:bis} and $a_{q,\chi_q}$ is given by \eqref{eq:defsymbolonRnaqchiq}, i.e., for $x\in\varphi_{q'}(U_{q'})$, $q'=0,p$, 
\beaa
\Opw(a)\psi(\varphi_{q'}^{-1}(x)) &=& \frac{1}{(2\pi)^4}\sum_{q=0,p}\widetilde{\chi}_q(\varphi_{q'}^{-1}(x))\\
&&\times\int_{\mathbb{R}^4}\int_{\mathbb{R}^4}e^{i(x_{q,q'}-y)\c\xi}a_{q,\chi_q}\left(\frac{x_{q,q'}+y}{2}, \xi\right)(\widetilde{\chi}_q\psi)\circ\varphi_q^{-1}(y)dy d\xi,
\eeaa
where $x_{q,q'}=\varphi_q\circ\varphi_{q'}^{-1}(x)$ if $\widetilde{\chi}_q(\varphi_{q'}^{-1}(x))\neq 0$.
\end{definition} 

\begin{remark}
Since $a_{q,\chi_q}\in \widetilde{S}^{m,N}(\mathbb{R}^4)$, Remark \ref{rmk:WeylquantizationofmixedsymbolsonRnlocalinxn} applies to $\Opw(a_{q,\chi_q})$. In view of \eqref{eq:defintionWeylquantizationmixedsymbolonMM}, we infer that $\Opw(a)$ for $a\in\widetilde{S}^{m,N}(\MM)$ is pseudo-differential on $H_r$ but differential w.r.t. $r$ so that it can be applied to functions that are defined on $\MM_{r_1, r_2}$ for $r_+(1-\dhor)\leq r_1<r_2<+\infty$. 
\end{remark}

\begin{remark}
Note that Definition \ref{def:weylquantizationforrhomsymbolsMM} is invariant modulo a smoothing operator under change of coordinates that preserve the isochore property of Lemma \ref{lemma:isochorecoordinates}, but not under general change of coordinates, see Remark 5.29 in \cite{MaSz24}.
\end{remark}

Next, we consider the properties of the Weyl quantization of symbols in $\widetilde{S}^{m,N}(\MM)$ w.r.t. composition and adjoint. 
\begin{proposition}
\label{prop:PDO:MM:Weylquan:mixedoperators}
The Weyl quantization satisfies the following properties for symbols in the class 
$\widetilde{S}^{m,N}(\MM)$: 
\begin{enumerate}[label=\arabic*)] 
\item For mixed symbols $a_1$ and $a_2$ of respective orders $(m_1,N_1)$ and $(m_2,N_2)$, we have
\bea\lab{eq:propWeylquantization:MM:composition:mixedsymbols}
\bsplit
[\Opw(a_1), \Opw(a_2)]=\Opw(a_3),\quad a_3=\frac{1}{i}\{a_1, a_2\} +\widetilde{S}^{m_1+m_2-3,N_1+N_2}(\MM),\\
\Opw(a_1)\circ\Opw(a_2) +\Opw(a_2)\circ\Opw(a_1)=\Opw(a_3),\quad a_3=2a_1a_2 + \widetilde{S}^{m_1+m_2-2,N_1+N_2}(\MM).
\end{split}
\eea

\item In the particular case where $a_1(\varphi_q^{-1}(x),\xi_\a dx^\a)=v_1(r)\xi_3^{N_1}$ for $x=(r,x')\in\varphi_q(U_q)$ and $\xi=(\xi', \xi_3)\in\mathbb{R}^4$, which is a mixed symbol of order $(m_1, N_1)$ with $m_1=N_1$, we have, with $a_2$ of order $(m_2,N_2)$
\bea\lab{eq:propWeylquantization:MM:composition:mixedsymbols:specialcase}
\bsplit
[\Opw(a_1), \Opw(a_2)]=\Opw(a_3),\quad a_3=\frac{1}{i}\{a_1, a_2\} +\tilde{a}_3,\\
\tilde{a}_3=0\quad\textrm{if}\quad \max(N_1, N_2)\leq 2, \qquad \tilde{a}_3\in\widetilde{S}^{m_1+m_2-3,N_1+N_2-3}(\MM)\quad\textrm{if}\quad \max(N_1, N_2)\geq 3,\\
\Opw(a_1)\circ\Opw(a_2) +\Opw(a_2)\circ\Opw(a_1)=\Opw(a_4),\quad a_4=2a_1a_2 + \tilde{a}_4,\\
\tilde{a}_4=0\quad\textrm{if}\quad \max(N_1, N_2)\leq 1, \qquad \tilde{a}_4\in\widetilde{S}^{m_1+m_2-2,N_1+N_2-2}(\MM)\quad\textrm{if}\quad \max(N_1, N_2)\geq 2.
\end{split}
\eea

\item In the particular case where $a_1$ and $a_2$ are mixed symbols of respective orders $(m_1, 1)$ and $(m_2, 1)$, and $f=f(x^n)$, we have
\bea\lab{eq:propWeylquantization:MM:composition:mixedsymbols:specialcase:1}
\bsplit
&[\Opw(a_1), \Opw(f(r)a_2)]=\Opw(a_3),\quad a_3=\frac{1}{i}\{a_1, f(r)a_2\} +\tilde{a}_3,\\
&\tilde{a}_3=f(r)\widetilde{S}^{m_1+m_2-3,2}(\MM)+\widetilde{S}^{m_1+m_2-3,1}(\MM),\\
&\Opw(a_1)\circ\Opw(f(r)a_2) +\Opw(f(r)a_2)\circ\Opw(a_1)=\Opw(a_4),\quad a_4=2f(r)a_1a_2 + \tilde{a}_4,\\
&\tilde{a}_4=f(r)\widetilde{S}^{m_1+m_2-2,2}(\MM)+\widetilde{S}^{m_1+m_2-2,1}(\MM).
\end{split}
\eea

\item For a mixed symbol $a(x,\xi)$, the adjoint, w.r.t. the Lebesgue measure $d\Vref$ in $(\tau, r, x^1, x^2)$ coordinates, of its Weyl quantization is given by 
\bea\lab{eq:propWeylquantization:MM:adjoint:mixedsymbols}
(\Opw(a))^{\star} = \Opw(\bar{a}).
\eea
In particular, the Weyl quantization of a real-valued symbol is a self-adjoint operator  w.r.t. the Lebesgue measure $d\Vref$ in $(\tau, r, x^1, x^2)$ coordinates.
\end{enumerate}
\end{proposition}

\begin{proof}
See Proposition 5.31 in \cite{MaSz24}.
\end{proof}

Also, we consider the action of the Weyl quantization of mixed symbols on Sobolev spaces.  
\begin{lemma}\lab{lemma:actionmixedsymbolsSobolevspaces:MM}
Let $m\in\mathbb{R}$, $N\in\mathbb{N}$, let $I$ be an interval of $[r_+(1-\dhor), +\infty)$, and let $a\in\widetilde{S}^{m,N}(\MM)$ be of the form, for all $q=0,p$, $x\in \varphi(U_q)$ and $\xi\in\mathbb{R}^4$,
\beaa
a(\varphi_q^{-1}(x),\xi_\a dx^\a)=v_{m-N}(\varphi_q^{-1}(x),\xi_i'(dx^i_q)_{|_{\varphi_q^{-1}(x)}})(\xi_3)^N,  \qquad v_{m-N}\in S^{m-N}_{tan}(\MM).
\eeaa 
Then we have for all $s\in\mathbb{R}$
\beaa
\|\Opw(a)\psi\|_{H^{s-m+N}(H_r)}&\les& \sum_{j=0}^N\|\pr_r^j\psi\|_{H^s(H_r)},\\
\|\Opw(a)\psi\|_{L^2_r(I, H^{s-m+N}(H_r))}&\les& \sum_{j=0}^N\|\pr_r^j\psi\|_{L^2_r(I, H^s(H_r))}.
\eeaa
\end{lemma}

\begin{proof}
See Lemma 5.32 in \cite{MaSz24}.
\end{proof}

Finally, we introduce a notation for the symbol of an operator corresponding to the Weyl quantization of a mixed symbol. 
\begin{definition}\lab{def:notationforthesymbolofWeylquantizationmixedsymbol}
If $A=\Opw(a)$ for some mixed symbol $a\in\widetilde{S}^{m,N}(\MM)$, then we denote the symbol $a$ of $A$ by $\sigma(A)$, i.e., 
\beaa
\sigma(A):=a, \qquad A=\Opw(a).
\eeaa
\end{definition}


\subsubsection{Relevant mixed symbols and operators on $\MM$}
\lab{sec:relevantmixedsymbolsonMM}


In this section, we list the symbols and operators introduced in \cite[Sections 5, 6 and 7]{MaSz24} that will appear in the microlocal energy-Morawetz estimates of Sections \ref{sect:microlocalenergyMorawetztensorialwaveequation} and \ref{sec:proofofth:main:intermediary}.  We start with the mixed symbols: 
\begin{enumerate}
\item\lab{items:symbolsandoperators:point1}  Recalling that mixed symbols are defined w.r.t. the $(\tau, r, x^1, x^2)$ coordinates systems of Section \ref{sec:isochorecoordinatesonHr}, we introduce the following symbols:
\bea\lab{eq:basicparticularsymbolsusedeverywhere}
\xit:=\langle \xi, \partial_{\tt}\rangle=\xi_0, \qquad {\xiphi}:= \langle \xi, \partial_{\tphi}\rangle=\frac{\partial x^a}{\partial \tphi} \xi_a, \qquad \xi_r:= \langle \xi, \partial_{r} \rangle = \xi_3,
\eea
where $\xit$, $\xiphi$ and $\xir$ satisfy 
\bea
\Opw(i\xit)=\pr_\tau, \qquad \Opw(i\xiphi)=\pr_{\tphi}, \qquad \Opw(i\xir)=\pr_r.
\eea

\item The mixed symbols $e_0$, $s_0$, $b_{\tphi}$, $b_{\tt}$, $r_{\trap}$ and $\Theta_i$, $i=-1,0, 1,\cdots, \iota$, are in $\widetilde{S}^{0,0}(\MM)$, where $\iota$ will be chosen as a large enough integer. 

\item The mixed symbols $e$, $\sigma_{\trap}$ and $\upsilon$ are in $\widetilde{S}^{1,0}(\MM)$.

\item The mixed symbols $\sigma_{\trap}$ and $\upsilon$ are given by 
\bea\lab{eq:precisedefinitionofsigmatrapandvarpi}
\sigma_{\trap}:=(r-r_{\trap})\upsilon,\qquad \upsilon:=\sqrt{1+\xi_0^2+\mathring{\ga}^{bc} \langle \xi, \pr_{x^b}\rangle \langle \xi, \pr_{x^c}\rangle}.
\eea

\item The mixed symbols $r_{\trap}$, $\upsilon$ and $\Theta_i$ satisfy 
\bea\label{eq:nodependenceonrforthesymbolesrtrapupsilonTheta}
\pr_r(r_{\trap})=0, \qquad \pr_r(\upsilon)=0, \qquad \pr_r(\Theta_i)=0, \quad i=-1,0, 1,\cdots, \iota.
\eea

\item The mixed symbols $e_0$, $s_0$, $b_{\tphi}$, $b_{\tt}$, $r_{\trap}$, $\Theta_j$, $j=-1,0, 1,\cdots, \iota$, $e$, $\sigma_{\trap}$ and $\upsilon$ are real-valued. 

\item The following pointwise estimates hold in $\Mtrap$\footnote{Note that $\lesssim$ appearing in the last  estimate of \eqref{eq:pointwisecontrolbyesetuptogetGardingtypeinequalities} depends on the choice of $b$. We do not indicate this dependence as it will be used in practice for a finite number of choices of $b\in\widetilde{S}^{1,0}(\MM)$.}  
\begin{equation}\lab{eq:pointwisecontrolbyesetuptogetGardingtypeinequalities}
\begin{split}
1+\big(|e_0|+|s_0| +|b_{\tphi}|+|b_{\tt}|\big)\upsilon \lesssim {}&e, \quad 1+\big(|b_{\tphi}|+|b_{\tt}|\big)\upsilon^2 \lesssim e^2,\\
1+(|\{b_{\tphi}, b\}|+|\{b_{\tt}, b\}|)\upsilon\lesssim {}&e\qquad \forall\, b\in\widetilde{S}^{1,0}(\MM).
\end{split}
\end{equation}
\end{enumerate}
In the above, item \eqref{items:symbolsandoperators:point1} follows from \cite[Section 5.2.6]{MaSz24}, and the remaining items are taken from Proposition 7.23 in \cite{MaSz24}.

In addition, following Proposition 7.23 in \cite{MaSz24}, we consider the following pseudodifferential operators
\bea
\lab{eq:definitionofXandEinMtrap:sect8.1.5}
X:=\Opw(is_0\mu\xi_r+ib_{\tphi}\xiphi+ib_{\tt}\xit)+A\pr_\tau, \qquad E:=\Opw(e_0),
\eea
where $A\geq 2$ is a large enough constant, and the following vectorfields
\bea
V_i:=\pr_\tau+d_i(r)\pr_{\tphi}, \qquad i=-1,0, 1,\cdots,\iota, 
\eea 
for some large integer $\iota$ and smooth real-valued functions $d_i(r)$ supported  in $r\leq 10m$.


\subsection{Microlocal energy-Morawetz norms and main estimates for \eqref{eq:ScalarizedWaveeq:general:Kerrpert}}



\subsubsection{Microlocal energy-Morawetz norms}
\lab{sec:definitionmicrolocalenergyMorawetznorms}


\begin{definition}[Microlocal energy-Morawetz norms]
\lab{def:microlocalenergyMorawetznorms}
Let $e,\, \sigma_{\trap}\in\widetilde{S}^{1,0}(\MM)$ be the mixed symbols introduced in Section \ref{sec:relevantmixedsymbolsonMM}. Then, for a scalar field $\psi$, we define  the microlocal Morawetz norm $\widetilde {\M}[\psi]$ by 
\bea\lab{eq:definitionofmicrolocalMorawetznormwidetildeM}
\widetilde{\M}[\psi]:={}{\M}[\psi](\Iti) +\int_{\MM_{r_+(1+2\dhor), 10m}}\Big(|\Opw(\sigma_{\trap})\psi|^2+|\Opw(e)\psi|^2\Big),
\eea
where we recall from Definition \ref{def:definitionofthetimetauR} that $\Iti=(\tmic, +\infty)$, and for a family of scalars $\psi_{ij}$, we define  the microlocal Morawetz norm $\widetilde {\M}[\pmb\psi]$ by\footnote{Recall that according to our slight abuse of notations introduced in Remark \ref{rmk:abusenotationbetweentensorandscalarizedversioninEMFnorms}, we use the notation $\widetilde {\M}[\pmb\psi]$ even though $\psi_{ij}$ do not necessarily come from the scalarization of a tensor $\pmb\psi\in\sk_2(\mathbb{C})$.} 
\bea\lab{eq:definitionofmicrolocalMorawetznormwidetildeM:tensor}
\widetilde{\M}[\pmb\psi]:={}{\M}[\pmb\psi](\Iti) +\sum_{i,j}\int_{\MM_{r_+(1+2\dhor), 10m}}\Big(|\Opw(\sigma_{\trap})\psi_{ij}|^2+\Opw(e)\psi_{ij}|^2\Big).
\eea
Also, for any $\de\in[0,\frac{1}{3}]$, we define $\widetilde{\M}_{\de}[\psi]$ and $\widetilde{\M}_{\de}[\pmb\psi]$ by replacing respectively ${\M}[\psi](\Iti)$ and ${\M}[\pmb\psi](\Iti)$ in \eqref{eq:definitionofmicrolocalMorawetznormwidetildeM} and \eqref{eq:definitionofmicrolocalMorawetznormwidetildeM:tensor} by ${\M}_{\de}[\psi](\Iti)$ and ${\M}_{\de}[\pmb\psi](\Iti)$. Finally, for any nonnegative integer $\reg$, we define the high-order microlocal Morawetz norm  $\widetilde{\M}^{(\reg)}[\pmb\psi]$ and $\widetilde{\M}^{(\reg)}_{\de}[\pmb\psi]$, as well as the combined norms $\widetilde{\EMF}^{(\reg)}_\de[\pmb\psi]$, $\widetilde{\EMF}^{(\reg)}[\pmb\psi]$, $\widetilde{\EM}^{(\reg)}_\de[\pmb\psi]$, $\widetilde{\EM}^{(\reg)}[\pmb\psi]$ and $\widetilde{\MF}^{(\reg)}[\pmb\psi]$, as in \eqref{eq:definitionofhigherordernormsFregEregMdeltaregMregNNdeltaregwidehatNreg} \eqref{eq:definitionofhigherordernormsFregEregMdeltaregMregNNdeltaregwidehatNreg:combinednorms}.
\end{definition}

Next, for families of scalars $\psi_{ij}, F_{ij}$, we introduce\footnote{Again, recall that according to our slight abuse of notations introduced in Remark \ref{rmk:abusenotationbetweentensorandscalarizedversioninEMFnorms}, we use the notation $\widetilde{\mathcal{N}}[\pmb\psi, \pmb F]$ even though $\psi_{ij}, F_{ij}$ do not necessarily come from the scalarization of tensors $\pmb\psi, \F\in\sk_2(\mathbb{C})$.}
\bea
\lab{def:NNtintermsofNNtMora:NNtEner:NNtaux:wavesystem:EMF} 
\NNt[\pmb \psi,\pmb F]:=\NNtmora[\pmb \psi, \pmb F]+\NNtener[\pmb \psi, \pmb F]
+\NNtaux[\pmb F],
\eea
with 
\bsub
\lab{def:NNtMora:NNtEner:NNtaux:wavesystem:EMF}
\begin{align}
\NNtmora[\pmb \psi, \pmb F]:=&\sum_{i,j}\bigg|\int_{\MM_{r_+(1+\dhor'),\Rmic}}\Re\Big(F_{ij} \ov{X\psi_{ij}}\Big)+\int_{\MM_{\Rmic,+\infty}(\Iti)}\Re\Big(F_{ij} \ov{({X} + w)\psi_{ij}}\Big)\bigg|,\\
\NNtener[\pmb \psi, \pmb F]:=&  \sup_{\tau\in\Reals}\sum_{n=-1}^{\iota}\sum_{i,j}\bigg|\int_{\Mtrap(\tmic,\tau)}\Re\Big(\ov{|q|^{-2}\Opw(\Theta_n)(\qs F_{ij})}V_n\Opw(\Theta_n)\psi_{ij}\Big)\bigg|\nn\\
&+\sum_{i,j}\sup_{\tau\geq\tmic}\bigg|\int_{\Mntrap(\tmic, \tau)}{\Re\Big(F_{ij}\ov{\pr_{\tau}\psi_{ij}}\Big)}\bigg|,\\
\NNtaux[\pmb F]:=&\sum_{i,j}\bigg(\int_{\Mntrap(\Iti)}|F_{ij}|^2 +(\ep+\dhor)\int_{\Mtrap(\Iti)}|F_{ij}|^2 \nn\\
& \qquad\quad + \ep\int_{\Mtrap(\Iti)}\tt^{-1-\dec}|F_{ij}|^2\bigg),
\end{align}
\esub
where: 
\begin{itemize}
\item the constants $\dhor'$ and $\Rmic$ are chosen as in Remark \ref{rmk:choiceofconstantRbymeanvalue} below,
\item in $\MM_{r_+(1+\dhor'),\Rmic}$, $X\in\Opw(\widetilde{S}^{1,1}(\MM))$ is defined as in Section \ref{sec:relevantmixedsymbolsonMM},
\item in $\MM_{\Rmic,+\infty}(\Iti)$, $X$ is a vectorfield satisfying $X=(1+O(r^{-1}))\pr^{\textrm{BL}}_r+A\pr_{\tt}$ and $w$ is a real-valued function satisfying $w=cr^{-1} + O(r^{-2})$,
\item in $\Mtrap$, for $n=-1,0,1,\cdots,\iota$, the mixed symbols $\Theta_n\in\widetilde{S}^{0,0}(\MM)$ and vectorfields $V_n$ have been introduced in Section \ref{sec:relevantmixedsymbolsonMM}.
\end{itemize}

\begin{remark}[Choice of constants $\dhor'$ an $\Rmic$]\label{rmk:choiceofconstantRbymeanvalue}
The constant $\dhor'\in [\dhor, 2\dhor]$ is chosen to verify
\bea
\lab{de:choiceofdhor'}
\nn&&\sum_{i,j=1}^3\int_{H_{r_+(1+\dhor')}(\Iti)}\big(|{\pr^{\leq 1}}\psi_{ij}|^2+|\square_\g\psi_{ij}|^2\big) d\tt dx^1dx^2\\ 
&\leq& \frac{1}{\dhor}\sum_{i,j=1}^3\int_{\MM_{r_+(1+\dhor), r_+(1+2\dhor)}(\Iti)}\big(|{\pr^{\leq 1}}\psi_{ij}|^2+|\square_\g\psi_{ij}|^2\big) d\Vref,
\eea
and the constant $\Rmic\in [{\Nmic}m, ({\Nmic}+1)m]$, with $\Nmic\geq 20$ a large enough integer, is chosen to verify 
\bea
\label{eq:choiceofRvalue:Kerr}
\nn&&\sum_{i,j=1}^3\int_{H_{\Rmic}(\Iti)}\big(|{\pr^{\leq 1}}\psi_{ij}|^2+|\square_\g\psi_{ij}|^2\big) d\tt dx^1dx^2\\
&\leq& \frac{1}{m}\sum_{i,j=1}^3\int_{\MM_{\Nmic m, (\Nmic +1)m}(\Iti)}\big(|{\pr^{\leq 1}}\psi_{ij}|^2+|\square_\g\psi_{ij}|^2\big) d\Vref,
\eea
where $d\Vref=d\tau d r dx^1 dx^2$ is given as in \eqref{def:dVref}.
This implies that for solutions to the coupled system of scalar wave equations \eqref{eq:ScalarizedWaveeq:general:Kerrpert}, we have
{\bea
\lab{de:choiceofdhor':extendedRWsystem}
\nn&& \sum_{i,j=1}^3\int_{H_{r_+(1+\dhor')}(\Reals)}\big(|{\pr^{\leq 1}}\psi_{ij}|^2+|\square_\g\psi_{ij}|^2\big) d\tt dx^1dx^2\\ 
&\les& \frac{1}{\dhor}\sum_{i,j=1}^3\int_{\MM_{r_+(1+\dhor), r_+(1+2\dhor)}(\Iti)}\big(|{\pr^{\leq 1}}\psi_{ij}|^2+|F_{ij}|^2\big) d\Vref,
\eea}
and 
{\bea
\label{eq:choiceofRvalue:Kerr:extendedRWsystem}
&& \sum_{i,j=1}^3\int_{H_{\Rmic}(\Reals)}\big(|{\pr^{\leq 1}}\psi_{ij}|^2+|\square_\g\psi_{ij}|^2\big) d\tt dx^1dx^2\nn\\
&\les& \frac{1}{m}\sum_{i,j=1}^3\int_{\MM_{{\Nmic}m, ({\Nmic}+1)m}(\Iti)}\big(|{\pr^{\leq 1}}\psi_{ij}|^2+|F_{ij}|^2\big) d\Vref.
\eea}
\end{remark}


\subsubsection{Main global energy-Morawetz estimates for \eqref{eq:ScalarizedWaveeq:general:Kerrpert}}
\label{sec:statemenmainresultsection8globalEMFscalarizedwave}


We are now ready to state the main result of Section \ref{sect:microlocalenergyMorawetztensorialwaveequation} on global energy-Morawetz estimates for solutions to \eqref{eq:ScalarizedWaveeq:general:Kerrpert} in perturbations of Kerr.

\begin{theorem}[Global energy-Morawetz estimates for \eqref{eq:ScalarizedWaveeq:general:Kerrpert}]\lab{th:mainenergymorawetzmicrolocal}
Let $(\MM, \g)$ satisfy the assumptions of Sections \ref{subsect:assumps:perturbednullframe}, \ref{subsubsect:assumps:perturbedmetric} and \ref{sec:regulartripletinperturbationsofKerr} with  $\g=\gam$ for $\tau\in\Reals\setminus(\tau_1, \tau_2)$, and let $\psi_{ij}$ be a solution to the system of scalarized wave equations \eqref{eq:ScalarizedWaveeq:general:Kerrpert}  with RHS $F_{ij}$, where the cut-off $\chi_{\tau_1, \tau_2}$ appearing in \eqref{eq:ScalarizedWaveeq:general:Kerrpert} satisfies \eqref{eq:propertieschitoextendmetricg}. Assume that $\psi_{ij}$ vanishes in $\MM(-\infty, \tmic)\cap\{r\leq (\Nmic+1)m\}$, with $\tmic$ and $N_0$ introduced in Definition \ref{def:definitionofthetimetauR}, and satisfies $\psi_{ij}=\pmb\psi(\Om_i, \Om_j)$ in $\MM(\tau_1+1, \tau_2-2)$ for a tensor $\pmb\psi\in\sk_2(\mathbb{C})$ and $\psi_{ij}=\pmb\psi((\Om_K)_i, (\Om_K)_j)$ in $\MM(\tau_2, +\infty)$ for a tensor $\pmb\psi\in\sk_{2,K}(\mathbb{C})$. Finally, assume that $F_{ij}$ are supported in $\MM(\tmic,\tau_2)$.
Then, we have 
\bea\lab{th:eq:mainenergymorawetzmicrolocal:tensorialwave:scalarized:eachpsisp}
\widetilde{\EMF}[\pmb \psi] \les \sup_{\tau\in[\tmic,\tau_1+2]}\E[\pmb\psi](\tau)+
\Errdefect[\pmb\psi]
+\A[\pmb \psi](\Iti)+\NNt[\pmb \psi,\pmb F],
\eea
where $\NNt[\pmb \psi,\pmb F]$ is given as in \eqref{def:NNtintermsofNNtMora:NNtEner:NNtaux:wavesystem:EMF} and where we have defined\footnote{While the control provided by $\A[\pmb \psi](\Iti)$ on $\tau\in[\tmic,\tau_1+2]$ and $\Errdefect[\pmb\psi]$ on $\tau\in[\tmic,\tau_1+1]$ is already included in $\sup_{\tau\in[\tmic,\tau_1+2]}\E[\pmb\psi](\tau)$, we nevertheless include the interval $\tau\in[\tmic,\tau_1+2]$ (respectively $\tau\in[\tmic,\tau_1+1]$) in their definition for convenience.}
\bea
\lab{def:EM-1norms:Reals}
\A[\pmb \psi](\Iti):=\sum_{i,j}\bigg(\int_{\MM(\Iti)}r^{-3}|\psi_{ij}|^2 +\sup_{\tau\in\Iti}\int_{\Sigma_{\tau}}r^{-2}|\psi_{ij}|^2+\int_{\II_+(\Iti)}r^{-2}|\psi_{ij}|^2\bigg)
\eea
and
\bea
\lab{def:Errdefectofpsi}
\Errdefect[\pmb\psi]:=\sup_{\tau\in[\tmic,\tau_1+1]\cup[\tau_2-2,\tau_2]}\sum_{j}\Big(\E[x^i\psi_{ij}](\tau)+\E[x^i\psi_{ji}](\tau)\Big).
\eea
\end{theorem}

The rest of Section \ref{sect:microlocalenergyMorawetztensorialwaveequation} is dedicated to the proof of Theorem \ref{th:mainenergymorawetzmicrolocal}. In Section 
\ref{sec:MorawetzestimatesinMMrplus1plusdhorprimRmicinhomscalarwavesystem}, we derive a Morawetz estimate on $\MM_{r_+(1+\dhor'),{\Rmic}}$ for the system of scalar wave equations \eqref{eq:ScalarizedWaveeq:general:Kerrpert}. Then, we derive energy-Morawetz estimates near infinity and redshift estimates for \eqref{eq:ScalarizedWaveeq:general:Kerrpert} in Section \ref{sec:energyMorawetznearinfinityandredshiftestimates:globalsystmscalarwave}. Finally, in Section \ref{subsect:globalEMFestimate:scalarizeeqfromtensorial:0}, we conclude the proof of Theorem \ref{th:mainenergymorawetzmicrolocal}.


\subsection{Morawetz estimate in $\MM_{r_+(1+\dhor'),{\Rmic}}$ for \eqref{eq:ScalarizedWaveeq:general:Kerrpert}}
\lab{sec:MorawetzestimatesinMMrplus1plusdhorprimRmicinhomscalarwavesystem}


The goal of this section is to prove Proposition \ref{prop:Morawetz:middleregion:globalextendedwavesystem} on microlocal Morawetz estimates for the coupled system of scalar wave equations \eqref{eq:ScalarizedWaveeq:general:Kerrpert} in $\MM_{r_+(1+\dhor'),{\Rmic}}$. To this end, we will rely on the following lemma providing a microlocal Morawetz estimate for the scalar wave equation.

\begin{lemma}\lab{lem:conditionaldegenerateMorawetzflux:pertKerrrp1pdhorpR}
 Let $\g$ satisfy the assumptions of Section \ref{subsubsect:assumps:perturbedmetric}, and let the scalar function $\psi$ be a solution to the inhomogeneous wave equation 
 \beaa
 \square_{\g}\psi=F.
 \eeaa
  Assume that $\psi$ vanishes in $\MM(-\infty, \tmic)\cap\{r\leq (\Nmic+1)m\}$,  that the metric $\g$ coincides with $\gam$ for $\tau\in\Reals\setminus(\tau_1, \tau_2)$, and that $F$ is supported in $\MM(\tmic,\tau_2)$. Let $\dhor'$ and ${\Rmic}$ be constants such that\footnote{Here, we do not require $\dhor'$ and ${\Rmic}$ to satisfy \eqref{de:choiceofdhor'} and \eqref{eq:choiceofRvalue:Kerr} respectively. This will be needed in Lemma \ref{lem:lowerorderterms:controlled:scalarizedwavefromtensorial} below in order to estimate the boundary terms on the timelike hypersurfaces $r=r_+(1+\dhor')$ and $r=\Rmic$.} $\dhor'\in[\dhor, 2\dhor]$ and $\Rmic\in [\Nmic m, (\Nmic+1)m]$. Then, we have
\bea
\lab{eq:intermediaryestimateformicrolocalconditiondegenerateMorrp1pdhorpR:robust}
\nn&& c\Bigg[\int_{\MM_{r_+(1+\dhor'),{\Rmic}}}\frac{\mu^2|\pr_r\psi|^2}{r^2} +\int_{\MM_{r_+(1+\dhor'),10m}}\Big(|\Opw(\sigma_{\trap})\psi|^2+|\Opw(e)\psi|^2\Big)\\
\nn&&+\int_{\Mntrap_{r_+(1+\dhor'),{\Rmic}}}\frac{|\pr_\tau\psi|^2+|\nab\psi|^2}{r^2}\Bigg]+\textbf{BDR}^{-}_{r={\Rmic}}[\psi](\Iti)+\int_{\MM_{r_+(1+\dhor'),{\Rmic}}}\Re\Big(\square_{\g}\psi\ov{(X+E)\psi}\Big)\nn\\
\nn&&\les_{\Rmic} (\ep+\dhor){\int_{\Mntrap}|\square_{\g}\psi|^2}+\ep\int_{\Mtrap}\tau^{-1-\dec}|\square_{\g}\psi|^2\\
&&+(\ep+\dhor^6)\int_{\Mtrap}\Big|\Opw(\widetilde{S}^{-1,0}(\MM))\square_{\g}\psi\Big|^2+\ep\EM[\psi](\Iti)+\dhor\M[\psi](\Iti)
\nn \\
&&+\frac{1}{\dhor^6}\int_{\MM_{r_+(1+\dhor'),{\Rmic}}}|\psi|^2+\bigg(\int_{H_{{\Rmic}}}\Big(|\pr\psi|^2+|\psi|^2\Big)\bigg)^{\frac{1}{2}}\bigg(\int_{H_{{\Rmic}}}|\psi|^2\bigg)^{\frac{1}{2}},
\eea
where $\textbf{BDR}^{-}_{r={\Rmic}}[\psi]$ denotes a boundary term\footnote{See the third line of \cite[Inequality (7.100)]{MaSz24} for its explicit form. We do not recall it here as it will be canceled, up to lower order terms, by the corresponding boundary term arising from integration by parts in $r\geq\Rmic$, see Section \ref{subsubsect:globalmicrolocalMora:scalarizeeqfromtensorial:0}.} on $H_{\Rmic}$, and where $\sigma_{\trap},\,  e\in\widetilde{S}^{1,0}(\MM)$, $X\in\Opw(\widetilde{S}^{1,1}(\MM))$ and $E\in \Opw(\widetilde{S}^{0,0}(\MM))$ are defined as in Section \ref{sec:relevantmixedsymbolsonMM}.
\end{lemma}

\begin{proof}
See \cite[Equation (7.100)]{MaSz24}.
\end{proof}


\subsubsection{Application of \eqref{eq:intermediaryestimateformicrolocalconditiondegenerateMorrp1pdhorpR:robust} to the components of the wave system \eqref{eq:ScalarizedWaveeq:general:Kerrpert}}
\lab{subsubsec:directapplication:scalsyst}


In view of the assumptions of Theorem \ref{th:mainenergymorawetzmicrolocal}, the metric $\g$ and the scalars $\psi_{ij}$ verify the assumptions required in Lemma \ref{lem:conditionaldegenerateMorawetzflux:pertKerrrp1pdhorpR}, but the assumption that $\widehat{F}_{ij}$ are compactly supported in $\MM(\tmic, \tau_2)$ is not satisfied. On the other hand, the requirement that $F$ is compactly supported in $\MM(\tmic, \tau_2)$ is only used in the proof of Lemma \ref{lem:conditionaldegenerateMorawetzflux:pertKerrrp1pdhorpR} in \cite[Section 7]{MaSz24} to ensure the square integrability of $\pr^{\leq 1}\psi$ in $\tau\in\Reals$ in the region $\MM_{r\leq \Rmic}$, which in turn allows to justify the various computations involving PDOs with mixed symbols on $\MM_{r\leq \Rmic}$. For our model problem \eqref{eq:ScalarizedWaveeq:general:Kerrpert} with assumptions made in Theorem \ref{th:mainenergymorawetzmicrolocal}, it thus suffices to prove that $\pr^{\leq 1}\psi_{ij}$ are square integrable in $\tau\in\Reals$ in the region $\MM_{r\leq \Rmic}$, which indeed holds in view of the following:
\begin{itemize} 
\item First, in view of the assumptions of Theorem \ref{th:mainenergymorawetzmicrolocal}, $\psi_{ij}$ vanishes in $\MM(-\infty, \tmic)\cap\{r\leq (\Nmic+1)m\}$ so that $\pr^{\leq 1}\psi_{ij}$ is square integrable in $\MM_{r\leq \Rmic}(-\infty, \tmic)$. 
\item Next, square integrability of $\pr^{\leq 1}\psi_{ij}$ in time on $\MM(\tmic, \tau_2)$ follows immediately from local energy estimates for \eqref{eq:ScalarizedWaveeq:general:Kerrpert}.
\item Finally, since, in view of the assumptions of Theorem \ref{th:mainenergymorawetzmicrolocal}, $\psi_{ij}=\pmb\psi((\Om_K)_i, (\Om_K)_j)$ in $\MM(\tau_2, +\infty)$ for a tensor $\pmb\psi\in\sk_{2,K}(\mathbb{C})$ and $F_{ij}$ are supported in $\MM(\tmic,\tau_2)$, we have from Lemma \ref{lemma:formoffirstordertermsinscalarazationtensorialwaveeq} that the tensor $\pmb\psi\in\sk_{2,K}(\mathbb{C})$ satisfies the tensorial wave equation \eqref{eq:basictensorialwaveequationinKerr:forbrevephi} in a subextremal Kerr spacetime $\MM(\tau_2, +\infty)$ which, by a direct generalization of the estimates in Theorem \ref{prop:weakMorawetzfortensorialwaveeqfromDRSR} to higher derivatives, implies that  $\pr^{\leq 1}\pmb\psi$ (and hence $\pr^{\leq 1}\psi_{ij}$) are square integrable in $\MM_{r\leq \Rmic}(\tau_2, +\infty)$. 
\end{itemize}

As a consequence of the above square integrability property together with the assumptions of Theorem \ref{th:mainenergymorawetzmicrolocal} that the metric $\g$ coincides with $\gam$ for $\tau\in\Reals\setminus(\tau_1, \tau_2)$ and that  $\psi_{ij}$ vanishes in $\MM(-\infty, \tmic)\cap\{r\leq (\Nmic+1)m\}$, we may now apply \eqref{eq:intermediaryestimateformicrolocalconditiondegenerateMorrp1pdhorpR:robust} to each wave equation of the system \eqref{eq:ScalarizedWaveeq:general:Kerrpert} and sum up over $i,j=1,2,3$, arriving at
\bea\lab{eq:microlocalMora:middle:scalarizedwave:1}
\nn&& \sum_{i,j}\Bigg(c\Bigg[\int_{\MM_{r_+(1+\dhor'),{\Rmic}}}\frac{\mu^2|\pr_r\psi_{ij}|^2}{r^2} +\int_{\Mntrap_{r_+(1+\dhor'),{\Rmic}}}\frac{|\pr_\tau\psi_{ij}|^2+|\nab\psi_{ij}|^2}{r^2}\\
\nn&&+\int_{\MM_{r_+(1+\dhor'),10m}}\Big(|\Opw(\sigma_{\trap})\psi_{ij}|^2+|\Opw(e)\psi_{ij}|^2\Big)\Bigg]\\
&&+\textbf{BDR}^{-}_{r={\Rmic}}[\psi_{ij}](\Iti)+\int_{\MM_{r_+(1+\dhor'),{\Rmic}}}\Re\Big(\square_{\g}\psi_{ij}\ov{(X+E)\psi_{ij}}\Big)\Bigg)\nn\\
\nn&\les_{\Rmic}& \sum_{i,j}\bigg((\ep+\dhor)\int_{\Mntrap(\Iti)}|\square_{\g}\psi_{ij}|^2+\ep\int_{\Mtrap(\Iti)}\tau^{-1-\dec}|\square_{\g}\psi_{ij}|^2\\
&&
+(\ep+\dhor^6)\int_{\Mtrap}\Big|\Opw(\widetilde{S}^{-1,0}(\MM))\square_{\g}\psi_{ij}\Big|^2+\frac{1}{\dhor^6}\int_{\MM_{r_+(1+\dhor'),\Rmic}(\Iti)}|\psi_{ij}|^2\bigg)\nn \\
\nn&&+\ep\EM[\pmb \psi](\Iti)+\dhor\M[\pmb \psi](\Iti)\\
&&+\sum_{i,j}\bigg(\int_{H_{{\Rmic}}}\big(|\pr\psi_{ij}|^2+|\psi_{ij}|^2\big)\bigg)^{\frac{1}{2}}\bigg(\int_{H_{{\Rmic}}}|\psi_{ij}|^2\bigg)^{\frac{1}{2}}.
\eea

In order to prove microlocal Morawetz estimates for  \eqref{eq:ScalarizedWaveeq:general:Kerrpert} in $\MM_{r_+(1+\dhor'),{\Rmic}}$, stated later in Proposition \ref{prop:Morawetz:middleregion:globalextendedwavesystem}, we start with estimating the last integral term on the LHS of \eqref{eq:microlocalMora:middle:scalarizedwave:1}. By the choice of the PDO $X$ in Section \ref{sec:relevantmixedsymbolsonMM}, it can be decomposed into
\bea\lab{eq:decompositionofX:Mora}
X=X_1 + X_2 +A\pr_\tau,\qquad X_1:=\Opw(is_0\mu\xi_r),\qquad X_2:=\Opw(ib_{\tphi}\xiphi+ib_{\tt}\xit),
\eea
hence we decompose the integrand of the last integral term on the LHS of \eqref{eq:microlocalMora:middle:scalarizedwave:1} into
\bea
\lab{def:errorterms:general:scalarizedwavefromtensorialwave}
\bsplit
\Err:=&\sum_{i,j}\Re\Big(\square_{\g}\psi_{ij} \ov{(X+E)\psi_{ij}}\Big)=\Err_1+\Err_2+\Err_3+\Err_4,\\
\Err_1:=&\sum_{i,j}\Re\Big(\square_{\g}\psi_{ij} \ov{X_1\psi_{ij}}\Big),\qquad
\Err_2:=\sum_{i,j}\Re\Big(\square_{\g}\psi_{ij} \ov{X_2\psi_{ij}}\Big),\\
\Err_3:=&\sum_{i,j}\Re\Big(\square_{\g}\psi_{ij} \ov{A\pr_{\tt}\psi_{ij}}\Big),\qquad
\Err_4:=\sum_{i,j}\Re\Big(\square_{\g}\psi_{ij} \ov{E\psi_{ij}}\Big),
\end{split}
\eea
such that the last integral term on the LHS of \eqref{eq:microlocalMora:middle:scalarizedwave:1} is given by
\beaa
\sum_{i,j}\int_{\MM_{r_+(1+\dhor'),{\Rmic}}}\Re\Big(\square_{\g}\psi_{ij} \ov{(X+E)\psi_{ij}}\Big)=\int_{\MM_{r_+(1+\dhor'),{\Rmic}}} \Err=\sum_{n=1}^4\int_{\MM_{r_+(1+\dhor'),{\Rmic}}} \Err_n. 
\eeaa


\subsubsection{Control provided by the microlocal Morawetz norm $\widetilde{\M}[\psi]$}


In order to control the error terms in \eqref{def:errorterms:general:scalarizedwavefromtensorialwave}, we will in particular rely on the following two lemmas that identify terms controlled by the microlocal Morawetz norm $\widetilde{\M}[\psi]$. 
\begin{lemma}\lab{lemma:gardinginequalitiesyieldcontroloftermsbywidetileM}
Let the mixed symbols $e_0$, $s_0$, $\upsilon$, $b_{\tphi}$, $b_{\tt}$ and $e$ be as in Section  \ref{sec:relevantmixedsymbolsonMM}, and let $\widetilde{\M}[\psi]$ be as in \eqref{eq:definitionofmicrolocalMorawetznormwidetildeM}. Then, we have for any scalar function $\psi$ vanishing on $\MM(-\infty, \tmic)\cap\{r\leq 10m\}$ 
\bea\lab{eq:gardinginequalitiesyieldcontroloftermsbywidetileM:1}
&&\int_{\MM_{r_+(1+\dhor'), 10m}}\Big(|\Opw(e_0\upsilon)\psi|^2+|\Opw(s_0\upsilon)\psi|^2\nn\\
&&\qquad \qquad\qquad \quad+|\Opw(b_{\tphi}\upsilon)\psi|^2+|\Opw(b_{\tt}\upsilon)\psi|^2\Big)\lesssim \widetilde{\M}[\psi]
\eea
and\footnote{Note that $\lesssim$ appearing in \eqref{eq:gardinginequalitiesyieldcontroloftermsbywidetileM:2} depends on the choices of $b_1$ and $b_0$. We do not indicate this dependence as it will be used in practice for a finite number of choices of $b_1\in\widetilde{S}^{1,0}(\MM)$ and $b_0\in\widetilde{S}^{0,0}(\MM)$.} 
\bsub
\lab{eq:gardinginequalitiesyieldcontroloftermsbywidetileM:2}
\bea\lab{eq:gardinginequalitiesyieldcontroloftermsbywidetileM:2:1}
&&\int_{\MM_{r_+(1+\dhor'), 10m}}\big(|\Opw(\{b_{\tphi}, b_1\}\upsilon)\psi|^2+|\Opw(\{b_{\tt}, b_1\}\upsilon)\psi|^2\big)\nn\\
&\lesssim& \widetilde{\M}[\psi],\quad \forall b_1\in\widetilde{S}^{1,0}(\MM),
\eea
\bea\lab{eq:gardinginequalitiesyieldcontroloftermsbywidetileM:2:2}
&&\int_{\MM_{r_+(1+\dhor'), 10m}}\big(|\Opw(\{b_{\tphi}, b_0\}\upsilon^2)\psi|^2+|\Opw(\{b_{\tt}, b_0\}\upsilon^2)\psi|^2\big)\nn\\
&\lesssim & \widetilde{\M}[\psi],\quad \forall b_0\in\widetilde{S}^{0,0}(\MM).
\eea
\esub
\end{lemma}

\begin{proof}
Recalling, in view of Section  \ref{sec:relevantmixedsymbolsonMM}, that  $e$ and $\upsilon$, are in $\widetilde{S}^{1,0}(\MM)$, that $e_0$,  $s_0$, $b_{\tphi}$ and $b_{\tt}$ are in $\widetilde{S}^{0,0}(\MM)$, and that  $e$, $\upsilon$, $e_0$, $s_0$, $b_{\tphi}$ and $b_{\tt}$ are real-valued, and relying on \eqref{eq:pointwisecontrolbyesetuptogetGardingtypeinequalities}, there exists a constant $c>0$ small enough such that 
\beaa
e_1:=\sqrt{e^2-c\big((e_0\upsilon)^2+(s_0\upsilon)^2+(b_{\tphi}\upsilon)^2+(b_{\tt}\upsilon)^2\big)}, \qquad e_1\in\widetilde{S}^{1,0}(\MM).
\eeaa 
Together with items 1 and 4 of Proposition \ref{prop:PDO:MM:Weylquan:mixedoperators}, we infer 
\beaa
&& c\Big((\Opw(e_0\upsilon))^{\star}\Opw(e_0\upsilon)+(\Opw(s_0\upsilon))^{\star}\Opw(s_0\upsilon)+(\Opw(b_{\tphi}\upsilon))^{\star}\Opw(b_{\tphi}\upsilon)\nn\\
&&+(\Opw(b_{\tt}\upsilon))^{\star}\Opw(b_{\tt}\upsilon)\Big)+(\Opw(e_1))^{\star}\Opw(e_1)\\ 
&=& (\Opw(e))^{\star}\Opw(e)+\Opw(\widetilde{S}^{0,0}(\MM))
\eeaa
and hence, for any scalar function $\psi$ vanishing on $\MM(-\infty, \tmic)\cap\{r\leq 10m\}$
\beaa
&&\int_{\MM_{r_+(1+\dhor'), 10m}}\Big(|\Opw(e_0\upsilon)\psi|^2+|\Opw(s_0\upsilon)\psi|^2+|\Opw(b_{\tphi}\upsilon)\psi|^2+|\Opw(b_{\tt}\upsilon)\psi|^2\Big)d\Vref\\
&\lesssim& \int_{\MM_{r_+(1+\dhor'), 10m}}|\Opw(e)\psi|^2d\Vref+\int_{\MM_{r_+(1+\dhor'), 10m}}|\Opw(\widetilde{S}^{0,0}(\MM))\psi|^2d\Vref \\
&\lesssim& \int_{\MM_{r_+(1+\dhor'), 10m}}|\Opw(e)\psi|^2d\Vref+\int_{\MM_{r_+(1+\dhor'), 10m}(\Iti)}|\psi|^2d\Vref 
\eeaa
where we used Lemma \ref{lemma:actionmixedsymbolsSobolevspaces:MM} in the last estimate. Together with Lemma \ref{lemma:spacetimevolumeformusingisochorecoordinates} and the definition \eqref{eq:definitionofmicrolocalMorawetznormwidetildeM} of $\widetilde{\M}[\psi]$, we deduce 
\beaa
\int_{\MM_{r_+(1+\dhor'), 10m}}\Big(|\Opw(e_0\upsilon)\psi|^2+|\Opw(s_0\upsilon)\psi|^2+|\Opw(b_{\tphi}\upsilon)\psi|^2+|\Opw(b_{\tt}\upsilon)\psi|^2\Big)\lesssim \widetilde{\M}[\psi]
\eeaa
as stated in \eqref{eq:gardinginequalitiesyieldcontroloftermsbywidetileM:1}. The estimate \eqref{eq:gardinginequalitiesyieldcontroloftermsbywidetileM:2:1} is derived similarly. To show \eqref{eq:gardinginequalitiesyieldcontroloftermsbywidetileM:2:2}, it suffices to notice
\beaa
\{b_{\tphi}, b_0\}\upsilon=\{b_{\tphi}, b_0\upsilon\} - b_0 \{b_{\tphi}, \upsilon\}
\eeaa
and to apply \eqref{eq:gardinginequalitiesyieldcontroloftermsbywidetileM:2:1}, using also \eqref{eq:propWeylquantization:MM:composition:mixedsymbols} and Lemma \ref{lemma:actionmixedsymbolsSobolevspaces:MM} to deal with the second term.
This concludes the proof of Lemma \ref{lemma:gardinginequalitiesyieldcontroloftermsbywidetileM}.
\end{proof}

\begin{lemma}\lab{lemma:integrationbypartsforitimesfirstordertimesfirstorder:general:microlocalversion}
Let $d_1\in\widetilde{S}^{1,1}(\MM)$ be a symbol equal to $s_0\mu\xi_r$, $b_\tau\xi_\tau$ or $b_{\tphi}\xi_{\tphi}$. Also, let $\prtan$ be given by  \eqref{def:tangentialderivativeonHr:Kerrpert}. Then, for a smooth real-valued function $f$ with bounded derivatives, we have for any scalar function $\psi$ vanishing on $\MM(-\infty, \tmic)\cap\{r\leq R_0\}$
\bea\lab{eq:controlofsumijklintMMReifprtanpsiijovOpwidpsikl}
\nn&&\left|\int_{\MM_{r_+(1+\dhor'),\Rmic}(\Iti)}\Re\big(if\prtan(\psi)\ov{\Opw(id_1)\psi}\big)\right|\\
\nn&\les&  \left(\int_{\MM_{r_+(1+\dhor'),\Rmic}(\Iti)}|\psi|^2\right)^{\frac{1}{2}}\left(\widetilde{\M}[\psi]\right)^{\frac{1}{2}}\\
&&+\left(\int_{H_{r_+(1+\dhor')}(\Iti)\cup H_{\Rmic}(\Iti)}|\psi|^2\right)^{\frac{1}{2}}\left(\int_{H_{r_+(1+\dhor')}(\Iti)\cup H_{\Rmic}(\Iti)}|\pr^{\leq 1}\psi|^2\right)^{\frac{1}{2}}.
\eea
Similarly, for $A^{ij}$ a family of  smooth real-valued scalars antisymmetric w.r.t. $(i,j)$ with bounded derivatives, we have for any family of scalar functions  $\psi_i$ vanishing on $\MM(-\infty, \tmic)\cap\{r\leq R_0\}$
\bea\lab{eq:controlofsumijklintMMReifprtanpsiijovOpwidpsikl:caseantisymmatrix}
\nn&&\left|\int_{\MM_{r_+(1+\dhor'),\Rmic}(\Iti)}\Re\big(A^{ij}\prtan(\psi_i)\ov{\Opw(id_1)\psi_j}\big)\right|\\
\nn&\les&  \left(\sum_{i=1}^3\int_{\MM_{r_+(1+\dhor'),\Rmic}(\Iti)}|\psi_i|^2\right)^{\frac{1}{2}}\left(\sum_{i=1}^3\widetilde{\M}[\psi_i]\right)^{\frac{1}{2}}\\
&+&\left(\sum_{i=1}^3\int_{H_{r_+(1+\dhor')}(\Iti)\cup H_{\Rmic}(\Iti)}|\psi_i|^2\right)^{\frac{1}{2}}\left(\sum_{i=1}^3\int_{H_{r_+(1+\dhor')}(\Iti)\cup H_{\Rmic}(\Iti)}|\pr^{\leq 1}\psi_i|^2\right)^{\frac{1}{2}}.\hspace{1cm}
\eea
Finally, under the above assumptions for $A^{ij}$, $d_1$ and $\psi_i$, and assuming in addition that $A^{ij}$ vanishes on $\Mtrap$, we have 
\bea\lab{eq:controlofsumijklintMMReifprtanpsiijovOpwidpsikl:caseantisymmatrix:bis}
\nn&&\left|\int_{\MM_{r_+(1+\dhor'),\Rmic}(\Iti)}\Re\big(A^{ij}\pr_r(\psi_i)\ov{\Opw(id_1)\psi_j}\big)\right|\\
\nn&\les&  \left(\sum_{i=1}^3\int_{\MM_{r_+(1+\dhor'),\Rmic}(\Iti)}|\psi_i|^2\right)^{\frac{1}{2}}\left(\sum_{i=1}^3\M[\psi_i](\Iti)\right)^{\frac{1}{2}}\\
\nn&+& \left(\sum_{i=1}^3\int_{\MM_{r_+(1+\dhor'),\Rmic}(\Iti)}|\psi_i|^2\right)^{\frac{1}{2}}\left(\sum_{i=1}^3\int_{\Mntrap_{r_+(1+\dhor'),\Rmic}(\Iti)}|\Opw(\widetilde{S}^{-1,0}(\MM))\mu\pr_r^2\psi_i|^2\right)^{\frac{1}{2}}\\
&+&\left(\sum_{i=1}^3\int_{H_{r_+(1+\dhor')}(\Iti)\cup H_{\Rmic}(\Iti)}|\psi_i|^2\right)^{\frac{1}{2}}\left(\sum_{i=1}^3\int_{H_{r_+(1+\dhor')}(\Iti)\cup H_{\Rmic}(\Iti)}|\pr^{\leq 1}\psi_i|^2\right)^{\frac{1}{2}}.
\eea
\end{lemma}

\begin{remark}
The estimates \eqref{eq:controlofsumijklintMMReifprtanpsiijovOpwidpsikl},  \eqref{eq:controlofsumijklintMMReifprtanpsiijovOpwidpsikl:caseantisymmatrix} and \eqref{eq:controlofsumijklintMMReifprtanpsiijovOpwidpsikl:caseantisymmatrix:bis} correspond to the integration of a microlocal analog of the following differential identities
\bea
\lab{eq:integrationbypartsforitimesfirstordertimesfirstorder:general}
2\Re(i f \pr_{\a} \psi \ov{\pr_{\b}\psi}) &=& \pr_{\a}\big(\Re(i f \psi \ov{\pr_{\b}\psi})\big)-\pr_{\b}\big(\Re(i f \psi \ov{\pr_{\a}\psi})\big) \nn\\
&-& \Re(i \pr_{\a}f \psi \ov{\pr_{\b}\psi})+\Re(i \pr_{\b}f \psi \ov{\pr_{\a}\psi}), \,\,\,\,\text{for any real-valued function } f
\eea
and 
\bea
\lab{eq:integrationbypartsforitimesfirstordertimesfirstorder:general:caseantisymmatrix}
2\Re(A^{ij}\pr_{\a} \psi_i \ov{\pr_{\b}\psi_j}) &=& \pr_{\a}\big(\Re(A^{ij}\psi_i \ov{\pr_{\b}\psi_j})\big)-\pr_{\b}\big(\Re(A^{ij} \psi_i \ov{\pr_{\a}\psi_j})\big) \nn\\
&-& \Re(\pr_{\a}(A^{ij}) \psi_i \ov{\pr_{\b}\psi_j})+\Re(\pr_{\b}(A^{ij})\psi_i\ov{\pr_{\a}\psi_j}), \,\,\,\,\text{for any antisymmetric}\nn\\
&&\qquad\qquad\qquad\qquad\qquad\quad\textrm{family of real-valued functions }A^{ij}
\eea
with the correspondance $\pr_\a\leftrightarrow (\prtan,\,\pr_r)$ and $\pr_\b\leftrightarrow\Opw(id_1)$.
\end{remark}

\begin{proof}
We start with the proof of \eqref{eq:controlofsumijklintMMReifprtanpsiijovOpwidpsikl}. We first integrate by parts w.r.t. $\prtan$, using also \eqref{def:dVref}, and obtain\footnote{Note that integration by parts w.r.t. $\prtan$ does not generate boundary terms on $\pr(\MM_{r_+(1+\dhor'),\Rmic}(\Iti))$ since $\prtan$ is tangent to the level hypersurfaces of $r$ and since $\psi$ and its derivatives vanish on $\{\tau=\tmic\}\cap\{r\leq R_0\}$.}
\beaa
&&\left|\int_{\MM_{r_+(1+\dhor'),\Rmic}(\Iti)}\Re\big(if\prtan(\psi)\ov{\Opw(id_1)\psi}\big) + \int_{\MM_{r_+(1+\dhor'),\Rmic}(\Iti)}\Re\big(if\psi\ov{\prtan\Opw(id_1)\psi}\big)\right|\\
&\les& \left|\int_{\MM_{r_+(1+\dhor'),\Rmic}(\Iti)}\Re\big(i\prtan(ff_0)\psi\ov{\Opw(id_1)\psi}\big)f_0^{-1}\right|\\
&\les& \left(\int_{\MM_{r_+(1+\dhor'),\Rmic}(\Iti)}|\psi|^2\right)^{\frac{1}{2}}\left(\int_{\MM_{r_+(1+\dhor'),\Rmic}(\Iti)}|
\Opw(id_1)\psi|^2\right)^{\frac{1}{2}}.
\eeaa
Since $d_1$ is equal to $s_0\mu\xi_r$, $b_\tau\xi_\tau$ or $b_{\tphi}\xi_{\tphi}$, we infer from Lemma \ref{lemma:gardinginequalitiesyieldcontroloftermsbywidetileM}
\beaa
&&\left|\int_{\MM_{r_+(1+\dhor'),\Rmic}(\Iti)}\Re\big(if\prtan(\psi)\ov{\Opw(id_1)\psi}\big) + \int_{\MM_{r_+(1+\dhor'),\Rmic}(\Iti)}\Re\big(if\psi\ov{\prtan\Opw(id_1)\psi}\big)\right|\\
&\les& \left(\int_{\MM_{r_+(1+\dhor'),\Rmic}(\Iti)}|\psi|^2\right)^{\frac{1}{2}}\left(\widetilde{\M}[\psi]\right)^{\frac{1}{2}}
\eeaa
and hence
\beaa
&&\left|\int_{\MM_{r_+(1+\dhor'),\Rmic}(\Iti)}\Re\big(if\prtan(\psi)\ov{\Opw(id_1)\psi}\big) +\int_{\MM_{r_+(1+\dhor'),\Rmic}(\Iti)}\Re\big(if\psi\ov{\Opw(id_1)\prtan\psi}\big)\right|\\
&\les& \left|\int_{\MM_{r_+(1+\dhor'),\Rmic}(\Iti)}\Re\big(if\psi\ov{[\prtan,\Opw(id_1)]\psi}\big)\right|+\left(\int_{\MM_{r_+(1+\dhor'),\Rmic}(\Iti)}|\psi|^2\right)^{\frac{1}{2}}\left(\widetilde{\M}[\psi]\right)^{\frac{1}{2}}.
\eeaa
Next, using \eqref{def:dVref}, the fourth item of Proposition \ref{prop:PDO:MM:Weylquan:mixedoperators} and the fact that $d_1$ is real-valued, we have
\beaa
&&\int_{\MM_{r_+(1+\dhor'),\Rmic}(\Iti)}\Re\big(if\psi\ov{\Opw(id_1)\prtan\psi}\big)\\ 
&=& \int_{\MM_{r_+(1+\dhor'),\Rmic}(\Iti)}\Re\big(if\psi\ov{\Opw(id_1)\prtan\psi}\big)f_0d\Vref\\
&=&  -\int_{\MM_{r_+(1+\dhor'),\Rmic}(\Iti)}\Re\big(i\Opw(id_1)(f_0f\psi)\ov{\prtan\psi}\big)d\Vref +\BB\\
&=&  -\int_{\MM_{r_+(1+\dhor'),\Rmic}(\Iti)}\Re\big(if\Opw(id_1)\psi\ov{\prtan\psi}\big) \\
&& -\int_{\MM_{r_+(1+\dhor'),\Rmic}(\Iti)}\Re\big(i[\Opw(id_1), f_0f]\psi)\ov{\prtan\psi}\big)d\Vref +\BB\\
&=&  \int_{\MM_{r_+(1+\dhor'),\Rmic}(\Iti)}\Re\big(if\prtan(\psi)\ov{\Opw(id_1)\psi}\big)\\
&&  -\int_{\MM_{r_+(1+\dhor'),\Rmic}(\Iti)}\Re\big(i[\Opw(id_1), f_0f]\psi)\ov{\prtan\psi}\big)f_0^{-1} +\BB,
\eeaa
where $\BB$ denotes boundary terms on $H_{r_+(1+\dhor')}\cup H_{\Rmic}(\Iti)$ generated by the integration by parts of $\Opw(id_1)$, which, using Cauchy-Schwarz and Lemma \ref{lemma:actionmixedsymbolsSobolevspaces:MM}, can be estimated by
\bea\lab{eq:intermediarystep:lemma:integrationbypartsforitimesfirstordertimesfirstorder:general:bourdarytermsIBP}
\nn|\BB| &\les& \int_{H_{r_+(1+\dhor')}(\Iti)\cup H_{\Rmic}(\Iti)}|\Opw(\widetilde{S}^{0,0}(\MM))\psi||\pr^{\leq 1}\psi|\\
&\les& \left(\int_{H_{r_+(1+\dhor')}(\Iti)\cup H_{\Rmic}(\Iti)}|\psi|^2\right)^{\frac{1}{2}}\left(\int_{H_{r_+(1+\dhor')}(\Iti)\cup H_{\Rmic}(\Iti)}|\pr^{\leq 1}\psi|^2\right)^{\frac{1}{2}}.
\eea
Plugging in the above, we deduce
\bea\lab{eq:intermediarystep:lemma:integrationbypartsforitimesfirstordertimesfirstorder:general:microlocalversion}
\nn&&\left|\int_{\MM_{r_+(1+\dhor'),\Rmic}(\Iti)}\Re\big(if\prtan(\psi)\ov{\Opw(id_1)\psi}\big)\right|\\
\nn&\les& \left|\int_{\MM_{r_+(1+\dhor'),\Rmic}(\Iti)}\Re\big(if\psi\ov{[\prtan,\Opw(id_1)]\psi}\big)\right|\\
\nn&&+\left|\int_{\MM_{r_+(1+\dhor'),\Rmic}(\Iti)}\Re\big(i[\Opw(id_1), f_0f]\psi)\ov{\prtan\psi}\big)f_0^{-1}\right|\\
&&+\left(\int_{\MM_{r_+(1+\dhor'),\Rmic}(\Iti)}|\psi|^2\right)^{\frac{1}{2}}\left(\widetilde{\M}[\psi]\right)^{\frac{1}{2}}+|\BB|.
\eea

Next, we estimate the first two terms on the RHS of \eqref{eq:intermediarystep:lemma:integrationbypartsforitimesfirstordertimesfirstorder:general:microlocalversion}. Since $d_1\in\widetilde{S}^{1,1}(\MM)$, using the first item of Proposition \ref{prop:PDO:MM:Weylquan:mixedoperators} and the fact that $\sigma(\prtan)\in\widetilde{S}^{1,0}(\MM)$, we obtain 
\beaa
\,[\prtan,\Opw(id_1)] &=& \Opw(\{\sigma(\prtan), d_1\})+\Opw(\widetilde{S}^{-1,1}(\MM)), \quad \\
\,[\Opw(id_1), f_0f] &=& \Opw(\{d_1, ff_0\})+\Opw(\widetilde{S}^{-2,1}(\MM)),
\eeaa
which in view of \eqref{eq:intermediarystep:lemma:integrationbypartsforitimesfirstordertimesfirstorder:general:microlocalversion} yields 
\bea\lab{eq:intermediarystep:lemma:integrationbypartsforitimesfirstordertimesfirstorder:general:microlocalversion:1}
\nn&&\left|\int_{\MM_{r_+(1+\dhor'),\Rmic}(\Iti)}\Re\big(if\prtan(\psi)\ov{\Opw(id_1)\psi}\big)\right|\\
\nn&\les& \left|\int_{\MM_{r_+(1+\dhor'),\Rmic}(\Iti)}\Re\big(if\psi\ov{\Opw(\{\sigma(\prtan), d_1\})\psi}\big)\right|\\
\nn&&+\left|\int_{\MM_{r_+(1+\dhor'),\Rmic}(\Iti)}\Re\big(i\Opw(\{d_1, ff_0\})\psi)\ov{\prtan\psi}\big)f_0^{-1}\right|\\
&&+\left(\int_{\MM_{r_+(1+\dhor'),\Rmic}(\Iti)}|\psi|^2\right)^{\frac{1}{2}}\left(\widetilde{\M}[\psi]\right)^{\frac{1}{2}}+|\BB|.
\eea

Next, we evaluate the Poisson brackets on the RHS of \eqref{eq:intermediarystep:lemma:integrationbypartsforitimesfirstordertimesfirstorder:general:microlocalversion:1}. We recall that $d_1$ is equal to $s_0\mu\xi_r$, $b_\tau\xi_\tau$ or $b_{\tphi}\xi_{\tphi}$, and we start with the case $d_1=s_0\mu\xi_r$ for which we have 
\beaa
\{\sigma(\prtan), s_0\mu\xi_r\} &=& \{\sigma(\prtan), s_0\}\mu\xi_r+\{\sigma(\prtan), \mu\xi_r\}s_0 = \widetilde{S}^{0,0}(\MM)\mu\xi_r+\widetilde{S}^{1,0}(\MM)s_0,\\
\{s_0\mu\xi_r, ff_0\} &=& \{s_0, ff_0\}\mu\xi_r+\{\mu\xi_r, ff_0\}s_0=\widetilde{S}^{-1,0}(\MM)\mu\xi_r+\widetilde{S}^{0,0}(\MM)s_0,
\eeaa
and hence, using the second item of Proposition \ref{prop:PDO:MM:Weylquan:mixedoperators}, we deduce 
\beaa
\Opw\big(\{\sigma(\prtan), s_0\mu\xi_r\}\big) &=& \Opw(\widetilde{S}^{0,0}(\MM))\mu\pr_r+\Opw(\widetilde{S}^{1,0}(\MM))\circ\Opw(s_0)+\Opw(\widetilde{S}^{0,1}(\MM)),\\
\Opw\big(\{s_0\mu\xi_r, ff_0\}\big) &=& \Opw(\widetilde{S}^{-1,0}(\MM))\mu\pr_r+\Opw(\widetilde{S}^{0,0}(\MM))\circ\Opw(s_0)\\
&& +\Opw(\widetilde{S}^{-1,1}(\MM)),
\eeaa
which, together with Lemma \ref{lemma:gardinginequalitiesyieldcontroloftermsbywidetileM}, implies 
\bea\lab{eq:controlofthe2termsinvolvingPoinssonbracketofsigmaprtanwithdanddwithff0:cased=s0muxir}
\nn&& \left|\int_{\MM_{r_+(1+\dhor'),\Rmic}(\Iti)}\Re\big(if\psi\ov{\Opw(\{\sigma(\prtan), d_1\})\psi}\big)\right|\\
\nn&&+\left|\int_{\MM_{r_+(1+\dhor'),\Rmic}(\Iti)}\Re\big(i\Opw(\{d_1, ff_0\})\psi)\ov{\prtan\psi}\big)f_0^{-1}\right|\\
&\les& \left(\int_{\MM_{r_+(1+\dhor'),\Rmic}(\Iti)}|\psi|^2\right)^{\frac{1}{2}}\left(\widetilde{\M}[\psi]\right)^{\frac{1}{2}}, \quad\textrm{if}\quad  d_1=s_0\mu\xi_r.
\eea

Next, we evaluate the Poisson brackets on the RHS of \eqref{eq:intermediarystep:lemma:integrationbypartsforitimesfirstordertimesfirstorder:general:microlocalversion:1} in the case $d_1=b_\tau\xi_\tau$ and $d=b_{\tphi}\xi_{\tphi}$. In this case, we have
\beaa
\{\sigma(\prtan), d_1\} &=& \{b_{\a}, b_1\}\widetilde{S}^{1,0}(\MM)+b_\a\widetilde{S}^{1,0}(\MM), \quad \a=\tau,\, \tphi,\quad b_1\in\widetilde{S}^{1,0}(\MM), \\
\{d_1, ff_0\} &=&  \{b_{\a}, b_0\}\widetilde{S}^{1,0}(\MM)+b_\a\widetilde{S}^{0,0}(\MM), \quad \a=\tau,\, \tphi, \quad b_0\in\widetilde{S}^{0,0}(\MM),
\eeaa
which together with Lemma \ref{lemma:gardinginequalitiesyieldcontroloftermsbywidetileM} yields
\bea\lab{eq:controlofthe2termsinvolvingPoinssonbracketofsigmaprtanwithdanddwithff0:cased=baxiafora=tauortphi}
\nn&& \left|\int_{\MM_{r_+(1+\dhor'),\Rmic}(\Iti)}\Re\big(if\psi\ov{\Opw(\{\sigma(\prtan), d_1\})\psi}\big)\right|\\
\nn&&+\left|\int_{\MM_{r_+(1+\dhor'),\Rmic}(\Iti)}\Re\big(i\Opw(\{d_1, ff_0\})\psi)\ov{\prtan\psi}\big)f_0^{-1}\right|\\
&\les& \left(\int_{\MM_{r_+(1+\dhor'),\Rmic}(\Iti)}|\psi|^2\right)^{\frac{1}{2}}\left(\widetilde{\M}_{r_+(1+\dhor'),\Rmic}[\psi]\right)^{\frac{1}{2}}, \,\,\,\textrm{if}\,\,\,  d_1=b_\tau\xi_\tau\,\,\,\textrm{or}\,\,\, d_1=b_{\tphi}\xi_{\tphi}.
\eea
Now, plugging \eqref{eq:controlofthe2termsinvolvingPoinssonbracketofsigmaprtanwithdanddwithff0:cased=s0muxir} \eqref{eq:controlofthe2termsinvolvingPoinssonbracketofsigmaprtanwithdanddwithff0:cased=baxiafora=tauortphi} on the RHS of \eqref{eq:intermediarystep:lemma:integrationbypartsforitimesfirstordertimesfirstorder:general:microlocalversion:1}, and using the control of $\BB$ in \eqref{eq:intermediarystep:lemma:integrationbypartsforitimesfirstordertimesfirstorder:general:bourdarytermsIBP}, we obtain 
\beaa
\nn&&\left|\int_{\MM_{r_+(1+\dhor'),\Rmic}(\Iti)}\Re\big(if\prtan(\psi)\ov{\Opw(id_1)\psi}\big)\right| \\
&\les&  \left(\int_{\MM_{r_+(1+\dhor'),\Rmic}(\Iti)}|\psi|^2\right)^{\frac{1}{2}}\left(\widetilde{\M}[\psi]\right)^{\frac{1}{2}}\\
&&+ \left(\int_{H_{r_+(1+\dhor')}(\Iti)\cup H_{\Rmic}(\Iti)}|\psi|^2\right)^{\frac{1}{2}}\left(\int_{H_{r_+(1+\dhor')}(\Iti)\cup H_{\Rmic}(\Iti)}|\pr^{\leq 1}\psi|^2\right)^{\frac{1}{2}}
\eeaa
as stated in \eqref{eq:controlofsumijklintMMReifprtanpsiijovOpwidpsikl}. 

The proof of \eqref{eq:controlofsumijklintMMReifprtanpsiijovOpwidpsikl:caseantisymmatrix} is completely analogous to \eqref{eq:controlofsumijklintMMReifprtanpsiijovOpwidpsikl} and left to the reader. Finally, we focus on the proof of \eqref{eq:controlofsumijklintMMReifprtanpsiijovOpwidpsikl:caseantisymmatrix:bis}. Proceeding similarly to the proof of \eqref{eq:controlofsumijklintMMReifprtanpsiijovOpwidpsikl} and using the fact that $A^{ij}$ vanishes on $\Mtrap$ (which makes the argument simpler), we obtain the following analog of \eqref{eq:intermediarystep:lemma:integrationbypartsforitimesfirstordertimesfirstorder:general:microlocalversion}
\bea\lab{eq:intermediarystep:lemma:integrationbypartsforitimesfirstordertimesfirstorder:general:microlocalversion:a}
\nn&&\left|\int_{\MM_{r_+(1+\dhor'),\Rmic}(\Iti)}\Re\big(A^{ij}\pr_r(\psi_i)\ov{\Opw(id_1)\psi_j}\big)\right|\\
\nn&\les& \left|\int_{\MM_{r_+(1+\dhor'),\Rmic}(\Iti)}\Re\big(A^{ij}\psi_i\ov{[\pr_r,\Opw(id_1)]\psi_j}\big)\right|\\
\nn&&+\left|\int_{\MM_{r_+(1+\dhor'),\Rmic}(\Iti)}\Re\big([\Opw(id_1), A^{ij}f_0]\psi_i\ov{\pr_r\psi_j}\big)f_0^{-1}\right|\\
&&+\left(\sum_{i=1}^3\int_{\MM_{r_+(1+\dhor'),\Rmic}(\Iti)}|\psi_i|^2\right)^{\frac{1}{2}}\left(\sum_{i=1}^3\M[\psi_i](\Iti)\right)^{\frac{1}{2}}+|\BB|,
\eea
where the boundary terms $\BB$ satisfy \eqref{eq:intermediarystep:lemma:integrationbypartsforitimesfirstordertimesfirstorder:general:bourdarytermsIBP}. Now, since $d_1$ is equal to $s_0\mu\xi_r$, $b_\tau\xi_\tau$ or $b_{\tphi}\xi_{\tphi}$, using the first two items of Proposition \ref{prop:PDO:MM:Weylquan:mixedoperators}, we have
\beaa
\,[\pr_r,\Opw(id_1)] &=& \Opw(\widetilde{S}^{1,1}(\MM)), \\ 
\,[\Opw(id_1), A^{ij}f_0] &=& [\Opw(s_0), A^{ij}f_0]\mu\pr_r+\Opw(\widetilde{S}^{0,0}(\MM)).
\eeaa
Plugging in \eqref{eq:intermediarystep:lemma:integrationbypartsforitimesfirstordertimesfirstorder:general:microlocalversion:a} and using Lemma \ref{lemma:actionmixedsymbolsSobolevspaces:MM}, we infer
\beaa
\nn&&\left|\int_{\MM_{r_+(1+\dhor'),\Rmic}(\Iti)}\Re\big(A^{ij}\pr_r(\psi_i)\ov{\Opw(id_1)\psi_j}\big)\right|\\
\nn&\les& \sum_{i,j=1}^3\left|\int_{\MM_{r_+(1+\dhor'),\Rmic}(\Iti)}\Re\big([\Opw(s_0), A^{ij}f_0]\mu\pr_r(\psi_i)\ov{\pr_r\psi_j}\big)f_0^{-1}\right|\\
&&+\left(\sum_{i=1}^3\int_{\MM_{r_+(1+\dhor'),\Rmic}(\Iti)}|\psi_i|^2\right)^{\frac{1}{2}}\left(\sum_{i=1}^3\M[\psi_i](\Iti)\right)^{\frac{1}{2}}+|\BB|.
\eeaa
Integrating by parts in $\pr_r$ in the first term on the RHS, which generates additional boundary terms of the type $\BB$, and applying Cauchy-Schwarz, we deduce
\beaa
\nn&&\left|\int_{\MM_{r_+(1+\dhor'),\Rmic}(\Iti)}\Re\big(A^{ij}\pr_r(\psi_i)\ov{\Opw(id_1)\psi_j}\big)\right|\\
\nn&\les& \left(\sum_{i=1}^3\int_{\MM_{r_+(1+\dhor'),\Rmic}(\Iti)}|\psi_i|^2\right)^{\frac{1}{2}}\left(\sum_{i=1}^3\int_{\MM_{r_+(1+\dhor'),\Rmic}(\Iti)}\big|[\Opw(s_0), A^{ij}f_0]\mu\pr_r^2(\psi_i)\big|^2\right)^{\frac{1}{2}}\\
&&+\left(\sum_{i=1}^3\int_{\MM_{r_+(1+\dhor'),\Rmic}(\Iti)}|\psi_i|^2\right)^{\frac{1}{2}}\left(\sum_{i=1}^3\M[\psi_i](\Iti)\right)^{\frac{1}{2}}+|\BB|.
\eeaa
Finally, since $[\Opw(s_0), A^{ij}f_0]$ vanishes on $\Mtrap$ and $[\Opw(s_0), A^{ij}f_0]=\Opw(\widetilde{S}^{-1,0}(\MM))$, and using the estimate \eqref{eq:intermediarystep:lemma:integrationbypartsforitimesfirstordertimesfirstorder:general:bourdarytermsIBP} for $\BB$, we deduce \eqref{eq:controlofsumijklintMMReifprtanpsiijovOpwidpsikl:caseantisymmatrix:bis}. This concludes the proof of Lemma \ref{lemma:integrationbypartsforitimesfirstordertimesfirstorder:general:microlocalversion}.
\end{proof}


\subsubsection{Notations for lower order terms and error terms}
\lab{subsec:definitionoflowerordertermsintheerrors:scalsyst}


As our microlocal energy-Morawetz estimates are derived on $\MM_{r_+(1+\dhor'), \Rmic}$, where the constants $\dhor'$ and $\Rmic$ are introduced in Remark \ref{rmk:choiceofconstantRbymeanvalue}, it is convenient to introduce the following notation $\Gac$ for error terms 
\bea\lab{eq:decaypropertiesofGac:microlocalregion}
|\dk^{\leq 15}\Gac|\les \ep\tau^{-1-\dec}\qquad\textrm{on}\,\,\MM_{r_+(1+\dhor'), \Rmic},
\eea
so that we may replace $(\Ga_g, \Ga_b)$  by $\Gac$ in $\MM_{r_+(1+\dhor'), \Rmic}$ in view of \eqref{eq:decaypropertiesofGabGag}. Also, we introduce a notation for all tangential derivatives to $H_r$ 
\bea
\lab{def:tangentialderivativeonHr:Kerrpert}
\prtan:=\pr\setminus\{\pr_r\},
\eea
which allows to decompose $\pr_r^2\psi$ as follows, see \cite[(7.91)]{MaSz24}, 
\bea\lab{eq:decompositionofpr2psiinfunctionwaveandprprtan:microlocal}
\left(\frac{\De}{|q|^2}+\Gac\right)\pr_r^2\psi = \square_\g\psi+\Big(O(1)+\Gac\Big)\prtan\pr\psi+\Big(O(1)+\dk^{\leq 1}(\Gac)\Big)\pr\psi\quad\textrm{on}\,\,\MM_{r_+(1+\dhor'), \Rmic}.
\eea

Next, we introduce the notation $\good[\psi]$ for lower order terms and error terms that will appear when controlling the error terms in \eqref{def:errorterms:general:scalarizedwavefromtensorialwave}.

\begin{definition}[Notation ${\good[\psi]}$ for lower order terms and error terms]\lab{def:lotnotationforlowerordertermsinENGMorawetz}
Assume that $\dhor'$ and $\Rmic$ are chosen as in Remark \ref{rmk:choiceofconstantRbymeanvalue}. Then, we denote by $\good[\psi]$ lower order terms and error terms of the following form
\bea\lab{eq:def:lotnotationforlowerordertermsinENGMorawetz}
\bsplit
\good[\psi] :=& \sum_{n=1}^6 \good^{(n)}[\psi],\\
\good^{(1)}[\psi] :=& \sum_{i,j,k,l=1}^3 \left[\int_{H_r}|\Opw(\widetilde{S}^{0,0}(\MM))\psi_{ij}||\pr^{\leq 1}\psi_{kl}|\right]_{r=r_+(1+\dhor')}^{r={\Rmic}},\\
\good^{(2)}[\psi] :=& \sum_{i,j,k,l=1}^3  \int_{\MM_{r_+(1+\dhor'), {\Rmic}}}|\Opw(\widetilde{S}^{0,0}(\MM))\psi_{ij}|\big(|\pr_r^{\leq 1}\psi_{kl}|+|\pr\psi_{kl}|\mathbf{1}_{\Mntrap}\big),\\
\good^{(3)}[\psi] :=& \sum_{i,j,k,l=1}^3\int_{\MM_{r_+(1+\dhor'), {\Rmic}}}|\psi_{ij}||\Opw(s_0)\Opw(\widetilde{S}^{1,1}(\MM))\psi_{kl}|,\\
\good^{(4)}[\psi] :=& \sum_{i,j,k,l=1}^3\int_{\MM_{r_+(1+\dhor'), {\Rmic}}}\widecheck{\Ga}|\Opw(\widetilde{S}^{1,1}(\MM))\psi_{ij}| |\Opw(\widetilde{S}^{1,1}(\MM))\psi_{kl}|,\\
\good^{(5)}[\psi] :=& \sum_{i,j,k,l=1}^3\bigg|\int_{\MM_{r_+(1+\dhor'), {\Rmic}}}h_0\Re\left(\ov{\pr_r\psi_{ij}}\Opw(\widetilde{S}^{-1,0}(\MM))\mu\pr_r\psi_{kl}\right)\bigg|,\\
\good^{(6)}[\psi] :=& \sum_{i,j,k,l=1}^3\int_{\MM_{r_+(1+\dhor'), {\Rmic}}}|\Opw(\widetilde{S}^{0,0}(\MM))\psi_{ij}|\big(|\Opw(b_{\tphi})\pr^{\leq 1}\psi_{kl}|+|\Opw(b_{\tt})\pr^{\leq 1}\psi_{kl}|\big),
\end{split}
\eea
where the symbols $s_0$, $b_{\tphi}$ and $b_{\tt}$ are as in Section \ref{sec:relevantmixedsymbolsonMM}, and
where $h_0$ is any smooth scalar function such that $h_0=\Gac$ on $\Mtrap$.
\end{definition}

In order to control the error term $\good^{(4)}[\psi]$ in Lemma \ref{lem:lowerorderterms:controlled:scalarizedwavefromtensorial} below, we will rely on the following lemma which is taken from \cite[Lemma 7.15]{MaSz24}.

\begin{lemma}\label{lem:gpert:MMtrap}
Let  $h$ be a scalar function in $\MM_{r_1, r_2}$ supported in $\tau\geq 1$ with $r_+(1+\dhor)\leq r_1<r_2<+\infty$, let $S\in \Opw(\widetilde{S}^{1,1}(\MM))$, and let $\psi$ be supported on $\{\tt\geq 1\}$ in $\MM_{r_1,r_2}$. Then, for any $\dec>0$, we have 
\bea
\int_{\MM_{r_1, r_2}} |h|\abs{S\psi}^2 \les_{r_2, \dec}\|\tt^{1+\dec}h\|_{L^\infty(\MM_{r_1, r_2})}{\EM}_{r_1,r_2}[\psi](\Reals).
\eea
\end{lemma} 

The following lemma justifies that $\good[\psi]$ indeed corresponds to lower order terms and error terms.

\begin{lemma}
\lab{lem:lowerorderterms:controlled:scalarizedwavefromtensorial}
Assume that $\dhor'$ and $\Rmic$ are chosen as in Remark \ref{rmk:choiceofconstantRbymeanvalue}, and that the scalars $\psi_{ij}$ vanish in $\MM(-\infty, \tmic)\cap\{r\leq (\Nmic+1)m\}$. Then, the terms $\good[\psi]$ are bounded by
\bea
\good[\psi] &\les& \ep\EM[\pmb\psi](\Iti)
+\frac{1}{\sqrt{\dhor}}\left(\sum_{i,j=1}^3\int_{\MM_{r_+(1+\dhor'),\Rmic}(\Iti)}|\psi_{ij}|^2\right)^{\frac{1}{4}}\nn\\
&& \qquad\qquad\qquad\qquad\qquad\quad\times\left(\widetilde{\M}[\pmb\psi]+\sum_{i,j=1}^3\int_{\Mntrap(\Iti)}|\square_\g\psi_{ij}|^2\right)^{\frac{3}{4}}.
\eea
\end{lemma}

\begin{proof}
For convenience, we drop the indices $(i,j,k,l)$ appearing in \eqref{eq:def:lotnotationforlowerordertermsinENGMorawetz} throughout the proof since they do not play a role. We start with the boundary terms on $r=\dhor'$ and $\Rmic$, i.e., with $\good^{(1)}[\psi]$. Using Lemma \ref{lemma:actionmixedsymbolsSobolevspaces:MM}  and a trace estimate on $H_{r_+(1+\dhor')}$ and $H_{\Rmic}$ yields 
\beaa
 &&\int_{H_{r_+(1+\dhor')}(\Iti)}|\Opw(\widetilde{S}^{0,0}(\MM))\psi|^2+\int_{H_{\Rmic}(\Iti)}|\Opw(\widetilde{S}^{0,0}(\MM))\psi|^2\\ 
 &\les&\int_{H_{r_+(1+\dhor')}(\Iti)}|\psi|^2+\int_{H_{\Rmic}(\Iti)}|\psi|^2\\ 
 &\les& \left(\int_{\MM_{r_+(1+\dhor'), \Rmic}(\Iti)}\big(|\pr_r\psi|^2+|\psi|^2\big)\right)^{\frac{1}{2}}\left(\int_{\MM_{r_+(1+\dhor'), \Rmic}(\Iti)}|\psi|^2\right)^{\frac{1}{2}}\\
 &\les& \bigg(\int_{\MM_{r_+(1+\dhor'), \Rmic}(\Iti)}|\psi|^2\bigg)^{\frac{1}{2}}\big(\M[\psi](\Iti)\big)^{\frac{1}{2}}
\eeaa
which together with Cauchy-Schwarz and \eqref{de:choiceofdhor'} \eqref{eq:choiceofRvalue:Kerr} implies
\bea\lab{eq:lowerorderterms:controlled:scalarizedwavefromtensorial:proof:1}
\good^{(1)}[\psi] &=& \left[\int_{H_r}|\Opw(\widetilde{S}^{0,0}(\MM))\psi||\pr^{\leq 1}\psi|\right]_{r=r_+(1+\dhor')}^{r={\Rmic}}\nn\\ 
\nn&\les& \bigg(\int_{H_{r_+(1+\dhor')}(\Iti)\cup H_{\Rmic}(\Iti)}|\Opw(\widetilde{S}^{0,0}(\MM))\psi|^2\bigg)^{\frac{1}{2}}\\
\nn&&\times\left(\int_{H_{r_+(1+\dhor')}(\Iti)\cup H_{\Rmic}(\Iti)}|\pr^{\leq 1}\psi|^2\right)^{\frac{1}{2}}\\
&\les& \frac{1}{\sqrt{\dhor}}\bigg(\int_{\MM_{r_+(1+\dhor'), \Rmic}(\Iti)}|\psi|^2\bigg)^{\frac{1}{4}}\left(\M[\psi](\Iti)+\int_{\Mntrap(\Iti)}|\square_\g\psi|^2\right)^{\frac{3}{4}}.
\eea

Next, using Cauchy-Schwarz, we have
\bea\lab{eq:lowerorderterms:controlled:scalarizedwavefromtensorial:proof:2}
\nn\good^{(2)}[\psi] &=&\int_{\MM_{r_+(1+\dhor'), {\Rmic}}}|\Opw(\widetilde{S}^{0,0}(\MM))\psi|\big(|\pr_r^{\leq 1}\psi| +|\pr\psi_{kl}|\mathbf{1}_{\Mntrap}\big)\\
\nn&\les& \bigg(\int_{\MM_{r_+(1+\dhor'),\Rmic}}|\Opw(\widetilde{S}^{0,0}(\MM))\psi|^2\bigg)^{\frac{1}{2}}\big(\M[\psi](\Iti)\big)^{\frac{1}{2}}\\
&\les& \bigg(\int_{\MM_{r_+(1+\dhor'),\Rmic}}|\psi|^2\bigg)^{\frac{1}{2}}\big(\M[\psi](\Iti)\big)^{\frac{1}{2}}.
\eea

Next, using again Cauchy-Schwarz, together with \eqref{eq:gardinginequalitiesyieldcontroloftermsbywidetileM:1}, we have
\bea\lab{eq:lowerorderterms:controlled:scalarizedwavefromtensorial:proof:3}
\nn \good^{(3)}[\psi] &=& \int_{\MM_{r_+(1+\dhor'), {\Rmic}}}|\psi||\Opw(s_0)\Opw(\widetilde{S}^{1,1}(\MM))\psi| \\
\nn&\les&  \bigg(\int_{\MM_{r_+(1+\dhor'),\Rmic}(\Iti)}|\psi|^2\bigg)^{\frac{1}{2}}\bigg(\int_{\MM_{r_+(1+\dhor'),\Rmic}}|\Opw(s_0)\Opw(\widetilde{S}^{1,1}(\MM))\psi|^2\bigg)^{\frac{1}{2}}\\
&\les& \bigg(\int_{\MM_{r_+(1+\dhor'),\Rmic}(\Iti)}|\psi|^2\bigg)^{\frac{1}{2}}\Big(\widetilde{\M}[\psi]\Big)^{\frac{1}{2}}
\eea
where we also used the fact that $\upsilon$ controls all tangential derivatives to $H_r$ in view of its definition in \eqref{eq:precisedefinitionofsigmatrapandvarpi}. 

Next, for $\good^{(4)}[\psi]$, we apply Lemma \ref{lem:gpert:MMtrap} to find 
\bea\lab{eq:lowerorderterms:controlled:scalarizedwavefromtensorial:proof:7}
\good^{(4)}[\psi]\les \ep\EM[\pmb\psi](\Iti).
\eea

Next, consider $\good^{(6)}[\psi]$. We control only the part containing $b_{\tphi}$, and the other part containing $b_{\tt}$ can be estimated in a similar manner. Using again Cauchy-Schwarz, together with \eqref{eq:gardinginequalitiesyieldcontroloftermsbywidetileM:1}, we have
\bea\lab{eq:lowerorderterms:controlled:scalarizedwavefromtensorial:proof:6}
\nn \good^{(6)}[\psi] &=& \int_{\MM_{r_+(1+\dhor'), {\Rmic}}}|\Opw(\widetilde{S}^{0,0}(\MM))\psi| |\Opw(b_{\tphi})\pr^{\leq 1}\psi|\\
\nn&\les&  \bigg(\int_{\MM_{r_+(1+\dhor'),\Rmic}(\Iti)}|\psi|^2\bigg)^{\frac{1}{2}}\bigg(\int_{\MM_{r_+(1+\dhor'),\Rmic}}|\Opw(b_0)\pr^{\leq 1}\psi|^2\bigg)^{\frac{1}{2}}\\
&\les& \bigg(\int_{\MM_{r_+(1+\dhor'),\Rmic}(\Iti)}|\psi|^2\bigg)^{\frac{1}{2}}\Big(\widetilde{\M}[\psi]\Big)^{\frac{1}{2}}
\eea
where we also used the fact that $\upsilon$ controls all tangential derivatives to $H_r$ in view of its definition in \eqref{eq:precisedefinitionofsigmatrapandvarpi}. 

Finally, we consider $\good^{(5)}[\psi]$. First, we rewrite it as 
\beaa
\good^{(5)}[\psi] &=&\left|\int_{\MM_{r_+(1+\dhor'), {\Rmic}}}h_0\Re\left(\ov{\pr_r\psi}\Opw(\widetilde{S}^{-1,0}(\MM))\mu\pr_r\psi\right)\right|\\ 
&\les& \left|\int_{\MM_{r_+(1+\dhor'), {\Rmic}}}h_0\Re\left(\ov{\pr_r\psi}\Opw(\widetilde{S}^{-1,0}(\MM))\left(\frac{|q|^2}{\R}\left(\frac{\De}{|q|^2}+\Gac\right)\pr_r\psi\right)\right)\right|\\
&&+\ep\M[\psi](\Iti),
\eeaa
where we have used the control \eqref{eq:decaypropertiesofGac:microlocalregion} for $\Gac$. Then, using integration by parts in $\pr_r$, the control of boundary terms in \eqref{eq:lowerorderterms:controlled:scalarizedwavefromtensorial:proof:1}, Cauchy-Schwarz, Lemma \ref{lemma:actionmixedsymbolsSobolevspaces:MM}, and the fact that $h_0=\Gac$ on $\Mtrap$, we obtain
\beaa
&&\good^{(5)}[\psi]\\ 
&\les& \bigg(\int_{\MM_{r_+(1+\dhor'),\Rmic}}|\psi|^2\bigg)^{\frac{1}{2}}\left(\int_{\Mntrap_{r_+(1+\dhor'), {\Rmic}}}\left|\Opw(\widetilde{S}^{-1,0}(\MM))\left(\left(\frac{\De}{|q|^2}+\Gac\right)\pr_r^2\psi\right)\right|^2\right)^{\frac{1}{2}}\\
&&+\frac{1}{\sqrt{\dhor}}\bigg(\int_{\MM_{r_+(1+\dhor'),\Rmic}(\Iti)}|\psi|^2\bigg)^{\frac{1}{4}}\left(\M[\psi](\Iti)+\int_{\Mntrap(\Iti)}|\square_\g\psi|^2\right)^{\frac{3}{4}}+\ep\M[\psi](\Iti).
\eeaa
Now, in view of \eqref{eq:decompositionofpr2psiinfunctionwaveandprprtan:microlocal}, using also Lemma \ref{lemma:actionmixedsymbolsSobolevspaces:MM}, we have
\bea\lab{eq:tobeusedlater:estimateMntrapOpwSminus10muplusGacpr2ofpsiij}
\nn&&\int_{\Mntrap_{r_+(1+\dhor'), {\Rmic}}}\left|\Opw(\widetilde{S}^{-1,0}(\MM))\left(\left(\frac{\De}{|q|^2}+\Gac\right)\pr_r^2\psi\right)\right|^2\\ 
&\les&\int_{\Mntrap(\Iti)}|\square_\g\psi|^2+\M[\psi](\Iti),
\eea
and hence
\beaa
\good^{(5)}[\psi] \les \bigg(\int_{\MM_{r_+(1+\dhor'),\Rmic}(\Iti)}|\psi|^2\bigg)^{\frac{1}{4}}\left(\M[\psi](\Iti)+\int_{\Mntrap(\Iti)}|\square_\g\psi|^2\right)^{\frac{3}{4}} +\ep\M[\psi](\Iti).
\eeaa
Together with \eqref{eq:lowerorderterms:controlled:scalarizedwavefromtensorial:proof:1}, \eqref{eq:lowerorderterms:controlled:scalarizedwavefromtensorial:proof:2}, \eqref{eq:lowerorderterms:controlled:scalarizedwavefromtensorial:proof:3}, \eqref{eq:lowerorderterms:controlled:scalarizedwavefromtensorial:proof:7} and \eqref{eq:lowerorderterms:controlled:scalarizedwavefromtensorial:proof:6}, this concludes the proof of Lemma \ref{lem:lowerorderterms:controlled:scalarizedwavefromtensorial}.
\end{proof}

Finally, we have the following corollary of Lemma \ref{lemma:integrationbypartsforitimesfirstordertimesfirstorder:general:microlocalversion}.
\begin{corollary}\lab{cor:integrationbypartsforitimesfirstordertimesfirstorder:general:microlocalversion}
Let $d_1\in\widetilde{S}^{1,1}(\MM)$ be a symbol equal to $s_0\mu\xi_r$, $b_\tau\xi_\tau$ or $b_{\tphi}\xi_{\tphi}$. Also, let $f$ be any smooth real-valued function with bounded derivatives, let $A^{ij}$ be any family of  smooth real-valued scalars antisymmetric w.r.t. $(i,j)$ with bounded derivatives, 
and let $\prtan$ be given by  \eqref{def:tangentialderivativeonHr:Kerrpert}. Then, under the assumptions of Lemma \ref{lem:lowerorderterms:controlled:scalarizedwavefromtensorial}, we have 
\beaa
\nn&&\sum_{i,j=1}^3\left|\int_{\MM_{r_+(1+\dhor'),\Rmic}(\Iti)}\Re\big(if\prtan(\psi_{ij})\ov{\Opw(id_1)\psi_{ij}}\big)\right|\\
\nn&&+\sum_{i,j=1}^3\left|\int_{\MM_{r_+(1+\dhor'),\Rmic}(\Iti)}\Re\big(A_i^k\prtan(\psi_{ij})\ov{\Opw(id_1)\psi_{kj}}\big)\right|\\
&\les&  \frac{1}{\sqrt{\dhor}}\bigg(\sum_{i,j=1}^3\int_{\MM_{r_+(1+\dhor'), \Rmic}(\Iti)}|\psi_{ij}|^2\bigg)^{\frac{1}{4}}\left(\widetilde{\M}[\pmb\psi]+\sum_{i,j=1}^3\int_{\Mntrap(\Iti)}|\square_\g\psi_{ij}|^2\right)^{\frac{3}{4}}.
\eeaa
Moreover, if in addition $A^{ij}$ vanishes on $\Mtrap$, we have 
\beaa
\nn&&\left|\int_{\MM_{r_+(1+\dhor'),\Rmic}(\Iti)}\Re\big(A^{ij}\pr_r(\psi_i)\ov{\Opw(id_1)\psi_j}\big)\right|\\
\nn&\les& \frac{1}{\sqrt{\dhor}}\bigg(\sum_{i,j=1}^3\int_{\MM_{r_+(1+\dhor'), \Rmic}(\Iti)}|\psi_{ij}|^2\bigg)^{\frac{1}{4}}\left(\M[\pmb\psi](\Iti)+\sum_{i,j=1}^3\int_{\Mntrap(\Iti)}|\square_\g\psi_{ij}|^2\right)^{\frac{3}{4}}\\
&&+\ep\M[\pmb\psi](\Iti).
\eeaa
\end{corollary}

\begin{proof}
We apply \eqref{eq:controlofsumijklintMMReifprtanpsiijovOpwidpsikl} with $\psi\to\psi_{ij}$ and \eqref{eq:controlofsumijklintMMReifprtanpsiijovOpwidpsikl:caseantisymmatrix} with $\psi_i\to\psi_{ij}$ and sum over $(i,j)$ to obtain 
\beaa
\nn&&\sum_{i,j=1}^3\left|\int_{\MM_{r_+(1+\dhor'),\Rmic}(\Iti)}\Re\big(if\prtan(\psi_{ij})\ov{\Opw(id_1)\psi_{ij}}\big)\right|\\
\nn&&+\sum_{i,j=1}^3\left|\int_{\MM_{r_+(1+\dhor'),\Rmic}(\Iti)}\Re\big(A_i^k\prtan(\psi_{ij})\ov{\Opw(id_1)\psi_{kj}}\big)\right|\\
\nn&\les&  \left(\sum_{i,j=1}^3\int_{\MM_{r_+(1+\dhor'),\Rmic}(\Iti)}|\psi_{ij}|^2\right)^{\frac{1}{2}}\left(\widetilde{\M}[\pmb\psi]\right)^{\frac{1}{2}}\\
&&+\left(\sum_{i,j=1}^3\int_{H_{r_+(1+\dhor')}(\Iti)\cup H_{\Rmic}(\Iti)}|\psi_{ij}|^2\right)^{\frac{1}{2}}\left(\sum_{i,j=1}^3\int_{H_{r_+(1+\dhor')}(\Iti)\cup H_{\Rmic}(\Iti)}|\pr^{\leq 1}\psi_{ij}|^2\right)^{\frac{1}{2}}.
\eeaa
We then rely on \eqref{eq:lowerorderterms:controlled:scalarizedwavefromtensorial:proof:1} to control the last line of the RHS and we infer
\beaa
\nn&&\sum_{i,j=1}^3\left|\int_{\MM_{r_+(1+\dhor'),\Rmic}(\Iti)}\Re\big(if\prtan(\psi_{ij})\ov{\Opw(id_1)\psi_{ij}}\big)\right|\\
\nn&&+\sum_{i,j=1}^3\left|\int_{\MM_{r_+(1+\dhor'),\Rmic}(\Iti)}\Re\big(A_i^k\prtan(\psi_{ij})\ov{\Opw(id_1)\psi_{kj}}\big)\right|\\
\nn&\les&  \frac{1}{\sqrt{\dhor}}\bigg(\sum_{i,j=1}^3\int_{\MM_{r_+(1+\dhor'), \Rmic}(\Iti)}|\psi_{ij}|^2\bigg)^{\frac{1}{4}}\left(\widetilde{\M}[\pmb\psi]+\sum_{i,j=1}^3\int_{\Mntrap(\Iti)}|\square_\g\psi_{ij}|^2\right)^{\frac{3}{4}}
\eeaa
as stated in the first estimate. 

To prove the second estimate, we apply \eqref{eq:controlofsumijklintMMReifprtanpsiijovOpwidpsikl:caseantisymmatrix:bis} with $\psi_i\to\psi_{ij}$ and sum over $(i,j)$. Using \eqref{eq:lowerorderterms:controlled:scalarizedwavefromtensorial:proof:1} as above to control the boundary terms, we obtain
\beaa
\nn&&\left|\int_{\MM_{r_+(1+\dhor'),\Rmic}(\Iti)}\Re\big(A^{ij}\pr_r(\psi_i)\ov{\Opw(id_1)\psi_j}\big)\right|\\
\nn&\les&  \frac{1}{\sqrt{\dhor}}\bigg(\sum_{i,j=1}^3\int_{\MM_{r_+(1+\dhor'), \Rmic}(\Iti)}|\psi_{ij}|^2\bigg)^{\frac{1}{4}}\left(\widetilde{\M}[\pmb\psi]+\sum_{i,j=1}^3\int_{\Mntrap(\Iti)}|\square_\g\psi_{ij}|^2\right)^{\frac{3}{4}}\\
\nn&+& \left(\sum_{i=1,j}^3\int_{\MM_{r_+(1+\dhor'),\Rmic}(\Iti)}|\psi_{ij}|^2\right)^{\frac{1}{2}}\left(\sum_{i,j=1}^3\int_{\Mntrap_{r_+(1+\dhor'),\Rmic}(\Iti)}|\Opw(\widetilde{S}^{-1,0}(\MM))\mu\pr_r^2\psi_{ij}|^2\right)^{\frac{1}{2}}.\eeaa
Using an analog\footnote{That is, we first multiply \eqref{eq:decompositionofpr2psiinfunctionwaveandprprtan:microlocal} by $\frac{|q|^2}{\R}$ to obtain $(\mu+\Gac)\pr_r^2\psi_{ij}$ on the RHS and then proceed as in the proof of  \eqref{eq:tobeusedlater:estimateMntrapOpwSminus10muplusGacpr2ofpsiij}.} of \eqref{eq:tobeusedlater:estimateMntrapOpwSminus10muplusGacpr2ofpsiij} to control the last term on the RHS, we infer
\beaa
\nn&&\left|\int_{\MM_{r_+(1+\dhor'),\Rmic}(\Iti)}\Re\big(A^{ij}\pr_r(\psi_i)\ov{\Opw(id_1)\psi_j}\big)\right|\\
\nn&\les& \frac{1}{\sqrt{\dhor}}\bigg(\sum_{i,j=1}^3\int_{\MM_{r_+(1+\dhor'), \Rmic}(\Iti)}|\psi_{ij}|^2\bigg)^{\frac{1}{4}}\left(\M[\pmb\psi](\Iti)+\sum_{i,j=1}^3\int_{\Mntrap(\Iti)}|\square_\g\psi_{ij}|^2\right)^{\frac{3}{4}}\\
&&+\ep\M[\pmb\psi](\Iti)
\eeaa
as stated in the second estimate. This concludes the proof of Corollary \ref{cor:integrationbypartsforitimesfirstordertimesfirstorder:general:microlocalversion}.
\end{proof}


\subsubsection{Control of the error term $\Err_1$}
\lab{subsubsect:controlerrorterms:Morawetz:scalartoscalarized}


From now on, we estimate the integrals of the error terms in \eqref{def:errorterms:general:scalarizedwavefromtensorialwave}.
In this section, we estimate the error term $\Err_1$ associated to the PDO $X_1= \Opw(is_0\mu\xi_{r})$.

In view of the wave equations \eqref{eq:ScalarizedWaveeq:general:Kerrpert} for $\psi_{ij}$, and recalling \eqref{eq:definitionwidehatSandwidehatQperturbationsofKerr}, we have
\bea
\lab{eq:ScalarizedWaveeq:general:Kerrpert:rewriteinsect8}
{\square}_{\g}\psi_{ij}&=& S_K(\psi)_{ij} +\chi_{\tau_1, \tau_2}\big({S}(\psi)_{ij}-{S}_K(\psi)_{ij}\big) +\frac{4ia\cos\th}{|q|^2} \pr_{\tt}\psi_{ij}\nn\\
&&
+\Big(\chi_{\tau_1, \tau_2}(\widehat{Q}\psi)_{ij}+(1-\chi_{\tau_1, \tau_2})\big((\widehat{Q}_K\psi)_{ij}+f_{D_0}\psi_{ij}\big)+D_0|q|^{-2}\psi_{ij} \Big)\nn\\
&&+F_{ij} ,\quad i,j=1,2,3, \quad\textrm{on}\quad\MM,
\eea
thus,  we further decompose $\Err_1$ into
\bea
\lab{eq:definitionofErr1:decomposition}
\bsplit
\Err_1={}&\Err_{1,1} +\Err_{1,2} +\Err_{1,3}  +\Err_{1,4}+\sum_{i,j}\Re\Big(F_{ij}\ov{X_1\psi_{ij}}\Big),\\
\Err_{1,1}:={}&\sum_{i,j}\Re\Big( S_K(\psi)_{ij} \ov{X_1\psi_{ij}}\Big),\\
 \Err_{1,2}:={}&\sum_{i,j}\Re\bigg(\frac{4ia\cos\th}{|q|^2} \pr_{\tt}(\psi_{ij})\ov{X_1\psi_{ij}}\bigg),\\
\Err_{1,3}:={}&\sum_{i,j}\Re\Big(\Big(\chi_{\tau_1, \tau_2}(\widehat{Q}\psi)_{ij}+(1-\chi_{\tau_1, \tau_2})\big((\widehat{Q}_K\psi)_{ij}+f_{D_0}\psi_{ij}\big)+D_0|q|^{-2}\psi_{ij} \Big)\ov{X_1\psi_{ij}}\Big),\\
\Err_{1,4}:={}& \sum_{i,j}\Re\Big(\chi_{\tau_1, \tau_2}\big({S}(\psi)_{ij}-{S}_K(\psi)_{ij}\big)\ov{X_1\psi_{ij}}\Big).
\end{split}
\eea

In the remainder of Section \ref{sec:MorawetzestimatesinMMrplus1plusdhorprimRmicinhomscalarwavesystem}, unless otherwise stated, the integrals over the region $\MM_{r_+(1+\dhor'),{\Rmic}}$ will always be compressed in $\int$ for convenience.\\

\noindent\textbf{Step 1. Control of $\Err_{1,4}$.}  In view of \eqref{SandV} and the assumption \eqref{eq:assumptionsonregulartripletinperturbationsofKerr:0}, and since we can replace $(\Ga_g, \Ga_b)$ by $\widecheck{\Ga}$ on $\MM_{r_+(1+\dhor'), \Rmic}$, we have
\beaa
{S}(\psi)_{ij}-{S}_K(\psi)_{ij}=\sum_{k,l}\widecheck{\Ga}\pr\psi_{kl}, \quad \textrm{on} \quad \MM_{r_+(1+\dhor'), \Rmic},
\eeaa
from which we infer, by applying Lemma \ref{lem:gpert:MMtrap} to the integral of $\Err_{1,4}$, 
\bea
\lab{eq:estiofErr14}
\int |\Err_{1,4}| \les \ep \EM[\pmb\psi](\Iti).
\eea

\noindent\textbf{Step 2. Control of $\Err_{1,3}$.} Next, we consider the integral of the term $\Err_{1,3}$. Using Lemma \ref{lemma:actionmixedsymbolsSobolevspaces:MM}, and in view of the fact that $X_1=\Opw(\widetilde{S}^{0,0}(\MM))\mu\pr_r+\Opw(\widetilde{S}^{0,0}(\MM))$ which follows from the second item of Proposition \ref{prop:PDO:MM:Weylquan:mixedoperators} and the form of $X_1$, we have
\bea
\lab{eq:EstimateforErr13:proof:MorawetzRW}
 \int|\Err_{1,3}| &\les& \bigg(\sum_{i,j}\int|\psi_{ij}|^2\bigg)^{\frac{1}{2}} \bigg(\sum_{i,j}\int \mu^2|\pr_r\psi_{ij}|^2\bigg)^{\frac{1}{2}}+ \sum_{i,j}\int|\psi_{ij}|^2.
\eea

\noindent\textbf{Step 3. Control of $\Err_{1,2}$.} Next, we consider the integral of the term $\Err_{1,2}$. Using Corollary \ref{cor:integrationbypartsforitimesfirstordertimesfirstorder:general:microlocalversion} with $d_1=s_0\mu\xi_r$ and $f=\frac{4a\cos\th}{|q|^2}$, we have, in view of the form of $X_1$,
\bea\lab{eq:Err12:rewriteform}
\nn&&\bigg|\int\Err_{1,2}\bigg|\\
\nn&\leq& \sum_{i,j}\bigg|\int\Re\bigg(\frac{4ia\cos\th}{|q|^2} \pr_{\tt}(\psi_{ij})\ov{X_1\psi_{ij}}\bigg)\bigg|= \sum_{i,j}\bigg|\int\Re\bigg(if\pr_{\tt}(\psi_{ij})\ov{\Opw(id_1)\psi_{ij}}\bigg)\bigg|\\
&\les&  \frac{1}{\sqrt{\dhor}}\bigg(\sum_{i,j=1}^3\int_{\MM_{r_+(1+\dhor'), \Rmic}(\Iti)}|\psi_{ij}|^2\bigg)^{\frac{1}{4}}\left(\widetilde{\M}[\pmb\psi]+\sum_{i,j=1}^3\int_{\Mntrap(\Iti)}|\square_\g\psi_{ij}|^2\right)^{\frac{3}{4}}\!\!\!.
\eea

\noindent\textbf{Step 4. Control of $\Err_{1,1}$.}   In view of the assumptions of Theorem \ref{th:mainenergymorawetzmicrolocal}, $\psi_{ij}=\pmb\psi(\Om_i,\Om_j)$ in $\MM(\tau_1+1, \tau_2-2)$ for a tensor $\pmb\psi\in\sk_2(\mathbb{C})$, $\psi_{ij}=\pmb\psi((\Om_K)_i, (\Om_K)_j)$ in $\MM(\tau_2, +\infty)$ for a tensor $\pmb\psi\in\sk_{2,K}(\mathbb{C})$, and $\psi_{ij}$ vanishes in $\MM_{r_+(1+\dhor'), \Rmic}(-\infty, \tmic)$. Hence it follows from the first part of Lemma \ref{lemma:backandforthbetweenhorizontaltensorsk2andscalars:complex} that 
\bea
\lab{eq:widetildepsiij:tau2-4totau2}
\textrm{supp}(x^i\psi_{ij})\subset [\tmic,\tau_1+1]\cup[\tau_2-2, \tau_2], \quad \textrm{supp}(x^j\psi_{ij})\subset [\tmic,\tau_1+1]\cup[\tau_2-2, \tau_2].
\eea
Next, in view of \eqref{SandV}, we have
\bea
\lab{eq:Err11:firstdecomp:pf}
\int \Err_{1,1}=\sum_{i,j}\int \Re\Big(2(M_K)_{i}^{k\a}\pr_\a(\psi_{kj})\ov{X_1(\psi_{ij})}\Big)
+\sum_{i,j}\int \Re\Big(2(M_K)_{j}^{k\a}
\pr_\a(\psi_{ik})\ov{X_1(\psi_{ij})}\Big).
\eea
In the following, we control only the first term on the RHS of \eqref{eq:Err11:firstdecomp:pf}, the control of the second term on the RHS being completely analogous.
 Using \eqref{formula:symmetricpartofMmatrices} to decompose $(M_K)_{i}^{k\a}$ as
\bea
\lab{eq:decompositionofMikalpha:symandantisym}
(M_K)_{i}^{k\a}=(M_{K,S})_{i}^{k\a} +(M_{K,A})_{i}^{k\a}= -\frac{1}{2}\pr^{\a} (x^i x^k) +(M_{K,A})_{i}^{k\a}
\eea
where $(M_{K,S})_{i}^{k\a}$ and $(M_{K,A})_{i}^{k\a}$ denote respectively the symmetric and antisymmetric part of $(M_K)_{i}^{k\a}$ w.r.t. $(i,k)$, we deduce for the first term on the RHS of \eqref{eq:Err11:firstdecomp:pf}
\bea
\lab{eq:controlofErr11:onepart}
\bigg|\sum_{i,j}\int \Re\Big(2(M_K)_{i}^{k\a}\pr_\a(\psi_{kj})\ov{X_1(\psi_{ij})}\Big)\bigg|
\les\sum_{n=1}^3\Err_{1,1}^{(n)}
\eea
where
\bsub
\lab{def:Err11n:all}
\bea
\Err_{1,1}^{(1)}&=&\bigg|\sum_{i,j}\int \Re\Big((M_{K,A})_{i}^{k\a}\pr_\a(\psi_{kj})\ov{X_1(\psi_{ij})}\Big)\bigg|,\\
\Err_{1,1}^{(2)}&=&\bigg|\sum_{i,j}\int \Re\Big(\pr^{\a} (x^k) \pr_\a(\psi_{kj})\ov{x^i X_1(\psi_{ij})}\Big)\bigg|,\\
\Err_{1,1}^{(3)}&=&\bigg|\sum_{i,j}\int \Re\Big(\pr^{\a} (x^i ) (x^k\pr_\a\psi_{kj})\ov{X_1(\psi_{ij})}\Big)\bigg|.
\eea
\esub

Consider first the term $\Err_{1,1}^{(1)}$. Using the fact that $(M_{K,A})_{i}^{k\a}$ are antisymmetric w.r.t. $(i,k)$ and real-valued, we may apply Corollary \ref{cor:integrationbypartsforitimesfirstordertimesfirstorder:general:microlocalversion} with $A_i^k=(M_{K,A})_{i}^{k\a}$ and $d_1=s_0\mu\xi_r$ which yields\footnote{For tangential derivatives, i.e., for $(M_{K,A})_{i}^{k\a}$ with $x^\a=\tau, x^1, x^2$, we use the first estimate of Corollary \ref{cor:integrationbypartsforitimesfirstordertimesfirstorder:general:microlocalversion}, while for $(M_{K,A})_{i}^{kr}$, we use the second estimate of Corollary \ref{cor:integrationbypartsforitimesfirstordertimesfirstorder:general:microlocalversion} recalling that $(M_K)_{i}^{kr}$ (and hence $(M_{K,A})_{i}^{kr}$) vanishes on $\Mtrap$ in view of Lemma \ref{lemma:computationoftheMialphajinKerr}.} 
\bea
\lab{eq:finalesti:Err110}
\nn\Err_{1,1}^{(1)}&\les& \frac{1}{\sqrt{\dhor}}\bigg(\sum_{i,j=1}^3\int_{\MM_{r_+(1+\dhor'), \Rmic}(\Iti)}|\psi_{ij}|^2\bigg)^{\frac{1}{4}}\left(\widetilde{\M}[\pmb\psi]+\sum_{i,j=1}^3\int_{\Mntrap(\Iti)}|\square_\g\psi_{ij}|^2\right)^{\frac{3}{4}}\\
&&+\ep\M[\pmb\psi](\Iti).
\eea

Next, we consider $\Err_{1,1}^{(2)}$. Noticing that $\{x^i, \xi_r\}=-\pr_r(x^i)=0$, we have
\bea\lab{eq:precisecommutatorxiX1neededforcontrolofErr112}
\nn [x^i, X_1] &=& \Opw(\{x^i, s_0\mu\xi_r\})+\Opw(\widetilde{S}^{-2,1}(\MM))=\Opw(\{x^i, s_0\mu\}\xi_r)+\Opw(\widetilde{S}^{-2,1}(\MM))\\
&=& \Opw(\widetilde{S}^{-1,0}(\MM))\pr_r+\Opw(\widetilde{S}^{-1,1}(\MM))=\Opw(\widetilde{S}^{-1,0}(\MM))\pr_r^{\leq 1},
\eea
and hence, using again the fact that $\pr_r(x^i)=0$,
\beaa
\Err_{1,1}^{(2)}
&\les& \bigg|\sum_{i,j}\int \Re\Big(\pr^{\a} (x^k) \pr_\a(\psi_{kj})\ov{X_1(x^i\psi_{ij})}\Big)\bigg|\\
&&+\bigg|\sum_{i,j}\int \Re\Big(\pr^{\a} (x^k)\pr_\a(\psi_{kj})\ov{\Opw(\widetilde{S}^{-1,0}(\MM))\pr_r^{\leq 1}(\psi_{ij})}\Big)\bigg|\nn\\
&\les& \bigg|\sum_{i,j}\sum_{\a\neq r}\int \Re\Big(\pr^{\a} (x^k) \pr_\a(\psi_{kj})\ov{X_1(x^i\psi_{ij})}\Big)\bigg|\\
&&+\bigg|\sum_{i,j}\sum_{\a\neq r}\int \Re\Big(\Opw(\widetilde{S}^{-1,0}(\MM))\big(f_0\pr^{\a} (x^k)\pr_\a(\psi_{kj})\big)\ov{\pr_r^{\leq 1}(\psi_{ij})}\Big)f_0^{-1}\bigg|
\eeaa
where we also used \eqref{eq:spacetimevolumeformusingisochorecoordinates} and the fourth item of Proposition \ref{prop:PDO:MM:Weylquan:mixedoperators}. Now, we have
\beaa
\sum_{\a\neq r}\Opw(\widetilde{S}^{-1,0}(\MM))\big(f_0\pr^{\a} (x^k)\pr_\a(\psi_{kj})\big) =\Opw(\widetilde{S}^{0,0}(\MM))\psi_{kj}
\eeaa
and hence
\beaa
\Err_{1,1}^{(2)} &\les& \bigg|\sum_{i,j}\sum_{\a\neq r}\int \Re\Big(\pr^{\a} (x^k) \pr_\a(\psi_{kj})\ov{X_1(x^i\psi_{ij})}\Big)\bigg|+\good^{(2)}[\psi].
\eeaa
Next, we integrate by parts the first term on the RHS, first in $\pr_{\a}$ and next in $X_1$, to deduce
\bea
\lab{eq:finalesti:Err112:middlestep}
&&\bigg|\sum_{i,j}\sum_{\a\neq r}\int \Re\Big(\pr^{\a} (x^k) \pr_\a(\psi_{kj})\ov{X_1(x^i\psi_{ij})}\Big)\bigg|\nn\\
&\les& \bigg|\sum_{i,j}\sum_{\a\neq r}\int \Re\Big(\pr^{\a} (x^k) X_1(\psi_{kj})\ov{\pr_{\a}(x^i\psi_{ij})}\Big)\bigg| +\good^{(1)}[\psi]+ \good^{(2)}[\psi]+ \good^{(3)}[\psi]\nn\\
&\les& \big(\M[\pmb\psi](\Iti)\big)^{\frac{1}{2}}
\big(\Errdefect[\pmb\psi]\big)^{\frac{1}{2}} + \good[\psi]
\eea
where we have applied Cauchy-Schwarz and used \eqref{eq:widetildepsiij:tau2-4totau2} as well as the definition \eqref{def:Errdefectofpsi} of $\Errdefect[\pmb\psi]$ in the last step. Plugging \eqref{eq:finalesti:Err112:middlestep} in the previous estimate yields
\bea
\lab{eq:finalesti:Err112}
\Err_{1,1}^{(2)}
\les \big(\M[\pmb\psi](\Iti)\big)^{\frac{1}{2}}
\big(\Errdefect[\pmb\psi]\big)^{\frac{1}{2}}+\good[\psi].
\eea

For the term $\Err_{1,1}^{(3)}$, commuting $x^k$ with $\pr_\a$, we have
\beaa
\Err_{1,1}^{(3)}&\les&\bigg|\sum_{i,j}\int \Re\Big(\pr^{\a}(x^i )\pr_\a(x^k\psi_{kj})\ov{X_1(\psi_{ij})}\Big)\bigg|+\good^{(2)}[\psi].
\eeaa
Hence, taking Cauchy-Schwarz and using \eqref{eq:widetildepsiij:tau2-4totau2} as well as the definition \eqref{def:Errdefectofpsi} of $\Errdefect[\pmb\psi]$, we infer 
\bea
\lab{eq:finalesti:Err113}
\Err_{1,1}^{(3)} \les \big(\M[\pmb\psi](\Iti)\big)^{\frac{1}{2}}
\big(\Errdefect[\pmb\psi]\big)^{\frac{1}{2}}+\good[\psi].
\eea

In view of the above estimates \eqref{eq:finalesti:Err110}, \eqref{eq:finalesti:Err112} and \eqref{eq:finalesti:Err113}, as well as \eqref{eq:controlofErr11:onepart}, we infer
\bea
\lab{eq:controlofErr11}
\nn&&\bigg|\sum_{i,j}\int \Re\Big(2(M_K)_{i}^{k\a}\pr_\a(\psi_{kj})\ov{X_1(\psi_{ij})}\Big)\bigg|\\
\nn&\les& \frac{1}{\sqrt{\dhor}}\bigg(\sum_{i,j=1}^3\int_{\MM_{r_+(1+\dhor'), \Rmic}(\Iti)}|\psi_{ij}|^2\bigg)^{\frac{1}{4}}\left(\widetilde{\M}[\pmb\psi]+\sum_{i,j=1}^3\int_{\Mntrap(\Iti)}|\square_\g\psi_{ij}|^2\right)^{\frac{3}{4}}\\
&&+\big(\M[\pmb\psi](\Iti)\big)^{\frac{1}{2}}
\big(\Errdefect[\pmb\psi]\big)^{\frac{1}{2}}+\good[\psi]+\ep\M[\pmb\psi](\Iti),
\eea
and hence, the same bound holds for $|\int \Err_{1,1}|$. Together with the estimates \eqref{eq:definitionofErr1:decomposition}, \eqref{eq:estiofErr14}, \eqref{eq:EstimateforErr13:proof:MorawetzRW} and \eqref{eq:Err12:rewriteform}, and applying Lemma \ref{lem:lowerorderterms:controlled:scalarizedwavefromtensorial} to control $\good[\psi]$, we infer
\bea
\lab{esti:Err1:MoraEstimiddleregion}
&&\bigg|\int \bigg(\Err_{1}-\sum_{i,j}\Re\Big(F_{ij}\ov{X_1\psi_{ij}}\Big)\bigg)\bigg|\nn\\
&\les_{\Rmic}& \ep\EM[\pmb\psi](\Iti) +\A[\pmb \psi](\Iti)+
\big(\M[\pmb\psi](\Iti)\big)^{\frac{1}{2}}
\big(\Errdefect[\pmb\psi]\big)^{\frac{1}{2}}\nn\\
&&+ \frac{1}{\sqrt{\dhor}}\big(\A[\pmb \psi](\Iti)\big)^{\frac{1}{4}}\bigg(\widetilde{\M}[\pmb\psi]+\sum_{i,j=1}^3\int_{\Mntrap(\Iti)}|F_{ij}|^2\bigg)^{\frac{3}{4}},
\eea
with $\A[\pmb \psi](\Iti)$ given by \eqref{def:EM-1norms:Reals}, where we have also used the following consequence of \eqref{eq:ScalarizedWaveeq:general:Kerrpert}
\bea\lab{eq:controlofL2spacetimenormofsquaregpsiijbyL2spacetimenormofFijandMorawetzusingwaveeqsec8}
\nn\sum_{i,j=1}^3\int_{\Mntrap(\Iti)}|\square_\g\psi_{ij}|^2&\les& \sum_{i,j=1}^3\int_{\Mntrap(\Iti)}\big(r^{-4}|\dk^{\leq 1}\psi_{ij}|^2+|F_{ij}|^2\big)\\
&\les& \sum_{i,j=1}^3\int_{\Mntrap(\Iti)}|F_{ij}|^2+\M[\pmb\psi](\Iti).
\eea


\subsubsection{Control of the error term  $\Err_2$}
\lab{subsubsect:controlerrorterms:Morawetz:scalartoscalarized:Err2}


Next, we estimate the error term $\Err_2$ introduced in \eqref{def:errorterms:general:scalarizedwavefromtensorialwave} and associated to the PDO $X_2= \Opw(ib_{\tphi}\xiphi+ib_{\tt}\xit)$, with the symbols $b_{\tphi}$ and $b_{\tt}$ given as in Section \ref{sec:relevantmixedsymbolsonMM}. We follow closely the analysis on the control of the error term $\Err_1$ in the previous Section \ref{subsubsect:controlerrorterms:Morawetz:scalartoscalarized}.

Similarly to \eqref{eq:definitionofErr1:decomposition}, we decompose $\Err_2$ as follows
\bea
\lab{eq:definitionofErr2:decomposition}
\bsplit
\Err_2={}&\Err_{2,1} +\Err_{2,2} +\Err_{2,3}  +\Err_{2,4}+\sum_{i,j}\Re\Big(F_{ij}\ov{X_2\psi_{ij}}\Big),\\
\Err_{2,1}:={}&\sum_{i,j}\Re\Big(S_K(\psi)_{ij} \ov{X_2\psi_{ij}}\Big),\\
 \Err_{2,2}:={}&\sum_{i,j}\Re\bigg(\frac{4ia\cos\th}{|q|^2} \pr_{\tt}(\psi_{ij})\ov{X_2\psi_{ij}}\bigg),\\
\Err_{2,3}:={}&\sum_{i,j}\Re\Big(\Big(\chi_{\tau_1, \tau_2}(\widehat{Q}\psi)_{ij}+(1-\chi_{\tau_1, \tau_2})\big((\widehat{Q}_K\psi)_{ij}+f_{D_0}\psi_{ij}\big)+D_0|q|^{-2}\psi_{ij} \Big)\ov{X_2\psi_{ij}}\Big),\\
\Err_{2,4}:={}& \sum_{i,j}\Re\Big(\chi_{\tau_1, \tau_2}\big({S}(\psi)_{ij}-{S}_K(\psi)_{ij}\big)\ov{X_2\psi_{ij}}\Big).
\end{split}
\eea
Next, proceeding as for the proof of \eqref{eq:estiofErr14}, we have
\bea
\lab{eq:estiofErr24}
 \bigg|\int \Err_{2,4}\bigg| \les \ep \EM[\pmb\psi](\Iti).
\eea
Also, in view of the definition of $\good^{(6)}[\psi]$ in \eqref{eq:def:lotnotationforlowerordertermsinENGMorawetz}, we have 
\bea
\lab{eq:estiofErr23}
\bigg|\int \Err_{2,3}\bigg|\les \good^{(6)}[\psi].
\eea 
Moreover, using Corollary \ref{cor:integrationbypartsforitimesfirstordertimesfirstorder:general:microlocalversion} with $d_1=b_{\tphi}\xiphi+b_{\tt}\xit$ and $f=\frac{4a\cos\th}{|q|^2}$, we have, in view of the fact that $X_2= \Opw(ib_{\tphi}\xiphi+ib_{\tt}\xit)$,
\bea\lab{eq:estiofErr22}
\nn&&\bigg|\int\Err_{2,2}\bigg|\\
\nn&\leq& \sum_{i,j}\bigg|\int\Re\bigg(\frac{4ia\cos\th}{|q|^2} \pr_{\tt}(\psi_{ij})\ov{X_2\psi_{ij}}\bigg)\bigg|= \sum_{i,j}\bigg|\int\Re\bigg(if\pr_{\tt}(\psi_{ij})\ov{\Opw(id_1)\psi_{ij}}\bigg)\bigg|\\
&\les&  \frac{1}{\sqrt{\dhor}}\bigg(\sum_{i,j=1}^3\int_{\MM_{r_+(1+\dhor'), \Rmic}(\Iti)}|\psi_{ij}|^2\bigg)^{\frac{1}{4}}\left(\widetilde{\M}[\pmb\psi]+\sum_{i,j=1}^3\int_{\Mntrap(\Iti)}|\square_\g\psi_{ij}|^2\right)^{\frac{3}{4}}\!\!\!.
\eea

Next, we estimate the integral of the error term $\Err_{2,1}$ which we decompose as in \eqref{eq:Err11:firstdecomp:pf} and \eqref{eq:controlofErr11:onepart}. This leads us to consider, analogously to \eqref{def:Err11n:all}, the terms $\{\Err_{2,1}^{(n)}\}_{n=0,1,2,3}$ which are defined as follows 
\bsub
\lab{def:Err21n:all}
\bea
\Err_{2,1}^{(1)}&=&\bigg|\sum_{i,j}\int \Re\Big((M_{K,A})_{i}^{k\a}\pr_\a(\psi_{kj})\ov{X_2(\psi_{ij})}\Big)\bigg|,\\
\Err_{2,1}^{(2)}&=&\bigg|\sum_{i,j}\int \Re\Big(\pr^{\a} (x^k) \pr_\a(\psi_{kj})\ov{x^i X_2(\psi_{ij})}\Big)\bigg|,\\
\Err_{2,1}^{(3)}&=&\bigg|\sum_{i,j}\int \Re\Big(\pr^{\a} (x^i ) (x^k\pr_\a\psi_{kj})\ov{X_2(\psi_{ij})}\Big)\bigg|.
\eea
\esub

As for the control of $\Err_{1,1}^{(1)}$, we may apply Corollary \ref{cor:integrationbypartsforitimesfirstordertimesfirstorder:general:microlocalversion} with $A_i^k=(M_{K,A})_{i}^{k\a}$, and this time with $d_1=b_{\tphi}\xiphi+b_{\tt}\xit$, which yields
\bea
\lab{eq:finalesti:Err210}
\nn\Err_{2,1}^{(1)}&\les& \frac{1}{\sqrt{\dhor}}\bigg(\sum_{i,j=1}^3\int_{\MM_{r_+(1+\dhor'), \Rmic}(\Iti)}|\psi_{ij}|^2\bigg)^{\frac{1}{4}}\left(\widetilde{\M}[\pmb\psi]+\sum_{i,j=1}^3\int_{\Mntrap(\Iti)}|\square_\g\psi_{ij}|^2\right)^{\frac{3}{4}}\\
&&+\ep\M[\pmb\psi](\Iti).
\eea
Next, in order to control $\Err_{2,1}^{(2)}$, we first derive the following analog of \eqref{eq:precisecommutatorxiX1neededforcontrolofErr112}
\beaa
\nn [x^i, X_2] &=& \Opw\big(\{x^i, b_{\tphi}\xiphi+b_{\tt}\xit\}\big)+\Opw(\widetilde{S}^{-2,0}(\MM))\\
&=& \Opw\big(\{x^i, \xiphi\}b_{\tphi}+\{x^i, \xit\}b_{\tt}\big)+\Opw\big(\{x^i, b_{\tphi}\}\xiphi+\{x^i, b_{\tt}\}\xit\big)+\Opw(\widetilde{S}^{-2,0}(\MM)),
\eeaa
and hence, commuting $x^i$ with $X_2$ in the formula \eqref{def:Err21n:all} for $\Err_{2,1}^{(2)}$ and using Lemma \ref{lemma:gardinginequalitiesyieldcontroloftermsbywidetileM}, we infer
\beaa
\Err_{2,1}^{(2)}&\les&\bigg|\sum_{i,j}\int \Re\Big(\pr^{\a} (x^k) \pr_\a(\psi_{kj})\ov{X_2(x^i\psi_{ij})}\Big)\bigg|+\bigg(\sum_{i,j=1}^3\int_{\MM_{r_+(1+\dhor'), \Rmic}(\Iti)}|\psi_{ij}|^2\bigg)^{\frac{1}{2}}\big(\widetilde{\M}[\pmb\psi]\big)^{\frac{1}{2}}.
\eeaa
From there, we argue as for the control of $\Err_{1,1}^{(2)}$, i.e., we integrate by parts first in $\pr_{\a}$ and then in $X_2$, and we then take Cauchy-Schwarz and use again Lemma \ref{lemma:gardinginequalitiesyieldcontroloftermsbywidetileM} as well as \eqref{eq:widetildepsiij:tau2-4totau2} and  \eqref{def:Errdefectofpsi} to obtain the following analog of \eqref{eq:finalesti:Err112}
\bea
\lab{eq:finalesti:Err212}
\Err_{2,1}^{(2)} \les \bigg(\sum_{i,j=1}^3\int_{\MM_{r_+(1+\dhor'), \Rmic}(\Iti)}|\psi_{ij}|^2\bigg)^{\frac{1}{2}}\big(\widetilde{\M}[\pmb\psi]\big)^{\frac{1}{2}}+\big(\widetilde{\M}[\pmb\psi]\big)^{\frac{1}{2}}
\big(\Errdefect[\pmb\psi]\big)^{\frac{1}{2}}.
\eea
Next, for $\Err_{2,1}^{(3)}$, we commute $x^k$ and $\pr_\a$ in the formula \eqref{def:Err21n:all} for $\Err_{2,1}^{(3)}$ and obtain, using again Lemma \ref{lemma:gardinginequalitiesyieldcontroloftermsbywidetileM}, \eqref{eq:widetildepsiij:tau2-4totau2} and  \eqref{def:Errdefectofpsi},  
\bea\lab{eq:finalesti:Err213}
\nn\Err_{2,1}^{(3)}&\les&\bigg|\sum_{i,j}\int \Re\Big(\pr^{\a} (x^i ) \pr_\a(x^k\psi_{kj})\ov{X_2(\psi_{ij})}\Big)\bigg|\\
\nn&&+\bigg(\sum_{i,j=1}^3\int_{\MM_{r_+(1+\dhor'), \Rmic}(\Iti)}|\psi_{ij}|^2\bigg)^{\frac{1}{2}}\big(\widetilde{\M}[\pmb\psi]\big)^{\frac{1}{2}}\\
&\les& \bigg(\sum_{i,j=1}^3\int_{\MM_{r_+(1+\dhor'), \Rmic}(\Iti)}|\psi_{ij}|^2\bigg)^{\frac{1}{2}}\big(\widetilde{\M}[\pmb\psi]\big)^{\frac{1}{2}}+\big(\widetilde{\M}[\pmb\psi]\big)^{\frac{1}{2}}
\big(\Errdefect[\pmb\psi]\big)^{\frac{1}{2}}.
\eea
In conclusion, \eqref{eq:finalesti:Err210}, \eqref{eq:finalesti:Err212} and \eqref{eq:finalesti:Err213} yield
\bea
\lab{eq:controlofErr21}
\nn\bigg|\int \Err_{2,1}\bigg|
&\les& \frac{1}{\sqrt{\dhor}}\bigg(\sum_{i,j=1}^3\int_{\MM_{r_+(1+\dhor'), \Rmic}(\Iti)}|\psi_{ij}|^2\bigg)^{\frac{1}{4}}\left(\widetilde{\M}[\pmb\psi]+\sum_{i,j=1}^3\int_{\Mntrap(\Iti)}|\square_\g\psi_{ij}|^2\right)^{\frac{3}{4}}\\
&&+\ep\M[\pmb\psi](\Iti) +\big(\widetilde{\M}[\pmb\psi]\big)^{\frac{1}{2}}
\big(\Errdefect[\pmb\psi]\big)^{\frac{1}{2}}.
\eea
In view of the above estimates \eqref{eq:estiofErr24}, \eqref{eq:estiofErr23}, \eqref{eq:estiofErr22} and \eqref{eq:controlofErr21} and the decomposition \eqref{eq:definitionofErr2:decomposition}, and using \eqref{eq:controlofL2spacetimenormofsquaregpsiijbyL2spacetimenormofFijandMorawetzusingwaveeqsec8} to control $\square_\g\psi_{ij}$, we deduce
\bea
\lab{esti:Err2:MoraEstimiddleregion}
&&\bigg|\int \bigg(\Err_{2}-\sum_{i,j}\Re\Big(F_{ij}\ov{X_2\psi_{ij}}\Big)\bigg)\bigg|\nn\\
&\les_{\Rmic}& \ep\M[\pmb\psi](\Iti) +\A[\pmb \psi](\Iti)+
\big(\widetilde{\M}[\pmb\psi]\big)^{\frac{1}{2}}
\big(\Errdefect[\pmb\psi]\big)^{\frac{1}{2}}\nn\\
&&+ \frac{1}{\sqrt{\dhor}}\big(\A[\pmb \psi](\Iti)\big)^{\frac{1}{4}}\bigg(\widetilde{\M}[\pmb\psi]+\sum_{i,j=1}^3\int_{\Mntrap(\Iti)}|F_{ij}|^2\bigg)^{\frac{3}{4}}.
\eea


\subsubsection{Control of the error term $\Err_4$}
\lab{subsubsect:controlerrorterms:Morawetz:scalartoscalarized:Err4}


Next, we estimate the error term $\Err_4$ introduced in \eqref{def:errorterms:general:scalarizedwavefromtensorialwave} and associated to the PDO $E=\Opw(e_0)$, with the symbol $e_0$ given as in Section \ref{sec:relevantmixedsymbolsonMM}.
Using the wave equations \eqref{eq:ScalarizedWaveeq:general:Kerrpert}  for $\psi_{ij}$ and the self-adjointness of the PDO $E\in \Opw(\widetilde{S}^{0,0}(\MM))$ for the measure $d\Vref$, we infer
\bea
\lab{esti:Err4:MoraEstimiddleregion}
\bigg|\int\Err_4- \sum_{i,j}\int \Re\Big(F_{ij}\ov{E\psi_{ij}}\Big)\bigg|&=&\bigg|\sum_{i,j,k,l}\int \Re\Big(E\big(O(1)\pr^{\leq 1}\psi_{kl}\big)\ov{\psi_{ij}}\Big)\bigg|\nn\\
&\les& \bigg(\sum_{i,j}\int |\psi_{ij}|^2\bigg)^{\frac{1}{2}} \bigg(\sum_{k,l}\int \Big|E\big(O(1)\pr^{\leq 1}\psi_{kl}\big)\Big|^2\bigg)^{\frac{1}{2}}\nn\\
&\les&\bigg(\sum_{i,j}\int |\psi_{ij}|^2\bigg)^{\frac{1}{2}}\Big(\widetilde{\M}[\pmb \psi]\Big)^{\frac{1}{2}}
\eea
where we have used \eqref{eq:gardinginequalitiesyieldcontroloftermsbywidetileM:1} in the last line.


\subsubsection{Control of the error term $\Err_3$}
\lab{subsubsect:controlerrorterms:Err3:Morawetz:scalartoscalarized}


We now estimate the remaining error term $\Err_3$ introduced in \eqref{def:errorterms:general:scalarizedwavefromtensorialwave} and associated to the vectorfield $A\pr_{\tau}$.  To this end, for $\pmb\psi\in \mathfrak{s}_k(\mathbb{C})$, $k=0,1,2$, satisfying the variational tensorial wave equation \eqref{eq:Gen.RW-general}, i.e., 
\beaa
\squared_k \pmb\psi-V\pmb\psi=\pmb{N},
\eeaa
where $V$ is a real potential, we recall the energy-momentum tensor \eqref{eq:definition-QQ-mu-nu}, i.e., 
 \beaa
 \QQ_{\mu\nu}[\pmb\psi]=\Re\Big(\Db_\mu  \pmb\psi \c \ov{\Db _\nu  \pmb\psi }\Big)
          -\frac 12 \g_{\mu\nu}  \LL[\pmb\psi],
 \eeaa
 where the Lagrangian $ \LL[\pmb\psi]$ is given by
\beaa
 \LL[\pmb\psi]= \Re\left(\Db_\la  \pmb\psi\c\ov{\Db^\la  \pmb\psi }\right)+ V |\pmb\psi |^2.
 \eeaa
Also, recall from Proposition \ref{prop-app:stadard-comp-Psi} that the $1$-form   $\PP_\mu[\pmb\psi](X, w)$, for a vectorfield $X$, a real scalar function $w$ and a tensor $\pmb\psi\in \mathfrak{s}_k(\mathbb{C})$, $k=0,1,2$, is given by 
 \beaa
\PP_\mu[\pmb\psi](X, w)=\QQ_{\mu\nu}[\pmb\psi] X^\nu +\frac 1 2  w \Re\big(\pmb\psi \c \ov{\Db_\mu \pmb\psi }\big)-\frac 1 4|\pmb\psi|^2   \pr_\mu w
  \eeaa
and satisfies
  \bea\lab{eq:DivofPPmu:tensor:RW}
\D^\mu  \PP_\mu[\pmb\psi](X, w)&=& \frac 1 2 \QQ[\pmb\psi]  \c\piX - \frac 1 2 X( V) |\pmb\psi|^2 +\frac{k}{2}{}^{(X)}A_\nu\Im\Big(\pmb\psi\c\ov{\Ddot^{\nu}\pmb\psi}\Big)+\frac 12  w \LL[\pmb\psi]\nn\\
  && -\frac 1 4|\pmb\psi|^2   \square_\g  w   +  \Re\bigg(\ov{\bigg(\nab_X\pmb\psi +\frac 1 2   w \pmb\psi\bigg)}\c \left(\squared_k \pmb\psi- V\pmb\psi\right)\bigg),
\eea
with ${}^{(X)}A_\nu$ the 1-form introduced in \eqref{eq:thespacetime1formXA}.

The following basic estimate for the $1$-form $\PP_\mu[\pmb\psi](X, w)$ will be useful.

\begin{lemma}
\lab{lem:estimatesforcurrentforprttandprrBL}
Let $\pmb\psi\in\sk_2(\mathbb{C})$  and let $\psi_{ij}$ be the corresponding scalars given by $\psi_{ij}:=\pmb\psi(\Om_i, \Om_j)$ for $i,j=1,2,3$ where $\Om_i$, $i=1,2,3$ is the regular triplet introduced in Section \ref{sec:regulartripletinperturbationsofKerr}. For a real-valued vectorfield $X=O(1)\pr_{\tt}+O(1)\pr_{r} + O(r^{-2})\pr_{x^a}$ and  a real-valued scalar function $w=O(r^{-1})$, we have
\bsub
\lab{eq:difference:PPmupmbpsiandpsiij}
\bea
\lab{eq:difference:PPmupmbpsiandpsiij:energy}
&&\int_{\Sigma(\tau)}\bigg|\bigg({\PP}_{\mu}[\pmb\psi](X, w)- \sum_{i,j}\PP_\mu[\psi_{ij}](X, w)\bigg)N_{\Sigma_{\tau}}^{\mu}\bigg|\nn\\
 &\les&\left(\int_{\Sigma(\tau)}r^{-2}|\pmb\psi|^2\right)^{\frac{1}{2}}\Big(\E[\pmb\psi](\tau)\Big)^{\frac{1}{2}},\quad \forall \tau\geq \tmic, \\
\lab{eq:difference:PPmupmbpsiandpsiij:scri}
&&\int_{\II_+(\tau',\tau'')}\bigg|\bigg({\PP}_{\mu}[\pmb\psi](X, w)- \sum_{i,j}\PP_\mu[\psi_{ij}](X, w)\bigg)N_{\II_+}^{\mu}\bigg| \nn\\
&\les& \left(\int_{\II_+(\tau', \tau'')}r^{-2}|\pmb\psi|^2\right)^{\frac{1}{2}}\Big(\F_{\II_+}[\pmb \psi](\tau',\tau'')\Big)^{\frac{1}{2}},\quad \forall \tmic\leq \tau'<\tau'',\\
\lab{eq:difference:PPmupmbpsiandpsiij:Hr}
&&\int_{H_{r}(\tau',\tau'')}\bigg|\bigg({\PP}_{\mu}[\pmb\psi](X, w)- \sum_{i,j}\PP_\mu[\psi_{ij}](X, w)\bigg)N_{H_r}^{\mu}\bigg| \nn\\
&\les& \int_{H_{r}(\tau',\tau'')}r^{-1}|\pmb\psi| |\pr^{\leq 1}\pmb\psi|, \quad \forall r_+(1-\dhor)\leq r <+\infty, \,\, \tmic\leq \tau'<\tau''.
\eea
\esub
\end{lemma}

\begin{proof}
To begin with, we first show the following identity for a real-valued vectorfield $X=O(1)\pr_{\tt}+O(1)\pr_{r} + O(r^{-2})\pr_{x^a}$ and a real-valued scalar function $w$
\bea
\lab{esti:currentforprttandprrBL}
&&\PP_\mu[\pmb\psi](X, w)- \sum_{i,j}\PP_\mu[\psi_{ij}](X, w)\nn\\
&=&\Re\Big(M_{i\mu}^j \psi\ov{X\psi}  + (O(r^{-3})+\Ga_b) \psi\ov{\pr_{\mu}\psi} +M_{i\mu}^j (O(r^{-3})+\Ga_b) \psi\ov{\psi} \Big)+w\Re\big(M_{i\mu}^j \psi \ov{\psi}\big)\nn\\
&&+\g_{\mu \a}X^{\a} \Re\Big(O(r^{-2})\psi \ov{(e_3)^{\leq 1}\psi }+O(r^{-1})\psi\ov{e_a\psi}+(O(r^{-2})+\Ga_b) \psi \ov{e_4 \psi}\Big),
\eea
where the schematic notations $\Re(\psi\ov{e_{\a}\psi})$ and $\Re(\psi\ov{\psi})$ denote respectively any term of the form $\Re(\psi_{ij}\ov{e_\a(\psi_{kl})})$ and $\Re(\psi_{ij}\ov{\psi_{kl}})$.
To this end, recall from Lemma \ref{lem:estimatesforMialphaj:Kerrpert} the following estimates for $M_{i\a}^j$:
 \begin{equation}
 \begin{split}
& M_{i4}^j=O(r^{-2}), \quad
 M_{ia}^j = O(r^{-1}), \quad  M_{i3}^j=O(r^{-2})+\Ga_b,\\
 &M_{i\a}^j(\pr_{\tau})^{\a}=O(r^{-3}) + \Ga_b, \quad M_{i\a}^j( \pr_{r}^{\text{BL}})^{\a}=\Ga_b,
 \end{split}
 \end{equation}
 and notice that
\bea\lab{eq:PPmuXw:X=prttprrBL:general}
 && \PP_\mu[\pmb\psi](X, w)- \sum_{i,j}\PP_\mu[\psi_{ij}](X, w)\nn\\
 &=&
 \QQ_{\mu\nu}[\pmb\psi] X^\nu-\sum_{i,j}\QQ_{\mu\nu}[\psi_{ij}] X^\nu+\frac 1 2  w \bigg(\Re\big(\pmb\psi \c \ov{\Db_\mu \pmb\psi }\big)-\sum_{i,j}\Re\big(\psi_{ij} \ov{\pr_\mu \psi_{ij} }\big)\bigg).
 \eea
Using Lemma \ref{lem:ucdotv:product}, we have
\beaa
&&\Re\Big(\Db_\mu  \pmb\psi \c \ov{\Db _\nu  \pmb\psi }\Big)(\pr_{\tt})^{\nu} - \sum_{i,j}\Re\Big(\pr_\mu\psi_{ij}   \ov{\pr _\nu  \psi_{ij} }\Big)(\pr_{\tt})^{\nu} \nn\\
&=&\Re\Big(M_{i\mu}^j \psi\ov{\pr_{\tt}\psi}  + (O(r^{-3})+\Ga_b) \psi\ov{\pr_{\mu}\psi} +M_{i\mu}^j (O(r^{-3})+\Ga_b) \psi\ov{\psi} \Big),\\
&&\Re\left(\Db_\la  \pmb\psi\c\ov{\Db^\la  \pmb\psi }\right)-\sum_{i,j}\Re\left(\pr_\la  \psi_{ij} \ov{\pr^\la  \psi_{ij}  }\right)\nn\\
&=&\Re\Big(O(r^{-2})\psi \ov{(e_3)^{\leq 1}\psi }+O(r^{-1})\psi\ov{e_a\psi}+ (O(r^{-2})+\Ga_b) \psi \ov{e_4 \psi}\Big)
\eeaa
and
\beaa
w\bigg(\Re\big(\pmb\psi \c \ov{\Db_\mu \pmb\psi }\big)-\sum_{i,j}\Re\big(\psi_{ij} \ov{\pr_\mu \psi_{ij} }\big)\bigg)
=w\Re\big(M_{i\mu}^j \psi \ov{\psi}\big),
\eeaa
and plugging these three estimates into \eqref{eq:PPmuXw:X=prttprrBL:general} implies the desired estimate \eqref{esti:currentforprttandprrBL} for $X=O(1)\pr_{\tt}$. The estimate \eqref{esti:currentforprttandprrBL} for $X=O(1)\pr_{r}$ and $X=O(r^{-2})\pr_{x^a}$ follows in the exact same manner.

Next, we rely on \eqref{esti:currentforprttandprrBL} to show the desired estimates \eqref{eq:difference:PPmupmbpsiandpsiij}. Using \eqref{eq:consequenceasymptoticKerrandassumptionsinverselinearizedmetric} as well as \eqref{expression:prtauIIplus:nullinf} for $N_{\II_+}$, we have 
\bsub\lab{eq:expressionforthenormalstoSigmatauIIplusandHr}
\bea
N_{\Sigma(\tau)} &=& -\g^{\a\tau}\pr_\a=O(1)\pr_r+O(r^{-2})\pr_\tau+O(r^{-2})\pr_{x^a},\\
N_{\II_+} &=& \pr_\tau^{\II_+}=\pr_\tau +O(1)\pr_r+O(r^{-1})\pr_{x^a},\\
N_{H_r} &=& \g^{\a r}\pr_\a= O(1)\pr_r+O(1)\pr_\tau+O(r^{-1})\pr_{x^a}, 
\eea
\esub
and, in view of \eqref{estimates:Mialphaj:Kerrperturbations} and \eqref{eq:relationsbetweennullframeandcoordinatesframe2:moreprecise}, we have
\bea\lab{estimates:Mialphaj:Kerrperturbations:slightvariationfocusingoncoordinatesvectorfield}
M_{i\a}^j(\pr_r)^{\a}=O(r^{-2}), \qquad M_{i\a}^j(\pr_{\tau})^{\a}=O(r^{-3}) + \Ga_b, \qquad M_{i\a}^j(\pr_{x^a})^{\a}=O(1).
\eea
Moreover, in view of \cite[Lemma 2.22]{MaSz24}, for any $\tau_1<\tau_2$ and any $\de>0$, we have for any scalar function $\psi$
\bea\lab{eq:lem:nullnfFluxBdedByEnergy:Lemma222scalarwavepaper}
\liminf_{\tauu\to +\infty}\int_{\tau_1}^{\tau_2}\int_{\mathbb{S}^2}{(1+\tau-\tau_1)^{-1-\de}}r^{-1}|\mathfrak{d}^{\leq 1}\psi|^2{(\tauu, \tau, \omega)}r^2d\mathring{\ga}d\tau   \les \sup_{\tau\in[\tau_1, \tau_2]}\E[\psi](\tau).
\eea
Hence, \eqref{esti:currentforprttandprrBL}, \eqref{eq:expressionforthenormalstoSigmatauIIplusandHr}, \eqref{estimates:Mialphaj:Kerrperturbations:slightvariationfocusingoncoordinatesvectorfield} and \eqref{eq:lem:nullnfFluxBdedByEnergy:Lemma222scalarwavepaper}, together with the estimate
\eqref{estimates:Mialphaj:Kerrperturbations} for ${M_{i\a}^j}$ and Cauchy-Schwarz, implies, for a real-valued vectorfield $X=O(1)\pr_{\tt}+O(1)\pr_{r} + O(r^{-2})\pr_{x^a}$ and  a real-valued scalar function $w=O(r^{-1})$, and for any $\tmic\leq \tau'<\tau''$ and any $r_+(1-\dhor)\leq r<+\infty$,
\begin{align*}
&\int_{\Sigma(\tau')}\bigg|\bigg({\PP}_{\mu}[\pmb\psi](X, w)- \sum_{i,j}\PP_\mu[\psi_{ij}](X, w)\bigg)N_{\Sigma_{\tau}}^{\mu}\bigg|\\
 \les& \int_{\Sigma(\tau')}r^{-2}|\pmb\psi||\dk^{\leq 1}\pmb\psi| \les \left(\int_{\Sigma(\tau')}r^{-2}|\pmb\psi|^2\right)^{\frac{1}{2}}\big(\E[\pmb\psi](\tau')\big)^{\frac{1}{2}},\\
&\int_{\II_+(\tau',\tau'')}\bigg|\bigg({\PP}_{\mu}[\pmb\psi](X, w)- \sum_{i,j}\PP_\mu[\psi_{ij}](X, w)\bigg)N_{\II_+}^{\mu}\bigg|\\
\les& \int_{\tau'}^{\tau''}\int_{\mathbb{S}^2}r^{-1}|\pmb\psi|\left(|\nab_{\pr_\tau^{\II_+}}\pmb\psi|+r^{-1}|\nabla_{\pr_{x^a}^{\II_+}}^{\leq 1}\pmb\psi|+|\Ga_b||\dk^{\leq 1}\pmb\psi|\right)(\tauu=+\infty, \tau, \om)r^2d\mathring{\ga}d\tau\\
\les&  \left(\int_{\II_+(\tau', \tau'')}r^{-2}|\pmb\psi|^2\right)^{\frac{1}{2}}\Big(\F_{\II_+}[\pmb \psi](\tau',\tau'')\Big)^{\frac{1}{2}},\\
&\int_{H_{r}(\tau',\tau'')} \bigg|\bigg({\PP}_{\mu}[\pmb\psi](X, w)- \sum_{i,j}\PP_\mu[\psi_{ij}](X, w)\bigg)N_{H_r}^{\mu}\bigg| \les \int_{H_{r}(\tau',\tau'')}r^{-1}|\pmb\psi| |\pr^{\leq 1}\pmb\psi|, 
\end{align*}
as desired. This concludes the proof of Lemma \ref{lem:estimatesforcurrentforprttandprrBL}.
\end{proof}

We now state a key lemma to control the integral of $\err_3$ on $\MM_{r_+(1+\dhor'), \Rmic}(\Iti)$. For convenience, we prove a statement that holds on more general domains $\MM_{r_1,r_2}(\tau',\tau'')$.

\begin{lemma}\lab{lemma:howtodealwithmultiplierAxitforcoupledsystemwavespsiij}
Under the same assumptions for the scalars $\psi_{ij}$, $F_{ij}$ and the spacetime $(\MM,\g)$ as in Theorem \ref{th:mainenergymorawetzmicrolocal}, let $\chi_n(\tau)$, $n=1,2,3,4$, be smooth nonnegative cut-off functions satisfying
\begin{equation}
\lab{def:cutoffsintime1234}
\begin{split}
&\sum_{n=1}^4\chi_n(\tau)=1\,\,\forall\tau\in\Reals, \qquad\mathrm{supp}(\chi_1)\subset(-\infty, \tau_1+2), \qquad \mathrm{supp}(\chi_2)\subset(\tau_1+1, \tau_2-2),\\ 
&\mathrm{supp}(\chi_3)\subset(\tau_2-3, \tau_2+1), \qquad \mathrm{supp}(\chi_4)\subset(\tau_2, +\infty),
\end{split}
\end{equation}
and let the 1-form $\BB_\mu[\pmb\psi]$ be given by
\bea
\lab{eq:defofBBmupsi}
\BB_\mu[\pmb\psi] := \big(\chi_1(\tau)+\chi_3(\tau)\big)\sum_{i,j}{\PP}_\mu[\psi_{ij}](\pr_{\tt},0) +\chi_2(\tau){}^{(\pr_{\tt})}\widetilde{\PP}_\mu[\pmb\psi]+\chi_4(\tau)\big({}^{(\pr_{\tt})}\widetilde{\PP}_\mu[\pmb\psi]\big)_{K}
\eea
with ${}^{(\pr_{\tt})}\widetilde{\PP}_\mu[\pmb\psi]$ defined as in Lemma \ref{cor:modifiedcurrentsforprtandprvphi} for $\pmb\psi\in\sk_2(\mathbb{C})$ and $V=D_0|q|^{-2}$, $\big({}^{(\pr_{\tt})}\widetilde{\PP}_\mu[\pmb\psi]\big)_{K}$ being the corresponding quantity in Kerr for $\pmb\psi\in\sk_2(\mathbb{C})$ and $V=\frac{4}{|q|^2}- \frac{4a^2\cos^2\th(|q|^2+6mr)}{|q|^6}$, and ${\PP}_\mu[\psi_{ij}](\pr_{\tt},0)$ given in 
\eqref{definitionofcurrentPPmuXw:generaltensor} for a scalar $\psi_{ij}\in\sk_0(\mathbb{C})$ and $V=D_0|q|^{-2}$. Then, $\BB_\mu[\pmb\psi]$ satisfies 
\bea
\lab{eq:1formBB_mu:equality}
\BB_\mu[\pmb\psi]=\sum_{i,j}\QQ_{\mu\nu}[\psi_{ij}]  (\pr_{\tt})^\nu
+H_{\mu}[\psi,\pr\psi]
\eea
and 
\bea
\lab{eq:1formBB_mu:divergence}
 \D^\mu\BB_\mu[\pmb\psi]
= \sum_{i,j}\Re\left(F_{ij} \ov{T_{\tau_2}(\psi)_{ij}}\right)+G[\psi,\pr\psi],
\eea
 where
\bea
 \lab{def:widehatTpsiij:equalsnabTpsiij}
T_{\tau_2}(\psi)_{ij} &:=& \pr_{\tt}(\psi_{ij})  -\chi_2\big(M_{i\tt}^k\psi_{kj} + M_{j\tt}^k\psi_{ik} + 2i\widetilde{w}\psi_{ij}\big)\nn\\
&& - \chi_4\big((M_K)_{i\tt}^k\psi_{kj} + (M_K)_{j\tt}^k\psi_{ik} + 2i\widetilde{w}\psi_{ij}\big),
\eea
where $H_{\mu}[\psi,\pr\psi]$ denotes any term that satisfies for any $\tmic\leq\tau'<\tau''$, 
\bsub
\lab{eq:1formBB_mu:equality:errortermestimates}
\begin{align}
&\int_{\Sigma(\tau')}\big|H_{\mu}[\psi,\pr\psi]N_{\Sigma(\tau)}^{\mu}\big| \les \left(\int_{\Sigma(\tau')}r^{-2}|\pmb\psi|^2\right)^{\frac{1}{2}}\big(\E[\pmb\psi](\tau')\big)^{\frac{1}{2}},  \\
&\int_{\II_+(\tau',\tau'')}\big|H_{\mu}[\psi,\pr\psi]N_{\II_+}^{\mu} \big|\les  \left(\int_{\II_+(\tau', \tau'')}r^{-2}|\pmb\psi|^2\right)^{\frac{1}{2}}\Big(\F_{\II_+}[\pmb \psi](\tau',\tau'')\Big)^{\frac{1}{2}},\\
&\int_{H_{r}(\tau',\tau'')}\big|H_{\mu}[\psi,\pr\psi]N_{H_r}^{\mu}\big| \les \int_{H_{r}(\tau',\tau'')}r^{-1}|\pmb\psi| |\pr^{\leq 1}\pmb\psi|, \quad\forall r_+(1-\dhor)\leq r<+\infty,
\end{align}
\esub
and where $G[\psi,\pr\psi]$ denotes any term that satisfies for any $r_+(1+\dhor)\leq r_1<r_2<+\infty$ and any $\tmic\leq\tau'<\tau''$,
\bsub
\lab{def:Gpsiprpsi:lowerordertermofDivofBmu}
\bea
\int_{\MM_{r_1,r_2}(\tau',\tau'')} |G[\psi,\pr\psi] |&\les& \NNtlocal[\pmb \psi](\tau',\tau'')+\sum_{r=r_1, r_2}\int_{H_r(\tau',\tau'')}  r^{-2} |\psi| |\pr_{\tt}\psi|+\ep\EM[\pmb \psi](\tau',\tau'')\nn\\
&&+\sup_{\tau\in[\tmic,\tau_1+2]}\E[\pmb \psi](\tau), 
\eea
and for any $r_+(1+\dhor)\leq r_1<+\infty$ and any $\tmic\leq\tau'<\tau''$,  
\bea
\int_{\MM_{r_1,+\infty}(\tau',\tau'')} |G[\psi,\pr\psi] |&\les& \NNtlocal[\pmb \psi](\tau',\tau'')+\int_{H_{r_1}(\tau',\tau'')}  r^{-2} |\psi| |\pr_{\tt}\psi|+\ep\EM[\pmb \psi](\tau',\tau'')\nn\\
&&+\sup_{\tau\in[\tmic,\tau_1+2]}\E[\pmb \psi](\tau),
\eea
\esub
with
\bea\lab{def:NNtlocalinr:NNtEner:NNtaux:wavesystem:EMF:proof}
 \NNtlocal[\pmb \psi](\tau',\tau'') &:=&\bigg(\sup_{\tau\in (\tau',\tau'')}\int_{\Sigma(\tau)}r^{-2}|\pmb \psi|^2
 +\int_{\MM(\tau',\tau'')}r^{-3}|\pmb\psi|^2\nn\\
&&+\int_{\II_+(\tau',\tau'')}r^{-2}|\pmb \psi|^2+\Errdefect[\pmb\psi]\bigg)^{\frac{1}{2}}\times\Big(\widetilde{\EMF}[\pmb \psi](\tau',\tau'')\Big)^{\frac{1}{2}}.
 \eea
\end{lemma}

\begin{proof}
In view of the definition of ${}^{(\pr_{\tt})}\widetilde{\PP}_\mu[\pmb\psi]$ in Lemma \ref{cor:modifiedcurrentsforprtandprvphi}, and using the fact that $\widetilde{w}=O(r^{-3})$, we have
\bea
\lab{eq:esti:widetildePmupmbpsi}
&&{}^{(\pr_{\tt})}\widetilde{\PP}_\mu[\pmb\psi]
-\sum_{i,j}{\PP}_\mu[\psi_{ij}](\pr_{\tt},0)\nn\\
&=&{\PP}_\mu[\pmb\psi](\pr_{\tt}, 0) -\sum_{i,j}{\PP}_\mu[\psi_{ij}](\pr_{\tt},0)+ 2\widetilde{w}\Im\left( \pmb\psi\c\ov{\Ddot_\mu\pmb\psi}\right)+(\pr_{\tt})_{\mu}\frac{2a\cos\th}{\qs}2\widetilde{w}|\pmb\psi|^2 \nn\\
&=& {\PP}_\mu[\pmb\psi](\pr_{\tt}, 0) -\sum_{i,j}{\PP}_\mu[\psi_{ij}](\pr_{\tt},0)+ O(r^{-3})\Im\left( \pmb\psi\c\ov{\Ddot_\mu\pmb\psi}\right)+(\pr_{\tt})_{\mu}O(r^{-5})|\pmb\psi|^2.
\eea
Plugging this into the expression \eqref{eq:defofBBmupsi} of $\BB_{\mu}[\pmb\psi]$, and using the fact that
\beaa
\left(\sum_{n=1}^4\chi_n(\tau)\right)\sum_{i,j}{\PP}_\mu[\psi_{ij}](\pr_{\tt},0)=\sum_{i,j}{\PP}_\mu[\psi_{ij}](\pr_{\tt},0)=\sum_{i,j}\QQ_{\mu\nu}[\psi_{ij}]  (\pr_{\tt})^\nu,
\eeaa
yields the identity \eqref{eq:1formBB_mu:equality} with $H_{\mu}[\psi,\pr\psi]$ given by
\beaa
&& H_{\mu}[\psi,\pr\psi]\\ 
&=& \chi_2(\tau)\bigg({\PP}_\mu[\pmb\psi](\pr_{\tt}, 0) -\sum_{i,j}{\PP}_\mu[\psi_{ij}](\pr_{\tt},0)+ O(r^{-3})\Im\left( \pmb\psi\c\ov{\Ddot_\mu\pmb\psi}\right)+(\pr_{\tt})_{\mu}O(r^{-5})|\pmb\psi|^2\bigg)\\
&&+\chi_4(\tau)\bigg({\PP}_\mu[\pmb\psi](\pr_{\tt}, 0) -\sum_{i,j}{\PP}_\mu[\psi_{ij}](\pr_{\tt},0)+ O(r^{-3})\Im\left( \pmb\psi\c\ov{\Ddot_\mu\pmb\psi}\right)+(\pr_{\tt})_{\mu}O(r^{-5})|\pmb\psi|^2\bigg)_K
\eeaa
which together with \eqref{eq:difference:PPmupmbpsiandpsiij} in the particular case $(X, w)=(\pr_{\tau}, 0)$ 
immediately implies \eqref{eq:1formBB_mu:equality:errortermestimates}.

Next, we compute the divergence of $\BB_\mu[\pmb\psi]$. We have, from the formula \eqref{eq:defofBBmupsi} of $\BB_\mu[\pmb\psi]$,
\bea
\lab{eq:computeDivergenceofBB:intermsofBulk:00}
 \D^\mu\BB_\mu[\pmb\psi]
 &=&\sum_{n=1}^4\Bulk{n}+\bigg(\big(\chi_1'(\tau)+\chi_3'(\tau)\big)\sum_{i,j}{\PP}_\mu[\psi_{ij}](\pr_{\tt},0)  \nn\\
 &&\qquad\qquad\qquad\qquad\quad+\chi_2'(\tau){}^{(\pr_{\tt})}\widetilde{\PP}_\mu[\pmb\psi]+\chi_4'(\tau)\big({}^{(\pr_{\tt})}\widetilde{\PP}_\mu[\pmb\psi]\big)_{K} \bigg)\D^\mu(\tau)\nn\\
 &=& \sum_{n=1}^4\Bulk{n}+\chi_2'(\tau)\D^\mu(\tau)\bigg({}^{(\pr_{\tt})}\widetilde{\PP}_\mu[\pmb\psi] - \sum_{i,j}{\PP}_\mu[\psi_{ij}](\pr_\tau, 0)\bigg)\nn\\
&& +\chi_4'(\tau)\D^\mu(\tau)\bigg(\big({}^{(\pr_{\tt})}\widetilde{\PP}_\mu[\pmb\psi]\big)_{K} - \sum_{i,j}{\PP}_\mu[\psi_{ij}](\pr_\tau, 0)\bigg)\nn\\
&=&\sum_{n=1}^4\Bulk{n} +\chi_2'(\tau)\D^\mu(\tau)H_{\mu}[\psi,\pr\psi]+\chi_4'(\tau)\D^\mu(\tau)H_{\mu}[\psi,\pr\psi]\nn\\
 &=& \sum_{n=1}^4\Bulk{n}+G[\psi,\pr\psi],
\eea
where we have defined in the first equality of \eqref{eq:computeDivergenceofBB:intermsofBulk:00}
\bsub
\begin{align}
\Bulk{1} :=& \chi_1(\tau)\sum_{i,j}\D^\mu{\PP}_\mu[\psi_{ij}](\pr_{\tt},0) , &\Bulk{2} :=&\chi_2(\tau)\D^\mu{}^{(\pr_{\tt})}\widetilde{\PP}_\mu[\pmb\psi],  \\
\Bulk{3} :=& \chi_3(\tau)\sum_{i,j}\D^\mu{\PP}_\mu[\psi_{ij}](\pr_{\tt},0),&
\Bulk{4} :=&\chi_4(\tau)\D_K^\mu\big({}^{(\pr_{\tt})}\widetilde{\PP}_\mu[\pmb\psi]\big)_{K},
\end{align}
\esub
used in the second equality of \eqref{eq:computeDivergenceofBB:intermsofBulk:00} the fact $\sum_{n=1}^4\chi_n'(\tau)=0$ which follows from $\sum_{n=1}^4\chi_n(\tau)=1$, used in the third equality of \eqref{eq:computeDivergenceofBB:intermsofBulk:00} the following consequence of Lemma \ref{lem:estimatesforcurrentforprttandprrBL}
\beaa
{}^{(\pr_{\tt})}\widetilde{\PP}_\mu[\pmb\psi]- \sum_{i,j}{\PP}_\mu[\psi_{ij}](\pr_{\tt},0)&=&H_{\mu}[\psi,\pr\psi], \\\big({}^{(\pr_{\tt})}\widetilde{\PP}_\mu[\pmb\psi]\big)_{K}- \sum_{i,j}{\PP}_\mu[\psi_{ij}](\pr_{\tt},0)&=&H_{\mu}[\psi,\pr\psi],
\eeaa
and used the first estimate in \eqref{eq:1formBB_mu:equality:errortermestimates} and the support properties of $\chi_2'(\tau)$ and $\chi_4'(\tau)$ in the last equality of \eqref{eq:computeDivergenceofBB:intermsofBulk:00}.

Next, we estimate the terms $\Bulk{n}$, $n=1,2,3,4$. We start with the term $\Bulk{4}$, which is supported on $\MM(\tau_2,+\infty)$. In view of the form of the wave equations \eqref{eq:ScalarizedWaveeq:general:Kerrpert}, and since $\psi_{ij}=\pmb\psi((\Om_K)_i, (\Om_K)_j)$ for a tensor $\pmb\psi\in\sk_{2,K}(\mathbb{C})$  in $\MM(\tau_2,+\infty)$,
we apply Corollary \ref{cor:energyestimatetensoriallevel:scalarizedversion} with $V=D_0|q|^{-2}+f_{D_0}=\frac{4}{|q|^2}- \frac{4a^2\cos^2\th(|q|^2+6mr)}{|q|^6}$ to deduce
\beaa
\D_K^\mu\big({}^{(\pr_{\tt})}\widetilde{\PP}_\mu[\pmb\psi]\big)_{K}=\Re\left(F_{ij} \ov{\big[\pr_{\tau}\psi_{ij}-\big((M_K)_{i\tt}^k\psi_{kj} + (M_K)_{j\tt}^k\psi_{ik} + 2i\widetilde{w}\psi_{ij}\big)\big]}\right), \quad \textrm{on}\quad \MM(\tau_2,+\infty).
\eeaa
and hence, by multiplying  on both  sides by $\chi_4$,
\bea
\lab{eq:1formBB_mu:divergence:geqtau2}
\Bulk{4}= \Re\left(F_{ij} \ov{\big[\chi_4\pr_{\tau}\psi_{ij}-\chi_4\big((M_K)_{i\tt}^k\psi_{kj} + (M_K)_{j\tt}^k\psi_{ik} + 2i\widetilde{w}\psi_{ij}\big)\big]}\right).
\eea

Similarly, in view of the form of the wave equations \eqref{eq:ScalarizedWaveeq:general:Kerrpert}, and since $\psi_{ij}=\pmb\psi(\Om_i, \Om_j)$ for a tensor $\pmb\psi\in\sk_2(\mathbb{C})$  in $\MM(\tau_1+1,\tau_2-2)$, we apply Corollary \ref{cor:energyestimatetensoriallevel:scalarizedversion} with $V=D_0|q|^{-2}$ and multiply  on both  sides by $\chi_2$ to deduce
\beaa
\Bulk{2}&=& \frac{1}{2}\chi_2\QQ[\pmb\psi]  \c {^{(\pr_{\tt})}}\pi+\chi_2\Div (\pr_{\tt})\frac{4a\cos\th}{\qs}\widetilde{w}|\pmb\psi|^2+\chi_2 r^{-1}\Ga_b\Im\Big(\pmb\psi\c\ov{\Ddot\pmb\psi}\Big)\nn\\
&&+\sum_{i,j} \Re\left(F_{ij} \ov{\big[\chi_2\pr_{\tt}(\psi_{ij})  -\chi_2\big(M_{i\tt}^k\psi_{kj} + M_{j\tt}^k\psi_{ik} + 2i\widetilde{w}\psi_{ij}\big)\big]}\right).
\eeaa
In view of the estimate
\bea
\lab{eq:estimatesofdeformationtensorofprtt:proof}
\int_{\MM_{r_1,r_2}(\tau',\tau'')}\big| \QQ[\psi_{ij}]\c{^{(\pr_{\tt})}}\pi\big|\les \ep\EM_{r_1,r_2}[\psi_{ij}](\tau',\tau'')
 \eea
which follows from  Lemma \ref{lemma:basiclemmaforcontrolNLterms:ter} and the estimates for the deformation tensor ${^{(\pr_{\tt})}}\pi$ in Lemma \ref{lemma:controlofdeformationtensorsforenergyMorawetz}, we infer
  \bea
\lab{eq:estimatesforQQtimesdeformationtensorofT}
&&\int_{\MM_{r_1,r_2}(\tau',\tau'')} \bigg(\left|  \QQ[\pmb\psi]  \c {^{(\pr_{\tt})}}\pi\right|+\left|\Div (\pr_{\tt})\frac{4a\cos\th}{\qs}\widetilde{w}|\pmb\psi|^2\right| +\left|r^{-1}\Ga_b\Im\Big(\pmb\psi\c\ov{\Ddot\pmb\psi}\Big)\right|\bigg) \nn\\
&\les& \ep\EM_{r_1,r_2}[\pmb \psi](\tau',\tau''),
 \eea
which then implies
\bea
\lab{eq:1formBB_mu:divergence:tau1+1totau2}
 \Bulk{2}=\sum_{i,j} \Re\left(F_{ij} \ov{\big[\chi_2\pr_{\tt}(\psi_{ij})  -\chi_2\big(M_{i\tt}^k\psi_{kj} + M_{j\tt}^k\psi_{ik} + 2i\widetilde{w}\psi_{ij}\big)\big]}\right)+G[\psi,\pr\psi].
\eea
 
Next, we consider the term $\Bulk{3}$. Recall that the scalars $\psi_{ij}$ satisfy the coupled system of wave equations \eqref{eq:ScalarizedWaveeq:general:Kerrpert}, which can be rewritten in the form of   \eqref{eq:ScalarizedWaveeq:general:Kerrpert:rewriteinsect8}, and then further decomposed as follows
\bea
\lab{eq:scalarizedextendedeqs:psiij:moraandener}
\square_{\g}(\psi_{ij})-D_0|q|^{-2}\psi_{ij}=F_{ij} +\sum_{n=1}^4G_{n,ij} ,
\eea
where 
\bsub
\lab{def:Gnij:energyesti:withpotential}
\bea
G_{1,ij}&=& S_K(\psi)_{ij} ,\\
 G_{2,ij}&=&\frac{4ia\cos\th}{|q|^2} \pr_{\tt}\psi_{ij},\\
  G_{3,ij}&=&\chi_{\tau_1, \tau_2}(\widehat{Q}\psi)_{ij}+(1-\chi_{\tau_1, \tau_2})\big((\widehat{Q}_K\psi)_{ij}+f_{D_0}\psi_{ij}\big) ,\\
    G_{4,ij}&=& \chi_{\tau_1, \tau_2}\big({S}(\psi)_{ij}-{S}_K(\psi)_{ij}\big).
\eea
\esub
Hence, by applying \eqref{eq:DivofPPmu:tensor:RW} to each scalar $\psi_{ij}\in\sk_0(\mathbb{C})$ with $(X,w)=(\pr_{\tt}, 0)$, $k=0$ and $V=D_0|q|^{-2}$, multiplying on both sides by $\chi_3$, and summing over $i,j=1,2,3$, we deduce
\bea
\lab{eq:divergenceenergyidentityforscalarwave:psiij:arbitraryr}
\Bulk{3}
=\sum_{i,j} \bigg(\Re\Big(F_{ij}\ov{\chi_3\pr_{\tt}\psi_{ij}}\Big)
+\frac{1}{2}\chi_3\QQ[\psi_{ij}]\c{^{(\pr_{\tt})}}\pi+\sum_{n=1}^4\Re\Big(\chi_3G_{n,ij}\ov{\pr_{\tt}\psi_{ij}}\Big)\bigg).
\eea
In view of the estimates \eqref{eq:assumptionsonregulartripletinperturbationsofKerr:0}  \eqref{estimates:Mialphaj:Kerrperturbations} for ${M_{i\a}^j}$, the form of ${S}(\psi)_{ij}$ and $({Q}\psi)_{ij}$ in \eqref{SandV}, and the form of $(\widehat{Q}\psi)_{ij}$ in \eqref{hatSandV:generalwave:Kerrpert}, we have
\begin{equation}
\lab{estimates:widehatSandQ:perturbationsofKerr}
S(\psi)_{ij}=O(r^{-2})\dk\psi +\Ga_g\dk\psi, \quad (\widehat{Q}\psi)_{ij}=O(r^{-2})\psi +\dk^{\leq 1}\Ga_g \psi, \quad  S(\psi)_{ij}-S_K(\psi)_{ij}=\Ga_g\dk\psi,
 \end{equation}
 and hence, noticing also $\Re(G_{2,ij}\ov{\pr_{\tt}\psi_{ij}})=0$, we infer
\beaa
\sum_{i,j}\sum_{n=2}^4\Big|\Re\Big(\chi_3 G_{n,ij}\ov{\pr_{\tt}\psi_{ij}}\Big)\Big|
 \les \ep \chi_3 r^{-2}|\dk^{\leq 1}\pmb\psi|^2 +\chi_3 r^{-2}|\pmb\psi| |\dk^{\leq 1}\pmb\psi|,
\eeaa
which together with \eqref{eq:estimatesofdeformationtensorofprtt:proof} and the support properties of $\chi_3$ yields 
\bea
\lab{eq:erroG234:prtau:arbitraryregion}
\sum_{i,j}\bigg(\frac{1}{2}\chi_3\QQ[\psi_{ij}]\c{^{(\pr_{\tt})}}\pi+\sum_{n=2}^4\Re\Big(\chi_3 G_{n,ij}\ov{\pr_{\tt}\psi_{ij}}\Big)\bigg)
=G[\psi,\pr\psi].
\eea

Next, we consider the term $\Re(\chi_3G_{1,ij}\ov{\pr_{\tt}\psi_{ij}})$. As $S_K(\psi)_{ij}=2(M_K)_{i}^{k\a}\pr_\a(\psi_{kj}) +2(M_K)_{j}^{k\a}\pr_\a(\psi_{ik})$, we have
\bea
\Re(\chi_3G_{1,ij}\ov{\pr_{\tt}\psi_{ij}})
=2\Re\big(\ov{\chi_3\pr_{\tt}\psi_{ij}}(M_K)_{i}^{k\a}\pr_\a(\psi_{kj}) \big)+2\Re\big(\ov{\chi_3\pr_{\tt}\psi_{ij}}(M_K)_{j}^{k\a}\pr_\a(\psi_{ik})\big),
\eea
and we integrate this identity over $\MM_{r_1, r_2}(\tau',\tau'')$.

Consider first the case that $r_2=+\infty$. We estimate only the integral of the part involving $(M_K)_{i}^{k\a}\pr_\a(\psi_{kj})$ over $\MM_{r_1, +\infty}(\tau',\tau'')$ with $r_+(1+\dhor)\leq r_1<+\infty$, the control of the integral of the other part involving $(M_K)_{j}^{k\a}\pr_\a(\psi_{ik})$ being estimated in the exact same way.  
Using the decomposition \eqref{eq:decompositionofMikalpha:symandantisym} of $(M_K)_{i}^{k\a}$
\beaa
(M_K)_{i}^{k\a}&=&(M_{K,S})_{i}^{k\a} +(M_{K,A})_{i}^{k\a}= -\frac{1}{2}\pr^{\a} (x^i x^k) +(M_{K,A})_{i}^{k\a}\nn\\
&=& -\frac{1}{2}\gam^{\a\b}\pr_{\b}(x^i x^k) +\gam^{\a\b}(M_{K,A})_{i\b}^{k},
\eeaa
and in view of the estimate for ${(M_K)_{i\a}^j}$ in Lemma \ref{lemma:computationoftheMialphajinKerr}, we obtain the following:
\begin{itemize}
\item By integrating the differential identity \eqref{eq:integrationbypartsforitimesfirstordertimesfirstorder:general:caseantisymmatrix} with $A_i^k=\chi_3\gam^{\a\b}(M_{K,A})_{i\b}^{k}$ and $x^\b=\tau$, we deduce that the integral involving $\gam^{\a\b}(M_{K,A})_{i\b}^{k}$ is bounded by\footnote{For the boundary term on $H_{r_1}$, we use $\g(N_{H_r}, e_a)=\g(\g^{\a r}\pr_\a, e_a)=e_a(r)=r\Ga_g$ which implies $\g(N_{H_r}, e_a)(M_{K,A})_{i}^{ka}=\Ga_g=O(r^{-2})$, and for the boundary term on $\II_+(\tau',\tau'')$, we use $\g(N_{\II_+}, \pr_{\tt})=O(1)$ and $\g(N_{\II_+}, e_a)=O(1)$ which imply $\g(N_{\II_+}, \pr_{\tt})\gam^{\a\b}(M_{K,A})_{i\b}^{k}\pr_{\a}=O(r^{-2})\dk$ and $\g(N_{\II_+}, e_{\a})\gam^{\a\b}(M_{K,A})_{i\b}^{k}=O(r^{-2})$.}  
\beaa
\sup_{\tau\in[\tau',\tau'']}\int_{\Sigma(\tau)} r^{-2} |\psi| |\dk^{\leq 1}\psi| +\int_{H_{r_1}(\tau',\tau'')} r^{-2} |\psi| |\pr_\tau\psi|\nn\\
+\bigg(\int_{\II_+(\tau',\tau'')}r^{-2} |\pmb\psi|^2\bigg)^{\frac{1}{2}} \Big(\F_{\II_+}[\pmb\psi](\tau',\tau'')\Big)^{\frac{1}{2}} &\les& \int_{\MM_{r_1,+\infty}(\tau',\tau'')}|G[\psi,\pr\psi]|
\eeaa
in view of the estimates \eqref{def:Gpsiprpsi:lowerordertermofDivofBmu}.

\item Since we have 
\begin{align*}
\gam^{\a\b}\pr_\b(x^ix^k)\pr_\a(\psi_{kj})\ov{\pr_\tau(\psi_{ij})} 
={}&\gam^{\a\b}\pr_\b(x^k)\pr_\a(\psi_{kj})\ov{\pr_\tau(x^i\psi_{ij})}+\gam^{\a\b}\pr_\b(x^i)\pr_\a(x^k\psi_{ik})\ov{\pr_\tau(\psi_{ij})}\\
& - \gam^{\a\b}\pr_\b(x^i)\pr_\a(x^k)\psi_{ik}\ov{\pr_\tau(\psi_{ij})},
\end{align*}
and in view of the fact that
\beaa
\gam^{\a\b}\pr_\b(x^k)\pr_\a=r^{-2}\dk, \quad \gam^{\a\b}\pr_\b(x^i)\pr_\a(x^k)=O(r^{-2}),
\eeaa
the integral involving $-\frac{1}{2}\gam^{\a\b}\pr_\b (x^i x^k)$ is, in view of \eqref{def:Errdefectofpsi} and the support properties of $\chi_3$, bounded by 
\beaa
\sup_{\tt\in[\tau',\tau'']}\big(\E[\pmb\psi](\tau)\big)^{\frac{1}{2}}\left(\Errdefect[\pmb\psi]+\int_{\Sigma(\tt)} r^{-2}|\pmb\psi|^2\right)^{\frac{1}{2}}.
\eeaa
\end{itemize}
Therefore, we deduce, in view of the estimates \eqref{def:Gpsiprpsi:lowerordertermofDivofBmu},
\bea
\lab{eq:erroG1:prtau:middleregion:after}
\bigg|\int_{\MM_{r_1, +\infty}(\tau',\tau'')}  \Re\big(\ov{\pr_{\tt}\psi_{ij}}\chi_3 G_{1,ij}\big)\bigg|
&\les& \int_{\MM_{r_1,+\infty}(\tau',\tau'')}|G[\psi,\pr\psi]|.
\eea
Plugging the two estimates \eqref{eq:erroG234:prtau:arbitraryregion} and \eqref{eq:erroG1:prtau:middleregion:after} into \eqref{eq:divergenceenergyidentityforscalarwave:psiij:arbitraryr}, we infer
\bea
\lab{eq:erroG1234:prtau:middleregion}
\Bulk{3}
=\sum_{i,j}\Re\Big(F_{ij}\ov{\chi_3\pr_{\tt}\psi_{ij}}\Big)
+G[\psi,\pr\psi].
\eea

Consider next the case that $r_2<+\infty$. This is in fact simpler than the above case in which $r_2=+\infty$, as the boundary term on $\II_+(\tau',\tau'')$ is now replaced by a boundary term on $H_{r_2}$ which can be bounded in the same manner as for the boundary term on $H_{r_1}$. By the same argument, we conclude the estimate \eqref{eq:erroG1234:prtau:middleregion} as well.

Next, we consider the term $\Bulk{1}$. We apply the energy identity \eqref{eq:DivofPPmu:tensor:RW:prop} to each equation \eqref{eq:scalarizedextendedeqs:psiij:moraandener} of $\psi_{ij}$ in $\MM(\tau',\tau'')$, multiply  on both  sides by $\chi_1$ and sum up over $i,j$ to deduce
\bea
\lab{eq:divergenceenergyidentityforscalarwave:psiij:arbitraryr:leqtau1+1}
\Bulk{1}
=\sum_{i,j} \bigg(\Re\Big(F_{ij}\ov{\chi_1\pr_{\tt}\psi_{ij}}\Big)
+\frac{1}{2}\chi_1\QQ[\psi_{ij}]\c{^{(\pr_{\tt})}}\pi+\sum_{n=1}^4\Re\Big(\chi_1G_{n,ij}\ov{\pr_{\tt}\psi_{ij}}\Big)\bigg).
\eea
 In view of the estimate
\eqref{estimates:widehatSandQ:perturbationsofKerr} and the fact that $\chi_1(\tau)$ is supported in $\tau\leq \tau_1+2$, we deduce
\bea
\lab{eq:erroG234:prtau:arbitraryregion:leqtau1+1}
\bigg|\int_{\MM_{r_1, r_2}(\tau',\tau'')} \sum_{n=1}^4  \Re\Big(\ov{\pr_{\tt}\psi_{ij}}\chi_1G_{n,ij}\Big)\bigg|
\les \sup_{\tt\in[\tmic,\tau_1+2]}\E[\pmb\psi](\tau),
\eea
and together with the estimate \eqref{eq:estimatesofdeformationtensorofprtt:proof}, we infer
\bea
\lab{eq:erroG1234:prtau:leqtau1+1}
\Bulk{1}
=\sum_{i,j} \Re\Big(F_{ij}\ov{\chi_1\pr_{\tt}\psi_{ij}}\Big)
+G[\psi,\pr\psi].
\eea

In the end, plugging the estimates \eqref{eq:1formBB_mu:divergence:geqtau2}, \eqref{eq:1formBB_mu:divergence:tau1+1totau2},  \eqref{eq:erroG1234:prtau:middleregion} and \eqref{eq:erroG1234:prtau:leqtau1+1} into \eqref{eq:computeDivergenceofBB:intermsofBulk:00}, and using $\sum_{n=1}^4\chi_n=1$, we then obtain
the identity \eqref{eq:1formBB_mu:divergence} with $T_{\tau_2}$ defined as in \eqref{def:widehatTpsiij:equalsnabTpsiij} and with $G[\psi,\pr\psi]$ satisfying \eqref{def:Gpsiprpsi:lowerordertermofDivofBmu}. This concludes the proof of Lemma \ref{lemma:howtodealwithmultiplierAxitforcoupledsystemwavespsiij}.
\end{proof}

We now estimate the error term $\Err_3$ defined in \eqref{def:errorterms:general:scalarizedwavefromtensorialwave}. Using the 1-form $\BB_\mu[\pmb\psi]$ introduced in Lemma \ref{lemma:howtodealwithmultiplierAxitforcoupledsystemwavespsiij}, we have 
\beaa
&&\bigg|\int\frac{1}{A}\Err_3-\int \D^\mu\BB_\mu[\pmb\psi]\bigg| \nn\\
&=&\bigg|\int\sum_{i,j}\Re\Big(\square_{\g}\psi_{ij} \ov{\pr_{\tt}\psi_{ij}}\Big)-\int \D^\mu\BB_\mu[\pmb\psi]\bigg|\nn\\
&\les&\bigg|\int \D^{\mu}\bigg( \sum_{i,j}\QQ_{\mu\nu}[\psi_{ij}]  (\pr_{\tt})^\nu\bigg) -\int \D^\mu\BB_\mu[\pmb\psi]\bigg|+\bigg| \int\sum_{i,j}\QQ[\psi_{ij}]{^{(\pr_{\tt})}}\pi\bigg|\\
&\les&\bigg|\int \D^{\mu}\bigg(\BB_\mu[\pmb\psi] - \sum_{i,j}\QQ_{\mu\nu}[\psi_{ij}](\pr_{\tt})^\nu\bigg)  \bigg| +\ep\EM[\pmb \psi](\Iti)
\eeaa
where we have used \eqref{eq:estimatesofdeformationtensorofprtt:proof}  in the last estimate. Hence, in view  of the identity \eqref{eq:1formBB_mu:equality} and the estimate \eqref{eq:1formBB_mu:equality:errortermestimates}, we infer
\beaa
\bigg|\int\frac{1}{A}\Err_3-\int \D^\mu\BB_\mu[\pmb\psi]\bigg| &\les&\bigg|\int \D^{\mu}H_{\mu}[\psi,\pr\psi]\bigg| +\ep\EM[\pmb \psi](\Iti) \nn\\
&\les&\sum_{i,j,k,l}\bigg(\int_{H_{r_+(1+\dhor')}(\Iti)}+\int_{H_{\Rmic}(\Iti)}\bigg) |\psi_{ij}| |\pr^{\leq 1}\psi_{kl}|+\ep\EM[\pmb \psi](\Iti)\nn\\
&\les& \good^{(1)}[\psi]+\ep\EM[\pmb \psi](\Iti).
\eeaa
This estimate, together with the identity \eqref{eq:1formBB_mu:divergence} and the estimate \eqref{def:Gpsiprpsi:lowerordertermofDivofBmu} with the choice $(r_1, r_2)=(r_+(1+\dhor'), \Rmic)$, yields
\bea
\lab{esti:Err3:MoraEstimiddleregion}
&&\bigg|\int\Err_3-\int\Re\left(F_{ij} \ov{AT_{\tau_2}(\psi)_{ij}}\right)\bigg|\nn\\
&\les & \good[\psi]+\ep\EM[\pmb \psi](\Iti)+\bigg|\int G[\psi,\pr\psi]\bigg|\nn\\
&\les & \good[\psi]+\ep\EM[\pmb \psi](\Iti)+ \NNtlocal[\pmb \psi](\Iti)+\sup_{\tau\in[\tmic,\tau_1+2]}\E[\pmb \psi](\tau)\nn\\
&\les &\sup_{\tau\in[\tmic,\tau_1+2]}\E[\pmb \psi](\tau)+ \ep\EM[\pmb \psi](\Iti)+ \NNtlocal[\pmb \psi](\Iti)\nn\\
&&+\frac{1}{\sqrt{\dhor}}\big(\A[\pmb \psi](\Iti)\big)^{\frac{1}{4}}\bigg(\widetilde{\M}[\pmb\psi]+\sum_{i,j=1}^3\int_{\Mntrap(\Iti)}|F_{ij}|^2\bigg)^{\frac{3}{4}},
\eea
where we have  controlled $\good[\psi]$ thanks to Lemma \ref{lem:lowerorderterms:controlled:scalarizedwavefromtensorial} and \eqref{eq:controlofL2spacetimenormofsquaregpsiijbyL2spacetimenormofFijandMorawetzusingwaveeqsec8} in the last line.


\subsubsection{Concluding the proof of a Morawetz 
estimate in $\MM_{r_+(1+\dhor'),{\Rmic}}$}


In view of the estimates \eqref{esti:Err1:MoraEstimiddleregion}, \eqref{esti:Err2:MoraEstimiddleregion} and \eqref{esti:Err4:MoraEstimiddleregion} for the error terms $\Err_1$, $\Err_2$ and $\Err_4$,  we conclude
\bea
\lab{esti:Err124:MoraEstimiddleregion}
&&\bigg|\int \left(\Err_1 + \Err_2 + \Err_4\right)- \sum_{i,j}\int\Re\left(F_{ij} \ov{(X_1+X_2+E)\psi_{ij}}\right)\bigg|\nn\\
&\les&\ep\EM[\pmb\psi](\Iti) +\A[\pmb \psi](\Iti)+\big(\widetilde{\M}[\pmb\psi]\big)^{\frac{1}{2}}
\big(\Errdefect[\pmb\psi]\big)^{\frac{1}{2}}\nn\\
&&+ \frac{1}{\sqrt{\dhor}}\big(\A[\pmb \psi](\Iti)\big)^{\frac{1}{4}}\bigg(\widetilde{\M}[\pmb\psi]+\sum_{i,j=1}^3\int_{\Mntrap(\Iti)}|F_{ij}|^2\bigg)^{\frac{3}{4}}.
\eea
Combining this with the  estimate \eqref{esti:Err3:MoraEstimiddleregion}  for $\int\Err_{3}$, and using \eqref{eq:decompositionofX:Mora} and \eqref{def:errorterms:general:scalarizedwavefromtensorialwave}, we infer
\bea\lab{esti:sumofErr1234:MoraEstimiddleregion} 
&&\bigg|\sum_{i,j}\int \Re\Big(\square_{\g}\psi_{ij} \ov{(X+E)\psi_{ij}}\Big)- \sum_{i,j}\int\Re\Big(F_{ij}\ov{\big[(X+E)\psi_{ij}+A(T_{\tau_2}(\psi)_{ij}-\pr_{\tt}\psi_{ij})\big]}\Big)\bigg|\nn\\
&=&\bigg|\int \sum_{n=1}^4\Err_n- \sum_{i,j}\int\Re\Big(F_{ij} \ov{\big[(X+E)\psi_{ij}+A(T_{\tau_2}(\psi)_{ij}-\pr_{\tt}\psi_{ij})\big]}\Big)\bigg|\nn\\
&\les&\sup_{\tau\in[\tmic,\tau_1+2]}\E[\pmb \psi](\tau)+ \ep\EM[\pmb \psi](\Iti) +\A[\pmb \psi](\Iti)+ \NNtlocal[\pmb \psi](\Iti)\nn\\
&&+\frac{1}{\sqrt{\dhor}}\big(\A[\pmb \psi](\Iti)\big)^{\frac{1}{4}}\bigg(\widetilde{\M}[\pmb\psi]+\sum_{i,j=1}^3\int_{\Mntrap(\Iti)}|F_{ij}|^2\bigg)^{\frac{3}{4}}.
\eea

We are now ready to state a global Morawetz estimate in $\MM_{r_+(1+\dhor'),R}$ for solutions to the system of wave equations \eqref{eq:ScalarizedWaveeq:general:Kerrpert}. 

\begin{proposition}
\lab{prop:Morawetz:middleregion:globalextendedwavesystem}
Under the same assumptions for the scalars $\psi_{ij}$, $F_{ij}$ and the spacetime $(\MM,\g)$ as in Theorem \ref{th:mainenergymorawetzmicrolocal}, there exists  a constant $c>0$ such that  the following Morawetz estimate in $\MM_{r_+(1+\dhor'),{\Rmic}}$ for solutions to the coupled system of wave equations \eqref{eq:ScalarizedWaveeq:general:Kerrpert} holds true
\bea
\lab{eq:microlocalMora:middle:scalarizedwave:withFij}
\nn&&\sum_{i,j}\Bigg(c\Bigg[\int_{\MM_{r_+(1+\dhor'),{\Rmic}}}\frac{\mu^2|\pr_r\psi_{ij}|^2}{r^2} +\int_{\Mntrap_{r_+(1+\dhor'),{\Rmic}}}\frac{|\pr_\tau\psi_{ij}|^2+|\nab\psi_{ij}|^2}{r^2}\\
\nn&&+\int_{\MM_{r_+(1+\dhor'),10m}}\Big(|\Opw(\sigma_{\trap})\psi_{ij}|^2+|\Opw(e)\psi_{ij}|^2\Big)\Bigg]\\
&&+{\textbf{BDR}^{-}_{r={\Rmic}}[\psi_{ij}](\Iti)}+\int_{\MM_{r_+(1+\dhor'),{\Rmic}}}\Re\Big(F_{ij} \ov{(X+ E)\psi_{ij}}\Big)\Bigg)\nn\\
&\les_{\Rmic}& \sup_{\tau\in[\tmic,\tau_1+1]}\E[\pmb \psi](\tau)+ (\ep+\dhor)\widetilde{\EM}[\pmb \psi]+\NNtlocal[\pmb \psi](\Iti)\nn\\
&& +\frac{1}{\dhor^6}\A[\pmb \psi](\Iti)+(\ep+\dhor)\int_{\MM(\Iti)}|\pmb F|^2,
\eea
where $\NNtlocal[\pmb \psi](\Iti)$ is given by  \eqref{def:NNtlocalinr:NNtEner:NNtaux:wavesystem:EMF:proof}, where $\textbf{BDR}^{-}_{r={\Rmic}}[\psi_{ij}](\Iti)$ denotes a boundary term\footnote{The boundary term $\textbf{BDR}^{-}_{r={\Rmic}}[\psi_{ij}](\Iti)$ is the same as the one appearing in \eqref{eq:microlocalMora:middle:scalarizedwave:1} which in turn comes from the one in Lemma \ref{lem:conditionaldegenerateMorawetzflux:pertKerrrp1pdhorpR} with the substitution $\psi\to\psi_{ij}$.} on $H_{\Rmic}(\Iti)$, and where the symbols $\sigma_{\trap},\,  e\in\widetilde{S}^{1,0}(\MM)$ and the PDOs $X\in\Opw(\widetilde{S}^{1,1}(\MM)), E\in \Opw(\widetilde{S}^{0,0}(\MM))$ are defined as in Section \ref{sec:relevantmixedsymbolsonMM}.
\end{proposition}

\begin{proof}
In view of 
\eqref{eq:ScalarizedWaveeq:general:Kerrpert} and \eqref{hatSandV:generalwave:Kerrpert}, we have 
\beaa
\square_{\g}(\psi_{ij}) =F_{ij} +\sum_{k,l}O(r^{-2})\dk^{\leq 1}\psi_{kl},
\eeaa
and hence
\beaa
&&\sum_{i,j}\bigg((\ep+\dhor)\int_{\Mntrap(\Iti)}|\square_{\g}\psi_{ij}|^2+\ep\int_{\Mtrap(\Iti)}\tau^{-1-\dec}|\square_{\g}\psi_{ij}|^2\\
&&\qquad\qquad\qquad\qquad\qquad +(\ep+\dhor^6)\int_{\Mtrap}\Big|\Opw(\widetilde{S}^{-1,0}(\MM))\square_{\g}\psi_{ij}\Big|^2\bigg)\nn \\
&\les& (\ep+\dhor) \EM[\pmb \psi](\Iti) 
 +\sum_{i,j}\bigg((\ep+\dhor)\int_{\Mntrap(\Iti)}|F_{ij}|^2\nn\\
 &&+\ep\int_{\Mtrap(\Iti)}\tau^{-1-\dec}|F_{ij}|^2+(\ep+\dhor^6)\int_{\Mtrap}\left|\Opw(\widetilde{S}^{-1,0}(\MM))F_{ij}\right|^2\bigg).
\eeaa
Using the above estimate to control the first two lines of the RHS of \eqref{eq:microlocalMora:middle:scalarizedwave:1},  bounding the last term on the RHS of \eqref{eq:microlocalMora:middle:scalarizedwave:1} by $\good^{(1)}[\psi]$, and using the control of $\good^{(1)}[\psi]$ provided by Lemma \ref{lem:lowerorderterms:controlled:scalarizedwavefromtensorial}, we obtain 
\beaa
\nn&& \sum_{i,j}\Bigg(c\Bigg[\int_{\MM_{r_+(1+\dhor'),{\Rmic}}}\frac{\mu^2|\pr_r\psi_{ij}|^2}{r^2} +\int_{\Mntrap_{r_+(1+\dhor'),{\Rmic}}}\frac{|\pr_\tau\psi_{ij}|^2+|\nab\psi_{ij}|^2}{r^2}\\
\nn&&+\int_{\MM_{r_+(1+\dhor'),10m}}\Big(|\Opw(\sigma_{\trap})\psi_{ij}|^2+|\Opw(e)\psi_{ij}|^2\Big)\Bigg]\\
&&+\textbf{BDR}^{-}_{r={\Rmic}}[\psi_{ij}](\Iti)+\int_{\MM_{r_+(1+\dhor'),{\Rmic}}}\Re\Big(\square_{\g}\psi_{ij}\ov{(X+E)\psi_{ij}}\Big)\Bigg)\nn\\
\nn&\les_{\Rmic}& (\ep+\dhor) \EM[\pmb \psi](\Iti) +\dhor^{-6}\int_{\MM_{r_+(1+\dhor'), {\Rmic}}}|\pmb \psi|^2 
 +(\ep+\dhor)\int_{\Mntrap(\Iti)}|\F|^2\nn\\
 &&+\ep\int_{\Mtrap(\Iti)}\tau^{-1-\dec}|\F|^2+(\ep+\dhor^6)\int_{\Mtrap}\left|\Opw(\widetilde{S}^{-1,0}(\MM))\F\right|^2\\
\nn&& +\frac{1}{\sqrt{\dhor}}\big(\A[\pmb \psi](\Iti)\big)^{\frac{1}{4}}\bigg(\widetilde{\M}[\pmb\psi]+\int_{\Mntrap(\Iti)}|\F|^2\bigg)^{\frac{3}{4}}\\
&\les_{\Rmic}& (\ep+\dhor) \EM[\pmb \psi](\Iti) +\dhor^{-6}\int_{\MM_{r_+(1+\dhor'), {\Rmic}}}|\pmb \psi|^2 
 +(\ep+\dhor)\int_{\MM(\Iti)}|\F|^2\nn\\
\nn&& +\frac{1}{\sqrt{\dhor}}\big(\A[\pmb \psi](\Iti)\big)^{\frac{1}{4}}\bigg(\widetilde{\M}[\pmb\psi]+\int_{\MM(\Iti)}|\F|^2\bigg)^{\frac{3}{4}}.
\eeaa
Then, relying on \eqref{esti:sumofErr1234:MoraEstimiddleregion} to control the last term on the LHS, we infer
\bea\lab{eq:microlocalMora:middle:scalarizedwave:withFij:proof}
\nn&&\sum_{i,j}\Bigg(c\Bigg[\int_{\MM_{r_+(1+\dhor'),{\Rmic}}}\frac{\mu^2|\pr_r\psi_{ij}|^2}{r^2} +\int_{\Mntrap_{r_+(1+\dhor'),{\Rmic}}}\frac{|\pr_\tau\psi_{ij}|^2+|\nab\psi_{ij}|^2}{r^2}\\
\nn&&+\int_{\MM_{r_+(1+\dhor'),10m}}\Big(|\Opw(\sigma_{\trap})\psi_{ij}|^2+|\Opw(e)\psi_{ij}|^2\Big)\Bigg]\\
&&+\textbf{BDR}^{-}_{r={\Rmic}}[\psi_{ij}](\Iti)+\int_{\MM_{r_+(1+\dhor'),{\Rmic}}}\Re\Big(F_{ij} \ov{\big[(X+E)\psi_{ij}+A(T_{\tau_2}(\psi)_{ij}-\pr_{\tt}\psi_{ij})\big]}\Big)\Bigg)\nn\\
&\les_{\Rmic}& \sup_{\tau\in[\tmic,\tau_1+2]}\E[\pmb \psi](\tau)+(\ep+\dhor) \EM[\pmb \psi](\Iti) +\dhor^{-6}\int_{\MM}r^{-3}|\pmb \psi|^2 
 +(\ep+\dhor)\int_{\MM(\Iti)}|\F|^2\nn\\
\nn&& +\frac{1}{\sqrt{\dhor}}\big(\A[\pmb \psi](\Iti)\big)^{\frac{1}{4}}\bigg(\widetilde{\M}[\pmb\psi]+\int_{\MM(\Iti)}|\F|^2\bigg)^{\frac{3}{4}}+\A[\pmb \psi](\Iti)+ \NNtlocal[\pmb \psi](\Iti)\\
&\les_{\Rmic}&\sup_{\tau\in[\tmic,\tau_1+2]}\E[\pmb \psi](\tau)+(\ep+\dhor)\widetilde{\EM}[\pmb \psi]+\NNtlocal[\pmb \psi](\Iti)\nn\\
&& +\frac{1}{\dhor^6}\A[\pmb \psi](\Iti)+(\ep+\dhor)\int_{\MM(\Iti)}|\pmb F|^2.
\eea
Now, in view of \eqref{def:widehatTpsiij:equalsnabTpsiij} and the estimate
\eqref{estimates:Mialphaj:Kerrperturbations}
 for ${M_{i\tau}^j}$, we have
\beaa
T_{\tau_2}(\psi)_{ij}
&=&\pr_{\tt}\psi_{ij} + \sum_{k,l}\big(O(r^{-3}) + \Ga_b\big)\psi_{kl}
\eeaa
which yields
\begin{align*}
&\sum_{i,j=1}^3\bigg|\int_{\MM_{r_+(1+\dhor'),{\Rmic}}}\Re\Big(F_{ij} \ov{A(T_{\tau_2}(\psi)_{ij}-\pr_{\tt}\psi_{ij})}\Big)\bigg|\nn\\
\les{}&\int_{\MM_{r_+(1+\dhor'),{\Rmic}}} |\pmb F| |\pmb \psi|\les\frac{1}{\dhor}\int_{\MM_{r_+(1+\dhor'), {\Rmic}}}|\pmb \psi|^2
+\dhor \int_{\MM_{r_+(1+\dhor'), {\Rmic}}}|\pmb F|^2,
\end{align*}
hence substituting this estimate into \eqref{eq:microlocalMora:middle:scalarizedwave:withFij:proof} yields the desired estimate \eqref{eq:microlocalMora:middle:scalarizedwave:withFij}. This concludes the proof of Proposition \ref{prop:Morawetz:middleregion:globalextendedwavesystem}.
\end{proof}


\subsection{Energy-Morawetz estimates near infinity and redshift estimates}
\lab{sec:energyMorawetznearinfinityandredshiftestimates:globalsystmscalarwave}


In this section, we derive energy-Morawetz estimates near infinity and redshift estimates for solutions to the coupled system of  wave equations \eqref{eq:ScalarizedWaveeq:general:Kerrpert}. We start with the derivation of a divergence identity.

\begin{proposition}\lab{prop-app:stadard-comp-Psi:extendedscalarizedRW}
Let $\psi_{ij}$ be a solution to the system of wave equations \eqref{eq:ScalarizedWaveeq:general:Kerrpert}. Under the same assumptions for the scalars $\psi_{ij}$, $F_{ij}$ and the spacetime $(\MM,\g)$ as in Theorem \ref{th:mainenergymorawetzmicrolocal}, let $X$ be  a real-valued vectorfield satisfying
 \bea
\lab{eq:assump:generalvectorfieldX:EMnearinf}
X=O(1)\pr_r + O(1)\pr_{\tau} + O(r^{-2})\pr_{x^a},
\eea
let $w$ be a real scalar function satisfying $w=O(r^{-1})$, let $\chi_n(\tau)$, $n=1,2,3,4$, be smooth nonnegative cut-off functions satisfying \eqref{def:cutoffsintime1234}, and define the following modified current
\bea
\lab{eq:DefofwidetilePPmu:extendedscalarizedRW}
{\PP}_{\mu,\tau_2}[\pmb\psi](X, w)&:=&
\big(\chi_1(\tau)+\chi_3(\tau)\big)\sum_{i,j}{\PP}_\mu[\psi_{ij}](X,w)\nn\\
&&+\chi_2(\tau)
{\PP}_{\mu}[\pmb\psi](X, w)+\chi_4(\tau)
\big({\PP}_{\mu}[\pmb\psi](X, w)\big)_{K}, 
\eea
with ${\PP}_\mu[\pmb\psi](X, w)$ given in \eqref{definitionofcurrentPPmuXw:generaltensor} for $\pmb\psi\in\sk_2(\mathbb{C})$ and $V=D_0|q|^{-2}$, $\big({\PP}_\mu[\pmb\psi](X, w)\big)_{K}$ being the corresponding quantity in Kerr for $\pmb\psi\in\sk_2(\mathbb{C})$ and $V=\frac{4}{|q|^2}- \frac{4a^2\cos^2\th(|q|^2+6mr)}{|q|^6}$, and ${\PP}_\mu[\psi_{ij}](X,w)$ given in 
\eqref{definitionofcurrentPPmuXw:generaltensor} for a scalar $\psi_{ij}\in\sk_0(\mathbb{C})$ and $V=D_0|q|^{-2}$. 
 Then, ${\PP}_{\mu,\tau_2}[\pmb\psi](X, w)$ satisfies
 \bea
\lab{eq:1formPP_mu:equality}
{\PP}_{\mu,\tau_2}[\pmb\psi](X, w)=\sum_{i,j}{\PP}_\mu[\psi_{ij}](X,w)
+H_{\mu}[\psi,\pr\psi]
\eea
where $H_{\mu}[\psi,\pr\psi]$ satisfies \eqref{eq:1formBB_mu:equality:errortermestimates} for any $\tmic\leq \tau'<\tau''$, and its divergence equals
\bea
\lab{eq:DivwidetilePPmu:extendedscalarizedRW}
&&\D^\mu {\PP}_{\mu,\tau_2}[\pmb\psi](X, w)\nn\\
&=&\sum_{i,j}\Re\bigg(F_{ij}\ov{\bigg(X_{\tau_2}(\psi)_{ij} + \frac{1}{2}w\psi_{ij}  \bigg)}\bigg) \nn\\
&&
+(\chi_1(\tau)+\chi_3(\tau))\frac{1}{2}\sum_{i,j}\Big(\QQ[\psi_{ij}]  \c\piX +w \LL[\psi_{ij}]  -X( V) |\psi_{ij}|^2-\frac{1}{2}|\psi_{ij}|^2   \square_\g  w\Big) \nn \\
&&+\chi_2(\tau)\frac 1 2 \Big(\QQ[\pmb\psi]  \c\piX +w \LL[\pmb\psi] -X( V) |\pmb\psi|^2-\frac{1}{2}|\pmb\psi|^2   \square_\g  w\Big)\nn\\
&&+\chi_4(\tau)\frac 1 2 \Big(\QQ[\pmb\psi]  \c\piX +w \LL[\pmb\psi] -X( V) |\pmb\psi|^2-\frac{1}{2}|\pmb\psi|^2   \square_\g  w\Big)_{K}\nn\\
&&+\Err_{\text{l.o.t.}}, \quad\textrm{for}\quad \tau\geq\tmic,
\eea
 where 
\bea
\lab{def:Xtau2psiij:generalforEM}
X_{\tau_2}(\psi)_{ij}&:=& 
X(\psi_{ij}) - \chi_2\big(X^{\a}M_{i\a}^k\psi_{kj}+ X^{\a}M_{j\a}^k\psi_{ik}\big)\nn\\
&& -\chi_4\big(X^{\a} (M_K)_{i\a}^k\psi_{kj}+ X^{\a} (M_K)_{j\a}^k\psi_{ik}\big),
\eea
and where $\Err_{\text{l.o.t.}}$ denotes terms satisfying the following bound for any $r_1 \geq 10m$ and $(\tau',\tau'')\subset (\tmic,+\infty)$, 
\bea
\lab{def:Errlowerorderterm:DivPPmuchi}
\bigg|\int_{\MM_{r_1, +\infty}(\tau',\tau'')}\Err_{\text{l.o.t.}}\bigg|
&\les&\sup_{\tau\in[\tmic,\tau_1+2]}\E[\pmb \psi](\tau)+\ep\sup_{\tt\in[\tau_2-3,\tau_2+1]}\E[\pmb\psi](\tau) + \NNtlocal[\pmb \psi](\tau',\tau'')\nn\\
&&+\int_{H_{r_1}(\tau',\tau'')}r^{-1} |\pmb\psi| |\pr^{\leq 1}\pmb\psi|
\eea
with $\NNtlocal[\pmb \psi](\tau',\tau'')$ as given in \eqref{def:NNtlocalinr:NNtEner:NNtaux:wavesystem:EMF:proof}.
\end{proposition}

\begin{proof}
The identity \eqref{eq:1formPP_mu:equality} with $H_{\mu}[\psi,\pr\psi]$ satisfying \eqref{eq:1formBB_mu:equality:errortermestimates} follows from the definition \eqref{eq:DefofwidetilePPmu:extendedscalarizedRW} for ${\PP}_{\mu,\tau_2}[\pmb\psi](X, w)$ together with the estimates \eqref{eq:difference:PPmupmbpsiandpsiij} for a real-valued vectorfield $X=O(1)\pr_{\tt}+O(1)\pr_{r} + O(r^{-2})\pr_{x^a}$ and   a real-valued scalar function $w=O(r^{-1})$.

Next, we compute the divergence of ${\PP}_{\mu,\tau_2}[\pmb\psi](X, w)$.
We have, from the formula \eqref{eq:DefofwidetilePPmu:extendedscalarizedRW} for ${\PP}_{\mu,\tau_2}[\pmb\psi](X, w)$,
\bea
\lab{eq:computeDivergenceofBB:intermsofBulk}
 \D^\mu{\PP}_{\mu,\tau_2}[\pmb\psi](X, w)
 &=&\sum_{n=1}^4\Bulkxw{n}+\bigg(\big(\chi_1'(\tau)+\chi_3'(\tau)\big)\sum_{i,j}{\PP}_\mu[\psi_{ij}](X,w)  \nn\\
 &&\qquad\qquad\qquad\qquad+\chi_2'(\tau){\PP}_{\mu}[\pmb\psi](X, w)+\chi_4'(\tau)\big({\PP}_{\mu}[\pmb\psi](X, w)\big)_{K}\bigg)\D^\mu(\tau)\nn\\
&=&\sum_{n=1}^4\Bulkxw{n}+\chi_2'(\tau)\D^\mu(\tau)\bigg({\PP}_{\mu}[\pmb\psi](X, w)- \sum_{i,j}{\PP}_\mu[\psi_{ij}](X, w)\bigg)\nn\\
&&+\chi_4'(\tau)\D^\mu(\tau)\bigg(\big({\PP}_{\mu}[\pmb\psi](X, w)\big)_{K}- \sum_{i,j}{\PP}_\mu[\psi_{ij}](X, w)\bigg)\nn\\
 &=& \sum_{n=1}^4\Bulkxw{n}+\Err_{\text{l.o.t.}},
\eea
where we have defined in the first equality of \eqref{eq:computeDivergenceofBB:intermsofBulk}
\bsub
\lab{def:Bulkxw:1234}
\begin{align}
\Bulkxw{1} :=& \chi_1(\tau)\sum_{i,j}\D^\mu{\PP}_\mu[\psi_{ij}](X,w) , &\Bulkxw{2} :=&\chi_2(\tau)\D^\mu{\PP}_{\mu}[\pmb\psi](X, w),  \\
\Bulkxw{3} :=& \chi_3(\tau)\sum_{i,j}\D^\mu{\PP}_\mu[\psi_{ij}](X,w),&
\Bulkxw{4} :=&\chi_4(\tau)\D_K^\mu\big({\PP}_{\mu}[\pmb\psi](X, w)\big)_{K},
\end{align}
\esub
used in the second equality of \eqref{eq:computeDivergenceofBB:intermsofBulk} the fact $\sum_{n=1}^4\chi_n'(\tau)=0$ which follows from $\sum_{n=1}^4\chi_n(\tau)=1$, and used in the last equality of \eqref{eq:computeDivergenceofBB:intermsofBulk} the first estimate in \eqref{eq:difference:PPmupmbpsiandpsiij} and the support properties of $\chi_2'(\tau)$ and $\chi_4'(\tau)$.  

Next, we estimate the terms $\Bulkxw{n}$, $n=1,2,3,4$. We start with the term $\Bulkxw{2}$, which is supported on $\MM(\tau_1+1,\tau_2-2)$. By  multiplying on both sides of the formula \eqref{eq:DivofPPmu:tensor:RW:prop} by $\chi_2$, the scalars $\psi_{ij}=\pmb\psi(\Om_i,\Om_j)$, with $\pmb\psi\in\sk_2(\mathbb{C})$, satisfy on $\MM(\tau_1+1,\tau_2-2)$
\beaa
\Bulkxw{2}
&=&\Re\bigg(\chi_2\ov{\bigg(\nab_X \pmb\psi +\frac 1 2   w \pmb\psi\bigg)}\c \big(\squared_2 \pmb\psi -V\pmb\psi\big)\bigg) + \chi_2{}^{(X)}A_\nu\Im\Big(\pmb\psi\c\ov{\Ddot^{\nu}\pmb\psi}\Big)\nn\\
&&+\frac 1 2 \chi_2\Big(\QQ[\pmb\psi]  \c\piX +w \LL[\pmb\psi] -X( V ) |\pmb\psi|^2-\frac{1}{2}|\pmb\psi|^2   \square_\g  w\Big).
\eeaa
 In view of the form of the system of wave equations \eqref{eq:ScalarizedWaveeq:general:Kerrpert} and using \eqref{def:Mialphaj}, we have on $\MM(\tau_1+1,\tau_2-2)$ 
\bsub
\bea
\squared_2 \pmb\psi-V\pmb\psi &=&\F+\frac{4ia\cos\th}{\qs}\nab_{\pr_{\tt}}\pmb\psi ,\\
\lab{eq:nabXpsi+wpsi:intermsofscalars}
\bigg(\nab_X \pmb\psi +\frac 1 2   w \pmb\psi\bigg)_{ij}&=&X(\psi_{ij}) - \big(X^{\a}M_{i\a}^k\psi_{kj}+ X^{\a}M_{j\a}^k\psi_{ik}\big) +\frac{1}{2}w\psi_{ij},
\eea
\esub
and in view of Corollary \ref{cor:asymtpticbehavioroftheRdottermintensorialenergyidentity:larger}, we have
\beaa
\Big|\chi_2{}^{(X)}A_\nu\Im\Big(\pmb\psi\c\ov{\Ddot^{\nu}\pmb\psi}\Big)\Big|\les \chi_2 r^{-3}\big(|\nab_3\pmb\psi|+|\nab_4\pmb\psi|\big)|\pmb\psi| + \chi_2 r^{-2} |\nab\pmb\psi||\pmb\psi|.
\eeaa
The above implies
\bea
\lab{eq:DivPPmuchi:general:proof:part1}
\Bulkxw{2} &=&\sum_{i,j}\Re\bigg(\chi_2\ov{\bigg(X(\psi_{ij}) - \big(X^{\a}M_{i\a}^k\psi_{kj}+ X^{\a}M_{j\a}^k\psi_{ik}\big) +\frac{1}{2}w\psi_{ij}\bigg)}F_{ij}\bigg)
\nn\\
&&+\frac 1 2\chi_2\Big(\QQ[\pmb\psi]  \c\piX +w \LL[\pmb\psi] -X( V ) |\pmb\psi|^2-\frac{1}{2}|\pmb\psi|^2   \square_\g  w\Big)\nn\\
&& +\chi_2 O(r^{-3})\big(|\nab_3\pmb\psi|+|\nab_4\pmb\psi|\big)|\pmb\psi| + \chi_2O(r^{-2}) |\nab\pmb\psi||\pmb\psi|\nn\\
&&+\Re\bigg(\ov{\chi_2\bigg(\nab_X \pmb\psi +\frac 1 2   w \pmb\psi\bigg)}\c \frac{4ia\cos\th}{\qs}\nab_{\pr_{\tt}}\pmb\psi\bigg) .
\eea
In view of the estimate \eqref{estimates:Mialphaj:Kerrperturbations} and the identity \eqref{eq:nabXpsi+wpsi:intermsofscalars}, and since the real-valued vectorfield $X$ and  the real-valued function $w$ satisfy respectively $X=O(1)\pr_r + O(1)\pr_{\tau} + O(r^{-2})\pr_{x^a}$ and $w=O(r^{-1})$, the term in the last line of the RHS of \eqref{eq:DivPPmuchi:general:proof:part1} equals
\beaa
\Re\bigg(\ov{\chi_2X^{\a}\pr_{\a}(\psi_{ij})}\c \frac{4ia\cos\th}{\qs}\pr_{\tt}(\psi_{ij})\bigg) +\chi_2 O(r^{-3})|\pmb\psi| |\pr^{\leq 1}\pmb\psi|.
\eeaa
Hence, applying the estimate \eqref{eq:integrationbypartsforitimesfirstordertimesfirstorder:general} with $f=\chi_2f_0X^{\a}\frac{2a\cos\th }{\qs}$ and $\pr_{\b}=\pr_{\tt}$ to  the integral of the first term over $\MM_{r_1,+\infty}(\tau',\tau'')$, and using also the fact that $\MM_{r_1,+\infty}(\tau',\tau'')\subset\Mntrap(\tau',\tau'')$, we deduce for the integral over $\MM_{r_1,+\infty}(\tau',\tau'')$ of the last two lines of \eqref{eq:DivPPmuchi:general:proof:part1} that
\beaa
&& \bigg|\int_{\MM_{r_1,+\infty}(\tau',\tau'')} \bigg(O(r^{-3})\big(|\nab_3\pmb\psi|+|\nab_4\pmb\psi|\big)|\pmb\psi| + O(r^{-2}) |\nab\pmb\psi||\pmb\psi|\nn\\
&&\qquad\qquad\qquad\quad+\Re\bigg(\ov{\chi_2\bigg(\nab_X \pmb\psi +\frac 1 2   w \pmb\psi\bigg)}\c \frac{4ia\cos\th}{\qs}\nab_{\pr_{\tt}}\pmb\psi\bigg)\bigg)\bigg|\nn\\
&\les&\NNtlocal[\pmb \psi](\tau',\tau'')+\int_{H_{r_1}(\tau',\tau'')}r^{-1}|\pmb\psi| |\pr^{\leq 1}\pmb\psi|,
\eeaa
which together with the identity \eqref{eq:DivPPmuchi:general:proof:part1} yields
\begin{align}
\lab{eq:Bulkxw:general:proof:part2}
\Bulkxw{2}
={}&\sum_{i,j}\Re\bigg(\ov{\chi_2\bigg(X(\psi_{ij}) - \big(X^{\a}M_{i\a}^k\psi_{kj}+ X^{\a}M_{j\a}^k\psi_{ik}\big) +\frac{1}{2}w\psi_{ij}\bigg)}F_{ij}\bigg)
\nn\\
&+\chi_2\frac 1 2 \Big(\QQ[\pmb\psi]  \c\piX +w \LL[\pmb\psi] -X( V ) |\pmb\psi|^2-\frac{1}{2}|\pmb\psi|^2   \square_\g  w\Big)+\Err_{\text{l.o.t.}} .
\end{align}

Proceeding in exactly the same manner, we infer, for $\Bulkxw{4}$, 
\begin{align}
\lab{eq:Bulkxw:general:proof:part4}
\Bulkxw{4}
={}&\sum_{i,j}\Re\bigg(\chi_4\ov{\bigg(X(\psi_{ij}) - \big(X^{\a}(M_K)_{i\a}^k\psi_{kj}+ X^{\a}(M_K)_{j\a}^k\psi_{ik}\big) +\frac{1}{2}w\psi_{ij}\bigg)}F_{ij}\bigg)
\nn\\
&+\chi_4\frac 1 2 \Big(\QQ[\pmb\psi]  \c\piX +w \LL[\pmb\psi] -X( V) |\pmb\psi|^2-\frac{1}{2}|\pmb\psi|^2   \square_\g  w\Big)_{K}+\Err_{\text{l.o.t.}} .
\end{align}

Next, we consider the term $\Bulkxw{3}$ which is supported on $\MM(\tau_2-3,\tau_2+1)$. Applying \eqref{eq:DivofPPmu:tensor:RW:prop} to the wave equations \eqref{eq:scalarizedextendedeqs:psiij:moraandener} for $\psi_{ij}$ with $\{G_{n,ij}\}_{n=1,2,3,4}$ given in \eqref{def:Gnij:energyesti:withpotential}, multiplying on both sides by $\chi_3$, and summing over $i,j$, we deduce
\bea
\lab{eq:energyidentity:generalvectorfieldX:witherrortermsexplicit}
\Bulkxw{3} 
&=&\sum_{i,j}\chi_3\D^\mu{\PP}_\mu[\psi_{ij}](X,w)\nn\\
&=&\sum_{i,j}\Re\bigg(F_{ij}\chi_3\ov{\bigg(X(\psi_{ij}) + \frac{1}{2}w\psi_{ij}  \bigg)}\bigg) \nn\\
&&
+ \chi_3\frac{1}{2}\sum_{i,j}\Big(\QQ[\psi_{ij}]  \c\piX +w \LL[\psi_{ij}]  -X( V) |\psi_{ij}|^2-\frac{1}{2}|\psi_{ij}|^2   \square_\g  w\Big) \nn \\
&&+\sum_{i,j}\sum_{n=1}^4\Re\bigg(\chi_3G_{n,ij}\ov{\bigg(X(\psi_{ij}) + \frac{1}{2}w\psi_{ij}  \bigg)}\bigg).
\eea
It remains to  estimate the term in the last line of  the above equation \eqref{eq:energyidentity:generalvectorfieldX:witherrortermsexplicit}. In view of \eqref{estimates:widehatSandQ:perturbationsofKerr} and the fact that $f_{D_0}=O(r^{-2})$ and $w=O(r^{-1})$, we have
\beaa
\sum_{i,j}\sum_{n=1}^4\Re\big(\chi_3G_{n,ij}\ov{w\psi_{ij}  }\big)=\chi_3O(r^{-2}) |\dk^{\leq 1}\pmb\psi| |\pmb\psi|,
\eeaa
and further, in view of the assumption \eqref{eq:assump:generalvectorfieldX:EMnearinf} for the vectorfield $X$, we have
\beaa
\sum_{i,j}\sum_{n=3}^4\Re\big(\chi_3G_{n,ij}\ov{X\psi_{ij} }\big)=\chi_3\Big(O(r^{-2}) |\dk^{\leq 1}\pmb\psi| |\pmb\psi| + \Ga_{g}|\dk\pmb\psi|^2\Big).
\eeaa

For the term with $G_{2,ij}=\frac{4ia\cos\th}{|q|^2} \pr_{\tt}\psi_{ij}$, we integrate the differential identity \eqref{eq:integrationbypartsforitimesfirstordertimesfirstorder:general} with $f=\chi_3\frac{2a\cos\th}{|q|^2}f_0X^\a$ and $x^\b=\tau$. Using \eqref{eq:spacetimevolumeformusingisochorecoordinates} and \eqref{def:dVref}, we deduce, for any $r_1 \geq 10m$, 
\beaa
&&\bigg|\int_{\MM_{r_1, +\infty}(\tau',\tau'')}\Re\big(\chi_3G_{2,ij}\ov{X\psi_{ij} }\big)\bigg|\nn\\
&\les& \int_{\MM_{r_1, +\infty}(\tau_2-3,\tau_2+1)}r^{-2} |\psi| |\pr_{\tt}\psi|+
\sup_{\tau\in [\tau_2-3,\tau_2+1]}\int_{\Si(\tau)}r^{-2} |\psi| |\pr^{\leq 1}\psi|\nn\\
&&+\int_{H_{r_1}(\tau',\tau'')} r^{-2} |\psi| |\pr_{\tt}\psi|
+\left(\int_{\II_+(\tau', \tau'')}r^{-2}|\pmb\psi|^2\right)^{\frac{1}{2}}\Big(\F_{\II_+}[\pmb \psi](\tau',\tau'')\Big)^{\frac{1}{2}}\nn\\
&\les& \NNtlocal[\pmb \psi](\tau',\tau'') + \int_{H_{r_1}(\tau',\tau'')} r^{-2} |\psi| |\pr_{\tt}\psi|.
\eeaa

For the term with $G_{1,ij}$, proceeding exactly as for the proof of \eqref{eq:erroG1:prtau:middleregion:after}, we obtain
\beaa
\bigg|\int_{\MM_{r_1, +\infty}(\tau',\tau'')}\Re\big(\chi_3G_{1,ij}\ov{X\psi_{ij}}\big)\bigg|
&\les& \sup_{\tau\in [\tau_2-3,\tau_2+1]}\big(\E[\pmb\psi](\tau)\big)^{\frac{1}{2}}\bigg(\int_{\Sigma(\tt)}r^{-2}|\pmb\psi|^2+\Errdefect[\pmb\psi]\bigg)^{\frac{1}{2}}\nn\\
&&+\int_{H_{r_1}(\tau',\tau'')} r^{-2} |\psi| |\pr\psi|.
\eeaa

In view of the above estimates, we infer for any $r_1 \geq 10m$, 
\beaa
&&\bigg|\int_{\MM_{r_1, +\infty}(\tau',\tau'')}  \sum_{i,j}\sum_{n=1}^4\Re\bigg(\chi_3G_{n,ij}\ov{\bigg(X(\psi_{ij}) + \frac{1}{2}w\psi_{ij}  \bigg)}\bigg)\bigg|\nn\\
&\les&\ep\sup_{\tau\in [\tau_2-3,\tau_2+1]}\E[\pmb\psi](\tau) +  \sup_{\tau\in [\tau_2-3,\tau_2+1]}\big(\E[\pmb\psi](\tau)\big)^{\frac{1}{2}}\bigg(\int_{\Sigma(\tt)}r^{-2}|\pmb\psi|^2+\Errdefect[\pmb\psi]\bigg)^{\frac{1}{2}}\nn\\
&&+\int_{H_{r_1}(\tau',\tau'')} r^{-2} |\pmb\psi| |\pr\pmb\psi|
\eeaa
which implies
\beaa
\sum_{i,j}\sum_{n=1}^4\Re\bigg(\chi_3G_{n,ij}\ov{\bigg(X(\psi_{ij}) + \frac{1}{2}w\psi_{ij}  \bigg)}\bigg)=\Err_{\text{l.o.t.}}
\eeaa
and hence, together with \eqref{eq:energyidentity:generalvectorfieldX:witherrortermsexplicit},
\bea
\lab{eq:Bulkxw:general:proof:part3}
\Bulkxw{3}
&=&\chi_3\frac{1}{2}\sum_{i,j}\Big(\QQ[\psi_{ij}]  \c\piX +w \LL[\psi_{ij}]  -X( V) |\psi_{ij}|^2-\frac{1}{2}|\psi_{ij}|^2   \square_\g  w\Big)
\nn\\
&&
+\sum_{i,j}\Re\bigg(F_{ij}\chi_3\ov{\bigg(X(\psi_{ij}) + \frac{1}{2}w\psi_{ij}  \bigg)}\bigg) +\Err_{\text{l.o.t.}} .
\eea

In the end, we consider the term $\Bulkxw{1}$ which is supported on $\MM(\tmic, \tau_1+2)$. We have the same identity \eqref{eq:energyidentity:generalvectorfieldX:witherrortermsexplicit} as in the previous case for $\Bulkxw{3}$, and it follows from \eqref{estimates:widehatSandQ:perturbationsofKerr} that 
\beaa
\bigg|\sum_{i,j}\sum_{n=1}^4\Re\bigg(\chi_1 G_{n,ij}\ov{\bigg(X\psi_{ij} + \frac{1}{2}w\psi_{ij}  \bigg)}\bigg)\bigg|\les \chi_1 r^{-2} |\dk^{\leq 1}\pmb\psi|^2,
\eeaa
and hence
\beaa
\int_{\MM_{r_1, +\infty}(\tau',\tau'')}\bigg|\sum_{i,j}\sum_{n=1}^4\Re\bigg(\chi_1G_{n,ij}\ov{\bigg(X\psi_{ij} + \frac{1}{2}w\psi_{ij}  \bigg)}\bigg)\bigg|\les \sup_{\tau\in[\tmic,\tau_1+2]}\E[\pmb \psi](\tau)
\eeaa
which yields
\bea
\lab{eq:Bulkxw:general:proof:part1}
\Bulkxw{1}
&=&\chi_1\frac{1}{2}\sum_{i,j}\Big(\QQ[\psi_{ij}]  \c\piX +w \LL[\psi_{ij}]  -X( V) |\psi_{ij}|^2-\frac{1}{2}|\psi_{ij}|^2   \square_\g  w\Big)
\nn\\
&&
+\sum_{i,j}\Re\bigg(F_{ij}\chi_1\ov{\bigg(X(\psi_{ij}) + \frac{1}{2}w\psi_{ij}  \bigg)}\bigg) +\Err_{\text{l.o.t.}} .
\eea

Plugging  the estimates \eqref{eq:Bulkxw:general:proof:part2}, \eqref{eq:Bulkxw:general:proof:part4}, \eqref{eq:Bulkxw:general:proof:part3} and \eqref{eq:Bulkxw:general:proof:part1} into \eqref{eq:computeDivergenceofBB:intermsofBulk} and using $\sum_{n=1}^4\chi_n=1$, we infer the identity  \eqref{eq:DivwidetilePPmu:extendedscalarizedRW} with $X_{\tau_2}$ defined as in \eqref{def:Xtau2psiij:generalforEM} and with $\Err_{\text{l.o.t.}}$ satisfying \eqref{def:Errlowerorderterm:DivPPmuchi}, which then concludes the proof of Proposition \ref{prop-app:stadard-comp-Psi:extendedscalarizedRW}.
\end{proof}

Next, we deduce energy and Morawetz estimates for the globally extended coupled system of wave equations \eqref{eq:ScalarizedWaveeq:general:Kerrpert} in a large radius region away from the trapping region. 

\begin{proposition}
\lab{prop:energymorawetznearinfinity:extendedRWsystem:kerrpert}
Let $\psi_{ij}$ be a solution to the system of wave equations \eqref{eq:ScalarizedWaveeq:general:Kerrpert}, and let $T_{\tau_2}(\psi)_{ij}$, $X_{\tau_2}(\psi)_{ij}$ and $\NNtlocal[\pmb \psi](\tau',\tau'')$ be given as in \eqref{def:widehatTpsiij:equalsnabTpsiij}, \eqref{def:Xtau2psiij:generalforEM} and 
\eqref{def:NNtlocalinr:NNtEner:NNtaux:wavesystem:EMF:proof}, respectively. 
 Then, under the assumptions of Theorem \ref{th:mainenergymorawetzmicrolocal} for the scalars $\psi_{ij}$, $F_{ij}$ and the spacetime $(\MM,\g)$:
 \begin{itemize}
 \item there exists a constant $c>0$ such that the following energy estimate holds, for  any  constant $r_1\geq 10m$,
\bea
\label{eq:energynearinf:extendedsystem:Kerrpert:1}
&&c\Big(\E_{r\geq r_1}[\pmb\psi](\tau'')
+\F_{\II_+}[\pmb\psi](\tau',\tau'')\Big)
-\int_{H_{r_1}(\tau',\tau'')}\sum_{i,j}\PP_\mu[\psi_{ij}](\pr_{\tt}, 0)N_{H_r}^{\mu}\nn\\
&&+\int_{\MM_{r\geq r_1}(\tau',\tau'')}\sum_{i,j}\Re\big(F_{ij}\ov{T_{\tau_2}(\psi)_{ij}}\big) 
\nn \\
&\les& \sup_{\tau\in[\tmic,\tau_1+2]}\E[\pmb \psi](\tau)+\E_{r\geq r_1}[\pmb\psi](\tau') +\int_{H_{r_1}(\tau',\tau'')} (r_1)^{-1}|\psi| |\pr^{\leq 1}\psi|  \nn\\
&&+ \NNtlocal[\pmb \psi](\tau',\tau'') + \ep\EM[\pmb\psi](\tau',\tau'');
\eea

 \item there exists a constant $c>0$ such that the following energy estimate holds, for  any  constant $r_1\geq 11m$,
\bea
\label{eq:energynearinf:extendedsystem:Kerrpert:2:cutoff}
&&c\Big(\E_{r\geq r_1}[\pmb\psi](\tau'')
+\F_{\II_+}[\pmb\psi](\tau',\tau'')\Big)
+\int_{\MM_{r\geq r_1}(\tau',\tau'')}\sum_{i,j}\Re\big(F_{ij}\ov{T_{\tau_2}(\psi)_{ij}}\big) 
\nn \\
&\les& \sup_{\tau\in[\tmic,\tau_1+2]}\E[\pmb \psi](\tau)+\E_{r\geq r_1-m}[\pmb\psi](\tau')  + \NNtlocal[\pmb \psi](\tau',\tau'') \nn\\
&&+ \ep\EM[\pmb\psi](\tau',\tau'')
+\M_{r_1-m, r_1}[\pmb\psi](\tau',\tau'')
+\int_{\MM_{r_1-m, r_1}(\tau',\tau'')}|\pmb F|^2;
\eea

\item there exists  a constant $c>0$ such that the following Morawetz estimate holds for  a suitably large constant $R_1\gg 12m$:
\bea
\label{eq:Morawetznearinf:extendedsystem:Kerrpert:1:00}
\nn&& c\M_{r\geq R_1}[\pmb\psi](\tau',\tau'') -\int_{H_{R_1}(\tau',\tau'')}\sum_{i,j}\PP_\mu[\psi_{ij}](X_1, w_1)N_{H_r}^{\mu}\\
&&+\int_{\MM_{r\geq R_1}(\tau',\tau'')}\sum_{i,j}\Re\left(F_{ij}\ov{\left((X_{1})_{\tau_2}(\psi)_{ij} +\frac{1}{2}w_1\psi_{ij}\right)}\right) 
\nn \\
&\les&\sup_{\tau\in[\tmic,\tau_1+2]}\E[\pmb \psi](\tau)+\E_{r\geq R_1}[\pmb\psi](\tau'')
+\E_{r\geq R_1}[\pmb\psi](\tau') 
+\F_{\II_+}[\pmb\psi](\tau',\tau'')
 \nn\\
&&
+\ep\sup_{\tt\in[\tau_2-3,\tau_2+1]}\E[\pmb\psi](\tau) +\int_{H_{R_1}(\tau',\tau'')}(R_1)^{-1}|\pmb\psi| |\pr^{\leq 1}\pmb\psi|+ \NNtlocal[\pmb \psi](\tau',\tau''),
\eea
where\footnote{Note that the value of the function $w$ here is twice the value of the function $w$ chosen in \cite[Lemma 3.10]{MaSz24} which is due to a different normalization in the definition of the current $\PP_{\mu}[\pmb\psi](X,w)$. This is also the case for the choice of the function $w_{\de}$ in \eqref{def:Xandw:improvedMorawetz:extendedRW:Kerrpert}.}
\beaa
X_1=2\mu(1-mr^{-1})\pr_{r}^{\text{BL}},\qquad w_1=4\mu r^{-1}(1-mr^{-1});
\eeaa

\item for any $\de\in(0,1]$, there exists  a constant $c>0$ such that the following Morawetz estimate holds for  a suitably large constant $R_1\gg 12m$:
\bea
\label{eq:Morawetznearinf:extendedRWsystem:Kerrpert:delta}
\nn&& c\de\M_{\de, r\geq R_1}[\pmb\psi](\tau',\tau'')+\int_{\MM_{r\geq R_1-m}(\tau',\tau'')}\sum_{i,j}\Re\left(F_{ij}\ov{\left((X_{\de})_{\tau_2}(\psi)_{ij} +\frac{1}{2}w_{\de}\psi_{ij}\right)}\right) 
\nn \\
&\les&\sup_{\tau\in[\tmic,\tau_1+2]}\E[\pmb \psi](\tau)+\E_{r\geq R_1-m}[\pmb\psi](\tau'')
+\E_{r\geq R_1-m}[\pmb\psi](\tau') 
+\F_{\II_+}[\pmb\psi](\tau',\tau'')
\nn\\
&&
+\ep\sup_{\tt\in[\tau_2-3,\tau_2+1]}\E[\pmb\psi](\tau) +\M_{R_1-m,R_1}[\pmb\psi](\tau',\tau'')+ \NNtlocal[\pmb \psi] (\tau',\tau''),
\eea
with
\bea
\lab{def:Xandw:improvedMorawetz:extendedRW:Kerrpert}
X_{\de}=2\mu f_{\de}{\pr}_r^{\text{BL}}, \quad w_{\de}=4\mu h_{\de}, \quad f_{\de}=\chi_{R_1}(1-m^{\de}r^{-\de}),\quad 
h_{\de}=\chi_{R_1} r^{-1} (1-m^{\de}r^{-\de}),
\eea
where $\chi_{R_1}=\chi_{R_1}(r)$ is a smooth cutoff function that equals $1$ for $r\geq R_1$ and vanishes for $r\leq R_1-m$.
 \end{itemize}
\end{proposition}

\begin{proof}
Integrating the divergence identity \eqref{eq:1formBB_mu:divergence} in $\MM_{r_1, +\infty}(\tau',\tau'')$, where $r_1\geq 10m$, and making use of the identity 
\eqref{eq:1formBB_mu:equality} for $\BB_{\mu}[\pmb\psi]$ and the estimate \eqref{eq:1formBB_mu:equality:errortermestimates}, as well as the estimate  \eqref{def:Gpsiprpsi:lowerordertermofDivofBmu} for $G[\psi,\pr\psi]$ which appears on the RHS of the identity \eqref{eq:1formBB_mu:divergence}, we deduce the following energy estimate
\beaa
&&c\Big(\E_{r\geq r_1}[\pmb\psi](\tau'')
+\F_{\II_+}[\pmb\psi](\tau',\tau'')\Big)
-\int_{H_{r_1}(\tau',\tau'')}\sum_{i,j}\PP_\mu[\psi_{ij}](\pr_{\tt}, 0)N_{H_r}^{\mu}\nn\\
&&+\int_{\MM_{r\geq r_1}(\tau',\tau'')}\sum_{i,j}\Re\big(F_{ij}\ov{T_{\tau_2}(\psi)_{ij}}\big) 
\nn \\
&\les& \sup_{\tau\in[\tmic,\tau_1+2]}\E[\pmb \psi](\tau)+\E_{r\geq r_1}[\pmb\psi](\tau') +\int_{H_{r_1}(\tau',\tau'')}(r_1)^{-1} |\psi| |\pr^{\leq 1}\psi|  \nn\\
&&+ \NNtlocal[\pmb \psi](\tau',\tau'') + \ep\EM[\pmb\psi](\tau',\tau'')
\eeaa
with $c>0$ a constant. This proves the energy estimate \eqref{eq:energynearinf:extendedsystem:Kerrpert:1}. 

Next, we consider the other energy estimate \eqref{eq:energynearinf:extendedsystem:Kerrpert:2:cutoff}. This follows easily from applying the energy estimate \eqref{eq:energynearinf:extendedsystem:Kerrpert:1} proven above to the wave equations for $\chi_{r_1}(r)\psi_{ij}$\footnote{Notice that the cutoff function $\chi_{r_1}(r)$  induces an additional term $\M_{r_1-m,r_1}[\pmb\psi](\tau',\tau'')$ on the RHS of \eqref{eq:energynearinf:extendedsystem:Kerrpert:2:cutoff} instead of a boundary term at $H_{r_1}$.}, where $\chi_{r_1}(r)$ is a smooth cut-off function in $r$ satisfying $\chi_{r_1}(r)=1$ on $r\geq r_1$ and $\chi_{r_1}(r)=0$ on $r\leq r_1-m$.

Next, we consider the Morawetz estimate. It is shown in the proof of Lemma 3.10 in \cite{MaSz24} that with the choice of
\bea
\lab{eq:choicesofX0w:Morawetzde=1:Kerrpert}
X_1=2\mu(1-mr^{-1})\pr_{r}^{\text{BL}},\qquad w_1=4\mu r^{-1}(1-mr^{-1}), 
\eea
we have for both $V=D_0|q|^{-2}$ and $V=\frac{4}{|q|^2}- \frac{4a^2\cos^2\th(|q|^2+6mr)}{|q|^6}$, and for $r\geq R_1$, with $R_1$ suitably large\footnote{The reason that we have $\frac{|\psi_{ij}|^2}{r^3}$ instead of $\frac{|\psi_{ij}|^2}{r^4}$ as in \cite[Lemma 3.10]{MaSz24}  lies in the fact that we have a positive potential $V$ which furthermore satisfies $-X_1(V)\gtrsim r^{-3}$ for $r$ large enough.},
\beaa
&&\bigg(\QQ[\psi_{ij}]  \c {^{(X_1)}}\pi +w_1 \LL[\psi_{ij}]  -X_1( V) |\psi_{ij}|^2- |\psi_{ij}|^2   \square_\g  w_1\bigg)\nn\\
&\gtrsim & \frac{|\pr_{\tt}\psi_{ij}|^2}{r^2}
+\frac{|\pr_{r}\psi_{ij}|^2}{r^2} +\frac{|\nab\psi_{ij}|^2}{r}+\frac{|\psi_{ij}|^2}{r^3},
\eeaa
and 
\begin{equation}
\lab{esti:fluxes:Morawetznearinfinity:control:scalar:Kerrpert}
\begin{split}
&\int_{\Sigma_{r\geq R_1}(\tau)}\big|{\PP}_\mu[\psi_{ij}](X_1, w_1)N_{\Sigma_{\tau}}^{\mu}\big|\les \E_{r\geq R_1}[\psi_{ij}](\tau),\\
&\int_{\II_+(\tau',\tau'')}\big|{\PP}_\mu[\psi_{ij}](X_1, w_1)N_{\II_+}^{\mu} \big|\les \F_{\II_+}[\psi_{ij}](\tau',\tau'').
\end{split}
\end{equation}
By the same argument, we have
\beaa
&&\bigg(\QQ[\pmb\psi]  \c {^{(X_1)}}\pi +w_1 \LL[\pmb\psi]  -X_1(V) |\pmb\psi|^2-\frac{1}{2}|\pmb\psi|^2   \square_\g  w_1\bigg)\nn\\
&\gtrsim& \frac{|\nab_{\pr_{\tt}}\pmb\psi|^2}{r^2}
+\frac{|\nab_{\pr_{r}}\pmb\psi|^2}{r^2} +\frac{|\nab\pmb\psi|^2}{r}+\frac{|\pmb\psi|^2}{r^3},
\eeaa
and hence, we deduce from the above that there is a  constant $c>0$ such that 
\bea
\lab{eq:Morawetznearinfinity:extendedsystem:controlofmainterm:Kerrpert}
&&\int_{\MM_{r\geq R_1}(\tau',\tau'')}\bigg\{\chi_2(\tau)\frac 1 2 \Big(\QQ[\pmb\psi]  \c{}^{(X_1)}\pi +w_1 \LL[\pmb\psi] -X_1( V) |\pmb\psi|^2-\frac{1}{2}|\pmb\psi|^2   \square_\g  w_1\Big)
 \nn \\
&&\qquad+(\chi_1(\tau)+\chi_3(\tau))\frac{1}{2}\sum_{i,j}\Big(\QQ[\psi_{ij}]  \c{}^{(X_1)}\pi +w_1 \LL[\psi_{ij}]  -X_1( V) |\psi_{ij}|^2-\frac{1}{2}|\psi_{ij}|^2   \square_\g  w_1\Big)\nn\\
&&\qquad+\chi_4(\tau)\frac 1 2 \Big(\QQ[\pmb\psi]  \c{}^{(X_1)}\pi +w_1 \LL[\pmb\psi] -X_1( V) |\pmb\psi|^2-\frac{1}{2}|\pmb\psi|^2   \square_\g  w_1\Big)_{K}\bigg\}\nn\\
&\geq &c\M_{r\geq R_1}[\pmb \psi](\tau',\tau'').
\eea

In view of Proposition \ref{prop-app:stadard-comp-Psi:extendedscalarizedRW}, for a real-valued vectorfield $X=O(1)\pr_{\tt}+O(1)\pr_{r} + O(r^{-2})\pr_{x^a}$ and  a real-valued scalar function $w=O(r^{-1})$, ${\PP}_{\mu,\tau_2}[\pmb\psi](X, w)$ defined in \eqref{eq:DefofwidetilePPmu:extendedscalarizedRW} satisfies the identity \eqref{eq:1formPP_mu:equality} with 
$H_{\mu}[\psi,\pr\psi]$ satisfying the  bound \eqref{eq:1formBB_mu:equality:errortermestimates} for any $\tmic\leq \tau'<\tau''$.  From \eqref{eq:choicesofX0w:Morawetzde=1:Kerrpert}, $X_1$ and $w_1$ verify the above conditions, so we combine this estimate with the estimate \eqref{esti:fluxes:Morawetznearinfinity:control:scalar:Kerrpert} and deduce
\bea
&&\bigg|\int_{\MM_{r\geq R_1}(\tau',\tau'')}\D^\mu {\PP}_{\mu,\tau_2}[\pmb\psi](X_1, w_1)-\int_{H_{R_1}(\tau',\tau'')}\sum_{i,j}\PP_\mu[\psi_{ij}](X_1, w_1)N_{H_r}^{\mu}\bigg| \nn\\
&\les &
\E_{r\geq R_1}[\pmb\psi](\tau'')
+\E_{r\geq R_1}[\pmb\psi](\tau') 
+\F_{\II_+}[\pmb\psi](\tau',\tau'')\nn\\
&&
+\int_{H_{R_1}(\tau',\tau'')}(R_1)^{-1}|\pmb\psi| |\pr^{\leq 1}\pmb\psi| + \NNtlocal[\pmb \psi](\tau',\tau'').
\eea
Plugging this estimate and the estimate \eqref{eq:Morawetznearinfinity:extendedsystem:controlofmainterm:Kerrpert} into the divergence identity 
\eqref{eq:DivwidetilePPmu:extendedscalarizedRW} integrated over $\MM_{R_1,+\infty}(\tau',\tau'')$ with $(X,w)=(X_1,w_1)$, and using the bound \eqref{def:Errlowerorderterm:DivPPmuchi} for $\Err_{\text{l.o.t.}}$,
we infer
\bea
\nn&& c\M_{r\geq R_1}[\pmb\psi](\tau',\tau'') -\int_{H_{R_1}(\tau',\tau'')}\sum_{i,j}\PP_\mu[\psi_{ij}](X_1, w_1)N_{H_r}^{\mu}\\
&&+\int_{\MM_{r\geq R_1}(\tau',\tau'')}\sum_{i,j}\Re\left(F_{ij}\ov{\left((X_1)_{\tau_2}(\psi)_{ij} +\frac{1}{2}w_1\psi_{ij}\right)}\right) 
\nn \\
&\les&\sup_{\tau\in[\tmic,\tau_1+2]}\E[\pmb \psi](\tau)+\E_{r\geq R_1}[\pmb\psi](\tau'')
+\E_{r\geq R_1}[\pmb\psi](\tau') 
+\F_{\II_+}[\pmb\psi](\tau',\tau'')
\nn\\
&&
+\ep\sup_{\tt\in[\tau_2-3,\tau_2+1]}\E[\pmb\psi](\tau) +\int_{H_{R_1}(\tau',\tau'')}(R_1)^{-1}|\pmb\psi| |\pr^{\leq 1}\pmb\psi|+ \NNtlocal[\pmb \psi](\tau',\tau''),
\eea
where $c>0$ is a constant. This proves the Morawetz estimate \eqref{eq:Morawetznearinf:extendedsystem:Kerrpert:1:00}.

We next consider the improved Morawetz estimate \eqref{eq:Morawetznearinf:extendedRWsystem:Kerrpert:delta}. The proof is identical to the one of the improved Morawetz estimate in \cite[Lemma 3.10]{MaSz24} which consists in making the choice $(X,w)=(X_\de, w_\de)$ with $(X_\de, w_\de)$ given by \eqref{def:Xandw:improvedMorawetz:extendedRW:Kerrpert} to obtain 
\beaa
&&-\bigg(\QQ[\psi_{ij}]  \c {^{(X_{\de})}}\pi +w_{\de} \LL[\psi_{ij}]  -X_{\de}(V ) |\psi_{ij}|^2-\frac{1}{2}|\psi_{ij}|^2   \square_\g  w_{\de}\bigg)\nn\\
&\gtrsim &\de\bigg(\frac{|\pr_{\tt}\psi_{ij}|^2}{r^{1+\de}}
+\frac{|\pr_{r}\psi_{ij}|^2}{r^{1+\de}}\bigg) +\frac{|\nab\psi_{ij}|^2}{r}+\frac{|\psi_{ij}|^2}{r^3}
\eeaa
and 
\begin{equation*}
\begin{split}
&\sum_{i,j}\int_{\Sigma_{r\geq R_1}(\tau)}\big|{\PP}_\mu[\psi_{ij}](X_{\de}, w_{\de})N_{\Sigma_{\tau}}^{\mu}\big|\les \E_{r\geq R_1}[\pmb\psi](\tau),\\
&\sum_{i,j}\int_{\II_+(\tau',\tau'')}\big|{\PP}_\mu[\psi_{ij}](X_{\de}, w_{\de})N_{\II_+}^{\mu} \big|\les \F_{\II_+}[\pmb\psi](\tau',\tau'').
\end{split}
\end{equation*}
The improved Morawetz estimate \eqref{eq:Morawetznearinf:extendedRWsystem:Kerrpert:delta} then follows in the same manner as proving the above Morawetz estimate \eqref{eq:Morawetznearinf:extendedsystem:Kerrpert:1:00}, noticing that the cutoff function $\chi_{R_1}$ appearing in the definition \eqref{def:Xandw:improvedMorawetz:extendedRW:Kerrpert} of $(X_{\de}, w_{\de})$ induces an additional term $\M_{R_1-m,R_1}[\pmb\psi](\tau',\tau'')$ on the RHS of \eqref{eq:Morawetznearinf:extendedRWsystem:Kerrpert:delta} instead of a boundary term at $H_{R_1}$. This concludes the proof of Proposition \ref{prop:energymorawetznearinfinity:extendedRWsystem:kerrpert}.
\end{proof}

Next, we derive redshift estimates near the event horizon  for the coupled system of wave equations \eqref{eq:ScalarizedWaveeq:general:Kerrpert}.

\begin{lemma}[Redshift estimates for the coupled system of wave equations \eqref{eq:ScalarizedWaveeq:general:Kerrpert}]
\lab{lem:redshiftestimates:globallyextendedscalarizedwavesystem}
Let $\psi_{ij}$ be a solution to the system of wave equations \eqref{eq:ScalarizedWaveeq:general:Kerrpert}. 
 Then, under the assumptions of Theorem \ref{th:mainenergymorawetzmicrolocal} for the scalars $\psi_{ij}$, $F_{ij}$ and the spacetime $(\MM,\g)$, for any $1\leq \tau'<\tau''$ and any $\reg\leq 14$,  we have the following redshift estimates
\bea
\lab{eq:redshiftestimates:globallyextendedscalarizedwavesystem}
\nn\sum_{i,j}\EMF^{(\reg)}_{r\leq r_+(1+\dred)}[\psi_{ij}](\tau', \tau'')&\les&\sum_{i,j}\bigg( \E^{(\reg)}[\psi_{ij}](\tau')+\dred^{-1}\M_{r_+(1+\dred), r_+(1+2\dred)}^{(\reg)}[\psi_{ij}](\tau', \tau'')\\
&&+\int_{\MM_{r\leq r_+(1+2\dred)}(\tau', \tau'')}|\pr^{\leq \reg}F_{ij}|^2\bigg).
\eea
\end{lemma}

\begin{proof}
In view of  \eqref{eq:ScalarizedWaveeq:general:Kerrpert}, \eqref{hatSandV:generalwave:Kerrpert} and 
Lemma \ref{lem:estimatesforMialphaj:Kerrpert}, $\psi_{ij}$ satisfies
\beaa
{\square}_{\g}\psi_{ij} =\widehat{S}_K(\psi)_{ij} +\sum_{k,l}O(\ep)\pr\psi_{kl}+  \sum_{k,l}O(1)\psi_{kl}+F_{ij} ,\quad i,j=1,2,3, \quad\textrm{on}\quad\MM_{r\leq 3m}.
\eeaa
Also, we have in view of Lemma \ref{lemma:computationoftheMialphajinKerr}
\beaa
(M_K)_{i4}^k=O(\De)=O(|r-r_+|), \qquad (M_K)_{i3}^k=(M_K)_{ia}^k=O(1)\quad\textrm{on}\quad\MM_{r\leq 3m},
\eeaa
and hence
\beaa
\widehat{S}_K(\psi)_{ij} 
&=&2(M_K)_{i}^{k\a}\pr_\a(\psi_{kj}) +2(M_K)_{j}^{k\a}\pr_\a(\psi_{ik})+\frac{4ia\cos\th}{|q|^2} \pr_{\tt}(\psi_{ij})\nn\\
&=&O(|r-r_+|)e_{3}\psi + O(1)(e_{4}\psi, e_{a}\psi, \pr_\tau\psi), \quad \text{in} \,\,\MM_{r\leq 3m}.
\eeaa
Moreover, in view of \eqref{eq:actionofingoingprincipalnullframeonnormalizedcoordinates} \eqref{eq:actionofingoingprincipalnullframeonnormalizedcoordinates:bis}, together with the fact that $e_\a=e_\a(x^\b)\pr_\b$, we have
\beaa
e_{3}&=&-(1+O(\ep))\pr_{r} + O(1)\sum_{\a\neq r}\pr_{\a}, \quad \text{in} \,\,\MM_{r\leq 3m},\\
e_{\a}&=& O(1)\sum_{\b\neq r}\pr_{\b} +\big(O(\ep)+O(|r-r_+|)\big)\pr_{r}, \,\, \forall \a\neq 3, \quad \text{in} \,\,\MM_{r\leq 3m}.
\eeaa
The above implies that the scalars $\psi_{ij}$ satisfy,  in the redshift region $r\leq r_+(1+2\dred)$, 
\beaa
\square_{\g}\psi_{ij}=\sum_{k,l}\Big(O(|r-r_+|)\pr_r\psi_{kl} + O(1)(\pr_{\tt}\psi_{kl}, \pr_{x^a}\psi_{kl}, \psi_{kl})\Big)  +F_{ij}+\sum_{k,l}O(\ep)\pr_r\psi_{kl}.
\eeaa
This system of wave equations for the scalars $\psi_{ij}$ can be put into the form of \eqref{eq:Redshift-gen.scalarwaveeqs} and, applying Lemma \ref{lemma:redshiftestimatesscalarwaveeqs:general}, we have,  for any $1\leq \tau'<\tau'' $ and $\reg\leq 14$ and for $\ep$ suitably small,  
\beaa
\nn\sum_{i,j}\EMF^{(\reg)}_{r\leq r_+(1+\dred)}[\psi_{ij}](\tau', \tau'')&\les&\sum_{i,j}\bigg( \E^{(\reg)}[\psi_{ij}](\tau')+\dred^{-1}\M_{r_+(1+\dred), r_+(1+2\dred)}^{(\reg)}[\psi_{ij}](\tau', \tau'')\\
&&+\int_{\MM_{r\leq r_+(1+2\dred)}(\tau', \tau'')}|\pr^{\leq \reg}F_{ij}|^2\bigg)
\eeaa
as desired. 
\end{proof}


\subsection{Proof of Theorem \ref{th:mainenergymorawetzmicrolocal} on 
global energy-Morawetz estimates}
\lab{subsect:globalEMFestimate:scalarizeeqfromtensorial:0}


In this section, we show a global microlocal Morawetz estimate in Section \ref{subsubsect:globalmicrolocalMora:scalarizeeqfromtensorial:0} and an energy estimate in Section \ref{subsect:globalmicrolocalEner:scalarizeeqfromtensorial}, and we finally conclude the proof of Theorem \ref{th:mainenergymorawetzmicrolocal} in Section \ref{subsubsect:concludeglobalmicrolocalEner:scalarizeeqfromtensorial}.


\subsubsection{Global microlocal Morawetz estimate}
\lab{subsubsect:globalmicrolocalMora:scalarizeeqfromtensorial:0}


Recall that $A$ is the large constant appearing in the choice of multiplier $X$ in the region $\MM_{r_+(1+\dhor'),{\Rmic}}$, see \eqref{eq:decompositionofX:Mora}. Now, in the region $r\geq\Rmic$ we multiply the energy estimate \eqref{eq:energynearinf:extendedsystem:Kerrpert:1} by the large constant $A$ and add it to the Morawetz estimate \eqref{eq:Morawetznearinf:extendedsystem:Kerrpert:1:00}. We deduce the following energy-Morawetz estimate near infinity for solutions to \eqref{eq:ScalarizedWaveeq:general:Kerrpert}, for any $\tmic\leq \tau'<\tau''$,
\bea
\lab{eq:largerMoraesti:scalarizedwave:withFij}
&&c\MF_{r\geq {\Rmic}}[\pmb\psi ](\tau', \tau'')+c\E_{r\geq {\Rmic}}[\pmb\psi](\tau'')+\sum_{i,j}\textbf{BDR}[\psi_{ij} ]\vert_{H_{\Rmic}(\tau',\tau'')}\nn\\
&&
+\sum_{i,j}\int_{\MM_{{\Rmic},+\infty}(\tau',\tau'')}\Re\Big(F_{ij}\ov{ ({X_1} +A{\pr_{\tt}}+ w_1)\psi_{ij}}\Big)\nn\\
&&
+\sum_{i,j}\int_{\MM_{{\Rmic},+\infty}(\tau',\tau'')}\Re\bigg(F_{ij}\ov{ \Big( (X_1)_{\tau_2}(\psi)_{ij} - X_1 (\psi_{ij} ) + AT_{\tau_2}(\psi)_{ij} - A \pr_{\tau}(\psi_{ij})\Big)}\bigg)\nn\\
&\les&\E[\pmb\psi](\tau')+\sup_{\tau\in[\tmic,\tau_1+2]}\E[\pmb \psi](\tau)+\ep\sup_{\tau\in[\tau_2-3,\tau_2+1]}\E[\pmb\psi](\tau)+\ep\EM[\pmb\psi](\tau',\tau'')\nn\\
&&+ \NNtlocal[\pmb \psi](\tau',\tau'')
+\int_{H_{{\Rmic}}(\tau',\tau'')}(R_0)^{-1}|\pmb\psi| |\pr^{\leq 1}\pmb\psi|.
\eea
We now take $\tau'=\tmic$ and $\tau''\to +\infty$ in \eqref{eq:largerMoraesti:scalarizedwave:withFij} and sum the resulting estimate on $\MM_{r\geq\Rmic}(\Iti)$ with the estimate \eqref{eq:microlocalMora:middle:scalarizedwave:withFij} on $\MM_{r_+(1+\dhor'),\Rmic}(\Iti)$. Noticing, by the same argument as in \cite[Proposition 7.20]{MaSz24}, that the boundary terms $\textbf{BDR}^{-}_{r={\Rmic}}[\psi_{ij}](\Iti)$ in \eqref{eq:microlocalMora:middle:scalarizedwave:withFij}  and $\textbf{BDR}[\psi_{ij} ]\vert_{H_{\Rmic}(\Iti)}$ in the above inequality cancel to each other up to lower order terms controlled by $\int_{H_{\Rmic}(\Iti)}|\pmb\psi| |\pr^{\leq 1}\pmb\psi|$, we infer
\bea
\nn&& \sum_{i,j}\Bigg[\int_{\MM_{r\geq r_+(1+\dhor')}(\Iti)}\frac{\mu^2|\pr_r\psi_{ij}|^2}{r^2} +\int_{\Mntrap_{r\geq r_+(1+\dhor')}(\Iti)}\bigg(\frac{|\pr_\tau\psi_{ij}|^2}{r^2}+\frac{|\pr_{x^a}^{\leq 1}\psi_{ij}|^2}{r^3}\bigg)\\
&&+\int_{\MM_{r_+(1+\dhor'),10m}}\big(|\Opw(\sigma_{\trap})\psi_{ij}|^2+|\Opw(e)\psi_{ij}|^2\big)\Bigg]
+\F_{\II_+}[\pmb \psi](\Iti)\nn\\
&\les& \sup_{\tau\in[\tmic,\tau_1+2]}\E[\pmb \psi](\tau)+ (\ep+\dhor) \widetilde{\EM}[\pmb \psi]
+\NNtmora[\pmb \psi, \pmb F]+\NNtlocal[\pmb \psi](\Iti)\nn\\
&&+\frac{1}{\dhor^6}\A[\pmb \psi](\Iti) +(\ep+\dhor)\int_{\MM(\Iti)}|\pmb F|^2,
\eea
with $\NNtmora[\pmb \psi, \pmb F]$ defined as in \eqref{def:NNtMora:NNtEner:NNtaux:wavesystem:EMF} and $\NNtlocal[\pmb\psi](\Iti)$ defined as in \eqref{def:NNtlocalinr:NNtEner:NNtaux:wavesystem:EMF:proof},
where we have used \eqref{eq:lowerorderterms:controlled:scalarizedwavefromtensorial:proof:1}, and further Cauchy-Schwarz, to control the integral over $H_{\Rmic}(\Iti)$ which belongs to $\good^{(1)}[\psi]$
and used
\beaa
&&\sum_{i,j}\bigg|\int_{\MM_{{\Rmic},+\infty}(\tmic,+\infty)}\Re\bigg(F_{ij}\ov{ \Big( (X_1)_{\tau_2}(\psi)_{ij} - X_1 (\psi_{ij} ) + AT_{\tau_2}(\psi)_{ij} - A \pr_{\tau}(\psi_{ij})\Big)}\bigg)\bigg| \nn\\
&\les&  \ep \sup_{\tau\in\Iti}{\E}[\pmb\psi](\tau)
+\dhor^{-1}\int_{\MM(\Iti)}\frac{|\pmb\psi|^2}{r^3}+\dhor\int_{\MM_{{\Rmic},+\infty}(\Iti)}|\pmb F|^2
\eeaa
in view of
\beaa
M_{i\a}^k \psi_{kj} (\pr_{r}^{\text{BL}})^{\a} =\Ga_b \psi,\quad M_{i\a}^k \psi_{kj} (\pr_{\tt})^{\a}=(O(r^{-3}) + \Ga_b)\psi
\eeaa
which follows from Lemma \ref{lem:estimatesforMialphaj:Kerrpert}.

Finally, we state global Morawetz estimates for solutions to the system of wave equations \eqref{eq:ScalarizedWaveeq:general:Kerrpert}.

\begin{proposition}[Global Morawetz estimates for the system of wave equations \eqref{eq:ScalarizedWaveeq:general:Kerrpert}]
\lab{prop:globalMorawetz:systemofwaveequations:Kerrandpert}
Let $\psi_{ij}$ be a solution to the system of wave equations \eqref{eq:ScalarizedWaveeq:general:Kerrpert}.  
 Then, under the assumptions of Theorem \ref{th:mainenergymorawetzmicrolocal} for the scalars $\psi_{ij}$, $F_{ij}$ and the spacetime $(\MM,\g)$, we have, for $\ep+\dhor\ll \dred^4\ll 1$, the following global Morawetz estimates
\bea
\label{eq:nondegMora:generalscalarizedwave:Kerrandpert}
&&\dred^{4}\bigg(\widetilde{\MF}[\pmb \psi] + \sup_{\tau\in\Iti}\E_{r\leq r_+(1+\dred)}[\pmb \psi](\tau) \bigg)\nn\\
&\les&\sup_{\tau\in[\tmic,\tau_1+2]}\E[\pmb \psi](\tau)+ (\ep+\dhor)\sup_{\tt\in\Iti}{\E}[\pmb \psi](\tt)+\dhor^{-6}\A[\pmb \psi](\Iti) \nn\\
&&+\NNtmora[\pmb \psi, \pmb F]+\NNtaux[\pmb F]
+\NNtlocal[\pmb \psi](\Iti),
\eea
with $\NNtmora[\pmb \psi,\pmb  F]$ and  $\NNtaux[\pmb F]$ defined as in \eqref{def:NNtMora:NNtEner:NNtaux:wavesystem:EMF} and $\NNtlocal[\pmb\psi](\Iti)$ defined as in \eqref{def:NNtlocalinr:NNtEner:NNtaux:wavesystem:EMF:proof}.
\end{proposition}

\begin{proof}
We consider the redshift estimate \eqref{eq:redshiftestimates:globallyextendedscalarizedwavesystem} with $\tau'=\tmic$ and $\tau''\to+\infty$ and multiply it with $\dred^4$. We then sum the resulting estimate with  \eqref{eq:largerMoraesti:scalarizedwave:withFij} which yields for $\dred$ small enough the following global Morawetz estimate 
\beaa
&&\dred^{4}\bigg(\widetilde{\MF}[\pmb \psi] + \sup_{\tau\in\Iti}\E_{r\leq r_+(1+\dred)}[\pmb \psi](\tau) \bigg)\nn\\
&\les&\sup_{\tau\in[\tmic,\tau_1+2]}\E[\pmb \psi](\tau)+ (\ep+\dhor) \widetilde{\EM}[\pmb \psi]+\dhor^{-6}\A[\pmb \psi](\Iti)+\NNtmora[\pmb \psi, \pmb F]
+\NNtlocal[\pmb \psi](\Iti)\nn\\
&&{+(\ep+\dhor)\int_{\MM_{r_+(1+\dhor'), {\Rmic}}(\Iti)}|\pmb F|^2+\int_{\MM_{{\Rmic},+\infty}(\Iti)}|\pmb F|^2+\dred^4\int_{\MM_{r\leq r_+(1+2\dred)}(\Iti)}|\pmb F|^2}.
\eeaa
Then, requiring $\ep+\dhor\ll \dred^4$, we infer
\beaa
&&\dred^{4}\bigg(\widetilde{\MF}[\pmb \psi] + \sup_{\tau\in\Iti}\E_{r\leq r_+(1+\dred)}[\pmb \psi](\tau) \bigg)\nn\\
&\les&\sup_{\tau\in[\tmic,\tau_1+2]}\E[\pmb \psi](\tau)+ (\ep+\dhor)\sup_{\tt\in\Iti}{\E}[\pmb \psi](\tt)+\dhor^{-6}\A[\pmb \psi](\Iti) \nn\\
&&+\NNtmora[\pmb \psi, \pmb F]+\NNtaux[\pmb F]
+\NNtlocal[\pmb \psi](\Iti)
\eeaa
as desired. This concludes the proof of Proposition \ref{prop:globalMorawetz:systemofwaveequations:Kerrandpert}.
\end{proof}


\subsubsection{Energy estimate}
\lab{subsect:globalmicrolocalEner:scalarizeeqfromtensorial}


In order to show an energy estimate for the system of coupled scalar wave equations \eqref{eq:ScalarizedWaveeq:general:Kerrpert}, we shall recall, for the scalar wave equation $\square_{\g}\psi = F$, the following energy estimate which is a consequence of the one proven in \cite[Proposition 7.22]{MaSz24}.

\begin{lemma}[Conditional energy estimate]
\lab{lem:nondegEnerand:scalarfield:Kerrandpert:copy}
Assuming that the scalar fields $\psi$, $F$ and the metric $\g$ satisfy the same assumptions as in Lemma \ref{lem:conditionaldegenerateMorawetzflux:pertKerrrp1pdhorpR} and that $\psi$ solves $\square_{\g}\psi=F$,  we have the following conditional energy estimate
\bea\lab{thm:eq:nondeg:EnerandMora:Kerrandpert}
\sup_{\tau\in\Iti}\widehat{\E}[\psi](\tau) &\les& \E[\psi](\tmic)+\widetilde{\M}[\psi]+\NNtener[\psi, F]+\int_{\Mntrap(\Iti)}|F|^2\nn\\
&&+\ep\int_{\Mtrap(\Iti)}\tau^{-1-\dec}|F|^2 +\ep\int_{\Mtrap}\left|\Opw(\widetilde{S}^{-1,0}(\MM))F\right|^2,
\eea
where
\bea
\sup_{\tau\in\Iti}\widehat{\E}[\psi](\tau)&:=&\sup_{\tau\in\Iti}\E[\psi](\tau)
+\sum_{n=-1}^{\iota}\sup_{\tau\in\Reals}\E_{\trap}[\Opw(\Theta_n)\psi](\tau),\\
\E_{\trap}[\psi](\tau)&:=&\int_{\Sigma(\tau)\cap\Mtrap}r^{-2}|\dk^{\leq 1}\psi|^2,\\
\lab{def:NNtener:scalar:conditionalenergyesti}
\NNtener[\psi, F]&:=&\sup_{\tau\geq \tmic}\bigg|\int_{\Mntrap(\tmic, \tau)}{\Re\Big(F\ov{\pr_{\tau}\psi}\Big)}\bigg|
\nn\\
&&
+\sum_{i=-1}^{\iota}\sup_{\tt\in\Reals}\bigg|\int_{\Mtrap(\tmic, \tau)}\Re\Big(\ov{|q|^{-2}\Opw(\Theta_i)(\qs F)}V_i\Opw(\Theta_i)\psi\Big)\bigg|,
\eea
and where the symbols $\Theta_i\in\widetilde{S}^{0,0}(\MM)$ and first-order differential operators $V_i$, $i=-1, 0, 1,2,\ldots, \iota$, have been introduced in Section \ref{sec:relevantmixedsymbolsonMM}. 
Moreover, we have the following alternative conditional energy estimates
\bea\lab{eq:energyestimate:scalar:Kerrandpert:lastsect:widehatM:esti1}
\sup_{\tau\in\Iti}\widehat{\E}[\psi](\tau)
 \les \E[\psi](\tmic)+\widehat{\M}[\psi](\Iti)+\int_{\MM(\Iti)}|F|^2+\sup_{\tau\geq \tmic}\bigg|\int_{\Mntrap(\tmic, \tau)}{\Re\Big(F\ov{\pr_{\tau}\psi}\Big)}\bigg|
\eea
and
\bea\lab{eq:energyestimate:scalar:Kerrandpert:lastsect:widehatM}
\sup_{\tau\in\Iti}\widehat{\E}[\psi](\tau)
 \les  \E[\psi](\tmic)+\widehat{\M}[\psi](\Iti)+\int_{\MM(\Iti)}r^2|F|^2,
\eea
where 
\bea
\widehat{\M}[\psi](\Iti):={\M}[\psi](\Iti)+\int_{\Mtrap(\Iti)}|\pr\psi|^2.
\eea
\end{lemma}

\begin{proof}
The estimate \eqref{thm:eq:nondeg:EnerandMora:Kerrandpert} is a direct consequence of the estimates in the statement of  \cite[Inequality (7.145)]{MaSz24}.

To prove the estimate \eqref{eq:energyestimate:scalar:Kerrandpert:lastsect:widehatM:esti1}, we first notice that 
\beaa
&&\sum_{i=-1}^{\iota}\sup_{\tt\in\Reals}\bigg|\int_{\Mtrap(\tmic, \tau)}\Re\Big(\ov{|q|^{-2}\Opw(\Theta_i)(\qs F)}V_i\Opw(\Theta_i)\psi\Big)\bigg|\\
&\les& \left(\int_{\Mtrap(\Iti)}|F|^2\right)^{\frac{1}{2}}\left(\int_{\Mtrap(\Iti)}|\pr^{\leq 1}\psi|^2\right)^{\frac{1}{2}}.
\eeaa
Using this bound to control $\NNtener[\psi, F]$ and then plugging it in the RHS of \eqref{thm:eq:nondeg:EnerandMora:Kerrandpert}, and using also the fact that 
\beaa
\widetilde{\M}[\psi]+\int_{\Mtrap(\Iti)}|\pr^{\leq 1}\psi|^2\les \widehat{\M}[\psi](\Iti),
\eeaa
we obtain \eqref{eq:energyestimate:scalar:Kerrandpert:lastsect:widehatM:esti1}.

Finally, the estimate \eqref{eq:energyestimate:scalar:Kerrandpert:lastsect:widehatM} follows by applying Cauchy-Schwarz to the last term on the RHS of \eqref{eq:energyestimate:scalar:Kerrandpert:lastsect:widehatM:esti1}. This concludes the proof of Lemma \ref{lem:nondegEnerand:scalarfield:Kerrandpert:copy}.
\end{proof}

We next state a preliminary estimate which will be useful in generalizing the above statement from a scalar wave equation to the coupled system of scalar wave equations \eqref{eq:ScalarizedWaveeq:general:Kerrpert}.

\begin{lemma}
\lab{lem:recovinghalfderivativesinmicrolocalMorawetz}
Define
\bea
\lab{def:Errproduct:zeroordertimesfirstorder}
\Errprod[\psi]:=\sum_{i,j,k,l}\sup_{\tau'<\tau}\bigg|\int_{\Mtrap(\tau', \tau)}\Re\Big(\ov{\Opw(\widetilde{S}^{0,0}(\MM))\psi_{kl}} \Opw(\widetilde{S}^{1,1}(\MM))\psi_{ij}\Big)\bigg|.
\eea
Then, we have
\bea
\lab{eq:recovinghalfderivativesinmicrolocalMorawetz}
\Errprod[\psi]\les \widetilde{\M}[\pmb \psi]+\bigg(\int_{\Mtrap}|\pmb \psi|^2\bigg)^{\frac{1}{2}} \Big(\EM_{\trap}[\pmb \psi](\Reals)\Big)^{\frac{1}{2}}.
\eea
\end{lemma}

\begin{proof}
Notice that to prove \eqref{eq:recovinghalfderivativesinmicrolocalMorawetz}, it suffices to show 
\bea
\lab{def:Errproduct:zeroordertimesfirstorder:reduction}
\sup_{\tau\in\Reals}\bigg|\int_{\Mtrap(-\infty, \tau)}\Re\Big(\ov{E\psi_{kl}} S\psi_{ij}\Big)\bigg|
\les \widetilde{\M}[\pmb \psi]+\bigg(\int_{\Mtrap}|\pmb \psi|^2\bigg)^{\frac{1}{2}} \Big(\EM_{\trap}[\pmb \psi](\Reals)\Big)^{\frac{1}{2}}
\eea
for any $S\in \Opw(\widetilde{S}^{1,1}(\MM))$, $E\in \Opw(\widetilde{S}^{0,0}(\MM))$ and $i,j,k,l$.

Given $S\in \Opw(\widetilde{S}^{1,1}(\MM))$,  we can decompose it into
$$
S=S_0\pr_r +\pr_rS_0+ S_1, \quad S_0\in \Opw(\widetilde{S}^{0,0}(\MM)),\quad S_1\in \Opw(\widetilde{S}^{1,0}(\MM)).
$$
Hence, for any $\tau'\in \Reals$ and for any $E\in \Opw(\widetilde{S}^{0,0}(\MM))$, we have
\bea\lab{eq:Mtrapfrom-inftytotau':zeroordertimesfirstorder:proof:00}
\nn&&\bigg|\int_{\Mtrap(-\infty, \tau')}\Re\Big(\ov{E\psi_{kl}}S\psi_{ij}\Big)
\bigg|\nn\\
\nn&\les &
\bigg|\int_{\Mtrap(-\infty, \tau')}\Re\Big(\ov{E\psi_{kl}}(S_0\pr_r+\pr_rS_0)\psi_{ij}\Big)\bigg|
+\bigg|\int_{\Mtrap(-\infty, \tau')}\Re\Big(\ov{E\psi_{kl}}S_1\psi_{ij}\Big)\bigg|\nn\\
&\les&\int_{\Mtrap}|\Opw(\widetilde{S}^{0,0}(\MM))\psi_{kl}||\pr^{\leq 1}_r\psi_{ij}|+\bigg|\int_{\Mtrap(-\infty, \tau')}\Re\Big(\ov{E\psi_{kl}}S_1\psi_{ij}\Big)\bigg|\nn\\
&\les& \M[\pmb \psi](\Reals)+\bigg|\int_{\Mtrap(-\infty, \tau')}\Re\Big(\ov{E\psi_{kl}}S_1\psi_{ij}\Big)\bigg|.
\eea

Next, we estimate the last term on the RHS of \eqref{eq:Mtrapfrom-inftytotau':zeroordertimesfirstorder:proof:00}. To this end, we introduce the smooth cut-off functions $\chi_{\tau',j}=\chi_{\tau',j}(\tau)$, $j=0,1,2$, such that
\begin{equation}
\lab{def:variouscutofffunctionsfortau'}
\begin{split}
&\textrm{supp}(\chi_{\tau',0})\subset(-\infty, \tau'+1), \quad \chi_{\tau',0}=1\,\,\,\textrm{on}\,\,(-\infty, \tau'), \quad 0\leq\chi_{\tau',j}\leq 1, \,\,j=0,1,\\
&\textrm{supp}(\chi_{\tau',1})\subset(\tau'-1, \tau'+2), \quad \chi_{\tau',1}=1\,\,\,\textrm{on}\,\,(\tau', \tau'+1).
\end{split}
\end{equation}
This allows us to estimate the last term on the RHS of \eqref{eq:Mtrapfrom-inftytotau':zeroordertimesfirstorder:proof:00} as follows
\bea
\lab{eq:Mtrapfrom-inftytotau':zeroordertimesfirstorder:proof}
&&\bigg|\int_{\Mtrap(-\infty, \tau')}\Re\Big(\ov{E\psi_{kl}}S_1\psi_{ij}\Big)\bigg|\nn\\
&\les&\bigg|\int_{\Mtrap}\Re\Big(\chi_{\tau',0}\ov{E\psi_{kl}}S_1\psi_{ij}\Big)\bigg|
+\int_{\Mtrap}\chi_{\tau',1}|E{\psi_{kl}}| |S_1\psi_{ij}|\nn\\
&\les&\bigg|\int_{\Mtrap}\Re\Big(\chi_{\tau',0}\ov{E\psi_{kl}}S_1\psi_{ij}\Big)\bigg|
+\int_{\Mtrap}|E{\psi_{kl}}| |S_1(\chi_{\tau',1}\psi_{ij})|+\int_{\Mtrap}|E{\psi_{kl}}| |[S_1, \chi_{\tau',1}]\psi_{ij}|\nn\\
&\les&\bigg|\int_{\Mtrap}\Re\Big(\chi_{\tau',0}\ov{E\psi_{kl}}S_1\psi_{ij}\Big)\bigg|
+\bigg(\int_{\Mtrap}|\pmb \psi|^2\bigg)^{\frac{1}{2}} \Big(\int_{\Mtrap}|\pr^{\leq 1}(\chi_{\tau',1}\psi_{ij})|^2+
\M_{\trap}[\pmb \psi]\Big)^{\frac{1}{2}}\nn\\
&\les&\bigg|\int_{\Mtrap}\Re\Big(\chi_{\tau',0}\ov{E\psi_{kl}}S_1\psi_{ij}\Big)\bigg|
+\bigg(\int_{\Mtrap}|\pmb \psi|^2\bigg)^{\frac{1}{2}} \Big(\EM_{\trap}[\pmb \psi]\Big)^{\frac{1}{2}},
\eea
where we used the size of the support of $\chi_{\tau',1}$ given by \eqref{def:variouscutofffunctionsfortau'} in the last inequality.

Next, we control the first term on the RHS of \eqref{eq:Mtrapfrom-inftytotau':zeroordertimesfirstorder:proof}. To this end, recall from \eqref{eq:definitionofmicrolocalMorawetznormwidetildeM} that 
\bea\lab{eq:Mtrapfrom-inftytotau':zeroordertimesfirstorder:proof:pogacar4}
\int_{\Mtrap}|\Opw(\sigma_{\trap})\psi|^2\les \widetilde{\M}[\psi],
\eea
where, as introduced in Section  \ref{sec:relevantmixedsymbolsonMM}, $\sigma_{\trap}=(r-r_{\trap})\upsilon$ with the mixed symbols $r_{\trap}\in\widetilde{S}^{0,0}(\MM)$ and $\upsilon\in\widetilde{S}^{1,0}(\MM)$. Also, recall that $\pr_r(r_{\trap})=0$ which implies $\pr_r(r-r_{\trap})=1$ and thus, for any scalar function $\psi$,
\beaa
\pr_r(\Opw(r-r_{\trap})\psi) &=& \Opw(r-r_{\trap})\pr_r\psi+\Opw(\pr_r(r-r_{\trap}))\psi\\
&=& \Opw(r-r_{\trap})\pr_r\psi+\psi.
\eeaa
Then, using a smooth cut-off function $\chi_0(r)$ supported in $(r_+(1-\dhor), 11m)$ with $\chi_0=1$ on the support of $\Mtrap$, we may estimate the first term on the RHS of \eqref{eq:Mtrapfrom-inftytotau':zeroordertimesfirstorder:proof} as follows
\beaa
\bigg|\int_{\Mtrap}\Re\Big(\chi_{\tau',0}\ov{E\psi_{kl}}S_1\psi_{ij}\Big)\bigg|
&\les& \bigg|\int_{\MM}\Re\Big(\chi_0(r)\chi_{\tau',0}\ov{E\psi_{kl}}S_1\psi_{ij}\Big)\bigg|+\M[\pmb \psi](\Reals)\\
&\les& \bigg|\int_{\MM}\Re\Big(\chi_0(r)\chi_{\tau',0}\ov{\Opw(r-r_{\trap})\pr_r(E\psi_{kl})}S_1\psi_{ij}\Big)\bigg|\\
&+&\bigg|\int_{\MM}\Re\Big(\chi_0(r)\chi_{\tau',0}\ov{\pr_r(\Opw(r-r_{\trap})E\psi_{kl})}S_1\psi_{ij}\Big)\bigg|+\M[\pmb \psi](\Reals),
\eeaa
and hence, after integration by parts in $r$ for the second term, taking the adjoint of $\Opw(r-r_{\trap})$ for the first term and the adjoint of $S_1$ for the second term, and using the fact that commutators generate lower order terms, we infer
\beaa
\bigg|\int_{\Mtrap}\Re\Big(\chi_{\tau',0}\ov{E\psi_{kl}}S_1\psi_{ij}\Big)\bigg| &\les& \int_{\Mtrap}\Big(|\pr_r^{\leq 1}\psi_{ij}|^2+|\pr^{\leq 1}_r\psi_{kl}|^2+|\Opw(\sigma_{\trap})\psi_{ij}|^2\\
&&\qquad\qquad\qquad\qquad\qquad\qquad +|\Opw(\sigma_{\trap})\psi_{kl}|^2\Big)+\M[\pmb \psi](\Reals)\\
&\les& \int_{\Mtrap}\Big(|\Opw(\sigma_{\trap})\psi_{ij}|^2+|\Opw(\sigma_{\trap})\psi_{kl}|^2\Big)+\M[\pmb \psi](\Reals).
\eeaa
In view of \eqref{eq:Mtrapfrom-inftytotau':zeroordertimesfirstorder:proof:pogacar4}, this yields
\beaa
\bigg|\int_{\Mtrap}\Re\Big(\chi_{\tau',0}\ov{E\psi_{kl}}S_1\psi_{ij}\Big)\bigg| &\les&\widetilde{\M}[\pmb \psi],
\eeaa
which, together with \eqref{eq:Mtrapfrom-inftytotau':zeroordertimesfirstorder:proof:00}
and \eqref{eq:Mtrapfrom-inftytotau':zeroordertimesfirstorder:proof}, proves \eqref{def:Errproduct:zeroordertimesfirstorder:reduction}, and hence \eqref{eq:recovinghalfderivativesinmicrolocalMorawetz}. This concludes the proof of Lemma \ref{lem:recovinghalfderivativesinmicrolocalMorawetz}.
\end{proof}

We are now ready to derive an energy estimate for the system of scalar wave equations \eqref{eq:ScalarizedWaveeq:general:Kerrpert}.

\begin{proposition}[Conditional energy estimate for the system of scalar wave equations \eqref{eq:ScalarizedWaveeq:general:Kerrpert}]
\lab{lem:nondegEnerand:scalarizedwave:Kerrandpert}
Let $\psi_{ij}$ be a solution to the system of wave equations \eqref{eq:ScalarizedWaveeq:general:Kerrpert}, and let $\NNtener[\pmb \psi, \pmb F]$, $\NNtaux[\pmb F]$, $\Errdefect[\pmb\psi]$ and $\NNtlocal[\pmb \psi](\Iti)$ be given as in \eqref{def:NNtMora:NNtEner:NNtaux:wavesystem:EMF}, 
\eqref{def:Errdefectofpsi} and \eqref{def:NNtlocalinr:NNtEner:NNtaux:wavesystem:EMF:proof}, respectively. Then, under the assumptions of Theorem \ref{th:mainenergymorawetzmicrolocal} for the scalars $\psi_{ij}$, $F_{ij}$ and the spacetime $(\MM,\g)$, we have the following conditional energy estimate 
\bea\lab{thm:eq:nondeg:EnerandMora:scalarized:Kerrandpert}
\sup_{\tau\in\Iti}\widehat{\E}[\pmb \psi](\tau) &\les& \sup_{\tau\in[\tmic,\tau_1+2]}\E[\pmb \psi](\tau) +\widetilde{\M}[\pmb \psi]
+\Errdefect[\pmb\psi]\nn\\
&&+\NNtener[\pmb \psi, \pmb F] +\NNtaux[\pmb F]+ \NNtlocal[\pmb \psi](\Iti)+\int_{\Mntrap(\Iti)}|\pmb F|^2.
\eea
Moreover, we have the following alternative conditional energy estimates
\bea\lab{eq:energyestimate:scalarized:Kerrandpert:lastsect:widehatM:esti1}
\sup_{\tau\in\Iti}\widehat{\E}[\pmb \psi](\tau)
& \les & \sup_{\tau\in[\tmic,\tau_1+2]}\E[\pmb \psi](\tau)+\widehat{\M}[\pmb \psi](\Iti)+\Errdefect[\pmb\psi]+ \NNtlocal[\pmb \psi](\Iti)\nn\\
 &&+\int_{\MM(\Iti)}|\pmb F|^2+\sup_{{\tau\geq \tmic}}\bigg|\int_{\Mntrap(\tmic, \tau)}\sum_{i,j}{\Re\Big(F_{ij}\ov{\pr_{\tau}\psi_{ij}}\Big)}\bigg|
\eea
and
\bea\lab{eq:energyestimate:scalarized:Kerrandpert:widehatM}
\sup_{\tau\in\Iti}\widehat{\E}[\pmb \psi](\tau)& \les&  \sup_{\tau\in[\tmic,\tau_1+2]}\E[\pmb \psi](\tau)+\widehat{\M}[\pmb \psi](\Iti)+\Errdefect[\pmb\psi]\nn\\
&&+ \NNtlocal[\pmb \psi](\Iti)+\int_{\MM(\Iti)}r^2|\pmb F|^2,
\eea
where we have defined 
\bea\lab{eq:definitionofsuptauinItiofwidehatEofpmbpsioftau}
\sup_{\tau\in\Iti}\widehat{\E}[\pmb \psi](\tau):=\sup_{\tau\in\Iti}\E[\pmb \psi](\tau)
+\sum_{n=-1}^{\iota}\sum_{i,j}\sup_{\tau\in\Reals}\E_{\trap}[\Opw(\Theta_n)\psi_{ij}](\tau)
\eea
with the symbols $\Theta_{n}\in\widetilde{S}^{0,0}(\MM)$, $n=-1,0,1,2,\ldots, \iota$, given as in Section \ref{sec:relevantmixedsymbolsonMM}.
\end{proposition}

\begin{proof}
We first control $\sup_{\tau\in\Iti}\widehat{\E}_{r\geq 11m}[\pmb \psi](\tau)$. To this end, we apply the energy estimate \eqref{eq:energynearinf:extendedsystem:Kerrpert:2:cutoff} with $r_1=11m$ to deduce 
\bea
&&\E_{r\geq 11m}[\pmb\psi](\tau'')
+\F_{\II_+}[\pmb\psi](\tau',\tau'')\nn\\
&\les&  \sup_{\tau\in[\tmic,\tau_1+2]}\E[\pmb \psi](\tau)+\E_{r\geq 10m}[\pmb\psi](\tau')   + \NNtlocal[\pmb \psi](\tau',\tau'')+ \ep\EM[\pmb\psi](\tau',\tau'')\nn\\
&&
+\int_{\MM_{10m,11m}(\tau',\tau'')} \big(|\pr^{\leq 1}\pmb\psi|^2+|\pmb F|^2\big)+\bigg|\int_{\MM_{r\geq 11m}(\tau',\tau'')}\sum_{i,j}\Re\big(F_{ij}\ov{T_{\tau_2}(\psi)_{ij}}\big)\bigg|.
\eea
By taking $\tau'=\tmic$ and  the supremum of $\tau''\in\Reals$, we infer
\bea
\lab{eq:energyestiinrgeqRregion:extendedRWsystem:Kerrpert:proof}
\nn\EF_{r\geq 11m}[\pmb\psi](\Iti) &\les&  \sup_{\tau\in[\tmic,\tau_1+2]}\E[\pmb \psi](\tau)+\M[\pmb \psi](\Iti) +  \NNtlocal[\pmb \psi](\Iti)+\int_{\Mntrap}|\pmb F|^2\\
&&+\ep\sup_{\tau\in\Iti}\E[\pmb\psi](\tau)+\sum_{i,j}\sup_{\tau\geq\tmic}\bigg|\int_{\Mntrap(\tmic, \tau)}{\Re\Big(F_{ij}\ov{\pr_{\tau}\psi_{ij}}\Big)}\bigg|.
\eea

It remains to estimate the energy in the region $r\leq 11m$. Recall from \eqref{eq:ScalarizedWaveeq:general:Kerrpert:rewriteinsect8} that the scalars $\psi_{ij}$ satisfy 
\bea
\lab{eq:waveofpsiij:inscalarform}
\square_{\g}(\psi_{ij})=F_{ij} +\sum_{p=1}^4G_{p,ij},
\eea
where 
\bsub
\lab{def:Gnij:energyesti}
\bea
G_{1,ij}&=& S_K(\psi)_{ij} ,\\
 G_{2,ij}&=&\frac{4ia\cos\th}{|q|^2} \pr_{\tt}\psi_{ij},\\
  G_{3,ij}&=& \chi_{\tau_1, \tau_2}\big({S}(\psi)_{ij}-{S}_K(\psi)_{ij}\big),\\
    G_{4,ij}&=&\chi_{\tau_1, \tau_2}(\widehat{Q}\psi)_{ij}+(1-\chi_{\tau_1, \tau_2})\big((\widehat{Q}_K\psi)_{ij}+f_{D_0}\psi_{ij}\big)+D_0|q|^{-2}\psi_{ij}.
\eea
\esub
This implies that the scalars $\chi_{12m}(r)\psi_{ij}$ satisfy
\bea
\lab{eq:waveofpsiij:inscalarform:cutoff}
\square_{\g}(\chi_{12m}(r)\psi_{ij})&=&\chi_{12m}(r)F_{ij} +\chi_{12m}(r)\sum_{p=1}^4G_{p,ij} \nn\\
&&+2\g^{r\b}\pr_{r}(\chi_{12m}(r))\pr_{\b}(\psi_{ij})+\square_{\g}(\chi_{12m}(r))\psi_{ij},
\eea
where $\chi_{12m}(r)$ is a smooth cut-off function in $r$ satisfying $\chi_{12m}(r)=1$ on $r\leq 11m$ and $\chi_{12m}(r)=0$ on $r\geq 12m$. Consequently, applying Lemma \ref{lem:nondegEnerand:scalarfield:Kerrandpert:copy} to the above wave equations \eqref{eq:waveofpsiij:inscalarform:cutoff} and summing over $i,j=1,2,3$, we deduce, using also the definition of $\NNtener[\pmb \psi, \pmb F]$ in \eqref{def:NNtMora:NNtEner:NNtaux:wavesystem:EMF}, 
\bea\lab{eq:nondeg:Ener:scalarapplytoscalarized:firststep}
\sup_{\tau\in\Iti}\widehat{\E}_{r\leq 11m}[\pmb \psi](\tau) &\les& \E[\pmb\psi](\tmic)+\widetilde{\M}[\pmb \psi]+\NNtener[\pmb \psi, \pmb F] +\NNtaux[\pmb F]\nn\\
&&+\ep\sup_{\tau\in\Iti}\E[\pmb\psi](\tau)+\int_{\Mntrap(\Iti)}|\pmb F|^2 +\sum_{p=1}^4\Err_{G_p},
\eea
where for $p=1,2,3,4$,
\beaa
\Err_{G_p}:=\sum_{n=-1}^{\iota}\sum_{i,j}\sup_{\tt\in\Reals}\bigg|\int_{\Mtrap(\tmic, \tau)}\Re\Big(\ov{|q|^{-2}\Opw(\Theta_n)(\qs G_{p,ij})}V_n\Opw(\Theta_n)\psi_{ij}\Big)\bigg|,
\eeaa
where we have used the following estimate that follows from the definition of $\NNtaux[\pmb F]$ in \eqref{def:NNtMora:NNtEner:NNtaux:wavesystem:EMF}
\beaa
\sum_{i,j}\Bigg(\ep\int_{\Mtrap(\Iti)}\tau^{-1-\dec}|F_{ij}|^2 +\ep\int_{\Mtrap}\left|\Opw(\widetilde{S}^{-1,0}(\MM))F_{ij}\right|^2\Bigg)\les \NNtaux[\pmb F],
\eeaa
and where we have applied Cauchy-Schwarz to control the integrals over $11m\leq r\leq 12m$ by 
\beaa
\M[\pmb\psi](\Iti) + \int_{\Mntrap(\Iti)}|\pmb F|^2  \les \widetilde{\M}[\pmb \psi]+ \int_{\Mntrap(\Iti)}|\pmb F|^2.
\eeaa

Next, we control the terms $\Err_{G_p}$, $p=1,2,3,4$, appearing on the RHS of \eqref{eq:nondeg:Ener:scalarapplytoscalarized:firststep}. First, in view of the expression of $G_{4,ij}$ in \eqref{def:Gnij:energyesti}, and using Lemma \ref{lem:recovinghalfderivativesinmicrolocalMorawetz}, we have
\bea
\lab{eq:NNtener:lowerorder:ErrG3:fromscalartoscalarized:proof}
\Err_{G_4}
&\les& \widetilde{\M}[\pmb \psi]+\bigg(\int_{\Mtrap(\Iti)}|\pmb \psi|^2\bigg)^{\frac{1}{2}} \Big(\EM_{\trap}[\pmb \psi](\Iti)\Big)^{\frac{1}{2}}.
\eea

Next, in view of \eqref{SandV} and the assumption \eqref{eq:assumptionsonregulartripletinperturbationsofKerr:0}, we have
${S}(\psi)_{ij}-{S}_K(\psi)_{ij}=\sum_{k,l}\widecheck{\Ga}\pr\psi_{kl}$ in $\Mtrap$. To control the term $\Err_{G_3}$, we introduce the smooth cut-off functions $\chi_{\tmic, \tau,j}=\chi_{\tmic, \tau,j}(\tau)$, $j=0,1$, satisfying
\begin{align*}
&\textrm{supp}(\chi_{\tmic,\tau,0})\subset(\tmic, \tau), \quad \chi_{\tmic,\tau,0}=1\,\,\,\textrm{on}\,\,(\tmic+1, \tau-1), \quad 0\leq\chi_{\tmic,\tau,j}\leq 1, \,\,j=0,1,\\
&\textrm{supp}(\chi_{\tmic,\tau,1})\subset(\tmic-1, \tmic+2)\cup(\tau-2, \tau+1), \quad \chi_{\tmic,\tau,1}=1\,\,\,\textrm{on}\,\,(\tmic, \tmic+1)\cup(\tau-1, \tau).
\end{align*}
Using the properties of the cut-offs $\chi_{\tmic, \tau,j}$, $j=0,1$, we infer
\beaa
\Err_{G_3}&\les&\sum_{n=-1}^{\iota}\sum_{i,j,k,l}\sup_{\tt\in\Reals}\bigg|\int_{\Mtrap(\tmic, \tau)}\Re\Big(\ov{|q|^{-2}\Opw(\Theta_n)(\qs \widecheck{\Ga}\pr\psi_{kl})}V_n\Opw(\Theta_n)\psi_{ij}\Big)\bigg|\nn\\
&\les&\sum_{n=-1}^{\iota}\sum_{i,j,k,l}\sup_{\tt\in\Reals}\bigg|\int_{\Mtrap}\chi_{\tmic,\tau,0}\Re\Big(\ov{|q|^{-2}\Opw(\Theta_n)(\qs \widecheck{\Ga}\pr\psi_{kl})}V_n\Opw(\Theta_n)\psi_{ij}\Big)\bigg|\nn\\
&&+\sum_{n=-1}^{\iota}\sum_{i,j,k,l}\sup_{\tt\in\Reals}\int_{\Mtrap}\chi_{\tmic,\tau,1}\Big||q|^{-2}\Opw(\Theta_n)(\qs \widecheck{\Ga}\pr\psi_{kl})\Big| \Big|V_n\Opw(\Theta_n)\psi_{ij}\Big|.
\eeaa
Taking the adjoint of $\Opw(\Th_n)$ in the before to last line, and using  
Proposition \ref{prop:PDO:MM:Weylquan:mixedoperators} and Lemmas \ref{lemma:actionmixedsymbolsSobolevspaces:MM} and \ref{lem:gpert:MMtrap}, we deduce
\bea
\lab{eq:NNtener:lowerorder:ErrG3perturbation:fromscalartoscalarized:proof}
\Err_{G_3}
&\les&\sum_{i,j,k,l}\bigg|\int_{\Mtrap}\widecheck{\Ga}\Re\Big(\ov{\pr\psi_{kl}}\Opw(\widetilde{S}^{1,0}(\MM))\psi_{ij}\Big)\bigg|\nn\\
&&+\sum_{n=-1}^{\iota}\sum_{i,j,k,l}\int_{\Mtrap}\Big||q|^{-2}\Opw(\Theta_n)(\qs \widecheck{\Ga}\pr\psi_{kl})\Big| \Big|\Opw(\widetilde{S}^{0,0}(\MM))\psi_{ij}\Big|\nn\\
&&+\sum_{n=-1}^{\iota}\sum_{i,j,k,l}\sup_{\tt\in\Reals}\int_{\Mtrap}\Big||q|^{-2}\Opw(\Theta_n)(\qs \widecheck{\Ga}\pr\psi_{kl})\Big| \Big|V_n\Opw(\Theta_n)(\chi_{\tmic,\tau,1}\psi_{ij})\Big|\nn\\
&\les&\ep \EM[\pmb\psi](\Iti).
\eea

Next, we estimate the term $\Err_{G_2}$. For any fixed $i,j,n$ and $\tau'\in\Reals$, we have
\beaa
\Err_{G_2} &=&\bigg|\int_{\Mtrap(\tmic, \tau')}\Re\Big(\ov{|q|^{-2}\Opw(\Theta_n)(\qs G_{2,ij})}V_n\Opw(\Theta_n)\psi_{ij}\Big)\bigg|\nn\\
&=&\bigg|\int_{\Mtrap(\tmic, \tau')}\Re\bigg(\ov{|q|^{-2}\Opw(\Theta_n)\big(4ia\cos\th\pr_{\tau}\psi_{ij}\big)}V_n\Opw(\Theta_n)\psi_{ij}\bigg)\bigg|\nn\\
&\les &\bigg|\int_{\Mtrap(\tmic, \tau')}\Re\Big(\ov{\Opw(\widetilde{S}^{0,0}(\MM))\psi_{kl}} \Opw(\widetilde{S}^{1,0}(\MM))\psi_{ij}\Big)\bigg|\nn\\
&&+\bigg|\int_{\Mtrap(\tmic, \tau')}\Re\bigg(\ov{\frac{4ia\cos\th}{\qs}\pr_{\tau}\Opw(\Theta_n)\psi_{ij}}V_n\Opw(\Theta_n)\psi_{ij}\bigg)\bigg|\nn\\
&\les &\widetilde{\M}[\pmb \psi]+\bigg(\int_{\Mtrap(\Iti)}|\pmb \psi|^2\bigg)^{\frac{1}{2}} \Big(\EM_{\trap}[\pmb \psi](\Iti)\Big)^{\frac{1}{2}}\\
&&+\bigg|\int_{\Mtrap(\tmic, \tau')}\Re\bigg(\ov{\frac{4ia\cos\th}{\qs}\pr_{\tau}\varphi_{ij}}V_n\varphi_{ij}\bigg)\bigg|,
\eeaa
where in the last step we have used Lemma \ref{lem:recovinghalfderivativesinmicrolocalMorawetz} and defined new scalars $\varphi_{ij}:=\Opw(\Theta_n)\psi_{ij}$. Then, in order to control the last term on the RHS, we integrate the differential identity \eqref{eq:integrationbypartsforitimesfirstordertimesfirstorder:general} on $\Mtrap(-\infty, \tau')$ respectively with $f=\frac{2a\cos\th}{\qs}f_0$, $x^\a=x^\b=\tau$ and $f=\frac{2a\cos\th}{\qs}f_0d_i(r)$, $x^\a=\tau$, $x^\b=\tphi$, where we recall that $V_i=\pr_\tau+d_i(r)\pr_{\tphi}$,  $i=-1,0,1,\cdots,\iota,$ and $f_0$ is given as in \eqref{eq:spacetimevolumeformusingisochorecoordinates}. Since we have, in view of \eqref{eq:spacetimevolumeformusingisochorecoordinates} and the definition of $f$, 
\beaa
f_0^{-1}\pr_\tau(f)=\dk\Ga_g, \qquad f_0^{-1}\pr_{\tphi}(f)=\dk\Ga_g,
\eeaa
we infer
\beaa
\Err_{G_2} &\les &\widetilde{\M}[\pmb \psi]+\bigg(\int_{\Mtrap(\Iti)}|\pmb \psi|^2\bigg)^{\frac{1}{2}} \Big(\EM_{\trap}[\pmb \psi](\Iti)\Big)^{\frac{1}{2}}\\
&&+\bigg(\int_{\Sigma(\tau')\cap\Mtrap}|\varphi_{ij}|^2\bigg)^{\frac{1}{2}}\Big(\E_{\trap}[\varphi_{ij}](\tau')\Big)^{\frac{1}{2}}\\
&&+\bigg(\int_{\Sigma(\tmic)\cap\Mtrap}|\varphi_{ij}|^2\bigg)^{\frac{1}{2}}\Big(\E_{\trap}[\varphi_{ij}](\tmic)\Big)^{\frac{1}{2}}+\ep\sup_{\tau\in\Reals}\E_{\trap}[\varphi_{ij}](\tau).
\eeaa
Using the definition \eqref{eq:definitionofsuptauinItiofwidehatEofpmbpsioftau} and the fact that $\varphi_{ij}=\Opw(\Theta_n)\psi_{ij}$, we deduce
\bea
\lab{eq:NNtener:itimesfirstorder:fromscalartoscalarized:proof}
\Err_{G_2} &\les &\widetilde{\M}[\pmb \psi]+\bigg(\int_{\Mtrap}|\pmb \psi|^2\bigg)^{\frac{1}{2}} \Big(\EM_{\trap}[\pmb \psi](\Iti)\Big)^{\frac{1}{2}}
+\ep\sup_{\tau\in\Iti}\widehat{\E}[\pmb\psi](\tau)\nn\\
&&+\bigg(\sup_{\tau\in\Iti}\widehat{\E}[\pmb\psi](\tau)\bigg)^{\frac{1}{2}}\bigg(\int_{\Sigma(\tau')\cap\Mtrap}|\varphi_{ij}|^2
+\int_{\Sigma(\tmic)\cap\Mtrap}|\varphi_{ij}|^2\bigg)^{\frac{1}{2}}.
\eea
Next, we estimate the last term on the RHS of \eqref{eq:NNtener:itimesfirstorder:fromscalartoscalarized:proof}. To this end, we rely on $\chi_{\tau'}=\chi_{\tau'}(\tau)$ such that
\beaa
&&\textrm{supp}(\chi_{\tau'})\subset(\tau'-2, \tau'+2), \quad \chi_{\tau'}=1\,\,\,\textrm{on}\,\,(\tau'-1, \tau'+1), \quad 0\leq\chi_{\tau'}\leq 1
\eeaa
and derive the following trace estimate for any $\tau'\in\Reals$
\bea
\lab{eq:energyofvarphiij:PDOS:trapregion:proof}
\int_{\Sigma(\tau')\cap\Mtrap}|\varphi_{ij}|^2&=&\int_{\Sigma(\tau')\cap\Mtrap}\chi_{\tau'}|\varphi_{ij}|^2\les  \int_{\Mtrap}|\pr_{\tau}\big(\chi_{\tau'}|\varphi_{ij}|^2\big)|\nn\\
&\les&
\int_{\Mtrap}|\psi_{ij}|^2 +\int_{\Mtrap}\Big|2\chi_{\tau'}\Re\Big(\varphi_{ij}\ov{\pr_{\tau}\varphi_{ij}}\Big)\Big|\nn\\
&\les&\int_{\Mtrap(\Iti)}|\pmb\psi|^2 +\bigg(\int_{\Mtrap(\Iti)}|\pmb \psi|^2\bigg)^{\frac{1}{2}} \Big(\sup_{\tau\in\Iti}\widehat{\E}[\psi](\tau)\Big)^{\frac{1}{2}}.
\eea
Using \eqref{eq:energyofvarphiij:PDOS:trapregion:proof} to control the integrals on $\Sigma(\tau')\cap\Mtrap$ and $\Sigma(\tmic)\cap\Mtrap$ in \eqref{eq:NNtener:itimesfirstorder:fromscalartoscalarized:proof}, we infer 
\bea
\lab{eq:NNtener:lowerorder:ErrG2:fromscalartoscalarized:proof}
\Err_{G_2}
&\les&\widetilde{\M}[\pmb \psi]+\bigg(\int_{\Mtrap(\Iti)}|\pmb \psi|^2\bigg)^{\frac{1}{2}} \Big(\EM_{\trap}[\pmb \psi](\Iti)\Big)^{\frac{1}{2}}
+\ep\sup_{\tau\in\Iti}\widehat{\E}[\pmb\psi](\tau)\nn\\
&&+\bigg(\sup_{\tau\in\Iti}\widehat{\E}[\pmb\psi](\tau)\bigg)^{\frac{3}{4}}\bigg(\int_{\Mtrap(\Iti)}|\pmb \psi|^2\bigg)^{\frac{1}{4}}.
\eea

Next, we estimate the term $\Err_{G_1}$. We have
\bea
\lab{eq:NNtener:firstordertimesfirstorder:fromscalartoscalarized:proof}
\Err_{G_1}&=&\sum_n\sum_{i,j}\bigg|\int_{\Mtrap(\tmic, \tau')}\Re\Big(\ov{|q|^{-2}\Opw(\Theta_n)(\qs G_{1,ij})}V_n\Opw(\Theta_n)\psi_{ij}\Big)\bigg|\nn\\
&=&\sum_n\sum_{i,j}\bigg|\int_{\Mtrap(\tmic, \tau')}\Re\Big(\ov{|q|^{-2}\Opw(\Theta_n)\big(\qs S_K(\psi)_{ij}\big)}V_n\Opw(\Theta_n)\psi_{ij}\Big)\bigg|\nn\\
&\les &\sum_n\sum_{i,j}\bigg(\bigg|\int_{\Mtrap(\tmic, \tau')}\Re\Big(\ov{\Opw(\widetilde{S}^{0,0}(\MM))\psi_{kl}} \Opw(\widetilde{S}^{1,0})\psi_{ij}\Big)\bigg|\nn\\
&+&\bigg|\int_{\Mtrap(\tmic, \tau')}\Re\Big(\ov{\big(2(M_K)_{i}^{k\a}\pr_\a(\Opw(\Theta_n)\psi_{kj}) +2(M_K)_{j}^{k\a}\pr_\a(\Opw(\Theta_n)\psi_{ik})\big)}\nn\\
&&\qquad\qquad\qquad\qquad\qquad\qquad\qquad\qquad\qquad\qquad\qquad\qquad\times V_n\Opw(\Theta_n)\psi_{ij}\Big)\bigg|\bigg)\nn\\
&\les &\widetilde{\M}[\pmb \psi]+\bigg(\int_{\Mtrap(\Iti)}|\pmb \psi|^2\bigg)^{\frac{1}{2}} \Big(\EM_{\trap}[\pmb \psi](\Iti)\Big)^{\frac{1}{2}}\nn\\
&+&\sum_n\sum_{i,j}\bigg|\int_{\Mtrap(\tmic, \tau')}\Re\Big(\ov{\big(2(M_K)_{i}^{k\a}\pr_\a(\varphi_{kj}) +2(M_K)_{j}^{k\a}\pr_\a(\varphi_{ik})\big)}V_n\varphi_{ij}\Big)\bigg|,
\eea
where in the last step we have used Lemma \ref{lem:recovinghalfderivativesinmicrolocalMorawetz} as well as  the definition of the scalars $\varphi_{ij}=\Opw(\Theta_n)\psi_{ij}$. It remains to control the last term on the RHS of \eqref{eq:NNtener:firstordertimesfirstorder:fromscalartoscalarized:proof}, and it suffices to estimate the term $|\int_{\Mtrap(\tmic, \tau')}\Re(\ov{2(M_K)_{i}^{k\a}\pr_\a(\varphi_{kj}) }V_n\varphi_{ij})|$ since the other term can be  controlled analogously. Using the decomposition \eqref{eq:decompositionofMikalpha:symandantisym} which reads
\beaa
(M_K)_{i}^{k\a}=(M_{K,S})_{i}^{k\a} +(M_{K,A})_{i}^{k\a}=-\frac{1}{2}\pr^{\a}(x^i x^k)+(M_{K,A})_{i}^{k\a},
\eeaa 
with $(M_{K,S})_{i}^{k\a}$ and $(M_{K,A})_{i}^{k\a} $ denoting respectively the symmetric and antisymmetric parts of $(M_K)_{i}^{k\a}$ w.r.t. $(i,k)$,  we obtain
\bea
\lab{eq:firstorderproductofThetanpsi}
&&\sum_n\sum_{i,j}\bigg|\int_{\Mtrap(\tmic, \tau')}\Re\Big(\ov{2(M_K)_{i}^{k\a}\pr_\a(\varphi_{kj}) }V_n\varphi_{ij}\Big)\bigg|\nn\\
&\les&\sum_n\sum_{i,j}\bigg(\bigg|\int_{\Mtrap(\tmic, \tau')}\Re\Big(\ov{\pr^{\a}(x^i x^k)\pr_\a(\varphi_{kj}) }V_n\varphi_{ij}\Big)\bigg|\nn\\
&&+\bigg|\int_{\Mtrap(\tmic, \tau')}\Re\Big(\ov{2(M_{K,A})_{i}^{k\a}\pr_\a(\varphi_{kj}) }V_n\varphi_{ij}\Big)\bigg|\bigg).
\eea
We first deal with the second term on the RHS of \eqref{eq:firstorderproductofThetanpsi}. To this end, we integrate the differential identity \eqref{eq:integrationbypartsforitimesfirstordertimesfirstorder:general:caseantisymmatrix} with the choice $A_i^k=(M_{K,A})_{i}^{k\a}f_0$ or $(M_{K,A})_{i}^{k\a}f_0d_i(r)$ and obtain\footnote{Since we have $(M_K)_{i}^{kr}=0$ in $\Mtrap$ in view of Lemma \ref{lemma:computationoftheMialphajinKerr}, it follows $(M_{K,A})_{i}^{kr}=0$ in $\Mtrap$ as well, so integrating the the differential identity \eqref{eq:integrationbypartsforitimesfirstordertimesfirstorder:general:caseantisymmatrix} does not generate boundary terms on the parts $\{r=r_+(1+2\dbl)\}$ and $\{r=10m\}$ of $\pr\Mtrap$.}
\beaa
&&\sum_n\sum_{i,j}\bigg|\int_{\Mtrap(\tmic, \tau')}\Re\Big(\ov{2(M_{K,A})_{i}^{k\a}\pr_\a(\varphi_{kj}) }V_n\varphi_{ij}\Big)\bigg|\nn\\
&\les& \sum_n\sum_{i,j}\bigg|\int_{\Mtrap(\tmic, \tau')}\Re\Big(\pr(A_i^k)\ov{\varphi_{kj}}\pr\varphi_{ij}\Big)\bigg|+\sum_n\sum_{i,j}\bigg|\int_{\Mtrap(\tmic, \tau')}\Re\Big(\pr(A_i^k)\ov{\pr\varphi_{kj}}\varphi_{ij}\Big)\bigg|\\
&&+ \sum_{i,j,k}\Bigg(\int_{\Sigma(\tau')\cap\Mtrap}
+\int_{\Sigma(\tmic)\cap\Mtrap}\bigg)\big(|\varphi_{kj}||\pr\varphi_{ij}|+|\pr\varphi_{kj}||\varphi_{ij}|\big)+\int_{\Mtrap(\Iti)}|\pmb \psi|^2\nn\\
&\les& \widetilde{\M}[\pmb\psi] +\bigg(\int_{\Mtrap(\Iti)}|\pmb \psi|^2\bigg)^{\frac{1}{2}} \Big(\EM_{\trap}[\pmb \psi](\Iti)\Big)^{\frac{1}{2}}\\
&&+\sum_{i,j,k,l}\left(\int_{\Sigma(\tau')\cap\Mtrap}|\varphi_{ij}|^2+\int_{\Sigma(\tmic)\cap\Mtrap}|\varphi_{ij}|^2\right)^{\frac{1}{2}}\left(\sup_{\tau\in\Iti}\widehat{\E}[\pmb{\psi}](\tau)\right)^{\frac{1}{2}}
\eeaa
where we used \eqref{eq:recovinghalfderivativesinmicrolocalMorawetz} in the last estimate. Together with \eqref{eq:energyofvarphiij:PDOS:trapregion:proof} and the fact that $\varphi_{ij}=\Opw(\Theta_n)\psi_{ij}$, we infer
\bea
\lab{eq:estiofzeroordertimesfirstorder:varphi:Kerrantisym}
&&\sum_n\sum_{i,j}\bigg|\int_{\Mtrap(\tmic, \tau')}\Re\Big(\ov{2(M_{K,A})_{i}^{k\a}\pr_\a(\varphi_{kj}) }V_n\varphi_{ij}\Big)\bigg|\nn\\
&\les& \widetilde{\M}[\pmb\psi] +\bigg(\int_{\Mtrap(\Iti)}|\pmb \psi|^2\bigg)^{\frac{1}{2}} \Big(\EM_{\trap}[\pmb \psi](\Iti)\Big)^{\frac{1}{2}}\nn\\
&&+ \bigg(\int_{\Mtrap(\Iti)}|\pmb \psi|^2\bigg)^{\frac{1}{4}}\Big(\sup_{\tau\in\Iti}\widehat{\E}[\pmb{\psi}](\tau)\Big)^{\frac{3}{4}}.
\eea
Plugging the estimate \eqref{eq:estiofzeroordertimesfirstorder:varphi:Kerrantisym} into \eqref{eq:firstorderproductofThetanpsi}, we infer
\bea
\lab{eq:firstorderproduct:varphiij:reduced}
&&\sum_n\sum_{i,j}\bigg|\int_{\Mtrap(\tmic, \tau')}\Re\Big(\ov{2(M_K)_{i}^{k\a}\pr_\a(\varphi_{kj}) }V_n\varphi_{ij}\Big)\bigg|\nn\\
&\les& \widetilde{\M}[\pmb\psi] +\bigg(\int_{\Mtrap(\Iti)}|\pmb \psi|^2\bigg)^{\frac{1}{2}} \Big(\EM_{\trap}[\pmb \psi](\Iti)\Big)^{\frac{1}{2}}\nn\\
&&+ \bigg(\int_{\Mtrap(\Iti)}|\pmb \psi|^2\bigg)^{\frac{1}{4}}\Big(\sup_{\tau\in\Iti}\widehat{\E}[\pmb{\psi}](\tau)\Big)^{\frac{3}{4}}\nn\\
&&+\sum_n\sum_{i,j}\bigg|\int_{\Mtrap(\tmic, \tau')}\Re\Big(\ov{\pr^{\a}(x^i x^k)\pr_\a(\Opw(\Th_n)\psi_{kj}) }V_n\Opw(\Th_n)\psi_{ij}\Big)\bigg|.
\eea
Next, we focus on the last term on the RHS of \eqref{eq:firstorderproduct:varphiij:reduced}. Since $\pr^{\a}(x^i x^k)= x^i \pr^{\a} (x^k) +x^k \pr^{\a} (x^i)$, we may commute $x^i$ with $V_n \Opw(\Th_n)$ and commute $x^k$ with $\pr_{\a}\Opw(\Th_n)$. Relying on Lemma \ref{lem:recovinghalfderivativesinmicrolocalMorawetz} to control the terms involving the commutators $[x^i, V_n \Opw(\Th_n)]$,  $[x^k, \pr_{\a}\Opw(\Th_n)]$, $[V_n, \Opw(\Th_n)]$ and $[\pr_{\a}, \Opw(\Th_n)]$ which are of the type $\Errprod[\psi]$, we infer
\bea
\lab{eq:firstorderproduct:symmetricpart:varphiij:reduced}
&&\sum_n\sum_{i,j}\bigg|\int_{\Mtrap(\tmic, \tau')}\Re\Big(\ov{\pr^{\a}(x^i x^k)\pr_\a(\Opw(\Th_n)\psi_{kj}) }V_n\Opw(\Th_n)\psi_{ij}\Big)\bigg|\nn\\
&\les&\widetilde{\M}[\pmb \psi]+\bigg(\int_{\Mtrap(\Iti)}|\pmb \psi|^2\bigg)^{\frac{1}{2}} \Big(\EM_{\trap(\Iti)}[\pmb \psi](\Iti)\Big)^{\frac{1}{2}}\nn\\
&& +\Err_{G_1,\textrm{main}, 1}+\Err_{G_1,\textrm{main}, 2},
\eea
where we have defined
\beaa
\Err_{G_1,\textrm{main}, 1}&:=&\sum_n\sum_{i,j,k}\bigg|\int_{\Mtrap(\tmic, \tau')}\Re\Big(\ov{\pr^{\a}(x^k)\Opw(\Th_n)\pr_\a\psi_{kj} }\Opw(\Th_n)V_n(x^i\psi_{ij})\Big)\bigg|,\\
\Err_{G_1,\textrm{main}, 2}&:=&\sum_n\sum_{i,j,k}\bigg|\int_{\Mtrap(\tmic, \tau')}\Re\Big(\ov{\pr^{\a}(x^i)\Opw(\Th_n) \pr_\a(x^k\psi_{kj}) }\Opw(\Th_n)V_n\psi_{ij}\Big)\bigg|.
\eeaa

Next, let $\chi_0=\chi_0(\tau)$ be a smooth cut-off function such that 
\beaa
\chi_0(\tau)=1\quad\textrm{on}\quad [\tmic,\tau_1+1]\cup[\tau_2-2, \tau_2], \qquad \textrm{supp}(\chi_0)\subset [\tmic-1,\tau_1+2]\cup[\tau_2-3, \tau_2+1].
\eeaa
Then, in view of  the support property \eqref{eq:widetildepsiij:tau2-4totau2} of $x^i\psi_{ij}$, we have $x^i\psi_{ij}=\chi_0(\tau)x^i\psi_{ij}$, and hence
\beaa
\Err_{G_1,\textrm{main}, 1}&=&\sum_n\sum_{i,j,k}\bigg|\int_{\Mtrap(\tmic, \tau')}\Re\Big(\ov{\pr^{\a}(x^k)\Opw(\Th_n)\pr_\a\psi_{kj} }\Opw(\Th_n)V_n(\chi_0 x^i\psi_{ij})\Big)\bigg|\\
&\les& \sum_n\sum_{i,j,k}\bigg|\int_{\Mtrap(\tmic, \tau')}\Re\Big(\ov{\pr^{\a}(x^k)\Opw(\Th_n)\pr_\a(\chi_0\psi_{kj})}\Opw(\Th_n)V_n(x^i\psi_{ij})\Big)\bigg|\\
&&+\Errprod[\psi]\\
&\les& \left(\sum_{ij}\int_{\Mtrap}|\pr(\chi_0\psi_{ij})|^2\right)^{\frac{1}{2}}\big(\Errdefect[\pmb\psi]\big)^{\frac{1}{2}}+\Errprod[\psi]\\
&\les& \Errprod[\psi]+\Big(\sup_{\tau\in\Iti}\E_{\trap}[\pmb\psi](\tau)\Big)^{\frac{1}{2}}\big(\Errdefect[\pmb\psi]\big)^{\frac{1}{2}},
\eeaa
where, to go from the first to the second line, we have used the fact that all the terms involving a commutator with $\chi_0$ are of the type $\Errprod[\psi]$. By the same argument, the other term $\Err_{G_1,\textrm{main}, 2}$ satisfies the same bound, 
which together with \eqref{eq:firstorderproduct:symmetricpart:varphiij:reduced} and \eqref{eq:firstorderproduct:varphiij:reduced} implies
\bea
\lab{eq:firstorderproduct:varphiij:reduced:final}
&&\sum_n\sum_{i,j}\bigg|\int_{\Mtrap(\tmic, \tau')}\Re\Big(\ov{2(M_K)_{i}^{k\a}\pr_\a(\varphi_{kj}) }V_n\varphi_{ij}\Big)\bigg|\nn\\
&\les& \widetilde{\M}[\pmb\psi] +\bigg(\Errdefect[\pmb\psi]+\int_{\Mtrap(\Iti)}|\pmb \psi|^2\bigg)^{\frac{1}{2}} \Big(\EM_{\trap}[\pmb \psi](\Iti)\Big)^{\frac{1}{2}}\nn\\
&&+ \bigg(\int_{\Mtrap(\Iti)}|\pmb \psi|^2\bigg)^{\frac{1}{4}}\Big(\sup_{\tau\in\Iti}\widehat{\E}[\pmb\psi](\tau)\Big)^{\frac{3}{4}},
\eea
where we have also used the estimate \eqref{eq:recovinghalfderivativesinmicrolocalMorawetz} to control $\Errprod[\psi]$. Plugging this estimate into \eqref{eq:NNtener:firstordertimesfirstorder:fromscalartoscalarized:proof}, we infer
\bea
\lab{eq:NNtener:lowerorder:ErrG1:fromscalartoscalarized:proof}
\Err_{G_1} &\les& \widetilde{\M}[\pmb\psi] +\bigg(\Errdefect[\pmb\psi]+\int_{\Mtrap(\Iti)}|\pmb \psi|^2\bigg)^{\frac{1}{2}} \Big(\EM_{\trap}[\pmb \psi](\Iti)\Big)^{\frac{1}{2}}\nn\\
&&+ \bigg(\int_{\Mtrap(\Iti)}|\pmb \psi|^2\bigg)^{\frac{1}{4}}\Big(\sup_{\tau\in\Iti}\widehat{\E}[\pmb\psi](\tau)\Big)^{\frac{3}{4}}.
\eea

Plugging the above estimates \eqref{eq:NNtener:lowerorder:ErrG3:fromscalartoscalarized:proof}, \eqref{eq:NNtener:lowerorder:ErrG3perturbation:fromscalartoscalarized:proof}, \eqref{eq:NNtener:lowerorder:ErrG2:fromscalartoscalarized:proof} and \eqref{eq:NNtener:lowerorder:ErrG1:fromscalartoscalarized:proof} for $\{\Err_{G_j}\}_{j=1,2,3,4}$ into the inequality \eqref{eq:nondeg:Ener:scalarapplytoscalarized:firststep}, we infer
\beaa
\sup_{\tau\in\Iti}\widehat{\E}_{r\leq 11m}[\pmb \psi](\tau)&\les& \E[\pmb\psi](\tmic)+\widetilde{\M}[\pmb \psi]+\NNtener[\pmb \psi, \pmb F] +\NNtaux[\pmb F]+\int_{\Mntrap(\Iti)}|\pmb F|^2\\
&&+\bigg(\Errdefect[\pmb\psi]+\int_{\Mtrap(\Iti)}|\pmb \psi|^2\bigg)^{\frac{1}{2}} \Big(\EM_{\trap}[\pmb \psi](\Iti)\Big)^{\frac{1}{2}}\\
&& +\ep\sup_{\tau\in\Iti}\widehat{\E}[\pmb\psi](\tau)+\bigg(\sup_{\tau\in\Iti}\widehat{\E}[\pmb\psi](\tau)\bigg)^{\frac{3}{4}}\bigg(\int_{\Mtrap}|\pmb \psi|^2\bigg)^{\frac{1}{4}}.
\eeaa
Combining this estimate for $\sup_{\tau\in\Iti}\widehat{\E}_{r\leq 11m}[\pmb \psi](\tau) $ with the estimate \eqref{eq:energyestiinrgeqRregion:extendedRWsystem:Kerrpert:proof}
for $\EF_{r\geq 11m}[\pmb \psi](\Iti)$, and taking $\ep$ small enough, we then conclude the desired estimate \eqref{thm:eq:nondeg:EnerandMora:scalarized:Kerrandpert}.

Similarly to the proof of Lemma \ref{lem:nondegEnerand:scalarfield:Kerrandpert:copy}, the estimate \eqref{eq:energyestimate:scalarized:Kerrandpert:lastsect:widehatM:esti1} follows from \eqref{thm:eq:nondeg:EnerandMora:scalarized:Kerrandpert} by applying Cauchy-Schwarz to the last term on the RHS of the expression \eqref{def:NNtener:scalar:conditionalenergyesti} of $\NNtener[\pmb \psi, \pmb F]$, and the estimate \eqref{eq:energyestimate:scalarized:Kerrandpert:widehatM} follows from applying Cauchy-Schwarz to each term in $\NNtener[\pmb \psi, \pmb F]$.
\end{proof}


\subsubsection{Proof of Theorem \ref{th:mainenergymorawetzmicrolocal}}
\lab{subsubsect:concludeglobalmicrolocalEner:scalarizeeqfromtensorial}


Combining the above  energy estimate \eqref{thm:eq:nondeg:EnerandMora:scalarized:Kerrandpert} with the Morawetz estimate \eqref{eq:nondegMora:generalscalarizedwave:Kerrandpert}, taking $\dhor$ and $\ep$ suitably small, and applying Cauchy-Schwarz to the term $\NNtlocal[\pmb \psi](\Iti)$ defined as in \eqref{def:NNtlocalinr:NNtEner:NNtaux:wavesystem:EMF:proof}, we deduce the desired energy-Morawetz estimate 
\eqref{th:eq:mainenergymorawetzmicrolocal:tensorialwave:scalarized:eachpsisp} and hence complete the proof of Theorem \ref{th:mainenergymorawetzmicrolocal}.


\section{EMF estimates for $\phis{p}$ and $\psis{p}$ at zeroth-order}
\lab{sec:proofofth:main:intermediary:0order}


Let $(\MM, \g)$ satisfy the assumptions of Sections \ref{subsect:assumps:perturbednullframe}, \ref{subsubsect:assumps:perturbedmetric} and \ref{sec:regulartripletinperturbationsofKerr}. By a slight abuse of notation, we denote throughout this section the metric $\g_{\chi_{\tau_1, \tau_2}}$ in \eqref{eq:extendedmetricgchitau1tau2} by $\g$. Also, let $\{\pmb\phi_s^{(p)}\}_{s=\pm 2, p=0,1,2}$ be a solution to the tensorial Teukolsky wave/transport systems \eqref{eq:TensorialTeuSysandlinearterms:rescaleRHScontaine2:general:Kerrperturbation}  \eqref{def:TensorialTeuScalars:wavesystem:Kerrperturbation}  in perturbations of Kerr, and let $\{\psis{p}\}_{s=\pm 2, p=0,1,2}$ be a solution to the following system of wave equations\footnote{In practice, $\{\psis{p}\}_{s=\pm 2, p=0,1,2}$ is a solution to \eqref{eq:waveeqwidetildepsi1} \eqref{def:tildef} so that $\widehat{F}_{total,s,ij}^{(p)}$ is  given by 
\bea\lab{eq:definitionofwidehatFtotalsij}
\widehat{F}_{total,s,ij}^{(p)}:=\chi_{\tau_1, \tau_2}^{(1)}N^{(p)}_{W,s,ij}+\underline{F}^{(p)}_{s,ij} +\breve{F}^{(p)}_{s,ij},
\eea
with $\underline{F}^{(p)}_{s,ij}$ and $\breve{F}^{(p)}_{s,ij}$ given respectively by \eqref{def:tildef0} and  \eqref{def:breveFpsij}. In this section, we do not specify the form of $\widehat{F}_{total,s,ij}^{(p)}$ appearing in \eqref{eq:definitionofwidehatFtotalsij} as we need a version that is stable under commutation w.r.t. higher order derivatives in view of our applications in Section \ref{sec:proofofth:main:intermediary}.}
\begin{align}
\lab{eq:waveeqwidetildepsi1:gtilde:Teu}
\bigg({\square}_{{\g}} -\frac{4-2\de_{p0}}{|q|^{2}}\bigg)\psi^{(p)}_{s,ij}
={}&\chi_{\tau_1, \tau_2}\big( \widehat{S}(\psi^{(p)}_s)_{ij} +(\widehat{Q}\psi^{(p)}_s)_{ij}\big) +(1-\chi_{\tau_1, \tau_2})\big(\widehat{S}_K(\psi^{(p)}_s)_{ij} +(\widehat{Q}_K\psi^{(p)}_s)_{ij} \big)\nn\\
&+(1-\chi_{\tau_1, \tau_2})f_p\psi^{(p)}_{s,ij}+\chi_{\tau_1, \tau_2}^{(1)}L^{(p)}_{s,ij}+\widehat{F}_{total,s,ij}^{(p)} 
\end{align}
on $\MM$ satisfying \eqref{eq:causlityrelationsforwidetildepsi1}. 

The goal of this section is to initiate the proof of the EMF estimates for unweighted derivatives of $\phis{p}$ and $\psis{p}$ stated in Theorem \ref{th:main:intermediary} by first proving the EMF estimates for $\phis{p}$ and $\psis{p}$ at zeroth-order stated below in Theorem \ref{thm:EMF:systemofTeuscalarized:order0:final}. Relying on Theorem \ref{thm:EMF:systemofTeuscalarized:order0:final}, Theorem \ref{th:main:intermediary} will then be proved in Section \ref{sec:proofofth:main:intermediary}.


\subsection{Statement of the EMF estimates for $\phis{p}$ and $\psis{p}$ at zeroth-order}
\lab{subsect:globalEMFestimate:loworder:mainstatement}


In this section, we state Theorem \ref{thm:EMF:systemofTeuscalarized:order0:final} on EMF estimates for $\phis{p}$ and $\psis{p}$ at zeroth-order. To this end, we start by introducing several notations.

For the sake of convenience, we define, for any $\tau'<\tau''$,\footnote{Note that $\A[\pmb\psi](\tau',\tau'')$ generalizes to an arbitrary interval $(\tau',\tau'')$ the notation $\A[\pmb \psi](\Iti)$ in \eqref{def:EM-1norms:Reals}.}
\bsub
\lab{def:AandAonorms:tau1tau2}
\begin{align}
\A[\pmb\psi](\tau',\tau''):={}&\sum_{i,j=1}^3\left(\int_{\MM(\tau',\tau'')}\frac{|\psi_{ij}|^2}{r^3} +\sup_{\tau\in(\tau',\tau'')}\int_{\Sigma_{\tau}}\frac{|\psi_{ij}|^2}{r^2}
+\int_{\II_+(\tau',\tau'')}\frac{|\psi_{ij}|^2}{r^2}\right), \\
\Ao[\pmb\psi](\tau',\tau''):={}&\sum_{i,j=1}^3\int_{\MM(\tau',\tau'')}\frac{|\psi_{ij}|^2}{r^3}.
\end{align}
\esub
We also define, for any $\tau'<\tau''$,
\bsub
\lab{def:AandAonorms:tau1tau2:phisandpsis}
\begin{align}
\A[\pmb\psi_s](\tau',\tau''):=&\sum_{i,j=1}^3\sum_{p=0}^2\bigg(\int_{\MM(\tau',\tau'')}\frac{|\psiss{ij}{p}|^2}{r^3}  +\sup_{\tau\in(\tau',\tau'')}\int_{\Sigma_{\tau}}\frac{|\psiss{ij}{p}|^2}{r^2}+\int_{\II_+(\tau',\tau'')}\frac{|\psiss{ij}{p}|^2}{r^2}\bigg), \\
\Ao[\pmb\psi_s](\tau',\tau''):=&\sum_{i,j=1}^3\sum_{p=0}^2\int_{\MM(\tau',\tau'')}\frac{|\psiss{ij}{p}|^2}{r^3},
\end{align}
\esub
and the corresponding quantities $\A[\pmb\phi_s](\tau',\tau'')$ and $\Ao[\pmb\phi_s](\tau',\tau'')$ with $\psiss{ij}{p}$ replaced by $\phiss{ij}{p}$ in the above formulas.

Next, we define, for any $\tau'<\tau''$, $\reg\in\mathbb{N}$ and $\de\in [0,1]$, the following EMF norm 
\bea
\lab{def:widehatEMFdenorm}
\bsplit
\widehat{\EMF}_{\de}^{(\reg)}[ \pmb\psi](\tau',\tau''):={}& \EF^{(\reg)}[ \pmb\psi](\tau',\tau'')+\widehat{\M}_{\de}^{(\reg)}[ \pmb\psi](\tau',\tau''),\\
\widehat{\M}_{\de}^{(\reg)}[ \pmb\psi](\tau',\tau''):={}& \M_{\de}^{(\reg)}[ \pmb\psi](\tau',\tau'')+\int_{\Mtrap(\tau',\tau'')} |\pr^{\leq \reg+1}\pmb\psi|^2\\
\simeq{}&\int_{\MM(\tau',\tau'')} \bigg(\frac{|\nab_{\pr_r} \dk^{\leq \reg}\pmb\psi|^2}{r^{1+\de}}
+\frac{|\nab_{\pr_\tt} \dk^{\leq \reg}\pmb\psi|^2}{r^{1+\de}}
+\frac{|(r\nab)^{\leq 1}\dk^{\leq \reg}\pmb\psi|^2}{r^{3}}\bigg),
\end{split}
\eea
in which the spacetime integrand is non-degenerate w.r.t. all derivatives in the trapping region. Also, denote for convenience  $\widehat{\EMF}_{1}^{(\reg)}[\pmb\psi](\tau',\tau'')$ by $\widehat{\EMF}^{(\reg)}[\pmb\psi](\tau',\tau'')$.

We define the following early time energy norm for $s=\pm 2$ that will be useful later 
\begin{align}
\lab{eq:definitionofinitialenergyofphiandpsi:IEterm}
\IE{\pmb\phi_s, \pmb\psi_s}:={}&\sum_{p=0,1,2}\left(\E[\phis{p}](\tau_1)+\sup_{\tau\in[\tmic, \tau_1+2]}\E[\psis{p}](\tau)\right).
\end{align}
Next, we define the following EMF norm for $s=\pm 2$: 
\begin{align}
\lab{def:EMFtotalps:pm2}
\EMFtotalp{s}:={}&\sum_{p=0,1,2}\widetilde{\EMF}[\psis{p}]
+\sum_{p=0,1}\widehat{\M}_{\de}[\psis{p}](\tau_1+1, \tau_2-3)
\nn\\
&+\sum_{p=0,1}\widehat{\EMF}_{\de}[\phis{p}](\tau_1,\tau_2)+ {\EMF}_{\de}[\phis{2}](\tau_1,\tau_2).
\end{align}
Also, we define, for $\pmb\psi\in\sk_k(\mathbb{C})$ and $\pmb H\in\sk_k(\mathbb{C})$, $k=1,2$, 
\bea
\lab{eq:defofNNhat'}
\widehat{\mathcal{N}}'[\pmb\psi, \H](\tau_1, \tau_2)
&:=& \sup_{\tau_1< \tau'<\tau''< \tau_2}\bigg|\int_{\Mntrap(\tau', \tau'')}\big(1+O(r^{-\de})\big)\nab_{\pr_\tau}\pmb\psi\c\ov{\H}\bigg| \nn\\
&&+\int_{\Mntrap(\tau_1, \tau_2)}r^{-1}|\dk^{\leq 1}\pmb\psi||\H|+\int_{\MM(\tau_1, \tau_2)}|\H|^2,
\eea
where the coefficient $1+O(r^{-\de})$ appearing on the RHS of \eqref{eq:defofNNhat'} is in practice either equal to $1$ or to the smooth function $f_\de=f_\de(r)$ introduced in \eqref{def:Xandw:improvedMorawetz:extendedRW:Kerrpert}. Finally, we define, for $s=\pm 2$,
\bea\lab{expression:NNttotalps:pm2}
\NNttotalp{s}
&:=&\sum_{p=0,1,2}\widetilde{\mathcal{N}}[\psis{p},\widehat{\F}_{total,s}^{(p)}]+\sum_{p=0,1,2}\widehat{\mathcal{N}}'[\phis{p}, \N_{W, s}^{(p)}](\tau_1, \tau_2)\nn\\
&&+\sum_{p=0,1}\widehat{\mathcal{N}}'[ \phis{p}, r^{-2}\N_{T, s}^{(p)}](\tau_1, \tau_2)+\sum_{p=0,1}\int_{\MM_{r\leq 12m}(\tau_1,\tau_2)}|\pr^{\leq 1}\N_{T, s}^{(p)}|^2\nn\\
&&+\sum_{p=0,1}\int_{\MM(\tau_1,\tau_2)} \Big(r^{-1+\de}|\N^{(p)}_{W,s}|+r^{-2+\de}|\pr\N^{(p)}_{T,s}|\Big)|\phis{p}|,
\eea
where $\widetilde{\mathcal{N}}[\c, \c]$ and $\widehat{\NN}'[\c, \c]$ are defined in \eqref{def:NNtintermsofNNtMora:NNtEner:NNtaux:wavesystem:EMF}-\eqref{def:NNtMora:NNtEner:NNtaux:wavesystem:EMF} and \eqref{eq:defofNNhat'}. Also,  we define $\EMFtotalhps{s}{\pr^{\leq \reg}}$ and $\NNttotalph{s}{\pr^{\leq \reg}}$ accordingly by making the replacements 
\beaa
(\phis{p}, \psis{p}, \widehat{\F}_{total,s}^{(p)}, \N_{W, s}^{(p)}, \N_{T, s}^{(p)}) \to  (\pr^{\leq \reg}\phis{p}, \pr^{\leq \reg}\psis{p}, \pr^{\leq \reg}\widehat{\F}_{total,s}^{(p)}, \pr^{\leq \reg}\N_{W, s}^{(p)}, \pr^{\leq \reg}\N_{T, s}^{(p)}) .
\eeaa

We are now ready to state our main EMF estimates for $\phis{p}$ and $\psis{p}$ at zeroth-order.

\begin{theorem}[EMF estimates for  $\phis{p}$ and $\psis{p}$ at zeroth-order]
\lab{thm:EMF:systemofTeuscalarized:order0:final}
Let $(\MM, \g)$ satisfy the assumptions of Sections \ref{subsect:assumps:perturbednullframe}, \ref{subsubsect:assumps:perturbedmetric} and \ref{sec:regulartripletinperturbationsofKerr}. Let $\{\pmb\phi_s^{(p)}\}_{s=\pm 2, p=0,1,2}$ be a solution to the tensorial Teukolsky wave/transport systems \eqref{eq:TensorialTeuSysandlinearterms:rescaleRHScontaine2:general:Kerrperturbation}  \eqref{def:TensorialTeuScalars:wavesystem:Kerrperturbation}  in perturbations of Kerr, 
and let $\{\psis{p}\}_{s=\pm 2, p=0,1,2}$ be a solution to \eqref{eq:waveeqwidetildepsi1:gtilde:Teu} satisfying \eqref{eq:causlityrelationsforwidetildepsi1}. Then, we have, for $s=\pm 2$ and any $\de\in(0,\frac{1}{3}]$, 
\bea\lab{eq:EMFtotalp:sumup:rweightscontrolled:pm2}
\EMFtotalp{s} &\les&  \IE{\pmb\phi_s, \pmb\psi_s}  + \NNttotalp{s}\nn\\
&&+\A[\pmb\psi_{s}](\Iti)+\A[\pmb\phi_{s}](\tau_1,\tau_2)+\sum_{p=0}^2\Errdefect[\psis{p}],
\eea
where $\EMFtotalp{s}$, $\IE{\pmb\phi_s, \pmb\psi_s}$, $\NNttotalp{s} $, $\Errdefect[\c]$ and $\A[\c](\c,\c)$  are given respectively as in \eqref{def:EMFtotalps:pm2}, \eqref{eq:definitionofinitialenergyofphiandpsi:IEterm},  \eqref{expression:NNttotalps:pm2}, \eqref{def:Errdefectofpsi} and \eqref{def:AandAonorms:tau1tau2:phisandpsis}.
\end{theorem}

The rest of Section \ref{sec:proofofth:main:intermediary:0order} is as follows. In Section \ref{subsect:globalEMFestimate:localenergyestis}, we derive local-in-time energy estimates for $\phis{p}$, and in Section \ref{subsect:globalEMFestimate:scalarizeeqfromtensorial}, we apply Theorem \ref{th:mainenergymorawetzmicrolocal} to the equations \eqref{eq:waveeqwidetildepsi1:gtilde:Teu} of $\psiss{ij}{p}$ to deduce preliminary EMF estimates that are conditional, in particular, on the control of the error terms arising from the linear coupling terms $\chi_{\tau_1, \tau_2}^{(1)}L^{(p)}_{s,ij}$. These error terms are in turn controlled in Sections \ref{subsect:controloferrorfromlinearcoupling} and \ref{subsect:improvedMora:nearinfinity:phisp01}. In the end, we prove the EMF estimates of Theorem \ref{thm:EMF:systemofTeuscalarized:order0:final} in Section \ref{subsect:proofofEMFmain:reg=0}.


\subsection{Local energy estimates for $\phis{p}$}
\lab{subsect:globalEMFestimate:localenergyestis}


Recall from \eqref{eq:causlityrelationsforwidetildepsi1} that the following holds
\bea
\lab{eq:psispequalsphispintau1+1tau2-1}
\psis{p}=\phis{p}, \quad \text{in } \MM(\tau_1+1,\tau_2-3), \quad s=\pm 2, \quad p=0,1,2.
\eea
For our later purpose of deriving global energy-Morawetz estimates for the globally extended coupled system of wave equations \eqref{eq:waveeqwidetildepsi1:gtilde:Teu}, we first need to derive  local (in time) energy estimates for $\phis{p}$ in the regions $\MM(\tau_1,\tau_1+1)$ and $\MM(\tau_2-3,\tau_2)$, where $\psis{p}$ differs from $\phis{p}$. 

Recall that the Teukolsky equations \eqref{eq:ScalarizedTeuSys:general:Kerrperturbation} satisfied by $\phiss{ij}{p}$ have the following form
\beaa
\square_{\g}\phiss{ij}{0}&=&\sum_{k,l}O(r^{-2})\dk^{\leq 1}\phiss{kl}{0} + O(r^{-3}) \phiss{ij}{1} + N_{W,s,ij}^{(0)},\\
\square_{\g}\phiss{ij}{1}&=&\sum_{k,l}O(r^{-2})\dk^{\leq 1}\phiss{kl}{1} + O(r^{-3}) \phiss{ij}{2} +\sum_{k,l}O(r^{-2})\dk^{\leq 1}\phiss{kl}{0}+ N_{W,s,ij}^{(1)},\\
\square_{\g}\phiss{ij}{2}&=&\sum_{k,l}O(r^{-2})\dk^{\leq 1}\phiss{kl}{2} +\sum_{k,l}\sum_{p=0,1}O(r^{-2})\dk^{\leq 1}\phiss{kl}{p}+ N_{W,s,ij}^{(2)},
\eeaa
so that they all solve a coupled system of scalar wave equations of the form  \eqref{eq:eqsforlocalenergyestimatelemma:general} with $D_1=0$. Hence, applying \eqref{eq:localenergyestimate:future} with $(\tau_0=\tau_2-3,q=3)$ and $(\tau_0=\tau_1, q=1)$  to this system of wave  equations, and using Cauchy-Schwarz to control the coupling terms between the three equations, we deduce, for any $\de\in(0,1]$,
\bsub
\label{eq:localenergyestimate:future:Teu:phis01}
\bea
\label{eq:localenergyestimate:future:Teu:phis0}
\EMF_{\de}[\phis{0}](\tau_2-3, \tau_2) & \les& \E[\phis{0}](\tau_2-3) + \NNtlede[\phis{0}, \N_{W,s}^{(0)}](\tau_2-3, \tau_2)\nn\\
&&+\A[\phis{1}](\tau_2-3, \tau_2),\\
\label{eq:localenergyestimate:future:Teu:phis1}
\EMF_{\de}[\phis{1}](\tau_2-3, \tau_2) & \les& \E[\phis{1}](\tau_2-3) + \NNtlede[\phis{1}, \N_{W,s}^{(1)}](\tau_2-3, \tau_2)\nn\\
&&+\A[\phis{2}](\tau_2-3, \tau_2)
+\sup_{\tau\in[\tau_2-3,\tau_2]} \E[\phis{0}](\tau),\\
\label{eq:localenergyestimate:future:Teu:phis2}
\EMF_{\de}[\phis{2}](\tau_2-3, \tau_2) & \les& \E[\phis{2}](\tau_2-3) + \NNtlede[\phis{2}, \N_{W,s}^{(2)}](\tau_2-3, \tau_2)\nn\\
&&+\sum_{p=0,1}\sup_{\tau\in[\tau_2-3,\tau_2]} \E[\phis{p}](\tau)
\eea
\esub
and
\bsub
\label{eq:localenergyestimate:past:Teu:phis01}
\bea
\label{eq:localenergyestimate:past:Teu:phis0}
\EMF_{\de}[\phis{0}](\tau_1, \tau_1+1)& \les& \E[\phis{0}](\tau_1) + \NNtlede[\phis{0}, \N_{W,s}^{(0)}](\tau_1, \tau_1+1)\nn\\
&&+\A[\phis{1}](\tau_1, \tau_1+1),\\
\label{eq:localenergyestimate:past:Teu:phis1}
\EMF_{\de}[\phis{1}](\tau_1, \tau_1+1) & \les& \E[\phis{1}](\tau_1)+ \NNtlede[\phis{1}, \N_{W,s}^{(1)}](\tau_1, \tau_1+1)\nn\\
&&
+\A[\phis{2}](\tau_1, \tau_1+1)+\sup_{\tau\in[\tau_1,\tau_1+1]} \E[\phis{0}](\tau),\\
\label{eq:localenergyestimate:past:Teu:phis2}
\EMF_{\de}[\phis{2}](\tau_1, \tau_1+1) & \les& \E[\phis{2}](\tau_1)
+ \NNtlede[\phis{2}, \N_{W,s}^{(2)}](\tau_1, \tau_1+1)\nn\\
&&
+\sum_{p=0,1}\sup_{\tau\in[\tau_1,\tau_1+1]} \E[\phis{p}](\tau),
\eea
\esub
where $\NNtlede$ is given by \eqref{def:NNtleinlocalenergyestimate}.


\subsection{Preliminary energy-Morawetz estimates for $\psis{p}$}
\lab{subsect:globalEMFestimate:scalarizeeqfromtensorial}


In this section, we show some preliminary energy-Morawetz estimates for $\psis{p}$, $s=\pm 2$, $p=0,1,2$, which are stated below in Proposition \ref{prop:nondegenerateEMF:psisp:couplingtermnotcontrolled}.

To begin with, an application of Theorem \ref{th:mainenergymorawetzmicrolocal} yields the following statement of global energy-Morawetz estimates for $\psis{p}$.

\begin{lemma}
For $s=\pm 2$ and $p=0,1,2$, we have the following energy-Morawetz estimate 
\bea\lab{th:eq:mainenergymorawetzmicrolocal:tensorialwave:scalarized}
\widetilde{\EMF}[\psis{p}]
&\les& \sup_{\tau\in[\tmic, \tau_1+2]}\E[\psis{p}](\tau)+\A[\psis{p}](\Iti)+\Errdefect[\psis{p}]\nn\\
&&
+\widetilde{\mathcal{N}}[\psis{p},\chi_{\tau_1, \tau_2}^{(1)}\L^{(p)}_{s}]+\widetilde{\mathcal{N}}[\psis{p}, \widehat{\F}_{total,s}^{(p)}].
\eea
\end{lemma}

\begin{proof}
Noticing that, for each $p=0,1,2$, the system of wave equations \eqref{eq:waveeqwidetildepsi1:gtilde:Teu} for $\psiss{ij}{p}$ corresponds to \eqref{eq:ScalarizedWaveeq:general:Kerrpert} with $D_0=4-2\de_{p0}$  and $F_{ij}=\chi_{\tau_1, \tau_2}^{(1)}L^{(p)}_{s,ij}+\widehat{F}_{total,s,ij}^{(p)}$, we may apply Theorem \ref{th:mainenergymorawetzmicrolocal}
to the system of equations \eqref{eq:waveeqwidetildepsi1:gtilde:Teu} of $\psiss{ij}{p}$ to infer, for $s=\pm 2$ and $p=0,1,2$,
\beaa
\widetilde{\EMF}[\psis{p}]
&\les&\sup_{\tau\in[\tmic, \tau_1+2]}\E[\psis{p}](\tau)+\A[\psis{p}](\Iti)+\Errdefect[\psis{p}]\nn\\
&&+ \widetilde{\mathcal{N}}[\psis{p}, \chi_{\tau_1, \tau_2}^{(1)}\L^{(p)}_{s}+\widehat{\F}_{total,s}^{(p)}]\nn\\
&\les& \sup_{\tau\in[\tmic, \tau_1+2]}\E[\psis{p}](\tau)
+\A[\psis{p}](\Iti)+\Errdefect[\psis{p}]\nn\\
&&+\widetilde{\mathcal{N}}[\psis{p},\chi_{\tau_1, \tau_2}^{(1)}\L^{(p)}_{s}]+\widetilde{\mathcal{N}}[\psis{p},\widehat{\F}_{total,s}^{(p)}]
\eeaa
as desired. This concludes the proof of the lemma.
\end{proof}

Next, we make use of the Teukolsky transport equations \eqref{eq:ScalarizedQuantitiesinTeuSystem:Kerrperturbation} to derive Morawetz estimates for $\{\psis{p}\}_{p=0,1}$ and $\{\phis{p}\}_{p=0,1}$ which do not contain degeneracy in the trapping region $\Mtrap(\tau_1,\tau_2)$. To this end, we define, for tensors $\pmb\psi\in\sk_2(\mathbb{C})$, $\F\in\sk_2(\mathbb{C})$ and any $\de\in(0,1]$, 
\begin{align}
\lab{def:NNtdemora}
&\NNtdemora[\pmb\psi,\F](\tau_1+1,\tau_2-3)\nn\\
:=&\Bigg|\int_{\MM_{r\geq 12m}(\tau_1+1,\tau_2-3)}\sum_{i,j}\Re\Bigg(F_{ij}\ov{\bigg(X_{\de}(\psi_{ij})  - \chi_2(X_\de)^{\a}M_{i\a}^k\psi_{kj}- \chi_2(X_\de)^{\a}M_{j\a}^k\psi_{ik} + \frac{1}{2}w_{\de}\psi_{ij}   \bigg)}\Bigg) \Bigg|,
\end{align}
where the vectorfield $X_\de$ and the scalar function $w_{\de}$ are given as in \eqref{def:Xandw:improvedMorawetz:extendedRW:Kerrpert} and where the cut-off function $\chi_2$,  introduced in \eqref{def:cutoffsintime1234}, is supported in $\tau\in(\tau_1+1, \tau_2-2)$ and satisfies $\chi_2=1$ on $\tau\in(\tau_1+2, \tau_2-3)$. Also, define $\{\widetilde{\mathcal{N}}_{W,\de}[\psis{p},\phis{p}]\}_{s=\pm2, p=0,1,2}$ by
\bea
\lab{def:NNtWdepsiphi}
\widetilde{\mathcal{N}}_{W,\de}[\psis{p},\phis{p}]&:={}&\widetilde{\mathcal{N}}[\psis{p},\widehat{\F}_{total,s}^{(p)}]+\NNtdemora[\psis{p}, \widehat{\F}_{total,s}^{(p)}](\tau_1+1,\tau_2-3)\nn\\
&&+ \NNtlede[\phis{p}, \N_{W,s}^{(p)}](\tau_1, \tau_1+1)+ \NNtlede[\phis{p}, \N_{W,s}^{(p)}](\tau_2-3, \tau_2),
\eea
with $\NNt[\pmb \psi,\pmb F]$, $\NNtdemora[\pmb\psi,\F](\tau',\tau'')$ and $\NNtlede[\pmb\psi, \F](\tau',\tau'')$ given respectively in \eqref{def:NNtintermsofNNtMora:NNtEner:NNtaux:wavesystem:EMF}, \eqref{def:NNtdemora} and \eqref{def:NNtleinlocalenergyestimate}.

\begin{proposition}[Preliminary EMF estimates for $\psis{p}$]
\lab{prop:nondegenerateEMF:psisp:couplingtermnotcontrolled}
For $s=\pm 2$ and $p=0,1,2$, and for any given $\de\in(0,1]$, we have the following energy-Morawetz estimates 
\bsub
\lab{th:eq:mainnondegeenergymorawetzmicrolocal:tensorialwave:scalarized}
\bea
\lab{th:eq:mainnondegeenergymorawetzmicrolocal:tensorialwave:scalarized:0}
&&\widetilde{\EMF}[\psis{0}]
+\widehat{\M}_{\de}[\psis{0}](\tau_1+1, \tau_2-3)
+\widehat{\EMF}_{\de}[\phis{0}](\tau_1,\tau_2)\nn\\
&\les& \IE{\pmb\phi_s, \pmb\psi_s}+ \widetilde{\mathcal{N}}[\psis{0},\chi_{\tau_1, \tau_2}^{(1)}\L^{(0)}_{s}]
+\NNtdemora[\psis{0},\L_{s}^{(0)}](\tau_1+1,\tau_2-3)+\widetilde{\mathcal{N}}_{W,\de}[\psis{0},\phis{0}]\nn\\
&&+\|\N^{(0)}_{W,s}\|^2_{L^2(\MM_{r\leq 12m}(\tau_1+1,\tau_2-3))}+ \|\N^{(0)}_{T,s}\|^2_{L^2(\MM_{r\leq 12m}(\tau_1+1,\tau_2-3))}\nn\\
&&+\A[\pmb\psi_s](\Iti)+\A[\pmb\phi_s](\tau_1,\tau_2)+\Errdefect[\psis{0}],\\
\lab{th:eq:mainnondegeenergymorawetzmicrolocal:tensorialwave:scalarized:1}
&&\widetilde{\EMF}[\psis{1}]
+\widehat{\M}_{\de}[\psis{1}](\tau_1+1, \tau_2-3)
+\widehat{\EMF}_{\de}[\phis{1}](\tau_1,\tau_2)\nn\\
&\les& \IE{\pmb\phi_s, \pmb\psi_s} +  \widetilde{\mathcal{N}}[\psis{1},\chi_{\tau_1, \tau_2}^{(1)}\L^{(1)}_{s}]+\NNtdemora[\psis{1},\L_{s}^{(1)}](\tau_1+1,\tau_2-3)+\widetilde{\mathcal{N}}_{W,\de}[\psis{1},\phis{1}]\nn\\
&&+\|\N^{(1)}_{W,s}\|^2_{L^2(\MM_{r\leq 12m}(\tau_1+1,\tau_2-3))}+ \|\N^{(1)}_{T,s}\|^2_{L^2(\MM_{r\leq 12m}(\tau_1+1,\tau_2-3))}
\nn\\
&&+\A[\pmb\psi_s](\Iti)+\A[\pmb\phi_s](\tau_1,\tau_2)+\Errdefect[\psis{1}] +\widehat{\EMF}[\phis{0}](\tau_1+1,\tau_2-3)\nn\\
&&+\sup_{\tau\in[\tau_1,\tau_1+1]\cup[\tau_2-3,\tau_2]}\E[\phis{0}](\tau),\\
\lab{th:eq:mainnondegeenergymorawetzmicrolocal:tensorialwave:scalarized:2}
&&\widetilde{\EMF}[\psis{2}]
+\M_{\de}[\psis{2}](\tau_1+1, \tau_2-3)+{\EMF}_{\de}[\phis{2}](\tau_1,\tau_2)\nn\\
&\les& \IE{\pmb\phi_s, \pmb\psi_s} +  \widetilde{\mathcal{N}}[\psis{2},\chi_{\tau_1, \tau_2}^{(1)}\L^{(2)}_{s}]+\NNtdemora[\psis{2},\L_{s}^{(2)}](\tau_1+1,\tau_2-3)+\widetilde{\mathcal{N}}_{W,\de}[\psis{2},\phis{2}]\nn\\
&&+\A[\pmb\psi_s](\Iti)+\Errdefect[\psis{2}] +\sum_{p=0,1}\sup_{\tau\in[\tau_1,\tau_1+1]\cup[\tau_2-3,\tau_2]}\E[\phis{p}](\tau),
\eea
\esub
with $\IE{\pmb\phi_s, \pmb\psi_s}$, $\NNt[\psis{p},\chi_{\tau_1,\tau_2}^{(1)}\L_{s}^{(p)}]$, $\NNtdemora[\psis{p},\L_{s}^{(p)}](\tau_1+1,\tau_2-3)$, $\widetilde{\mathcal{N}}_{W,\de}[\psis{p},\phis{p}]$, $\Errdefect[\psis{p}]$ and $\A[\cdot](\cdot,\cdot)$ given as in \eqref{eq:definitionofinitialenergyofphiandpsi:IEterm},  \eqref{def:NNtintermsofNNtMora:NNtEner:NNtaux:wavesystem:EMF}, \eqref{def:NNtdemora},  \eqref{def:NNtWdepsiphi}, \eqref{def:Errdefectofpsi} and \eqref{def:AandAonorms:tau1tau2:phisandpsis}
respectively.
\end{proposition}

\begin{proof}
Using the transport equations \eqref{def:TensorialTeuScalars:wavesystem:Kerrperturbation} in the trapping region $\Mtrap(\tau_1+1,\tau_2-3)$ as well as the fact that
\bea
e_4=\frac{\R}{|q|^2}\widehat{T} +O(1)\pr_r + O(\ep) \pr_{\a},\quad e_3=\frac{\R}{\De}\widehat{T} - \pr_r + O(\ep) \pr_{\a},  \quad  \text{in }\Mtrap,
\eea
where $\widehat{T}:=\pr_{\tt} + \frac{a}{r^2+a^2}\pr_{\tphi}$,  and noticing that the square of the $L^2$ norm of $\nab_{\pr_r}^{\leq 1}\pmb\phi^{(p)}_{s}$ on $\Mtrap(\tau_1+1,\tau_2-3)$ is bounded by $\M[\phis{p}](\tau_1+1,\tau_2-3)$, we deduce, for $p=0,1$, 
\bea
\lab{eq:controlofhatTderivativeofpsisijp:proofofEMF:pre}
\|\nab_{\widehat{T}}\phis{p}\|^2_{L^2(\Mtrap(\tau_1+1,\tau_2-3))}&\les &\|\phis{p+1}\|^2_{L^2(\Mtrap(\tau_1+1,\tau_2-3))}+ \|\N^{(p)}_{T,s}\|^2_{L^2(\Mtrap(\tau_1+1,\tau_2-3))}\nn\\
&&+\M[\phis{p}](\tau_1+1,\tau_2-3)+\ep\widehat{ \M}[\phis{p}](\tau_1+1,\tau_2-3).
\eea

Next, we use the following decomposition of $\squared_2$ which is a non-sharp consequence of Lemma 4.7.2 in \cite{GKS22}
\beaa
|q|^2\squared_2 = -\frac{(\R)^2}{\De}\nab_{\widehat{T}}^2+|q|^2\De_2 +O(1)\nab_{\pr_r}^2+O(\ep)\pr^2+O(1)\pr \quad  \text{in }\,\,\MM_{r_+(1+\dbl), 12m}.
\eeaa
Together with the Teukolsky wave equations \eqref{eq:TensorialTeuSys:rescaleRHScontaine2:general:Kerrperturbation}, we infer, for $s=\pm 2$, $p=0,1,2$,
\beaa
&&\bigg(\Delta_2-\frac{4-2\de_{p0}}{\qs}\bigg)\pmb\phi^{(p)}_{s} +O(\dbl^{-1})\nab_{\widehat{T}}^2\pmb\phi^{(p)}_{s}+O(\dbl^{-1})\nab_{\pr_r}^2\pmb\phi^{(p)}_{s}+O(1)\pr\pmb\phi^{(p)}_{s}+O(\ep)\pr^2\pmb\phi^{(p)}_{s}\\
&=& \L_{s}^{(p)}[\pmb\phi_{s}]+\N_{W,s}^{(p)}\quad \text{in}\,\,\MM_{r_+(1+\dbl), 12m}(\tau_1+1,\tau_2-3).
\eeaa
We then contract this equation with $\ov{-\chi^2(r)\pmb\phi^{(p)}_{s}}$, where $\chi(r)$ is a smooth cut-off function in $r$ which equals $1$ on $\Mtrap=\MM_{r_+(1+2\dbl), 10m}$ and vanishes outside $\MM_{r_+(1+\dbl), 12m}$, take the real part, integrate over $\Mtrap(\tau_1+1,\tau_2-3)$, and integrate second order derivatives  by parts, arriving at\footnote{Recall from Section \ref{sec:smallnesconstants} that $\les$ may in particular depend on $O(\dbl^{-1})$ factors.}
\beaa
&&\|\nab\phis{p}\|^2_{L^2(\Mtrap(\tau_1+1,\tau_2-3))}\nn\\
&\les &\|\nab_{\widehat{T}}\phis{p}\|^2_{L^2(\MM_{r\leq 12m}(\tau_1+1,\tau_2-3))}+\M[\phis{p}](\tau_1+1,\tau_2-3)
\nn\\
&&+\Big(\widehat{\EMF}[\phis{p}](\tau_1+1,\tau_2-3)
+\|\L^{(p)}_{s}\|^2_{L^2(\MM_{r\leq 12m}(\tau_1+1,\tau_2-3))}\Big)^{\frac{1}{2}}\Big(\A[\pmb\phi_s](\tau_1+1,\tau_2-3)\Big)^{\frac{1}{2}}\nn\\
&&
+\ep\widehat{\M}[\phis{p}](\tau_1+1,\tau_2-3)+\|\N^{(p)}_{W,s}\|^2_{L^2(\MM_{r\leq 12m}(\tau_1+1,\tau_2-3))}
,\qquad p=0,1.
\eeaa
In view of the form \eqref{eq:tensor:Lsn:onlye_2present:general:Kerrperturbation} of the linear coupling terms $\L^{(p)}_{s}$, and using the estimate \eqref{eq:controlofhatTderivativeofpsisijp:proofofEMF:pre} to control $\|\nab_{\widehat{T}}\phis{p}\|^2_{L^2(\Mtrap(\tau_1+1,\tau_2-3))}$, we deduce
\beaa
&&\|\nab\phis{0}\|^2_{L^2(\Mtrap(\tau_1+1,\tau_2-3))}\nn\\
&\les& \M[\phis{0}](\tau_1+1,\tau_2-3)
+\ep\widehat{\M}[\phis{0}](\tau_1+1,\tau_2-3)+\|\N^{(0)}_{W,s}\|^2_{L^2(\MM_{r\leq 12m}(\tau_1+1,\tau_2-3))}\nn\\
&& +\Big(\widehat{\EMF}[\phis{0}](\tau_1+1,\tau_2-3)
\Big)^{\frac{1}{2}}\Big(\A[\pmb\phi_s](\tau_1+1,\tau_2-3)\Big)^{\frac{1}{2}}\nn\\
&&
+ \|\N^{(0)}_{T,s}\|^2_{L^2(\MM_{r\leq 12m}(\tau_1+1,\tau_2-3))}
+\A[\pmb\phi_s](\tau_1+1,\tau_2-3)
\eeaa
and
\beaa
&&\|\nab\phis{1}\|^2_{L^2(\Mtrap(\tau_1+1,\tau_2-3))}\nn\\
&\les & \M[\phis{1}](\tau_1+1,\tau_2-3)
+\ep\widehat{\M}[\phis{1}](\tau_1+1,\tau_2-3)+\|\N^{(1)}_{W,s}\|^2_{L^2(\MM_{r\leq 12m}(\tau_1+1,\tau_2-3))}\nn\\
&&+\Big(\widehat{\M}[\phis{0}](\tau_1+1,\tau_2-3)+\widehat{\EMF}[\phis{1}](\tau_1+1,\tau_2-3)\Big)^{\frac{1}{2}}\Big(\A[\pmb\phi_s](\tau_1+1,\tau_2-3)\Big)^{\frac{1}{2}}\nn\\
&& + \|\N^{(1)}_{T,s}\|^2_{L^2(\MM_{r\leq 12m}(\tau_1+1,\tau_2-3))}
+\A[\pmb\phi_s](\tau_1+1,\tau_2-3).
\eeaa
 Taking $\ep$ suitably small, and combining with the estimate \eqref{eq:controlofhatTderivativeofpsisijp:proofofEMF:pre} which controls $\widehat{T}$ derivative, this yields
\beaa
&&\widehat{\EMF}[\phis{0}](\tau_1+1,\tau_2-3)\nn\\
&\les& {\EMF}[\phis{0}](\tau_1+1,\tau_2-3)+\A[\pmb\phi_s](\tau_1+1,\tau_2-3)\nn\\
&& +\|\N^{(0)}_{W,s}\|^2_{L^2(\MM_{r\leq 12m}(\tau_1+1,\tau_2-3))}+ \|\N^{(0)}_{T,s}\|^2_{L^2(\MM_{r\leq 12m}(\tau_1+1,\tau_2-3))}
,\\
&&\widehat{\EMF}[\phis{1}](\tau_1+1,\tau_2-3)\nn\\
&\les& {\EMF}[\phis{1}](\tau_1+1,\tau_2-3)+\A[\pmb\phi_s](\tau_1+1,\tau_2-3)
+\|\N^{(1)}_{W,s}\|^2_{L^2(\MM_{r\leq 12m}(\tau_1+1,\tau_2-3))}\nn\\
&&+ \|\N^{(1)}_{T,s}\|^2_{L^2(\MM_{r\leq 12m}(\tau_1+1,\tau_2-3))}+\Big(\widehat{\EMF}[\phis{0}](\tau_1+1,\tau_2-3)\Big)^{\frac{1}{2}}\Big(\A[\pmb\phi_s](\tau_1+1,\tau_2-3)\Big)^{\frac{1}{2}},
\eeaa
and together with  the microlocal EMF estimate \eqref{th:eq:mainenergymorawetzmicrolocal:tensorialwave:scalarized} and using the fact \eqref{eq:psispequalsphispintau1+1tau2-1} that $\psis{p}=\phis{p}$ in $\MM(\tau_1+1,\tau_2-3)$, we infer
\bsub
\lab{eq:mainnondegeenergymorawetzmicrolocal:1:tensorialwave:scalarized}
\bea
&&\widetilde{\EMF}[\psis{0}]+\widehat{\EMF}[\phis{0}](\tau_1+1,\tau_2-3)\nn\\
&\les&\sup_{\tau\in[\tmic, \tau_1+2]}\E[\psis{0}](\tau)+ \widetilde{\mathcal{N}}[\psis{0},\chi_{\tau_1, \tau_2}^{(1)}\L^{(0)}_{s}]+\widetilde{\mathcal{N}}[\psis{0},\widehat{\F}_{total,s}^{(0)}]+\A[\pmb\psi_s](\Iti)\nn\\
&&+{\Errdefect[\psis{0}]}+\|\N^{(0)}_{W,s}\|^2_{L^2(\MM_{r\leq 12m}(\tau_1+1,\tau_2-3))}+ \|\N^{(0)}_{T,s}\|^2_{L^2(\MM_{r\leq 12m}(\tau_1+1,\tau_2-3))},\\
&&\widetilde{\EMF}[\psis{1}]+\widehat{\EMF}[\phis{1}](\tau_1+1,\tau_2-3)
\nn\\
&\les&\sup_{\tau\in[\tmic, \tau_1+2]}\E[\psis{1}](\tau)+\widetilde{\mathcal{N}}[\psis{1},\chi_{\tau_1, \tau_2}^{(1)}\L^{(1)}_{s}]+\widetilde{\mathcal{N}}[\psis{1},\widehat{\F}_{total,s}^{(1)}]+\A[\pmb\psi_s](\Iti)\nn\\
&&+\|\N^{(1)}_{W,s}\|^2_{L^2(\MM_{r\leq 12m}(\tau_1+1,\tau_2-3))}+ \|\N^{(1)}_{T,s}\|^2_{L^2(\MM_{r\leq 12m}(\tau_1+1,\tau_2-3))}\nn\\
&&+{\Errdefect[\psis{1}]}+\Big(\widehat{\EMF}[\psis{0}](\tau_1+1,\tau_2-3)\Big)^{\frac{1}{2}}\Big(\A[\pmb\psi_s](\tau_1+1,\tau_2-3)\Big)^{\frac{1}{2}},\\
&&\widetilde{\EMF}[\psis{2}]\nn\\
&\les&\sup_{\tau\in[\tmic, \tau_1+2]}\E[\psis{2}](\tau)+ \widetilde{\mathcal{N}}[\psis{2},\chi_{\tau_1, \tau_2}^{(1)}\L^{(2)}_{s}]+\widetilde{\mathcal{N}}[\psis{2},\widehat{\F}_{total,s}^{(2)}]\nn\\
&&
+{\Errdefect[\psis{2}]}+\A[\pmb\psi_s](\Iti).
\eea
\esub
Combining the above estimates \eqref{eq:mainnondegeenergymorawetzmicrolocal:1:tensorialwave:scalarized} with the local-in-time energy estimates \eqref{eq:localenergyestimate:future:Teu:phis01} and \eqref{eq:localenergyestimate:past:Teu:phis01} for $\{\phis{p}\}_{p=0,1,2}$ in $\MM(\tau_2-3,\tau_2)\cup \MM(\tau_1,\tau_1+1)$, and in view of 
the fact that $\psis{p}=\phis{p}$ in $\MM(\tau_1+1,\tau_2-3)$ from \eqref{eq:psispequalsphispintau1+1tau2-1},
we infer the following EMF estimates 
\bsub
\lab{eq:EMFpsisandphisp:012:secondpreliminary}
\bea
&&\widetilde{\EMF}[\psis{0}]+\widehat{\EMF}[\phis{0}](\tau_1,\tau_2)+\widehat{\M}_{\de}[\phis{0}](\tau_1, \tau_1+1)+\widehat{\M}_{\de}[\phis{0}](\tau_2-3, \tau_2)\nn\\
&\les& \IE{\pmb\phi_s, \pmb\psi_s} + \widetilde{\mathcal{N}}[\psis{0},\chi_{\tau_1, \tau_2}^{(1)}\L^{(0)}_{s}]+\widetilde{\mathcal{N}}_{W,\de}[\psis{0},\phis{0}]+\|\N^{(0)}_{W,s}\|^2_{L^2(\MM_{r\leq 12m}(\tau_1+1,\tau_2-3))}\nn\\
&&+ \|\N^{(0)}_{T,s}\|^2_{L^2(\MM_{r\leq 12m}(\tau_1+1,\tau_2-3))}+{\Errdefect[\psis{0}]}+\A[\pmb\psi_s](\Iti)+\A[\pmb\phi_s](\tau_1,\tau_2),\\
&&\widetilde{\EMF}[\psis{1}]
+\widehat{\EMF}[\phis{1}](\tau_1,\tau_2)
+\widehat{\M}_{\de}[\phis{1}](\tau_1, \tau_1+1)+\widehat{\M}_{\de}[\phis{1}](\tau_2-3, \tau_2)\nn\\
&\les& \IE{\pmb\phi_s, \pmb\psi_s} +\widetilde{\mathcal{N}}[\psis{1},\chi_{\tau_1, \tau_2}^{(1)}\L^{(1)}_{s}]+\widetilde{\mathcal{N}}_{W,\de}[\psis{1},\phis{1}]+\|\N^{(1)}_{W,s}\|^2_{L^2(\MM_{r\leq 12m}(\tau_1+1,\tau_2-3))}
\nn\\
&&
+ \|\N^{(1)}_{T,s}\|^2_{L^2(\MM_{r\leq 12m}(\tau_1+1,\tau_2-3))}+\A[\pmb\psi_s](\Iti)+\A[\pmb\phi_s](\tau_1,\tau_2)+{\Errdefect[\psis{1}]}
\nn\\
&&+\Big(\widehat{\EMF}[\psis{0}](\tau_1+1,\tau_2-3)\Big)^{\frac{1}{2}}\Big(\A[\pmb\psi_s](\tau_1+1,\tau_2-3)\Big)^{\frac{1}{2}} +\sup_{\tau\in[\tau_1,\tau_1+1]\cup[\tau_2-3,\tau_2]}\E[\phis{0}](\tau)\nn\\
&\les& \IE{\pmb\phi_s, \pmb\psi_s} +\widetilde{\mathcal{N}}[\psis{1},\chi_{\tau_1, \tau_2}^{(1)}\L^{(1)}_{s}]+\widetilde{\mathcal{N}}_{W,\de}[\psis{1},\phis{1}]+\|\N^{(1)}_{W,s}\|^2_{L^2(\MM_{r\leq 12m}(\tau_1+1,\tau_2-3))}
\nn\\
&&
+ \|\N^{(1)}_{T,s}\|^2_{L^2(\MM_{r\leq 12m}(\tau_1+1,\tau_2-3))}+\A[\pmb\psi_s](\Iti)+\A[\pmb\phi_s](\tau_1,\tau_2)+{\Errdefect[\psis{1}]}
\nn\\
&&+\widehat{\EMF}[\phis{0}](\tau_1+1,\tau_2-3) +\sup_{\tau\in[\tau_1,\tau_1+1]\cup[\tau_2-3,\tau_2]}\E[\phis{0}](\tau),\\
&&\widetilde{\EMF}[\psis{2}]
+{\EMF}[\phis{2}](\tau_1, \tau_2)
+{\M}_{\de}[\phis{2}](\tau_1, \tau_1+1)+\M_{\de}[\phis{2}](\tau_2-3, \tau_2)\nn\\
&\les& \IE{\pmb\phi_s, \pmb\psi_s} + \widetilde{\mathcal{N}}[\psis{2},\chi_{\tau_1, \tau_2}^{(1)}\L^{(2)}_{s}]
+\widetilde{\mathcal{N}}_{W,\de}[\psis{2},\phis{2}]+\A[\pmb\psi_s](\Iti)+{\Errdefect[\psis{2}]}\nn\\
&& +\sum_{p=0,1}\sup_{\tau\in[\tau_1,\tau_1+1]\cup[\tau_2-3,\tau_2]}\E[\phis{p}](\tau).
\eea
\esub

Finally, we improve the Morawetz estimates in $\MM(\tau_1+1,\tau_2-3)$ near infinity for $\psis{p}$, $p=0,1,2$. Noticing that, for each $p=0,1,2$, the system of wave equations \eqref{eq:waveeqwidetildepsi1:gtilde:Teu} for $\psiss{ij}{p}$ corresponds to \eqref{eq:ScalarizedWaveeq:general:Kerrpert} with $D_0=4-2\de_{p0}$  and $F_{ij}=\chi_{\tau_1, \tau_2}^{(1)}L^{(p)}_{s,ij}+\widehat{F}_{total,s,ij}^{(p)}$, we may apply the improved Morawetz estimate \eqref{eq:Morawetznearinf:extendedRWsystem:Kerrpert:delta} with $R_1=12m$ to the system of equations \eqref{eq:waveeqwidetildepsi1:gtilde:Teu} of $\psiss{ij}{p}$ to infer, for $\de\in(0,1]$ and $p=0,1,2$,\footnote{Since we are integrating over $\tau\in(\tau_1+1,\tau_2-3)$ which contains neither $[\tmic, \tau_1+1]$ nor $[\tau_2-3,\tau_2+1]$, the terms $\sup_{\tt\in[\tmic,\tau_1+1]}\E[\pmb\psi](\tau)$  and $\ep\sup_{\tt\in[\tau_2-3,\tau_2+1]}\E[\pmb\psi](\tau)$ on the RHS of \eqref{eq:Morawetznearinf:extendedRWsystem:Kerrpert:delta} do not show up in the resulting estimate. Also, note that the term involving $F_{ij}$ on the LHS of  \eqref{eq:Morawetznearinf:extendedRWsystem:Kerrpert:delta} for $(\tau', \tau'')=(\tau_1+1, \tau_2-3)$ is consistent with the definition of $\NNtdemora[\pmb\psi,\F](\tau_1+1,\tau_2-3)$ in view of the definition of $X_{\tau_2}(\psi)_{ij}$ in \eqref{def:Xtau2psiij:generalforEM}.} 
 \beaa
\nn &&\M_{\de, r\geq 12m}[\psis{p}](\tau_1+1,\tau_2-3)+\M_{\de, r\geq 12m}[\phis{p}](\tau_1+1,\tau_2-3)\\
&\les& \EMF[\psis{p}](\tau_1+1,\tau_2-3) + \NNtlocal[\psis{p}](\tau_1+1,\tau_2-3)\nn\\
&&+\NNtdemora[\psis{p},\L_{s}^{(p)}](\tau_1+1,\tau_2-3)+\NNtdemora[\psis{p},\widehat{\F}_{total,s}^{(p)}](\tau_1+1,\tau_2-3),
\eeaa
where we have added the term $\M_{\de, r\geq 12m}[\phis{p}](\tau_1+1,\tau_2-3)$ to the LHS in view of 
the fact that $\psis{p}=\phis{p}$ in $\MM(\tau_1+1,\tau_2-3)$ from \eqref{eq:psispequalsphispintau1+1tau2-1}. Adding this improved Morawetz estimate to the above estimates \eqref{eq:EMFpsisandphisp:012:secondpreliminary} for $p=0,1,2$ respectively, and using the definition \eqref{def:NNtlocalinr:NNtEner:NNtaux:wavesystem:EMF:proof} of $\NNtlocal[\c](\tau',\tau'')$, 
we then conclude the desired EMF estimates \eqref{th:eq:mainnondegeenergymorawetzmicrolocal:tensorialwave:scalarized}
 for $\psis{p}$ and $\phis{p}$. This concludes the proof of Proposition \ref{prop:nondegenerateEMF:psisp:couplingtermnotcontrolled}.
\end{proof}


\subsection{Control of the error terms arising from the linear coupling terms}
\lab{subsect:controloferrorfromlinearcoupling}


In this section, we provide estimates for the error terms arising from the linear coupling terms, that is, for the terms $\{\widetilde{\mathcal{N}}[\psis{p},\chi_{\tau_1, \tau_2}^{(1)}\L^{(p)}_{s}]\}_{s=\pm2, p=0,1,2}$ and $\{\NNtdemora[\psis{p},\L_{s}^{(p)}](\tau_1+1,\tau_2-3)\}_{s=\pm2, p=0,1,2}$ on the RHS of the estimates \eqref{th:eq:mainnondegeenergymorawetzmicrolocal:tensorialwave:scalarized}. To this end, for $\de>0$, we define a bulk term for a tensor $\pmb\psi\in \sk_2(\mathbb{C})$ 
\bea
\lab{def:ExtrabulknotcontrolledbystandardMora}
\Extra[\pmb\psi](\tau',\tau''):= \int_{\MM(\tau',\tau'')}r^{-3+\de}|(\nab_{\pr_{\tphi}+ a\pr_{\tt}})^{\leq 1}\pmb\psi|^2.
\eea

The estimates for these error terms are contained in the following proposition.

\begin{proposition}[Control of the error terms arising from the linear coupling terms]
\lab{prop:naivecontroloferrorfromlinearcoupling:scalarizedTeu}
For $s=\pm 2$, we have
\bsub
\lab{subeq:naivecontroloferrorfromlinearcoupling:scalarizedTeu}
\bea
&&\widetilde{\mathcal{N}}[\psis{0},\chi_{\tau_1, \tau_2}^{(1)}\L^{(0)}_{s}]
+\NNtdemora[\psis{0},\L_{s}^{(0)}](\tau_1+1,\tau_2-3)\nn\\
&\les& \left(\A[\pmb\phi_s](\tau_1,\tau_2-2)+\Ao[r\nab\phis{0}](\tau_1,\tau_2-2) +\ep\widehat{\M}[\phis{0}](\tau_1,\tau_2-2)\right)^{\frac{1}{2}}
\nn\\
&&\times\Big(\widehat{\M}_{\de}[\psis{0}](\tau_1+1,\tau_2-3)+\widetilde{\EM}[\psis{0}] \Big)^{\frac{1}{2}} +\Ao[r\nab\phis{0}](\tau_1,\tau_2-2)+\A[\pmb\phi_s](\tau_1,\tau_2-2)\nn\\
&& +\ep \widehat{\M}[\phis{0}](\tau_1,\tau_2-2)\\
&&\widetilde{\mathcal{N}}[\psis{1},\chi_{\tau_1, \tau_2}^{(1)}\L^{(1)}_{s}]
+\NNtdemora[\psis{1},\L_{s}^{(1)}](\tau_1+1,\tau_2-3)\nn\\
&\les& \Bigg(\Extra[\phis{0}](\tau_1,\tau_2-2)+\A[\pmb\phi_s](\tau_1,\tau_2-2)
+\Ao[r\nab\phis{1}](\tau_1,\tau_2-2) +\ep\widehat{\M}[\phis{1}](\tau_1,\tau_2-2) \nn\\ 
&& +\sup_{\tau\in[\tau_1, \tau_1+1]\cup[\tau_2-3,\tau_2-2]}\E[\phis{0}](\tau)\Bigg)^{\frac{1}{2}}\Big(\widehat{\M}_{\de}[\psis{1}](\tau_1+1,\tau_2-3)+\widetilde{\EM}[\psis{1}] \Big)^{\frac{1}{2}}\nn\\
&&+\widehat{\M}[\phis{0}](\tau_1,\tau_2-2)+ \Ao[r\nab\phis{1}](\tau_1,\tau_2-2)+\A[\pmb\phi_s](\tau_1,\tau_2-2)\nn\\
&&+\ep \widehat{\M}[\phis{1}](\tau_1,\tau_2-2),\\
&&\widetilde{\mathcal{N}}[\psis{2},\chi_{\tau_1, \tau_2}^{(1)}\L^{(2)}_{s}]+\NNtdemora[\psis{2},\L_{s}^{(2)}](\tau_1+1,\tau_2-3)\nn\\
&\les& \bigg(\sum_{p=0,1}\widehat{\EM}[\phis{p}](\tau_1,\tau_2-2)+\sum_{p=0,1}\Extra[\phis{p}](\tau_1,\tau_2-2)+\A[\pmb\phi_s](\tau_1,\tau_2-2)\bigg)^{\frac{1}{2}} \nn\\
&&\times\Big(\widetilde{\EM}[\psis{2}]+\EM[\phis{2}](\tau_1,\tau_2-2)
+{\M}_{\de}[\psis{2}](\tau_1+1,\tau_2-3)\Big)^{\frac{1}{2}} +\A[\pmb\phi_s](\tau_1,\tau_2-2)\nn\\
&&+\int_{\MM_{r\leq 12m}(\tau_1,\tau_2-2)}|\pr^{\leq 1}\N_{T, s}^{(1)}|^2+  \widehat{\M}[\phis{0}](\tau_1,\tau_2-2)+\widehat{\M}[\phis{1}](\tau_1,\tau_2-2),
\eea
\esub
where $\A[\cdot](\cdot,\cdot)$, $\Ao[\cdot](\cdot,\cdot)$ and $\Extra[\cdot](\cdot,\cdot)$ are given as in \eqref{def:AandAonorms:tau1tau2:phisandpsis} and \eqref{def:ExtrabulknotcontrolledbystandardMora}.
\end{proposition}

The rest of this section is devoted to proving this proposition.
To this end, we first
recall from \eqref{def:NNtintermsofNNtMora:NNtEner:NNtaux:wavesystem:EMF} and \eqref{def:NNtMora:NNtEner:NNtaux:wavesystem:EMF} that
\bea
\lab{def:NNtnorms:copy}
\NNt[\pmb\psi,\pmb F]=\NNtmora[\pmb\psi, \pmb F]+\NNtener[\pmb\psi, \pmb F]+\NNtaux[\pmb F],
\eea
with
\bsub
\lab{def:NNtmoraeneraux:copy}
\begin{align}
\lab{def:NNtmoracopy}
\NNtmora[\pmb \psi, \pmb F]=&\sum_{i,j}\bigg|\int_{\MM_{r_+(1+\dhor'),\Rmic}}\Re\Big(F_{ij} \ov{X\psi_{ij}}\Big)+\int_{\MM_{\Rmic,+\infty}(\Iti)}\Re\Big(F_{ij} \ov{({X} + w)\psi_{ij}}\Big)\bigg|,\\
\lab{def:NNtenercopy}
\NNtener[\pmb \psi, \pmb F]=&\sum_{i,j}\sup_{\tau\geq\tmic}\bigg|\int_{\Mntrap(\tmic, \tau)}{\Re\Big(F_{ij}\ov{\pr_{\tau}\psi_{ij}}\Big)}\bigg|
\nn\\
&+ \sup_{\tau\in\Reals}\sum_{n=-1}^{\iota}\sum_{i,j}\bigg|\int_{\Mtrap(\tmic,\tau)}\Re\Big(\ov{|q|^{-2}\Opw(\Theta_n)(\qs F_{ij})}V_n\Opw(\Theta_n)\psi_{ij}\Big)\bigg|,\\
\lab{def:NNtauxcopy}
\NNtaux[\pmb F]=&\sum_{i,j}\bigg(\int_{\Mntrap(\Iti)}\frac{|F_{ij}|^2}{r} +(\ep+\dhor)\int_{\Mtrap(\Iti)}\frac{|F_{ij}|^2}{r}  + \ep\int_{\MM}\tt^{-1-\dec}|F_{ij}|^2\bigg),
\end{align}
\esub
where $\dhor'$, $\Rmic$, $X, w, \Th_n$ and $V_n$ satisfy the properties given in  Section \ref{sec:definitionmicrolocalenergyMorawetznorms}.

Recall from \eqref{def:NNtdemora} the formula of the terms $\{\NNtdemora[\psis{p},\L_{s}^{(p)}](\tau_1+1,\tau_2-3)\}_{s=\pm2, p=0,1,2}$, where the linear coupling terms ${L_{s,ij}^{(p)}}=(\L_{s}^{(p)}[\pmb\phi_s])_{ij}$ are given, in view of \eqref{eq:linearterms:ScalarizedTeuSys:general:Kerrperturbation}, by 
\bsub
\lab{eq:expression:Lsijp:copy}
\begin{align}
{L_{s,ij}^{(0)}}={}& (2sr^{-3} +O(mr^{-4}))\phiss{ij}{1} + O(mr^{-3}){\Xcal_s}\phiss{ij}{0}\nn\\
&+\sum_{k,l=1,2,3}O(mr^{-3})\phiss{kl}{0},\\
{L_{s,ij}^{(1)}}={}&(sr^{-3} +O(mr^{-4}))\phiss{ij}{2} + O(mr^{-3}){\Xcal_s}\phiss{ij}{1}+O(mr^{-2})(\pr_{\tphi} +a\pr_\tau)\phiss{ij}{0}\nn\\
  &+\sum_{k,l=1,2,3}\Big(O(mr^{-3})\phiss{kl}{1}+O(mr^{-2})\phiss{kl}{0}\Big),\\
{L_{s,ij}^{(2)}}={}&h^{(2,2)} \phiss{ij}{2}+h^{(2,1)}(\pr_{\tphi} +a\pr_\tau)\phiss{ij}{1} +O(m^2r^{-2})\phiss{ij}{0}\nn\\
&+\sum_{k,l=1,2,3}O(mr^{-2})\phiss{kl}{1},
\end{align}
\esub
with $\Xcal_s$, $s=\pm 2$, being the regular vectorfields introduced in \eqref{eq:formofregularhorizontalvectorfieldmathcalXs:Kerrperturbation}, with the coefficients $h^{(2,2)}$ and $h^{(2,1)}$ being real-valued scalar functions satisfying 
\beaa
h^{(2,2)}=O(mr^{-3}), \qquad h^{(2,1)}=O(mr^{-2}),
\eeaa
and with all the coefficients in the first line\footnote{The coefficients on the second line of the three equations in \eqref{eq:expression:Lsijp:copy} involve $M_{i\a}^j$ so that they are independent of $\tau$ but depend on $\tphi$.} of the three equations in \eqref{eq:expression:Lsijp:copy} being independent of coordinates $\tau$ and $\tphi$.

We now estimate the error terms arising from the linear coupling terms starting with the control of $\NNtaux[\chi_{\tau_1, \tau_2}^{(1)}\L^{(p)}_{s}]$.

\begin{lemma}[Bounds for the error terms $\NNtaux$ coming from linear coupling terms]
\lab{lem:NNtaux:scalarizedTeukolsky:estimates}
For the terms $\NNtaux[\chi_{\tau_1, \tau_2}^{(1)}\L^{(p)}_{s}]$, $s=\pm2$, $p=0,1,2$, we have
 \bsub
\lab{eq:esti:linearcouplingterm:NNtaux}
 \begin{align}
\NNtaux[\chi_{\tau_1, \tau_2}^{(1)}\L^{(0)}_{s}] \les{}& \Ao[(r\nab)^{\leq 1}\phis{0}](\tau_1,\tau_2-2)+\A[\phis{1}](\tau_1,\tau_2-2)\nn\\
& +\ep \widehat{\M}[\phis{0}](\tau_1,\tau_2-2),\\
\NNtaux[\chi_{\tau_1, \tau_2}^{(1)}\L^{(1)}_{s}]\les{}& \widehat{\M}[\phis{0}](\tau_1,\tau_2-2)+ \Ao[(r\nab)^{\leq 1}\phis{1}](\tau_1,\tau_2-2) +\A[\phis{2}](\tau_1,\tau_2-2)\nn\\
& +\ep \widehat{\M}[\phis{1}](\tau_1,\tau_2-2),\\
\NNtaux[\chi_{\tau_1, \tau_2}^{(1)}\L^{(2)}_{s}]\les{}& \widehat{\M}[\phis{0}](\tau_1,\tau_2-2)+\widehat{\M}[\phis{1}](\tau_1,\tau_2-2)+\A[\phis{2}](\tau_1,\tau_2-2).
\end{align}
 \esub
 \end{lemma}
 
 \begin{proof}
In view of the definition of $\NNtaux[\F]$, see \eqref{def:NNtauxcopy}, and the fact that $\textrm{supp}(\chi_{\tau_1, \tau_2}^{(1)})\subset [\tau_1,\tau_2-2]$, we have
\beaa
\NNtaux[\chi_{\tau_1, \tau_2}^{(1)}\L^{(p)}_{s}]\les \sum_{i,j}\int_{\MM}|\chi_{\tau_1, \tau_2}^{(1)}L^{(p)}_{s,ij}|^2\les \sum_{i,j}\int_{\MM(\tau_1,\tau_2-2)}|L^{(p)}_{s,ij}|^2. 
\eeaa 
Then, the estimates \eqref{eq:esti:linearcouplingterm:NNtaux} follow immediately from the above form \eqref{eq:expression:Lsijp:copy} of $L_{s,ij}^{(p)}$ and the fact that
\bea\lab{eq:formofregularvectorfieldmathcalXs:Kerr:copy:horizontal:consequence}
|\mathcal{X}_s\phiss{ij}{p}|\lesssim |r\nab\phiss{ij}{p}|+\ep|\dk\phiss{ij}{p}|,\quad p=0,1,
\eea
in view of \eqref{eq:formofregularvectorfieldmathcalXs:Kerr:copy:horizontal}.
\end{proof}

Next, we estimate in the following lemma the error terms $\NNtdemora[\psis{p},\L_{s}^{(p)}](\tau_1+1,\tau_2-3)$  and  $\NNtmora[\psis{p},\chi_{\tau_1, \tau_2}^{(1)}\L^{(p)}_{s}]$, both of which  arise from deriving the Morawetz estimates.

\begin{lemma}[Bounds for the error terms $\NNtmora$ and $\NNtdemora$ coming from linear coupling terms]
\lab{lem:NNtmora:scalarizedTeukolsky:estimates}
For the terms $\NNtmora[\psis{p},\chi_{\tau_1, \tau_2}^{(1)}\L^{(p)}_{s}]$ and $\NNtdemora[\psis{p},\L_{s}^{(p)}](\tau_1+1,\tau_2-3)$, $s=\pm2$, $p=0,1,2$, we have 
\bsub
\lab{eq:esti:NNtmora:psisp:012:nocontrolonangularderis}
\begin{align}
\lab{eq:esti:NNtmora:psisp:0:nocontrolonangularderis}
&\NNtmora[\psis{0},\chi_{\tau_1, \tau_2}^{(1)}\L^{(0)}_{s}]
+\NNtdemora[\psis{0},\L_{s}^{(0)}](\tau_1+1,\tau_2-3)\nn\\
\les {}&\left(\sum_{p=0,1}\A[\phis{p}](\tau_1,\tau_2-2)+\Ao[r\nab\phis{0}](\tau_1,\tau_2-2) +\ep\widehat{\M}[\phis{0}](\tau_1,\tau_2-2)\right)^{\frac{1}{2}}
\nn\\
&\times\Big(\widehat{\M}_{\de}[\psis{0}](\tau_1+1,\tau_2-3)+\widetilde{\EM}[\psis{0}] \Big)^{\frac{1}{2}},\\
\lab{eq:esti:NNtmora:psisp:1:nocontrolonangularderis}
&\NNtmora[\psis{1},\chi_{\tau_1, \tau_2}^{(1)}\L^{(1)}_{s}]
+\NNtdemora[\psis{1},\L_{s}^{(1)}](\tau_1+1,\tau_2-3)\nn\\
\les{}& \Bigg(\Extra[\phis{0}](\tau_1,\tau_2-2)+\sum_{p=1,2}\A[\phis{p}](\tau_1,\tau_2-2)
+\Ao[r\nab\phis{1}](\tau_1,\tau_2-2)\nn\\
& +\ep\widehat{\M}[\phis{1}](\tau_1,\tau_2-2) +\sup_{\tau\in[\tau_1, \tau_1+1]\cup[\tau_2-3,\tau_2-2]}\E[\phis{0}](\tau)\Bigg)^{\frac{1}{2}}\nn\\
&\times \Big(\widehat{\M}_{\de}[\psis{1}](\tau_1+1,\tau_2-3)+\widetilde{\EM}[\psis{1}] \Big)^{\frac{1}{2}},\\
\lab{eq:esti:NNtmora:psisp:2:nocontrolonangularderis}
&\NNtmora[\psis{2},\chi_{\tau_1, \tau_2}^{(1)}\L^{(2)}_{s}]
+\NNtdemora[\psis{2},\L_{s}^{(2)}](\tau_1+1,\tau_2-3)
\nn\\
\les{}&\bigg(\sum_{p=0,1}\widehat{\EM}[\phis{p}](\tau_1,\tau_2-2)+\sum_{p=0,1}\Extra[\phis{p}](\tau_1,\tau_2-2)+\A[\pmb\phi_s](\tau_1,\tau_2-2)\bigg)^{\frac{1}{2}} \nn\\
&\times\Big(\widetilde{\EM}[\psis{2}]+\EM[\phis{2}](\tau_1,\tau_2-2)
+{\M}_{\de}[\psis{2}](\tau_1+1,\tau_2-3)\Big)^{\frac{1}{2}}\nn\\
&+\widehat{\M}[\phis{1}](\tau_1,\tau_2-2)+\int_{\MM_{r\leq 12m}(\tau_1,\tau_2-2)}|\pr^{\leq 1}\N_{T, s}^{(1)}|^2.
\end{align}
\esub
\end{lemma}

\begin{proof}
  In view of \eqref{eq:definitionofXandEinMtrap:sect8.1.5} and Section  \ref{sec:definitionmicrolocalenergyMorawetznorms}, we have\footnote{As $R_0\geq 20m$ in view of Remark \ref{rmk:choiceofconstantRbymeanvalue}, due to our choice of normalized coordinates, note that we have
\beaa
\pr_{r}^{\text{BL}} &=& \left(-\mu^{-1}+\frac{m^2}{r^2}\right)\pr_\tau+\pr_r-\frac{a}{\De}\pr_{\tphi}=-\pr_\tau+O(r^{-1})\dk\quad\textrm{on}\,\,r\geq R_0.
\eeaa} 
\beaa
X&=&A\pr_{\tau} + \Opw(is_0\mu\xi_r+ib_{\tphi}\xiphi+ib_{\tt}\xit), \quad \textrm{in} \quad \MM_{r_+(1+\dhor'),\Rmic}, \\
 X+w&=&A\pr_{\tau} + (1+O(r^{-1}))\pr_{r}^{\text{BL}} +O(r^{-1})\nn\\
 &=& A\pr_{\tau}  -\pr_\tau+ O(r^{-1})\dk^{\leq 1}
  \quad \textrm{in} \quad \MM_{\Rmic,+\infty},
\eeaa
with the symbols $s_0$, $b_{\tphi}$ and $b_{\tau}$ given in Section \ref{sec:relevantmixedsymbolsonMM}. Hence, by \eqref{def:NNtmoracopy}, we infer, for $p=0,1,2$,
\bea
\lab{eq:NNtmora:bound:by1234<pf}
\NNtmora[\psis{p},\chi_{\tau_1, \tau_2}^{(1)}\L^{(p)}_{s}]\les {\bf{I}}_p+{\bf{II}}_p+{\bf{III}}_p
\eea
with
\beaa
{\bf{I}}_p&:=&\sum_{i,j}\bigg|\int_{\MM_{r_+(1+\dhor'),\Rmic}}\Re\bigg(\chi_{\tau_1, \tau_2}^{(1)}L^{(p)}_{s,ij}\ov{\Opw(is_0\mu\xi_r+ib_{\tphi}\xiphi+ib_{\tt}\xit){\psi}^{(p)}_{s,ij}}\bigg)\bigg|,\\
{\bf{II}}_p&:=&\sum_{i,j}\bigg|\int_{\MM_{r_+(1+\dhor'), \Rmic}(\Iti)}\Re\Big(\pr_{\tau}{\psi}^{(p)}_{s,ij}\chi_{\tau_1, \tau_2}^{(1)}\ov{L^{(p)}_{s,ij}}\Big)\bigg|,\\
{\bf{III}}_p&:=&\sum_{i,j}\bigg|\int_{\MM_{\Rmic, +\infty}(\Iti)}\Re\Big(\big(\pr_\tau, O(r^{-1})\dk^{\leq 1}\big)\psi^{(p)}_{s,ij}\chi_{\tau_1, \tau_2}^{(1)}\ov{L^{(p)}_{s,ij}}\Big)\bigg|,
\eeaa
in which we have dropped the positive constant $A$ which has been fixed. In the following, we will estimate these three terms one by one, and we will constantly use the following fact  
\bea
\lab{eq:supportpropertyofchi1tau1tau2:pf}
\textrm{supp}(\chi_{\tau_1, \tau_2}^{(1)})\subset [\tau_1,\tau_2-2].
\eea

Applying Cauchy-Schwarz and using the estimate \eqref{eq:gardinginequalitiesyieldcontroloftermsbywidetileM:1},  we deduce, for $p=0,1,2$,
\beaa
{\bf{I}}_p&\les& \bigg(\int_{\MM_{r_+(1+\dhor'),\Rmic}}\big|\chi_{\tau_1, \tau_2}^{(1)}\L^{(p)}_{s}\big|^2 \bigg)^{\frac{1}{2}}\Big(\widetilde{\M}[\psis{p}] \Big)^{\frac{1}{2}}\\
&\les& \bigg(\int_{\MM_{r_+(1+\dhor'),\Rmic}(\tau_1,\tau_2-2)}\big|\L^{(p)}_{s}\big|^2 \bigg)^{\frac{1}{2}}\Big(\widetilde{\M}[\psis{p}] \Big)^{\frac{1}{2}}
\eeaa
where we have used \eqref{eq:supportpropertyofchi1tau1tau2:pf} in the last inequality. Now, in view of the form \eqref{eq:expression:Lsijp:copy} of $\L^{(p)}_{s}$, we have 
\bsub\lab{eq:basiccontrolofL2spacetimenormLsponrleqR0tau1tau2minus2}
\begin{align}
\int_{\MM_{r_+(1+\dhor'),\Rmic}(\tau_1,\tau_2-2)}\big|\L^{(0)}_{s}\big|^2 \les{}& \sum_{p=0,1}\A[\phis{p}](\tau_1,\tau_2-2)+\Ao[r^{-\frac{1}{2}}\Xcal_s\phis{0}](\tau_1,\tau_2-2),\\
\nn\int_{\MM_{r_+(1+\dhor'),\Rmic}(\tau_1,\tau_2-2)}\big|\L^{(1)}_{s}\big|^2 \les{}& \Extra[\phis{0}](\tau_1,\tau_2-2)+\sum_{p=1,2}\A[\phis{p}](\tau_1,\tau_2-2)\\
&+\Ao[r^{-\frac{1}{2}}\Xcal_s\phis{1}](\tau_1,\tau_2-2),\\
\int_{\MM_{r_+(1+\dhor'),\Rmic}(\tau_1,\tau_2-2)}\big|\L^{(2)}_{s}\big|^2 \les{}& \Extra[\phis{1}](\tau_1,\tau_2-2)+\sum_{p=0,2}\A[\phis{p}](\tau_1,\tau_2-2),
\end{align}
\esub
which, together with \eqref{eq:formofregularvectorfieldmathcalXs:Kerr:copy:horizontal:consequence}, implies that ${\bf{I}}_p$ is bounded by the RHS of \eqref{eq:esti:NNtmora:psisp:012:nocontrolonangularderis} for $p=0,1,2$.

Next, we consider the integral term ${\bf{III}}_p$. In view of \eqref{eq:supportpropertyofchi1tau1tau2:pf}, and decomposing the integration over $\MM_{r\geq \Rmic}(\tau_1,\tau_1+1)$, $\MM_{r\geq \Rmic}(\tau_1+1,\tau_2-3)$ and $\MM_{r\geq \Rmic}(\tau_2-3,\tau_2-2)$, we infer
\beaa
{\bf{III}}_p&\les& \bigg(\int_{\MM_{r\geq \Rmic}(\tau_1+1,\tau_2-3)}r^{1+\de}\big|\L^{(p)}_{s}\big|^2 \bigg)^{\frac{1}{2}}\Big({\M}_{\de}[\psis{p}](\tau_1+1,\tau_2-3)\Big)^{\frac{1}{2}}\nn\\
&&+\bigg(\int_{\MM_{r\geq \Rmic}(\tau_2-3,\tau_2-2)}r^2\big|\L^{(p)}_{s}\big|^2 \bigg)^{\frac{1}{2}}\bigg(\sup_{\tau\in (\tau_2-3,\tau_2-2)}\E[\psis{p}](\tau)\bigg)^{\frac{1}{2}}\nn\\
&&+\bigg(\int_{\MM_{r\geq \Rmic}(\tau_1,\tau_1+1)}r^2\big|\L^{(p)}_{s}\big|^2 \bigg)^{\frac{1}{2}}\bigg(\sup_{\tau\in (\tau_1,\tau_1+1)}\E[\psis{p}](\tau)\bigg)^{\frac{1}{2}}.
\eeaa
Now, we have, for any $(\tau_0, \tau_0+q)\subset(\tau_1, \tau_2)$,
\bsub
\lab{eq:controlr2L2square:largeradiusregionlocalintimeregion}
\begin{align}
\int_{\MM_{r\geq \Rmic}(\tau_0,\tau_0+q)}r^2\big|\L^{(0)}_{s}\big|^2 \les{}& \Ao[r^{-\frac{1}{2}}\Xcal_s\phis{0}](\tau_0,\tau_0+q)+\sum_{p=0,1}\A[\phis{p}](\tau_0,\tau_0+q),\\
\int_{\MM_{r\geq \Rmic}(\tau_0,\tau_0+q)}r^2\big|\L^{(1)}_{s}\big|^2 \les{}& q\sup_{\tau\in(\tau_0,\tau_0+q)}\E[\phis{0}](\tau)
+\Ao[r^{-\frac{1}{2}}(\Xcal_s)^{\leq 1}\phis{1}](\tau_0,\tau_0+q)\nn\\
&
+\A[\phis{2}](\tau_0,\tau_0+q),\\
\int_{\MM_{r\geq \Rmic}(\tau_0,\tau_0+q)}r^2\big|\L^{(2)}_{s}\big|^2 \les{}& \sum_{p=0,1}q\sup_{\tau\in(\tau_0,\tau_0+q)}\E[\phis{p}](\tau)+\A[\phis{2}](\tau_0,\tau_0+q),
\end{align}
\esub
and we also have
\bsub
\lab{eq:controlr1+deL2square:largeradiusregion}
\begin{align}
\int_{\MM_{r\geq \Rmic}(\tau_1,\tau_2-2)}r^{1+\de}\big|\L^{(0)}_{s}\big|^2 \les{}& \Ao[r^{-\frac{1}{2}}\Xcal_s\phis{0}](\tau_1,\tau_2-2)+\sum_{p=0,1}\A[\phis{p}](\tau_1,\tau_2-2),\\
\int_{\MM_{r\geq \Rmic}(\tau_1,\tau_2-2)}r^{1+\de}\big|\L^{(1)}_{s}\big|^2 \les{}& \Extra[\phis{0}](\tau_1,\tau_2-2)
+\Ao[r^{-\frac{1}{2}}(\Xcal_s)^{\leq 1}\phis{1}](\tau_1,\tau_2-2)\nn\\
&
+\A[\phis{2}](\tau_1,\tau_2-2),\\
\int_{\MM_{r\geq \Rmic}(\tau_1,\tau_2-2)}r^{1+\de}\big|\L^{(2)}_{s}\big|^2 \les{}& \sum_{p=0,1}\Extra[\phis{p}](\tau_1,\tau_2-2)
+\A[\phis{2}](\tau_1,\tau_2-2),
\end{align}
\esub
from which, together with \eqref{eq:formofregularvectorfieldmathcalXs:Kerr:copy:horizontal:consequence}, we infer that ${\bf{III}}_p$ is bounded by the RHS of \eqref{eq:esti:NNtmora:psisp:012:nocontrolonangularderis} for $p=0,1,2$.

Next, we use Cauchy-Schwarz to find
\beaa
{\bf{II}}_p &\les& \bigg(\int_{\MM_{r_+(1+\dhor'),\Rmic}(\tau_1,\tau_2-2)}\big|\L^{(p)}_{s}\big|^2 \bigg)^{\frac{1}{2}}\\
&&\times\left(\sup_{\tau\in[\tau_1, \tau_1+1]\cup[\tau_2-3,\tau_2-2]}\E[\psis{p}](\tau)+\widehat{\M}[\psis{p}](\tau_1+1, \tau_2-3) \right)^{\frac{1}{2}},\qquad p=0,1.
\eeaa
Using the estimates \eqref{eq:basiccontrolofL2spacetimenormLsponrleqR0tau1tau2minus2}, this yields that the term ${\bf{II}}_p$, $p=0,1$, is bounded by the RHS of \eqref{eq:esti:NNtmora:psisp:012:nocontrolonangularderis}.

It remains to show that the term ${\bf{II}}_2$ is bounded by the RHS of \eqref{eq:esti:NNtmora:psisp:2:nocontrolonangularderis}. First, we have
\beaa
{\bf{II}}_2 &\les& \sum_{i,j}\bigg|\int_{\MM_{r_+(1+\dhor'), \Rmic}(\Iti)}\Re\Big(\pr_{\tau}{\phi}^{(2)}_{s,ij}\chi_{\tau_1, \tau_2}^{(1)}\ov{L^{(2)}_{s,ij}}\Big)\bigg|\\
&&+\sum_{i,j}\int_{\MM_{r_+(1+\dhor'), \Rmic}(\tau_1, \tau_2-2)}\chi_{\tau_1, \tau_2}^{(1)}|\pr_\tau({\psi}^{(2)}_{s,ij}-{\phi}^{(2)}_{s,ij})||L^{(2)}_{s,ij}|
\eeaa
which together with \eqref{eq:psispequalsphispintau1+1tau2-1} and \eqref{eq:supportpropertyofchi1tau1tau2:pf} implies
\bea\lab{eq:intermediaryestimateforthecontrolofII2tofocusontermwithonlyphiinsteadofpsitodoIBP}
\nn {\bf{II}}_2 &\les& \sum_{i,j}\bigg|\int_{\MM_{r_+(1+\dhor'), \Rmic}(\Iti)}\Re\Big(\pr_{\tau}{\phi}^{(2)}_{s,ij}\chi_{\tau_1, \tau_2}^{(1)}\ov{L^{(2)}_{s,ij}}\Big)\bigg|\\
\nn&&+\left(\sup_{\tau\in[\tau_1, \tau_1+1]\cup[\tau_2-3,\tau_2-2]}\Big(\E[\psis{2}](\tau)+\E[\phis{2}](\tau)\Big)\right)^{\frac{1}{2}}\\
&&\times\left(\sum_{p=0,1}\sup_{\tau\in[\tau_1, \tau_1+1]\cup[\tau_2-3,\tau_2-2]}\E[\phis{p}](\tau)+\A[\phis{2}](\tau_1,\tau_2-2)\right)^{\frac{1}{2}}.
\eea
Next, we focus on estimating the first term on the RHS of \eqref{eq:intermediaryestimateforthecontrolofII2tofocusontermwithonlyphiinsteadofpsitodoIBP}. To this end, we proceed in the same manner as in \cite{Ma}. Given the form of $L^{(2)}_{s,ij}$ in \eqref{eq:expression:Lsijp:copy}, we have
\beaa
&&\bigg|\int_{\MM_{r_+(1+\dhor'), \Rmic}(\Iti)}\Re\Big(\pr_{\tau}{\phi}^{(2)}_{s,ij}\chi_{\tau_1, \tau_2}^{(1)}\ov{L^{(2)}_{s,ij}}\Big)\bigg|\nn\\
&\les& \bigg|\int_{\MM_{r_+(1+\dhor'), \Rmic}(\Iti)}\Re\Big(\pr_{\tau}{\phi}^{(2)}_{s,ij}\chi_{\tau_1, \tau_2}^{(1)}\ov{h^{(2,1)}(\pr_{\tphi}+a\pr_{\tt})\phi^{(1)}_{s,ij}}\Big)\bigg|\nn\\
&&+\bigg|\int_{\MM_{r_+(1+\dhor'), \Rmic}(\Iti)}\Re\Big(\pr_{\tau}{\phi}^{(2)}_{s,ij}\chi_{\tau_1, \tau_2}^{(1)}\ov{\Big(h^{(2,2)}\phi_{s,ij}^{(2)}\Big)}\Big)\bigg|\nn\\
&&+\bigg|\int_{\MM_{r_+(1+\dhor'), \Rmic}(\Iti)}\Re\bigg(\pr_{\tau}{\phi}^{(2)}_{s,ij}\chi_{\tau_1, \tau_2}^{(1)}\ov{\bigg(\sum_{k,l=1,2,3}O(mr^{-2})\phiss{kl}{1}+O(m^2r^{-2})\phiss{ij}{0}\bigg)}\bigg)\bigg|.
\eeaa
Using integration by parts in $\pr_{\tau}$ for the last two integrals, and in view of the fact that the coefficients $h^{(2,2)}=O(mr^{-3})$, $O(mr^{-2})$ and $O(m^2r^{-2})$ are all real-valued functions, we deduce
\bea
\lab{eq:errortermfromprtau:psi2:wholepart:prestep:proof}
&&\bigg|\int_{\MM_{r_+(1+\dhor'), \Rmic}(\Iti)}\Re\Big(\pr_{\tau}{\phi}^{(2)}_{s,ij}\chi_{\tau_1, \tau_2}^{(1)}\ov{L^{(2)}_{s,ij}}\Big)\bigg|\nn\\
&\les&\Bigg(\sum_{p=0,1}\widehat{\M}[\phis{p}](\tau_1,\tau_2-2)
+\A[\pmb\phi_s](\tau_1, \tau_2-2)\Bigg)^{\frac{1}{2}}\bigg(\M[\phis{2}](\tau_1,\tau_2-2)\bigg)^{\frac{1}{2}}\nn\\
&&
+\bigg|\int_{\MM(\Iti)}\Re\Big(\chi_{\dbl}(r)\pr_{\tau}{\phi}^{(2)}_{s,ij}\chi_{\tau_1, \tau_2}^{(1)}\ov{h^{(2,1)}(\pr_{\tphi}+a\pr_{\tt})\phi^{(1)}_{s,ij}}\Big)\bigg|,
\eea
where $\chi_{\dbl}(r)$ is a smooth cut-off function in $r$ such that
\begin{equation}
\lab{eq:defofchidbl}
\begin{split}
\chi_{\dbl}(r)={}&1 \quad \text{for } r_+(1+3\dbl/2)\leq r\leq 11m, \\
\chi_{\dbl}(r)={}&0 \quad \text{for } r\in[r_+(1-\dhor), r_+(1+\dbl)]\cup [12m,+\infty).
\end{split}
\end{equation}
In view of \eqref{eq:psispequalsphispintau1+1tau2-1} and the fact that $h^{(2,1)}$ is a real-valued function, we have, for the last term in \eqref{eq:errortermfromprtau:psi2:wholepart:prestep:proof},
\bea
\lab{eq:errortermfromprtau:psi2:onepart:proof}
&&\bigg|\int_{\MM(\Iti)}\Re\Big(\chi_{\dbl}(r)\pr_{\tau} {\phi}^{(2)}_{s,ij}\chi_{\tau_1, \tau_2}^{(1)}h^{(2,1)}\ov{{(\pr_{\tphi}+a\pr_{\tt})}\phi^{(1)}_{s,ij}}\Big)\bigg|\nn\\
&\les& \widehat{\M}[\phis{1}](\tau_1,\tau_2-2)+\int_{\MM_{r\leq 12m}(\tau_1,\tau_2-2)}|\pr^{\leq 1}\N_{T, s}^{(1)}|^2\nn\\
&&+\bigg|\int_{\MM(\Iti)}\Re\Big(\chi_{\dbl}(r)\pr_{\tau}Y_s{\phi}^{(1)}_{s,ij}\chi_{\tau_1, \tau_2}^{(1)}\tilde{h}^{(2,1)}\ov{{(\pr_{\tphi}+a\pr_{\tt})}\phi^{(1)}_{s,ij}}\Big)\bigg|,
\eea
with $\tilde{h}^{(2,1)}$ a real-valued scalar function and 
with $Y_s$ a first-order differential operator defined by
\bea
\lab{def:Ysdifferentialoperator}
Y_{+2}:=e_3, \qquad Y_{-2}:=\frac{|q|^2}{\De}e_4,
\eea
 where in the last step we have used 
the transport equation\footnote{In the case $s=-2$, the use of the transport equation \eqref{eq:ScalarizedQuantitiesinTeuSystem:Kerrperturbation} generates a factor $\De^{-1}$, and in turn a factor $\dbl^{-1}$ on the support of $\chi_{\dbl}(r)$. This factor of $\dbl^{-1}$ is absorbed in $\les$ in view of our convention in Section \ref{sec:smallnesconstants}.}  \eqref{eq:ScalarizedQuantitiesinTeuSystem:Kerrperturbation} for $p=1$ on the support of $\chi_{\dbl}(r)\chi_{\tau_1, \tau_2}^{(1)}(\tau)$. For the last term in \eqref{eq:errortermfromprtau:psi2:onepart:proof}, we utilize the following identity for any scalar field $\varphi$ and any real-valued function $f$
\begin{align}
\lab{eq:firstordertimessecondorder:reducetolowerorder}
2\Re( f\pr_{\a}\varphi \ov{\pr_{\b}\pr_{\ga}\varphi})={}&\Re\Big(\pr_{\b}(f\pr_{\a}\varphi \ov{\pr_{\ga}\varphi})-\pr_{\a}(f\pr_{\b}\varphi \ov{\pr_{\ga}\varphi})
+\pr_{\ga}(f\pr_{\b}\varphi \ov{\pr_{\a}\varphi})\nn\\
&-\pr_{\b}(f)\pr_{\a}\varphi \ov{\pr_{\ga}\varphi}+\pr_{\a}(f)\pr_{\b}\varphi \ov{\pr_{\ga}\varphi}
-\pr_{\ga}(f)\pr_{\b}\varphi \ov{\pr_{\a}\varphi}\Big)
\end{align}
to deduce\footnote{In the case $s=-2$, the factor $\De^{-1}$ in the definition of $Y_{-2}$ generates a factor $\dbl^{-1}$ on the support of $\chi_{\dbl}(r)$. This factor of $\dbl^{-1}$ is absorbed in $\les$ in view of our convention in Section \ref{sec:smallnesconstants}.}
\beaa
\sum_{i,j}\bigg|\int_{\MM(\Iti)}\Re\Big(\chi_{\dbl}(r)\pr_{\tau}Y_s{\phi}^{(1)}_{s,ij}\chi_{\tau_1, \tau_2}^{(1)}\tilde{h}^{(2,1)}\ov{{(\pr_{\tphi}+a\pr_{\tt})}(\phi^{(1)}_{s,ij})}\Big)\bigg|\les \widehat{\M}[\phis{1}](\tau_1,\tau_2-2),
\eeaa
which together with \eqref{eq:intermediaryestimateforthecontrolofII2tofocusontermwithonlyphiinsteadofpsitodoIBP},  \eqref{eq:errortermfromprtau:psi2:wholepart:prestep:proof} and \eqref{eq:errortermfromprtau:psi2:onepart:proof}  implies 
\beaa
{\bf{II}}_2&=&\sum_{i,j}\bigg|\int_{\MM_{r_+(1+\dhor'), \Rmic}(\Iti)}\Re\Big(\pr_{\tau}{\psi}^{(2)}_{s,ij}\chi_{\tau_1, \tau_2}^{(1)}\ov{L^{(2)}_{s,ij}}\Big)\bigg|\nn\\
&\les& \Bigg(\sum_{p=0,1}\widehat{\EM}[\phis{p}](\tau_1,\tau_2-2)
+\A[\pmb\phi_s](\tau_1, \tau_2-2)\Bigg)^{\frac{1}{2}}\\
&&\times\bigg(\EM[\psis{2}](\tau_1,\tau_2-2)+\EM[\phis{2}](\tau_1,\tau_2-2)\bigg)^{\frac{1}{2}}\nn\\
&&+\widehat{\M}[\phis{1}](\tau_1,\tau_2-2)+\int_{\MM_{r\leq 12m}(\tau_1,\tau_2-2)}|\pr^{\leq 1}\N_{T, s}^{(1)}|^2
\eeaa
which is in turn bounded by the RHS of inequality \eqref{eq:esti:NNtmora:psisp:2:nocontrolonangularderis}. 

In view of the bound \eqref{eq:NNtmora:bound:by1234<pf} and the above estimates for the terms ${\bf{I}}_p$, ${\bf{II}}_p$ and  ${\bf{III}}_p$, this concludes the proof of the estimates \eqref{eq:esti:NNtmora:psisp:012:nocontrolonangularderis} for the terms $\NNtmora[\psis{p},\chi_{\tau_1, \tau_2}^{(1)}\L^{(p)}_{s}]$, $p=0,1,2$.

The estimates \eqref{eq:esti:NNtmora:psisp:0:nocontrolonangularderis}, \eqref{eq:esti:NNtmora:psisp:1:nocontrolonangularderis} and \eqref{eq:esti:NNtmora:psisp:2:nocontrolonangularderis} for the terms $\NNtdemora[\psis{p},\L_{s}^{(p)}](\tau_1+1,\tau_2-1)$ with $p=0,1,2$ follow in the same manner as the proof of the above bounds for the term ${\bf{III}}_p$. This concludes the proof of 
Lemma \ref{lem:NNtmora:scalarizedTeukolsky:estimates}.
\end{proof}

Next, we estimate in the following lemma the error terms $\NNtener[\psis{p},\chi_{\tau_1, \tau_2}^{(1)}\L^{(p)}_{s}]$. 

\begin{lemma}[Bounds for the error terms $\NNtener$ coming from linear coupling terms]
\lab{lem:NNtener:scalarizedTeukolsky:estimates}
For the terms $\NNtener[\psis{p},\chi_{\tau_1, \tau_2}^{(1)}\L^{(p)}_{s}]$,  $s=\pm 2$, $p=0,1,2$, we have 
\bsub
\lab{eq:esti:NNtener:psisp:012:nocontrolonangularderis}
\begin{align}
\lab{eq:esti:NNtener:psisp:0:nocontrolonangularderis}
&\NNtener[\psis{0},\chi_{\tau_1, \tau_2}^{(1)}\L^{(0)}_{s}]\nn\\
\les {}&\left(\sum_{p=0,1}\A[\phis{p}](\tau_1,\tau_2-2)+\Ao[r\nab\phis{0}](\tau_1,\tau_2-2) +\ep\widehat{\M}[\phis{0}](\tau_1,\tau_2-2)\right)^{\frac{1}{2}}
\nn\\
&\times\Big(\widehat{\M}_{\de}[\psis{0}](\tau_1+1,\tau_2-3)+\widetilde{\EM}[\psis{0}] \Big)^{\frac{1}{2}},\\
\lab{eq:esti:NNtener:psisp:1:nocontrolonangularderis}
&\NNtener[\psis{1},\chi_{\tau_1, \tau_2}^{(1)}\L^{(1)}_{s}]\nn\\
\les{}& \Bigg(\Extra[\phis{0}](\tau_1,\tau_2-2)+\sum_{p=1,2}\A[\phis{p}](\tau_1,\tau_2-2)
+\Ao[r\nab\phis{1}](\tau_1,\tau_2-2)\nn\\
& +\ep\widehat{\M}[\phis{1}](\tau_1,\tau_2-2) +\sup_{\tau\in[\tau_1, \tau_1+1]\cup[\tau_2-3,\tau_2-2]}\E[\phis{0}](\tau)\Bigg)^{\frac{1}{2}}\nn\\
&\times \Big(\widehat{\M}_{\de}[\psis{1}](\tau_1+1,\tau_2-3)+\widetilde{\EM}[\psis{1}] \Big)^{\frac{1}{2}},\\
\lab{eq:esti:NNtener:psisp:2:nocontrolonangularderis}
&\NNtener[\psis{2},\chi_{\tau_1, \tau_2}^{(1)}\L^{(2)}_{s}]\nn\\
\les{}&\bigg(\sum_{p=0,1}\widehat{\EM}[\phis{p}](\tau_1,\tau_2-2)+\sum_{p=0,1}\Extra[\phis{p}](\tau_1,\tau_2-2)+\A[\pmb\phi_s](\tau_1,\tau_2-2)\bigg)^{\frac{1}{2}} \nn\\
&\times\Big(\widetilde{\EM}[\psis{2}]+\EM[\phis{2}](\tau_1,\tau_2-2)
+{\M}_{\de}[\psis{2}](\tau_1+1,\tau_2-3)\Big)^{\frac{1}{2}}\nn\\
&+\widehat{\M}[\phis{1}](\tau_1,\tau_2-2)+\int_{\MM_{r\leq 12m}(\tau_1,\tau_2-2)}|\pr^{\leq 1}\N_{T, s}^{(1)}|^2.
\end{align}
\esub
\end{lemma}

\begin{proof}
For convenience, we may assume that the smooth cut-off function $\chi^{(1)}_{\tau_1, \tau_2}$ satisfying  \eqref{eq:propertieschi:thisonetodealwithRHSofTeukolsky} has been chosen such that $(\chi^{(1)}_{\tau_1, \tau_2})^{\frac{1}{4}}$ satisfies the same properties, i.e., 
\beaa
\chi^{(1)}_{\tau_1, \tau_2}=(\chi^{(2)}_{\tau_1, \tau_2})^4
\eeaa 
where $\chi^{(2)}_{\tau_1, \tau_2}=\chi^{(2)}_{\tau_1, \tau_2}(\tau)$ is a smooth cut-off function satisfying 
\bea\lab{eq:propertieschi:thisonetodealwithRHSofTeukolsky:forchi2=chi11over4}
\chi^{(2)}_{\tau_1, \tau_2}(\tau)=0\,\,\,\textrm{on}\,\,\mathbb{R}\setminus(\tau_1, \tau_2-2), \,\,\,\chi^{(2)}_{\tau_1, \tau_2}(\tau)=1\,\,\textrm{on}\,\, [\tau_1+1, \tau_2-3], \,\,\, \|\chi^{(2)}_{\tau_1, \tau_2}\|_{W^{{15},+\infty}(\mathbb{R})}\les 1.
\eea
In particular, to prove the lemma, it thus suffices to control $\NNtener[\psis{p}, (\chi_{\tau_1, \tau_2}^{(2)})^4\L^{(p)}_{s}]$.

In view of the formula of $\NNtener[\c, \c]$  in \eqref{def:NNtenercopy},  we have 
\begin{align}
\lab{eq:NNtener:scalarized:Teu:firststep:control:pf}
&\NNtener[\psis{p}, (\chi_{\tau_1, \tau_2}^{(2)})^4\L^{(p)}_{s}]\nn\\
\les{}&\sup_{\tau\geq \tmic}\sum_{i,j}\bigg|\int_{\Mntrap(\tmic, \tau)}{\Re\Big((\chi_{\tau_1, \tau_2}^{(2)})^4 L^{(p)}_{s,ij}\ov{\pr_{\tau}\psi_{s,ij}^{(p)}}\Big)}\bigg|
\nn\\
&+
\sup_{\tau\in\Reals}\sum_{n=-1}^{\iota}\sum_{i,j}\bigg|\int_{\Mtrap(\tmic, \tau)}\Re\Big(\ov{|q|^{-2}\Opw(\Theta_n)\big(\qs (\chi_{\tau_1, \tau_2}^{(2)})^4L^{(p)}_{s,ij}\big)}V_{n}\Opw(\Theta_{n})\psi_{s,ij}^{(p)}\Big)\bigg|.
\end{align}
In a similar manner as estimating the term ${\bf{III}_p}$ in the proof of Lemma  \ref{lem:NNtmora:scalarizedTeukolsky:estimates}, we infer that  for $p=0,1,2$, the first term on the RHS of  \eqref{eq:NNtener:scalarized:Teu:firststep:control:pf}
is bounded by the RHS of \eqref{eq:esti:NNtmora:psisp:012:nocontrolonangularderis}, and hence by the RHS of \eqref{eq:esti:NNtener:psisp:012:nocontrolonangularderis} since the RHS of \eqref{eq:esti:NNtener:psisp:012:nocontrolonangularderis} coincides with the RHS of \eqref{eq:esti:NNtmora:psisp:012:nocontrolonangularderis}. Therefore, it remains to control the last term  on the RHS of \eqref{eq:NNtener:scalarized:Teu:firststep:control:pf} by the RHS of \eqref{eq:esti:NNtener:psisp:012:nocontrolonangularderis} for $p=0,1,2$ respectively.

For the last term on the RHS of \eqref{eq:NNtener:scalarized:Teu:firststep:control:pf}, we have, for each $n=-1,0,1,2,\ldots, \iota$,  and $i,j=1,2,3$, 
\bea
\lab{eq:lastterminNNtener:precontrol}
&&\sup_{\tau\in\Reals}\bigg|\int_{\Mtrap(\tmic, \tau)}\Re\Big(\ov{|q|^{-2}\Opw(\Theta_{n})\big(\qs (\chi_{\tau_1, \tau_2}^{(2)})^4L^{(p)}_{s,ij}\big)}V_{n}\Opw(\Theta_{n})\psi_{s,ij}^{(p)}\Big)\bigg|\nn\\
&\les&\sup_{\tau\in\Reals}\bigg|\int_{\Mtrap(\tmic, \tau)}\Re\Big(\ov{\Opw(\widetilde{S}^{-1,0}(\MM))(\chi_{\tau_1, \tau_2}^{(2)})^2L^{(p)}_{s,ij}}V_{n}\Opw(\Theta_{n})\psi_{s,ij}^{(p)}\Big)\bigg|\nn\\
&&
+\sup_{\tau\in\Reals}\bigg|\int_{\Mtrap(\tmic, \tau)}\Re\Big(\ov{\Opw(\Theta_{n})(\chi_{\tau_1, \tau_2}^{(2)})^2L^{(p)}_{s,ij}}\Opw(\widetilde{S}^{0,0}(\MM))\psi_{s,ij}^{(p)}\Big)\bigg|\nn\\
&&
+\sup_{\tau\in\Reals}\bigg|\int_{\Mtrap(\tmic, \tau)}\Re\Big(\ov{\Opw(\Theta_{n})(\chi_{\tau_1, \tau_2}^{(2)})^2L^{(p)}_{s,ij}}V_{n}\Opw(\Theta_{n})\big((\chi_{\tau_1, \tau_2}^{(2)})^2\psi_{s,ij}^{(p)}\big)\Big)\bigg|\nn\\
&\les&{\bf{I}_{s}^{(p)}} + {\bf{II}_{s}^{(p)}} +\bigg(\int_{\Mtrap}|\psis{p}|^2\bigg)^{\frac{1}{2}}\bigg(\int_{\MM(\tau_1,\tau_2-2)}|\L^{(p)}_{s}|^2\bigg)^{\frac{1}{2}},
\eea
where we have used in the first step
\beaa
[\Opw(\Theta_{n}), \qs (\chi_{\tau_1, \tau_2}^{(2)})^2]=\Opw(\widetilde{S}^{-1,0}(\MM)),\qquad
[(\chi_{\tau_1, \tau_2}^{(2)})^2, V_{n}\Opw(\Theta_{n})]=\Opw(\widetilde{S}^{0,0}(\MM))
\eeaa
and where we have introduced in the last step ${\bf{I}_{s}^{(p)}}$ and ${\bf{II}_{s}^{(p)}}$ defined by
\bsub
\begin{align}
\lab{expressionofIsp:estimateNNtener:error}
{\bf{I}_{s}^{(p)}}:={}&\sup_{\tau\in\Reals}\bigg|\int_{\Mtrap(-\infty, \tau)}\Re\Big(\ov{\Opw(\widetilde{S}^{-1,0}(\MM))(\chi_{\tau_1, \tau_2}^{(2)})^2L^{(p)}_{s,ij}}V_{n}\Opw(\Theta_{n})\psi_{s,ij}^{(p)}\Big)\bigg|,\\
\lab{expressionofIIsp:estimateNNtener:error}
{\bf{II}_{s}^{(p)}}:={}&\sup_{\tau\in\Reals}\bigg|\int_{\Mtrap(-\infty, \tau)}\Re\Big(\ov{\Opw(\Theta_{n})(\chi_{\tau_1, \tau_2}^{(2)})^2L^{(p)}_{s,ij}}V_{n}\Opw(\Theta_{n})\big((\chi_{\tau_1, \tau_2}^{(2)})^2\psi_{s,ij}^{(p)}\big)\Big)\bigg|.
\end{align}
\esub

In order to estimate the term ${\bf{I}_{s}^{(p)}}$, we introduce the smooth cut-off functions $\chi_{\tau,j}$, $j=0,1$, on $\Reals$ as defined in  \eqref{def:variouscutofffunctionsfortau'}, i.e., 
\begin{equation}
\lab{def:variouscutofffunctionsfortau':tautau'case}
\begin{split}
&\textrm{supp}(\chi_{\tau,0})\subset(-\infty, \tau+1), \quad \chi_{\tau,0}=1\,\,\,\textrm{on}\,\,(-\infty, \tau), \quad 0\leq\chi_{\tau,j}\leq 1, \,\,j=0,1,\\
&\textrm{supp}(\chi_{\tau,1})\subset(\tau-1, \tau+2), \quad \chi_{\tau,1}=1\,\,\,\textrm{on}\,\,(\tau, \tau+1).
\end{split}
\end{equation}
Then, we have
\bea
\lab{eq:estimatesofIsp:inNNtener:errorterms}
{\bf{I}_{s}^{(p)}}
&\les&\sup_{\tau\in\Reals}\bigg(\bigg|\int_{\Mtrap}\Re\Big(\chi_{\tau,0}\ov{\Opw(\widetilde{S}^{-1,0}(\MM))(\chi_{\tau_1, \tau_2}^{(2)})^2L^{(p)}_{s,ij}}V_{n}\Opw(\Theta_{n})\psi_{s,ij}^{(p)}\Big)\bigg|\nn\\
&&+\int_{\Mtrap}\chi_{\tau,1}\left|{\Opw(\widetilde{S}^{-1,0}(\MM))(\chi_{\tau_1, \tau_2}^{(2)})^2L^{(p)}_{s,ij}}\right| \left|V_{n}\Opw(\Theta_{n})\psi_{s,ij}^{(p)}\right|\bigg)\nn\\
&\les&\sup_{\tau\in\Reals}\bigg|\int_{\Mtrap}\Re\Big(\ov{V_{n}\Big(\chi_{\tau,0}\Opw(\widetilde{S}^{-1,0}(\MM))(\chi_{\tau_1, \tau_2}^{(2)})^2L^{(p)}_{s,ij}\Big)}\Opw(\Theta_{n})\psi_{s,ij}^{(p)}\Big)\bigg|\nn\\
&&
+\bigg(\int_{\MM(\tau_1,\tau_2-2)}|\L^{(p)}_{s}|^2\bigg)^{\frac{1}{2}}\bigg(\EM[\psis{p}](\Iti)\bigg)^{\frac{1}{2}}\nn\\
&\les&\bigg(\int_{\MM(\tau_1,\tau_2-2)}|\L^{(p)}_{s}|^2\bigg)^{\frac{1}{2}} \Big(\EM[\psis{p}](\Iti)\Big)^{\frac{1}{2}},
\eea
where,  in the second step, we have used integration by parts in $V_{n}$ and the size of the support of $\chi_{\tau,1}$, and where, in the last step, we have used $V_{n}(\chi_{\tau,0}\Opw(\widetilde{S}^{-1,0}(\MM)))=\Opw(\widetilde{S}^{0,0}(\MM))$. This hence implies that for the last term on the RHS of  \eqref{eq:NNtener:scalarized:Teu:firststep:control:pf}, we have
\bea\lab{eq:auxiliaryinequlalitywhereeverythingisfineexpcetthetermIIsp=2whichremainstobedone:000}
&&\sup_{\tau\in\Reals}\bigg|\int_{\Mtrap(\tmic, \tau)}\Re\Big(\ov{|q|^{-2}\Opw(\Theta_{n})\big(\qs (\chi_{\tau_1, \tau_2}^{(2)})^4L^{(p)}_{s,ij}\big)}V_{n}\Opw(\Theta_{n})\psi_{s,ij}^{(p)}\Big)\bigg|\nn\\
&\les& {\bf{II}_{s}^{(p)}} +\bigg(\int_{\MM(\tau_1,\tau_2-2)}|\L^{(p)}_{s}|^2\bigg)^{\frac{1}{2}}\Big(\EM[\psis{p}](\Iti)\Big)^{\frac{1}{2}}.
\eea
Also, using Cauchy-Schwarz, we have
\beaa
{\bf{II}_{s}^{(p)}} &\les& \bigg(\int_{\MM(\tau_1,\tau_2-2)}\big|\L^{(p)}_{s}\big|^2 \bigg)^{\frac{1}{2}}\\
&&\times\left(\sup_{\tau\in[\tau_1, \tau_1+1]\cup[\tau_2-3,\tau_2-2]}\E[\psis{p}](\tau)+\widehat{\M}[\psis{p}](\tau_1+1, \tau_2-3) \right)^{\frac{1}{2}},\qquad p=0,1.
\eeaa
Together with the estimates \eqref{eq:basiccontrolofL2spacetimenormLsponrleqR0tau1tau2minus2} \eqref{eq:controlr1+deL2square:largeradiusregion} for the linear coupling terms $L^{(p)}_{s,ij}$, and \eqref{eq:formofregularvectorfieldmathcalXs:Kerr:copy:horizontal:consequence}, we infer that all the terms on the RHS of \eqref{eq:auxiliaryinequlalitywhereeverythingisfineexpcetthetermIIsp=2whichremainstobedone:000}, except ${\bf{II}_{s}^{(2)}}$, are bounded by the RHS of \eqref{eq:esti:NNtener:psisp:012:nocontrolonangularderis}.

It remains to show that the term ${\bf{II}_{s}^{(2)}}$ is bounded by the RHS of \eqref{eq:esti:NNtener:psisp:2:nocontrolonangularderis}. First, we have
\beaa
{\bf{II}_{s}^{(2)}} &\les& \sup_{\tau\in\Reals}\bigg|\int_{\Mtrap(-\infty, \tau)}\Re\Big(\ov{\Opw(\Theta_{n})(\chi_{\tau_1, \tau_2}^{(2)})^2L^{(2)}_{s,ij}}V_{n}\Opw(\Theta_{n})\big((\chi_{\tau_1, \tau_2}^{(2)})^2\phi_{s,ij}^{(2)}\big)\Big)\bigg|\\
&&+\int_{\Mtrap}|\Opw(\Theta_{n})(\chi_{\tau_1, \tau_2}^{(2)})^2L^{(2)}_{s,ij}||V_{n}\Opw(\Theta_{n})\big((\chi_{\tau_1, \tau_2}^{(2)})^2(\psi_{s,ij}^{(2)}-\phi_{s,ij}^{(2)}\big)|
\eeaa
which together with \eqref{eq:psispequalsphispintau1+1tau2-1},  \eqref{eq:propertieschi:thisonetodealwithRHSofTeukolsky:forchi2=chi11over4} and the form \eqref{eq:expression:Lsijp:copy} of the linear coupling terms $L^{(2)}_{s,ij}$ implies
\bea\lab{eq:intermediaryestimateforthecontrolofII2tofocusontermwithonlyphiinsteadofpsitodoIBP}
\nn {\bf{II}_{s}^{(p)}} &\les& \sup_{\tau\in\Reals}\bigg|\int_{\Mtrap(-\infty, \tau)}\Re\Big(\ov{\Opw(\Theta_{n})(\chi_{\tau_1, \tau_2}^{(2)})^2L^{(2)}_{s,ij}}V_{n}\Opw(\Theta_{n})\big((\chi_{\tau_1, \tau_2}^{(2)})^2\phi_{s,ij}^{(2)}\big)\Big)\bigg|\\
\nn&&+\left(\sup_{\tau\in[\tau_1, \tau_1+1]\cup[\tau_2-3,\tau_2-2]}\Big(\E[\psis{2}](\tau)+\E[\phis{2}](\tau)\Big)\right)^{\frac{1}{2}}\\
&&\times\left(\sum_{p=0,1}\widehat{\EM}[\phis{p}](\tau_1,\tau_2-2)+\A[\phis{2}](\tau_1,\tau_2-2)\right)^{\frac{1}{2}}.
\eea
Substituting the form  \eqref{eq:expression:Lsijp:copy} of the linear coupling terms $L^{(2)}_{s,ij}$  into the first term on the RHS of \eqref{eq:intermediaryestimateforthecontrolofII2tofocusontermwithonlyphiinsteadofpsitodoIBP}, we deduce
\bea
\lab{eq:NNtener:scalarized:Teu:firststep:microlocalterm}
{\bf{II}_{s}^{(2)}}
&\les&
{\bf{II}_{s,1}^{(2)}}+{\bf{II}_{s,2}^{(2)}} +\left(\EM[\psis{2}](\tau_1, \tau_2-2)+\EM[\phis{2}](\tau_1, \tau_2-2)\right)^{\frac{1}{2}}\nn\\
&&\times\left(\sum_{p=0,1}\widehat{\EM}[\phis{p}](\tau_1,\tau_2-2)+\A[\phis{2}](\tau_1,\tau_2-2)\right)^{\frac{1}{2}},
\eea
where we have defined
\begin{align*}
{\bf{II}_{s,1}^{(2)}}:=&\sup_{\tau\in\Reals}\bigg|\int_{\Mtrap(-\infty, \tau)}\Re\Big(\ov{\Opw(\Theta_{n})(\chi_{\tau_1, \tau_2}^{(2)})^2h^{(2,2)}\phi^{(2)}_{s,ij}}V_{n}\Opw(\Theta_{n})(\chi_{\tau_1, \tau_2}^{(2)})^2\phi_{s,ij}^{(2)}\Big)\bigg|,\nn\\
{\bf{II}_{s,2}^{(2)}}:=&\sup_{\tau\in\Reals}\bigg|\int_{\Mtrap(-\infty, \tau)}\Re\Big(\ov{\Opw(\Theta_{n})(\chi_{\tau_1, \tau_2}^{(2)})^2h^{(2,1)}{(\pr_{\tphi}+a\pr_{\tt})}\phi^{(1)}_{s,ij}}V_{n}\Opw(\Theta_{n})(\chi_{\tau_1, \tau_2}^{(2)})^2\phi_{s,ij}^{(2)}\Big)\bigg|
\end{align*}
and where, to estimate the terms arising from the part $\sum_{k,l=1,2,3}O(mr^{-2})\phiss{kl}{1}+O(m^2r^{-2})\phiss{ij}{0}$ of $L^{(2)}_{s,ij}$, we have relied on the cut-off functions $\chi_{\tau, j}$, $j=0,1$, as in the proof of  \eqref{eq:estimatesofIsp:inNNtener:errorterms}, and then integrated by parts in $V_n$ in the term containing $\chi_{\tau, 0}$.

For ${\bf{II}_{s,1}^{(2)}}$, we rely again on the cut-off functions $\chi_{\tau, j}$, $j=0,1$, integrate by parts in $V_{n}$ and make use of the fact that $h^{(2,2)}$ is a real-valued function, which implies
\bea
\lab{eq:NNtener:scalarized:Teu:firststep:microlocaltermI}
{\bf{II}_{s,1}^{(2)}}\les \Big(\A[\phis{2}](\tau_1,\tau_2-2)\Big)^{\frac{1}{2}} \Big(\EM[\phis{2}](\tau_1,\tau_2-2)\Big)^{\frac{1}{2}} +\A[\phis{2}](\tau_1,\tau_2-2).
\eea

For ${\bf{II}_{s,2}^{(2)}}$, we rely again on the cut-off functions $\chi_{\tau, j}$, $j=0,1$, and use also the cut-off $\chi_{\dbl}(r)$ introduced in \eqref{eq:defofchidbl} to obtain
\begin{align*}
{\bf{II}_{s,2}^{(2)}}
\les&\sup_{\tau\in\Reals}\bigg|\int_{\MM}\Re\Big(\chi_{\dbl}(r)\chi_{\tau,0}h^{(2,1)}\ov{\Opw(\Theta_{n})(\chi_{\tau_1, \tau_2}^{(2)})^2{(\pr_{\tphi}+a\pr_{\tt})}\phi^{(1)}_{s,ij}}\nn\\
&\times V_{n}\Opw(\Theta_{n})(\chi_{\tau_1, \tau_2}^{(2)})^2\phi_{s,ij}^{(2)}\Big)\bigg|+ \Big(\widehat{\M}[\phis{1}](\tau_1,\tau_2-2)\Big)^{\frac{1}{2}}\Big(\EM[\phis{2}](\tau_1,\tau_2-2)\Big)^{\frac{1}{2}}.
\end{align*}
Next, using the transport equation \eqref{eq:ScalarizedQuantitiesinTeuSystem:Kerrperturbation} with $p=1$ to rewrite $\phi_{s,ij}^{(2)}$, we infer\footnote{In the case $s=-2$, the use of the transport equation \eqref{eq:ScalarizedQuantitiesinTeuSystem:Kerrperturbation} generates a factor $\De^{-1}$, and in turn a factor $\dbl^{-1}$ on the support of $\chi_{\dbl}(r)$. This factor of $\dbl^{-1}$ is absorbed in $\les$ in view of our convention in Section \ref{sec:smallnesconstants}.}
\begin{align*}
{\bf{II}_{s,2}^{(2)}}
\les{}& \widehat{\M}[\phis{1}](\tau_1,\tau_2-2)+\int_{\MM_{r\leq 12m}(\tau_1,\tau_2-2)} |\pr^{\leq 1}\N_{T,s}^{(1)}|^2\nn\\
&+ \Big(\widehat{\M}[\phis{1}](\tau_1,\tau_2-2)\Big)^{\frac{1}{2}}\Big(\EM[\phis{2}](\tau_1,\tau_2-2)\Big)^{\frac{1}{2}}
\nn\\
&+\sup_{\tau\in\Reals}\bigg|\int_{\MM}\Re\Big(\chi_{\dbl}(r)\chi_{\tau,0}\tilde{h}^{(2,1)}\ov{\Opw(\Theta_{n})(\chi_{\tau_1, \tau_2}^{(2)})^2{(\pr_{\tphi}+a\pr_{\tt})}\phi^{(1)}_{s,ij}}\nn\\
&\times V_{n}\Opw(\Theta_{n})(\chi_{\tau_1, \tau_2}^{(2)})^2Y_s\phi_{s,ij}^{(1)}\Big)\bigg|,
\end{align*}
where $Y_s$ are defined in \eqref{def:Ysdifferentialoperator} and where $\tilde{h}^{(2,1)}$ is a real-valued function. By  denoting the last term on the RHS of the previous estimate in the following form
\bea
\lab{def:IIs212term:IItener:error}
{\bf{II}_{s,2,1}^{(2)}}&:=&\sup_{\tau\in\Reals}\bigg|\int_{\MM}\Re\Big(\chi_{\dbl}(r)\chi_{\tau,0}\tilde{h}\ov{\Opw(\Theta_{n})(\chi_{\tau_1, \tau_2}^{(2)})^2\pr_{\text{tan},1}\phi^{(1)}_{s,ij}}\nn\\
&&\qquad \qquad \qquad \qquad \qquad\times \pr_{\text{tan},2}\Opw(\Theta_{n})(\chi_{\tau_1, \tau_2}^{(2)})^2\pr\phi_{s,ij}^{(1)}\Big)\bigg|
\eea
up to terms that are bounded by $\widehat{\M}[\phis{1}](\tau_1,\tau_2-2)$, where $\tilde{h}$ is a real-valued function and where $\pr_{\text{tan},i}\in\{\pr_{\tt}, \pr_{x^a}\}$, $i=1,2$, are derivatives tangential to the constant-$r$ hypersurfaces, we infer\footnote{In the case $s=-2$, the factor $\De^{-1}$ in the definition of $Y_{-2}$ generates a factor $\dbl^{-1}$ on the support of $\chi_{\dbl}(r)$. This factor of $\dbl^{-1}$ is absorbed in $\les$ in view of our convention in Section \ref{sec:smallnesconstants}.}
\begin{align}
\lab{eq:IIs22:boundedbyIIs212:NNtener:error}
{\bf{II}_{s,2}^{(2)}}\les{}& \widehat{\M}[\phis{1}](\tau_1,\tau_2-2)+\int_{\MM_{r\leq 12m}(\tau_1,\tau_2-2)} |\pr^{\leq 1}\N_{T,s}^{(1)}|^2\nn\\
&+ \Big(\widehat{\M}[\phis{1}](\tau_1,\tau_2-2)\Big)^{\frac{1}{2}}\Big(\EM[\phis{2}](\tau_1,\tau_2-2)\Big)^{\frac{1}{2}}+{\bf{II}_{s,2,1}^{(2)}}.
\end{align}
The last term ${\bf{II}_{s,2,1}^{(2)}}$,  as defined in \eqref{def:IIs212term:IItener:error}, can be controlled in a similar manner as estimating the last term in \eqref{eq:errortermfromprtau:psi2:onepart:proof} in the proof of Lemma \ref{lem:NNtmora:scalarizedTeukolsky:estimates}, by relying on a variant of the identity \eqref{eq:firstordertimessecondorder:reducetolowerorder}. Specifically, for  real-valued functions $f$ and $\chi$, a scalar $\varphi$ and a PDO $S^0\in \Opw(\widetilde{S}^{(0,0)}(\MM))$ such that $[\pr_r, S^0]=0$, we have the following generalization of \eqref{eq:firstordertimessecondorder:reducetolowerorder}:
\begin{align}\lab{eq:firstordertimessecondorder:reducetolowerorder:secondtimeenergyestmicrolocversion}
&2\Re( fS^0\chi^2\pr_{\text{tan},1}\varphi \ov{\pr_{\text{tan},2}S^0\chi^2\pr_{\ga}\varphi})\nn\\
={}&\Re\Big(\pr_{\text{tan},2}(fS^0\chi^2\pr_{\text{tan},1}\varphi \ov{S^0\chi^2\pr_{\ga}\varphi})-\pr_{\text{tan},1}(fS^0\chi^2\pr_{\text{tan},2}\varphi \ov{S^0\chi^2\pr_{\ga}\varphi})
\nn\\
&+\pr_{\ga}(fS^0\chi^2\pr_{\text{tan},2}\varphi \ov{S^0\chi^2\pr_{\text{tan},1}\varphi})\Big)\nn\\
&+\Re\Big(\big(f, \pr f\big)\Opw(\widetilde{S}^{(0,0)}(\MM))\chi \pr\varphi \ov{\Opw(\widetilde{S}^{(0,0)}(\MM))\chi \pr\varphi}\Big),
\end{align}
where we have used 
\beaa
[\prtan, S^0\chi^2]\prtan\varphi&=& [\prtan, S^0]\chi^2\prtan\varphi + 2S^0\prtan(\chi)\chi\prtan\varphi\nn\\
&=&\Opw(\widetilde{S}^{(0,0)}(\MM))\chi\pr\varphi 
\eeaa
as well as 
\beaa
[\pr_r, S^0\chi^2]\pr\varphi=[\pr_r, S^0]\chi^2\pr\varphi+2S^0\pr_r(\chi)\chi\pr\varphi=\Opw(\widetilde{S}^{(0,0)}(\MM))\chi\pr\varphi
\eeaa
in view of the assumption $[\pr_r, S^0]=0$. Applying \eqref{eq:firstordertimessecondorder:reducetolowerorder:secondtimeenergyestmicrolocversion} with\footnote{Note that $[\pr_r, \Opw(\Th_n)]=0$ in view of \eqref{eq:nodependenceonrforthesymbolesrtrapupsilonTheta} and
 \eqref{eq:propWeylquantization:MM:composition:mixedsymbols:specialcase}.} $S^0=\Opw(\Th_n)$, $f=\chi_{\dbl}(r)\chi_{\tau,0}\tilde{h}$, $\chi=\chi_{\tau_1, \tau_2}^{(2)}$ and $\varphi=\phiss{ij}{1}$ to the integrand of 
${\bf{II}_{s,2,1}^{(2)}}$ in \eqref{def:IIs212term:IItener:error}, we infer
\beaa
\lab{eq:estimatesforIIs212term:onederivativetimestwoderivativesusingIBP}
{\bf{II}_{s,2,1}^{(2)}}\les\widehat{\M}[\phis{1}](\tau_1,\tau_2-2),
\eeaa
which together with \eqref{eq:IIs22:boundedbyIIs212:NNtener:error} implies
\bea
\lab{eq:IIs22:bound:NNtener:error}
{\bf{II}_{s,2}^{(2)}}&\les& \Big(\widehat{\M}[\phis{1}](\tau_1,\tau_2-2)\Big)^{\frac{1}{2}}\Big(\EM[\phis{2}](\tau_1,\tau_2-2)\Big)^{\frac{1}{2}}
\nn\\
&&+\widehat{\M}[\phis{1}](\tau_1,\tau_2-2)+\int_{\MM_{r\leq 12m}(\tau_1,\tau_2-2)} |\pr^{\leq 1}\N_{T,s}^{(1)}|^2.
\eea
Combining this estimate for ${\bf{II}_{s,2}^{(2)}}$ with the estimate \eqref{eq:NNtener:scalarized:Teu:firststep:microlocalterm} for ${\bf{II}_{s}^{(2)}}$ and the estimate \eqref{eq:NNtener:scalarized:Teu:firststep:microlocaltermI} for ${\bf{II}_{s,1}^{(2)}}$, we infer
\beaa
{\bf{II}_{s}^{(2)}}
&\les& \left(\EM[\psis{2}](\tau_1, \tau_2-2)+\EM[\phis{2}](\tau_1, \tau_2-2)\right)^{\frac{1}{2}}\nn\\
&&\times\left(\sum_{p=0,1}\widehat{\EM}[\phis{p}](\tau_1,\tau_2-2)+\A[\phis{2}](\tau_1,\tau_2-2)\right)^{\frac{1}{2}}\\
&&+\widehat{\M}[\phis{1}](\tau_1,\tau_2-2)+\int_{\MM_{r\leq 12m}(\tau_1,\tau_2-2)} |\pr^{\leq 1}\N_{T,s}^{(1)}|^2
\eeaa
which is bounded by the RHS of \eqref{eq:esti:NNtener:psisp:2:nocontrolonangularderis} as desired. This concludes the proof of Lemma \ref{lem:NNtener:scalarizedTeukolsky:estimates}.
\end{proof}

\begin{proof}[Proof of Proposition \ref{prop:naivecontroloferrorfromlinearcoupling:scalarizedTeu}]
In view of \eqref{def:NNtnorms:copy}, combining the estimates in Lemmas \ref{lem:NNtaux:scalarizedTeukolsky:estimates},  \ref{lem:NNtmora:scalarizedTeukolsky:estimates} and \ref{lem:NNtener:scalarizedTeukolsky:estimates}, we deduce the desired bounds \eqref{subeq:naivecontroloferrorfromlinearcoupling:scalarizedTeu} for $\widetilde{\mathcal{N}}[\psis{p},\chi_{\tau_1, \tau_2}^{(1)}\L^{(p)}_{s}]+\NNtdemora[\psis{p},\L_{s}^{(p)}](\tau_1+1,\tau_2-3)$ with $p=0,1,2$. This concludes the proof of Proposition \ref{prop:naivecontroloferrorfromlinearcoupling:scalarizedTeu}.
\end{proof}


\subsection{Improved estimates for $\{\nab\phis{p}\}_{p=0,1}$}
\lab{subsect:improvedMora:nearinfinity:phisp01}


By combining the EMF estimates \eqref{th:eq:mainnondegeenergymorawetzmicrolocal:tensorialwave:scalarized} with the estimates \eqref{subeq:naivecontroloferrorfromlinearcoupling:scalarizedTeu} for the error terms arising from the linear coupling terms, one finds that there is a lower-triangular structure\footnote{Such a lower-triangular structure allows us to first get estimates for $p=0$, then use the estimates for $p=0$ to get estimates for $p=1$, and eventually use the estimates for $p=0,1$ to get the estimates for $p=2$.} in the ensuing energy-Morawetz estimates for $(\psis{p},\phis{p})$ conditional on the control of: 
\begin{itemize}
\item $\A[\pmb\phi_s](\tau_1,\tau_2)$ defined as in \eqref{def:AandAonorms:tau1tau2:phisandpsis},

\item the early time energy norm  $\IE{\pmb\phi_s, \pmb\psi_s}$ defined as in \eqref{eq:definitionofinitialenergyofphiandpsi:IEterm} and the supremum of $\E[\pmb\phi_s^{(p)}]$ on $\tau\in[\tau_1,\tau_1+1]\cup[\tau_2-3, \tau_2-2]$,

\item  $\{\Errdefect[\psis{p}]\}_{p=0,1,2}$ defined as in \eqref{def:Errdefectofpsi}, 

\item $\{\Ao[r\nab\phis{p}](\tau_1,\tau_2)\}_{p=0,1}$ which are defined in \eqref{def:AandAonorms:tau1tau2:phisandpsis} as spacetime integrals of the derivatives $\{\nab\phis{p}\}_{p=0,1}$, 

\item and  $\{\Extra[\phis{p}](\tau_1,\tau_2)\}_{p=0,1}$ which, as defined in \eqref{def:ExtrabulknotcontrolledbystandardMora},  require further control  of the derivatives $\{\nab_{\pr_{\tphi}+ a\pr_{\tt}}^{\leq 1}\phis{p}\}_{p=0,1}$ with $r$-weights near infinity.
\end{itemize}

In this section, we derive improved estimates for $\Ao[r\nab\phis{p}](\tau_1,\tau_2)$, $p=0,1$, as well as for $\Extra[\phis{p}](\tau_1,\tau_2)$, $p=0,1$, which will allow us in Section \ref{subsect:proofofEMFmain:reg=0}  to derive energy-Morawetz estimates for $(\psis{p},\phis{p})$ conditional only on the control of the first three items listed above. These improved estimates are obtained by exploiting the fact that the principal part of the Teukolsky wave equations \eqref{eq:TensorialTeuSys:rescaleRHScontaine2:general:Kerrperturbation} can be rewritten as $\De_2\phis{p}$ plus a null derivative of $\phis{p+1}$ in view of the Teukolsky transport equations \eqref{def:TensorialTeuScalars:wavesystem:Kerrperturbation}. 

\begin{proposition}
\lab{prop:ellipticestimates:phis01:angularderivatives}
For any $0<\de\leq \frac{1}{3}$, we have the following estimates for $\{\nab\phis{p}\}_{s=\pm 2, p=0,1}$ 
\bsub
\lab{eq:improvedestiforangularderiofkerrpert:prop}
\bea
\Ao[r^{\frac{\de}{2}}(r\nab)^{\leq 1}\phis{0}](\tau_1,\tau_2)
&\les &\Big(\widehat{\EMF}[\phis{0}](\tau_1,\tau_2)\Big)^{\frac{1}{2}}\big(\A[\pmb\phi_{s}](\tau_1,\tau_2)\big)^{\frac{1}{2}} +\A[\pmb\phi_{s}](\tau_1,\tau_2)\nn\\
&&
+\int_{\MM(\tau_1,\tau_2)}\Big(r^{-1+\de}|\N_{W,s}^{(0)}|+r^{-2+\de}|\nab_{X_s}\N^{(0)}_{T,s}|\Big)|\phis{0}|,\\
\Ao[r^{\frac{\de}{2}}(r\nab)^{\leq 1}\phis{1}](\tau_1,\tau_2)
&\les &\bigg(\sum_{p=0,1}\widehat{\EMF}[\phis{p}](\tau_1,\tau_2)\bigg)^{\frac{1}{2}}\big(\A[\pmb\phi_{s}](\tau_1,\tau_2)\big)^{\frac{1}{2}} +\A[\pmb\phi_{s}](\tau_1,\tau_2)\nn\\
&&
+\int_{\MM(\tau_1,\tau_2)}\Big(r^{-1+\de}|\N_{W,s}^{(1)}|+r^{-2+\de}|\nab_{X_s}\N^{(1)}_{T,s}|\Big)|\phis{1}|\nn\\
&& +\Ao[r^{\frac{\de}{2}}(r\nab)^{\leq 1}\phis{0}](\tau_1,\tau_2),
\eea
\esub
and the following estimates for $\{\nab_{\pr_{\tphi}+ a\pr_{\tt}}^{\leq 1}\phis{p}\}_{s=\pm2, p=0,1}$
\bea\lab{eq:improvedesti:prtphi+aprtau:ofkerrpert:prop}
\Extra[\phis{p}](\tau_1,\tau_2)\les \Ao[r^{\frac{\de}{2}}(r\nab)^{\leq 1}\phis{p}](\tau_1,\tau_2)+\EM[\phis{p}](\tau_1,\tau_2),
\eea
 where $X_s$ are null vectorfields given by
\bea
\lab{def:vectorfieldXs}
X_{+2}:=e_4, \qquad X_{-2}:=e_3.
\eea
\end{proposition}

\begin{proof}
Recall from Lemma \ref{lemma:expression-wave-operator}  that
\bsub
\bea
\lab{eq:tensorialwave:intermsofnab4nab3:kerrpert}
\squared_2\pmb\psi&=&-\nab_4\nab_3\pmb\psi {+\left(2\om -\frac{1}{2}\tr\chi\right)\nab_3\pmb\psi  -\frac{1}{2}\tr\chib\nab_4\pmb\psi}\nn\\
&&
+\De_2\pmb\psi + {2\etab\c\nab\pmb\psi + 2i (\dual\rho -\eta\wedge\etab)\pmb\psi}+(\Ga_b\c\Ga_g)\pmb\psi,\\
\lab{eq:tensorialwave:intermsofnab3nab4:kerrpert}
\squared_2\pmb\psi&=&-\nab_3\nab_4\pmb\psi {+\left(2\omb -\frac{1}{2}\tr\chib\right)\nab_4\pmb\psi  -\frac{1}{2}\tr\chi\nab_3\pmb\psi}\nn\\
&&
+\De_2\pmb\psi + {2\eta\c\nab\pmb\psi - 2i (\dual\rho -\eta\wedge\etab)\pmb\psi}+(\Ga_b\c\Ga_g)\pmb\psi.
\eea
\esub
For a tensor $\pmb\psi\in\sk_2(\mathbb{C})$ and a complex-valued scalar function $f_{-2}^{(p)}=f_{-2}^{(p)}(r,\cos\th)$, we have
\bea
-\nab_3\nab_4(f_{-2}^{(p)}\pmb\psi)=- e_3(e_4(f_{-2}^{(p)}))\pmb\psi - e_4(f_{-2}^{(p)}) 
\nab_3\pmb\psi-e_3(f_{-2}^{(p)})\nab_4\pmb\psi - f_{-2}^{(p)}\nab_3\nab_4\pmb\psi.
\eea
Choosing
\beaa
\pmb\psi=\phiminus{p}, \qquad f_{-2}^{(p)}=\frac{rq}{\bar{q}}\left(\frac{{r^2}}{|q|^2}\right)^{p-2}=r+2ia\cos\th +O(r^{-1})
\eeaa 
in the above equation, and using
\begin{align*}
e_3(f_{-2}^{(p)})={}&-1+O(r^{-2}) + r\Ga_b, \qquad e_4(f_{-2}^{(p)})=\frac{\De}{\qs} +O(r^{-2}) + \Ga_g, \\
e_3(e_4 (f_{-2}^{(p)}))={}&O(r^{-2})+\dk^{\leq 1}\Ga_g
\end{align*}
which follow from Definitions \ref{def:renormalizationofallnonsmallquantitiesinPGstructurebyKerrvalue} and \ref{definition.Ga_gGa_b}, we obtain
\beaa
-\nab_3\nab_4(f_{-2}^{(p)}\phiminus{p})&=& - f_{-2}^{(p)} \nab_3\nab_4\phiminus{p}+\big(1+O(r^{-2}) + r\Ga_b\big)\nab_4\phiminus{p}\nn\\
&&-\bigg(\frac{\De}{\qs} +O(r^{-2}) + \Ga_g\bigg)\nab_3\phiminus{p}+\big(O(r^{-2})+\dk^{\leq 1}\Ga_g\big)\phiminus{p},
\eeaa
hence, we deduce
\beaa
&&-\nab_3\nab_4\phiminus{p}{+\left(2\omb -\frac{1}{2}\tr\chib\right)\nab_4\phiminus{p}  -\frac{1}{2}\tr\chi\nab_3\phiminus{p}}\nn\\
&=&- \big(r^{-1}+O(r^{-2})\big)\nab_3\nab_4(f_{-2}^{(p)}\phiminus{p})
+\big(O(r^{-2}) +\Ga_b\big)\nab_4\phiminus{p}\nn\\
&&+\big(O(r^{-2}) + \Ga_g\big)\nab_3\phiminus{p}+\big(O(r^{-3})+r^{-1}\dk^{\leq 1}\Ga_g\big)\phiminus{p}
\eeaa
in view of Definitions \ref{def:renormalizationofallnonsmallquantitiesinPGstructurebyKerrvalue} and \ref{definition.Ga_gGa_b}.
Further,  in view of Definitions \ref{def:renormalizationofallnonsmallquantitiesinPGstructurebyKerrvalue} and \ref{definition.Ga_gGa_b}, we have, for $\pmb\psi\in\sk_2(\mathbb{C})$,
\beaa
&&\etab\c\nab\pmb\psi=(O(r^{-2}) +\Ga_g)\c\nab\pmb\psi, \quad \eta\c\nab\pmb\psi=(O(r^{-2}) +\Ga_b)\c\nab\pmb\psi,\\
&&-\dual\rho +\eta\wedge\etab=O(r^{-4}) +r^{-1}\Ga_g + \Ga_b\wedge\Ga_g,
\eeaa
and substituting these into the form \eqref{eq:tensorialwave:intermsofnab3nab4:kerrpert} of the tensorial wave equation, we infer
\bea
\squared_2\phiminus{p}&=&\De_2\phiminus{p}- \big(r^{-1}+O(r^{-2})\big)\nab_3\nab_4(f_{-2}^{(p)}\phiminus{p})
+\big(O(r^{-2}) + \Ga_b\big)\nab_4\phiminus{p}\nn\\
&&+\big(O(r^{-2}) + \Ga_g\big)\nab_3\phiminus{p}
+\big(O(r^{-2}) +\Ga_b\big)\c\nab\phiminus{p}\nn\\
&&+\Big(O(r^{-3}) + r^{-1}\dk^{\leq 1}\Ga_g+  \Ga_b\c\Ga_g\Big)\phiminus{p}.
\eea
Next, in view of the tensorial Teukolsky wave equation \eqref{eq:TensorialTeuSys:rescaleRHScontaine2:general:Kerrperturbation} for $s=-2$, using \eqref{eq:relationsbetweennullframeandcoordinatesframe2:moreprecise} to decompose the $\nab_{\pr_{\tt}}$ derivative term on the LHS of \eqref{eq:TensorialTeuSys:rescaleRHScontaine2:general:Kerrperturbation}, and using the Teukolsky transport equation \eqref{def:TensorialTeuScalars:wavesystem:Kerrperturbation:-2} in the following form
\beaa
\nab_4(f_{-2}^{(p)}\phiminus{p})=f_{-2}^{(p)}\De |q|^{-4} \phiminus{p+1}+\N_{T,-2}^{(p)}, \quad p=0,1,
\eeaa
 we infer
\bea
\lab{eq:Laplacianofphiprewrittenbyderivativeofphip+1:improveddecayforangular:pfofs=-2}
\De_2\phiminus{p} - \frac{4-2\de_{p0}}{\qs}\phiminus{p} &=& \big(r^{-1}+O(r^{-2})\big)\nab_3(|q|^{-4} \De f_{-2}^{(p)} \phiminus{p+1})\nn\\
&&+\big(O(r^{-2}) + \Ga_b\big)\nab_4\phiminus{p}+\big(O(r^{-2}) + \Ga_g\big)\nab_3\phiminus{p}
\nn\\
&&+\big(O(r^{-2}) +\Ga_b\big)\c\nab\phiminus{p}+\big(O(r^{-3}) + r^{-1}\dk^{\leq 1}\Ga_g+  \Ga_b\c\Ga_g\big)\phiminus{p}  \nn\\
&&+ \L_{-2}^{(p)}[\pmb\phi_{-2}] + \N_{W,-2}^{(p)}+O(r^{-1})\nab_3\N_{T,-2}^{(p)}, \quad p=0,1.
\eea
By multiplying on both sides of this equation \eqref{eq:Laplacianofphiprewrittenbyderivativeofphip+1:improveddecayforangular:pfofs=-2} by $r^{-1+\de}\ov{\phiminus{p}}$, taking the real part, integrating over $\MM(\tau_1,\tau_2)$, applying integration by parts in $\nab_3$ for the product between $r^{-1+\de}\ov{\phiminus{p}}$ and the term  in  the first line on the RHS of \eqref{eq:Laplacianofphiprewrittenbyderivativeofphip+1:improveddecayforangular:pfofs=-2}, and in view of the expression of $\L_{-2}^{(p)}[\pmb\phi_{-2}]$ in \eqref{eq:tensor:Lsn:onlye_2present:general:Kerrperturbation}, we deduce, for $0<\de\leq \frac{1}{3}$,
\bsub
\lab{eq:improvedestiforangularderiofkerrpert:proof:s=-2}
\bea
\Ao[r^{\frac{\de}{2}}(r\nab)^{\leq 1}\phiminus{0}](\tau_1,\tau_2)
&\les &\Big(\widehat{\EMF}[\phiminus{0}](\tau_1,\tau_2)\Big)^{\frac{1}{2}}\big(\A[\pmb\phi_{-2}](\tau_1,\tau_2)\big)^{\frac{1}{2}}+\A[\pmb\phi_{-2}](\tau_1,\tau_2)\nn\\
&&+\int_{\MM(\tau_1,\tau_2)}\Big(r^{-1+\de}|\N_{W,-2}^{(0)}|+r^{-2+\de}|\nab_{3}\N^{(0)}_{T,-2}|\Big)|\phiminus{0}|,\\
\Ao[r^{\frac{\de}{2}}(r\nab)^{\leq 1}\phiminus{1}](\tau_1,\tau_2)
&\les& \bigg(\sum_{p=0,1}\widehat{\EMF}[\phiminus{p}](\tau_1,\tau_2)\bigg)^{\frac{1}{2}}\big(\A[\pmb\phi_{-2}](\tau_1,\tau_2)\big)^{\frac{1}{2}}
+\A[\pmb\phi_{-2}](\tau_1,\tau_2)\nn\\
&&+\int_{\MM(\tau_1,\tau_2)}\Big(r^{-1+\de}|\N_{W,-2}^{(1)}|+r^{-2+\de}|\nab_{3}\N^{(1)}_{T,-2}|\Big)|\phiminus{1}|\nn\\
&&+\int_{\MM(\tau_1,\tau_2)} r^{-3+\de} |\nab_{\pr_{\tphi}+a\pr_{\tt}}^{\leq 1}\phiminus{0}| | \phiminus{1}|\nn\\
&\les& \bigg(\sum_{p=0,1}\widehat{\EMF}[\phiminus{p}](\tau_1,\tau_2)\bigg)^{\frac{1}{2}}\big(\A[\pmb\phi_{-2}](\tau_1,\tau_2)\big)^{\frac{1}{2}}
+\A[\pmb\phi_{-2}](\tau_1,\tau_2)\nn\\
&&+\int_{\MM(\tau_1,\tau_2)}\Big(r^{-1+\de}|\N_{W,-2}^{(1)}|+r^{-2+\de}|\nab_{3}\N^{(1)}_{T,-2}|\Big)|\phiminus{1}|\nn\\
&&+\Big(\Ao[r^{\frac{\de}{2}}(r\nab)^{\leq 1}\phiminus{0}](\tau_1,\tau_2)\Ao[r^{\frac{\de}{2}}\phiminus{1}](\tau_1,\tau_2)\Big)^{\frac{1}{2}},
\eea
\esub
where in the last step we have used the estimate \eqref{eq:formofprphi+aprtt:horizontal}. Then, the desired estimate \eqref{eq:improvedestiforangularderiofkerrpert:prop} for $s=-2$ follows from \eqref{eq:improvedestiforangularderiofkerrpert:proof:s=-2} after applying Cauchy-Schwarz to the last term of the second equation.

In the same manner, in view of the Teukolsky wave equation \eqref{eq:TensorialTeuSys:rescaleRHScontaine2:general:Kerrperturbation} for $s=+2$ and using the Teukolsky transport equation \eqref{def:TensorialTeuScalars:wavesystem:Kerrperturbation:+2}, we infer, relying this time on \eqref{eq:tensorialwave:intermsofnab4nab3:kerrpert},
\bea
\lab{eq:Laplacianofphiprewrittenbyderivativeofphip+1:improveddecayforangular:pfofs=+2}
\De_2\phiplus{p} - \frac{4-2\de_{p0}}{\qs}\phiplus{p} &=& \big(r^{-1}+O(r^{-2})\big)\nab_4(|q|^{-2}  f_{+2}^{(p)} \phiplus{p+1})\nn\\
&&+\big(O(r^{-2}) + \Ga_b\big)\nab_4\phiplus{p}+\big(O(r^{-2}) + \Ga_g\big)\nab_3\phiplus{p}
\nn\\
&&+\big(O(r^{-2}) +\Ga_g\big)\c\nab\phiplus{p}+\big(O(r^{-3}) + r^{-1}\dk^{\leq 1}\Ga_g+  \Ga_b\c\Ga_g\big)\phiplus{p}\nn\\
&& + \L_{+2}^{(p)}[\pmb\phi_{+2}]+ \N_{W,+2}^{(p)}+O(r^{-1})\nab_4\N_{T,+2}^{(p)}, \quad p=0,1,
\eea
where we have used Definitions \ref{def:renormalizationofallnonsmallquantitiesinPGstructurebyKerrvalue} and \ref{definition.Ga_gGa_b} and chosen the function $f_{+2}^{(p)}$ as
\beaa
f_{+2}^{(p)}=\frac{r\bar{q}}{q}\left(\frac{r^2}{|q|^2}\right)^{p-2}=r-2ia\cos\th +O(r^{-1}).
\eeaa
By multiplying on both sides of   \eqref{eq:Laplacianofphiprewrittenbyderivativeofphip+1:improveddecayforangular:pfofs=+2} by $r^{-1+\de}\ov{\phiplus{p}}$, taking the real part, integrating over $\MM(\tau_1,\tau_2)$, applying integration by parts in $\nab_4$ for the product with the term in  the first line on the RHS of \eqref{eq:Laplacianofphiprewrittenbyderivativeofphip+1:improveddecayforangular:pfofs=+2}, and in view of the expression of $\L_{+2}^{(p)}[\pmb\phi_{+2}]$ in \eqref{eq:tensor:Lsn:onlye_2present:general:Kerrperturbation}, we deduce the same estimates as \eqref{eq:improvedestiforangularderiofkerrpert:proof:s=-2} with $s=-2$ replaced by $s=+2$ and with $\nab_{3}\N^{(1)}_{T,-2}$ replaced by $\nab_{4}\N^{(1)}_{T,+2}$. This yields the desired estimates \eqref{eq:improvedestiforangularderiofkerrpert:prop} for $s=\pm 2$.

Next, consider the estimates for $\{\nab_{\pr_{\tphi}+ a\pr_{\tt}}^{\leq 1}\phis{p}\}_{s=\pm2, p=0,1}$. Recall from \eqref{eq:formofprphi+aprtt:horizontal} that 
\beaa
\pr_{\tphi} + a\pr_{\tt} =\mathcal{Y} +r^2\Ga_g \dk + O(|a|)\pr_{\tt}, \qquad {\mathcal{Y}}\in\OO(\MM), \qquad |{\mathcal{Y}}|\les r.
\eeaa
Hence, we infer, for $s=\pm2$,  $p=0,1$, and any $0<\de\leq \frac{1}{3}$,
\beaa
\Extra[\phis{p}](\tau_1,\tau_2)
&=&\int_{\MM(\tau_1,\tau_2)}r^{-3+\de}|(\nab_{\pr_{\tphi}+ a\pr_{\tt}})^{\leq 1}\phis{p}|^2\nn\\
&\les&\int_{\MM(\tau_1,\tau_2)}r^{-3+\de}|(r\nab)^{\leq 1}\phis{p}|^2
+\ep\sup_{\tau\in[\tau_1,\tau_2]}\E[\phis{p}](\tau)
+\M[\phis{p}](\tau_1,\tau_2)\nn\\
&\les&\Ao[r^{\frac{\de}{2}}(r\nab)^{\leq 1}\phis{p}](\tau_1,\tau_2)+\EM[\phis{p}](\tau_1,\tau_2),
\eeaa
which yields the desired estimates \eqref{eq:improvedesti:prtphi+aprtau:ofkerrpert:prop}. This concludes the proof of Proposition \ref{prop:ellipticestimates:phis01:angularderivatives}.
\end{proof}

We also show the following higher-order regularity analog of \eqref{eq:improvedestiforangularderiofkerrpert:prop}.

\begin{proposition}
\lab{prop:ellipticestimates:phis01:angularderivatives:highorder}
For all $\reg\leq 14$ and any $0<\de\leq \frac{1}{3}$, we have the following estimates for $\{\nab\phis{p}\}_{s=\pm 2, p=0,1}$ 
\bsub
\lab{eq:improvedestiforangularderiofkerrpert:prop:highorder}
\bea
\lab{eq:improvedestiforangularderiofkerrpert:prop:highorder:p=0}
&&\Ao[{r^{\frac{\de}{2}}(r\nab)^{\leq 1}}\dk^{\leq \reg}\phis{0}](\tau_1,\tau_2)\nn\\
&\les &\Big(\widehat{\EMF}^{(\reg)}_{\de}[\phis{0}](\tau_1,\tau_2)\Big)^{\frac{1}{2}}\big(\A[r^{-\frac{3\de}{2}}\dk^{\leq \reg} \pmb\phi_{s}](\tau_1,\tau_2)\big)^{\frac{1}{2}} +\A[r^{-\frac{3\de}{2}}\dk^{\leq \reg} \pmb\phi_{s}](\tau_1,\tau_2)\nn\\
&&
+\int_{\MM(\tau_1,\tau_2)}\Big(r^{-1+\de}|\dk^{\leq \reg} \N_{W,s}^{(0)}|+r^{-2+\de}|\dk^{\leq \reg} \nab_{X_s}\N^{(0)}_{T,s}|\Big)|\dk^{\leq \reg} \phis{0}|,\\
\lab{eq:improvedestiforangularderiofkerrpert:prop:highorder:p=1}
&&\Ao[r^{\frac{\de}{2}}(r\nab)^{\leq 1}\dk^{\leq \reg} \phis{1}](\tau_1,\tau_2)\nn\\
&\les &\bigg(\sum_{p=0,1}\widehat{\EMF}^{(\reg)}_{\de}[\phis{p}](\tau_1,\tau_2)\bigg)^{\frac{1}{2}}\big(\A[r^{-\frac{3\de}{2}}\dk^{\leq \reg} \pmb\phi_{s}](\tau_1,\tau_2)\big)^{\frac{1}{2}} +\A[r^{-\frac{3\de}{2}}\dk^{\leq \reg} \pmb\phi_{s}](\tau_1,\tau_2)\nn\\
&&
+\int_{\MM(\tau_1,\tau_2)}\Big(r^{-1+\de}|\dk^{\leq \reg} \N_{W,s}^{(1)}|+r^{-2+\de}|\dk^{\leq \reg} \nab_{X_s}\N^{(1)}_{T,s}|\Big)|\dk^{\leq \reg} \phis{1}|\nn\\
&&+\Ao[r^{\frac{\de}{2}}(r\nab)^{\leq 1}\dk^{\leq \reg} \phis{0}](\tau_1,\tau_2),
\eea
\esub
 where $X_s$ are the null vectorfields given in \eqref{def:vectorfieldXs}.
 
Also, for all $\reg\leq 14$, any $0<\de\leq \frac{1}{3}$ and any $R\geq 20m$, we have the following high-order estimates for $\{\nab\phis{p}\}_{s=\pm 2, p=0,1}$ 
\bsub
\lab{eq:improvedestiforangularderiofkerrpert:prop:highorder:largeR}
\bea
\lab{eq:improvedestiforangularderiofkerrpert:prop:highorder:largeR:p=0}
&&\Ao_{r\geq R}[{r^{\frac{\de}{2}}(r\nab)^{\leq 1}}\dk^{\leq \reg}\phis{0}](\tau_1,\tau_2)\nn\\
&\les &\Big({\EMF}^{(\reg)}_{\de, r\geq R/2}[\phis{0}](\tau_1,\tau_2)\Big)^{\frac{1}{2}}\big(\A_{r\geq R/2}[r^{-\frac{3\de}{2}}\dk^{\leq \reg} \pmb\phi_{s}](\tau_1,\tau_2)\big)^{\frac{1}{2}} +\A_{r\geq R/2}[r^{-\frac{3\de}{2}}\dk^{\leq \reg} \pmb\phi_{s}](\tau_1,\tau_2)\nn\\
&&
+\int_{\MM_{r\geq R/2}(\tau_1,\tau_2)}\Big(r^{-1+\de}|\dk^{\leq \reg} \N_{W,s}^{(0)}|+r^{-2+\de}|\dk^{\leq \reg} \nab_{X_s}\N^{(0)}_{T,s}|\Big)|\dk^{\leq \reg} \phis{0}|,\\
\lab{eq:improvedestiforangularderiofkerrpert:prop:highorder:largeR:p=1}
&&\Ao_{r\geq R}[r^{\frac{\de}{2}}(r\nab)^{\leq 1}\dk^{\leq \reg} \phis{1}](\tau_1,\tau_2)\nn\\
&\les &\bigg(\sum_{p=0,1}{\EMF}_{\de, r\geq R/2}^{(\reg)}[\phis{p}](\tau_1,\tau_2)\bigg)^{\frac{1}{2}}\big(\A_{r\geq R/2}[r^{-\frac{3\de}{2}}\dk^{\leq \reg} \pmb\phi_{s}](\tau_1,\tau_2)\big)^{\frac{1}{2}}\nn\\
&& +\A_{r\geq R/2}[r^{-\frac{3\de}{2}}\dk^{\leq \reg} \pmb\phi_{s}](\tau_1,\tau_2)\nn\\
&&
+\int_{\MM_{r\geq R/2}(\tau_1,\tau_2)}\Big(r^{-1+\de}|\dk^{\leq \reg} \N_{W,s}^{(1)}|+r^{-2+\de}|\dk^{\leq \reg} \nab_{X_s}\N^{(1)}_{T,s}|\Big)|\dk^{\leq \reg} \phis{1}|\nn\\
&&+\Ao_{r\geq R/2}[r^{\frac{\de}{2}}(r\nab)^{\leq 1}\dk^{\leq \reg} \phis{0}](\tau_1,\tau_2).
\eea
\esub
\end{proposition}

\begin{proof}
We start with the proof of \eqref{eq:improvedestiforangularderiofkerrpert:prop:highorder}. The proof of the $\reg=0$ cas is analogous to the one of \eqref{eq:improvedestiforangularderiofkerrpert:prop} noticing in addition that for $0<\de\leq \frac{1}{3}$ and horizontal tensors $\pmb\phi_1$ and $\pmb\phi_2$, we have
\bea\lab{eq:usefulgeneralestimatefortermsokforproof:prop:ellipticestimates:phis01:angularderivatives:highorder}
&&\int_{\MM(\tau_1,\tau_2)}r^{-3+\de}|\pmb\phi_1|\Big(|\nab_3\pmb\phi_2|+r^{-1}|\dk^{\leq 1}\pmb\phi_2|+\tau^{-1-\dec}|\dk^{\leq 1}\pmb\phi_2|\Big)\nn\\
&\les& \left(A[r^{-\frac{3\de}{2}}\pmb\phi_1](\tau_1, \tau_2)\right)^{\frac{1}{2}}\left(\widehat{\EMF}_\de[\pmb\phi_2](\tau_1, \tau_2)\right)^{\frac{1}{2}}.
\eea

Next, we focus on proving \eqref{eq:improvedestiforangularderiofkerrpert:prop:highorder} in the case $1\leq\reg\leq 14$. Multiplying on both sides of \eqref{eq:Laplacianofphiprewrittenbyderivativeofphip+1:improveddecayforangular:pfofs=-2} and \eqref{eq:Laplacianofphiprewrittenbyderivativeofphip+1:improveddecayforangular:pfofs=+2} by $\qs$, we deduce, for $s=\pm2, p=0,1$, 
\beaa
\qs \De_2\phis{p}- (4-2\de_{p0})\phis{p} &=& O(1)\nab_{X_s}\phis{p+1} +(O(r^{-1})+\Ga_b)\phis{p+1}\\
&& +O(1)\big(\nab_3\phis{p}, r^{-1}\dk^{\leq 1}\phis{p}\big)+r^2\Ga_g\nab_3\phis{p} +r\Ga_b\dk^{\leq 1}\phis{p}\\
&& +\qs \L_{s}^{(p)}[\pmb\phi_{s}] + \qs \N_{W,s}^{(p)}+O(r)\nab_{X_s}\N_{T,s}^{(p)}.
\eeaa
Differentiating this equation w.r.t. $\dk^{\leq \reg}$, and using the commutators in Corollary \ref{cor:corollaryofLemmacomm} and Definition \ref{definition.Ga_gGa_b}, we inductively show, for $s=\pm2$, $p=0,1$, $1\leq\reg\leq 14$,
\bea
\lab{eq:Laplacianofphiprewrittenbyderivativeofphip+1:pm2:weightedderi}
&&\qs \De_2\dk^{\reg}\phis{p}- (4-2\de_{p0})\dk^{\reg}\phis{p}\nn\\
&=&\big(O(1)+r\dk^{\leq \reg}\Ga_b\big)\nab_{X_s}\dk^{\leq \reg}\phis{p+1}
+\big(O(r^{-1}) + r\dk^{\leq \reg+1}\Ga_g\big)\dk^{\leq \reg}\phis{p+1}\nn\\
&& +O(1)\big(\nab_3\dk^{\leq \reg}\phis{p}, r^{-1}\dk^{\leq k+1}\phis{p}\big) +r^2\dk^{\leq \reg+1}\Ga_g \nab_3\dk^{\leq \reg}\phis{p} +r\dk^{\leq \reg+1}\Ga_b \dk^{\leq \reg+1}\phis{p} \nn\\
&&+O(1)(r\nab)^{\leq 1}\dk^{\leq \reg-1}\phis{p}+ \big(O(r^2)+r^3\dk^{\leq \reg}\Ga_g\big)\dk^{\leq \reg}\L_{s}^{(p)}[\pmb\phi_{s}] + \big(O(r^2)+r^3\dk^{\leq \reg}\Ga_g\big)\dk^{\leq \reg}\N_{W,s}^{(p)}\nn\\
&&+\big(O(r)+r^2\dk^{\leq \reg}\Ga_g\big)\dk^{\leq \reg}\nab_{X_s}\N_{T,s}^{(p)}.
\eea

Multiplying both sides of \eqref{eq:Laplacianofphiprewrittenbyderivativeofphip+1:pm2:weightedderi} by $|q|^{-2}r^{-1+\de}\ov{\dk^{\reg}\phis{p}}$, taking the real part, integrating over $\MM(\tau_1,\tau_2)$, applying integration by parts in $\nab_{X_s}$ for the product between $|q|^{-2}r^{-1+\de}\ov{\dk^{\reg}\phis{p}}$ and the first term  in  the first line on the RHS of \eqref{eq:Laplacianofphiprewrittenbyderivativeofphip+1:pm2:weightedderi}, and in view of the expression of $\L_{s}^{(p)}[\pmb\phi_{s}]$ in \eqref{eq:tensor:Lsn:onlye_2present:general:Kerrperturbation}, we deduce, for $1\leq \reg\leq 14$ and $0<\de\leq \frac{1}{3}$, relying again on the estimate \eqref{eq:usefulgeneralestimatefortermsokforproof:prop:ellipticestimates:phis01:angularderivatives:highorder},
\bsub
\lab{eq:improvedestiforangularderiofkerrpert:pm2:highorder}
\bea
&&{\Ao[r^{\frac{\de}{2}}(r\nab)^{\leq 1}\dk^{\reg}\phis{0}](\tau_1,\tau_2)}\nn\\
&\les &\Big(\widehat{\EMF}^{(\reg)}_{\de}[\phis{0}](\tau_1,\tau_2)\Big)^{\frac{1}{2}}\big(\A[r^{-\frac{3\de}{2}}\dk^{\leq \reg}\pmb\phi_{s}](\tau_1,\tau_2)\big)^{\frac{1}{2}}+\A[r^{-\frac{3\de}{2}}\dk^{\leq \reg}\pmb\phi_{s}](\tau_1,\tau_2)\nn\\
&&+\int_{\MM(\tau_1,\tau_2)}\Big(r^{-1{+\de}}|\dk^{\leq \reg}\N_{W,s}^{(0)}|+r^{-2{+\de}}|\dk^{\leq \reg}\nab_{X_s}\N^{(0)}_{T,s}|\Big)|\dk^{\reg}\phis{0}|\nn\\
&&+\Big(\Ao[r^{\frac{\de}{2}}(r\nab)^{\leq 1}\dk^{\leq \reg-1}\phis{0}](\tau_1,\tau_2)\Big)^{\frac{1}{2}}\Big(\Ao[r^{\frac{\de}{2}}(r\nab)^{\leq 1}\dk^{\leq \reg}\phis{0}](\tau_1,\tau_2)\Big)^{\frac{1}{2}},\\
&&{\Ao[r^{\frac{\de}{2}}(r\nab)^{\leq 1}\dk^{\leq \reg}\phis{1}](\tau_1,\tau_2)}\nn\\
&\les& \bigg(\sum_{p=0,1}\widehat{\EMF}^{(\reg)}_{\de}[\phis{p}](\tau_1,\tau_2)\bigg)^{\frac{1}{2}}\big(\A[r^{-\frac{3\de}{2}}\dk^{\leq \reg}\pmb\phi_{s}](\tau_1,\tau_2)\big)^{\frac{1}{2}}
+\A[r^{-\frac{3\de}{2}}\dk^{\leq \reg}\pmb\phi_{s}](\tau_1,\tau_2)\nn\\
&&+\int_{\MM(\tau_1,\tau_2)}\Big(r^{-1{+\de}}|\dk^{\leq \reg}\N_{W,s}^{(1)}|+r^{-2{+\de}}|\dk^{\leq \reg}\nab_{X_s}\N^{(1)}_{T,s}|\Big)|\dk^{\reg}\phis{1}|\nn\\
&&+\Big(\Ao[r^{\frac{\de}{2}}(r\nab)^{\leq 1}\dk^{\leq \reg-1}\phis{1}](\tau_1,\tau_2)\Big)^{\frac{1}{2}}\Big(\Ao[r^{\frac{\de}{2}}(r\nab)^{\leq 1}\dk^{\leq \reg}\phis{1}](\tau_1,\tau_2)\Big)^{\frac{1}{2}}\nn\\
&&+\int_{\MM(\tau_1,\tau_2)} r^{-3+\de} |\dk^{\leq \reg}\nab_{\pr_{\tphi}+a\pr_{\tt}}^{\leq 1}\phis{0}| | \dk^{\reg} \phis{1}|\nn\\
&\les& \bigg(\sum_{p=0,1}\widehat{\EMF}^{(\reg)}_{\de}[\phis{p}](\tau_1,\tau_2)\bigg)^{\frac{1}{2}}\big(\A[r^{-\frac{3\de}{2}}\dk^{\leq \reg}\pmb\phi_{s}](\tau_1,\tau_2)\big)^{\frac{1}{2}}
+\A[r^{-\frac{3\de}{2}}\dk^{\leq \reg}\pmb\phi_{s}](\tau_1,\tau_2)\nn\\
&&+\int_{\MM(\tau_1,\tau_2)}\Big(r^{-1{+\de}}|\dk^{\leq \reg}\N_{W,s}^{(1)}|+r^{-2{+\de}}|\dk^{\leq \reg}\nab_{X_s}\N^{(1)}_{T,s}|\Big)|\dk^{\reg}\phis{1}|\nn\\
&&+\Big(\Ao[r^{\frac{\de}{2}}(r\nab)^{\leq 1}\dk^{\leq \reg}\phis{0}](\tau_1,\tau_2)+\Ao[r^{\frac{\de}{2}}(r\nab)^{\leq 1}\dk^{\leq \reg-1}\phis{1}](\tau_1,\tau_2)\Big)^{\frac{1}{2}}\nn\\
&&\times\Big(\Ao[r^{\frac{\de}{2}}(r\nab)^{\leq 1}\dk^{\leq \reg}\phis{1}](\tau_1,\tau_2)\Big)^{\frac{1}{2}},
\eea
\esub
where we have used \eqref{eq:formofprphi+aprtt:horizontal} in the last step. The desired estimates \eqref{eq:improvedestiforangularderiofkerrpert:prop:highorder} then follow from the case $\reg=0$, from Cauchy-Schwarz, and from summing up the estimates \eqref{eq:improvedestiforangularderiofkerrpert:pm2:highorder} over $1\leq\reg\leq 14$. 

To prove the estimates \eqref{eq:improvedestiforangularderiofkerrpert:prop:highorder:largeR}, we multiply both sides of equation \eqref{eq:Laplacianofphiprewrittenbyderivativeofphip+1:pm2:weightedderi} by $\chi_R(r)r^{-3+\de}\ov{\dk^{\reg}\phis{p}}$, where $\chi_R(r)=\chi(r/R)$ with $\chi(r)$ a smooth nonnegative cut-off function that equals $1$ for $r\geq 1$ and vanishes for $r\leq \frac{1}{2}$, and the rest of the proof follows in the same manner as for \eqref{eq:improvedestiforangularderiofkerrpert:prop:highorder}. This concludes the proof of Proposition \ref{prop:ellipticestimates:phis01:angularderivatives:highorder}.
\end{proof}


\subsection{Proof of Theorem \ref{thm:EMF:systemofTeuscalarized:order0:final}}
\lab{subsect:proofofEMFmain:reg=0}


In this section, we provide the proof for Theorem \ref{thm:EMF:systemofTeuscalarized:order0:final}. 
To this end, let $B_1\gg 1$ be a large enough constant that will be fixed later, and define, for $s=\pm 2$ and $\de\in(0,\frac{1}{3}]$,  
\begin{align}\lab{def:EMFtotalppps:pm2}
\EMFtotalpp{s}:={}&\sum_{p=0,1,2}(B_1)^{2-p}\widetilde{\EMF}[\psis{p}]+ {\EMF}_{\de}[\phis{2}](\tau_1,\tau_2)
+ {\M}_{\de}[\psis{2}](\tau_1+1,\tau_2-3)
\nn\\
&+\sum_{p=0,1}(B_1)^{2-p}\Big(\widehat{\EMF}_{\de}[\phis{p}](\tau_1,\tau_2)+\widehat{\M}_{\de}[\psis{p}](\tau_1+1, \tau_2-3)\Big)
\end{align}
and
\begin{align}
\lab{def:NNttotalppps:pm2}
\NNttotalpp{s}:={}&
\sum_{p=0,1,2}(B_1)^{2-p}\widetilde{\mathcal{N}}_{W,\de}[\psis{p},\phis{p}]+\sum_{p=0,1}(B_1)^{2-p}\int_{\MM_{r\leq 12m}(\tau_1,\tau_2)}\Big(|\pr^{\leq 1}\N_{T, s}^{(p)}|^2+|\N^{(p)}_{W,s}|^2\Big)\nn\\
&+\sum_{p=0,1}(B_1)^{2-p}\int_{\MM(\tau_1,\tau_2)}\Big(r^{-1+\de}|\N_{W,s}^{(p)}|+r^{-2+\de}|\nab_{X_s}\N^{(p)}_{T,s}|\Big)|\phis{p}|,
\end{align}
with $\widetilde{\mathcal{N}}_{W,\de}[\c,\c]$ as given in \eqref{def:NNtWdepsiphi} and vectorfields $X_s$ as introduced in \eqref{def:vectorfieldXs}.

We then deduce the following conditional EMF estimates.

\begin{proposition}
\lab{prop:EMF:systemofTeuscalarized:order0}
For $B_1\gg 1$ sufficiently large, we have 
\bea\lab{eq:EMFtotal:sumup:final:rweightscontrolled:pm2}
\EMFtotalpp{s} &\les& \NNttotalpp{s} +(B_1)^{2}\Big(\A[\pmb\psi_s](\Iti) +\A[\pmb\phi_{s}](\tau_1,\tau_2) +\IE{\pmb\phi_s, \pmb\psi_s}\Big)\nn\\
&&+\sum_{p=0}^2 B_1^{2-p} \Errdefect[\psis{p}].
\eea
\end{proposition}

\begin{proof}
Adding $(B_1)^{2}$ multiple of the estimate \eqref{th:eq:mainnondegeenergymorawetzmicrolocal:tensorialwave:scalarized:0} for $\psis{0}$ together with $B_1$ multiple of the estimate \eqref{th:eq:mainnondegeenergymorawetzmicrolocal:tensorialwave:scalarized:1} for $\psis{1}$ to the estimate \eqref{th:eq:mainnondegeenergymorawetzmicrolocal:tensorialwave:scalarized:2} for $\psis{2}$,  we obtain
\begin{align}
\lab{eq:EMFtotals:esti:v1}
&\EMFtotalpp{s}\nn\\
\les{}& \NNttotalpp{s} +(B_1)^{2}\Big(\A[\pmb\psi_s](\Iti) +\A[\pmb\phi_{s}](\tau_1,\tau_2) +\IE{\pmb\phi_s, \pmb\psi_s}\Big)
+\sum_{p=0}^2 B_1^{2-p} \Errdefect[\psis{p}]
\nn\\
&+\sum_{p=0,1,2}(B_1)^{2-p}\Big(\widetilde{\mathcal{N}}[\psis{p},\chi_{\tau_1, \tau_2}^{(1)}\L^{(p)}_{s}]+\NNtdemora[\psis{p},\L_{s}^{(p)}](\tau_1+1,\tau_2-3)\Big),
\end{align}
where we have absorbed the terms 
\beaa
&&B_1\bigg(\widehat{\EMF}[\psis{0}](\tau_1+1,\tau_2-3)+\sup_{\tau\in[\tau_1,\tau_1+1]\cup[\tau_2-3,\tau_2]}\E[\phis{0}](\tau)\bigg)\\
&& +\sum_{p=0,1}\sup_{\tau\in[\tau_1,\tau_1+1]\cup[\tau_2-3,\tau_2]}\E[\phis{p}](\tau)
\eeaa
by the LHS by taking $B_1\gg 1$ large enough. 

Next, we rely on \eqref{subeq:naivecontroloferrorfromlinearcoupling:scalarizedTeu} to estimate the last line of  \eqref{eq:EMFtotals:esti:v1} as follows
\beaa
&&\sum_{p=0,1,2}(B_1)^{2-p}\Big(\widetilde{\mathcal{N}}[\psis{p},\chi_{\tau_1, \tau_2}^{(1)}\L^{(p)}_{s}]+\NNtdemora[\psis{p},\L_{s}^{(p)}](\tau_1+1,\tau_2-3)\Big)\nn\\
&\les&\sum_{p=0,1}(B_1)^{1-p}\widehat{\M}[\phis{p}](\tau_1,\tau_2-2)+(B_1)^{2}\A[\pmb\phi_{s}](\tau_1,\tau_2-2)\nn\\
&&
+\sum_{p=0,1}(B_1)^{2-p}\Big(\Ao[r\nab\phis{p}](\tau_1,\tau_2-2) +\ep \widehat{\M}[\phis{p}](\tau_1,\tau_2-2)\Big) +\int_{\MM_{r\leq 12m}(\tau_1,\tau_2)}|\pr^{\leq 1}\N_{T, s}^{(1)}|^2
\nn\\
&&+(B_1)^{2}\Big(\widehat{\M}_{\de}[\psis{0}](\tau_1+1,\tau_2-3)+\widetilde{\EM}[\psis{0}] \Big)^{\frac{1}{2}}\nn\\
&&\times\left(\A[\pmb\phi_s](\tau_1,\tau_2-2)+\Ao[r\nab\phis{0}](\tau_1,\tau_2-2) +\ep\widehat{\M}[\phis{0}](\tau_1,\tau_2-2)\right)^{\frac{1}{2}}\nn\\
&&+B_1\Bigg(\Extra[\phis{0}](\tau_1,\tau_2-2)+\A[\pmb\phi_s](\tau_1,\tau_2-2)
+\Ao[r\nab\phis{1}](\tau_1,\tau_2-2) +\ep\widehat{\M}[\phis{1}](\tau_1,\tau_2-2) \nn\\ 
&& +\sup_{\tau\in[\tau_1, \tau_1+1]\cup[\tau_2-3,\tau_2-2]}\E[\phis{0}](\tau)\Bigg)^{\frac{1}{2}}\Big(\widehat{\M}_{\de}[\psis{1}](\tau_1+1,\tau_2-3)+\widetilde{\EM}[\psis{1}] \Big)^{\frac{1}{2}}\nn\\
&& +\bigg(\sum_{p=0,1}\widehat{\EM}[\phis{p}](\tau_1,\tau_2-2)+\sum_{p=0,1}\Extra[\phis{p}](\tau_1,\tau_2-2)+\A[\pmb\phi_s](\tau_1,\tau_2-2)\bigg)^{\frac{1}{2}} \nn\\
&&\times\Big(\widetilde{\EM}[\psis{2}]+\EM[\phis{2}](\tau_1,\tau_2-2)
+{\M}_{\de}[\psis{2}](\tau_1+1,\tau_2-3)\Big)^{\frac{1}{2}},
\eeaa
and, substituting this back into inequality \eqref{eq:EMFtotals:esti:v1} and applying H\"older's inequality to the product terms, we deduce, by taking $B_1\gg 1$ large enough and $\ep$ small enough,
\bea
\lab{eq:estimatesof:EMFtotal''s}
&&\EMFtotalpp{s}\nn\\
&\les& \NNttotalpp{s} +(B_1)^{2}\Big(\A[\pmb\psi_s](\Iti) +\A[\pmb\phi_{s}](\tau_1,\tau_2) +\IE{\pmb\phi_s, \pmb\psi_s}\Big) +\sum_{p=0}^2(B_1)^{2-p}\Errdefect[\psis{p}]
\nn\\
&&+\sum_{p=0,1}(B_1)^{1-p}\Extra[\phis{p}](\tau_1,\tau_2)+\sum_{p=0,1}(B_1)^{2-p}\Ao[r\nab\phis{p}](\tau_1,\tau_2).
\eea

Next, we use the estimates for  $\Ao[r^{\frac{\de}{2}}(r\nab)^{\leq 1}\phis{p}](\tau_1,\tau_2)$ and $\Extra[\phis{p}](\tau_1,\tau_2)$ in Proposition \ref{prop:ellipticestimates:phis01:angularderivatives} to control the terms in the last line of \eqref{eq:estimatesof:EMFtotal''s} by
\bea
&&\sum_{p=0,1}(B_1)^{1-p}\Extra[\phis{p}](\tau_1,\tau_2)+\sum_{p=0,1}(B_1)^{2-p}\Ao[r\nab\phis{p}](\tau_1,\tau_2)\nn\\
&\les&\sum_{p=0,1}(B_1)^{2-p}\Ao[r^{\frac{\de}{2}}(r\nab)^{\leq 1}\phis{p}](\tau_1,\tau_2) +\sum_{p=0,1}(B_1)^{1-p}\EM[\phis{p}](\tau_1,\tau_2)\nn\\
&\les&\sum_{p=0,1}(B_1)^{1-p}\EM[\phis{p}](\tau_1,\tau_2)+\sum_{p=0,1}(B_1)^{2-p}\bigg(\Big(\widehat{\EMF}[\phis{p}](\tau_1,\tau_2)\Big)^{\frac{1}{2}}\big(\A[\pmb\phi_{s}](\tau_1,\tau_2)\big)^{\frac{1}{2}} \nn\\
&&
+\A[\pmb\phi_{s}](\tau_1,\tau_2)+\int_{\MM(\tau_1,\tau_2)}\Big(r^{-1+\de}|\N_{W,s}^{(p)}|+r^{-2+\de}|\nab_{X_s}\N^{(p)}_{T,s}|\Big)|\phis{p}|\bigg).
\eea
Combining this with the estimate \eqref{eq:estimatesof:EMFtotal''s}, 
 we take $B_1$ large enough to infer
\beaa
\EMFtotalpp{s} &\les& \NNttotalpp{s} +(B_1)^{2}\Big(\A[\pmb\psi_s](\Iti) +\A[\pmb\phi_{s}](\tau_1,\tau_2) +\IE{\pmb\phi_s, \pmb\psi_s}\Big)\nn\\
&&+\sum_{p=0}^2 (B_1)^{2-p} \Errdefect[\psis{p}]
\eeaa
as desired. This concludes the proof of Proposition \ref{prop:EMF:systemofTeuscalarized:order0}.
\end{proof}

We now fix the large positive constant $B_1$ and we are in position to prove Theorem  \ref{thm:EMF:systemofTeuscalarized:order0:final}. 

\begin{proof}[Proof of Theorem \ref{thm:EMF:systemofTeuscalarized:order0:final}]
Since the constant $B_1$ is fixed, it follows from \eqref{def:EMFtotalps:pm2} and \eqref{def:EMFtotalppps:pm2} that
\bea
\lab{equivalence:twotypesofEMFnorms:zeroorder}
\EMFtotalp{s}\simeq\EMFtotalpp{s}, \quad s=\pm2.
\eea
This together with \eqref{eq:EMFtotal:sumup:final:rweightscontrolled:pm2} implies 
\bea\lab{eq:EMFtotal:sumup:final:rweightscontrolled:pm2:rewrite}
\EMFtotalp{s} \les  \IE{\pmb\phi_s, \pmb\psi_s}  + \NNttotalpp{s}+\A[\pmb\psi_{s}](\Iti)+\A[\pmb\phi_{s}](\tau_1,\tau_2)+\sum_{p=0}^2\Errdefect[\psis{p}].
\eea

Next, we estimate the term $\NNttotalpp{s}$, $s=\pm 2$. We have from \eqref{def:NNttotalppps:pm2} and \eqref{def:NNtWdepsiphi} that
\bea
\lab{expression:NNttotalppps:pm2}
\NNttotalpp{s}
&\simeq&\sum_{p=0,1,2}\Big(\widetilde{\mathcal{N}}[\psis{p},\widehat{\F}_{total,s}^{(p)}]+\NNtdemora[\psis{p}, \widehat{\F}_{total,s}^{(p)}](\tau_1+1,\tau_2-3) \Big)\nn\\
&&+\sum_{p=0,1,2}\Big(\NNtlede[\phis{p}, \N_{W,s}^{(p)}](\tau_2-3, \tau_2)+\NNtlede[\phis{p}, \N_{W,s}^{(p)}](\tau_1, \tau_1+1)
 \Big)\nn\\
&&+\sum_{p=0,1}\int_{\MM(\tau_1,\tau_2)}\Big(r^{-1+\de}|\N_{W,s}^{(p)}|+r^{-2+\de}|\nab_{X_s}\N^{(p)}_{T,s}|\Big)|\phis{p}|\nn\\
&&+\sum_{p=0,1}\int_{\MM_{r\leq 12m}(\tau_1,\tau_2)}\Big(|\pr^{\leq 1}\N_{T, s}^{(p)}|^2+|\N^{(p)}_{W,s}|^2\Big).
\eea
In view of \eqref{def:NNtleinlocalenergyestimate}
which defines $\NNtlede[\pmb\psi,\F](\tau', \tau'')$, \eqref{def:NNtdemora} which defines $\NNtdemora[\pmb\psi,\F](\tau',\tau'')$, and \eqref{eq:causlityrelationsforwidetildepsi1} which yields $\widehat{F}_{total,s,ij}^{(p)}=N^{(p)}_{W,s,ij}$ in $\MM(\tau_1+1,\tau_2-3)$ since $\chi_{\tau_1, \tau_2}=\chi^{(1)}_{\tau_1, \tau_2}=1$ there, we infer 
\beaa
&&\sum_{p=0,1,2}\NNtdemora[\psis{p},\widehat{\F}_{total,s}^{(p)}](\tau_1+1,\tau_2-3) \nn\\
&&+\sum_{p=0,1,2}\Big(\NNtlede[\phis{p}, \N_{W,s}^{(p)}](\tau_2-3, \tau_2)+\NNtlede[\phis{p}, \N_{W,s}^{(p)}](\tau_1, \tau_1+1)\Big)
\nn\\
&\les&\sum_{p=0,1,2}\Bigg(\sum_{i,j}\sup_{\tau_1< \tau'<\tau''< \tau_2}\bigg|\int_{\Mntrap(\tau', \tau'')} \Re\Big(\ov{N_{W,s,ij}^{(p)}}\big(1+O(r^{-\de})\big)\pr_{\tt}\phiss{ij}{p}\Big)\bigg|
\nn\\
&&+\sum_{i,j}\int_{\Mntrap(\tau_1, \tau_2)}  r^{-1}|\dk^{\leq 1}\phiss{ij}{p}| |N_{W,s,ij}^{(p)}| +\int_{\MM(\tau_1,\tau_2)}|\N_{W,s}^{(p)}|^2\Bigg)
\nn\\
&\les& \sum_{p=0,1,2}\widehat{\mathcal{N}}'[ \phis{p}, \N_{W,s}^{(p)}](\tau_1, \tau_2),
 \eeaa
 where we used the definition of $\widehat{\mathcal{N}}'[\pmb\psi, \H](\tau_1, \tau_2)$ from
\eqref{eq:defofNNhat'}. 
Also, using  the formula \eqref{expression:NNttotalps:pm2} of $\NNttotalp{s}$ and the formula \eqref{eq:defofNNhat'} of $\widehat{\NN}'[\c,\c](\c,\c)$, we have 
\beaa
&&\sum_{p=0,1}\int_{\MM(\tau_1,\tau_2)}\Big(r^{-1+\de}|\N_{W,s}^{(p)}|+r^{-2+\de}|\nab_{X_s}\N^{(p)}_{T,s}|\Big)|\phis{p}|\nn\\
&&+\sum_{p=0,1}\int_{\MM_{r\leq 12m}(\tau_1,\tau_2)}\Big(|\pr^{\leq 1}\N_{T, s}^{(p)}|^2+|\N^{(p)}_{W,s}|^2\Big)
\les\NNttotalp{s}.
\eeaa
Hence, the above together yields
\beaa
\NNttotalpp{s} \les \NNttotalp{s},
\eeaa
and plugging this estimate into \eqref{eq:EMFtotal:sumup:final:rweightscontrolled:pm2:rewrite}, we infer
\beaa
\EMFtotalp{s} \les  \IE{\pmb\phi_s, \pmb\psi_s}  + \NNttotalp{s}+\A[\pmb\psi_{s}](\Iti)+\A[\pmb\phi_{s}](\tau_1,\tau_2)+\sum_{p=0}^2\Errdefect[\psis{p}].
\eeaa
This proves the desired estimate \eqref{eq:EMFtotalp:sumup:rweightscontrolled:pm2} and hence concludes the proof of Theorem \ref{thm:EMF:systemofTeuscalarized:order0:final}.
\end{proof}


\section{Proof of Theorem \ref{th:main:intermediary}}
\lab{sec:proofofth:main:intermediary}


In this section, we provide the proof for Theorem \ref{th:main:intermediary} which proves high order EMF estimates for unweighted derivatives. We derive first in Sections \ref{subsect:unweightedEMFintermsofprtauprtphihat}-\ref{subsect:highorderunweightedEMF:conditional} high order EMF estimates for unweighted derivatives conditional on the control of the terms $\mathbf{IE}$, $\Errdefect$ and $\NNt_{s, \de,\text{total}}$.  Next, in Sections \ref{subsect:controlofIEtotalandErrdefect}-\ref{subsect:controlofNNtsdetotal}, we provide the control of the terms $\Errdefect$ and $\NNt_{s, \de,\text{total}}$. The proof of Theorem \ref{th:main:intermediary} is finally concluded in Section \ref{sec:endfotheproofofth:main:intermediary}.


\subsection{Estimates for unweighted derivatives conditional on $(\pr_\tau, \chi_0(r)\pr_{\tphi})$ derivatives}
\lab{subsect:unweightedEMFintermsofprtauprtphihat}


We start by defining, for $\pmb\psi\in\sk_k(\mathbb{C})$ and $\pmb H\in\sk_k(\mathbb{C})$, $k=1,2$, 
\bea
\lab{eq:defofNNhat''}
\widehat{\mathcal{N}}''[\pmb\psi, \H](\tau_1, \tau_2)
&:=&\sup_{\tau_1< \tau'<\tau''< \tau_2}\bigg|\int_{\Mntrap(\tau', \tau'')}\nab_{\pr_\tau}\pmb\psi\c\ov{\H}\bigg| \nn\\
&&+\int_{\Mntrap(\tau_1, \tau_2)}r^{-1}|\dk^{\leq 1}\pmb\psi||\H|+\int_{\MM(\tau_1, \tau_2)}|\H|^2.
\eea
The following lemma allows to derive conditional EMF estimates for all high order (unweighted) derivatives conditional on the control of high order $(\pr_\tau, \chi_0(r)\pr_{\tphi})$ derivatives. 

\begin{lemma}\lab{lemma:higherorderenergyMorawetzestimates}
Let $\g$ satisfy the assumptions of Section  \ref{subsubsect:assumps:perturbedmetric}, let $\tau_2>\tau_1$ satisfy \eqref{eq:tau2geqtau1plus10conditionisthemaincase}, and let $\reg$ be a positive integer satisfying $1\leq \reg \leq 14$. Let $(\psi)_{ij}$ and $(F)_{ij}$ be families of complex-valued scalars satisfying the following coupled system of wave equations
\bea
\lab{eq:eqsforconditionalhighorderestimatelemma:general}
\square_{\g}\psi_{ij}=N_{ij}, \qquad N_{ij}:=D_1r^{-1}\pr_{\tt}\psi_{ij}+O(r^{-2})\dk^{\leq 1}\psi_{kl} +F_{ij},
\eea
where $D_1\geq0$ is a constant and the coefficients of the terms in $N_{ij}$ may depend on $\dk^{\leq 1}\Ga_b$ and $\dk^{\leq 1}\Ga_g$.
Assume that $(\psi)_{ij}$ satisfy the following redshift estimate:
\bea
\lab{eq:assum:redshift:highorderunweightedEMF:generalwave}
\nn\EMF_{r\leq r_+(1+\dred)}[\pr^{\leq \reg}\pmb\psi](\tau_1, \tau_2)&\les& \E[\pr^{\leq \reg}\pmb\psi](\tau_1)+\dred^{-1}\M_{r_+(1+\dred), r_+(1+2\dred)}[\pr^{\leq \reg}\pmb\psi](\tau_1, \tau_2)\\
&&+\int_{\MM_{r\leq r_+(1+2\dred)}(\tau_1, \tau_2)}|\pr^{\leq \reg}{\F}|^2.
\eea
Let $\chi_0=\chi_0(r)$ be a smooth cut-off function such that
\bea
\lab{def:cutoffcuntionchi0}
0\leq\chi_0\leq 1, \qquad \chi_0=1 \qquad \text{for } \,\, r\leq 11m, \quad \chi_0=0 \quad \text{for }\,\, r\geq 12m.
\eea
Then, for  any $1\leq\tau_1<\tau_2 <+\infty$, we have  the improved estimate
\bea\lab{eq:higherorderenergymorawetzforlargerconditionalonsprtau}
\nn&&\EMF_{r\geq 11m}[\pr^{\leq \reg}\pmb\psi](\tau_1, \tau_2) \nn\\
&\les& \EMF[\pr_\tau^{\leq \reg}\pmb\psi](\tau_1, \tau_2)+\big(\EMF[\pr^{\leq \reg-1}\pmb\psi](\tau_1, \tau_2)\big)^{\frac{1}{2}}\left(\EMF[\pr^{\leq \reg}\pmb\psi](\tau_1, \tau_2)\right)^{\frac{1}{2}}\nn\\
&&+\sup_{\tau\in[\tau_1,\tau_2-2]}\widehat{\mathcal{N}}_{r\geq {10m}}''[\pr^{\leq \reg}\pmb\psi, \pr^{\leq \reg}\F](\tau, \tau+2)+\int_{\MM_{r\geq 10m}(\tau_1,\tau_2)}|\pr^{\leq  \reg}\F|^2,
\eea
and
\bea\lab{eq:higherorderenergymorawetzconditionalonsprtausprphis}
\nn\EMF[\pr^{\leq \reg}\pmb\psi](\tau_1, \tau_2) &\les& \E[\pr^{\leq \reg}\pmb\psi](\tau_1)+\EMF[(\pr_\tau, \chi_0\pr_{\tphi})^{\leq \reg}\pmb\psi](\tau_1, \tau_2)+\int_{\MM(\tau_1, \tau_2)}|\pr^{\leq \reg}\F|^2\\
&&+\sup_{\tau\in[\tau_1,\tau_2-2]}\widehat{\mathcal{N}}_{r\geq {10m}}''[\pr^{\leq \reg}\pmb\psi, \pr^{\leq \reg}\F](\tau, \tau+2)
\eea
as well as
\bea\lab{eq:higherorderenergymorawetzconditionalonsprtausprphis:nontrapped}
\nn\widehat{\EMF}[\pr^{\leq \reg}\pmb\psi](\tau_1, \tau_2) &\les& \E[\pr^{\leq \reg}\pmb\psi](\tau_1)+\widehat{\EMF}[(\pr_\tau, \chi_0\pr_{\tphi})^{\leq \reg}\pmb\psi](\tau_1, \tau_2)+\int_{\MM(\tau_1, \tau_2)}|\pr^{\leq \reg}\F|^2\\
&&+\sup_{\tau\in[\tau_1,\tau_2-2]}\widehat{\mathcal{N}}_{r\geq {10m}}''[\pr^{\leq \reg}\pmb\psi, \pr^{\leq \reg}\F](\tau, \tau+2),
\eea
where $\widehat{\EMF}$ and $\widehat{\NN}$ are given as in \eqref{def:widehatEMFdenorm} and \eqref{eq:defintionwidehatmathcalNfpsinormRHS}, respectively.
Furthermore, we have, for any $0\leq\de\leq 1$, 
\bea\lab{eq:improvedhigherorderenergymorawetzforlargerconditionalonsprtau}
\nn&&\EMF_{\de, r\geq 11m}[\pr^{\leq \reg}\pmb\psi](\tau_1, \tau_2) \nn\\
&\les& \EMF_\de [\pr_\tau^{\leq \reg}\pmb\psi](\tau_1, \tau_2)+\big(\EMF [\pr^{\leq \reg-1}\pmb\psi](\tau_1, \tau_2)\big)^{\frac{1}{2}}\left(\EMF[\pr^{\leq \reg}\pmb\psi](\tau_1, \tau_2)\right)^{\frac{1}{2}}\nn\\
&&+\sup_{\tau\in[\tau_1,\tau_2-2]}\widehat{\mathcal{N}}_{r\geq {10m}}''[\pr^{\leq \reg}\pmb\psi, \pr^{\leq \reg}\F](\tau, \tau+2)+\int_{\MM_{r\geq 10m}(\tau_1,\tau_2)}|\pr^{\leq  \reg}\F|^2
\eea
and
\bea\lab{eq:imporvedhigherorderenergymorawetzconditionalonlowerrweigth}
\nn\EMF_{\de}[\pr^{\leq \reg}\pmb\psi](\tau_1, \tau_2) &\les& \E[\pr^{\leq \reg}\pmb\psi](\tau_1)+\EMF_{\de}[(\pr_\tau, \chi_0\pr_{\tphi})^{\leq \reg}\pmb\psi](\tau_1, \tau_2)+\int_{\MM(\tau_1, \tau_2)}|\pr^{\leq \reg}\F|^2\\
&&+\sup_{\tau\in[\tau_1,\tau_2-2]}\widehat{\mathcal{N}}_{r\geq {10m}}''[\pr^{\leq \reg}\pmb\psi, \pr^{\leq \reg}\F](\tau, \tau+2)
\eea
as well as
\bea\lab{eq:imporvedhigherorderenergymorawetzconditionalonlowerrweigth:nontrapped}
\widehat{\EMF}_{\de}[\pr^{\leq \reg}\pmb\psi](\tau_1, \tau_2) &\les& \E[\pr^{\leq \reg}\pmb\psi](\tau_1)+\widehat{\EMF}_{\de}[(\pr_\tau, \chi_0\pr_{\tphi})^{\leq \reg}\pmb\psi](\tau_1, \tau_2)+\int_{\MM(\tau_1, \tau_2)}|\pr^{\leq \reg}\F|^2\nn\\
&&+\sup_{\tau\in[\tau_1,\tau_2-2]}\widehat{\mathcal{N}}_{r\geq {10m}}''[\pr^{\leq \reg}\pmb\psi, \pr^{\leq \reg}\F](\tau, \tau+2).
\eea
\end{lemma}

\begin{proof}
The proof is an adaptation of the corresponding one for the scalar field in \cite[Lemma 3.13]{MaSz24}. Note from the assumption \eqref{eq:tau2geqtau1plus10conditionisthemaincase} that we have in particular $\tau_2>\tau_1+2$.

\noindent{\bf Step 1.} \textit{Proof of  \eqref{eq:higherorderenergymorawetzforlargerconditionalonsprtau}: Morawetz part}. We have from \cite[inequality (3.28)]{MaSz24} and the formula \eqref{eq:eqsforconditionalhighorderestimatelemma:general} of $N_{ij}$  that
\beaa
\M_{r\geq 10.5m}[\pr^{\leq 1}\pmb\psi](\tau_1,\tau_2) 
&\les& \M_{r\geq 10m}[\pr_{\tau}^{\leq 1}\pmb\psi](\tau_1,\tau_2)
+\int_{\MM_{{r\geq 10m}}(\tau_1,\tau_2)} r^{-1}  |\N|^2\nn\\
&& +\sqrt{\EM[\pmb\psi](\tau_1,\tau_2)}\sqrt{\EM[\pr^{\leq 1}\pmb\psi](\tau_1,\tau_2)}\nn\\
&\les&  \M_{r\geq 10m}[\pr_{\tau}^{\leq 1}\pmb\psi](\tau_1,\tau_2)
+\int_{\MM_{{r\geq 10m}}(\tau_1,\tau_2)} r^{-1}  |\F|^2\nn\\
&& +\sqrt{\EM[\pmb\psi](\tau_1,\tau_2)}\sqrt{\EM[\pr^{\leq 1}\pmb\psi](\tau_1,\tau_2)}.
\eeaa
In fact, by examining the proof for \cite[inequality (3.28)]{MaSz24}, the above estimates hold true as well if  replacing $10.5m$ and $10m$ by $10m+n_1m$ and $10m+n_2m$ respectively, for any $0\leq n_2<n_1$, i.e.,
\bea
\lab{eq:EM:AsympKerr:highMora:v1:mid}
\M_{r\geq 10m+n_1m}[\pr^{\leq 1}\pmb\psi](\tau_1,\tau_2) 
&\les&  \M_{r\geq 10m+n_2m}[\pr_{\tau}^{\leq 1}\pmb\psi](\tau_1,\tau_2)
+\int_{\MM_{{r\geq 10m+n_2m}}(\tau_1,\tau_2)} r^{-1}  |\F|^2\nn\\
&& +\sqrt{\EM[\pmb\psi](\tau_1,\tau_2)}\sqrt{\EM[\pr^{\leq 1}\pmb\psi](\tau_1,\tau_2)}.
\eea

Commuting $\pr_{\tt}^{\reg}$ with the wave equation \eqref{eq:eqsforconditionalhighorderestimatelemma:general}, and in view of the following commutator formula 
\bea
\label{esti:commutatorBoxgandT:general:moreexplicit}
 [ \pr_{\tau}, \square_{\g}]\psi = \dk\Ga_b\c \dk \pr\psi+
\dk^{\leq 2}\Ga_g\c\dk\psi
\eea
which follows from \eqref{esti:commutatorBoxgandT:general:0} together with the assumption \eqref{eq:controloflinearizedinversemetriccoefficients}, we deduce
\begin{equation}
\lab{eq:eqsforconditionalhighorderestimatelemma:general:prttreg}
\bsplit
\square_{\g}\pr_{\tt}^{\reg}\psi_{ij}={}&D_1r^{-1}\pr_{\tt} \pr_{\tt}^{\reg}\psi_{ij}+O(r^{-2})\dk^{\leq 1}\pr_{\tt}^{\leq \reg}\psi_{kl} +F_{ij}^{(0,\reg)},\\
F_{ij}^{(0,\reg)}={}&
\pr_{\tt}^{\reg}F_{ij}
+\dk^{\leq \reg+1}\Ga_b \c \dk \pr^{\leq \reg}\psi_{kl}.
\end{split}
\end{equation}
Further, commuting $(\pr_r, r^{-1}\pr_{x^a})$ with the above wave equation and using the commutator relation \eqref{eq:localwavecommutators:withfirstordergoodderis:0}, we deduce
\begin{equation}
\lab{eq:eqsforconditionalhighorderestimatelemma:general:prreg}
\bsplit
\square_{\g}(\pr_r, r^{-1}\pr_{x^a})^{\reg_1}\pr_{\tt}^{\reg_2}\psi_{ij}={}&D_1r^{-1}\pr_{\tt}(\pr_r, r^{-1}\pr_{x^a})^{\reg_1}\pr_{\tt}^{\reg_2}\psi_{ij}\\
&
+O(r^{-2})\dk^{\leq 1}\pr^{\leq \reg_1}\pr_{\tt}^{\leq \reg_2}\psi_{kl}+{F}_{ij}^{(\reg_1,\reg_2)},\\
{F}_{ij}^{(\reg_1,\reg_2)}={}&
(\pr_r, r^{-1}\pr_{x^a})^{\reg_1}\pr_{\tt}^{\reg_2}F_{ij}
+\dk^{\leq \reg_1+\reg_2+1}\Ga_b\c \dk \pr^{\leq \reg_1+\reg_2}\psi_{kl}.
\end{split}
\end{equation}
Combining this equation with the wave equation \eqref{eq:eqsforconditionalhighorderestimatelemma:general:prttreg}, we infer
\begin{equation}
\lab{eq:eqsforconditionalhighorderestimatelemma:general:allprreg}
\bsplit
\square_{\g}(\pr^{\reg_1}\pr_{\tt}^{\reg_2}\psi_{ij})={}&D_1r^{-1}\pr_{\tt}\pr^{\reg_1}\pr_{\tt}^{\reg_2}\psi_{ij}+O(r^{-2})\dk^{\leq 1}\pr^{\leq \reg_1}\pr_{\tt}^{\leq \reg_2}\psi_{kl}+{F}_{ij}^{(\reg_1,\reg_2)},\\
{F}_{ij}^{(\reg_1,\reg_2)}={}&
\pr^{\reg_1}\pr_{\tt}^{\reg_2}F_{ij}
+\dk^{\leq \reg_1+\reg_2+1}\Ga_b\c\dk \pr^{\leq \reg_1+\reg_2}\psi_{kl}.
\end{split}
\end{equation}

We view the above wave equations for $(\pr_{\tt}^{\leq \reg-1}\psi_{ij})$, i.e., with $(\reg_1=0, \reg_2\leq \reg-1)$ in \eqref{eq:eqsforconditionalhighorderestimatelemma:general:allprreg}, as a coupled wave system and apply \eqref{eq:EM:AsympKerr:highMora:v1:mid} with $(n_1,n_2)=(10m+(2\reg)^{-1}m, 10m)$ to this wave system for the scalars $(\pr_{\tt}^{\leq \reg-1}\psi_{ij})$ to deduce
\bea
\lab{eq:EM:AsympKerr:highMora:v1:middle:generalk:prttk}
\M_{r\geq 10m+(2\reg)^{-1}m}[\pr^{\leq 1}\pr_{\tt}^{\leq \reg-1}\pmb\psi](\tau_1,\tau_2) 
&\les&\M_{r\geq 10m}[\pr_{\tau}^{\leq \reg}\pmb\psi](\tau_1,\tau_2)
+\int_{\MM_{{r\geq 10m}}(\tau_1,\tau_2)} r^{-1}  |\pr^{\leq \reg-1}\F|^2\nn\\
&& +\sqrt{\EM[\pr^{\leq \reg-1}\pmb\psi](\tau_1,\tau_2)}\sqrt{\EM[\pr^{\leq \reg}\pmb\psi](\tau_1,\tau_2)},
\eea
where we used the fact that, for $\reg\leq 14$,
\bea
\lab{eq:highorderrecovering:errorforprreg}
\sum_{k,l}\int_{\MM_{{r\geq 10m}}(\tau_1,\tau_2)} r^{-1}  | \dk^{\leq\reg}\Ga_b\c\dk\pr^{\leq \reg-1}\psi_{kl}|^2\les \ep \sup_{\tau\in(\tau_1,\tau_2)}\E[\pr^{\leq \reg-1}\pmb\psi](\tau).
\eea

Next, we consider the following induction hypothesis, for any integer $n$ with $1\leq n\leq \reg-1$,
\bea
\lab{eq:EM:AsympKerr:highMora:v1:middle:generalk:prnprttk-n:inductionhypothesis}
&&\M_{r\geq 10m+n(2\reg)^{-1}m}[\pr^{\leq n}\pr_{\tt}^{\leq \reg-n}\pmb\psi](\tau_1,\tau_2) \nn\\
&\les&\M_{r\geq 10m}[\pr_{\tt}^{\leq \reg}\pmb\psi](\tau_1,\tau_2)
+\int_{\MM_{r\geq 10m}(\tau_1,\tau_2)} r^{-1}  |\pr^{\leq \reg-1}\F|^2\nn\\
&& +\sqrt{\EM[\pr^{\leq \reg-1}\pmb\psi](\tau_1,\tau_2)}\sqrt{\EM[\pr^{\leq \reg}\pmb\psi](\tau_1,\tau_2)}.
\eea
In view of \eqref{eq:EM:AsympKerr:highMora:v1:middle:generalk:prttk}, the induction hypothesis  \eqref{eq:EM:AsympKerr:highMora:v1:middle:generalk:prnprttk-n:inductionhypothesis} holds true for $n=1$. We now assume that the induction hypothesis  \eqref{eq:EM:AsympKerr:highMora:v1:middle:generalk:prnprttk-n:inductionhypothesis} holds true for some $1\leq n\leq \reg-2$ and our goal is to show that it also holds for $n+1$. To this end, we view the above wave equations for $(\pr^{\leq n}\pr_{\tt}^{\leq \reg-1-n}\psi_{ij})$, i.e., with $(\reg_1\leq n, \reg_2\leq \reg-1-n)$ in \eqref{eq:eqsforconditionalhighorderestimatelemma:general:allprreg}, as a coupled wave system and apply \eqref{eq:EM:AsympKerr:highMora:v1:mid} with $(n_1,n_2)=(10m+(n+1)(2\reg)^{-1}m, 10m+n(2\reg)^{-1}m)$ to this wave system for $(\pr^{\leq n}\pr_{\tt}^{\leq \reg-1-n}\psi_{ij})$ to deduce
\bea
\lab{eq:EM:AsympKerr:highMora:v1:middle:generalk:prnprttk-n:toproveinductionwithnplus1}
&&\M_{r\geq 10m+(n+1)(2\reg)^{-1}m}[\pr^{\leq n+1}\pr_{\tt}^{\leq \reg-(n+1)}\pmb\psi](\tau_1,\tau_2) \nn\\
&\les&\M_{r\geq 10m+n(2\reg)^{-1}m}[\pr^{\leq n}\pr_{\tt}^{\leq \reg-n}\pmb\psi](\tau_1,\tau_2)
+\int_{\MM_{r\geq 10m}(\tau_1,\tau_2)} r^{-1}  |\pr^{\leq\reg -1}\F|^2\nn\\
&& +\sqrt{\EM[\pr^{\leq \reg-1}\psi](\tau_1,\tau_2)}\sqrt{\EM[\pr^{\leq \reg}\pmb\psi](\tau_1,\tau_2)},
\eea
where we again used the estimate \eqref{eq:highorderrecovering:errorforprreg}. Using the induction hypothesis \eqref{eq:EM:AsympKerr:highMora:v1:middle:generalk:prnprttk-n:inductionhypothesis} to control the first term on the RHS of \eqref{eq:EM:AsympKerr:highMora:v1:middle:generalk:prnprttk-n:toproveinductionwithnplus1}, we infer
\beaa
&&\M_{r\geq 10m+(n+1)(2\reg)^{-1}m}[\pr^{\leq n+1}\pr_{\tt}^{\leq \reg-(n+1)}\pmb\psi](\tau_1,\tau_2) \nn\\
&\les&\M_{r\geq 10m}[\pr_{\tt}^{\leq \reg}\pmb\psi](\tau_1,\tau_2)
+\int_{\MM_{r\geq 10m}(\tau_1,\tau_2)} r^{-1}  |\pr^{\leq\reg -1}\F|^2\nn\\
&& +\sqrt{\EM[\pr^{\leq \reg-1}\pmb\psi](\tau_1,\tau_2)}\sqrt{\EM[\pr^{\leq \reg}\pmb\psi](\tau_1,\tau_2)},
\eeaa
which proves the induction hypothesis \eqref{eq:EM:AsympKerr:highMora:v1:middle:generalk:prnprttk-n:inductionhypothesis} with $n$ replaced by $n+1$. We thus deduce that  \eqref{eq:EM:AsympKerr:highMora:v1:middle:generalk:prnprttk-n:inductionhypothesis} holds for all integers $1\leq n\leq\reg -1$, and hence, choosing $n=\reg-1$, we deduce, for all $\reg\leq 14$,
\bea
\lab{eq:EM:AsympKerr:highMora:v1:middle:generalk:prk:final}
\M_{r\geq 10.5m}[\pr^{\leq \reg}\pmb\psi](\tau_1,\tau_2) 
&\les&\M_{r\geq 10m}[\pr_{\tt}^{\leq \reg}\pmb\psi](\tau_1,\tau_2)
+\int_{\MM_{{r\geq 10m}}(\tau_1,\tau_2)} r^{-1} | \pr^{\leq \reg-1}\F|^2\nn\\
&& +\sqrt{\EM[\pr^{\leq \reg-1}\pmb\psi](\tau_1,\tau_2)}\sqrt{\EM[\pr^{\leq \reg}\pmb\psi](\tau_1,\tau_2)}.
\eea

\noindent{\bf Step 2.} \textit{Proof of  \eqref{eq:higherorderenergymorawetzforlargerconditionalonsprtau}: energy part}. Let $\chi_2=\chi_2(r)$ be a smooth cut-off function such that  
\beaa
\chi_2(r)=1 \,\,\text{for } r\geq 11m, \qquad  \chi_2(r)=0 \,\,\text{for }r\leq 10.5m. 
\eeaa
Then, we have from \cite[inequality (3.29)]{MaSz24}  that 
\beaa
\nn&& \sup_{\tau\in[\tau_1,\tau_2]}\E[\pr^{\leq 1}(\chi_2\pmb\psi)](\tau) +D_1\sup_{\tau\in[\tau_1,\tau_2-2]}\int_{\MM(\tau, \tau+2)}r^{-1}|\pr_\tau\pr^{\leq 1}(\chi_2\pmb\psi)|^2\\
&\les& \EM_{r\geq 10m} [\pr_{\tau}^{\leq 1}\pmb\psi](\tau_1,\tau_2) +\sqrt{\EMF[\pmb\psi](\tau_1,\tau_2)}\sqrt{\EMF[\pr^{\leq 1}\pmb\psi](\tau_1,\tau_2)}\nn\\
&&+\sum_{i,j}\sup_{\tau\in[\tau_1,\tau_2-2]}\bigg|\int_{\MM(\tau,\tau+2)}\Re\Big(\ov{\pr_{\tt}\pr^{\leq 1}(\chi_2\psi_{ij})}\pr^{\leq 1}\big(\chi_2(N_{ij}-D_1r^{-1}\pr_{\tt}\psi_{ij})\big)\Big)\bigg|\nn\\
&&+{\sup_{\tau\in[\tau_1,\tau_2-2]} \int_{\MM_{r\geq {10.5 m}}(\tau,\tau+2)}}|\N|^2\nn\\
&\les&  \EM_{r\geq 10m} [\pr_{\tau}^{\leq 1}\pmb\psi](\tau_1,\tau_2) +\sqrt{\EMF[\pmb\psi](\tau_1,\tau_2)}\sqrt{\EMF[\pr^{\leq 1}\pmb\psi](\tau_1,\tau_2)}\nn\\
&&+\Big( \EM_{r\geq 10.5m} [\pr_{\tau}^{\leq 1}\pmb\psi](\tau_1,\tau_2)\Big)^{\frac{1}{2}}\Big(\EM[\pr^{\leq 1}(\chi_2\pmb\psi)](\tau_1,\tau_2)\Big)^{\frac{1}{2}}\nn\\
&&
+\sup_{\tau\in[\tau_1,\tau_2-2]}\widehat{\mathcal{N}}_{r\geq {10.5m}}''[\pr^{\leq 1}\pmb\psi, \pr^{\leq 1}\F](\tau, \tau+2),
\eeaa
where in the second step we have plugged in the formula \eqref{eq:eqsforconditionalhighorderestimatelemma:general}  of $N_{ij}$ and used, in view of the definition \eqref{eq:defintionwidehatmathcalNfpsinormRHS} of $\widehat{\NN}[\cdot, \cdot](\cdot, \cdot)$, the following bound
\beaa
{\sup_{\tau\in[\tau_1,\tau_2-2]} \int_{\MM_{r\geq {10.5 m}}(\tau,\tau+2)}}|\F|^2\les \sup_{\tau\in[\tau_1,\tau_2-2]}\widehat{\mathcal{N}}_{r\geq {10.5m}}''[\pr^{\leq 1}\pmb\psi, \pr^{\leq 1}\F](\tau, \tau+2).
\eeaa
This yields
\beaa
\nn&&\sup_{\tau\in[\tau_1,\tau_2]}\E_{r\geq 11 m}[\pr^{\leq 1}\pmb\psi](\tau) +D_1\sup_{\tau\in[\tau_1,\tau_2-2]}\int_{\MM_{r\geq 11m}(\tau, \tau+2)}r^{-1}|\pr_\tau\pr^{\leq 1}\pmb\psi|^2\\
&\les&  \EM_{r\geq 10m} [\pr_{\tau}^{\leq 1}\pmb\psi](\tau_1,\tau_2) +\sqrt{\EMF[\pmb\psi](\tau_1,\tau_2)}\sqrt{\EMF[\pr^{\leq 1}\pmb\psi](\tau_1,\tau_2)}\nn\\
&&+\Big( \EM_{r\geq 10.5m} [\pr_{\tau}^{\leq 1}\pmb\psi](\tau_1,\tau_2)\Big)^{\frac{1}{2}}\Big(\M_{r\geq 10.5m}[\pr^{\leq 1}\pmb\psi](\tau_1,\tau_2)\Big)^{\frac{1}{2}}\nn\\
&&
+\sup_{\tau\in[\tau_1,\tau_2-2]}\widehat{\mathcal{N}}_{r\geq {10.5m}}''[\pr^{\leq 1}\pmb\psi, \pr^{\leq 1}\F](\tau, \tau+2),
\eeaa
and, using the estimate \eqref{eq:EM:AsympKerr:highMora:v1:middle:generalk:prk:final} to control the term $\M_{r\geq 10.5 m} [\pr^{\leq 1}\psi](\tau_1,\tau_2)$, we infer
\bea
\lab{eq:highorderEMF:energylarger}
\nn&&\sup_{\tau\in[\tau_1,\tau_2]}\E_{r\geq 11 m}[\pr^{\leq 1}\pmb\psi](\tau) +D_1\sup_{\tau\in[\tau_1,\tau_2-2]}\int_{\MM_{r\geq 11m}(\tau, \tau+2)}r^{-1}|\pr_\tau\pr^{\leq 1}\pmb\psi|^2\\
&\les&  \EM_{r\geq 10m} [\pr_{\tau}^{\leq 1}\pmb\psi](\tau_1,\tau_2) +\sqrt{\EMF[\pmb\psi](\tau_1,\tau_2)}\sqrt{\EMF[\pr^{\leq 1}\pmb\psi](\tau_1,\tau_2)}\nn\\
&&
+\sup_{\tau\in[\tau_1,\tau_2-2]}\widehat{\mathcal{N}}_{r\geq {10.5m}}''[\pr^{\leq 1}\pmb\psi, \pr^{\leq 1}\F](\tau, \tau+2)
+\int_{\MM_{{r\geq 10m}}(\tau_1,\tau_2)}  r^{-1}|\F|^2.
\eea

Next, in view of the commutators \eqref{esti:commutatorBoxgandT:general:0} \eqref{eq:localwavecommutators:withfirstordergoodderis:0}, as well as \eqref{eq:controloflinearizedinversemetriccoefficients}, we notice that the term $F_{ij}^{(\reg_1,\reg_2)}$ appearing in \eqref{eq:eqsforconditionalhighorderestimatelemma:general:allprreg} has in fact the following more precise structure 
\bea\lab{eq:eqsforconditionalhighorderestimatelemma:general:allprreg:moreprecisestructure}
\nn F_{ij}^{(\reg_1,\reg_2)}&=&
(\pr_r, r^{-1}\pr_{x^a})^{\reg_1}\pr_{\tt}^{\reg_2}F_{ij}+\dk^{\leq \reg_1+\reg_2}\Ga_b\c\pr_r(r\pr_r, \pr_{x^a}) \pr^{\leq \reg_1+\reg_2-1}\psi_{kl}\\
&&+\dk^{\leq \reg_1+\reg_2+1}\Ga_g\c \dk \pr^{\leq \reg_1+\reg_2}\psi_{kl}.
\eea
Thus, applying \eqref{eq:highorderEMF:energylarger} to the wave equations \eqref{eq:eqsforconditionalhighorderestimatelemma:general:allprreg}, proceeding by induction as in Step 1, and taking \eqref{eq:eqsforconditionalhighorderestimatelemma:general:allprreg:moreprecisestructure} into account,  we deduce, for all $\reg\leq 14$,
\begin{align}\lab{eq:highorderEMF:energylarger:prreg:notyetdone}
&\sup_{\tau\in[\tau_1,\tau_2]}\E_{r\geq 11m}[\pr^{\leq \reg}\pmb\psi](\tau) +D_1\sup_{\tau\in[\tau_1,\tau_2-2]}\int_{\MM_{r\geq 11m}(\tau, \tau+2)}r^{-1}|\pr_\tau\pr^{\leq\reg}\pmb\psi|^2\nn\\
\les{}&   \EM_{r\geq 10m} [\pr_{\tau}^{\leq \reg}\pmb\psi](\tau_1,\tau_2) 
+\sqrt{\EMF[\pr^{\leq \reg-1}\pmb\psi](\tau_1,\tau_2)}\sqrt{\EMF[\pr^{\leq \reg}\pmb\psi](\tau_1,\tau_2)}\nn\\
&+\sup_{\tau\in[\tau_1,\tau_2-2]}\widehat{\mathcal{N}}_{r\geq {10.5m}}''\Big[\pr^{\leq \reg}\pmb\psi, \pr^{\leq \reg}\Ga_b\c\pr_r(r\pr_r, \pr_{x^a}) \pr^{\leq \reg-1}\pmb\psi+\dk^{\leq\reg+1}\Ga_g\c\dk \pr^{\leq \reg}\pmb\psi\Big](\tau, \tau+2)\nn\\
&+\sup_{\tau\in[\tau_1,\tau_2-2]}\widehat{\mathcal{N}}_{r\geq {10.5m}}''[\pr^{\leq \reg}\pmb\psi, \pr^{\leq \reg}\F](\tau, \tau+2)
+\int_{\MM_{{r\geq 10m}}(\tau_1,\tau_2)}   r^{-1}| \pr^{\leq \reg-1}\F|^2,
\end{align}
where we have used again the estimate \eqref{eq:highorderrecovering:errorforprreg}. Then, we estimate the last term on the RHS of \eqref{eq:highorderEMF:energylarger:prreg:notyetdone} as follows, with $\chi_2=\chi_2(r)$ as above, 
\beaa
&&\sup_{\tau\in[\tau_1,\tau_2-2]}\widehat{\mathcal{N}}_{r\geq {10.5m}}''\Big[\pr^{\leq \reg}\pmb\psi, \pr^{\leq\reg}\Ga_b\c\pr_r(r\pr_r, \pr_{x^a}) \pr^{\leq \reg-1}\pmb\psi+\dk^{\leq\reg+1}\Ga_g\c\dk \pr^{\leq \reg}\pmb\psi\Big](\tau, \tau+2)\\
&\les& \left|\int_{\MM(\tau, \tau+2)}\chi_2(r)\pr_\tau\pr^{\leq \reg}\pmb\psi \pr^{\leq \reg}\Ga_b\c\pr_r(r\pr_r, \pr_{x^a}) \pr^{\leq \reg-1}\pmb\psi\right|\\
&&+\ep\M_{10.5m,11m} [\pr^{\leq \reg}\pmb\psi](\tau_1,\tau_2) +\ep\EM_{r\geq 11m} [\pr^{\leq \reg}\pmb\psi](\tau_1,\tau_2) 
\eeaa
and integrating by parts first in $\pr_r$ and then in $\pr_\tau$, we infer
\bea\lab{eq:controlsuptauintau1tau2ofwidehatNNrgeq10point5mofNLtermgeneratedcommutationswaveeq:zz}
\nn&&\sup_{\tau\in[\tau_1,\tau_2-2]}\widehat{\mathcal{N}}_{r\geq {10.5m}}''\Big[\pr^{\leq \reg}\pmb\psi, \pr^{\leq\reg}\Ga_b\c\pr_r(r\pr_r, \pr_{x^a}) \pr^{\leq \reg-1}\pmb\psi+\dk^{\leq\reg+1}\Ga_g\c\dk \pr^{\leq \reg}\pmb\psi\Big](\tau, \tau+2)\\
\nn&\les& \left|\int_{\MM(\tau, \tau+2)}\chi_2(r)\pr_r\pr^{\leq \reg}\pmb\psi \pr^{\leq \reg}\Ga_b\c(r\pr_r, \pr_{x^a})\pr^{\leq \reg}\pmb\psi\right|\\
\nn&&+\ep\M_{10.5m,11m} [\pr^{\leq \reg}\pmb\psi](\tau_1,\tau_2) +\ep\EM_{r\geq 11m} [\pr^{\leq \reg}\pmb\psi](\tau_1,\tau_2)\\
&\les& \ep\M_{10.5m,11m} [\pr^{\leq \reg}\pmb\psi](\tau_1,\tau_2) +\ep\EM_{r\geq 11m} [\pr^{\leq \reg}\pmb\psi](\tau_1,\tau_2).  
\eea
Plugging in \eqref{eq:highorderEMF:energylarger:prreg:notyetdone}, we deduce, for $\ep$ small enough, 
\beaa
&&\sup_{\tau\in[\tau_1,\tau_2]}\E_{r\geq 11m}[\pr^{\leq \reg}\pmb\psi](\tau) +D_1\sup_{\tau\in[\tau_1,\tau_2-2]}\int_{\MM_{r\geq 11m}(\tau, \tau+2)}r^{-1}|\pr_\tau\pr^{\leq\reg}\pmb\psi|^2\\
&\les&   \EM_{r\geq 10m} [\pr_{\tau}^{\leq \reg}\pmb\psi](\tau_1,\tau_2) 
+\sqrt{\EMF[\pr^{\leq \reg-1}\pmb\psi](\tau_1,\tau_2)}\sqrt{\EMF[\pr^{\leq \reg}\pmb\psi](\tau_1,\tau_2)}\nn\\
\nn&&
+\sup_{\tau\in[\tau_1,\tau_2-2]}\widehat{\mathcal{N}}_{r\geq {10.5m}}''[\pr^{\leq \reg}\pmb\psi, \pr^{\leq \reg}\F](\tau, \tau+2)
+\int_{\MM_{{r\geq 10m}}(\tau_1,\tau_2)}   r^{-1}| \pr^{\leq \reg-1}\F|^2\\
\nn&&+\ep\M_{r\geq 10.5m} [\pr^{\leq \reg}\pmb\psi](\tau_1,\tau_2),
\eeaa
and, using the estimate \eqref{eq:EM:AsympKerr:highMora:v1:middle:generalk:prk:final} to control the term $\M_{r\geq 10.5 m} [\pr^{\leq\reg}\pmb\psi](\tau_1,\tau_2)$, we infer
\bea
\lab{eq:highorderEMF:energylarger:prreg}
&&\sup_{\tau\in[\tau_1,\tau_2]}\E_{r\geq 11m}[\pr^{\leq \reg}\pmb\psi](\tau) +D_1\sup_{\tau\in[\tau_1,\tau_2-2]}\int_{\MM_{r\geq 11m}(\tau, \tau+2)}r^{-1}|\pr_\tau\pr^{\leq\reg}\pmb\psi|^2\nn\\
&\les&   \EM_{r\geq 10m} [\pr_{\tau}^{\leq \reg}\pmb\psi](\tau_1,\tau_2) 
+\sqrt{\EMF[\pr^{\leq \reg-1}\pmb\psi](\tau_1,\tau_2)}\sqrt{\EMF[\pr^{\leq \reg}\pmb\psi](\tau_1,\tau_2)}\nn\\
&&
+\sup_{\tau\in[\tau_1,\tau_2-2]}\widehat{\mathcal{N}}_{r\geq {10.5m}}''[\pr^{\leq \reg}\pmb\psi, \pr^{\leq \reg}\F](\tau, \tau+2)
+\int_{\MM_{{r\geq 10m}}(\tau_1,\tau_2)}   r^{-1}| \pr^{\leq \reg-1}\F|^2.
\eea

\noindent{\bf Step 3.} \textit{Proof of  \eqref{eq:higherorderenergymorawetzforlargerconditionalonsprtau}: flux part}. We have from \cite[inequality (3.31)]{MaSz24} that 
\bea\lab{eq:highorderEMF:fluxnullinf:final}
&&\F_{\II_+}[\pr^{\leq 1}\pmb\psi](\tau_1, \tau_2)\nn\\
&\les& \EMF_{r\geq 10m}[\pr_{\tau}^{\leq 1}\pmb\psi](\tau_1,\tau_2)+\sqrt{\EMF[\pmb\psi](\tau_1,\tau_2)}\sqrt{\EMF[\pr^{\leq 1}\pmb\psi](\tau_1,\tau_2)}+\int_{\II_+(\tau_1,\tau_2)}|\N|^2 \nn\\
&&+\sup_{\tau\in[\tau_1,\tau_2-2]}\widehat{\mathcal{N}}_{r\geq {10.5m}}''[\pr^{\leq1}\pmb\psi, \pr^{\leq1 }\N](\tau, \tau+2)\nn\\
&\les& \EMF_{r\geq 10m}[\pr_{\tau}^{\leq 1}\pmb\psi](\tau_1,\tau_2)+\sqrt{\EMF[\pmb\psi](\tau_1,\tau_2)}\sqrt{\EMF[\pr^{\leq 1}\pmb\psi](\tau_1,\tau_2)}+\int_{\II_+(\tau_1,\tau_2)}|\F|^2 \nn\\
&&+\sup_{\tau\in[\tau_1,\tau_2-2]}\widehat{\mathcal{N}}_{r\geq {10.5m}}''[\pr^{\leq1}\pmb\psi, \pr^{\leq1 }\F](\tau, \tau+2) +D_1\sup_{\tau\in[\tau_1,\tau_2-2]}\int_{\MM_{r\geq 11m}(\tau, \tau+2)}r^{-1}|\pr_\tau\pr^{\leq 1}\pmb\psi|^2\nn\\
&&+\M_{10.5m,11m} [\pr^{\leq 1}\pmb\psi](\tau_1,\tau_2) +\EM_{r\geq 11m} [\pr^{\leq 1}\pmb\psi](\tau_1,\tau_2)\nn\\
&\les& \EMF_{r\geq 10m}[\pr_{\tau}^{\leq 1}\pmb\psi](\tau_1,\tau_2)+\sqrt{\EMF[\pmb\psi](\tau_1,\tau_2)}\sqrt{\EMF[\pr^{\leq 1}\pmb\psi](\tau_1,\tau_2)} +\int_{\II_+(\tau_1,\tau_2)}|\F|^2\nn\\
&& +\int_{\MM_{r\geq 10m}(\tau_1, \tau_2)}r^{-1}|\F|^2 +\sup_{\tau\in[\tau_1,\tau_2-2]}\widehat{\mathcal{N}}_{r\geq {10.5m}}''[\pr^{\leq1}\pmb\psi, \pr^{\leq1 }\F](\tau, \tau+2),
\eea
where in the second step we have used the formula   \eqref{eq:eqsforconditionalhighorderestimatelemma:general} of $N_{ij}$, and where in the last step we have used   \eqref{eq:EM:AsympKerr:highMora:v1:middle:generalk:prk:final} with $\reg=1$ and \eqref{eq:highorderEMF:energylarger}. Again, applying \eqref{eq:highorderEMF:fluxnullinf:final} to the wave equations \eqref{eq:eqsforconditionalhighorderestimatelemma:general:allprreg}, proceeding by induction as in Step 1, and taking \eqref{eq:controlsuptauintau1tau2ofwidehatNNrgeq10point5mofNLtermgeneratedcommutationswaveeq:zz} into account,  we deduce, for all $\reg\leq 14$,
\beaa
\F_{\II_+}[\pr^{\leq \reg}\pmb\psi](\tau_1, \tau_2)
&\les& \EMF_{r\geq 10m}[\pr_{\tau}^{\leq \reg}\pmb\psi](\tau_1,\tau_2)
+ \int_{\II_+(\tau_1,\tau_2)}|\pr^{\leq\reg-1}\F|^2\nn\\
&&+\sqrt{\EMF[\pr^{\leq \reg-1}\pmb\psi](\tau_1,\tau_2)}\sqrt{\EMF[\pr^{\leq \reg}\pmb\psi](\tau_1,\tau_2)} \nn\\
&&+\sup_{\tau\in[\tau_1,\tau_2-2]}\widehat{\mathcal{N}}_{r\geq {10.5m}}''[\pr^{\leq \reg}\pmb\psi, \pr^{\leq \reg}\F](\tau, \tau+2)\nn\\
&&
+\int_{\MM_{r\geq 10m}(\tau_1,\tau_2)}r^{-1}|\pr^{\leq\reg-1}\F|^2\\
&& +\ep\M_{10.5m,11m} [\pr^{\leq \reg}\pmb\psi](\tau_1,\tau_2) +\ep\EM_{r\geq 11m} [\pr^{\leq \reg}\pmb\psi](\tau_1,\tau_2),
\eeaa
where we have used 
\beaa
\sum_{k,l}\int_{\II_+(\tau_1,\tau_2)}  | \dk^{\leq \reg}\Ga_b\c\dk\pr^{\leq \reg-1}\psi_{kl}|^2
&\les&\sum_{k,l}\int_{\MM_{r\geq 11m}(\tau_1,\tau_2)}  | \pr_r^{\leq 1}\big(\dk^{\leq \reg}\Ga_b\c\dk\pr^{\leq \reg-1}\psi_{kl}\big)|^2
\nn\\
&\les&\sum_{k,l}\int_{\MM_{r\geq 11m}(\tau_1,\tau_2)}  | \dk^{\leq \reg+1}\Ga_b\c\dk\pr^{\leq \reg}\psi_{kl}|^2
\nn\\
&\les& \ep \sup_{\tau\in[\tau_1,\tau_2]}\E_{r\geq 11m}[\pr^{\leq \reg}\pmb\psi](\tau),
\eeaa
in which we have applied a trace estimate in the first step. Together with \eqref{eq:EM:AsympKerr:highMora:v1:middle:generalk:prk:final} and \eqref{eq:highorderEMF:energylarger:prreg}, and using a trace estimate to control the integral of $|\pr^{\leq\reg-1}\F|^2$ on $\II_+$, we infer, for all $\reg\leq 14$,
\bea\lab{eq:highorderEMF:fluxnullinf:final:prreg}
\F_{\II_+}[\pr^{\leq \reg}\pmb\psi](\tau_1, \tau_2)
&\les& \EMF_{r\geq 10m}[\pr_{\tau}^{\leq \reg}\pmb\psi](\tau_1,\tau_2) +\int_{\MM_{r\geq 10m}(\tau_1,\tau_2)}|\pr^{\leq \reg}\F|^2\nn\\
&&+\sqrt{\EMF[\pr^{\leq \reg-1}\pmb\psi](\tau_1,\tau_2)}\sqrt{\EMF[\pr^{\leq \reg}\pmb\psi](\tau_1,\tau_2)} \nn\\
&&+\sup_{\tau\in[\tau_1,\tau_2-2]}\widehat{\mathcal{N}}_{r\geq {10m}}''[\pr^{\leq \reg}\pmb\psi, \pr^{\leq \reg}\F](\tau, \tau+2).
\eea

{\noindent{\bf Step 4.} \textit{End of the proof of  \eqref{eq:higherorderenergymorawetzforlargerconditionalonsprtau}}}. Combining the estimates \eqref{eq:EM:AsympKerr:highMora:v1:middle:generalk:prk:final}, \eqref{eq:highorderEMF:energylarger:prreg} and \eqref{eq:highorderEMF:fluxnullinf:final:prreg},  we deduce
\beaa
\nn&&\EMF_{r\geq 11m}[\pr^{\leq \reg}\pmb\psi](\tau_1, \tau_2) \nn\\
&\les& \EMF[\pr_\tau^{\leq \reg}\pmb\psi](\tau_1, \tau_2)+\big(\EMF[\pr^{\leq \reg-1}\pmb\psi](\tau_1, \tau_2)\big)^{\frac{1}{2}}\left(\EMF[\pr^{\leq \reg}\pmb\psi](\tau_1, \tau_2)\right)^{\frac{1}{2}}\nn\\
&&+\sup_{\tau\in[\tau_1,\tau_2-2]}\widehat{\mathcal{N}}_{r\geq {10m}}''[\pr^{\leq \reg}\pmb\psi, \pr^{\leq \reg}\F](\tau, \tau+2)+\int_{\MM_{r\geq 10m}(\tau_1,\tau_2)}|\pr^{\leq  \reg}\F|^2,
\eeaa
which proves
the desired estimate \eqref{eq:higherorderenergymorawetzforlargerconditionalonsprtau}.

\noindent{\bf Step 5.} \textit{Proof of  \eqref{eq:higherorderenergymorawetzconditionalonsprtausprphis}: Morawetz part and flux part on $\AA$}. We have from \cite[equation (3.33)]{MaSz24}  that 
\bea\lab{eq:elliptic:waveform:outsidehorizon}
\Big(\De\pr_r^2+\mathring{\ga}^{ab}\pr_{x^a}\pr_{x^b}\Big)\psi_{ij} = |q|^2N_{ij} +O(1)(\pr_\tau, \pr_{\tphi})(\pr_\tau, \pr_{\tphi},\pr_r)\psi_{ij} +O(1)\pr\psi+O(\ep)\pr^2\psi_{ij} .
\eea
As in the discussions below  \cite[equation (3.33)]{MaSz24}, we multiply both sides of \eqref{eq:elliptic:waveform:outsidehorizon} by $\chi_4^2\ov{\pr^2_r\psi}$ where
\beaa
\chi_4(r)=1 \,\,\text{for }r_+(1+\dred)\leq r\leq 11m, \qquad  \chi_4(r)=0 \,\,\text{for }r\geq 12m\,\,\textrm{and }r\leq r_+(1+\dred/2),
\eeaa
take the real part and integrate over $\MM(\tau_1, \tau_2)$ to deduce, after performing integration by parts on the RHS, 
\bea
&&\sum_{i,j}\int_{\MM(\tau_1, \tau_2)}\Re\Big(\chi_4^2\big(\De\pr_r^2+\mathring{\ga}^{ab}\pr_{x^a}\pr_{x^b}\big)\psi_{ij}\ov{\pr^2_r\psi_{ij}}\Big)\nn\\
&\les& \left(\M[\pr_\tau^{\leq 1}\pmb\psi](\tau_1, \tau_2)+\M_{r\leq 12m}[\pr_{\tphi}\pmb\psi](\tau_1, \tau_2)\right)^{\frac{1}{2}}\Big(\M[\chi_4\pr_r\pmb\psi](\tau_1, \tau_2)\Big)^{\frac{1}{2}}\nn\\
&&+\sqrt{\EMF[\pmb\psi](\tau_1, \tau_2)}\sqrt{\EMF[\pr^{\leq 1}\pmb\psi](\tau_1, \tau_2)}+\ep\M[\pr^{\leq 1}\pmb\psi](\tau_1, \tau_2)\nn\\
&& +\sum_{i,j}\left|\int_{\MM(\tau_1, \tau_2)}\Re(\chi_4^2 \qs N_{ij} \ov{\pr^2_r\psi_{ij}})\right|,
\eea
an estimate similar to \cite[inequality (3.34)]{MaSz24}, the only difference being that we keep the integral of $\Re(\chi_4^2 \qs N_{ij} \ov{\pr^2_r\psi_{ij}})$ as it is. Then, integrating by parts the LHS, we deduce the following estimate similar to \cite[inequality (3.35)]{MaSz24}:
\bea
\lab{eq:Mnormofpr_rpsi:infiniterleq11m:pre}
\M[\chi_4\pr_r\pmb\psi](\tau_1, \tau_2)
&\les&\M[\pr_\tau^{\leq 1}\pmb\psi](\tau_1, \tau_2)+\M_{r\leq 11m}[\pr_{\tphi}\pmb\psi](\tau_1, \tau_2)+\ep\M[\pr^{\leq 1}\pmb\psi](\tau_1, \tau_2)\nn\\
&&+\M_{r\geq 11m}[\pr^{\leq 1}\pmb\psi](\tau_1, \tau_2)+\sqrt{\EMF[\pmb\psi](\tau_1, \tau_2)}\sqrt{\EMF[\pr^{\leq 1}\pmb\psi](\tau_1, \tau_2)}\nn\\
&& +\sum_{i,j}\left|\int_{\MM(\tau_1, \tau_2)}\Re(\chi_4^2 \qs N_{ij} \ov{\pr^2_r\psi_{ij}})\right|\nn\\
&\les&\M[(\pr_\tau, \chi_0\pr_{\tphi})^{\leq 1}\pmb\psi](\tau_1, \tau_2)
+\M_{r\geq 11m}[\pr^{\leq 1}\pmb\psi](\tau_1, \tau_2)+\ep\M[\pr^{\leq 1}\pmb\psi](\tau_1, \tau_2)\nn\\
&&+\sqrt{\EMF[\pmb\psi](\tau_1, \tau_2)}\sqrt{\EMF[\pr^{\leq 1}\pmb\psi](\tau_1, \tau_2)}
\nn\\
&& +\sum_{i,j}\left|\int_{\MM(\tau_1, \tau_2)}\Re(\chi_4^2 \qs N_{ij} \ov{\pr^2_r\psi_{ij}})\right|.
\eea
In view of the formula of $N_{ij}$ in \eqref{eq:eqsforconditionalhighorderestimatelemma:general}, we have
\beaa
 \sum_{i,j}\left|\int_{\MM(\tau_1, \tau_2)}\Re(\chi_4^2 \qs N_{ij} \ov{\pr^2_r\psi_{ij}})\right| &\les& \left(\int_{\MM_{r\leq 12m}(\tau_1, \tau_2)}|\F|^2\right)^{\frac{1}{2}}\Big(\M[\chi_4\pr_r\pmb\psi](\tau_1, \tau_2)\Big)^{\frac{1}{2}}\\
 && + \sqrt{\M[\pmb\psi](\tau_1, \tau_2)}\sqrt{\M[\pr^{\leq 1}\pmb\psi](\tau_1, \tau_2)},
 \eeaa
 where we have used integration by parts in $\pr_r$ for the second  term on the RHS, and plugging this  into \eqref{eq:Mnormofpr_rpsi:infiniterleq11m:pre}, we infer
 \begin{align}
\lab{eq:Mnormofpr_rpsi:infiniterleq11m}
 \M_{r_+(1+\dred), 11m}[\pr_r\pmb\psi](\tau_1, \tau_2)
\les{}& \int_{\MM_{r\leq 12m}(\tau_1, \tau_2)}|\F|^2 +\ep\M[\pr^{\leq 1}\pmb\psi](\tau_1, \tau_2)
\nn\\
&+\M[(\pr_\tau, \chi_0\pr_{\tphi})^{\leq 1}\pmb\psi](\tau_1, \tau_2)+\M_{r\geq 11m}[\pr^{\leq 1}\pmb\psi](\tau_1, \tau_2)\nn\\
&+\sqrt{\EMF[\pmb\psi](\tau_1, \tau_2)}\sqrt{\EMF[\pr^{\leq 1}\pmb\psi](\tau_1, \tau_2)}.
\end{align}

Next, we have from the first step in the derivation of \cite[inequality (3.36)]{MaSz24}  that
\beaa
&& \M_{r_+(1+\dred), 11m}[\pr^{\leq 1}\pmb\psi](\tau_1, \tau_2)\nn\\
&\les& \int_{\Mntrap_{r_+(1+\dred), 11m}(\tau_1, \tau_2)}|\N|^2+\M[(\pr_\tau, \chi_0\pr_{\tphi})^{\leq 1}\pmb\psi](\tau_1, \tau_2)
+\M_{r_+(1+\dred), 11m}[\pr_r\pmb\psi](\tau_1, \tau_2)\nn\\
&&+\M_{r\geq 11m}[\pr^{\leq 1}\pmb\psi](\tau_1, \tau_2) +\sqrt{\EMF[\pmb\psi](\tau_1, \tau_2)}\sqrt{\EMF[\pr^{\leq 1}\pmb\psi](\tau_1, \tau_2)},
\eeaa
and using the estimate \eqref{eq:Mnormofpr_rpsi:infiniterleq11m} and the formula  of \eqref{eq:eqsforconditionalhighorderestimatelemma:general} $N_{ij}$, we infer
\bea
\lab{eq:recoveringderivatives:firstorder:Mora:middlestep}
&& \M_{r_+(1+\dred), 11m}[\pr^{\leq 1}\pmb\psi](\tau_1, \tau_2)\nn\\
&\les& \int_{\MM_{r\leq 12m}(\tau_1, \tau_2)}|\F|^2  
+\M[(\pr_\tau, \chi_0\pr_{\tphi})^{\leq 1}\pmb\psi](\tau_1, \tau_2)+\M_{r\geq 11m}[\pr^{\leq 1}\pmb\psi](\tau_1, \tau_2)\nn\\
&&+\ep\M_{r\leq r_+(1+\dred)}[\pr^{\leq 1}\pmb\psi](\tau_1, \tau_2)+\sqrt{\EMF[\pmb\psi](\tau_1, \tau_2)}\sqrt{\EMF[\pr^{\leq 1}\pmb\psi](\tau_1, \tau_2)}.
\eea
 
Next, arguing as for the proof of \eqref{eq:eqsforconditionalhighorderestimatelemma:general:allprreg}, we have
\begin{equation}
\lab{eq:eqsforconditionalhighorderestimatelemma:general:allprreg:withprttprtphi}
\bsplit
&\square_{\g}(\pr^{\reg_1}(\pr_{\tt},\chi_0\pr_{\tphi})^{\reg_2}\psi_{ij})=D_1r^{-1}\pr_{\tt}\pr^{\leq \reg_1}(\pr_{\tt},\chi_0\pr_{\tphi})^{\leq \reg_2}\psi_{ij}\\
&\qquad \qquad \qquad\qquad \qquad \quad\,\,+O(r^{-2})\dk^{\leq 1}\pr^{\leq \reg_1}(\pr_{\tt},\chi_0\pr_{\tphi})^{\leq \reg_2}\psi_{kl} +\widehat{F}_{ij}^{(\reg_1,\reg_2)},\\
&\widehat{F}_{ij}^{(\reg_1,\reg_2)}=
\pr^{\reg_1}\pr_{\tt}^{\reg_2}F_{ij}
+\dk^{\leq \reg_1+\reg_2+1}\Ga_b\c\dk \pr^{\leq \reg_1+\reg_2}\psi_{kl}
+ O(1)\mathbf{1}_{11m\leq r\leq 12m}\pr^{\reg_1+\reg_2+1}\psi_{kl},
\end{split}
\end{equation}
where we also used the commutator relation
\beaa
[ \chi_0\pr_{\tphi}, \square_{\g}]\psi = \mathbf{1}_{r\leq 12m}r\dk^{\leq 2}\Ga_g\c \dk\pr^{\leq 1}\psi
+O(1)\mathbf{1}_{11m\leq r\leq 12m}\pr^{\leq 2}\psi_{kl}.
\eeaa
Thus, applying the estimate \eqref{eq:recoveringderivatives:firstorder:Mora:middlestep} to the wave equations \eqref{eq:eqsforconditionalhighorderestimatelemma:general:allprreg:withprttprtphi} and proceeding by induction as in Step 1, we infer, for $\reg\leq 14$,
\bea\lab{eq:controlofMorawetzrplus1plusdehto11mforregderivativesofpmbpsiconditionalonprtauandchi0prphit}
&&\M_{r_+(1+\dred), 11m}[\pr^{\leq \reg}\pmb\psi](\tau_1, \tau_2)\\
&\les& \int_{\MM_{r\leq 12m}(\tau_1, \tau_2)}|\pr^{\leq \reg-1}\F|^2 
+\M[(\pr_\tau, \chi_0\pr_{\tphi})^{\leq \reg}\pmb\psi](\tau_1, \tau_2)+\M_{r\geq 11m}[\pr^{\leq \reg}\pmb\psi](\tau_1, \tau_2)\nn\\
\nn&&+\ep\M_{r\leq r_+(1+\dred)}[\pr^{\leq \reg}\pmb\psi](\tau_1, \tau_2)+\sqrt{\EMF[\pr^{\leq \reg-1}\pmb\psi](\tau_1, \tau_2)}\sqrt{\EMF[\pr^{\leq \reg}\pmb\psi](\tau_1, \tau_2)}.
\eea
Together with the assumed red-shift estimate \eqref{eq:assum:redshift:highorderunweightedEMF:generalwave},  we deduce, for $\ep$ small enough, an analogous estimate to \cite[inequality (3.37)]{MaSz24}:
\beaa
&&\sup_{\tau\in[\tau_1,\tau_2]}\E_{r\leq r_+ (1+\dred)}[\pr^{\leq \reg}\pmb\psi](\tau)
+\F_{\AA}[\pr^{\leq \reg}\pmb\psi](\tau_1, \tau_2)+ \M_{r\leq 11m}[\pr^{\leq \reg}\pmb\psi](\tau_1, \tau_2)\nn\\
&\les& \E[\pr^{\leq \reg}\pmb\psi](\tau_1)+\M[(\pr_\tau, \chi_0\pr_{\tphi})^{\leq \reg}\pmb\psi](\tau_1, \tau_2)+\M_{r\geq 11m}[\pr^{\leq \reg}\pmb\psi](\tau_1, \tau_2)\nn\\
&&+\sqrt{\EMF[\pr^{\leq \reg-1}\pmb\psi](\tau_1, \tau_2)}\sqrt{\EMF[\pr^{\leq \reg}\pmb\psi](\tau_1, \tau_2)}+\int_{\MM(\tau_1, \tau_2)}|\pr^{\leq \reg}\F|^2.
\eeaa
Plugging the control for $\M_{r\geq 11m}[\pr^{\leq \reg}\pmb\psi](\tau_1, \tau_2)$ in \eqref{eq:EM:AsympKerr:highMora:v1:middle:generalk:prk:final}, we deduce 
\bea\lab{eq:highorderMora:improveMora:v3}
&&\sup_{\tau\in[\tau_1,\tau_2]}\E_{r\leq r_+ (1+\dred)}[\pr^{\leq \reg}\pmb\psi](\tau)
+\F_{\AA}[\pr^{\leq \reg}\pmb\psi](\tau_1, \tau_2)+ \M[\pr^{\leq \reg}\pmb\psi](\tau_1, \tau_2)\nn\\
&\les& \E[\pr^{\leq \reg}\pmb\psi](\tau_1)+\M[(\pr_\tau, \chi_0\pr_{\tphi})^{\leq \reg}\pmb\psi](\tau_1, \tau_2)+\int_{\MM(\tau_1, \tau_2)}|\pr^{\leq \reg}\F|^2
\nn\\
&&+\sqrt{\EMF[\pr^{\leq \reg-1}\pmb\psi](\tau_1, \tau_2)}\sqrt{\EMF[\pr^{\leq \reg}\pmb\psi](\tau_1, \tau_2)}.
\eea

\noindent{\bf Step 6.} \textit{Proof of  \eqref{eq:higherorderenergymorawetzconditionalonsprtausprphis}: energy and Morawetz parts}. Applying the estimate (3.39) in \cite{MaSz24}, we have
\beaa
\E_{r_+(1+\dred), 11m}[\pr^{\leq 1}\pmb\psi](\tau) &\les& {\int_{\Si(\tau)}|\N|^2}+\E[\pr_\tau\pmb\psi](\tau)+\E_{r_+(1+\dred/2), 11m}[\pr_{\tphi}\pmb\psi](\tau)+\E[\pmb\psi](\tau)\nn\\
&&+\E_{r\geq 11m}[\pr^{\leq 1}\pmb\psi](\tau)+\ep^2\E[\pr^{\leq 1}\pmb\psi](\tau)+\sqrt{\E[\pmb\psi](\tau)}\sqrt{\E[\pr^{\leq 1}\pmb\psi](\tau)}\nn\\
&\les& {\int_{\Si(\tau)}|\F|^2}+\E[(\pr_\tau,\chi_0\pr_{\tphi})^{\leq 1}\pmb\psi](\tau)+\E_{r\geq 11m}[\pr^{\leq 1}\pmb\psi](\tau)\nn\\
&&+\ep^2\E[\pr^{\leq 1}\pmb\psi](\tau)+\sqrt{\E[\psi](\tau)}\sqrt{\E[\pr^{\leq 1}\pmb\psi](\tau)},
\eeaa
where in the last step we have used the formula \eqref{eq:eqsforconditionalhighorderestimatelemma:general} of $N_{ij}$. 
Next, by inductively applying the above estimate to the wave equations \eqref{eq:eqsforconditionalhighorderestimatelemma:general:allprreg:withprttprtphi} as in Step 1, we infer
\beaa
\E_{r_+(1+\dred), 11m}[\pr^{\leq \reg}\pmb\psi](\tau) 
&\les& {\int_{\Si(\tau)}|\pr^{\leq \reg-1}\F|^2}+\E[(\pr_\tau,\chi_0\pr_{\tphi})^{\leq \reg}\pmb\psi](\tau)+\ep^2\E[\pr^{\leq \reg}\pmb\psi](\tau)\nn\\
&&
+\E_{r\geq 11m}[\pr^{\leq \reg}\pmb\psi](\tau)+\sqrt{\E[\pr^{\leq \reg-1}\pmb\psi](\tau)}\sqrt{\E[\pr^{\leq \reg}\pmb\psi](\tau)}.
\eeaa
Taking the supremum of $\tau\in[\tau_1, \tau_2]$, and using also \eqref{eq:highorderEMF:energylarger:prreg} and \eqref{eq:highorderMora:improveMora:v3}, 
 we infer, for $\ep$ small enough such that $\ep^2\E[\pr^{\leq \reg}\pmb\psi](\tau)$ is absorbed by the LHS, 
 \bea\lab{eq:highorder:energyneartrap:total}
\nn&&\EM[\pr^{\leq \reg}\pmb\psi](\tau_1, \tau_2)
+\F_{\AA}[\pr^{\leq \reg}\pmb\psi](\tau_1, \tau_2)\\
\nn&\les& \E[\pr^{\leq \reg}\pmb\psi](\tau_1)+\EM[(\pr_\tau, \chi_0\pr_{\tphi})^{\leq \reg}\pmb\psi](\tau_1, \tau_2)+\sup_{\tau\in[\tau_1,\tau_2-2]}\widehat{\mathcal{N}}_{r\geq {10m}}''[\pr^{\leq \reg}\pmb\psi, \pr^{\leq \reg}\F](\tau, \tau+2)\\
&&+\sqrt{\EMF[\pr^{\leq \reg-1}\pmb\psi](\tau_1, \tau_2)}\sqrt{\EMF[\pr^{\leq \reg}\pmb\psi](\tau_1, \tau_2)}+ \int_{\MM(\tau_1, \tau_2)}|\pr^{\leq \reg}\F|^2,
\eea
where we used a trace estimate to control the integral of $|\pr^{\leq \reg-1}\F|^2$ on $\Si(\tau)$.

{\noindent{\bf Step 7.} \textit{End of the proof of  \eqref{eq:higherorderenergymorawetzconditionalonsprtausprphis} and \eqref{eq:higherorderenergymorawetzconditionalonsprtausprphis:nontrapped}}}. 
Adding the flux estimate \eqref{eq:highorderEMF:fluxnullinf:final:prreg} and the energy-Morawetz estimate \eqref{eq:highorder:energyneartrap:total} together yields the desired estimate \eqref{eq:higherorderenergymorawetzconditionalonsprtausprphis}. 

In order to prove \eqref{eq:higherorderenergymorawetzconditionalonsprtausprphis:nontrapped}, we need to derive the analog of \eqref{eq:controlofMorawetzrplus1plusdehto11mforregderivativesofpmbpsiconditionalonprtauandchi0prphit} for $\widehat{\M}_{r_+(1+\dred), 11m}[\pr^{\leq \reg}\pmb\psi](\tau_1, \tau_2)$. To this end, we take the square of the modulus on both sides of \eqref{eq:elliptic:waveform:outsidehorizon}, multiply by $\chi_4^2$ with $\chi_4$ chosen as in Step 5, integrate over $\MM(\tau_1, \tau_2)$ and sum over $i,j$ to deduce, after performing integration by parts on the RHS, the following analog of \eqref{eq:recoveringderivatives:firstorder:Mora:middlestep}
\bea\lab{eq:recoveringderivatives:firstorder:Mora:middlestep:analogforwidehatcase}
&& \widehat{\M}_{r_+(1+\dred), 11m}[\pr^{\leq 1}\pmb\psi](\tau_1, \tau_2)\nn\\
&\les& \int_{\MM_{r\leq 12m}(\tau_1, \tau_2)}|\F|^2  
+\widehat{\M}[(\pr_\tau, \chi_0\pr_{\tphi})^{\leq 1}\pmb\psi](\tau_1, \tau_2)+\M_{r\geq 11m}[\pr^{\leq 1}\pmb\psi](\tau_1, \tau_2)\nn\\
&&+\ep\M_{r\leq r_+(1+\dred)}[\pr^{\leq 1}\pmb\psi](\tau_1, \tau_2)+\sqrt{\widehat{\EMF}[\pmb\psi](\tau_1, \tau_2)}\sqrt{\widehat{\EMF}[\pr^{\leq 1}\pmb\psi](\tau_1, \tau_2)}.
\eea
Then, relying on \eqref{eq:recoveringderivatives:firstorder:Mora:middlestep:analogforwidehatcase} and \eqref{eq:eqsforconditionalhighorderestimatelemma:general:allprreg:withprttprtphi}, and  proceeding by induction  as in Step 1, we obtain the following analog of \eqref{eq:controlofMorawetzrplus1plusdehto11mforregderivativesofpmbpsiconditionalonprtauandchi0prphit}
\beaa
&&\widehat{\M}_{r_+(1+\dred), 11m}[\pr^{\leq \reg}\pmb\psi](\tau_1, \tau_2)\\
&\les& \int_{\MM_{r\leq 12m}(\tau_1, \tau_2)}|\pr^{\leq \reg-1}\F|^2 
+\widehat{\M}[(\pr_\tau, \chi_0\pr_{\tphi})^{\leq \reg}\pmb\psi](\tau_1, \tau_2)+\M_{r\geq 11m}[\pr^{\leq \reg}\pmb\psi](\tau_1, \tau_2)\nn\\
\nn&&+\ep\M_{r\leq r_+(1+\dred)}[\pr^{\leq \reg}\pmb\psi](\tau_1, \tau_2)+\sqrt{\widehat{\EMF}[\pr^{\leq \reg-1}\pmb\psi](\tau_1, \tau_2)}\sqrt{\widehat{\EMF}[\pr^{\leq \reg}\pmb\psi](\tau_1, \tau_2)}.
\eeaa
The rest of the proof of \eqref{eq:higherorderenergymorawetzconditionalonsprtausprphis:nontrapped} is then identical to the one of \eqref{eq:higherorderenergymorawetzconditionalonsprtausprphis}.

\noindent{\bf Step 8.} \textit{Proof of  \eqref{eq:improvedhigherorderenergymorawetzforlargerconditionalonsprtau}, \eqref{eq:imporvedhigherorderenergymorawetzconditionalonlowerrweigth} and \eqref{eq:imporvedhigherorderenergymorawetzconditionalonlowerrweigth:nontrapped}}.  In view of \cite[Equation (3.41)]{MaSz24}, we have the following equation
\beaa
\pr^2_r\psi_{ij} = N_{ij}+O(1)\nab\pr\psi+O(1)\pr_r\pr_\tau\psi_{ij}
+O(\ep)\pr_r^2\psi_{ij}+O(r^{-1})\pr\pr^{\leq 1}\psi_{ij}, \quad \text{for }\,\, r\geq 10m.
\eeaa
We take the square of the modulus of both sides of the above identity, multiply both squares by $r^{-1-\de}$, sum over $i,j=1,2,3$, and integrate over $\MM_{r\geq 11m}(\tau_1, \tau_2)$ which yields
\beaa
\int_{\MM_{r\geq 11m}(\tau_1, \tau_2)}\frac{|\pr_r^2\pmb\psi|^2}{r^{1+\de}} 
&\les &\int_{\MM_{r\geq 11m}(\tau_1, \tau_2)}\frac{|\N|^2}{r^{1+\de}}+\M_{r\geq 11m}[\pr^{\leq 1}\pmb\psi](\tau_1, \tau_2)\\
&&+\M_{\de, r\geq 11m}[\pr_\tau^{\leq 1}\pmb\psi](\tau_1, \tau_2)+\ep^2\M_{\de, r\geq 11m}[\pr^{\leq 1}\pmb\psi](\tau_1, \tau_2)\nn\\
&\les &\int_{\MM_{r\geq 11m}(\tau_1, \tau_2)}\frac{|\F|^2}{r^{1+\de}}+\M_{r\geq 11m}[\pr^{\leq 1}\pmb\psi](\tau_1, \tau_2)\\
&&+\M_{\de, r\geq 11m}[\pr_\tau^{\leq 1}\pmb\psi](\tau_1, \tau_2)+\ep^2\M_{\de, r\geq 11m}[\pr^{\leq 1}\pmb\psi](\tau_1, \tau_2),
\eeaa
where in the last step we have used the formula  of $N_{ij}$ in \eqref{eq:eqsforconditionalhighorderestimatelemma:general}.
Since 
\beaa
&&\M_{\de, r\geq 11m}[\pr^{\leq 1}\pmb\psi](\tau_1, \tau_2)\nn\\
& \les &\M_{r\geq 11m}[\pr^{\leq 1}\pmb\psi](\tau_1, \tau_2)+\int_{\MM_{r\geq 11m}(\tau_1, \tau_2)}\frac{|\pr_r^2\pmb\psi|^2}{r^{1+\de}}+\M_{\de, r\geq 11m}[\pr_\tau^{\leq 1}\pmb\psi](\tau_1, \tau_2),
\eeaa
we infer, for $\ep$ small enough, 
\bea\lab{eq:highorderMora:improveMora:byTpsi}
&&\M_{\de, r\geq 11m}[\pr^{\leq 1}\pmb\psi](\tau_1, \tau_2)\nn\\
& \les& \int_{\MM_{r\geq 11m}(\tau_1, \tau_2)}\frac{|\F|^2}{r^{1+\de}}+\M_{r\geq 11m}[\pr^{\leq 1}\pmb\psi](\tau_1, \tau_2)+\M_{\de, r\geq 11m}[\pr_\tau^{\leq 1}\pmb\psi](\tau_1, \tau_2).
\eea
Then, relying on \eqref{eq:highorderMora:improveMora:byTpsi} and \eqref{eq:eqsforconditionalhighorderestimatelemma:general:allprreg}, and  proceeding by induction as in Step 1, we obtain
\beaa
\M_{\de, r\geq 11m}[\pr^{\leq \reg}\pmb\psi](\tau_1, \tau_2) 
&\les& \int_{\MM_{r\geq 11m}(\tau_1, \tau_2)}\frac{|\pr^{\leq \reg-1}\F|^2}{r^{1+\de}}+\M_{r\geq 11m}[\pr^{\leq \reg}\pmb\psi](\tau_1, \tau_2)\nn\\
&&
+\M_{\de, r\geq 11m}[\pr_\tau^{\leq \reg}\pmb\psi](\tau_1, \tau_2)+\ep {\M_{\de, r\geq 11m}[\pr^{\leq \reg}\pmb\psi](\tau_1,\tau_2)},
\eeaa
and hence for $\ep$ small enough 
\bea\lab{eq:highorderMora:improveMora:byTpsi:highorder}
\M_{\de, r\geq 11m}[\pr^{\leq \reg}\pmb\psi](\tau_1, \tau_2)
&\les& \int_{\MM_{r\geq 11m}(\tau_1, \tau_2)}|\pr^{\leq \reg-1}\F|^2+\M_{r\geq 11m}[\pr^{\leq \reg}\pmb\psi](\tau_1, \tau_2)\nn\\
&&+\M_{\de, r\geq 11m}[\pr_\tau^{\leq \reg}\pmb\psi](\tau_1, \tau_2),
\eea
which together with the estimate \eqref{eq:higherorderenergymorawetzforlargerconditionalonsprtau} yields the desired estimate \eqref{eq:improvedhigherorderenergymorawetzforlargerconditionalonsprtau}, together with the  estimate \eqref{eq:higherorderenergymorawetzconditionalonsprtausprphis} yields the desired estimate \eqref{eq:imporvedhigherorderenergymorawetzconditionalonlowerrweigth}, and  together with the estimate \eqref{eq:higherorderenergymorawetzconditionalonsprtausprphis:nontrapped} yields the desired estimate \eqref{eq:imporvedhigherorderenergymorawetzconditionalonlowerrweigth:nontrapped}.
This concludes the proof of Lemma \ref{lemma:higherorderenergyMorawetzestimates}.
\end{proof}

Next, we introduce the following EMF norms for $s=\pm 2$ and any $\de\in(0,\frac{1}{3}]$
\begin{align}
\lab{def:EMFtotalpsnotilde:pm2}
\EMFtotalps{s}:={}&\sum_{p=0,1,2}{\EMF}[\psis{p}](\Iti)
+\sum_{p=0,1}\widehat{\M}_{\de}[\psis{p}](\tau_1+1, \tau_2-3)
\nn\\
&+\sum_{p=0,1}\widehat{\EMF}_{\de}[\phis{p}](\tau_1,\tau_2)+ {\EMF}_{\de}[\phis{2}](\tau_1,\tau_2),
\end{align}
and the following early time EMF norms, for $s=\pm 2$ and any $\de\in[0,1]$,
\bea\lab{eq:definitionofinitialenergyofphiandpsi:IEtotalterm:delta}
\IEde{\pr^{\leq \reg}\pmb\phi_s, \pr^{\leq \reg}\pmb\psi_s} :=\IE{\pr^{\leq \reg}\pmb\phi_s, \pr^{\leq \reg}\pmb\psi_s} + \sum_{p=0}^2{\MF}_{\de}[\pr^{\leq \reg}\psis{p}](\tmic,\tau_1+1),
\eea
where $\IE{\pr^{\leq \reg}\pmb\phi_s, \pr^{\leq \reg}\pmb\psi_s}$ is given by \eqref{eq:definitionofinitialenergyofphiandpsi:IEterm}. Then, we apply Lemma \ref{lemma:higherorderenergyMorawetzestimates} to the coupled Teukolsky wave system and deduce the following high order EMF estimates which recover the control of all unweighted derivatives from the control of high order $(\pr_{\tt}, \chi_0\pr_{\tphi})$ derivatives.

\begin{lemma}\lab{lemma:higherorderenergyMorawetzestimates:Teu}
Under the same assumptions as in Theorem \ref{th:main:intermediary}, we have, for  $s=\pm 2$, $1\leq \reg\leq 14$ and $\de\in(0,\frac{1}{3}]$, 
\bea\lab{eq:higherorderenergymorawetzforlargerconditionalonsprtau:Teu}
&&\EMF_{s, \de, \text{total}, r\geq 11m}[\pr^{\leq \reg}\pmb\phi_{s}] \nn\\
&\les&\EMFtotalhps{s}{\nab_{\pr_{\tt}}^{\leq \reg}}+\IEde{\pr^{\leq \reg}\pmb\phi_s, \pr^{\leq \reg}\pmb\psi_s} 
+\NNttotalph{s}{\pr^{\leq \reg}}\nn\\
&&+\sum_{p=0}^{2}\widehat{\mathcal{N}}_{r\geq 10m}''[\pr^{\leq \reg}\psis{p}, \pr^{\leq \reg}\widehat{\F}_{total,s}^{(p)}](\tau_2-3, \tau_2)\nn\\
&&
+\Big(\EMFtotalhps{s}{\pr^{\leq \reg-1}}\Big)^{\frac{1}{2}}\Big(\EMFtotalhps{s}{\pr^{\leq \reg}}\Big)^{\frac{1}{2}},
\eea
with $\EMF_{s,\de, \text{total}}$, $\IEde{\pr^{\leq \reg}\pmb\phi_s, \pr^{\leq \reg}\pmb\psi_s}$, $\widetilde{\NN}_{s, \de, \text{total}}$ and $\widehat{\NN}''[\c,\c](\c,\c)$ given as in \eqref{def:EMFtotalpsnotilde:pm2}, \eqref{eq:definitionofinitialenergyofphiandpsi:IEtotalterm:delta},  \eqref{expression:NNttotalps:pm2} and \eqref{eq:defofNNhat''}, respectively.

Also, under the same assumptions  as in Theorem \ref{th:main:intermediary}, we have,  for  $s=\pm 2$, $1\leq \reg\leq 14$ and $\de\in(0,\frac{1}{3}]$, 
\bea
\lab{eq:imporvedhigherorderenergymorawetzconditionalonlowerrweigth:Teu}
&&\EMFtotalhps{s}{\pr^{\leq \reg}}\nn\\
&\les&\EMFtotalhps{s}{\nab_{(\pr_{\tt},\chi_{0}\pr_{\tilde{\phi}})}^{\leq \reg}}+\IEde{\pr^{\leq \reg}\pmb\phi_s, \pr^{\leq \reg}\pmb\psi_s} 
+\NNttotalph{s}{\pr^{\leq \reg}}\nn\\
&&+\sum_{p=0}^{2}\widehat{\mathcal{N}}''[\pr^{\leq \reg}\psis{p}, \pr^{\leq \reg}\widehat{\F}_{total,s}^{(p)}](\tau_2-3, \tau_2),
\eea
with the cut-off function $\chi_0$ given as in \eqref{def:cutoffcuntionchi0}.
\end{lemma}

\begin{proof}
Recall from \eqref{eq:waveeqwidetildepsi1:gtilde:Teu} the system of coupled wave equations for $\psiss{ij}{p}$
\bea
\lab{eq:waveeqwidetildepsi1:gtilde:Teu:paste}
{\square}_{{\g}_{\chi_{\tau_1, \tau_2}}} \psi^{(p)}_{s,ij}
&=&\chi_{\tau_1, \tau_2}\big( \widehat{S}(\psi^{(p)}_s)_{ij} +(\widehat{Q}\psi^{(p)}_s)_{ij}\big) +(1-\chi_{\tau_1, \tau_2})\big(\widehat{S}_K(\psi^{(p)}_s)_{ij} +(\widehat{Q}_K\psi^{(p)}_s)_{ij} \big)\nn\\
&&+\bigg((1-\chi_{\tau_1, \tau_2})f_p+\frac{4-2\de_{p0}}{|q|^{2}}\bigg)\psi^{(p)}_{s,ij}+\chi_{\tau_1, \tau_2}^{(1)}L^{(p)}_{s,ij}+\widehat{F}_{total,s,ij}^{(p)},
\eea
and recall from \eqref{eq:ScalarizedTeuSys:general:Kerrperturbation}-\eqref{eq:defofwidehatsquaregoperator} that
\bea
\lab{eq:waveeqphisijp:gtilde:Teu:paste}
{\square}_{{\g}}\phiss{ij}{p} 
=(4-2\de_{p0})|q|^{-2}\phiss{ij}{p}+\big( \widehat{S}(\pmb\phi^{(p)}_s)_{ij} +(\widehat{Q}\pmb\phi^{(p)}_s)_{ij}\big)+L^{(p)}_{s,ij}+N_{W,s,ij}^{(p)}.
\eea
Now, the above system of coupled wave equations \eqref{eq:waveeqwidetildepsi1:gtilde:Teu:paste}-\eqref{eq:waveeqphisijp:gtilde:Teu:paste} for the family of scalars $(\psi)_{ij}=(\psiss{ij}{p}, \phiss{ij}{p})$ is in the form \eqref{eq:eqsforconditionalhighorderestimatelemma:general} with both of the two metrics $\g_{\chi_{\tau_1,\tau_2}}$ and $\g$ satisfying the assumptions of Section \ref{subsubsect:assumps:perturbedmetric}, with $D_1=0$ and with $(F)_{ij}=(\widehat{F}_{total,s,ij}^{(p)} , N_{W,s,ij}^{(p)})$. Also, in view of \eqref{eq:expansionofSandQterms:neareventhorizon} and 
\eqref{eq:expansionofealpha:neareventhorizon}, this system can be put into the form \eqref{eq:Redshift-gen.scalarwaveeqs} with both metrics $\g_{\chi_{\tau_1,\tau_2}}$ and $\g$ satisfying the assumptions of Section \ref{subsubsect:assumps:perturbedmetric}, with $D_1=0$ and with $(F)_{ij}=(\widehat{F}_{total,s,ij}^{(p)} , N_{W,s,ij}^{(p)})$, so that applying  Lemma \ref{lemma:redshiftestimatesscalarwaveeqs:general} to this system yields the redshift estimates \eqref{eq:assum:redshift:highorderunweightedEMF:generalwave}. Thus, we can apply Lemma \ref{lemma:higherorderenergyMorawetzestimates}.

First, applying the estimate \eqref{eq:imporvedhigherorderenergymorawetzconditionalonlowerrweigth} to the system of coupled wave equations \eqref{eq:waveeqphisijp:gtilde:Teu:paste} for $p=0,1,2$ on $\MM(\tau_1, \tau_2)$, and the estimate \eqref{eq:imporvedhigherorderenergymorawetzconditionalonlowerrweigth:nontrapped} to the system of coupled wave equations \eqref{eq:waveeqphisijp:gtilde:Teu:paste} for $p=0,1$ on $\MM(\tau_1, \tau_2)$, and since $\psis{p}=\phis{p}$ on $\MM(\tau_1+1,\tau_2-3)$ in view of \eqref{eq:psispequalsphispintau1+1tau2-1}, we deduce 
\beaa
&&\sum_{p=0}^2{\EMF}_{\de}[\pr^{\leq \reg}\psis{p}](\tau_1+1,\tau_2-3)+\sum_{p=0}^1\widehat{\M}_{\de}[\pr^{\leq \reg}\psis{p}](\tau_1+1, \tau_2-3)\nn\\
&&+\sum_{p=0}^{1}\widehat{\EMF}_{\de}[\pr^{\leq \reg}\phis{p}](\tau_1,\tau_2)
+ {\EMF}_{\de}[\pr^{\leq \reg}\phis{2}](\tau_1,\tau_2)\nn\\
 &\les& \sum_{p=0}^{2}\E[\pr^{\leq \reg}\phis{p}](\tau_1) 
+\sum_{p=0}^1\widehat{\EMF}_{\de}[\nab_{(\pr_{\tt},\chi_{0}\pr_{\tilde{\phi}})}^{\leq \reg}\phis{p}](\tau_1, \tau_2)
+{\EMF}_{\de}[\nab_{(\pr_{\tt},\chi_{0}\pr_{\tilde{\phi}})}^{\leq \reg}\phis{2}](\tau_1, \tau_2)\nn\\
&& +\sum_{p=0}^{2}\sum_{i,j}\widehat{\mathcal{N}}_{r\geq 10m}''[\pr^{\leq \reg}\phiss{ij}{p}, \pr^{\leq \reg}N_{W,s,ij}^{(p)}](\tau_1, \tau_2) +\sum_{p=0}^{2}\sum_{i,j}\int_{\MM(\tau_1, \tau_2)}|\pr^{\leq\reg}N_{W,s,ij}|^2,
\eeaa
and hence
\bea
\lab{eq:imporvedhigherorderEMFfromprttprphi:tau1tau2:proof:incomplete}
&&\sum_{p=0}^2{\EMF}_{\de}[\pr^{\leq \reg}\psis{p}](\tmic,\tau_2-3)+\sum_{p=0}^1\widehat{\M}_{\de}[\pr^{\leq \reg}\psis{p}](\tau_1+1, \tau_2-3)\nn\\
&&+\sum_{p=0}^{1}\widehat{\EMF}_{\de}[\pr^{\leq \reg}\phis{p}](\tau_1,\tau_2)
+ {\EMF}_{\de}[\pr^{\leq \reg}\phis{2}](\tau_1,\tau_2)\nn\\
&\les& \IEde{\pr^{\leq \reg}\pmb\phi_s, \pr^{\leq \reg}\pmb\psi_s}  +\EMFtotalhps{s}{\nab_{(\pr_{\tt},\chi_{0}\pr_{\tilde{\phi}})}^{\leq \reg}}+\NNttotalph{s}{\pr^{\leq \reg}},
\eea
where we have used, in view of \eqref{eq:definitionofinitialenergyofphiandpsi:IEterm} and \eqref{eq:definitionofinitialenergyofphiandpsi:IEtotalterm:delta},
\beaa
 \sum_{p=0}^{2}\E[\pr^{\leq \reg}\phis{p}](\tau_1) +\sum_{p=0}^2{\EMF}_{\de}[\pr^{\leq \reg}\psis{p}](\tmic,\tau_1+1)\leq \IEde{\pr^{\leq \reg}\pmb\phi_s, \pr^{\leq \reg}\pmb\psi_s},
\eeaa
as well as the following consequence of \eqref{eq:defofNNhat'}, \eqref{expression:NNttotalps:pm2} and \eqref{eq:defofNNhat''}
\bea
\lab{eq:widehatNNrgeq10m:phisp:highorderprderi:bdbyNNtstotal}
\nn&&\sum_{p=0}^{2}\sum_{i,j}\widehat{\mathcal{N}}_{r\geq 10m}''[\pr^{\leq \reg}\phiss{ij}{p}, \pr^{\leq \reg}N_{W,s,ij}^{(p)}](\tau_1, \tau_2) +\sum_{p=0}^{2}\sum_{i,j}\int_{\MM(\tau_1, \tau_2)}|\pr^{\leq\reg}N_{W,s,ij}|^2\\
&\les&\sum_{p=0}^{2}\sum_{i,j}\widehat{\mathcal{N}}'[\pr^{\leq \reg}\phiss{ij}{p}, \pr^{\leq \reg}N_{W,s,ij}^{(p)}](\tau_1, \tau_2)\nn\\
&\les& \NNttotalph{s}{\pr^{\leq \reg}}.
\eea

Applying the estimate \eqref{eq:higherorderenergymorawetzconditionalonsprtausprphis} to the system of coupled wave equations \eqref{eq:waveeqwidetildepsi1:gtilde:Teu:paste}-\eqref{eq:waveeqphisijp:gtilde:Teu:paste} for the family of scalars $(\psi)_{ij}=(\psiss{ij}{p}, \phiss{ij}{p})$, $p=0,1,2$, on $\MM(\tau_2-3,\tau_2)$, and using  \eqref{eq:psispequalsphispintau1+1tau2-1}, we deduce
\beaa
&&\sum_{p=0}^2{\EMF}[\pr^{\leq \reg}\psis{p}](\tau_2-3,\tau_2)+\sum_{p=0}^{2}{\EMF}[\pr^{\leq \reg}\phis{p}](\tau_2-3,\tau_2)
\nn\\
 &\les& \sum_{p=0}^{2}\Big(\E[\pr^{\leq \reg}\phis{p}](\tau_2-3)
 +
 \EMF[\nab_{(\pr_\tau, \chi_0\pr_{\tphi})}^{\leq \reg}\psis{p}](\tau_2-3, \tau_2)\\
 &&+ \EMF[\nab_{(\pr_\tau, \chi_0\pr_{\tphi})}^{\leq \reg}\phis{p}](\tau_2-3, \tau_2)\Big)
+\sum_{p=0}^{2}\widehat{\mathcal{N}}_{r\geq 10m}''[\pr^{\leq \reg}\phis{p}, \pr^{\leq \reg}\N_{W,s}^{(p)}](\tau_2-3, \tau_2)\\
&&+\sum_{p=0}^{2}\widehat{\mathcal{N}}_{r\geq 10m}''[\pr^{\leq \reg}\psis{p}, \pr^{\leq \reg}\widehat{\F}_{total,s}^{(p)}](\tau_2-3, \tau_2)\\
&&+\int_{\MM(\tau_2-3, \tau_2)}|\pr^{\leq\reg}\N_{W,s}^{(p)}|^2+\int_{\MM(\tau_2-3, \tau_2)}|\pr^{\leq\reg}\widehat{\F}_{total,s}^{(p)}|^2,
\eeaa
which together with \eqref{eq:imporvedhigherorderEMFfromprttprphi:tau1tau2:proof:incomplete} and \eqref{eq:widehatNNrgeq10m:phisp:highorderprderi:bdbyNNtstotal} implies
\bea
\lab{eq:imporvedhigherorderEMFfromprttprphi:tau1tau2:proof}
&&\sum_{p=0}^2{\EMF}_{\de}[\pr^{\leq \reg}\psis{p}](\tmic,\tau_2-3)+\sum_{p=0}^2{\EMF}[\pr^{\leq \reg}\psis{p}](\tau_2-3,\tau_2)\nn\\
&&+\sum_{p=0}^1\widehat{\M}_{\de}[\pr^{\leq \reg}\psis{p}](\tau_1+1, \tau_2-3)+\sum_{p=0}^{1}\widehat{\EMF}_{\de}[\pr^{\leq \reg}\phis{p}](\tau_1,\tau_2)
+ {\EMF}_{\de}[\pr^{\leq \reg}\phis{2}](\tau_1,\tau_2)\nn\\
&\les& \IEde{\pr^{\leq \reg}\pmb\phi_s, \pr^{\leq \reg}\pmb\psi_s}  +\EMFtotalhps{s}{\nab_{(\pr_{\tt},\chi_{0}\pr_{\tilde{\phi}})}^{\leq \reg}}+\NNttotalph{s}{\pr^{\leq \reg}}\nn\\
&&+\sum_{p=0}^{2}\widehat{\mathcal{N}}''[\pr^{\leq \reg}\psis{p}, \pr^{\leq \reg}\widehat{\F}_{total,s}^{(p)}](\tau_2-3, \tau_2).
\eea

Next, since $\psiss{ij}{p}$ solves 
\bea
\lab{eq:waveeqwidetildepsi1:gtilde:Teu:paste:tau2toinfty}
{\square}_{\gam} \psi^{(p)}_{s,ij}
&=&\widehat{S}_K(\psi^{(p)}_s)_{ij} +(\widehat{Q}_K\psi^{(p)}_s)_{ij} +\bigg(\frac{4}{|q|^2}- \frac{4a^2\cos^2\th(|q|^2+6mr)}{|q|^6}\bigg)\psi^{(p)}_{s,ij}\nn\\
&=&\sum_{k,l}O(r^{-2})\dk^{\leq 1}\psiss{kl}{p}, \qquad \text{on} \quad \MM(\tau_2,+\infty),
\eea
an application of the estimate \eqref{eq:higherorderenergymorawetzconditionalonsprtausprphis} yields
\bea
\lab{eq:EMFde:improvedhighorderfromprttandchiprphi:tau2toinfty:proof}
\sum_{p=0}^2{\EMF}[\pr^{\leq \reg}\psis{p}](\tau_2,+\infty) 
\les \sum_{p=0}^2\Big(\E[\pr^{\leq \reg}\psis{p}](\tau_2)+\EMF[\nab_{(\pr_{\tt},\chi_{0}\pr_{\tilde{\phi}})}^{\leq \reg}\psis{p}](\tau_2,+\infty)\Big).
\eea

Now, we add  the estimates \eqref{eq:imporvedhigherorderEMFfromprttprphi:tau1tau2:proof} and  \eqref{eq:EMFde:improvedhighorderfromprttandchiprphi:tau2toinfty:proof}  together to deduce 
\beaa
\EMFtotalhps{s}{\pr^{\leq \reg}}&\les&\EMFtotalhps{s}{\nab_{(\pr_{\tt},\chi_{0}\pr_{\tilde{\phi}})}^{\leq \reg}}+\IEde{\pr^{\leq \reg}\pmb\phi_s, \pr^{\leq \reg}\pmb\psi_s} \nn\\
&&+\NNttotalph{s}{\pr^{\leq \reg}}+\sum_{p=0}^{2}\widehat{\mathcal{N}}''[\pr^{\leq \reg}\psis{p}, \pr^{\leq \reg}\widehat{\F}_{total,s}^{(p)}](\tau_2-3, \tau_2),
\eeaa
which is the desired estimate \eqref{eq:imporvedhigherorderenergymorawetzconditionalonlowerrweigth:Teu}.

The proof of the estimate \eqref{eq:higherorderenergymorawetzforlargerconditionalonsprtau:Teu} follows in a similar manner, the only difference being that we use the estimate \eqref{eq:improvedhigherorderenergymorawetzforlargerconditionalonsprtau} instead of the estimates \eqref{eq:imporvedhigherorderenergymorawetzconditionalonlowerrweigth}--\eqref{eq:imporvedhigherorderenergymorawetzconditionalonlowerrweigth:nontrapped}, and the estimate \eqref{eq:higherorderenergymorawetzforlargerconditionalonsprtau} instead of the estimate \eqref{eq:higherorderenergymorawetzconditionalonsprtausprphis}. This concludes the proof of Lemma \ref{lemma:higherorderenergyMorawetzestimates:Teu}.
\end{proof}


\subsection{Estimates for unweighted derivatives of $(\phis{p}, \psis{p})$}
\lab{subsect:highorderunweightedEMF:conditional}


The goal of this section is to prove the following proposition on the control of EMF estimates for high order unweighted derivatives of $(\phis{p}, \psis{p})$, which is the analog of Theorem \ref{thm:EMF:systemofTeuscalarized:order0:final} in the high order regularity case.

\begin{proposition}[Conditional EMF estimates for unweighted derivatives of $(\phis{p}, \psis{p})$]
\lab{prop:EMFtotalp:finalestimate:pm2:prhighorder}
Under the assumptions of Theorem \ref{th:main:intermediary}, we have, for $s=\pm2$, all $\reg\leq 14$ and $\de\in (0,\frac{1}{3}]$,
\bea
\lab{eq:EMFtotalp:finalestimate:pm2:prhighorder} 
&&\EMFtotalh{s}{\nab_{(\pr_{\tt},\chi_{0}\pr_{\tilde{\phi}})}^{\leq \reg}}+\EMFtotalhps{s}{\pr^{\leq \reg}} \nn\\
&\les& \IEde{\pr^{\leq  \reg}\pmb\phi_s, \pr^{\leq  \reg}\pmb\psi_s} +\sum_{p=0}^2\Errdefect[\nab_{(\pr_{\tt},\chi_{0}\prtphihat)}^{\leq \reg}\psis{p}]+\NNttotalph{s}{\pr^{\leq \reg}}
  \nn\\
&&+\A[\pmb\psi_{s}](\Iti)+\A[\pmb\phi_{s}](\tau_1,\tau_2)
+\sum_{p=0}^{2}\widehat{\mathcal{N}}''[\pr^{\leq \reg}\psis{p}, \pr^{\leq \reg}\widehat{\F}_{total,s}^{(p)}](\tau_2-3, \tau_2),
\eea
with $\widetilde{\EMF}_{s,\de, \text{total}}$, $\EMF_{s,\de,\text{total}}$, $\IEde{\pr^{\leq  \reg}\pmb\phi_s, \pr^{\leq  \reg}\pmb\psi_s}$, $\Errdefect[\c]$, $\widetilde{\NN}_{s, \de, \text{total}}$, $\A$ and $\widehat{\mathcal{N}}''$ given in \eqref{def:EMFtotalps:pm2}, \eqref{def:EMFtotalpsnotilde:pm2}, \eqref{eq:definitionofinitialenergyofphiandpsi:IEtotalterm:delta}, \eqref{def:Errdefectofpsi}, \eqref{expression:NNttotalps:pm2}, \eqref{def:AandAonorms:tau1tau2:phisandpsis} and \eqref{eq:defofNNhat''}, respectively, and with the cut-off function $\chi_0$ introduced in \eqref{def:cutoffcuntionchi0}.
\end{proposition}

We first prove EMF estimates for $\pr_{\tau}$ derivatives in Section \ref{subsubsect:commutationwithprtt:conditionalEMF}, then prove EMF estimates for $\chi_0\prtphihat$ derivatives in Section \ref{subsubsect:commutationwithprtphihat:conditionalEMF}, and in the end prove Proposition \ref{prop:EMFtotalp:finalestimate:pm2:prhighorder}  in Section \ref{subsubsect:commutationwithprttprtphi:conditionalEMF:high:v1}.


\subsubsection{Commutation with $\pr_{\tt}$}
\lab{subsubsect:commutationwithprtt:conditionalEMF}


Commuting $\pr_{\tt}$ with the scalarized Teukolsky transport equations \eqref{eq:ScalarizedQuantitiesinTeuSystem:Kerrperturbation}, we deduce
\bsub
\lab{eq:ScalarizedQuantitiesinTeuSystem:Kerrperturbation:pr_tt}
\bea
&& e_3\bigg(\frac{r\bar{q}}{q}\bigg({\frac{r^2}{|q|^2}}\bigg)^{p-2}\phiplussh{ij}{p}{\pr_{\tt}} \bigg) - \frac{r\bar{q}}{q}\bigg({\frac{r^2}{|q|^2}}\bigg)^{p-2}\Big(M_{i3}^k \phiplussh{kj}{p}{\pr_{\tt}} + M_{j3}^k \phiplussh{ik}{p}{\pr_{\tt}}\Big)\nn\\
&=&\frac{\bar{q}}{rq}\bigg({\frac{r^2}{|q|^2}}\bigg)^{p-1}\phiplussh{ij}{p+1}{\pr_{\tt}}+N_{T,+2,ij}^{(p),\pr_{\tt}}
\eea
and 
\bea
&& e_4\bigg(\frac{rq}{\bar{q}}\bigg(\frac{r^2}{|q|^2}\bigg)^{p-2}\phiminussh{ij}{p}{\pr_{\tt}}\bigg) - \frac{rq}{\bar{q}}\bigg(\frac{r^2}{|q|^2}\bigg)^{p-2}\Big(M_{i4}^k \phiminussh{kj}{p}{\pr_{\tt}} + M_{j4}^k \phiminussh{ik}{p}{\pr_{\tt}}\Big)\nn\\
&=& \frac{q}{r\bar{q}}\left(\frac{r^2}{|q|^2}\right)^{p-1}\frac{\De}{\qs}\phiminussh{ij}{p+1}{\pr_{\tt}} +N_{T,-2,ij}^{(p),\pr_{\tt}},
\eea
\esub
where the scalars $\phissh{ij}{p}{\pr_{\tt}}$, $s=\pm 2$, $i,j=1,2,3$,  $p=0,1,2$, are defined by
\beaa
\phissh{ij}{p}{\pr_{\tt}}:=\pr_{\tt}\phiss{ij}{p}
\eeaa
and the complex-valued scalars $N_{T,s,ij}^{(p),\pr_{\tt}}$,  $s=\pm 2$, $i,j=1,2,3$,  $p=0,1,2$,  are given by 
\bsub
\lab{def:Nt+-2ijp:pr_tt}
\bea
\lab{def:Nt+2ijp:pr_tt}
N_{T,+2,ij}^{(p), \pr_{\tt}}&=&\pr_{\tt}N_{T,+2,ij}^{(p)}+[e_3, \pr_{\tt}]\bigg(\frac{r\bar{q}}{q}\bigg({\frac{r^2}{|q|^2}}\bigg)^{p-2}\phipluss{ij}{p} \bigg)\nn\\
&&+\frac{r\bar{q}}{q}\bigg({\frac{r^2}{|q|^2}}\bigg)^{p-2}\Big(\pr_{\tt}(M_{i3}^k) \phipluss{kj}{p} + \pr_{\tt}(M_{j3}^k) \phipluss{ik}{p}\Big)\nn\\
&=&\pr_{\tt}N_{T,+2,ij}^{(p)}+\sum_{k,l}\Big(r\dk^{\leq 1}\Ga_b \dk^{\leq 1}\phipluss{kl}{p} +r^2\dk^{\leq 1}\Ga_g\pr_{\tt}\phipluss{kl}{p}\Big)
\eea
and
\bea
\lab{def:Nt-2ijp:pr_tt}
N_{T,-2,ij}^{(p), \pr_{\tt}}&=&\pr_{\tt}N_{T,-2,ij}^{(p)}+[e_4, \pr_{\tt}]\bigg(\frac{rq}{\ov{q}}\bigg({\frac{r^2}{|q|^2}}\bigg)^{p-2}\phiminuss{ij}{p} \bigg)\nn\\
&&+\frac{rq}{\ov{q}}\bigg({\frac{r^2}{|q|^2}}\bigg)^{p-2}\Big(\pr_{\tt}(M_{i4}^k) \phiminuss{kj}{p} + \pr_{\tt}(M_{j4}^k) \phiminuss{ik}{p}\Big)\nn\\
&=&\pr_{\tt}N_{T,-2,ij}^{(p)}+\sum_{k,l}r\dk^{\leq 1}\Ga_g\dk^{\leq 1} \phiminuss{kl}{p}
\eea
\esub
in view of the assumption \eqref{eq:assumptionsonregulartripletinperturbationsofKerr:0} and 
\beaa
&&[e_3, \pr_{\tt}]=[\widecheck{e_{3}(x^{\b})} \pr_{\b},\pr_{\tt}]=-\pr_{\tt}(\widecheck{e_{3}(x^{\b})})\pr_{\b} =\dk^{\leq 1}\Ga_b\dk+r\dk^{\leq 1}\Ga_g \pr_{\tt},\\
&&[e_4, \pr_{\tt}]=[\widecheck{e_{4}(x^{\b})} \pr_{\b},\pr_{\tt}]=-\pr_{\tt}(\widecheck{e_{4}(x^{\b})})\pr_{\b}=\dk^{\leq 1}\Ga_g\dk
\eeaa
which follow from Definition \ref{definition.Ga_gGa_b}.

Next, commuting $\pr_{\tt}$ with the scalarized Teukolsky wave equations \eqref{eq:ScalarizedTeuSys:general:Kerrperturbation}, and using \eqref{eq:defofwidehatsquaregoperator} \eqref{eq:definitionwidehatSandwidehatQperturbationsofKerr}   \eqref{eq:linearterms:ScalarizedTeuSys:general:Kerrperturbation} to expand out the wave operator and the linear coupling terms, we obtain scalarized wave equations for $\phissh{ij}{p}{\pr_{\tt}}$ in the following schematic form
\bea
\lab{eq:ScalarizedTeuSys:general:Kerrperturbation:pr_tt} 
\square_\g(\phissh{ij}{p}{\pr_{\tt}}) -\frac{4-2\de_{p0}}{\qs} \phissh{ij}{p} {\pr_{\tt}}=\widehat{S}(\phi_s^{(p),\pr_{\tau}})_{ij} +(\widehat{Q}\phi_s^{(p),\pr_{\tau}})_{ij} + {L_{s,ij}^{(p),\pr_{\tt}}}+N_{W,s,ij}^{(p),\pr_{\tt}}
\eea
for $s=\pm 2$, $i,j=1,2,3$,  $p=0,1,2$,
where $\widehat{S}(\phi_s^{(p),\pr_{\tau}})_{ij}$, $(\widehat{Q}\phi_s^{(p),\pr_{\tau}})_{ij}$ and ${L_{s,ij}^{(p),\pr_{\tt}}}$ are obtained by replacing each $\phiss{kl}{p}$ with $\phissh{kl}{p} {\pr_{\tt}}$ in the expanded forms of  $\widehat{S}(\phi_s^{(p)})_{ij}$, $(\widehat{Q}\phi_s^{(p)})_{ij}$ and ${L_{s,ij}^{(p)}}$ in \eqref{eq:definitionwidehatSandwidehatQperturbationsofKerr}-\eqref{eq:linearterms:ScalarizedTeuSys:general:Kerrperturbation}, respectively, and where 
\bea
\lab{def:Nwsijp:pr_tt}
N_{W,s,ij}^{(p),\pr_{\tt}}&=&\pr_{\tt}N_{W,s,ij}^{(p)} + [\square_{\g}, \pr_{\tt}]\phiss{ij}{p}-\frac{4ia\cos\th}{|q|^2}\big(\pr_{\tt}M_{i\tau}^l \phiss{lj}{p}+\pr_{\tt}M_{j\tau}^l \phiss{il}{p}\big)\nn\\
&&+2\pr_{\tt}M_{i}^{k\a}\pr_\a(\psi_{kj}) +2\pr_{\tt}M_{j}^{k\a}\pr_\a(\psi_{ik})
+\pr_{\tt}(\Ddot^\a M_{i\a}^k)\psi_{kj}
+\pr_{\tt}(\Ddot^\a M_{j\a}^k)\psi_{ik} \nn\\
&&-\pr_{\tt}(M_{i\a}^kM_k^{l\a})\psi_{lj}-2\pr_{\tt}(M_{i\a}^kM_{j}^{l\a})\psi_{kl}-\pr_{\tt}(M_{j\a}^kM_k^{l\a})\psi_{il}\nn\\
&=&\pr_{\tt}N_{W,s,ij}^{(p)} -\pr_{\tt}(\widecheck{\g}^{\a\b})\pr_{\a}\pr_{\b}\phiss{ij}{p} + \sum_{k,l}\dk^{\leq 2}\Ga_g \dk^{\leq 1}\phiss{kl}{p}
\eea
in the derivation of which we have used \eqref{eq:assumptionsonregulartripletinperturbationsofKerr:0}, Lemma \ref{lem:estimatesforMialphaj:Kerrpert},  the formula \eqref{esti:commutatorBoxgandT:general:0} 
and the fact from Lemma \ref{lem:scalarizedTeukolskywavetransportsysteminKerrperturbation:Omi} that all the coefficients in the formulas \eqref{eq:linearterms:ScalarizedTeuSys:general:Kerrperturbation} for $L_{s,ij}^{(p)}$ are $(\tau,\tphi)$-independent. 
Next, commuting $\pr_{\tt}$ with the globally extended system of coupled wave equations \eqref{eq:waveeqwidetildepsi1:gtilde:Teu}, we deduce
\bea\lab{eq:waveeqwidetildepsi1:gtilde:Teu:pr_tt}
&&\bigg({\square}_{{\g}_{\chi_{\tau_1,\tau_2}}} -\frac{4-2\de_{p0}}{|q|^{2}}\bigg)\psi^{(p),\pr_{\tt}}_{s,ij}\nn\\
&=&\chi_{\tau_1, \tau_2}\big(\widehat{S}(\psi_s^{(p),\pr_{\tau}})_{ij} +(\widehat{Q}\psi_s^{(p),\pr_{\tau}})_{ij}\big)
+(1-\chi_{\tau_1, \tau_2})\big(\widehat{S}_K(\psi_s^{(p),\pr_{\tau}})_{ij} +(\widehat{Q}_K\psi_s^{(p),\pr_{\tau}})_{ij} \big)\nn\\
&&
+(1-\chi_{\tau_1, \tau_2})f_p\psi^{(p),\pr_{\tt}}_{s,ij}
+\chi_{\tau_1, \tau_2}^{(1)}L^{(p),\pr_{\tt}}_{s,ij}+\widehat{F}_{total,s,ij}^{(p),\pr_{\tt}} 
\eea
on $\MM$, where the scalars $\psi^{(p),\pr_{\tt}}_{s,ij}$, $s=\pm 2, i,j=1,2,3,  p=0,1,2$,  are defined by
\beaa
\psissh{ij}{p}{\pr_{\tt}}:=\pr_{\tt}\psiss{ij}{p},
\eeaa
where, as above, $\widehat{S}(\psi_s^{(p),\pr_{\tau}})_{ij}$ and $(\widehat{Q}\psi_s^{(p),\pr_{\tau}})_{ij}$ are obtained by replacing each $\psiss{ij}{p}$ with $\psissh{ij}{p}{\pr_{\tt}}$ in the expanded forms of $\widehat{S}(\psi_s^{(p)})_{ij}$ and $(\widehat{Q}\psi_s^{(p)})_{ij}$, where $L^{(p),\pr_{\tt}}_{s,ij}$ is defined as above, and where
\bea
\lab{def:Ftotalsijp:pr_tt}
\widehat{F}_{total,s,ij}^{(p),\pr_{\tt}}&=&\pr_{\tt}\widehat{F}_{total,s,ij}^{(p)}
-\pr_{\tt}(\widecheck{\g}^{\a\b})\pr_{\a}\pr_{\b}\psiss{ij}{p} 
+ \sum_{k,l}\dk^{\leq 2}\Ga_g \dk^{\leq 1}\psiss{kl}{p} \nn\\
&&
-\pr_{\tt}(\chi_{\tau_1, \tau_2})f_p\psi^{(p)}_{s,ij}+\pr_{\tt}(\chi_{\tau_1, \tau_2}^{(1)})L^{(p)}_{s,ij}
\eea
which is obtained in a similar manner as \eqref{def:Nwsijp:pr_tt}, using in addition the fact that 
\beaa
{S}(\psi)_{ij}-{S}_K(\psi)_{ij}=\sum_{k,l}\Ga_g\dk\psi_{kl}, \qquad  (\widehat{Q}\psi)_{ij}-({\widehat{Q}}_K\psi)_{ij}=\sum_{k,l}\dk^{\leq 1}\Ga_g\psi_{kl},
\eeaa
in view of \eqref{SandV} and the assumption \eqref{eq:assumptionsonregulartripletinperturbationsofKerr:0}.

Relying on Remark \ref{rem:prtauandprtphihatpreservesk2C}, we infer the existence, for $s=\pm 2$, $p=0,1,2$, of tensors $\pmb\phi_s^{(p),\pr_\tau}\in\sk_2(\mathbb{C})$ such that $\phissh{ij}{p}{\pr_{\tt}}=\pmb\phi_s^{(p),\pr_\tau}(\Om_i, \Om_j)$. In view of the Teukolsky wave/transport systems \eqref{eq:ScalarizedTeuSys:general:Kerrperturbation:pr_tt}  \eqref{eq:ScalarizedQuantitiesinTeuSystem:Kerrperturbation:pr_tt} for  $\phissh{ij}{p}{\pr_{\tt}}$, we deduce from Lemma \ref{lemma:formoffirstordertermsinscalarazationtensorialwaveeq} that $\pmb\phi_s^{(p),\pr_\tau}$ are solutions to the tensorial Teukolsky wave/transport systems \eqref{eq:TensorialTeuSysandlinearterms:rescaleRHScontaine2:general:Kerrperturbation}  \eqref{def:TensorialTeuScalars:wavesystem:Kerrperturbation} with $\N_{W}^{(p)}$ and $\N_{T}^{(p)}$ being replaced by $\N_{W}^{(p),\pr_{\tt}}$ and $\N_{T}^{(p),\pr_{\tt}}$ respectively, where $\N_{W}^{(p),\pr_{\tt}}, \N_{T}^{(p),\pr_{\tt}}\in\sk_2(\mathbb{C})$ are such that $N_{W,s,ij}^{(p),\pr_{\tt}}=\N_{W}^{(p),\pr_{\tt}}(\Om_i, \Om_j)$ and $N_{T,s,ij}^{(p),\pr_{\tt}}=\N_{T}^{(p),\pr_{\tt}}(\Om_i, \Om_j)$. Furthermore, $\psi^{(p),\pr_{\tt}}_{s,ij}$ satisfies the wave system \eqref{eq:waveeqwidetildepsi1:gtilde:Teu:pr_tt}, and since $\psi^{(p),\pr_{\tt}}_{s,ij}=\pr_{\tt}\psiss{ij}{p}$ and $\phi^{(p),\pr_{\tt}}_{s,ij}=\pr_{\tt}\phiss{ij}{p}$, $\psi^{(p),\pr_{\tt}}_{s,ij}$ satisfies \eqref{eq:causlityrelationsforwidetildepsi1} with $(\phi^{(p)}_{s,ij}, \psi^{(p)}_{s,ij})$ being replaced by $(\phi^{(p),\pr_{\tt}}_{s,ij}, \psi^{(p),\pr_{\tt}}_{s,ij})$. We may thus apply Theorem \ref{thm:EMF:systemofTeuscalarized:order0:final} to infer, 
for $s=\pm 2$ and $\de\in(0, \frac{1}{3}]$,
\bea
\lab{eq:EMFtotalp:sumup:rweightscontrolled:pm2:prtt}
\EMFtotalh{s}{\nab_{\pr_{\tt}}^{\leq 1}} 
&\les&\IE{\nab_{\pr_{\tt}}^{\leq 1}\pmb\phi_s, \nab_{\pr_{\tt}}^{\leq 1}\pmb\psi_s} +\A[\nab_{\pr_{\tt}}^{\leq 1}\pmb\psi_{s}](\Iti)+\A[\nab_{\pr_{\tt}}^{\leq 1}\pmb\phi_{s}](\tau_1,\tau_2)\nn\\
&&+\sum_{p=0}^2\Errdefect[\nab_{\pr_{\tt}}^{\leq 1}\psis{p}]  +\NNttotalp{s}
+\NNttotalh{s}{\pr_{\tt}},
\eea
where $\NNttotalh{s}{\pr_{\tt}}$ is defined by the same formulas as \eqref{expression:NNttotalps:pm2} but with the following replacements
\beaa
(\phiss{ij}{p}, \psiss{ij}{p}, N_{W,s,ij}^{(p)}, \widehat{F}_{total,s,ij}^{(p)}, N_{T,s,ij}^{(p)})\to (\phissh{ij}{p}{\pr_{\tt}}, \psissh{ij}{p}{\pr_{\tt}}, N_{W,s,ij}^{(p), \pr_{\tt}}, \widehat{F}_{total,s,ij}^{(p),\pr_{\tt}}, N_{T,s,ij}^{(p), \pr_{\tt}}),
\eeaa
with $N_{W, s,ij}^{(p), \pr_{\tt}}$, $\widehat{F}_{total,s,ij}^{(p),\pr_{\tt}}$ and $N_{T,s,ij}^{(p), \pr_{\tt}}$ given as in \eqref{def:Nwsijp:pr_tt},  \eqref{def:Ftotalsijp:pr_tt} and \eqref{def:Nt+-2ijp:pr_tt} respectively. In view of \eqref{def:Nwsijp:pr_tt},  \eqref{def:Ftotalsijp:pr_tt} and \eqref{def:Nt+-2ijp:pr_tt}, the formulas for $N_{W, s,ij}^{(p), \pr_{\tt}}$, $\widehat{F}_{total,s,ij}^{(p),\pr_{\tt}}$ and $N_{T,s,ij}^{(p), \pr_{\tt}}$ can alternatively be given by
\bsub
\lab{def:Nt+-2ijp:NwFtotal:pr_tt:copy}
\bea
\lab{def:Nwsijp:pr_tt:copy}
N_{W,s,ij}^{(p),\pr_{\tt}}
&=&\pr_{\tt}N_{W,s,ij}^{(p)} +\widetilde{N_{W,s,ij}^{(p),\pr_{\tt}}},\\
\lab{def:Ftotalsijp:pr_tt:copy}
\widehat{F}_{total,s,ij}^{(p),\pr_{\tt}}&=&\pr_{\tt}\widehat{F}_{total,s,ij}^{(p)}+\widetilde{\widehat{F}_{total,s,ij}^{(p),\pr_{\tt}}},\\
\lab{def:Nt+2ijp:pr_tt:copy}
N_{T,+2,ij}^{(p), \pr_{\tt}}
&=&\pr_{\tt}N_{T,+2,ij}^{(p)}+\widetilde{N_{T,+2,ij}^{(p), \pr_{\tt}}},\\
\lab{def:Nt-2ijp:pr_tt:copy}
N_{T,-2,ij}^{(p), \pr_{\tt}}&=&\pr_{\tt}N_{T,-2,ij}^{(p)}+\widetilde{N_{T,-2,ij}^{(p), \pr_{\tt}}},
\eea
\esub
where 
\bsub
\lab{def:Nt+-2ijp:NwFtotal:pr_tt:difference}
\bea
\widetilde{N_{W,s,ij}^{(p),\pr_{\tt}}}&=&-\pr_{\tt}(\widecheck{\g}^{\a\b})\pr_{\a}\pr_{\b}\phiss{ij}{p} +\sum_{k.l} \dk^{\leq 2}\Ga_g \dk^{\leq 1}\phiss{kl}{p},\\
\widetilde{\widehat{F}_{total,s,ij}^{(p),\pr_{\tt}}}&=&-\pr_{\tt}(\widecheck{\g}^{\a\b})\pr_{\a}\pr_{\b}\psiss{ij}{p} + \sum_{k.l}\dk^{\leq 2}\Ga_g \dk^{\leq 1}\psiss{kl}{p} \nn\\
&&-\pr_{\tt}(\chi_{\tau_1, \tau_2})f_p\psi^{(p)}_{s,ij}+\pr_{\tt}(\chi_{\tau_1, \tau_2}^{(1)})L^{(p)}_{s,ij},\\
\widetilde{N_{T,+2,ij}^{(p), \pr_{\tt}}}&=&\sum_{k.l}r\dk^{\leq 1}\Ga_b \dk^{\leq 1}\phipluss{kl}{p} +\sum_{k.l}r^2\dk^{\leq 1}\Ga_g\pr_{\tt}\phipluss{kl}{p},\\
\widetilde{N_{T,-2,ij}^{(p), \pr_{\tt}}}&=&\sum_{k.l}r\dk^{\leq 1}\Ga_g\dk^{\leq 1} \phiminuss{kl}{p}.
\eea
\esub

It then follows that the last two terms on the RHS of \eqref{eq:EMFtotalp:sumup:rweightscontrolled:pm2:prtt} satisfy
\bea
\lab{eq:NNttotalpsandhs:prtt:intermsofwidetildeNNt}
\NNttotalp{s}
+\NNttotalh{s}{\pr_{\tt}}
\les\NNttotalph{s}{\nab_{\pr_{\tt}}^{\leq 1}}+\widetilde{\NNttotalh{s}{\pr_{\tt}}}
\eea
where, in view of \eqref{def:Nt+-2ijp:NwFtotal:pr_tt:copy} and the definition \eqref{expression:NNttotalps:pm2} of $\NNttotalp{s}$,
\begin{align}
\lab{NNtpm2total:pr_tt:diff:expression}
\widetilde{\NNttotalh{s}{\pr_{\tt}}}
={}&\sum_{p=0}^2\widetilde{\mathcal{N}}[\psish{p}{\pr_{\tt}},\widetilde{\widehat{F}_{total,s}^{(p),{\pr_{\tt}}}}]\nn\\
&+\sum_{p=0}^2\widehat{\mathcal{N}}'[\phish{p}{\pr_{\tt}}, \widetilde{\N_{W, s}^{(p),{\pr_{\tt}}}}](\tau_1, \tau_2)+\sum_{p=0}^1\widehat{\mathcal{N}}'[ \phish{p}{\pr_{\tt}}, r^{-2}\widetilde{\N_{T, s}^{(p),{\pr_{\tt}}}}](\tau_1, \tau_2)\nn\\
&+\sum_{p=0}^1\int_{\MM_{r\leq 12m}(\tau_1,\tau_2)}|\pr^{\leq 1}\widetilde{\N_{T, s}^{(p),{\pr_{\tt}}}}|^2\nn\\
&+\sum_{p=0}^1\int_{\MM(\tau_1,\tau_2)} \Big(r^{-1+\de}|\widetilde{\N_{W, s}^{(p),{\pr_{\tt}}}}|+r^{-2+\de}|\pr\widetilde{\N_{T, s}^{(p),{\pr_{\tt}}}}|\Big)|\phish{p}{\pr_{\tt}}|.
\end{align}
In view of the formula \eqref{eq:defofNNhat'} for $\widehat{\NN}'[\c, \c](\c,\c)$ and the above formulas \eqref{def:Nt+-2ijp:NwFtotal:pr_tt:difference}, and using Lemma \ref{lemma:basiclemmaforcontrolNLterms:bis} and \eqref{eq:controloflinearizedinversemetriccoefficients},  we infer that the last three lines of \eqref{NNtpm2total:pr_tt:diff:expression} are bounded by\footnote{Note also that $\int_{\Mntrap(\tau_1, \tau_2)}r^{-1-\de}|\phissh{ij}{p}{\pr_{\tt}}|^2\les\M_\de[\pmb\phi_s^{(p)}](\tau_1, \tau_2)$ since $\phissh{ij}{p}{\pr_{\tt}}=\pr_{\tt}\phiss{ij}{p}$.} \footnote{While $\de\leq\frac{1}{2}$ is enough for most estimates in this paper, the stronger constraint $\de\leq\frac{1}{3}$ is needed to bound the last line of \eqref{NNtpm2total:pr_tt:diff:expression} by $\ep\EM_\de[\pr^{\leq 1}\pmb\phi_s](\tau_1, \tau_2)$.}
 $\ep\EM_\de[\pr^{\leq 1}\pmb\phi_s](\tau_1, \tau_2)$, which together with \eqref{eq:EMFtotalp:sumup:rweightscontrolled:pm2:prtt} and \eqref{eq:NNttotalpsandhs:prtt:intermsofwidetildeNNt} yields 
\beaa
&&\EMFtotalh{s}{\nab_{\pr_{\tt}}^{\leq 1}} \nn\\
&\les&\IE{\nab_{\pr_{\tt}}^{\leq 1}\pmb\phi_s, \nab_{\pr_{\tt}}^{\leq 1}\pmb\psi_s}+\NNttotalph{s}{\nab_{\pr_{\tt}}^{\leq 1}}+\A[\nab_{\pr_{\tt}}^{\leq 1}\pmb\psi_{s}](\Iti)+\A[\nab_{\pr_{\tt}}^{\leq 1}\pmb\phi_{s}](\tau_1,\tau_2)\nn\\
&& 
+\sum_{p=0}^2\Errdefect[\nab_{\pr_{\tau}}^{\leq 1}\psis{p}]+\ep\EM_\de[\pr^{\leq 1}\pmb\phi_s](\tau_1, \tau_2)
+\sum_{p=0}^2\widetilde{\mathcal{N}}[\psish{p}{\pr_{\tt}},\widetilde{\widehat{F}_{total,s}^{(p),{\pr_{\tt}}}}].
\eeaa
Using the formula of $\widetilde{\widehat{F}_{total,s,ij}^{(p),\pr_{\tt}}}$ in \eqref{def:Nt+-2ijp:NwFtotal:pr_tt:difference} and the following estimate 
\beaa
\sum_{p=0}^2\sum_{i,j}\widetilde{\mathcal{N}}\bigg[\psissh{ij}{p}{\pr_{\tt}},-\pr_{\tt}(\widecheck{\g}^{\a\b})\pr_{\a}\pr_{\b}\psiss{ij}{p} + \sum_{k,l}\dk^{\leq 2}\Ga_g \dk^{\leq 1}\psiss{kl}{p}\bigg]
\les \ep \EMFtotalhps{s}{\pr^{\leq 1}}
\eeaa
which follows from \eqref{def:NNtintermsofNNtMora:NNtEner:NNtaux:wavesystem:EMF}-\eqref{def:NNtMora:NNtEner:NNtaux:wavesystem:EMF}, \eqref{eq:controloflinearizedinversemetriccoefficients} and Lemmas \ref{lemma:basiclemmaforcontrolNLterms:bis} and \ref{lem:gpert:MMtrap},  we infer
\bea
\lab{eq:EMFtotalp:sumup:middlestep:pm2:prtt}
&&\EMFtotalh{s}{\nab_{\pr_{\tt}}^{\leq 1}} \nn\\
&\les&\IE{\nab_{\pr_{\tt}}^{\leq 1}\pmb\phi_s, \nab_{\pr_{\tt}}^{\leq 1}\pmb\psi_s} +\NNttotalph{s}{\nab_{\pr_{\tt}}^{\leq 1}}+\A[\nab_{\pr_{\tt}}^{\leq 1}\pmb\psi_{s}](\Iti)+\A[\nab_{\pr_{\tt}}^{\leq 1}\pmb\phi_{s}](\tau_1,\tau_2)\nn\\
&& 
+\sum_{p=0}^2\Errdefect[\nab_{\pr_{\tau}}^{\leq 1}\psis{p}]+\ep \EMFtotalhps{s}{\pr^{\leq 1}} 
\nn\\
&&+\sum_{p=0}^2\widetilde{\mathcal{N}}\big[\psish{p}{\pr_{\tt}},\pr_{\tt}\chi_{\tau_1,\tau_2}^{(1)}\L_{s}^{(p)}\big]+\sum_{p=0}^2\widetilde{\mathcal{N}}\big[\psish{p}{\pr_{\tt}},-\pr_{\tt}(\chi_{\tau_1, \tau_2})f_p\psis{p}\big].
\eea

Next, we estimate the last two terms on the RHS of \eqref{eq:EMFtotalp:sumup:middlestep:pm2:prtt}. By \cite[inequality (7.148)]{MaSz24}, we have, for $\tau\geq\tmic$,
\beaa
&&\bigg|\int_{\Mtrap({\tmic},\tau)}\Re\Big(\ov{{|q|^{-2}}\Opw(\Theta_n){(\qs F)}_{ij}}V_n\Opw(\Theta_n)\psi_{ij}\Big)\bigg|\nn\\
&\les&\left(\int_{\Mtrap}|\pr^{\leq 1} F_{ij}|^2\right)^{\frac{1}{2}}\left(\int_{\Mtrap(\Iti)} |\psi_{ij}|^2\right)^{\frac{1}{2}}
+\bigg(\int_{\Mtrap}|F_{ij}|^2\bigg)^{\frac{1}{2}}\bigg(\sup_{\tau\in\Iti}\E[\psi_{ij}](\tau)\bigg)^{\frac{1}{2}}, 
\eeaa
 and hence, together with \eqref{def:NNtintermsofNNtMora:NNtEner:NNtaux:wavesystem:EMF}-\eqref{def:NNtMora:NNtEner:NNtaux:wavesystem:EMF}, we infer
 \bea
 \lab{esti:controlofNNtener:byprF}
 \NNt[\pmb \psi, \pmb F]&\les&\sum_{i,j}\sup_{\tau\geq\tmic}\bigg|\int_{\Mntrap(\tmic, \tau)}{\Re\Big(F_{ij}\ov{\pr_{\tau}\psi_{ij}}\Big)}\bigg|+\int_{\Mntrap(\Iti)}r^{-1}|\F||\dk^{\leq 1}\pmb\psi|\nn\\
 &&+\bigg(\int_{\Mtrap}|\pmb F|^2\bigg)^{\frac{1}{2}}\bigg(\EM[\pmb \psi](\Iti)\bigg)^{\frac{1}{2}} +\int_{\MM(\Iti)}|\F|^2\nn\\
 &&+ \left(\int_{\Mtrap}|\pr^{\leq 1} \pmb F|^2\right)^{\frac{1}{2}}\left(\int_{\Mtrap(\Iti)} |\pmb \psi|^2\right)^{\frac{1}{2}}.
\eea
Now, using \eqref{esti:controlofNNtener:byprF} together with \eqref{eq:expression:Lsijp:copy} which yields
\bea
\lab{eq:roughformofLspij:sect10pf}
L_{s,ij}^{(p)}=O(r^{-2}) \Big(\dk^{\leq 1}\phiss{kl}{0}, \dk^{\leq 1}\phiss{kl}{1}, \phiss{ij}{2}\Big),\qquad \forall \,p=0,1,2,
\eea
and the fact from \eqref{eq:propertieschi:thisonetodealwithRHSofTeukolsky} that $\pr_{\tt}\chi_{\tau_1,\tau_2}^{(1)}$ is supported in $[\tau_1,\tau_1+1]\cup[\tau_2-3,\tau_2-2]$, we infer
\begin{align}
\lab{control:NNtterm:commuteprtau:Lsp}
\sum_{p=0}^2\widetilde{\mathcal{N}}[\psish{p}{\pr_{\tt}},\pr_{\tt}\chi_{\tau_1,\tau_2}^{(1)}\L_{s}^{(p)}]
\les{}& \Bigg(\sum_{p=0}^2\EMF[\nab_{\pr_{\tt}}^{\leq 1}\psis{p}](\Iti)\Bigg)^{\frac{1}{2}}\big(\EMFtotalps{s}\big)^{\frac{1}{2}}
+\EMFtotalps{s}
\nn\\
&+\big(\EMFtotalhps{s}{\pr^{\leq 1}}\big)^{\frac{1}{2}} \big(\Ao[\nab_{\pr_{\tt}}^{\leq 1}\pmb\psi_{s}](\Iti)\big)^{\frac{1}{2}}\nn\\
\les {}&\Big(\EMFtotalhps{s}{\nab_{\pr_{\tt}}^{\leq 1}}\Big)^{\frac{1}{2}}\big(\EMFtotalps{s}\big)^{\frac{1}{2}}\nn\\
&+\big(\EMFtotalhps{s}{\pr^{\leq 1}}\big)^{\frac{1}{2}} \big(\Ao[\nab_{\pr_{\tt}}^{\leq 1}\pmb\psi_{s}](\Iti)\big)^{\frac{1}{2}},
 \end{align}
 where we have used in the second step the following fact  which follows from \eqref{def:EMFtotalps:pm2}
  \beaa
 \sum_{p=0}^2\EMF[\nab_{\pr_{\tt}}^{\leq 1}\psis{p}](\Iti) \les \EMFtotalhps{s}{\nab_{\pr_{\tt}}^{\leq 1}}. 
 \eeaa
 In a similar manner, we deduce
 \begin{align}
\lab{control:NNtterm:commuteprtau:LOT}
\sum_{p=0}^2\widetilde{\mathcal{N}}[\psish{p}{\pr_{\tt}},-\pr_{\tt}(\chi_{\tau_1, \tau_2})f_p\psis{p}]
\les{} & \Big(\EMFtotalhps{s}{\nab_{\pr_{\tt}}^{\leq 1}}\Big)^{\frac{1}{2}}\big(\EMFtotalps{s}\big)^{\frac{1}{2}}\nn\\
&+\big(\EMFtotalhps{s}{\pr^{\leq 1}}\big)^{\frac{1}{2}} \big(\Ao[\nab_{\pr_{\tt}}^{\leq 1}\pmb\psi_{s}](\Iti)\big)^{\frac{1}{2}}.
 \end{align}

Substituting \eqref{control:NNtterm:commuteprtau:Lsp} and \eqref{control:NNtterm:commuteprtau:LOT} into \eqref{eq:EMFtotalp:sumup:middlestep:pm2:prtt} to control the last line of \eqref{eq:EMFtotalp:sumup:middlestep:pm2:prtt}, 
 we infer
 \beaa
\EMFtotalh{s}{\nab_{\pr_{\tt}}^{\leq 1}} 
&\les&\IE{\nab_{\pr_{\tt}}^{\leq 1}\pmb\phi_s, \nab_{\pr_{\tt}}^{\leq 1}\pmb\psi_s} +\NNttotalph{s}{\nab_{\pr_{\tt}}^{\leq 1}}+\sum_{p=0}^2\Errdefect[\nab_{\pr_{\tau}}^{\leq 1}\psis{p}]\nn\\
&& 
+\A[\nab_{\pr_{\tt}}^{\leq 1}\pmb\psi_{s}](\Iti)+\A[\nab_{\pr_{\tt}}^{\leq 1}\pmb\phi_{s}](\tau_1,\tau_2)+\ep \EMFtotalhps{s}{\pr^{\leq 1}}
\nn\\
&&+\Big(\EMFtotalhps{s}{\nab_{\pr_{\tt}}^{\leq 1}}\Big)^{\frac{1}{2}}\big(\EMFtotalps{s}\big)^{\frac{1}{2}}\\
&&+\big(\EMFtotalhps{s}{\pr^{\leq 1}}\big)^{\frac{1}{2}} \big(\Ao[\nab_{\pr_{\tt}}^{\leq 1}\pmb\psi_{s}](\Iti)\big)^{\frac{1}{2}},
\eeaa
and hence
\begin{align}\lab{eq:EMFtotalp:sumup:finalestimate:pm2:prttandprtphi:sofaronlyprtt}
\EMFtotalh{s}{\nab_{\pr_{\tt}}^{\leq 1}}
\les{}&\IE{\nab_{\pr_{\tt}}^{\leq 1}\pmb\phi_s, \nab_{\pr_{\tt}}^{\leq 1}\pmb\psi_s} +\NNttotalph{s}{\nab_{\pr_{\tt}}^{\leq 1}}+\ep \EMFtotalhps{s}{\pr^{\leq 1}} \nn\\
&+\big(\EMFtotalhps{s}{\pr^{\leq 1}}\big)^{\frac{1}{2}} \big(\EMFtotalps{s}+\Ao[\nab_{\pr_{\tt}}^{\leq 1}\pmb\psi_{s}](\Iti)\big)^{\frac{1}{2}}\nn\\
&
+\sum_{p=0}^2\Errdefect[\nab_{\pr_{\tau}}^{\leq 1}\psis{p}]+\A[\nab_{\pr_{\tt}}^{\leq 1}\pmb\psi_{s}](\Iti)+\A[\nab_{\pr_{\tt}}^{\leq 1}\pmb\phi_{s}](\tau_1,\tau_2).
\end{align}

Commuting further with $\pr_{\tt}$, and by an inductive argument, we infer that, for any $\reg\leq 14$,
\begin{align}
\lab{eq:EMFtotalp:sumup:finalestimate:pm2:prtt}
\EMFtotalh{s}{\nab_{\pr_{\tt}}^{\leq \reg}} 
\les{}&\IE{\pr^{\leq \reg}\pmb\phi_s, \pr^{\leq \reg}\pmb\psi_s} +\NNttotalph{s}{\nab_{\pr_{\tt}}^{\leq \reg}}+\ep \EMFtotalhps{s}{\pr^{\leq \reg}}\nn\\
&+\big(\EMFtotalhps{s}{\pr^{\leq\reg}}\big)^{\frac{1}{2}} \big(\EMFtotalhps{s}{\pr^{\leq\reg-1}}+\Ao[\nab_{\pr_{\tt}}^{\leq\reg}\pmb\psi_{s}](\Iti)\big)^{\frac{1}{2}}\nn\\
&  +\sum_{p=0}^2\Errdefect[\nab_{\pr_{\tau}}^{\leq \reg}\psis{p}] +\A[\nab_{\pr_{\tt}}^{\leq\reg}\pmb\psi_{s}](\Iti)+\A[\nab_{\pr_{\tt}}^{\leq\reg}\pmb\phi_{s}](\tau_1,\tau_2).
\end{align}


\subsubsection{Commutation with $\chi_0\prtphihat$}
\lab{subsubsect:commutationwithprtphihat:conditionalEMF}


Since $\chi_0=\chi_0(r)$ given in \eqref{def:cutoffcuntionchi0} is supported in $r\leq 12m$, we introduce, as in  Section \ref{subsec:definitionoflowerordertermsintheerrors:scalsyst}, the notation $\Gac$ for error terms satisfying 
\bea\lab{eq:decaypropertiesofGac:microlocalregion:bisforcommprphi}
|\dk^{\leq 15}\Gac|\les \ep\tau^{-1-\dec}\qquad\textrm{on}\,\,\MM_{r\leq 12m}.
\eea
Now, recalling Definition \ref{definition:hor-Lie-derivative} for the horizontal Lie derivative $\Lieb$, we have, for any $U\in\sk_2$, 
\bea\lab{eq:formofcommutatorswithLiebprtphionrgeq12m}
\bsplit
&[\Lieb_{\pr_{\tphi}}, \squared_2]U=\pr^{\leq 2}(\Gac\c U)\quad\textrm{on}\quad\MM_{r\leq 12m}(\tau_1, \tau_2), \\
&[\Lieb_{\pr_{\tphi}}, \nab_4]U,\,\,  [\Lieb_{\pr_{\tphi}}, \nab_3]U,\,\,  [\Lieb_{\pr_{\tphi}}, \nab]U=\pr^{\leq 1}(\Gac\c U)\quad\textrm{on}\quad\MM_{r\leq 12m}(\tau_1, \tau_2),
\end{split}
\eea
which follows immediately from Proposition 4.3.4 and Lemma C.5.2 in \cite{GKS22}. Thus, commuting the tensorial Teukolsky wave-transport system \eqref{eq:TensorialTeuSysandlinearterms:rescaleRHScontaine2:general:Kerrperturbation} \eqref{def:TensorialTeuScalars:wavesystem:Kerrperturbation} with $\chi_0\Lieb_{\pr_{\tphi}}$ and using \eqref{eq:formofcommutatorswithLiebprtphionrgeq12m}, we infer in $\MM(\tau_1, \tau_2)$, for $s=\pm 2$ and $p=0,1,2$,
\bea\lab{eq:TensorialTeuSysandlinearterms:rescaleRHScontaine2:general:Kerrperturbation:commLiebprtphi}
\bigg(\squared_2 -\frac{4ia\cos\th}{|q|^2}\nab_{\pr_{\tt}}- \frac{4-2\de_{p0}}{\qs}\bigg)\chi_0\Lieb_{\pr_{\tphi}}\phis{p} = \L_{s}^{(p)}[\chi_0\Lieb_{\pr_{\tphi}}\pmb\phi_{s}]+\N_{W,s}^{(p),\chi_0\Lieb_{\pr_{\tphi}}}, 
\eea
and in $\MM(\tau_1, \tau_2)$, for $p=0,1$,
\bsub\lab{def:TensorialTeuScalars:wavesystem:Kerrperturbation:Kerrperturbation:commLiebprtphi}
\begin{align}
\nab_3 \left(\frac{r\bar{q}}{q}\left(\frac{r^2}{|q|^2}\right)^{p-2}\chi_0\Lieb_{\pr_{\tphi}}\pmb\phi_{+2}^{(p)}\right) =& \frac{\bar{q}}{rq}\left(\frac{r^2}{|q|^2}\right)^{p-1}\chi_0\Lieb_{\pr_{\tphi}}\pmb\phi_{+2}^{(p+1)}+\N_{T,+2}^{(p),\chi_0\Lieb_{\pr_{\tphi}}},\\
\nab_4\left(\frac{rq}{\bar{q}}\left(\frac{{r^2}}{|q|^2}\right)^{p-2}\chi_0\Lieb_{\pr_{\tphi}}\pmb\phi_{-2}^{(p)}\right) =& \frac{q}{r\bar{q}}\left(\frac{r^2}{|q|^2}\right)^{p-1}\frac{\De}{\qs}\chi_0\Lieb_{\pr_{\tphi}}\pmb\phi_{-2}^{(p+1)}+\N_{T,-2}^{(p),\chi_0\Lieb_{\pr_{\tphi}}},
\end{align}
\esub
where, with $\pmb\phi_s$ denoting $(\pmb\phi_s^{(p)})_{p=0,1,2}$,
\bsub\lab{eq:formofRHSNWandNTintensorialTeukolskytau1tau2aftercommutationchi0Liebtphi}
\bea
\N_{W,s}^{(p),\chi_0\Lieb_{\pr_{\tphi}}} &=& \chi_0\Lieb_{\pr_{\tphi}}\N_{W,s}^{(p)}+\mathbf{1}_{r\leq 12m}\pr^{\leq 2}(\Gac\c\pmb\phi_s)+O(1)\mathbf{1}_{11m\leq r\leq 12m}\pr^{\leq 2}\pmb\phi_s,\\
\N_{T,s}^{(p),\chi_0\Lieb_{\pr_{\tphi}}} &=& \chi_0\Lieb_{\pr_{\tphi}}\N_{T,s}^{(p)}+\mathbf{1}_{r\leq 12m}\pr^{\leq 1}(\Gac\c\pmb\phi_s)+O(1)\mathbf{1}_{11m\leq r\leq 12m}\pr^{\leq 1}\pmb\phi_s.
\eea
\esub
In the above derivation of \eqref{eq:TensorialTeuSysandlinearterms:rescaleRHScontaine2:general:Kerrperturbation:commLiebprtphi}, note that we have used the fact that all the coefficients on the RHS of \eqref{eq:tensor:Lsn:onlye_2present:general:Kerrperturbation} are independent of the coordinates $\tau$ and 
$\tphi$.

Next, recalling that $\prtphihat$ is given in Definition \ref{def:widehatprtphi}, we have in view of Lemma \ref{lemma:differencebetweenwidehatprtphiandLiebprtrphiisanerrorterm},
\bea\lab{eq:differencebetweenwidehatprtphiandLiebprtrphiisanerrorterm:onMMtau1tau2rleq12m}
\widehat{\pr}_{\tphi}(\phi_s)_{ij}=\Lieb_{\pr_{\tphi}}\pmb\phi_s(\Om_i, \Om_j)+\Gac\c\pmb\phi_s, \quad \forall\, i,j,  \quad\textrm{on}\quad\MM_{r\leq 12m}(\tau_1,\tau_2). 
\eea
Thus, commuting the transport equations \eqref{eq:ScalarizedQuantitiesinTeuSystem:Kerrperturbation}, the scalarized Teukolsky wave equations \eqref{eq:ScalarizedTeuSys:general:Kerrperturbation}, and the globally extended system of coupled wave equations \eqref{eq:waveeqwidetildepsi1:gtilde:Teu} with $\chi_0\prtphihat$, and taking \eqref{eq:TensorialTeuSysandlinearterms:rescaleRHScontaine2:general:Kerrperturbation:commLiebprtphi} \eqref{def:TensorialTeuScalars:wavesystem:Kerrperturbation:Kerrperturbation:commLiebprtphi} \eqref{eq:formofRHSNWandNTintensorialTeukolskytau1tau2aftercommutationchi0Liebtphi} and \eqref{eq:differencebetweenwidehatprtphiandLiebprtrphiisanerrorterm:onMMtau1tau2rleq12m} into account, we get the transport equations  for $p=0,1$
\bsub
\lab{eq:ScalarizedQuantitiesinTeuSystem:Kerrperturbation:pr_tphi}
\bea
&& e_3\bigg(\frac{r\bar{q}}{q}\bigg({\frac{r^2}{|q|^2}}\bigg)^{p-2}\phiplussh{ij}{p}{\chi_{0}\prtphihat} \bigg) - \frac{r\bar{q}}{q}\bigg({\frac{r^2}{|q|^2}}\bigg)^{p-2}\Big(M_{i3}^k \phiplussh{kj}{p}{\chi_{0}\prtphihat} + M_{j3}^k \phiplussh{ik}{p}{\chi_{0}\prtphihat}\Big)\nn\\
&=&\frac{\bar{q}}{rq}\bigg({\frac{r^2}{|q|^2}}\bigg)^{p-1}\phiplussh{ij}{p+1}{\chi_{0}\prtphihat}+\chi_{0}\prtphihat N_{T,+2,ij}^{(p)}+\widetilde{N_{T,+2,ij}^{(p), \chi_{0}\prtphihat}},\\
&& e_4\bigg(\frac{rq}{\bar{q}}\bigg(\frac{r^2}{|q|^2}\bigg)^{p-2}\phiminussh{ij}{p}{\chi_{0}\prtphihat}\bigg) - \frac{rq}{\bar{q}}\bigg(\frac{r^2}{|q|^2}\bigg)^{p-2}\Big(M_{i4}^k \phiminussh{kj}{p}{\chi_{0}\prtphihat} + M_{j4}^k \phiminussh{ik}{p}{\chi_{0}\prtphihat}\Big)\nn\\
&=& \frac{q}{r\bar{q}}\left(\frac{r^2}{|q|^2}\right)^{p-1}\frac{\De}{\qs}\phiminussh{ij}{p+1}{\chi_{0}\prtphihat} +\chi_{0}\prtphihat N_{T,-2,ij}^{(p)}+\widetilde{N_{T,-2,ij}^{(p), \chi_{0}\prtphihat}},
\eea
\esub
the scalarized Teukolsky wave equations for $p=0,1,2$
\bea
\lab{eq:ScalarizedTeuSys:general:Kerrperturbation:pr_tphi} 
&&\square_\g\phissh{ij}{p} {\chi_{0}\prtphihat} -\frac{4-2\de_{p0}}{\qs} \phissh{ij}{p}{\chi_{0}\prtphihat}\nn\\
&=&\big(\widehat{S}(\phi_s^{(p),\chi_0\prtphihat})_{ij} +(\widehat{Q}\phi_s^{(p),\chi_0\prtphihat})_{ij}\big)+ {L_{s,ij}^{(p),\chi_{0}\prtphihat}}+\chi_{0}\prtphihat N_{W,s,ij}^{(p)}+\widetilde{N_{W,s,ij}^{(p),\chi_{0}\prtphihat}}, 
\eea
and the globally extended system of coupled wave equations  for $p=0,1,2$
\begin{align}\lab{eq:waveeqwidetildepsi1:gtilde:Teu:pr_tphi}
&{\square}_{{\g}}\psi^{(p),\chi_{0}\prtphihat}_{s,ij} -\frac{4-2\de_{p0}}{\qs}\psi^{(p),\chi_{0}\prtphihat}_{s,ij}\nn\\
={}&\chi_{\tau_1, \tau_2}\big(\widehat{S}(\psi^{(p),\chi_{0}\prtphihat}_{s})_{ij} +(\widehat{Q}\psi^{(p),\chi_{0}\prtphihat}_{s})_{ij}\big)
+(1-\chi_{\tau_1, \tau_2})\big(\widehat{S}_K(\psi^{(p),\chi_{0}\prtphihat}_{s})_{ij} +(\widehat{Q}_K\psi^{(p),\chi_{0}\prtphihat}_{s})_{ij}\big)\nn\\
&
+(1-\chi_{\tau_1, \tau_2})f_p\psi^{(p),\chi_{0}\prtphihat}_{s}+\chi_{\tau_1, \tau_2}^{(1)}L^{(p),\chi_0\prtphihat}_{s,ij}+\chi_{0}\prtphihat\widehat{F}_{total,s,ij}^{(p)} +\widetilde{\widehat{F}_{total,s,ij}^{(p),\chi_{0}\prtphihat}},
\end{align}
where
\bsub
\lab{def:Nt+-2ijp:NwFtotal:pr_tphi:difference}
\bea
\phissh{ij}{p}{\chi_{0}\prtphihat}:=\chi_{0}\prtphihat(\phi_s^{(p)})_{ij}, \quad \psissh{ij}{p}{\chi_{0}\prtphihat}:=\chi_{0}\prtphihat(\psi_s^{(p)})_{ij}
\eea
and, with $\phi_{s,kl}$ denoting $(\phi_{s,kl}^{(p)})_{p=0,1,2}$ and $\psi_{s,kl}$ denoting $(\psi_{s,kl}^{(p)})_{p=0,1,2}$, 
\begin{align}
\widetilde{N_{W,s,ij}^{(p),\chi_{0}\prtphihat}}={}& \sum_{k,l}\Big(\mathbf{1}_{r\leq 12m}\pr^{\leq 2}(\Gac\,\phi_{s,kl})+ \mathbf{1}_{11m\leq r\leq 12m}O(1)\pr^{\leq 2}\phi_{s,kl}\Big),\\
\widetilde{\widehat{F}_{total,s,ij}^{(p),\chi_{0}\prtphihat}}={}& \sum_{k,l}\Big(\mathbf{1}_{r\leq 12m}\dk^{\leq 2}(\Gac\,\psiss{kl}{p})+\mathbf{1}_{r\leq 12m}\chi_{\tau_1, \tau_2}^{(1)}\dk^{\leq 2}(\Gac\,\phi_{s,kl})+\mathbf{1}_{11m\leq r\leq 12m}O(1)\pr^{\leq 2}\psiss{kl}{p}\nn\\
&+\mathbf{1}_{\tau\in[\tmic, \tau_1+1]\cup[\tau_2-3, \tau_2]}\mathbf{1}_{r\leq 12m}O(1)\pr^{\leq 1}\psi_{s,kl}\nn\\
&+\mathbf{1}_{\tau\in[\tau_1, \tau_1+1]\cup[\tau_2-3, \tau_2-2]}\mathbf{1}_{r\leq 12m}O(1)\Big(\phiss{kl}{0}, \phiss{kl}{1}\Big),\\
\widetilde{N_{T,s,ij}^{(p), \chi_{0}\prtphihat}}={}& \sum_{k,l}\Big(\mathbf{1}_{r\leq 12m}\pr^{\leq 1}(\Gac\,\phi_{s,kl}) + \mathbf{1}_{11m\leq r\leq 12m}O(1)\pr^{\leq 1}\phi_{s,kl}\Big).
\end{align}
\esub
In the above derivations, we have in particular used the fact that all the coefficients in the formulas \eqref{eq:linearterms:ScalarizedTeuSys:general:Kerrperturbation} for $L_{s,ij}^{(p)}$ are $(\tau,\tphi)$-independent. Additionally, we have used \eqref{eq:causlityrelationsforwidetildepsi1} for the structure of $\widetilde{\widehat{F}_{total,s,ij}^{(p),\chi_{0}\prtphihat}}$ on $\MM(\tau_1+1, \tau_2-3)$, and the fact that $\psi_{s,ij}^{(p)}$ corresponds to the scalarization of tensors $\breve{\pmb\phi}_s^{(p)}\in\sk_2(\mathbb{C})$ satisfying decoupled tensorial wave equations on Kerr in view of Steps 4 and 5 of Proposition \ref{prop:extensionprocedureoftheTeukolskywaveequations} for the structure of $\widetilde{\widehat{F}_{total,s,ij}^{(p),\chi_{0}\prtphihat}}$ on $\MM(\tau_2, +\infty)$.

Relying on Remark \ref{rem:prtauandprtphihatpreservesk2C}, we infer the existence, for $s=\pm 2$, $p=0,1,2$, of tensors $\pmb\phi_s^{(p),\chi_0\widehat{\pr}_{\tphi}}\in\sk_2(\mathbb{C})$ such that $\phissh{ij}{p}{\chi_0\widehat{\pr}_{\tphi}}=\pmb\phi_s^{(p),\chi_0\widehat{\pr}_{\tphi}}(\Om_i, \Om_j)$. In view of the Teukolsky wave/transport systems \eqref{eq:ScalarizedTeuSys:general:Kerrperturbation:pr_tphi}  \eqref{eq:ScalarizedQuantitiesinTeuSystem:Kerrperturbation:pr_tphi} for  $\phissh{ij}{p}{\chi_0\widehat{\pr}_{\tphi}}$, we deduce from Lemma \ref{lemma:formoffirstordertermsinscalarazationtensorialwaveeq} that $\pmb\phi_s^{(p),\chi_0\widehat{\pr}_{\tphi}}$ are solutions to the tensorial Teukolsky wave/transport systems \eqref{eq:TensorialTeuSysandlinearterms:rescaleRHScontaine2:general:Kerrperturbation}  \eqref{def:TensorialTeuScalars:wavesystem:Kerrperturbation} with $\N_{W}^{(p)}$ and $\N_{T}^{(p)}$ being replaced by $\N_{W}^{(p),\chi_0\widehat{\pr}_{\tphi}}$ and $\N_{T}^{(p),\chi_0\widehat{\pr}_{\tphi}}$ respectively, where $\N_{W}^{(p),\chi_0\widehat{\pr}_{\tphi}}, \N_{T}^{(p),\chi_0\widehat{\pr}_{\tphi}}\in\sk_2(\mathbb{C})$ are such that $N_{W,s,ij}^{(p),\chi_0\widehat{\pr}_{\tphi}}=\N_{W}^{(p),\chi_0\widehat{\pr}_{\tphi}}(\Om_i, \Om_j)$ and $N_{T,s,ij}^{(p),\chi_0\widehat{\pr}_{\tphi}}=\N_{T}^{(p),\chi_0\widehat{\pr}_{\tphi}}(\Om_i, \Om_j)$. Furthermore, $\psi^{(p),\chi_0\widehat{\pr}_{\tphi}}_{s,ij}$ satisfies the wave system \eqref{eq:waveeqwidetildepsi1:gtilde:Teu:pr_tphi}, and since $\psi^{(p),\chi_0\widehat{\pr}_{\tphi}}_{s,ij}=\chi_0\widehat{\pr}_{\tphi}\psiss{ij}{p}$ and $\phi^{(p),\chi_0\widehat{\pr}_{\tphi}}_{s,ij}=\chi_0\widehat{\pr}_{\tphi}\phiss{ij}{p}$, $\psi^{(p),\chi_0\widehat{\pr}_{\tphi}}_{s,ij}$ satisfies \eqref{eq:causlityrelationsforwidetildepsi1} with $(\phi^{(p)}_{s,ij}, \psi^{(p)}_{s,ij})$ being replaced by $(\phi^{(p),\chi_0\widehat{\pr}_{\tphi}}_{s,ij}, \psi^{(p),\chi_0\widehat{\pr}_{\tphi}}_{s,ij})$. We may thus apply Theorem \ref{thm:EMF:systemofTeuscalarized:order0:final} to infer, 
for $s=\pm 2$ and $\de\in(0, \frac{1}{3}]$,
\begin{align}
\lab{eq:EMFtotalp:sumup:rweightscontrolled:pm2:prtphi}
\EMFtotalh{s}{\nab_{\chi_{0}\prtphihat}^{\leq 1}} 
\les{}&\IE{\nab_{\chi_{0}\prtphihat}^{\leq 1}\pmb\phi_s, \nab_{\chi_{0}\prtphihat}^{\leq 1}\pmb\psi_s} +\A[\nab_{\chi_{0}\prtphihat}^{\leq 1}\pmb\psi_{s}](\Iti)+\A[\nab_{\chi_{0}\prtphihat}^{\leq 1}\pmb\phi_{s}](\tau_1,\tau_2)\nn\\
& 
+\NNttotalph{s}{\nab_{\chi_{0}\prtphihat}^{\leq 1}}+\sum_{p=0}^2\Errdefect[\nab_{\chi_0\prtphihat}^{\leq 1}\psis{p}]+\widetilde{\NNttotalh{s}{\chi_{0}\prtphihat}},
\end{align}
where $\widetilde{\NNttotalh{s}{\chi_{0}\prtphihat}}$ is given by
\begin{align}
\lab{NNtpm2total:prtphihat:diff:expression}
\widetilde{\NNttotalh{s}{\chi_{0}\prtphihat}}
={}&\sum_{p=0}^2\widetilde{\mathcal{N}}[\psish{p}{\chi_{0}\prtphihat},\widetilde{\widehat{F}_{total,s}^{(p),{\chi_{0}\prtphihat}}}]\nn\\
&+\sum_{p=0}^2\widehat{\mathcal{N}}'[\phish{p}{\chi_{0}\prtphihat}, \widetilde{\N_{W, s}^{(p),{\chi_{0}\prtphihat}}}](\tau_1, \tau_2)+\sum_{p=0}^1\widehat{\mathcal{N}}'[ \phish{p}{\chi_{0}\prtphihat}, r^{-2}\widetilde{\N_{T, s}^{(p),{\chi_{0}\prtphihat}}}](\tau_1, \tau_2)\nn\\
&+\sum_{p=0}^1\int_{\MM_{r\leq 12m}(\tau_1,\tau_2)}\Big|\pr^{\leq 1}\widetilde{\N_{T, s}^{(p),{\chi_{0}\prtphihat}}}\Big|^2\nn\\
&+\sum_{p=0}^1\int_{\MM(\tau_1,\tau_2)} \bigg(r^{-1+\de}\Big|\widetilde{\N_{W, s}^{(p),{\chi_{0}\prtphihat}}}\Big|+r^{-2+\de}\Big|\pr\widetilde{\N_{T, s}^{(p),{\chi_{0}\prtphihat}}}\Big|\bigg)\Big|\phish{p}{\chi_{0}\prtphihat}\Big|.
\end{align}

Next, we control the terms on the RHS of \eqref{NNtpm2total:prtphihat:diff:expression}. Notice that for $H$ and $\psi$ supported in $\MM_{r\leq 12m}(\Iti)$, we have the following analog of \eqref{esti:controlofNNtener:byprF}
\bea\lab{eq:controlRHSNNtpm2total:prtphihat:diff:expression:aux1}
\widetilde{\mathcal{N}}[\psi, H]  &\les& \bigg(\int_{\MM_{r\leq 12m}}|H|^2\bigg)^{\frac{1}{2}}\bigg(\EM_{r\leq 12m}[\psi](\Iti)\bigg)^{\frac{1}{2}} +\int_{\MM_{r\leq 12m}(\Iti)}|H|^2\nn\\
 &&+ \left(\int_{\Mtrap}|\pr^{\leq 1}H|^2\right)^{\frac{1}{2}}\left(\int_{\Mtrap(\Iti)} |\psi|^2\right)^{\frac{1}{2}},
\eea
for $H$ supported in supported in $\MM_{11m, 12m}(\Iti)$ and $\psi$ supported in $\MM_{r\leq 12m}(\Iti)$, we have\footnote{This estimate is significantly easier than \eqref{eq:controlRHSNNtpm2total:prtphihat:diff:expression:aux1} as $H$ vanishes identically on $\Mtrap$ in \eqref{eq:controlRHSNNtpm2total:prtphihat:diff:expression:aux1:biscasewhichiseasier}.}
\bea\lab{eq:controlRHSNNtpm2total:prtphihat:diff:expression:aux1:biscasewhichiseasier}
\widetilde{\mathcal{N}}[\psi, H]  \les \bigg(\int_{\MM_{11m, 12m}}|H|^2\bigg)^{\frac{1}{2}}\bigg(\EM_{11m,12m}[\psi](\Iti)\bigg)^{\frac{1}{2}} +\int_{\MM_{11m, 12m}(\Iti)}|H|^2,
\eea
and from \eqref{eq:defofNNhat'} that
\begin{align}\lab{eq:controlRHSNNtpm2total:prtphihat:diff:expression:aux2}
\widehat{\mathcal{N}}'[\psi, H](\tau_1, \tau_2)  \les&  \left(\int_{\Mntrap_{r\leq 12m}(\tau_1, \tau_2)}|\pr^{\leq 1}\psi|^2\right)^{\frac{1}{2}}\left(\int_{\MM_{r\leq 12m}(\tau_1, \tau_2)}|H|^2\right)^{\frac{1}{2}}+\int_{\MM_{r\leq 12m}(\tau_1, \tau_2)}|H|^2\nn\\
\les& \Big(\M_{r\leq 12m}[\psi](\tau_1, \tau_2)\Big)^{\frac{1}{2}}\left(\int_{\MM_{r\leq 12m}(\tau_1, \tau_2)}|H|^2\right)^{\frac{1}{2}}+\int_{\MM_{r\leq 12m}(\tau_1, \tau_2)}|H|^2.
\end{align}
Applying \eqref{eq:controlRHSNNtpm2total:prtphihat:diff:expression:aux1} and \eqref{eq:controlRHSNNtpm2total:prtphihat:diff:expression:aux1:biscasewhichiseasier}, we infer, for $\varphi_j$, $j=1,2,3$, and $\psi$ supported in $\MM_{r\leq 12m}(\Iti)$,
\begin{align}
\lab{eq:controlRHSNNtpm2total:prtphihat:diff:expression:aux2:result}
&\NNt\Big[\mathbf{1}_{11m\leq r\leq 12m}\pr^{\leq 2}\varphi_1+\mathbf{1}_{\tau\in[\tmic, \tau_1+1]\cup[\tau_2-3, \tau_2]}\pr^{\leq 1}\varphi_2
+\mathbf{1}_{\tau\in[\tau_1, \tau_1+1]\cup[\tau_2-3, \tau_2-2]}\pr^{\leq 1}\varphi_3, \psi\Big]\nn\\
\les& \,\bigg(\M_{{11m,12m}}[\pr^{\leq 1}\varphi_1](\Iti)+\sup_{\tau\in\Iti}\E_{r\leq 12m}[\varphi_2](\tau)+\sup_{\tau\in [\tau_1,\tau_2]}\E_{r\leq 12m}[\varphi_3](\tau)\bigg)^{\frac{1}{2}}\Big(\EM_{r\leq 12m}[\psi](\Iti)\Big)^{\frac{1}{2}}
\nn\\
&\,+\M_{{11m,12m}}[\pr^{\leq 1}\varphi_1](\Iti)+\sup_{\tau\in\Iti}\E_{r\leq 12m}[\varphi_2](\tau)+\sup_{\tau\in [\tau_1,\tau_2]}\E_{r\leq 12m}[\varphi_3](\tau).
\end{align}
Also, for $\varphi_4$ and $\psi$ supported in $\MM_{r\leq 12m}(\Iti)$, we have, using in particular Lemma \ref{lem:gpert:MMtrap}, 
\begin{align}\lab{eq:controlRHSNNtpm2total:prtphihat:diff:expression:aux3}
\NNt[\pr^{\leq 2}(\Gac\varphi_4),\psi]
\les \bigg(\ep^2\sup_{\tau\in\Iti}\E_{r\leq 12m}[\pr^{\leq 1}\varphi_4](\tau)\bigg)^{\frac{1}{2}}\Big(\EM_{r\leq 12m}[\psi](\Iti)\Big)^{\frac{1}{2}}.
\end{align}
 In view of \eqref{eq:controlRHSNNtpm2total:prtphihat:diff:expression:aux2}, \eqref{eq:controlRHSNNtpm2total:prtphihat:diff:expression:aux2:result},  \eqref{eq:controlRHSNNtpm2total:prtphihat:diff:expression:aux3}, \eqref{def:Nt+-2ijp:NwFtotal:pr_tphi:difference} and \eqref{NNtpm2total:prtphihat:diff:expression}, we infer 
\beaa
\widetilde{\NNttotalh{s}{\chi_{0}\prtphihat}} &\les& \left(\EM_{r\leq 12m}[\pmb\phi_s^{\chi_{0}\prtphihat}](\tau_1, \tau_2)+\EM_{r\leq 12m}[\pmb\psi_s^{\chi_{0}\prtphihat}](\Iti)\right)^{\frac{1}{2}}\\
&&\times\Bigg(\ep^2\sup_{\tau\in[\tau_1, \tau_2]}\E_{r\leq 12m}[\pr^{\leq 1}\pmb\phi_s](\tau)+\ep^2\sup_{\tau\in\Iti}\E_{r\leq 12m}[\pr^{\leq 1}\pmb\psi_s](\tau)\\
&&+\M_{11m,12m}[\pr^{\leq 1}\pmb\phi_s](\tau_1, \tau_2)+\M_{11m,12m}[\pr^{\leq 1}\pmb\psi_s](\Iti)\\
&&+\sup_{\tau\in[\tau_1, \tau_2]}\E_{r\leq 12m}[\pmb\phi_s](\tau)+\sup_{\tau\in\Iti}\E_{r\leq 12m}[\pmb\psi_s](\tau)\Bigg)^{\frac{1}{2}}\\
&&+\ep^2\sup_{\tau\in[\tau_1, \tau_2]}\E_{r\leq 12m}[\pr^{\leq 1}\pmb\phi_s](\tau)+\ep^2\sup_{\tau\in\Iti}\E_{r\leq 12m}[\pr^{\leq 1}\pmb\psi_s](\tau)\\
&&+\M_{11m,12m}[\pr^{\leq 1}\pmb\phi_s](\tau_1, \tau_2)+\M_{11m,12m}[\pr^{\leq 1}\pmb\psi_s](\Iti)\\
&&+\sup_{\tau\in[\tau_1, \tau_2]}\E_{r\leq 12m}[\pmb\phi_s](\tau)+\sup_{\tau\in\Iti}\E_{r\leq 12m}[\pmb\psi_{s}](\tau)\\
&&+\big(\EM_{r\leq 12m}[\pr^{\leq 1}\pmb\psi_s](\Iti)+\EM_{r\leq 12m}[\pr^{\leq 1}\pmb\psi_s](\tau_1, \tau_2)\big)^{\frac{1}{2}} \big(\Ao[\nab_{\chi_{0}\pr_{\tphi}}^{\leq 1}\pmb\psi_{s}](\Iti)\big)^{\frac{1}{2}}
\eeaa
and hence 
\begin{align*}
\widetilde{\NNttotalh{s}{\chi_{0}\prtphihat}} \les& \Big(\EMFtotalps{s}+\ep^2\EMFtotalhps{s}{\pr^{\leq 1}} +\EMF_{s,\de, \text{total}, 11m\leq r\leq 12m}[\pr^{\leq 1}\pmb\phi_{s}]\Big)^{\frac{1}{2}}\\
&\times\Big(\EMFtotalhps{s}{\nab_{\chi_{0}\prtphihat}^{\leq 1}}\Big)^{\frac{1}{2}}+\EMFtotalps{s}+\ep^2\EMFtotalhps{s}{\pr^{\leq 1}}\\
&+\big(\EMFtotalhps{s}{\pr^{\leq 1}}\big)^{\frac{1}{2}} \big(\Ao[\nab_{\chi_{0}\pr_{\tphi}}^{\leq 1}\pmb\psi_{s}](\Iti)\big)^{\frac{1}{2}} +\EMF_{s,\de, \text{total}, 11m\leq r\leq 12m}[\pr^{\leq 1}\pmb\phi_{s}].
\end{align*}
Plugging into \eqref{eq:EMFtotalp:sumup:rweightscontrolled:pm2:prtphi} to control the last term on the RHS of \eqref{eq:EMFtotalp:sumup:rweightscontrolled:pm2:prtphi}, we infer
\begin{align*}
&\EMFtotalh{s}{\nab_{\chi_{0}\prtphihat}^{\leq 1}}\\ 
\les{}&\IE{\nab_{\chi_{0}\prtphihat}^{\leq 1}\pmb\phi_s, \nab_{\chi_{0}\prtphihat}^{\leq 1}\pmb\psi_s} +\A[\nab_{\chi_{0}\prtphihat}^{\leq 1}\pmb\psi_{s}](\Iti)+\A[\nab_{\chi_{0}\prtphihat}^{\leq 1}\pmb\phi_{s}](\tau_1,\tau_2)\nn\\
& 
+\NNttotalph{s}{\nab_{\chi_{0}\prtphihat}^{\leq 1}}+\sum_{p=0}^2\Errdefect[\nab_{\chi_0\prtphihat}^{\leq 1}\psis{p}] +\EMF_{s,\de, \text{total}, 11m\leq r\leq 12m}[\pr^{\leq 1}\pmb\phi_{s}]\\
&+\Big(\EMFtotalps{s}+\ep^2\EMFtotalhps{s}{\pr^{\leq 1}} +\EMF_{s,\de, \text{total}, 11m\leq r\leq 12m}[\pr^{\leq 1}\pmb\phi_{s}]\Big)^{\frac{1}{2}}\\
&\times\Big(\EMFtotalhps{s}{\nab_{\chi_{0}\prtphihat}^{\leq 1}}\Big)^{\frac{1}{2}}+\EMFtotalps{s}+\ep^2\EMFtotalhps{s}{\pr^{\leq 1}}\\
&+\big(\EMFtotalhps{s}{\pr^{\leq 1}}\big)^{\frac{1}{2}} \big(\Ao[\nab_{\chi_{0}\pr_{\tphi}}^{\leq 1}\pmb\psi_{s}](\Iti)\big)^{\frac{1}{2}}
\end{align*}
and hence  
\begin{align*}
&\EMFtotalh{s}{\nab_{\chi_{0}\prtphihat}^{\leq 1}}\nn\\ 
\les{}& \IE{\pr^{\leq 1}\pmb\phi_s, \pr^{\leq 1}\pmb\psi_s} +\A[\nab_{\chi_{0}\prtphihat}^{\leq 1}\pmb\psi_{s}](\Iti)+\A[\nab_{\chi_{0}\prtphihat}^{\leq 1}\pmb\phi_{s}](\tau_1,\tau_2)\nn\\
& 
+\NNttotalph{s}{\nab_{\chi_{0}\prtphihat}^{\leq 1}}+\sum_{p=0}^2\Errdefect[\nab_{\chi_0\prtphihat}^{\leq 1}\psis{p}] +\EMF_{s,\de, \text{total}, 11m\leq r\leq 12m}[\pr^{\leq 1}\pmb\phi_{s}]\nn\\
& +\big(\EMFtotalhps{s}{\pr^{\leq 1}}\big)^{\frac{1}{2}} \big(\EMFtotalps{s}+\Ao[\nab_{\chi_{0}\pr_{\tphi}}^{\leq 1}\pmb\psi_{s}](\Iti)\big)^{\frac{1}{2}}+\ep^2\EMFtotalhps{s}{\pr^{\leq 1}}.
\end{align*}
Together with \eqref{eq:EMFtotalp:sumup:finalestimate:pm2:prttandprtphi:sofaronlyprtt}, this yields
\begin{align}\lab{eq:EMFtotalp:sumup:finalestimate:pm2:prttandprtphi} 
&\EMFtotalh{s}{\nab_{(\pr_{\tt}, \chi_{0}\prtphihat)}^{\leq 1}}\nn\\ 
\les{}& \IE{\pr^{\leq 1}\pmb\phi_s, \pr^{\leq 1}\pmb\psi_s} +\A[\nab_{(\pr_{\tt}, \chi_{0}\prtphihat)}^{\leq 1}\pmb\psi_{s}](\Iti)+\A[\nab_{(\pr_{\tt}, \chi_{0}\prtphihat)}^{\leq 1}\pmb\phi_{s}](\tau_1,\tau_2) +\ep\EMFtotalhps{s}{\pr^{\leq 1}}\nn\\
& 
+\NNttotalph{s}{{\nab_{\chi_{0}\prtphihat}^{\leq 1}}}+\sum_{p=0}^2\Errdefect[\nab_{(\pr_{\tt}, \chi_0\prtphihat)}^{\leq 1}\psis{p}] +\EMF_{s,\de, \text{total}, 11m\leq r\leq 12m}[\pr^{\leq 1}\pmb\phi_{s}]\nn\\
& +\big(\EMFtotalhps{s}{\pr^{\leq 1}}\big)^{\frac{1}{2}} \big(\EMFtotalps{s}+\Ao[\nab_{(\pr_{\tt}, \chi_{0}\pr_{\tphi})}^{\leq 1}\pmb\psi_{s}](\Iti)\big)^{\frac{1}{2}}.
\end{align}

We commute further with $(\pr_{\tt}, \chi_0\prtphihat)$, and by an inductive argument, we infer the following EMF estimates for $(\pr_\tau, \chi_{0}\prtphihat)^{\leq \reg}$-derivatives of $(\pmb\phi_s, \pmb\psi_s)$, with $\reg\leq 14$, 
 \begin{align}
\lab{eq:EMFtotalp:sumup:finalestimate:pm2:prttandprtphi:highorder}
&\EMFtotalh{s}{\nab_{(\pr_{\tt},\chi_{0}\prtphihat)}^{\leq \reg}} \nn\\
\les&\IE{\pr^{\leq  \reg}\pmb\phi_s, \pr^{\leq  \reg}\pmb\psi_s} +\A[\nab_{(\pr_{\tt},\chi_{0}\prtphihat)}^{\leq  \reg}\pmb\psi_{s}](\Iti)+\A[\nab_{(\pr_{\tt},\chi_{0}\prtphihat)}^{\leq  \reg}\pmb\phi_{s}](\tau_1,\tau_2)+\ep\EMFtotalhps{s}{\pr^{\leq  \reg}}\nn\\
& 
+\sum_{p=0}^2\Errdefect[\nab_{(\pr_{\tt},\chi_{0}\prtphihat)}^{\leq \reg}\psis{p}]+\NNttotalph{s}{\nab_{(\pr_{\tt},\chi_{0}\prtphihat)}^{\leq \reg}} +\EMF_{s,\de, \text{total}, 11m\leq r\leq 12m}[\pr^{\leq \reg}\pmb\phi_{s}]\nn\\
&+\big(\EMFtotalhps{s}{\pr^{\leq\reg}}\big)^{\frac{1}{2}} \big(\EMFtotalhps{s}{\pr^{\leq\reg-1}}+\Ao[\nab_{(\pr_{\tt}, \chi_{0}\pr_{\tphi})}^{\leq\reg}\pmb\psi_{s}](\Iti)\big)^{\frac{1}{2}}.
 \end{align}


\subsubsection{Proof of Proposition \ref{prop:EMFtotalp:finalestimate:pm2:prhighorder}}
\lab{subsubsect:commutationwithprttprtphi:conditionalEMF:high:v1}


Next, we recover the control of all the $\pr^{\leq \reg}$ derivatives from the above control of $(\pr_{\tt},\chi_{0}\prtphihat)^{\leq \reg}$ derivatives. First, in view of \eqref{eq:higherorderenergymorawetzforlargerconditionalonsprtau:Teu} and \eqref{eq:EMFtotalp:sumup:finalestimate:pm2:prtt}, we have, for $\reg\leq 14$,
\bea\lab{eq::higherorderenergymorawetzforlargerconditionalonsprtau:Teu:plusEMFestimatehigherorderprtt}
&&\EMF_{s, \de, \text{total}, r\geq 11m}[\pr^{\leq \reg}\pmb\phi_{s}] \nn\\
&\les&  \IE{\pr^{\leq \reg}\pmb\phi_s, \pr^{\leq \reg}\pmb\psi_s} +\NNttotalph{s}{\pr^{\leq \reg}}+\ep \EMFtotalhps{s}{\pr^{\leq \reg}}\nn\\
&& +\sum_{p=0}^2\Errdefect[\nab_{\pr_{\tau}}^{\leq \reg}\psis{p}]+\ov{\A}[\nab_{\pr_{\tt}}^{\leq \reg}\pmb\psi_{s}](\Iti)+\ov{\A}[\nab_{\pr_{\tt}}^{\leq \reg}\pmb\phi_{s}](\tau_1,\tau_2)\nn\\
&& +\A[\pmb\psi_{s}](\Iti)+\A[\pmb\phi_{s}](\tau_1,\tau_2)+\sum_{p=0}^{2}\widehat{\mathcal{N}}_{r\geq 10m}''[\pr^{\leq \reg}\psis{p}, \pr^{\leq \reg}\widehat{\F}_{total,s}^{(p)}](\tau_2-3, \tau_2)\nn\\
&&
+\Big(\EMFtotalhps{s}{\pr^{\leq \reg-1}} +\ov{\A}[\nab_{\pr_{\tt}}^{\leq \reg}\pmb\psi_{s}](\Iti)\Big)^{\frac{1}{2}}\Big(\EMFtotalhps{s}{\pr^{\leq \reg}}\Big)^{\frac{1}{2}},
\eea
where we have used the following consequence of  \eqref{def:EMFtotalps:pm2} and \eqref{def:EMFtotalpsnotilde:pm2}
\beaa
\EMFtotalhps{s}{\nab_{\pr_{\tt}}^{\leq \reg}}\les \EMFtotalh{s}{\nab_{\pr_{\tt}}^{\leq \reg}},
\eeaa
and where we used the following estimate, for any $\tau'<\tau''$ and any $\pmb\varphi\in\sk_2(\mathbb{C})$,  
\bea
\A[\nab_{\pr_{\tt}}\pmb\varphi](\tau', \tau'')\les \ov{\A}[\nab_{\pr_{\tt}}\pmb\varphi](\tau', \tau'')+ \EF[\pmb\varphi](\tau', \tau'')
\eea
which follows from \eqref{def:AandAonorms:tau1tau2}. Also, in view of the estimates \eqref{eq:imporvedhigherorderenergymorawetzconditionalonlowerrweigth:Teu} and \eqref{eq:EMFtotalp:sumup:finalestimate:pm2:prttandprtphi:highorder}, we have, for $\ep$ small enough and for $\reg\leq 14$,
\beaa
&&\EMFtotalh{s}{\nab_{(\pr_{\tt},\chi_{0}\prtphihat)}^{\leq \reg}}+\EMFtotalhps{s}{\pr^{\leq \reg}}\nn\\
&\les& \IEde{\pr^{\leq \reg}\pmb\phi_s, \pr^{\leq \reg}\pmb\psi_s} +\Ao[\nab_{(\pr_{\tt},\chi_{0}\prtphihat)}^{\leq  \reg}\pmb\psi_{s}](\Iti)+\Ao[\nab_{(\pr_{\tt},\chi_{0}\prtphihat)}^{\leq  \reg}\pmb\phi_{s}](\tau_1,\tau_2)\nn\\
&& 
+\sum_{p=0}^2\Errdefect[\nab_{(\pr_{\tt},\chi_{0}\prtphihat)}^{\leq \reg}\psis{p}]+\NNttotalph{s}{\pr^{\leq \reg}} +\EMF_{s,\de, \text{total}, 11m\leq r\leq 12m}[\pr^{\leq\reg}\pmb\phi_{s}]\nn\\
&& +\A[\pmb\psi_{s}](\Iti)+\A[\pmb\phi_{s}](\tau_1,\tau_2)  +\sum_{p=0}^{2}\widehat{\mathcal{N}}''[\pr^{\leq \reg}\psis{p}, \pr^{\leq \reg}\widehat{\F}_{total,s}^{(p)}](\tau_2-3, \tau_2)\nn\\
&& +\EMFtotalhps{s}{\pr^{\leq \reg-1}},
\eeaa
where we used the following estimate, for any $\tau'<\tau''$ and any $\pmb\varphi\in\sk_2(\mathbb{C})$,  
\bea
\A[\nab_{(\pr_{\tt}\chi_{0}\prtphihat)}\pmb\varphi](\tau', \tau'')\les \ov{\A}[\nab_{(\pr_{\tt}\chi_{0}\prtphihat)}\pmb\varphi](\tau', \tau'')+ \EF[\pmb\varphi](\tau', \tau'')
\eea
which follows from \eqref{def:AandAonorms:tau1tau2}. Thus, using \eqref{eq::higherorderenergymorawetzforlargerconditionalonsprtau:Teu:plusEMFestimatehigherorderprtt} to control $\EMF_{s,\de, \text{total}, 11m\leq r\leq 12m}[\pr^{\leq\reg}\pmb\phi_{s}]$, we infer, for $\ep$ small enough,
\bea\lab{eq:EMFtotalp:sumup:finalestimate:pm2:prhighorder:notthereyetbutalmost} 
&&\EMFtotalh{s}{\nab_{(\pr_{\tt},\chi_{0}\prtphihat)}^{\leq \reg}}+\EMFtotalhps{s}{\pr^{\leq \reg}}\nn\\
&\les& \IEde{\pr^{\leq \reg}\pmb\phi_s, \pr^{\leq \reg}\pmb\psi_s} +\Ao[\nab_{(\pr_{\tt},\chi_{0}\prtphihat)}^{\leq  \reg}\pmb\psi_{s}](\Iti)+\Ao[\nab_{(\pr_{\tt},\chi_{0}\prtphihat)}^{\leq  \reg}\pmb\phi_{s}](\tau_1,\tau_2)\nn\\
&& +\sum_{p=0}^2\Errdefect[\nab_{(\pr_{\tt},\chi_{0}\prtphihat)}^{\leq \reg}\psis{p}]+\NNttotalph{s}{\pr^{\leq \reg}}  +\A[\pmb\psi_{s}](\Iti)+\A[\pmb\phi_{s}](\tau_1,\tau_2) \nn\\
&& +\sum_{p=0}^{2}\widehat{\mathcal{N}}''[\pr^{\leq \reg}\psis{p}, \pr^{\leq \reg}\widehat{\F}_{total,s}^{(p)}](\tau_2-3, \tau_2)+\EMFtotalhps{s}{\pr^{\leq \reg-1}},
\eea
where we have used the following trivial bound which follows from \eqref{def:EMFtotalps:pm2} and \eqref{def:EMFtotalpsnotilde:pm2}
\beaa
\EMFtotalhps{s}{\nab_{(\pr_{\tt},\chi_{0}\pr_{\tilde{\phi}})}^{\leq \reg}}\les \EMFtotalh{s}{\nab_{(\pr_{\tt},\chi_{0}\pr_{\tphi})}^{\leq \reg}}\les \EMFtotalh{s}{\nab_{(\pr_{\tt},\chi_{0}\prtphihat)}^{\leq \reg}}.
\eeaa
Starting from \eqref{eq:EMFtotalp:sumup:rweightscontrolled:pm2}, we argue by induction to remove the term $\EMFtotalhps{s}{\pr^{\leq \reg-1}}$ on the RHS of \eqref{eq:EMFtotalp:sumup:finalestimate:pm2:prhighorder:notthereyetbutalmost}. This yields, for $\reg\leq 14$,
\bea\lab{eq:EMFtotalp:sumup:finalestimate:pm2:prhighorder} 
&&\EMFtotalh{s}{\nab_{(\pr_{\tt},\chi_{0}\prtphihat)}^{\leq \reg}}+\EMFtotalhps{s}{\pr^{\leq \reg}}\nn\\
&\les& \IEde{\pr^{\leq \reg}\pmb\phi_s, \pr^{\leq \reg}\pmb\psi_s} +\Ao[\nab_{(\pr_{\tt},\chi_{0}\prtphihat)}^{\leq  \reg}\pmb\psi_{s}](\Iti)+\Ao[\nab_{(\pr_{\tt},\chi_{0}\prtphihat)}^{\leq  \reg}\pmb\phi_{s}](\tau_1,\tau_2)\nn\\
&& +\sum_{p=0}^2\Errdefect[\nab_{(\pr_{\tt},\chi_{0}\prtphihat)}^{\leq \reg}\psis{p}]+\NNttotalph{s}{\pr^{\leq \reg}}  +\A[\pmb\psi_{s}](\Iti)+\A[\pmb\phi_{s}](\tau_1,\tau_2) \nn\\
&& +\sum_{p=0}^{2}\widehat{\mathcal{N}}''[\pr^{\leq \reg}\psis{p}, \pr^{\leq \reg}\widehat{\F}_{total,s}^{(p)}](\tau_2-3, \tau_2).
\eea

In order to use the above estimate \eqref{eq:EMFtotalp:sumup:finalestimate:pm2:prhighorder}  to prove Proposition \ref{prop:EMFtotalp:finalestimate:pm2:prhighorder}, it remains to control the term $\Ao[\nab_{(\pr_{\tt},\chi_{0}\pr_{\tphi})}^{\leq  \reg}\pmb\psi_{s}](\Iti)+\Ao[\nab_{(\pr_{\tt},\chi_{0}\pr_{\tphi})}^{\leq  \reg}\pmb\phi_{s}](\tau_1,\tau_2)$. This is the focus of the following two lemmas.

\begin{lemma}\lab{lemma:controllowerordertermusingtrickprrrmrtrap}
For a scalar $\psi$ that vanishes in $\MM_{r_+(1+\dhor'), 11m}(-\infty, \tmic)$, we have 
\beaa
\int_{\Mtrap}|\nab_{\pr_\tau}\psi|^2 &\les& \Big(\M_{r\leq 11m}[\nab_{\pr_\tau}\psi](\Iti)\Big)^{\frac{1}{2}}\left(\widetilde{\M}[\psi]\right)^{\frac{1}{2}}, \\ 
\int_{\Mtrap}|\nab_{\pr_{\tphi}}\psi|^2 &\les& \Big(\M_{r\leq 11m}[\nab_{\pr_{\tphi}}\psi](\Iti)\Big)^{\frac{1}{2}}\left(\widetilde{\M}[\psi]\right)^{\frac{1}{2}}.
\eeaa
\end{lemma}

\begin{proof}
Straightforward consequence of Lemma 6.3 in \cite{MaSz24}.
\end{proof}

\begin{lemma}
For any $1\leq \reg\leq 14$, we have
\bea
\lab{eq:controlofAohighorderbyEMF:prderivatives}
&&\Ao[\nab_{(\pr_{\tt},\chi_{0}\pr_{\tilde{\phi}})}^{\leq \reg}\pmb\psi_{s}](\Iti)+\Ao[\nab_{(\pr_{\tt},\chi_{0}\pr_{\tilde{\phi}})}^{\leq \reg}\pmb\phi_{s}](\tau_1,\tau_2)\nn\\
&\les&\Big(\EMFtotalhps{s}{\nab_{(\pr_{\tt},\chi_{0}\pr_{\tilde{\phi}})}^{\leq \reg}}
\Big)^{\frac{1}{2}} \Big(\EMFtotalh{s}{\nab_{(\pr_{\tt},\chi_{0}\pr_{\tilde{\phi}})}^{\leq \reg-1}} \Big)^{\frac{1}{2}} .
\eea
\end{lemma}

\begin{proof}
Since $\psis{p}=\phis{p}$ in $\MM(\tau_1+1, \tau_2-3)$, we have
\beaa
&&\Ao[\nab_{(\pr_{\tt},\chi_{0}\pr_{\tilde{\phi}})}^{\leq \reg}\pmb\psi_{s}](\Iti)+\Ao[\nab_{(\pr_{\tt},\chi_{0}\pr_{\tilde{\phi}})}^{\leq \reg}\pmb\phi_{s}](\tau_1,\tau_2)\nn\\
&\les& \EMFtotalhps{s}{\nab_{(\pr_{\tt},\chi_{0}\pr_{\tilde{\phi}})}^{\leq \reg-1}} + \sum_{i,j=1}^3\sum_{p=0}^2 \int_{\Mtrap} |(\pr_{\tt},\chi_{0}\pr_{\tilde{\phi}})^{\reg}\psiss{ij}{p}|^2\nn\\
&& +\int_{\MM(\tau_1, \tau_1+1)\cup\MM(\tau_2-3, \tau_2)}r^{-3}|\nab_{(\pr_{\tt},\chi_{0}\pr_{\tilde{\phi}})}^{\leq \reg}\pmb\phi_{s}|^2\nn\\
&\les& \EMFtotalhps{s}{\nab_{(\pr_{\tt},\chi_{0}\pr_{\tilde{\phi}})}^{\leq \reg-1}}+ \sum_{i,j=1}^3\sum_{p=0}^2\int_{\Mtrap} |(\pr_{\tt},\chi_{0}\pr_{\tilde{\phi}})^{\reg}\psiss{ij}{p}|^2 
\eeaa
which, together with the following estimate
\beaa
&&\sum_{i,j=1}^3\sum_{p=0}^2\int_{\Mtrap} |(\pr_{\tt},\chi_{0}\pr_{\tilde{\phi}})^{\reg}\psiss{ij}{p}|^2 \nn\\
&\les&\sum_{i,j=1}^3\sum_{p=0}^2\Big(\M_{r\leq 11m}[ (\pr_{\tt},\chi_{0}\pr_{\tilde{\phi}})^{\reg}\psiss{ij}{p}](\Iti)\Big)^{\frac{1}{2}}\left(\widetilde{\M}[(\pr_{\tt},\chi_{0}\pr_{\tilde{\phi}})^{\reg-1}\psiss{ij}{p}]\right)^{\frac{1}{2}}\nn\\
&\les&\Big(\EMFtotalhps{s}{\nab_{(\pr_{\tt},\chi_{0}\pr_{\tilde{\phi}})}^{\leq \reg}}
\Big)^{\frac{1}{2}} \Big(\EMFtotalh{s}{\nab_{(\pr_{\tt},\chi_{0}\pr_{\tilde{\phi}})}^{\leq \reg-1}} \Big)^{\frac{1}{2}} 
\eeaa
that follows from Lemma \ref{lemma:controllowerordertermusingtrickprrrmrtrap}, \eqref{def:EMFtotalps:pm2} and \eqref{def:EMFtotalpsnotilde:pm2}, yields the desired estimate.
\end{proof}

Now, based on \eqref{eq:EMFtotalp:sumup:finalestimate:pm2:prhighorder} and \eqref{eq:controlofAohighorderbyEMF:prderivatives}, we can conclude the proof of Proposition \ref{prop:EMFtotalp:finalestimate:pm2:prhighorder}.

\begin{proof}[Proof of Proposition \ref{prop:EMFtotalp:finalestimate:pm2:prhighorder}]
Using the estimate \eqref{eq:controlofAohighorderbyEMF:prderivatives} to control the terms $\Ao[\nab_{(\pr_{\tt},\chi_{0}\pr_{\tphi})}^{\leq  \reg}\pmb\psi_{s}](\Iti)+\Ao[\nab_{(\pr_{\tt},\chi_{0}\pr_{\tphi})}^{\leq  \reg}\pmb\phi_{s}](\tau_1,\tau_2)$ appearing on the RHS of \eqref{eq:EMFtotalp:sumup:finalestimate:pm2:prhighorder}, we deduce
\beaa
&&\EMFtotalh{s}{\nab_{(\pr_{\tt},\chi_{0}\pr_{\tilde{\phi}})}^{\leq \reg}}+\EMFtotalhps{s}{\pr^{\leq \reg}} \nn\\
&\les& \IEde{\pr^{\leq \reg}\pmb\phi_s, \pr^{\leq \reg}\pmb\psi_s} +\sum_{p=0}^2\Errdefect[\nab_{(\pr_{\tt},\chi_{0}\prtphihat)}^{\leq \reg}\psis{p}]+\NNttotalph{s}{\pr^{\leq \reg}}
  \nn\\
&& + \Big(\EMFtotalhps{s}{\pr^{\leq \reg}}
\Big)^{\frac{1}{2}} \Big(\EMFtotalh{s}{\nab_{(\pr_{\tt},\chi_{0}\pr_{\tilde{\phi}})}^{\leq \reg-1}} \Big)^{\frac{1}{2}} \nn\\
&&+\A[\pmb\psi_{s}](\Iti)+\A[\pmb\phi_{s}](\tau_1,\tau_2)
+\sum_{p=0}^{2}\widehat{\mathcal{N}}''[\pr^{\leq \reg}\psis{p}, \pr^{\leq \reg}\widehat{\F}_{total,s}^{(p)}](\tau_2-3, \tau_2),
\eeaa
which then yields by induction the desired estimate \eqref{eq:EMFtotalp:finalestimate:pm2:prhighorder} and hence concludes the proof of Proposition \ref{prop:EMFtotalp:finalestimate:pm2:prhighorder}.
\end{proof}


\subsection{Control of the terms $\mathbf{IE}_{\de}$ and $\Errdefect$}
\lab{subsect:controlofIEtotalandErrdefect}


In this section, we control the terms $\mathbf{IE}_{\de}$ and $\Errdefect$ appearing on the RHS of \eqref{eq:EMFtotalp:finalestimate:pm2:prhighorder}.


\subsubsection{Control of the term $\mathbf{IE}_{\de}$}
\lab{subsubsect:controlIEtotal:highorder}


We begin with controlling the term $\IEde{\pr^{\leq  \reg}\pmb\phi_{s}, \pr^{\leq  \reg}\pmb\psi_{s}}$. In view of the estimate \eqref{eq:localenergyforpsiontau1minus1tau1plus1}, we have for $s=\pm2$, all $\reg\leq 14$ and $\de\in (0,\frac{1}{3}]$,
\bea
\lab{esti:initialIEnorms:tmictotau1+2:prleqreg}
&&\sum_{p=0}^2\EMF_\de[\pr^{\leq \reg}\pmb\psi^{(p)}_{s}](\tmic, \tau_1+2)\nn\\
&\les&\sum_{p=0}^2\bigg(\E[\pr^{\leq \reg}\phis{p}](\tau_1)+\int_{\MM(\tau_1, \tau_1+2)}r^{1+\de} |\pr^{\leq \reg}\N_{W,s}^{(p)}|^2\bigg).
\eea
Plugging this estimate into the definition of $\IEde{\pr^{\leq  \reg}\pmb\phi_s, \pr^{\leq  \reg}\pmb\psi_s}$ in \eqref{eq:definitionofinitialenergyofphiandpsi:IEtotalterm:delta}, we infer
\begin{align}
\lab{eq:localenergyforpsiontau1minus1tau1plus2:initialEMF}
\IEde{\pr^{\leq  \reg}\pmb\phi_{s}, \pr^{\leq  \reg}\pmb\psi_{s}}\les\sum_{p=0}^2\bigg(\E[\pr^{\leq \reg}\phis{p}](\tau_1)+\int_{\MM(\tau_1, \tau_1+2)}r^{1+\de} |\pr^{\leq \reg}\N_{W,s}^{(p)}|^2\bigg).
\end{align}


\subsubsection{Control of the term $\Errdefect$}
\lab{subsubsect:controlErrdefect:highorder}


Next, we consider the second term $\sum\limits_{p=0}^2\Errdefect[\nab_{(\pr_{\tt},\chi_{0}\prtphihat)}^{\leq \reg}\psis{p}]$ on the RHS of  the high order EMF estimate \eqref{eq:EMFtotalp:finalestimate:pm2:prhighorder}. In view of the definition \eqref{def:Errdefectofpsi} of $\Errdefect$, we have for $s=\pm2$, all $\reg\leq 14$ and $\de\in (0,\frac{1}{3}]$,
\bea
\lab{esti:Errdefectterm:highorderprderi}
&&\sum_{p=0}^2\Errdefect[\nab_{(\pr_{\tt},\chi_{0}\prtphihat)}^{\leq \reg}\psis{p}]\nn\\
&\les&\sup_{\tau\in[\tmic,\tau_1+1]\cup[\tau_2-2,\tau_2]}\sum_{p=0}^2\sum_{j}\Big(\E[(\pr_{\tau}, \chi_0\pr_{\tphi})^{\leq \reg}(x^i\psiss{ij}{p})](\tau)+\E[(\pr_{\tau}, \chi_0\pr_{\tphi})^{\leq \reg}(x^i\psiss{ji}{p})](\tau)\Big)\nn\\
&\les&\sup_{\tau\in[\tmic, \tau_1+1]}\sum_{p=0}^2\E[\pr^{\leq \reg}\pmb\psi^{(p)}_{s}](\tau)+ \sum_{p=0}^2\EMF[\pr^{\leq \reg}\err_{\textrm{TDefect}}[\psi^{(p)}_{s}]](\tau_2-2, \tau_2)\nn\\
&\les&\sum_{p=0}^2\Big(\E[\pr^{\leq \reg}\phis{p}](\tau_1)+\ep^2\E[\pr^{\leq \reg}\phis{p}](\tau_2-3)\Big)\nn\\
&&+\sum_{p=0}^2\int_{\MM(\tau_1, \tau_1+2)}r^{1+\de} |\pr^{\leq \reg}\N_{W,s}^{(p)}|^2
+\ep^2\sum_{p=0}^2\int_{\MM(\tau_2-3, \tau_2-2)}r^{1+\de} |\pr^{\leq \reg}\N_{W,s}^{(p)}|^2,
\eea
where in the second step we have used Definition \ref{eq:definitionofthenotationerrforthescalarizationdefect} for  $\err_{\textrm{TDefect}}$ and where in the last step we have used the estimates
\eqref{esti:initialIEnorms:tmictotau1+2:prleqreg} and \eqref{eq:propertyofscalarizationdefect:prop2}.


\subsection{Control of the term $\NNt_{s, \de, \text{total}}$}
\lab{subsect:controlofNNtsdetotal}


In this section, we estimate the term $\NNttotalph{s}{\pr^{\leq \reg}}$ which appears on the RHS of  the high order EMF estimates \eqref{eq:EMFtotalp:finalestimate:pm2:prhighorder}. Our main estimates are contained in the following proposition.

\begin{proposition}[Control of $\NNttotalph{s}{\pr^{\leq \reg}}$]
\lab{prop:controlofhighorderNNtspin-2:lastcontrol:pfofthm7.5}
Under the assumptions of Theorem \ref{th:main:intermediary}, we have, for any $0<\la \leq 1$, all $\reg\leq 14$ and $\de\in (0,\frac{1}{3}]$, 
\bea
\lab{eq:controlofhighorderNNtspin-2:lastcontrol:pfofthm7.5}
&&\NNttotalph{s}{\pr^{\leq \reg}}\nn\\
&\les&
(\la+\ep)\EMFtotalhps{s}{\pr^{\leq \reg}} +\la^{-3}\sum_{p=0}^2\E[\pr^{\leq \reg}\phis{p}](\tau_1)
\nn\\
&&+\la^{-3}\sum_{p=0}^2{\mathcal{N}}_{\de}[\pr^{\leq \reg}\phis{p}, \pr^{\leq \reg}\N_{W, s}^{(p)}](\tau_1, \tau_2)+\la^{-1}\sum_{p=0}^1\int_{\MM(\tau_1,\tau_2)}r^{-3+\de}\big|\pr^{\leq \reg+1}\N_{T,s}^{(p)}\big|^2
\nn\\
&&
+\sum_{p=0}^1\int_{\MM_{r\geq 10m}(\tau_1,\tau_2)} \Big(r^{-1+\de}\big|\pr^{\leq \reg}\N^{(p)}_{W,s}\big|+r^{-2+\de}\big|\pr^{\leq \reg+1}\N^{(p)}_{T,s}\big|\Big)\big|\pr^{\leq \reg}\phis{p}\big|.
\eea
\end{proposition}

The rest of this section is dedicated to the proof of Proposition \ref{prop:controlofhighorderNNtspin-2:lastcontrol:pfofthm7.5}. In view of \eqref{expression:NNttotalps:pm2}, we have 
\begin{align}\lab{expression:NNttotalps:pm2:pasteinproof:s=pm2}
\NNttotalph{s}{\pr^{\leq \reg}}
=&\sum_{p=0}^2\widetilde{\mathcal{N}}[\pr^{\leq \reg}\psis{p},\pr^{\leq \reg}\widehat{\F}_{total,s}^{(p)}]+\sum_{p=0}^2\widehat{\mathcal{N}}'[\pr^{\leq \reg}\phis{p}, \pr^{\leq \reg}\N_{W, s}^{(p)}](\tau_1, \tau_2)\nn\\
&+\sum_{p=0}^1\widehat{\mathcal{N}}'[ \pr^{\leq \reg}\phis{p}, r^{-2}\pr^{\leq \reg}\N_{T, s}^{(p)}](\tau_1, \tau_2)+\sum_{p=0}^1\int_{\MM_{r\leq 12m}(\tau_1,\tau_2)}\big|\pr^{\leq \reg+1}\N_{T, s}^{(p)}\big|^2\nn\\
&+\sum_{p=0}^1\int_{\MM(\tau_1,\tau_2)} \Big(r^{-1+\de}\big|\pr^{\leq \reg}\N^{(p)}_{W,s}\big|+r^{-2+\de}\big|\pr\pr^{\leq \reg}\N^{(p)}_{T,s}\big|\Big)\big|\pr^{\leq \reg}\phis{p}\big|.
\end{align}
We start with estimating the term $\sum_{p=0,1,2}\widetilde{\mathcal{N}}[\pr^{\leq \reg}\psis{p},\pr^{\leq \reg}\widehat{\F}_{total,s}^{(p)}]$ appearing first on the RHS of \eqref{expression:NNttotalps:pm2:pasteinproof:s=pm2}.  Using the formula \eqref{eq:definitionofwidehatFtotalsij}, we deduce, for $p=0,1,2$,
\bea
\lab{decomp:NNtofhatFtotal-2:highorderpr}
&&\sum\limits_{p=0}^2\NNt[\pr^{\leq \reg}\psis{p},\pr^{\leq \reg}\widehat{\F}_{total,s}^{(p)}]\nn\\
&\les& 
\sum\limits_{p=0}^2\NNt[\pr^{\leq \reg}\psis{p},\pr^{\leq \reg}(\chi_{\tau_1, \tau_2}^{(1)}\N^{(p)}_{W,s}+\underline{\F}^{(p)}_{s})]
+\sum\limits_{p=0}^2\NNt[\pr^{\leq \reg}\psis{p},\pr^{\leq \reg}\breve{\F}^{(p)}_{s}].
\eea


\subsubsection{Comparison of $\NNt$ with $\widehat\NN$}


We have the following comparison of $\widetilde{\mathcal{N}}[\pmb \psi, \pmb F]$ with $\widehat{\mathcal{N}}[\pmb \psi, \pmb F]$.

\begin{lemma}\lab{lemma:additivepropertyofwidetildeN}
Given $N\geq 2$, consider any partition of $\Iti$ in intervals of the following form
\beaa
\Iti=(\tmic, \tau^{(1)}]\cup\bigcup_{j=1}^{N-1}( \tau^{(j)},  \tau^{(j+1)}]\cup( \tau^{(N)}, +\infty), \qquad \tmic+1< \tau^{(1)}<\cdots< \tau^{(N)}<+\infty.
\eeaa
Then, for families of scalars $\psi_{ij}, F_{ij}$ such that $\psi_{ij}$ and $F_{ij}$ vanish in $\Mtrap(-\infty, \tmic)$, $\widetilde{\mathcal{N}}[\pmb \psi, \pmb F]$, defined as in \eqref{def:NNtintermsofNNtMora:NNtEner:NNtaux:wavesystem:EMF}-\eqref{def:NNtMora:NNtEner:NNtaux:wavesystem:EMF}, satisfies the following bound, for any $0<\la\leq 1$,
\bea
\lab{eq:boundofNNtbyNNhat}
\widetilde{\mathcal{N}}[\pmb \psi, \pmb F] &\les_N& \la^{-1}\widehat{\mathcal{N}}[\pmb \psi, \pmb F](\tmic, \tau^{(1)}+1)+\la^{-1}\sum_{j=1}^{N-1}\widehat{\mathcal{N}}[\pmb \psi, \pmb F](\tau^{(j)}-1, \tau^{(j+1)}+1)\nn\\
&&+\la^{-1}\widehat{\mathcal{N}}[\pmb \psi, \pmb F](\tau^{(N)}-1, +\infty) +\la\EM[\pmb\psi](\Iti)\nn\\ 
&&+\left(\int_{\Mtrap(\Iti)}|\pmb F|^2\right)^{\frac{1}{2}}\left(\int_{\Mtrap(\Iti)}|\pmb \psi|^2\right)^{\frac{1}{2}},
\eea
where $\widehat{\mathcal{N}}[\pmb \psi, \pmb F](\tau', \tau'')$, for any $\tau'<\tau''$, is given as in \eqref{eq:defintionwidehatmathcalNfpsinormRHS}.
\end{lemma}

\begin{proof}
In view of the formulas \eqref{def:NNtMora:NNtEner:NNtaux:wavesystem:EMF} for $\NNtmora[\pmb{\psi},\F]$, $\NNtener[\pmb \psi, \pmb F]$ and $\NNtaux[\pmb F]$, and the assumption that $\psi_{ij}$ and $F_{ij}$ vanish in $\Mtrap(-\infty, \tmic)$, we have, for any $0<\lambda \leq 1$,
\beaa
\NNtaux[\pmb F]&\les& \int_{\MM(\Iti)} |\pmb F|^2,\\
\NNtmora[\pmb{\psi},\F]
&\les&\sum_{i,j}\bigg|\int_{\Mntrap(\Iti)}\Re\big({\pr_\tau}\psi_{ij}\ov{ F_{ij}}\big)\bigg|+\int_{\Mntrap(\Iti)}r^{-1}|\dk^{\leq 1}\pmb\psi| |\pmb F|\nn\\
&&
+ \sum_{i,j}\int_{\Mtrap(\Iti)}|F_{ij}||P^1\psi_{ij}|+\la^{-1}\int_{\Mntrap_{r_+(1+\dhor'),\Rmic}}|\pmb F|^2\nn\\
&&
+\la\M_{r_+(1+\dhor'),\Rmic}[\pmb\psi](\Iti)
\eeaa
and
\beaa
\NNtener[\pmb \psi, \pmb F]&\les&\sum_{i,j}\sup_{\tau\geq\tmic}\bigg|\int_{\Mntrap(\tmic, \tau)}{\Re\Big(F_{ij}\ov{\pr_{\tau}\psi_{ij}}\Big)}\bigg|+\sum_{n=-1}^{\iota}\sum_{i,j}\sup_{\tt\in\Reals}\NNt_{n,i,j,\tau}[\pmb\psi, \pmb F],
\eeaa
where $P^1$ is a PDO satisfying $P^1\in\Opw(\widetilde{S}^{1,0}(\MM))$, and where 
\bea\lab{eq:defofquantityNNtnijtaupfpmbpsipmbF:usefullater}
\NNt_{n,i,j,\tau}[\pmb\psi, \pmb F]:=\bigg|\int_{\Mtrap(-\infty,\tau)}\Re\Big(\ov{|q|^{-2}\Opw(\Theta_n)(\qs F_{ij})}V_n\Opw(\Theta_n)\psi_{ij}\Big)\bigg|
\eea
with the vectorfields $V_n$ and the symbols $\Theta_n\in\widetilde{S}^{0,0}(\MM)$, $n=-1,0,1,\ldots, \iota$, given as in Section \ref{sec:relevantmixedsymbolsonMM}. By summing up these estimates, we infer in view of \eqref{def:NNtintermsofNNtMora:NNtEner:NNtaux:wavesystem:EMF}, for any $0<\lambda \leq 1$,
\bea
\lab{eq:preesti:NNtnorm:proof}
\NNt[\pmb\psi, \pmb F]&\les&\sup_{\tau\geq\tmic}\sum_{i,j}\bigg|\int_{\Mntrap(\tmic, \tau)}{\Re\big(F_{ij}\ov{\pr_{\tau}\psi_{ij}}\big)}\bigg|+\int_{\Mntrap(\Iti)}r^{-1}|\dk^{\leq 1}\pmb\psi| |\pmb F|\nn\\
&&
+\la^{-1}\int_{\MM(\Iti)}|\pmb F|^2+\la\M[\pmb\psi](\Iti) + \sum_{i,j}\int_{\Mtrap(\Iti)}|F_{ij}||P^1\psi_{ij}|\nn\\
&&+\sum_{n=-1}^{\iota}\sum_{i,j}\sup_{\tt\in\Reals}\NNt_{n,i,j,\tau}[\pmb\psi, \pmb F].
\eea

Next, we estimate the last term $\NNt_{n,i,j,\tau}[\pmb\psi, \pmb F]$ on the RHS of \eqref{eq:preesti:NNtnorm:proof}. By introducing the smooth cut-off functions $\chi_{\tau,j}=\chi_{\tau,j}(\tau)$, $j=0,1$,  given as in \eqref{def:variouscutofffunctionsfortau':tautau'case}, we infer
\begin{align*}
\NNt_{n,i,j,\tau}[\pmb\psi, \pmb F]\les&\bigg|\int_{\Mtrap}\chi_{\tau,0}\Re\Big(\ov{|q|^{-2}\Opw(\Theta_n)(\qs F_{ij})}V_n\Opw(\Theta_n)\psi_{ij}\Big)\bigg|\nn\\
&+\int_{\Mtrap}\chi_{\tau,1}|\Opw(\widetilde{S}^{0,0}(\MM))F_{ij}||\Opw(\widetilde{S}^{1,0}(\MM))\psi_{ij}|\nn\\
\les&\int_{\Mtrap}|F_{ij}| |P^1\psi_{ij}|+\bigg(\int_{\Mtrap(\Iti)}|\pmb F|^2\bigg)^{\frac{1}{2}}\left(\int_{\Mtrap}\chi_{\tau,1}^2|\Opw(\widetilde{S}^{1,0}(\MM))\psi_{ij}|^2\right)^{\frac{1}{2}}.
\end{align*}
with $P^1\in\Opw(\widetilde{S}^{1,0}(\MM))$, where in the last step we have used the self-adjointness of $\Opw(\Theta_n)$. Now, in view of the support properties of $\chi_{\tau,1}$ in \eqref{def:variouscutofffunctionsfortau':tautau'case}, we have
\beaa
&&\int_{\Mtrap}\chi_{\tau,1}^2|\Opw(\widetilde{S}^{1,0}(\MM))\psi_{ij}|^2\\ 
&\les& \int_{\Mtrap}|\Opw(\widetilde{S}^{0,0}(\MM))\psi_{ij}|^2+\int_{\Mtrap}|\Opw(\widetilde{S}^{1,0}(\MM))\chi_{\tau,1}\psi_{ij}|^2\\
&\les& \int_{\Mtrap}|\psi_{ij}|^2+\int_{\Mtrap}|\pr^{\leq 1}(\chi_{\tau,1}\psi_{ij})|^2\les\EM_{\trap}[\pmb \psi](\Iti)
\eeaa
and hence
 \beaa
\NNt_{n,i,j,\tau}[\pmb\psi, \pmb F] &\les&\int_{\Mtrap}|F_{ij}| |P^1\psi_{ij}|+\bigg(\int_{\Mtrap(\Iti)}|\pmb F|^2\bigg)^{\frac{1}{2}} \Big(\EM_{\trap}[\pmb \psi](\Iti)\Big)^{\frac{1}{2}}\nn\\
&\les&\int_{\Mtrap(\Iti)}|F_{ij}| |P^1\psi_{ij}|+\la^{-1}\int_{\MM(\Iti)}|\pmb F|^2+ \la \EM_{\trap}[\pmb \psi](\Iti).
\eeaa
Plugging this estimate back into \eqref{eq:preesti:NNtnorm:proof}, we deduce
\bea
\lab{eq:preesti:NNtnorm:proof:v2}
\NNt[\pmb\psi, \pmb F]&\les&\sup_{\tau\geq\tmic}\sum_{i,j}\bigg|\int_{\Mntrap(\tmic, \tau)}{\Re\big(F_{ij}\ov{\pr_{\tau}\psi_{ij}}\big)}\bigg|+\int_{\Mntrap(\Iti)}r^{-1}|\dk^{\leq 1}\pmb\psi| |\pmb F|\nn\\
&&
+\la^{-1}\int_{\MM(\Iti)}|\pmb F|^2+\la\EM[\pmb\psi](\Iti) + \sum_{i,j}\int_{\Mtrap(\Iti)}|F_{ij}||P^1\psi_{ij}|,
\eea
where $P^1$ denotes a general PDO satisfying $P^1\in\Opw(\widetilde{S}^{1,0}(\MM))$.

 Let us denote the intervals $I_j$ by
\beaa
I_0:=(\tmic, \tau^{(1)}], \qquad I_j:=(\tau^{(j)}, \tau^{(j+1)}], \quad j=1,\cdots,N-1, \qquad I_N:=(\tau^{(N)}, +\infty)
\eeaa
and define the intervals $I_j^1$ by
\begin{align}
\lab{def:intervals:Ij1's}
I_0^1:=(\tmic, \tau^{(1)}+1), \,\,\,\, I_j^1:=(\tau^{(j)}-1, \tau^{(j+1)}+1), \,\, j=1,\cdots,N-1, \,\,\,\, I_N^1:=(\tau^{(N)}-1, +\infty).
\end{align}
Since $\cup_{j=0}^N I_j^1 = \cup_{j=0}^NI_j=\Iti$, we have, from the estimate \eqref{eq:preesti:NNtnorm:proof:v2} and the definition \eqref{eq:defintionwidehatmathcalNfpsinormRHS} of $\widehat{\NN}[\pmb{\psi}, \F](\tau_1,\tau_2)$,
\bea
\lab{eq:preesti:NNtnorm:proof:v3}
\NNt[\pmb\psi, \pmb F]\les\la^{-1} \sum_{j=0}^N\widehat{\mathcal{N}}[\pmb\psi, \pmb F](I_j^1)+\la\EM[\pmb\psi](\Iti) + \sum_{j=0}^N\sum_{k,l}\int_{\Mtrap(I_j)}|F_{kl}||P^1\psi_{kl}|.
\eea
We then follow the proof of Lemma 6.2 in \cite{MaSz24} to control the last term in \eqref{eq:preesti:NNtnorm:proof:v3}.
In view of the before to last estimate of the proof of Lemma 6.2 in \cite{MaSz24}, we have for any scalars $\psi$ and $F$, any $P^1\in \Opw(\widetilde{S}^{1,0}(\MM))$ and any $0<\la\leq 1$, assuming in addition that $\psi$ and $F$ vanish\footnote{The fact that $\psi$ and $F$ vanish in $\Mtrap(-\infty, \tmic)$ allows in particular to choose $I_0^1:=(\tmic, \tau^{(1)}+1)$ instead of $I_0^1:=(\tmic-\de, \tau^{(1)}+1)$ for some $\de>0$.} in $\Mtrap(-\infty, \tmic)$, 
\begin{align}\lab{eq:beforetolastestimateofproofLemma6dot2inMaSz24}
\int_{\Mtrap(I_j)}|F||P^1\psi| \les{}& \la^{-1}\widehat{\mathcal{N}}[\psi, F](I_j^1)+\la\EM[\psi](I^1_j)\nn\\
&+\bigg(\int_{\Mtrap(\Iti)}|F|^2\bigg)^{\frac{1}{2}}\bigg(\int_{\Mtrap(\Iti)}|\psi|^2\bigg)^{\frac{1}{2}}, \,\, j=0,\cdots,N.
\end{align}
 We now apply \eqref{eq:beforetolastestimateofproofLemma6dot2inMaSz24} with $(\psi=\psi_{kl}, F=F_{kl})$ and, together with the estimate \eqref{eq:preesti:NNtnorm:proof:v3}, we deduce, for any $0<\la\leq 1$,
\begin{align*}
\widetilde{\mathcal{N}}[\pmb\psi, \F]
{}&\les\la\EM[\pmb\psi](\Iti) +\la^{-1}\sum_{j=0}^N \widehat{\mathcal{N}}[\pmb \psi, \pmb F](I_j^1)
+ \sum_{j=0}^N\sum_{k,l}\int_{\Mtrap(I_j)}|F_{kl}||P^1\psi_{kl}|\\
{}&\les_N \sum_{j=0}^N\big( \la^{-1}\widehat{\mathcal{N}}[\pmb\psi, \F](I_j^1)+\la\EM[\pmb\psi](I^1_j)\big)+\left(\int_{\Mtrap(\Iti)}|\F|^2\right)^{\frac{1}{2}}\left(\int_{\Mtrap(\Iti)}|\pmb\psi|^2\right)^{\frac{1}{2}}\\
{}&\les_N \la^{-1}\sum_{j=0}^N\widehat{\mathcal{N}}[\pmb\psi, \F](I_j^1)+\la\EM[\pmb\psi](\Iti)+\left(\int_{\Mtrap(\Iti)}|\F|^2\right)^{\frac{1}{2}}\left(\int_{\Mtrap(\Iti)}|\pmb\psi|^2\right)^{\frac{1}{2}}
\end{align*}
which, in view of the above definition \eqref{def:intervals:Ij1's} of $\{I_j^1\}_{j=0,1,\ldots, N}$, is the stated estimate \eqref{eq:boundofNNtbyNNhat}. This concludes the proof of Lemma \ref{lemma:additivepropertyofwidetildeN}.
\end{proof}


\subsubsection{Control of the first term on the RHS of \eqref{decomp:NNtofhatFtotal-2:highorderpr}}


In order to estimate the first term on the RHS of \eqref{expression:NNttotalps:pm2:pasteinproof:s=pm2}, it suffices to control the RHS of \eqref{decomp:NNtofhatFtotal-2:highorderpr}, and the goal of this section is to control the first term on the RHS of \eqref{decomp:NNtofhatFtotal-2:highorderpr}. To this end, we use the comparison of $\NNt[\c,\c]$ with $\widehat{\NN}[\c,\c]$ provided by Lemma \ref{lemma:additivepropertyofwidetildeN} with $N=3$, $\tau^{(1)}=\tau_1+2$, $\tau^{(2)}=\tau_2-4$, and $\tau^{(3)}=\tau_2+1$ to obtain\footnote{Note that $\tau_2-4>\tau_1+2$ in view of the assumption \eqref{eq:tau2geqtau1plus10conditionisthemaincase}.}
\beaa
&&\widetilde{\mathcal{N}}[\pr^{\leq \reg}\psis{p},\pr^{\leq \reg}(\chi_{\tau_1, \tau_2}^{(1)}\N^{(p)}_{W,s}+\underline{\F}^{(p)}_{s})]\nn\\
&\les&\la^{-1}\widehat{\mathcal{N}}[\pr^{\leq \reg}\psis{p},\pr^{\leq \reg}(\chi_{\tau_1, \tau_2}^{(1)}\N^{(p)}_{W,s}+\underline{\F}^{(p)}_{s})](\tmic, \tau_1+3)\nn\\
&&+\la^{-1}\widehat{\mathcal{N}}[\pr^{\leq \reg}\psis{p},\pr^{\leq \reg}(\chi_{\tau_1, \tau_2}^{(1)}\N^{(p)}_{W,s}+\underline{\F}^{(p)}_{s})](\tau_1+1, \tau_2-3)\nn\\
&&+\la^{-1}\widehat{\mathcal{N}}[\pr^{\leq \reg}\psis{p},\pr^{\leq \reg}(\chi_{\tau_1, \tau_2}^{(1)}\N^{(p)}_{W,s}+\underline{\F}^{(p)}_{s})](\tau_2-5, \tau_2+2)\nn\\ 
&&+\la^{-1}\widehat{\mathcal{N}}[\pr^{\leq \reg}\psis{p},\pr^{\leq \reg}(\chi_{\tau_1, \tau_2}^{(1)}\N^{(p)}_{W,s}+\underline{\F}^{(p)}_{s})](\tau_2, +\infty) +\la\EM[\pr^{\leq \reg}\psis{p}](\Iti)\nn\\ 
&&+\bigg(\int_{\Mtrap(\Iti)}|\pr^{\leq \reg}(\chi_{\tau_1, \tau_2}^{(1)}\N^{(p)}_{W,s}+\underline{\F}^{(p)}_{s})|^2\bigg)^{\frac{1}{2}}\bigg(\int_{\Mtrap(\Iti)}|\pr^{\leq \reg}\psis{p}|^2\bigg)^{\frac{1}{2}}.
\eeaa  
In view of  the support property of $\underline{F}^{(p)}_{s,ij}$ from \eqref{eq:structureoftildef0:support}, the definition of the cutoff function $\chi_{\tau_1,\tau_2}^{(1)}$ from \eqref{eq:propertieschi:thisonetodealwithRHSofTeukolsky} and the fact from \eqref{eq:causlityrelationsforwidetildepsi1} that $\pmb\phi_{s}^{(p)}=\pmb\psi_{s}^{(p)}$ on $\MM(\tau_1+1,\tau_2-3)$,
we deduce 
\bea
\lab{eq:NNtofhighorderpr:spin-2:NwandFbar}
&&\sum_{p=0}^2\widetilde{\mathcal{N}}[\pr^{\leq \reg}\psis{p},\pr^{\leq \reg}(\chi_{\tau_1, \tau_2}^{(1)}\N^{(p)}_{W,s}+\underline{\F}^{(p)}_{s})]\nn\\
&\les&\sum_{p=0}^2\bigg[\la^{-1}\widehat{\mathcal{N}}[\pr^{\leq \reg}\psis{p},\pr^{\leq \reg}(\chi_{\tau_1, \tau_2}^{(1)}\N^{(p)}_{W,s})](\tau_1, \tau_1+3)\nn\\
&&+\la^{-1}\widehat{\mathcal{N}}[\pr^{\leq \reg}\psis{p},\pr^{\leq \reg}\underline{\F}^{(p)}_{s}](\tmic, \tau_1)+\la^{-1}\widehat{\mathcal{N}}[\pr^{\leq \reg}\phis{p},\pr^{\leq \reg}\N^{(p)}_{W,s}](\tau_1+1, \tau_2-3)\nn\\
&&+\la^{-1}\widehat{\mathcal{N}}[\pr^{\leq \reg}\psis{p},\pr^{\leq \reg}(\chi_{\tau_1, \tau_2}^{(1)}\N^{(p)}_{W,s})](\tau_2-5, \tau_2-2)+\la\EM[\pr^{\leq \reg}\psis{p}](\Iti)\nn\\
&&+\bigg(\int_{\Mtrap(\tau_1,\tau_2-2)}\big|\pr^{\leq \reg}(\chi_{\tau_1, \tau_2}^{(1)}\N^{(p)}_{W,s})\big|^2+
\int_{\Mtrap(\tmic, \tau_1)}\big|\pr^{\leq \reg}\underline{\F}^{(p)}_{s}\big|^2\bigg)^{\frac{1}{2}}\nn\\
&&\times\bigg(\int_{\Mtrap(\Iti)}|\pr^{\leq \reg}\psis{p}|^2\bigg)^{\frac{1}{2}}\bigg]\nn\\
&\les& \sum_{p=0}^2\Big(\la^{-1}\widehat{\mathcal{N}}[\pr^{\leq \reg}\psis{p},\pr^{\leq \reg}(\chi_{\tau_1, \tau_2}^{(1)}\N^{(p)}_{W,s})](\tau_1, \tau_1+3)\nn\\
&&+\la^{-1}\widehat{\mathcal{N}}[\pr^{\leq \reg}\psis{p},\pr^{\leq \reg}\underline{\F}^{(p)}_{s}](\tmic, \tau_1)+\la^{-1}\widehat{\mathcal{N}}[\pr^{\leq \reg}\phis{p},\pr^{\leq \reg}\N^{(p)}_{W,s}](\tau_1+1, \tau_2-3)\nn\\
&& +\la^{-1}\widehat{\mathcal{N}}[\pr^{\leq \reg}\psis{p},\pr^{\leq \reg}(\chi_{\tau_1, \tau_2}^{(1)}\N^{(p)}_{W, s})](\tau_2-5, \tau_2-2)+\la\EM[\pr^{\leq \reg}\psis{p}](\Iti)\Big),
\eea
where we have used in the last step the formula \eqref{eq:defintionwidehatmathcalNfpsinormRHS}. Applying Cauchy-Schwarz to the first two terms on the RHS of \eqref{eq:NNtofhighorderpr:spin-2:NwandFbar}, we have 
\begin{align}
\lab{eq:NNtofhighorderpr:spin-2:NwandFbar:s1}
&\sum_{p=0}^2\Big(\la^{-1}\widehat{\mathcal{N}}[\pr^{\leq \reg}\psis{p},\pr^{\leq \reg}(\chi_{\tau_1, \tau_2}^{(1)}\N^{(p)}_{W,s})](\tau_1, \tau_1+3)+\la^{-1}\widehat{\mathcal{N}}[\pr^{\leq \reg}\psis{p},\pr^{\leq \reg}\underline{\F}^{(p)}_{s}](\tmic, \tau_1)\Big)\nn\\
\les{}& \la^{-1}\bigg(\int_{\MM(\tmic, \tau_1)}r^{1+\de}|\pr^{\leq\reg}\underline{\F}_{s}|^2+\int_{\MM(\tau_1,\tau_1+3)} r^{1+\de}|\pr^{\leq \reg}\N_{W,s}|^2\bigg)^{\frac{1}{2}}\bigg(\EM_\de[\pr^{\leq \reg}\pmb\psi](\tmic, \tau_1+3)\bigg)^{\frac{1}{2}}\nn\\
& +\la^{-1}\bigg(\int_{\MM(\tmic, \tau_1)}|\pr^{\leq\reg}\underline{\F}_{s}|^2+\int_{\MM(\tau_1,\tau_1+3)}|\pr^{\leq \reg}\N_{W,s}|^2\bigg)\nn\\
\les{}&\la \EMFtotalhps{s}{\pr^{\leq \reg}} +\la^{-3}\sum_{p=0}^2\bigg(\E[\pr^{\leq \reg}\phis{p}](\tau_1)+\int_{\MM(\tau_1, \tau_2)}r^{1+\de}{|\pr^{\leq \reg}\N_{W,s}^{(p)}|^2}\bigg),
\end{align}
where we have used in the last step the estimates  \eqref{eq:localenergyforpsiontau1minus1tau1plus1} and \eqref{eq:localenergyestimateforunderlineFsijponSigma}.

For the first term in the last line of the RHS of \eqref{eq:NNtofhighorderpr:spin-2:NwandFbar}, we use the fact that it is integrated over a finite interval of time and the support property of $\chi_{\tau_1,\tau_2}^{(1)}$ from \eqref{eq:propertieschi:thisonetodealwithRHSofTeukolsky} to bound it by 
\bea
\lab{eq:NNtofhighorderpr:spin-2:NwandFbar:s2}
&&\la^{-1}\sum_{p=0}^2\widehat{\mathcal{N}}[\pr^{\leq \reg}\psis{p},\pr^{\leq \reg}(\chi_{\tau_1, \tau_2}^{(1)}\N^{(p)}_{W,s})](\tau_2-5, \tau_2-2)\nn\\
&\les&\la^{-1}\sum_{p=0}^2\Big(\EM_{\de}[\pr^{\leq \reg}\psis{p}](\tau_2-3,\tau_2-2)\Big)^{\frac{1}{2}}\bigg(\int_{\MM(\tau_2-5,\tau_2-2)} r^{1+\de}|\pr^{\leq \reg}\N_{W,s}^{(p)}|^2\bigg)^{\frac{1}{2}}\nn\\
&&+\la^{-1}\sum_{p=0}^2\bigg(\int_{\MM(\tau_2-5,\tau_2-2)} r^{1+\de}|\pr^{\leq \reg}\N_{W,s}^{(p)}|^2
+\widehat{\mathcal{N}}[\pr^{\leq \reg}\phis{p}, \pr^{\leq \reg}\N_{W, s}^{(p)}](\tau_2-5,\tau_2-3)\bigg)\nn\\
&\les&\la \sum_{p=0}^2\E[\pr^{\leq \reg}\phis{p}](\tau_2-3) +\la^{-3}\sum_{p=0}^2\int_{\MM(\tau_1,\tau_2)} r^{1+\de}|\pr^{\leq \reg}\N_{W,s}^{(p)}|^2\nn\\
&&
+\la^{-1}\sum_{p=0}^2\widehat{\mathcal{N}}[\pr^{\leq \reg}\phis{p}, \pr^{\leq \reg}\N_{W, s}^{(p)}](\tau_1,\tau_2),
\eea
where in the last step we have used the estimate \eqref{eq:localenergyestimateforpsisijponSigmatau2minus3tau2:partialderivatives}.
Plugging the estimates \eqref{eq:NNtofhighorderpr:spin-2:NwandFbar:s1} and \eqref{eq:NNtofhighorderpr:spin-2:NwandFbar:s2} into \eqref{eq:NNtofhighorderpr:spin-2:NwandFbar}, we deduce, for any $0<\la\leq 1$,
\bea
\lab{eq:NNtofhighorderpr:spin-2:NwandFbar:conclusion}
&&\sum_{p=0}^2\widetilde{\mathcal{N}}[\pr^{\leq \reg}\psis{p},\pr^{\leq \reg}(\chi_{\tau_1, \tau_2}^{(1)}\N^{(p)}_{W,s}+\underline{\F}^{(p)}_{s})]\nn\\
&\les&\la^{-3}\sum_{p=0}^2\bigg(\E[\pr^{\leq \reg}\phis{p}](\tau_1)+\int_{\MM(\tau_1, \tau_2)}r^{1+\de}{|\pr^{\leq \reg}\N_{W,s}^{(p)}|^2}\bigg)\nn\\
&& +\la^{-1}\sum_{p=0}^2\widehat{\mathcal{N}}[\pr^{\leq \reg}\phis{p},\pr^{\leq \reg}\N^{(p)}_{W,s}](\tau_1, \tau_2)+\la \EMFtotalhps{s}{\pr^{\leq \reg}}.
\eea


\subsubsection{Control of the last term on the RHS of \eqref{decomp:NNtofhatFtotal-2:highorderpr}}


Next, we consider the last term on the RHS of \eqref{decomp:NNtofhatFtotal-2:highorderpr}, that is, the term $\sum\limits_{p=0}^2\NNt[\pr^{\leq \reg}\psis{p},\pr^{\leq \reg}\breve{\F}^{(p)}_{s}]$. 

Recall from \eqref{def:breveFpsij} that $\breve{F}^{(p)}_{s,ij}$ is supported in $\MM(\tau_2-1,\tau_2)$, and recall from \eqref{def:breveFpsij} and \eqref{eq:structureofbreveFpsij} that we have on $\MM(\tau_2-1, \tau_2)$
\begin{align}
\lab{eq:origianldefinitionofbreveFwithoutindices}
\breve{F}^{(p)}_{s,ij} ={} &
\Big({\square}_{\gam}-(4-2\de_{p0})|q|^{-2} -\big( \widehat{S}_K+\widehat{Q}_K +f_p\big)\Big)(\chi_{\tau_2}\widetilde{\phi}^{(p)}_{s,ij}+(1-\chi_{\tau_2})\breve{\phi}^{(p)}_{s,ij})\nn\\
={}&{\square}_{\gam}(\chi_{\tau_2}\widetilde{\phi}^{(p)}_{s,ij}+(1-\chi_{\tau_2})\breve{\phi}^{(p)}_{s,ij}) +\sum_{k,l}O(r^{-2})\dk^{\leq 1}(\chi_{\tau_2}\widetilde{\phi}^{(p)}_{s,kl}+(1-\chi_{\tau_2})\breve{\phi}^{(p)}_{s,kl})
\end{align}
and
\bea
\lab{eq:schematicdecompositionofbreveFwithoutindices}
\breve{F}^{(p)}_{s,ij} =  -2\chi_{\tau_2}'(\tau)r^{-1}\pr_r(r(\widetilde{\phi}-\breve{\phi})^{(p)}_{s,ij}) +\sum_{k,l}O(r^{-2})\Big(\chi_{\tau_2}''(\tau), \chi_{\tau_2}'(\tau)\Big)\dk^{\leq 1}(\widetilde{\phi}-\breve{\phi})^{(p)}_{s,kl},
\eea
where the smooth cut-off function $\chi_{\tau_2}=\chi_{\tau_2}(\tau)$ is such that $\chi_{\tau_2}=1$ for $\tau\leq \tau_2-2/3$ and 
$\chi_{\tau_2}=0$ for $\tau\geq \tau_2-1/3$.

We start with proving the following estimate for the most sensitive term in $\NNt[\pr^{\leq \reg}\psis{p},\pr^{\leq \reg}\breve{\F}^{(p)}_{s}]$. This term is the only one which requires integration by parts, while all other terms will be estimated directly.

\begin{lemma}
\lab{lem:esti:breveF:withprtau:highorderpr}
We have, for any $\tau_2-1\leq \tau'<\tau''\leq \tau_2$ and for $\reg\leq 14$, 
\bea
\lab{eq:esti:breveF:withprtau:highorderpr}
\nn&&\sup_{\tau_2-1\leq \tau'<\tau''\leq \tau_2}\sum_{p=0}^2\sum_{i,j}\left|\int_{\MM(\tau',\tau'')}\Re\Big(\pr^{\leq \reg}\breve{F}_{s,ij}^{(p)}\ov{\pr_\tau\pr^{\leq \reg}\psi_{s,ij}^{(p)}}\Big)\right|\\
&\les&\ep \sum_{p=0}^2\bigg(\E[\pr^{\leq \reg}\pmb\phi^{(p)}_{s}](\tau_2-3)+\int_{\MM(\tau_2-3, \tau_2-2)}r^{1+\de} |\pr^{\leq \reg}\N_{W,s}^{(p)}|^2\bigg). 
\eea
\end{lemma}

\begin{proof}
By \eqref{eq:psispij:intermsoftildephibrevephiandphiaux}, we have
\beaa
&&\sum_{p=0}^2\sum_{i,j}\left|\int_{\MM(\tau',\tau'')}\Re\Big(\pr^{\leq \reg}\breve{F}_{s,ij}^{(p)}\ov{\pr_\tau\pr^{\leq \reg}\psi_{s,ij}^{(p)}}\Big)\right|\nn\\
&=&\sum_{p=0}^2\sum_{i,j}\left|\int_{\MM(\tau',\tau'')}\Re\Big(\pr^{\leq \reg}\breve{F}_{s,ij}^{(p)}\ov{\pr_\tau\pr^{\leq \reg}(\chi_{\tau_2}\widetilde{\phi}_{s,ij}^{(p)}+(1-\chi_{\tau_2})\breve{\phi}_{s,ij}^{(p)})}\Big)\right|.
\eeaa
The control of the RHS follows in a similar manner as the proof of \cite[Lemma 4.3]{MaSz24}. To ease the notations in Steps 1--4 below,  we drop the lower index $s,i,j$ as well as the upper index $p$, and denote $\chi_{\tau_2}$ by $\chi$.

\noindent{\bf Step 1.} First, to control the boundary terms on $\II_+$ that will appear in integration by parts, we need an estimate on $\II_+(\tau_2-1, \tau_2)$ for $\widetilde{\phi}-\breve{\phi}$. To this end, we rely on Lemma 2.22 in \cite{MaSz24} which yields
\beaa
\int_{\II_+(\tau_2-1, \tau_2)}r^{-1}|\dk^{\leq 1}\pr^{\leq \reg-1}(\widetilde{\phi}-\breve{\phi})|^2\les\sup_{\tau\in[\tau_2-1, \tau_2]}\E[\pr^{\leq \reg-1}(\widetilde{\phi}-\breve{\phi})](\tau)
\eeaa
so that 
\beaa
\int_{\II_+(\tau_2-1, \tau_2)}|\pr^{\leq \reg}(\widetilde{\phi}-\breve{\phi})|^2 &\les& \int_{\II_+(\tau_2-1, \tau_2)}|\nab_{\pr_\tau}^{\leq \reg}(\widetilde{\phi}-\breve{\phi})|^2+\int_{\II_+(\tau_2-1, \tau_2)}r^{-1}|\dk\pr^{\leq \reg-1}(\widetilde{\phi}-\breve{\phi})|^2\\
&\les& \int_{\II_+(\tau_2-1, \tau_2)}|\widetilde{\phi}-\breve{\phi}|^2+ \EF[\pr^{\leq \reg-1}(\widetilde{\phi}-\breve{\phi})](\tau_2-1, \tau_2).
\eeaa
Also, we have, in view of \eqref{eq:psispij:intermsoftildephibrevephiandphiaux} and \eqref{eq:propertyofscalarizationdefect:prop1},
\beaa
\err_{\textrm{TDefect}}[\widetilde{\phi}]=\err_{\textrm{TDefect}}[\psi]\,\,\,\,\textrm{on}\,\,\,\,\MM(\tau_1, \tau_2-1), \qquad \err_{\textrm{TDefect}}[\psi]=0\,\,\,\,\textrm{on}\,\,\,\,\MM(\tau_1+1, \tau_2-2),
\eeaa
and, in view of \eqref{eq:waveequationdefiningfirstextensionbrevepsi} and Lemma \ref{lemma:computationerrorscalarizationdeffect},
\beaa
|\widetilde{\phi}-\breve{\phi}| =|\widetilde{\phi} - \Pi_2[\widetilde{\phi}]|\les |\err_{\textrm{TDefect}}[\widetilde{\phi}]|\quad\textrm{on}\quad\Si(\tau_2-1), 
\eeaa
which implies 
\beaa
 \int_{\II_+(\tau_2-1, \tau_2)}|\widetilde{\phi}-\breve{\phi}|^2 &\les&  \int_{\II_+(\tau_2-1, \tau_2)}|\nab_{\pr_{\tau}}(\widetilde{\phi}-\breve{\phi})|^2+\int_{\II_+(\tau_2-2, \tau_2-1)}|\nab_{\pr_{\tau}}(\err_{\textrm{TDefect}}[\psi])|^2\\
 &\les& \F_{\II_+}[\widetilde{\phi}-\breve{\phi}](\tau_2-1, \tau_2)+\F_{\II_+}[\err_{\textrm{TDefect}}[\psi]](\tau_2-2, \tau_2-1).
\eeaa
Plugging in the above, we deduce 
\bea\lab{eq:controlofboundarytermsonIIplusforwidetildephiminusbrevephi}
&&\int_{\II_+(\tau_2-1, \tau_2)}|\pr^{\leq \reg}(\widetilde{\phi}-\breve{\phi})|^2\nn\\ 
&\les& \int_{\II_+(\tau_2-1, \tau_2)}|\widetilde{\phi}-\breve{\phi}|^2+ \EF[\pr^{\leq \reg-1}(\widetilde{\phi}-\breve{\phi})](\tau_2-1, \tau_2)\nn\\
&\les& \EF[\pr^{\leq \reg-1}(\widetilde{\phi}-\breve{\phi})](\tau_2-1, \tau_2)+\F_{\II_+}[\err_{\textrm{TDefect}}[\psi]](\tau_2-2, \tau_2-1).
\eea

\noindent{\bf Step 2.} Next, applying Cauchy-Schwarz, we have, in view of \eqref{eq:schematicdecompositionofbreveFwithoutindices},
\beaa
&&\left|\int_{\MM(\tau',\tau'')}\Re\bigg(\pr^{\leq \reg}\Big(\breve{F}+2\chi'(\tau)r^{-1}\pr_r(r(\widetilde{\phi}-\breve{\phi}))\Big)\ov{\pr_\tau\pr^{\leq \reg}(\chi\widetilde{\phi}+(1-\chi)\breve{\phi})}\bigg)\right|\\
&\les& \sqrt{\EM[\pr^{\leq \reg}(\widetilde{\phi}-\breve{\phi})](\tau_2-1, \tau_2)}\sqrt{\EM[\pr^{\leq \reg}(\chi\widetilde{\phi}+(1-\chi)\breve{\phi})](\tau_2-1, \tau_2)}
\eeaa
and hence 
\bea\lab{eq:auxilliaryestimatetocontrolebreveFterminenergyestimate}
\nn&&\left|\int_{\MM(\tau',\tau'')}\Re\Big(\pr^{\leq \reg}\breve{F}\ov{\pr_\tau\pr^{\leq \reg}(\chi\widetilde{\phi}+(1-\chi)\breve{\phi})}\Big)\right|\\
&\les& \mathcal{J}+\sqrt{\EM[\pr^{\leq \reg}(\widetilde{\phi}-\breve{\phi})](\tau_2-1, \tau_2)}\sqrt{\EM[\pr^{\leq \reg}(\chi\widetilde{\phi}+(1-\chi)\breve{\phi})](\tau_2-1, \tau_2)}
\eea
where 
\beaa
 \mathcal{J} &:=& \left|\Re\bigg(\int_{\MM(\tau',\tau'')}\frac{1}{r^2}\pr_r\pr^{\leq \reg}\Big(\chi'(\tau)(r(\widetilde{\phi}-\breve{\phi}))\Big)\ov{\pr_\tau\pr^{\leq \reg}\Big(r(\chi\widetilde{\phi}+(1-\chi)\breve{\phi})\Big)}\bigg)\right|.
\eeaa

\noindent{\bf Step 3.} Next, we integrate by parts in $\pr_r$ in $\mathcal{J}$ and obtain, using \eqref{eq:controlofboundarytermsonIIplusforwidetildephiminusbrevephi} to deal with the boundary terms at $\II_+$, 
\beaa
\mathcal{J}&\les&\left|\int_{\MM(\tau',\tau'')}\Re\bigg(\pr^{\leq \reg}\Big(\chi'(\tau)(r(\widetilde{\phi}-\breve{\phi}))\Big)\ov{\frac{1}{r^2}\pr_r\pr_\tau\pr^{\leq \reg}(r(\chi\widetilde{\phi}+(1-\chi)\breve{\phi}))}\bigg)\right|\\
&&+\sqrt{\EMF[\pr^{\leq \reg}(\widetilde{\phi}-\breve{\phi})](\tau_2-1, \tau_2)+\EMF[\pr^{\leq \reg}(\err_{\textrm{TDefect}}[\psi])](\tau_2-2, \tau_2-1)}\\
&&\times\sqrt{\EMF[\pr^{\leq \reg}(\chi\widetilde{\phi}+(1-\chi)\breve{\phi})](\tau_2-1, \tau_2)}\nn\\
&\les&\left|\int_{\MM(\tau',\tau'')}\Re\bigg(\pr^{\leq \reg}\Big(\chi'(\tau)(\widetilde{\phi}-\breve{\phi})\Big)\ov{\pr^{\leq \reg}\bigg(\frac{1}{r}\pr_r\pr_\tau(r(\chi\widetilde{\phi}+(1-\chi)\breve{\phi}))\bigg)}\bigg)\right|\\
&&+\sqrt{\EMF[\pr^{\leq \reg}(\widetilde{\phi}-\breve{\phi})](\tau_2-1, \tau_2)+\EMF[\pr^{\leq \reg}(\err_{\textrm{TDefect}}[\psi])](\tau_2-2, \tau_2-1)}\\
&&\times\sqrt{\EMF[\pr^{\leq \reg}(\chi\widetilde{\phi}+(1-\chi)\breve{\phi})](\tau_2-1, \tau_2)}.
\eeaa
In view of the following formula from \cite[Equation (4.27)]{MaSz24}
\beaa
\square_{\gam}\phi &=& -2r^{-1}\pr_r(r\pr_\tau\phi) +O(1)\big(\pr_r^2, r^{-1}\pr_{x^a}\pr_r, r^{-2}\pr_{x^a}\pr_{x^b}\big)\phi\nn\\
&&+O(r^{-1})\big(\pr_r, r^{-1}\pr_{x^a}, r^{-1}\pr_\tau\big)\pr^{\leq 1}\phi,
\eeaa
we infer, integrating by parts some of the terms and using again  \eqref{eq:controlofboundarytermsonIIplusforwidetildephiminusbrevephi} to deal with the corresponding boundary terms at $\II_+$, 
\beaa
\mathcal{J}&\les& \left|\int_{\MM(\tau',\tau'')}\Re\Big(\pr^{\leq \reg}\left(\chi'(\tau)(\widetilde{\phi}-\breve{\phi})\right)\ov{\pr^{\leq \reg}\left(\square_{\gam}(\chi\widetilde{\phi}+(1-\chi)\breve{\phi})\right)}\Big)\right|\\
&&+\sqrt{\EMF[\pr^{\leq \reg}(\widetilde{\phi}-\breve{\phi})](\tau_2-1, \tau_2)+\EMF[\pr^{\leq \reg}(\err_{\textrm{TDefect}}[\psi])](\tau_2-2, \tau_2-1)}\\
&&\times\sqrt{\EMF[\pr^{\leq \reg}(\chi\widetilde{\phi}+(1-\chi)\breve{\phi})](\tau_2-1, \tau_2)}
\eeaa
which together with \eqref{eq:origianldefinitionofbreveFwithoutindices} yields
\beaa
\mathcal{J}&\les& \left|\int_{\MM(\tau',\tau'')}\Re\Big(\pr^{\leq \reg}\left(\chi'(\tau)(\widetilde{\phi}-\breve{\phi})\right)\ov{\pr^{\leq \reg}\breve{F}}\Big)\right|\\
&&+\sqrt{\EMF[\pr^{\leq \reg}(\widetilde{\phi}-\breve{\phi})](\tau_2-1, \tau_2)+\EMF[\pr^{\leq \reg}(\err_{\textrm{TDefect}}[\psi])](\tau_2-2, \tau_2-1)}\\
&&\times\sqrt{\EMF[\pr^{\leq \reg}(\chi\widetilde{\phi}+(1-\chi)\breve{\phi})](\tau_2-1, \tau_2)}.
\eeaa

\noindent{\bf Step 4.} Next, using again \eqref{eq:schematicdecompositionofbreveFwithoutindices}, we infer
\beaa
\mathcal{J}&\les& \left|\int_{\MM(\tau',\tau'')}\Re\Big(\pr^{\leq \reg}\left(\chi'(\tau)r(\widetilde{\phi}-\breve{\phi})\right)\ov{r^{-2}\pr_r\pr^{\leq \reg}\left(\chi'(\tau)r(\widetilde{\phi}-\breve{\phi})\right)}\Big)\right|\\
&&+\sqrt{\EMF[\pr^{\leq \reg}(\widetilde{\phi}-\breve{\phi})](\tau_2-1, \tau_2)+\EMF[\pr^{\leq \reg}(\err_{\textrm{TDefect}}[\psi])](\tau_2-2, \tau_2-1)}\\
&&\times\sqrt{\EMF[\pr^{\leq \reg}(\widetilde{\phi})](\tau_2-1, \tau_2)+\EMF[\pr^{\leq \reg}(\breve{\phi})](\tau_2-1, \tau_2)}
\eeaa
and hence
\beaa
\mathcal{J}&\les& \left|\int_{\MM(\tau',\tau'')}\frac{1}{r^2}\pr_r\Big(\big|\pr^{\leq \reg}\big(r\chi'(\tau)(\widetilde{\phi}-\breve{\phi})\big)\big|^2\Big)\right|\\
&&+\sqrt{\EMF[\pr^{\leq \reg}(\widetilde{\phi}-\breve{\phi})](\tau_2-1, \tau_2)+\EMF[\pr^{\leq \reg}(\err_{\textrm{TDefect}}[\psi])](\tau_2-2, \tau_2-1)}\\
&&\times\sqrt{\EMF[\pr^{\leq \reg}(\widetilde{\phi})](\tau_2-1, \tau_2)+\EMF[\pr^{\leq \reg}(\breve{\phi})](\tau_2-1, \tau_2)}.
\eeaa
Integrating by parts in $r$ and using again \eqref{eq:controlofboundarytermsonIIplusforwidetildephiminusbrevephi} to deal with the boundary terms at $\II_+$, we infer
\beaa
\mathcal{J}&\les& \sqrt{\EMF[\pr^{\leq \reg}(\widetilde{\phi}-\breve{\phi})](\tau_2-1, \tau_2)+\EMF[\pr^{\leq \reg}(\err_{\textrm{TDefect}}[\psi])](\tau_2-2, \tau_2-1)}\\
&&\times\sqrt{\EMF[\pr^{\leq \reg}(\widetilde{\phi})](\tau_2-1, \tau_2)+\EMF[\pr^{\leq \reg}(\breve{\phi})](\tau_2-1, \tau_2)}.
\eeaa
Plugging this estimate into \eqref{eq:auxilliaryestimatetocontrolebreveFterminenergyestimate}, we deduce
\bea
\lab{eq:esti:breveF:prtau:highorderpr:middlestep}
\nn&&\left|\int_{\MM(\tau',\tau'')}\Re\Big(\pr^{\leq \reg}\breve{F}\ov{\pr_\tau\pr^{\leq \reg}(\chi\widetilde{\phi}+(1-\chi)\breve{\phi})}\Big)\right|\\
&\les& \mathcal{J}+\sqrt{\EM[\pr^{\leq \reg}(\widetilde{\phi}-\breve{\phi})](\tau_2-1, \tau_2)}\sqrt{\EM[\pr^{\leq \reg}(\chi\widetilde{\phi}+(1-\chi)\breve{\phi})](\tau_2-1, \tau_2)}\nn\\
&\les& \sqrt{\EMF[\pr^{\leq \reg}(\widetilde{\phi}-\breve{\phi})](\tau_2-1, \tau_2)+\EMF[\pr^{\leq \reg}(\err_{\textrm{TDefect}}[\psi])](\tau_2-2, \tau_2-1)}\nn\\
&&\times\sqrt{\EMF[\pr^{\leq \reg}(\widetilde{\phi})](\tau_2-1, \tau_2)+\EMF[\pr^{\leq \reg}(\breve{\phi})](\tau_2-1, \tau_2)}.
\eea

\noindent{\bf Step 5.} Finally, we control the terms appearing on the RHS of the estimate \eqref{eq:esti:breveF:prtau:highorderpr:middlestep}. From \eqref{eq:controlofwidetildephiminusbreevephiontau2minus1tau2} and \eqref{eq:propertyofscalarizationdefect:prop2}, we have, for $\reg\leq 14$,
\bea
\lab{esti:controloferrorofbreveF:pf:eq1}
\nn&&\sum_{p=0}^2\EMF_\de[\pr^{\leq \reg}(\widetilde{\phi}_{s}^{(p)}-\breve{\phi}_{s}^{(p)})](\tau_2-1, \tau_2)+\sum_{p=0}^2\EMF_\de[\pr^{\leq \reg}\err_{\textrm{TDefect}}[\psi^{(p)}_{s}]](\tau_2-2, \tau_2)\\ &\les& \ep^2 \sum_{p=0}^2\bigg(\E[\pr^{\leq \reg}\pmb\phi^{(p)}_{s}](\tau_2-3)+\int_{\MM(\tau_2-3, \tau_2-2)}r^{1+\de} |\pr^{\leq \reg}\N_{W,s}^{(p)}|^2\bigg).
\eea
In addition, we have from the estimates \eqref{eq:localenergyestimateforwidetildephisijponSigmatau2minus3tau2} and \eqref{eq:localenergyestimateforbrevephisijponSigmatau2minus1tau2},  for $\reg\leq 14$,
\bea
\lab{esti:controloferrorofbreveF:pf:eq2}
\nn&&\sum_{i,j=1}^3\sum_{p=0}^2\Big(\EMF_\de[\pr^{\leq \reg}\breve{\phi}^{(p)}_{s,ij}](\tau_2-1, \tau_2)+\EMF_\de[\pr^{\leq \reg}\widetilde{\phi}^{(p)}_{s,ij}](\tau_2-1, \tau_2)\Big)\\ 
&\les& \sum_{p=0}^2\bigg(\E[\pr^{\leq \reg}\pmb\phi^{(p)}_{s}](\tau_2-3)+\int_{\MM(\tau_2-3, \tau_2-2)}r^{1+\de} |\pr^{\leq \reg}\N_{W,s}^{(p)}|^2\bigg).
\eea
Plugging the above two estimates \eqref{esti:controloferrorofbreveF:pf:eq1} and \eqref{esti:controloferrorofbreveF:pf:eq2} into inequality \eqref{eq:esti:breveF:prtau:highorderpr:middlestep} yields the desired estimate \eqref{eq:esti:breveF:withprtau:highorderpr}.
This concludes the proof of Lemma \ref{lem:esti:breveF:withprtau:highorderpr}.
\end{proof}

Next, we control the term $\sum\limits_{p=0}^2\NNt[\pr^{\leq \reg}\psis{p},\pr^{\leq \reg}\breve{\F}^{(p)}_{s}]$ using the above Lemma \ref{lem:esti:breveF:withprtau:highorderpr}. Recall the definition of $\NNt[\c,\c]$ from \eqref{def:NNtintermsofNNtMora:NNtEner:NNtaux:wavesystem:EMF}-\eqref{def:NNtMora:NNtEner:NNtaux:wavesystem:EMF}, as well as the formula \eqref{eq:schematicdecompositionofbreveFwithoutindices} for 
 $\breve{F}^{(p)}_{s,ij}$. First, it follows from \eqref{eq:schematicdecompositionofbreveFwithoutindices} that
\bea
\lab{eq:esti:NNt:breveF:aux:highorderpr}
\sum_{p=0}^2\NNtaux[\pr^{\leq \reg}\breve{\F}^{(p)}_{s}]\les 
\sum_{p=0}^2\EMF[\pr^{\leq \reg}(\widetilde{\phi}_{s}^{(p)}-\breve{\phi}_{s}^{(p)})](\tau_2-1, \tau_2).
\eea
Second, recalling from the context below \eqref{def:NNtMora:NNtEner:NNtaux:wavesystem:EMF} that in $r\geq \Rmic$,  $X$ is a vectorfield satisfying $X=(1+O(r^{-1}))\pr^{\textrm{BL}}_r+A\pr_{\tt}=(A - 1)\pr_{\tau} + O(r^{-1})\dk$ and $w$ is a real-valued function satisfying $w=cr^{-1} + O(r^{-2})$, and in view of the formula \eqref{eq:schematicdecompositionofbreveFwithoutindices} of
 $\breve{F}^{(p)}_{s,ij}$, we have
\bea\lab{eq:esti:NNt:breveF:mora:highorderpr}
&&\sum_{p=0}^2\NNtmora[\pr^{\leq \reg}\psis{p},\pr^{\leq \reg}\breve{\F}^{(p)}_{s}]\nn\\
&\les& 
\Bigg(\sum_{p=0}^2\EMF[\pr^{\leq \reg}(\widetilde{\phi}_{s}^{(p)}-\breve{\phi}_{s}^{(p)})](\tau_2-1, \tau_2)\Bigg)^{\frac{1}{2}}
\Bigg(\sum_{p=0}^2\EMF[\pr^{\leq \reg}\psis{p}](\Iti)\Bigg)^{\frac{1}{2}}\nn\\
&&
+\sum_{p=0}^2\sum_{i,j}\left|\int_{\MM(\tau_1,\tau_2)}\Re\Big(\pr^{\leq \reg}\breve{F}_{s,ij}^{(p)}\ov{\pr_\tau\pr^{\leq \reg}\psi_{s,ij}^{(p)}}\Big)\right|.
\eea
Third, in order to control the $\NNtener$ part, we rely on the smooth cut-off functions $\chi_{\tau,j}$, $j=0,1$,  defined in  \eqref{def:variouscutofffunctionsfortau':tautau'case} to estimate 
\beaa
&&\sup_{\tau\in\Reals}\sum_{n=-1}^{\iota}\bigg|\int_{\Mtrap(\tmic,\tau)}\Re\Big(\ov{|q|^{-2}\Opw(\Theta_n)(\qs \pr^{\leq\reg}\breve{F}^{(p)}_{s,ij})}V_n\Opw(\Theta_n)\pr^{\leq\reg}\psi_{s,ij}^{(p)}\Big)\bigg|\\
&\les& \sup_{\tau\in\Reals}\bigg|\int_{\Mtrap}\Re\Big(\chi_{\tau,0}\ov{\Opw(\widetilde{S}^{0,0}(\MM))\pr^{\leq\reg}\breve{F}^{(p)}_{s,ij}}\Opw(\widetilde{S}^{1,0}(\MM)\pr^{\leq\reg}\psi_{s,ij}^{(p)}\Big)\bigg|\\
&&+\sup_{\tau\in\Reals}\bigg|\int_{\Mtrap}|\chi_{\tau,1}||\Opw(\widetilde{S}^{0,0}(\MM))\pr^{\leq\reg}\breve{F}^{(p)}_{s,ij}||\Opw(\widetilde{S}^{1,0}(\MM)\pr^{\leq\reg}\psi_{s,ij}^{(p)}|\\
&\les& \sup_{\tau\in\Reals}\bigg|\int_{\Mtrap}\Re\Big(\chi_{\tau,0}\ov{\Opw(\widetilde{S}^{0,0}(\MM))\pr^{\leq\reg}\breve{F}^{(p)}_{s,ij}}\Opw(\widetilde{S}^{1,0}(\MM)\pr^{\leq\reg}\psi_{s,ij}^{(p)}\Big)\bigg|\\
&&+\bigg(\int_{\Mtrap}|\pr^{\leq \reg}\breve{\F}^{(p)}_{s}|^2\bigg)^{\frac{1}{2}}\bigg(\sup_{\tau\in\Iti}\E[\pr^{\leq \reg}\psis{p}](\tau)\bigg)^{\frac{1}{2}}
\eeaa
where we used the support property of $\chi_{\tau,1}$ in the last line. New, we introduce a smooth cut-off function $\chi^{(1)}_{\tau_2}=\chi^{(1)}_{\tau_2}(\tau)$ supported in $(\tau_2-1, \tau_2)$ and such that $\chi^{(1)}_{\tau_2}=1$ on the support of $\chi_{\tau_2}$ so that $\breve{F}^{(p)}_{s,ij}=\chi^{(1)}_{\tau_2}\breve{F}^{(p)}_{s,ij}$ in view of \eqref{eq:schematicdecompositionofbreveFwithoutindices}. Then, together with the above, we deduce 
\beaa
&&\sup_{\tau\in\Reals}\sum_{n=-1}^{\iota}\bigg|\int_{\Mtrap(\tmic,\tau)}\Re\Big(\ov{|q|^{-2}\Opw(\Theta_n)(\qs \pr^{\leq\reg}\breve{F}^{(p)}_{s,ij})}V_n\Opw(\Theta_n)\pr^{\leq\reg}\psi_{s,ij}^{(p)}\Big)\bigg|\\
&\les& \sup_{\tau\in\Reals}\bigg|\int_{\Mtrap}\Re\Big(\chi_{\tau,0}\ov{\Opw(\widetilde{S}^{0,0}(\MM))\chi^{(1)}_{\tau_2}\pr^{\leq\reg}\breve{F}^{(p)}_{s,ij}}\Opw(\widetilde{S}^{1,0}(\MM)\pr^{\leq\reg}\psi_{s,ij}^{(p)}\Big)\bigg|\\
&&+\bigg(\int_{\Mtrap}|\pr^{\leq \reg}\breve{\F}^{(p)}_{s}|^2\bigg)^{\frac{1}{2}}\bigg(\sup_{\tau\in\Iti}\E[\pr^{\leq \reg}\psis{p}](\tau)\bigg)^{\frac{1}{2}}\\
&\les& \sup_{\tau\in\Reals}\bigg|\int_{\Mtrap}\Re\Big(\chi_{\tau,0}\ov{\Opw(\widetilde{S}^{0,0}(\MM))\pr^{\leq\reg}\breve{F}^{(p)}_{s,ij}}\Opw(\widetilde{S}^{1,0}(\MM)\chi^{(1)}_{\tau_2}\pr^{\leq\reg}\psi_{s,ij}^{(p)}\Big)\bigg|\\
&&+\bigg(\int_{\Mtrap}|\pr^{\leq \reg}\breve{\F}^{(p)}_{s}|^2\bigg)^{\frac{1}{2}}\bigg(\EM[\pr^{\leq \reg}\psis{p}](\Iti)\bigg)^{\frac{1}{2}}\\
&\les& \bigg(\int_{\Mtrap}|\pr^{\leq \reg}\breve{\F}^{(p)}_{s}|^2\bigg)^{\frac{1}{2}}\bigg(\EM[\pr^{\leq \reg}\psis{p}](\Iti)\bigg)^{\frac{1}{2}}
\eeaa
which yields
\beaa
&&\sum_{p=0}^2\NNtener[\pr^{\leq \reg}\psis{p},\pr^{\leq \reg}\breve{\F}^{(p)}_{s}]\nn\\
&\les&\sum_{p=0}^2\sum_{i,j}\sup_{\tau\geq\tmic}\bigg|\int_{\Mntrap(\tmic, \tau)}{\Re\Big(\pr^{\leq \reg}\breve{F}_{s,ij}^{(p)}\ov{\pr_{\tau}\pr^{\leq \reg}\psi_{s,ij}^{(p)}}\Big)}\bigg|\nn\\
&& +\sum_{p=0}^2\bigg(\int_{\Mtrap}|\pr^{\leq \reg}\breve{\F}^{(p)}_{s}|^2\bigg)^{\frac{1}{2}}\bigg(\EM[\pr^{\leq \reg}\psis{p}](\Iti)\bigg)^{\frac{1}{2}}.
\eeaa
Together with the formula \eqref{eq:schematicdecompositionofbreveFwithoutindices} for
 $\breve{F}^{(p)}_{s,ij}$, we infer
\bea\lab{eq:esti:NNt:breveF:ener:highorderpr}
&&\sum_{p=0}^2\NNtener[\pr^{\leq \reg}\psis{p},\pr^{\leq \reg}\breve{\F}^{(p)}_{s}]\nn\\
&\les& 
\Bigg(\sum_{p=0}^2\EMF[\pr^{\leq \reg}(\widetilde{\phi}_{s}^{(p)}-\breve{\phi}_{s}^{(p)})](\tau_2-1, \tau_2)\Bigg)^{\frac{1}{2}}
\Bigg(\sum_{p=0}^2\EMF[\pr^{\leq \reg}\psis{p}](\Iti)\Bigg)^{\frac{1}{2}}\nn\\
&&
+\sup_{\tau_2-1\leq \tau'<\tau''\leq \tau_2}\sum_{p=0}^2\sum_{i,j}\left|\int_{\MM(\tau',\tau'')}\Re\Big(\pr^{\leq \reg}\breve{F}_{s,ij}^{(p)}\ov{\pr_\tau\pr^{\leq \reg}\psi_{s,ij}^{(p)}}\Big)\right|.
\eea

Combining the above three estimates \eqref{eq:esti:NNt:breveF:aux:highorderpr}, \eqref{eq:esti:NNt:breveF:mora:highorderpr} and \eqref{eq:esti:NNt:breveF:ener:highorderpr}, and applying Lemma \ref{lem:esti:breveF:withprtau:highorderpr} and the estimate \eqref{esti:controloferrorofbreveF:pf:eq1}, we infer
\begin{align}\lab{eq:esti:NNt:breveF:sum:highorderpr}
&\sum_{p=0}^2\NNt[\pr^{\leq \reg}\psis{p},\pr^{\leq \reg}\breve{\F}^{(p)}_{s}]\nn\\
\les{}& 
\Bigg(\sum_{p=0}^2\EMF[\pr^{\leq \reg}(\widetilde{\phi}_{s}^{(p)}-\breve{\phi}_{s}^{(p)})](\tau_2-1, \tau_2)\Bigg)^{\frac{1}{2}}
\Bigg(\sum_{p=0}^2\EMF[\pr^{\leq \reg}\psis{p}](\Iti)\Bigg)^{\frac{1}{2}}\nn\\
&
+\sum_{p=0}^2\EMF[\pr^{\leq \reg}(\widetilde{\phi}_{s}^{(p)}-\breve{\phi}_{s}^{(p)})](\tau_2-1, \tau_2)\nn\\
&+\sup_{\tau_2-1\leq \tau'<\tau''\leq \tau_2}\sum_{p=0}^2\sum_{i,j}\left|\int_{\MM(\tau',\tau'')}\Re\Big(\pr^{\leq \reg}\breve{F}_{s,ij}^{(p)}\ov{\pr_\tau\pr^{\leq \reg}\psi_{s,ij}^{(p)}}\Big)\right|\nn\\
\les{}&\ep \sum_{p=0}^2\bigg(\E[\pr^{\leq \reg}\pmb\phi^{(p)}_{s}](\tau_2-3)+\EMF[\pr^{\leq \reg}\psis{p}](\Iti)+\int_{\MM(\tau_2-3, \tau_2-2)}r^{1+\de} |\pr^{\leq \reg}\N_{W,s}^{(p)}|^2\bigg).
\end{align}


\subsubsection{Concluding the proof of Proposition \ref{prop:controlofhighorderNNtspin-2:lastcontrol:pfofthm7.5}}
\lab{sec:endofproofofprop:controlofhighorderNNtspin-2:lastcontrol:pfofthm7.5}


Plugging the estimates \eqref{eq:esti:NNt:breveF:sum:highorderpr} and \eqref{eq:NNtofhighorderpr:spin-2:NwandFbar:conclusion} into \eqref{decomp:NNtofhatFtotal-2:highorderpr}, we deduce, for any $0<\la\leq 1$,
\bea
\lab{eq:controlofallhighorderNNtterms:s=-2:pfofthm7.5}
&&\sum\limits_{p=0}^2\NNt[\pr^{\leq \reg}\psis{p},\pr^{\leq \reg}\widehat{\F}_{total,s}^{(p)}]\nn\\
&\les&\ep \sum_{p=0}^2\bigg(\E[\pr^{\leq \reg}\pmb\phi^{(p)}_{s}](\tau_2-3)+\EMF[\pr^{\leq \reg}\psis{p}](\Iti)+\int_{\MM(\tau_2-3, \tau_2-2)}r^{1+\de} |\pr^{\leq \reg}\N_{W,s}^{(p)}|^2\bigg)\nn\\
&&+\la^{-3}\sum_{p=0}^2\bigg(\E[\pr^{\leq \reg}\phis{p}](\tau_1)+\int_{\MM(\tau_1, \tau_2)}r^{1+\de}{|\pr^{\leq \reg}\N_{W,s}^{(p)}|^2}\bigg)\nn\\
&& +
\la^{-1}\sum_{p=0}^2\widehat{\mathcal{N}}[\pr^{\leq \reg}\phis{p},\pr^{\leq \reg}\N^{(p)}_{W,s}](\tau_1, \tau_2)+\la \EMFtotalhps{s}{\pr^{\leq \reg}}.
\eea

We are now ready to prove Proposition \ref{prop:controlofhighorderNNtspin-2:lastcontrol:pfofthm7.5}.

\begin{proof}[Proof of Proposition \ref{prop:controlofhighorderNNtspin-2:lastcontrol:pfofthm7.5}] Recall the formula \eqref{expression:NNttotalps:pm2:pasteinproof:s=pm2} for $\NNttotalph{s}{\pr^{\leq \reg}}$. We use the bound for $\widehat{\NN}[\pmb\psi, \pmb H](\tau_1, \tau_2)$ in Remark \ref{rmk:controlofwidehatNfpsibyNfpsi}, as well as the following bound for $\widehat{\NN}'[\pmb\psi, \pmb H](\tau_1, \tau_2)$
\beaa
\widehat{\mathcal{N}}'[\pmb\psi, \H](\tau_1, \tau_2)\les \Big(\M_\de[\pmb\psi](\tau_1,\tau_2)\Big)^{\frac{1}{2}}\left(\int_{\MM(\tau_1, \tau_2)}r^{1+\de}|\H|^2\right)^{\frac{1}{2}}+\int_{\MM(\tau_1, \tau_2)}r^{1+\de}|\H|^2,
\eeaa
to control the $\widehat{\NN}'$ and $\widehat{\NN}$ terms appearing on the RHS of \eqref{expression:NNttotalps:pm2:pasteinproof:s=pm2} and \eqref{eq:controlofallhighorderNNtterms:s=-2:pfofthm7.5} by
\bea
\lab{esti:controlofNNhat'andNNhat:pfofThm7.5:-2}
&&\sum_{p=0}^2\widehat{\mathcal{N}}'[\pr^{\leq \reg}\phis{p}, \pr^{\leq \reg}\N_{W, s}^{(p)}](\tau_1, \tau_2)
+\sum_{p=0}^1\widehat{\mathcal{N}}'[ \pr^{\leq \reg}\phis{p}, r^{-2}\pr^{\leq \reg}\N_{T,  s}^{(p)}](\tau_1, \tau_2)\nn\\
&&+\la^{-1}\sum_{p=0}^2\widehat{\mathcal{N}}[\pr^{\leq \reg}\phis{p},\pr^{\leq \reg}\N^{(p)}_{W, s}](\tau_1, \tau_2)\nn\\
&\les&
\Bigg(\la^{-2}\sum_{p=0}^2{\mathcal{N}}_{\de}[\pr^{\leq \reg}\phis{p}, \pr^{\leq \reg}\N_{W,  s}^{(p)}](\tau_1, \tau_2)
+\sum_{p=0}^1\int_{\MM(\tau_1,\tau_2)}r^{-3+\de}\big|\pr^{\leq \reg}\N_{T,  s}^{(p)}\big|^2
\Bigg)^{\frac{1}{2}}\nn\\
&&\times \Big(\EMFtotalhps{s}{\pr^{\leq \reg}}\Big)^{\frac{1}{2}}+\la^{-1}\sum_{p=0}^2{\mathcal{N}}_{\de}[\pr^{\leq \reg}\phis{p}, \pr^{\leq \reg}\N_{W,  s}^{(p)}](\tau_1, \tau_2)\nn\\
&&
+\sum_{p=0}^1\int_{\MM(\tau_1,\tau_2)}r^{-3+\de}\big|\pr^{\leq \reg}\N_{T,  s}^{(p)}\big|^2.
\eea
For the terms in the last line of \eqref{expression:NNttotalps:pm2:pasteinproof:s=pm2}, we separate the integration onto $\MM_{r\leq 10m}(\tau_1,\tau_2)$ and $\MM_{r\geq 10m}(\tau_1,\tau_2)$ to deduce
\bea
\lab{eq:estimateofproductofNwsandNtswithphisp:p=01}
&&\sum_{p=0}^1\int_{\MM(\tau_1,\tau_2)} \Big(r^{-1+\de}\big|\pr^{\leq \reg}\N^{(p)}_{W,s}\big|+r^{-2+\de}\big|\pr\pr^{\leq \reg}\N^{(p)}_{T,s}\big|\Big)\big|\pr^{\leq \reg}\phis{p}\big|\nn\\
&\les&\la \EMFtotalhps{s}{\pr^{\leq \reg}}
+\la^{-1}\sum_{p=0}^1\int_{\MM_{r\leq 10m}(\tau_1,\tau_2)}\Big(r^{-3+\de}\big|\pr^{\leq \reg+1}\N_{T, s}^{(p)}\big|^2
+r^{1+\de}\big|\pr^{\leq \reg}\N_{W,s}^{(p)}\big|^2\Big)
\nn\\
&&+\sum_{p=0}^1\int_{\MM_{r\geq 10m}(\tau_1,\tau_2)} \Big(r^{-1+\de}\big|\pr^{\leq \reg}\N^{(p)}_{W,s}\big|+r^{-2+\de}\big|\pr^{\leq \reg+1}\N^{(p)}_{T,s}\big|\Big)\big|\pr^{\leq \reg}\phis{p}\big|.
\eea

Finally, plugging  \eqref{eq:controlofallhighorderNNtterms:s=-2:pfofthm7.5}, \eqref{esti:controlofNNhat'andNNhat:pfofThm7.5:-2} and \eqref{eq:estimateofproductofNwsandNtswithphisp:p=01} into \eqref{expression:NNttotalps:pm2:pasteinproof:s=pm2}, and using the following estimate which follows from the definition \eqref{eq:defmathcalNpsif} of $\NN_{\de}[\c,\c](\c,\c)$
\bea
\lab{eq:spacetimeofNwcontrolledbyNNdenorms}
\int_{\MM(\tau',\tau'')} r^{1+\de}\big|\pr^{\leq \reg}\N_{W,s}^{(p)}\big|^2\les {\mathcal{N}}_{\de}[\pr^{\leq \reg}\phis{p}, \pr^{\leq \reg}\N_{W, s}^{(p)}](\tau', \tau''), \qquad \forall \,\tau'<\tau'',
\eea
we infer the desired estimates \eqref{eq:controlofhighorderNNtspin-2:lastcontrol:pfofthm7.5}.
This concludes the proof of Proposition \ref{prop:controlofhighorderNNtspin-2:lastcontrol:pfofthm7.5}.
\end{proof}


\subsection{End of the proof of Theorem \ref{th:main:intermediary}}
\lab{sec:endfotheproofofth:main:intermediary}


In this section, we conclude the proof for Theorem \ref{th:main:intermediary}. 
Substituting the estimates \eqref{eq:localenergyforpsiontau1minus1tau1plus2:initialEMF}, \eqref{esti:Errdefectterm:highorderprderi} and \eqref{eq:controlofhighorderNNtspin-2:lastcontrol:pfofthm7.5} back into \eqref{eq:EMFtotalp:finalestimate:pm2:prhighorder}, taking both of $\la$ and $\ep$ suitably small such that  the term $(\la+\ep)\EMFtotalhps{s}{\pr^{\leq \reg}}$ in \eqref{eq:controlofhighorderNNtspin-2:lastcontrol:pfofthm7.5} is absorbed by the LHS of \eqref{eq:EMFtotalp:finalestimate:pm2:prhighorder},  and noticing from the definition \eqref{def:AandAonorms:tau1tau2:phisandpsis} of the norm $\A[\cdot, \cdot](\cdot, \cdot)$ that for $s=\pm 2$ we have
\beaa
\A[\pmb\psi_{s}](\Iti)+\A[\pmb\phi_{s}](\tau_1,\tau_2)
\les\sum_{p=0}^2\Big(\EMF[\psis{p}](\Iti)
+\EMF[\phis{p}](\tau_1, \tau_2)\Big),
\eeaa
we deduce
\bea
\lab{eq:EMFtotalpspm2:finalestimatetoprovemaintheom:proof} 
&&\EMFtotalh{s}{\nab_{(\pr_{\tt},\chi_{0}\pr_{\tilde{\phi}})}^{\leq \reg}}+\EMFtotalhps{s}{\pr^{\leq \reg}}  \nn\\
&\les&{\sum_{p=0}^2}\E[\pr^{\leq \reg}\phis{p}](\tau_1)  +\sum_{p=0}^2\Big(\EMF[\psis{p}](\Iti)
+\EMF[\phis{p}](\tau_1, \tau_2)\Big)
\nn\\
&&+\sum_{p=0}^2{\mathcal{N}}_{\de}[\pr^{\leq \reg}\phis{p}, \pr^{\leq \reg}\N_{W, s}^{(p)}](\tau_1, \tau_2)
+\sum_{p=0}^1\int_{\MM(\tau_1,\tau_2)}r^{-3+\de}|\pr^{\leq \reg+1}\N_{T,s}^{(p)}|^2
\nn\\
&&+\sum_{p=0}^1\int_{\MM_{r\geq 10m}(\tau_1,\tau_2)} \Big(r^{-1+\de}|\pr^{\leq \reg}\N^{(p)}_{W,s}|+r^{-2+\de}|\pr^{\leq \reg+1}\N^{(p)}_{T,s}|\Big)|\pr^{\leq \reg}\phis{p}|\nn\\
&&
+\sum_{p=0}^{2}\widehat{\mathcal{N}}''[\pr^{\leq \reg}\psis{p}, \pr^{\leq \reg}\widehat{\F}_{total,s}^{(p)}](\tau_2-3, \tau_2).
\eea

Next, we estimate the last term on the RHS of \eqref{eq:EMFtotalpspm2:finalestimatetoprovemaintheom:proof}. In view of the formula \eqref{eq:definitionofwidehatFtotalsij} for $\widehat{F}_{total,s,ij}^{(p)}$, we have
\beaa
\widehat{F}_{total,s,ij}^{(p)}=\chi_{\tau_1, \tau_2}^{(1)}N^{(p)}_{W,s,ij}+\breve{F}^{(p)}_{s,ij}, \quad \text{on} \quad \MM(\tau_2-3,\tau_2),
\eeaa
and applying Cauchy-Schwarz to $\widehat{\mathcal{N}}''$ as given in \eqref{eq:defofNNhat''}, we infer
\beaa
&&\sum_{p=0}^{2}\widehat{\mathcal{N}}''[\pr^{\leq \reg}\psis{p}, \pr^{\leq \reg}\widehat{\F}_{total,s}^{(p)}](\tau_2-3, \tau_2)\nn\\
&\les& \sum_{p=0}^2\Big(\EMF_\de[\pr^{\leq \reg}\psis{p}](\tau_2-3,\tau_2)\Big)^{\frac{1}{2}}\bigg(\int_{\MM(\tau_2-3,\tau_2)}r^{1+\de}|\pr^{\leq \reg}(\chi_{\tau_1, \tau_2}^{(1)}\N^{(p)}_{W,s})|^2\bigg)^{\frac{1}{2}}\nn\\
&&+ \sum_{p=0}^2\int_{\MM(\tau_2-3,\tau_2)}|\pr^{\leq \reg}(\chi_{\tau_1, \tau_2}^{(1)}\N^{(p)}_{W,s})|^2
+\sum_{p=0}^{2}\widehat{\mathcal{N}}''[\pr^{\leq \reg}\psis{p}, \pr^{\leq \reg}\breve{\F}^{(p)}_{s}](\tau_2-3, \tau_2)\nn\\
&\les&\sum_{p=0}^2\Big(\EMF_\de[\pr^{\leq \reg}\psis{p}](\tau_2-3,\tau_2)\Big)^{\frac{1}{2}}\bigg(\int_{\MM(\tau_2-3,\tau_2)}r^{1+\de}|\pr^{\leq \reg}\N^{(p)}_{W,s}|^2\bigg)^{\frac{1}{2}}\nn\\
&&+\sum_{p=0}^2{\mathcal{N}}_{\de}[\pr^{\leq \reg}\phis{p}, \pr^{\leq \reg}\N_{W, s}^{(p)}](\tau_2-3, \tau_2)
+\sum_{p=0}^{2}\widehat{\mathcal{N}}''[\pr^{\leq \reg}\psis{p}, \pr^{\leq \reg}\breve{\F}^{(p)}_{s}](\tau_2-3, \tau_2),
\eeaa
where in the last step we have used the estimate \eqref{eq:spacetimeofNwcontrolledbyNNdenorms}. 
Plugging the formula \eqref{eq:schematicdecompositionofbreveFwithoutindices} for $\breve{F}^{(p)}_{s,ij}$, using the estimates \eqref{eq:esti:breveF:withprtau:highorderpr} and \eqref{esti:controloferrorofbreveF:pf:eq1}, and in view of the definition \eqref{eq:defofNNhat''} of $\widehat{\mathcal{N}}''$, we infer 
\beaa
&&\sum_{p=0}^{2}\widehat{\mathcal{N}}''[\pr^{\leq \reg}\psis{p}, \pr^{\leq \reg}\breve{\F}^{(p)}_{s}](\tau_2-3, \tau_2)\nn\\
&\les&\ep \sum_{p=0}^2\bigg(\E[\pr^{\leq \reg}\pmb\phi^{(p)}_{s}](\tau_2-3)+\int_{\MM(\tau_2-3, \tau_2-2)}r^{1+\de} |\pr^{\leq \reg}\N_{W,s}^{(p)}|^2
+\EMF[\pr^{\leq \reg}\psis{p}](\tau_2-3,\tau_2)\bigg),
\eeaa
which together with the previous estimate yields
\bea
\lab{eq:controlofNNtrgeq10m:Fhattotal:tau2-3tau2}
&&\sum_{p=0}^{2}\widehat{\mathcal{N}}''[\pr^{\leq \reg}\psis{p}, \pr^{\leq \reg}\widehat{\F}_{total,s}^{(p)}](\tau_2-3, \tau_2)\nn\\
&\les&\sum_{p=0}^2\Big(\EMF_\de[\pr^{\leq \reg}\psis{p}](\tau_2-3,\tau_2)\Big)^{\frac{1}{2}}\bigg(\int_{\MM(\tau_2-3,\tau_2)}r^{1+\de}|\pr^{\leq \reg}\N^{(p)}_{W,s}|^2\bigg)^{\frac{1}{2}}\nn\\
&&+\ep \sum_{p=0}^2\Big(\EMF[\pr^{\leq \reg}\psis{p}](\tau_2-3,\tau_2)+\E[\pr^{\leq \reg}\pmb\phi^{(p)}_{s}](\tau_2-3)\Big)\nn\\
&&+ \sum_{p=0}^2{\mathcal{N}}_{\de}[\pr^{\leq \reg}\phis{p}, \pr^{\leq \reg}\N_{W, s}^{(p)}](\tau_2-3, \tau_2).
\eea

Next, since $\psis{p}$ is a solution in a subextremal Kerr background to the tensorial wave equation \eqref{eq:basictensorialwaveequationinKerr:forbrevephi} in $\MM(\tau\geq \tau_2)$,  we can apply a direct generalization of the estimates in Theorem \ref{prop:weakMorawetzfortensorialwaveeqfromDRSR} to high-order  unweighted derivatives, and improve $\M$ to $\M_\de$, to deduce
\beaa
\sum_{i,j=1}^3\sum_{p=0}^2\EMF_\de[\pr^{\leq \reg}\psi^{(p)}_{s,ij}](\tau_2, +\infty)
\les \sum_{i,j=1}^3\sum_{p=0}^2\E[\pr^{\leq \reg}\psi^{(p)}_{s,ij}](\tau_2),
\eeaa
which together with \eqref{eq:localenergyestimateforpsisijponSigmatau2minus3tau2:partialderivatives} implies
\bea
\lab{esti:EMFdeofpsisp:tau2-3toinfty:prderi}
\nn&&\sum_{i,j=1}^3\sum_{p=0}^2\EMF_\de[\pr^{\leq \reg}\psi^{(p)}_{s,ij}](\tau_2-3,+\infty)\\ 
&\les& \sum_{i,j=1}^3\sum_{p=0}^2\bigg(\E[\pr^{\leq \reg}\phi^{(p)}_{s,ij}](\tau_2-3)+\int_{\MM(\tau_2-3, \tau_2-2)}r^{1+\de} |\pr^{\leq \reg}N_{W,s,ij}^{(p)}|^2\bigg).
\eea

Combining the estimates \eqref{eq:localenergyforpsiontau1minus1tau1plus1},  \eqref{eq:EMFtotalpspm2:finalestimatetoprovemaintheom:proof}, \eqref{eq:controlofNNtrgeq10m:Fhattotal:tau2-3tau2}  and \eqref{esti:EMFdeofpsisp:tau2-3toinfty:prderi}, noticing that the sum of the LHS of \eqref{eq:localenergyforpsiontau1minus1tau1plus1},  \eqref{eq:EMFtotalpspm2:finalestimatetoprovemaintheom:proof}  and \eqref{esti:EMFdeofpsisp:tau2-3toinfty:prderi} controls the LHS of \eqref{eq:mainhighorderunweightedEMF:maintheorem:pm2}, and using the estimate \eqref{eq:spacetimeofNwcontrolledbyNNdenorms},  the desired estimate \eqref{eq:mainhighorderunweightedEMF:maintheorem:pm2} then follows.  This concludes the proof of Theorem \ref{th:main:intermediary}.


\section{EMF estimates near infinity for Teukolsky in perturbations of Kerr}
\label{sec:proofofprop:EnergyMorawetznearinfinitytensorialTeuk}


The goal of this section is to prove Proposition \ref{prop:EnergyMorawetznearinfinitytensorialTeuk} establishing EMF estimates near infinity for tensorial Teukolsky equations in perturbations of Kerr.
To this end, let us first state an EMF estimate near infinity for an inhomogeneous tensorial wave equation.

\begin{lemma}
\lab{lem:zeroorderEMF:inhomogeneoustensorialwave}
Let $(\MM, \g)$ satisfy the assumptions of Sections \ref{subsect:assumps:perturbednullframe} and  \ref{subsubsect:assumps:perturbedmetric}. Assume that $\pmb{\psi}\in\sk_k(\mathbb{C})$, $k=1,2$, satisfies
\bea
\lab{eq:weightedEMnearinf:inhomotensorialwave}
\squared_k\pmb{\psi} - D_0r^{-2}\pmb{\psi} &=& D_1r^{-1}\nab_{\pr_{\tt}}\pmb{\psi} -D_2r^{-2}\pmb{\psi}+O(r^{-2})\nab_3\pmb\psi+O(r^{-3})\dk^{\leq 1}\pmb\psi\nn\\
&& +O(\ep r^{-2}\tau^{-\frac{1}{2}-\dec})\dk^{\leq 1}\pmb{\psi} +D_1O(\ep r^{-1})\nab_{\pr_\tau}\pmb{\psi} + \N,
\eea
with constants $D_0> 0$ and $D_1, D_2\geq 0$. If $D_1>0$, then we have, for any $\tau_1<\tau_2$ and $R\geq 20m$ large enough,
\bea
\lab{eq:weightedEMFnearinf:D1>0:0order}
&&\EMF_{0,r\geq R}[\pmb{\psi}](\tau_1,\tau_2) \nn\\
&\les& \E_{r\geq R/2}[\pmb{\psi}](\tau_1) + R^2 \M_{R/2, R}[\pmb{\psi}](\tau_1,\tau_2) +D_2\sup_{\tau\in[\tau_1, \tau_2]}\int_{\Si(\tau)}r^{-2}|\pmb\psi|^2\nn\\
&&+D_2\int_{\MM(\tau_1, \tau_2)}r^{-3}|\pmb\psi|^2 +D_2\int_{\II_+(\tau_1, \tau_2)}r^{-2}|\pmb\psi|^2 + \int_{\MM_{r\geq R/2}(\tau_1,\tau_2)} r |\N|^2.
\eea
If $D_1=0$, then we have, for any $\tau_1<\tau_2$, $\de\in(0, \frac{1}{3}]$ and $R\geq 20m$ large enough,
\bea
\lab{eq:weightedEMFnearinf:D1=0:0order}
&&\EMF_{\de,r\geq R}[\pmb{\psi}](\tau_1,\tau_2)\nn\\
 &\les& \E_{r\geq R/2}[\pmb{\psi}](\tau_1) + R^2 \M_{R/2, R}[\pmb{\psi}](\tau_1,\tau_2) +D_2\sup_{\tau\in[\tau_1, \tau_2]}\int_{\Si(\tau)}r^{-2}|\pmb\psi|^2\nn\\
&& +D_2\int_{\MM(\tau_1, \tau_2)}r^{-3}|\pmb\psi|^2 +D_2\int_{\II_+(\tau_1, \tau_2)}r^{-2}|\pmb\psi|^2+ \int_{\MM_{r\geq R/2}(\tau_1,\tau_2)} r^{1+\de} |\N|^2.
\eea
\end{lemma}

\begin{proof}
We first rewrite \eqref{eq:weightedEMnearinf:inhomotensorialwave} as 
\bea
\lab{eq:weightedEMnearinf:inhomotensorialwave:simpifiedform}
\squared_k\pmb{\psi} - D_0r^{-2}\pmb{\psi} = D_1r^{-1}\nab_{\pr_{\tt}}\pmb{\psi} -D_2r^{-2}\pmb{\psi}+D_1O(\ep r^{-1})\nab_{\pr_\tau}\pmb{\psi}+\F
\eea
where $\F$ is given by
\beaa
\F:=O(r^{-2})\nab_3\pmb\psi+O(r^{-3})\dk^{\leq 1}\pmb\psi+O(\ep r^{-2}\tau^{-\frac{1}{2}-\dec})\dk^{\leq 1}\pmb{\psi}  + \N
\eeaa 
and satisfies for any $\de_0\in[0,\frac{1}{3}]$,
\bea\lab{eq:weightedEMnearinf:inhomotensorialwave:simpifiedform:estF}
\int_{\MM_{r\geq R}(\tau_1, \tau_2)}r^{1+\de_0}|\F|^2 &\les& R^{-1+\de_0}\M_{r\geq R}[\pmb\psi](\tau_1, \tau_2)+\ep\EM_{r\geq R}[\pmb\psi](\tau_1, \tau_2)\nn\\
&&+\int_{\MM_{r\geq R}(\tau_1, \tau_2)}r^{1+\de_0}|\N|^2.
\eea
Also, integrating by parts, we have
\bea\lab{eq:weightedEMnearinf:inhomotensorialwave:simpifiedform:estF1}
&&\left|\int_{\MM(\tau_1, \tau_2)}\Re\left(\ov{r^{-2}\pmb\psi}\c\Big(f_1(r)\nab_{\tau},\, f_2(r)\nab_{\pr_r}, f_3(r, \cos\th)\nab\Big)\pmb\psi\right)\right|\nn\\
&\les& \sup_{\tau\in[\tau_1, \tau_2]}\int_{\Si(\tau)}r^{-2}|\pmb\psi|^2+\int_{\MM(\tau_1, \tau_2)}r^{-3}|\pmb\psi|^2 +\int_{\II_+(\tau_1, \tau_2)}r^{-2}|\pmb\psi|^2
\eea
for any real-valued scalar functions $f_1$, $f_2$, $f_3$  supported in $r\geq R/2$ such that 
\bea\lab{eq:weightedEMnearinf:inhomotensorialwave:simpifiedform:condfjj=123}
f_1(r),\, f_2(r),\, f_3(r,\cos\th)=O(1), \qquad f_1'(r), \, f_2'(r),\, \big(\pr_r, r^{-1}\pr_{\cos\th}\big)f_3(r,\cos\th)=O(r^{-1}).
\eea

Next, consider first the case $D_1>0$. 
Applying Proposition \ref{prop-app:stadard-comp-Psi} with $(X=\chi_{R}\pr_{\tt}, w=0)$ to equation \eqref{eq:weightedEMnearinf:inhomotensorialwave}, where $\chi_R$ is a smooth cutoff function satisfying
\bea
\lab{def:cutoffcuntionbetweenRandR/2}
\chi_{R}=1 \quad \text{for} \quad r\geq R, \qquad  \chi_{ R}=0\quad \text{for} \quad r\leq R/2,
\eea
and noticing from \eqref{eq:estimatesforQQtimesdeformationtensorofT} and Corollary \ref{cor:computationofthecomponentsofthetensorAforexactenergyconservationKerr} that
\beaa
\bigg|\int_{\MM_{r\geq R}(\tau_1,\tau_2)}\frac{1}{2}\QQ[\pmb{\psi}]  \c{^{(\pr_{\tt})}}\pi\bigg|
&\les& \ep \EM_{r\geq R}[\pmb{\psi}](\tau_1,\tau_2)  ,\\
\frac{k}{2}{}^{(\pr_{\tt})}A_\nu \Im\Big(\pmb\psi\c\ov{\Ddot^{\nu}\pmb\psi}\Big)&=&\big(O(r^{-4}) +r^{-1}\Ga_b\big)\Im\Big(\pmb\psi\c\ov{\pr\pmb\psi }\Big),
\eeaa
 we infer, for $R\geq 20m$,
\beaa
&&\sup_{\tau\in(\tau_1, \tau_2)}\int_{\Si_{R/2,R}(\tau)}r^{-2}\chi_R|\dk^{\leq 1}\pmb\psi|^2+\EF_{r\geq R}[\pmb{\psi}](\tau_1,\tau_2)+D_1\int_{\MM_{r\geq R}(\tau_1,\tau_2)} r^{-1} \big|\nab_{\pr_{\tt}}\pmb{\psi}\big|^2\nn\\
&\les&\E_{r\geq R/2} [\pmb{\psi}](\tau_1) +\ep\EM_{r\geq R}[\pmb{\psi}](\tau_1,\tau_2) +R^{-1}\M_{r\geq R}[\pmb{\psi}](\tau_1,\tau_2) + R^2 \M_{ R/2, R}[\pmb{\psi}](\tau_1,\tau_2)\\
&&+\bigg|\int_{\MM_{r\geq  R/2}(\tau_1,\tau_2)} \Re\Big(\ov{\big(-D_2r^{-2}\pmb{\psi}+D_1O(\ep r^{-1})\nab_{\pr_\tau}\pmb{\psi}+\F\big)}\chi_{R}\nab_{\pr_{\tt}}\pmb{\psi}\Big)\bigg|,
\eeaa
which together with \eqref{eq:weightedEMnearinf:inhomotensorialwave:simpifiedform:estF} for $\de_0=0$ and \eqref{eq:weightedEMnearinf:inhomotensorialwave:simpifiedform:estF1} yields
\beaa
&&\sup_{\tau\in(\tau_1, \tau_2)}\int_{\Si_{R/2,R}(\tau)}r^{-2}\chi_R|\dk^{\leq 1}\pmb\psi|^2+\EF_{r\geq R}[\pmb{\psi}](\tau_1,\tau_2)+D_1\int_{\MM_{r\geq R}(\tau_1,\tau_2)} r^{-1} \big|\nab_{\pr_{\tt}}\pmb{\psi}\big|^2\nn\\
&\les&\E_{r\geq R/2} [\pmb{\psi}](\tau_1) +\ep\EM_{r\geq R}[\pmb{\psi}](\tau_1,\tau_2) +R^{-1}\M_{r\geq R}[\pmb{\psi}](\tau_1,\tau_2) + R^2 \M_{ R/2, R}[\pmb{\psi}](\tau_1,\tau_2)\\
&&+ D_2\sup_{\tau\in[\tau_1, \tau_2]}\int_{\Si(\tau)}r^{-2}|\pmb\psi|^2+D_2\int_{\MM(\tau_1, \tau_2)}r^{-3}|\pmb\psi|^2 +D_2\int_{\II_+(\tau_1, \tau_2)}r^{-2}|\pmb\psi|^2\nn\\
&&
+D_1\ep\int_{\MM_{r\geq R}(\tau_1,\tau_2)} r^{-1} \big|\nab_{\pr_{\tt}}\pmb{\psi}\big|^2+\int_{\MM_{r\geq R}(\tau_1, \tau_2)}r|\N|^2,
\eeaa
and hence, taking $\ep$ small enough, we deduce
\bea
\lab{eq:EFestimate:case1D1>0:0order:proof}
&& \sup_{\tau\in(\tau_1, \tau_2)}\int_{\Si_{R/2,R}(\tau)}r^{-2}\chi_R|\dk^{\leq 1}\pmb\psi|^2+\EF_{r\geq R}[\pmb{\psi}](\tau_1,\tau_2)+D_1\int_{\MM_{r\geq R}(\tau_1,\tau_2)} r^{-1} \big|\nab_{\pr_{\tt}}\pmb{\psi}\big|^2\nn\\
\nn&\les&\E_{r\geq R/2} [\pmb{\psi}](\tau_1) +\big(\ep+R^{-1}\big)\M_{r\geq R}[\pmb{\psi}](\tau_1,\tau_2) + R^2 \M_{ R/2, R}[\pmb{\psi}](\tau_1,\tau_2)\\
&& + D_2\sup_{\tau\in[\tau_1, \tau_2]}\int_{\Si(\tau)}r^{-2}|\pmb\psi|^2+D_2\int_{\MM(\tau_1, \tau_2)}r^{-3}|\pmb\psi|^2\nn\\
&& +D_2\int_{\II_+(\tau_1, \tau_2)}r^{-2}|\pmb\psi|^2 + \int_{\MM_{ r\geq R/2}(\tau_1,\tau_2)} r |\N|^2.
\eea

Next, applying  Proposition \ref{prop-app:stadard-comp-Psi} with $(X=0, w=\chi_{R}r^{-1})$, and noticing in view of \eqref{eq:assymptiticpropmetricKerrintaurxacoord:volumeform}, \eqref{eq:assymptiticpropmetricKerrintaurxacoord:1}, \eqref{eq:controloflinearizedinversemetriccoefficients} and Lemma \ref{lemma:computationofthederiveativeofsrqtg} that 
\beaa
\square_\g(r^{-1})=O(mr^{-4})+r^{-2}\dk^{\leq 1}\Ga_b, 
\eeaa
and
\beaa
&&\left|\int_{\MM_{r\geq R}(\tau_1,\tau_2)}r^{-1}\LL[\pmb\psi] - \int_{\MM_{r\geq R}(\tau_1,\tau_2)}r^{-1}\big(|\nab_{\pr_r}\pmb\psi|^2+r^{-2}|\nab_{\pr_{x^a}}\pmb\psi|^2\big)\right|\\
&\les& \int_{\MM_{r\geq R}(\tau_1,\tau_2)}r^{-1}\Big(|\nab_{\pr_\tau}\pmb\psi|^2+|\nab_{\pr_\tau}\psi|\big(|\nab_{\pr_r}\pmb\psi|+r^{-1}|\nab_{\pr_{x^a}}\pmb\psi|\big)\Big)\\
&&+\ep \EM_{r\geq R}[\pmb{\psi}](\tau_1,\tau_2)+R^{-1}\M_{r\geq R}[\pmb{\psi}](\tau_1,\tau_2),
\eeaa
 we deduce, for $R\geq 20m$ large enough, 
\beaa
\M_{0,r\geq R}[\pmb{\psi}](\tau_1,\tau_2)
&\les& \sup_{\tau\in(\tau_1, \tau_2)}\int_{\Si_{R/2,R}(\tau)}r^{-2}\chi_R|\dk^{\leq 1}\pmb\psi|^2+ \EF_{r\geq R}[\pmb{\psi}](\tau_1,\tau_2) \\
&&+\ep \EM_{r\geq R}[\pmb{\psi}](\tau_1,\tau_2) + R^2 \M_{ R/2, R}[\pmb{\psi}](\tau_1,\tau_2)+\int_{\MM_{r\geq R}(\tau_1,\tau_2)} r^{-1} \big|\nab_{\pr_{\tt}}\pmb{\psi}\big|^2\\
&&+\bigg|\int_{\MM_{r\geq R/2}(\tau_1,\tau_2)} \Re\Big(\ov{\big(-D_2r^{-2}\pmb{\psi}+D_1O(\ep r^{-1})\nab_{\pr_\tau}\pmb{\psi}+\F\big)}\chi_{R}r^{-1}\pmb{\psi}\Big)\bigg|
\eeaa
and hence, using \eqref{eq:weightedEMnearinf:inhomotensorialwave:simpifiedform:estF} for $\de_0=0$, and taking $\ep$ small enough and $R\geq 20m$ large enough, we deduce
\bea\lab{eq:Moraestimate:case1D1>0:0order:proof}
\M_{0,r\geq R}[\pmb{\psi}](\tau_1,\tau_2) &\les& \sup_{\tau\in(\tau_1, \tau_2)}\int_{\Si_{R/2,R}(\tau)}r^{-2}\chi_R|\dk^{\leq 1}\pmb\psi|^2+\EF_{r\geq R}[\pmb{\psi}](\tau_1,\tau_2)\nn\\
&&+ R^2 \M_{ R/2, R}[\pmb{\psi}](\tau_1,\tau_2) +\int_{\MM_{r\geq R}(\tau_1,\tau_2)} r^{-1} \big|\nab_{\pr_{\tt}}\pmb{\psi}\big|^2\nn\\
&& +D_2\int_{\MM(\tau_1, \tau_2)}r^{-3}|\pmb\psi|^2 + \int_{\MM_{ r\geq R/2}(\tau_1,\tau_2)} r |\N|^2.
\eea
In view of \eqref{eq:EFestimate:case1D1>0:0order:proof} and \eqref{eq:Moraestimate:case1D1>0:0order:proof}, and taking $\ep$ small enough and $R\geq 20m$ large enough, we deduce the desired estimate \eqref{eq:weightedEMFnearinf:D1>0:0order}.

Consider next the case $D_1=0$ and assume $\de\in(0, \frac{1}{3}]$. Similarly to the above discussions in the case $D_1>0$, we apply Proposition \ref{prop-app:stadard-comp-Psi} with $(X=\chi_{R}\pr_{\tt}, w=0)$ to equation \eqref{eq:weightedEMnearinf:inhomotensorialwave}, where $\chi_R$ is as given in \eqref{def:cutoffcuntionbetweenRandR/2},
 and we infer, for $R\geq 20m$ large enough and $\ep$ small enough, 
\bea
\lab{eq:EFestimate:case1D1=0:0order:proof}
\nn&& \sup_{\tau\in(\tau_1, \tau_2)}\int_{\Si_{R/2,R}(\tau)}r^{-2}\chi_R|\dk^{\leq 1}\pmb\psi|^2+ \EF_{r\geq R}[\pmb{\psi}](\tau_1,\tau_2)\\
&\les&\E_{r\geq R/2} [\pmb{\psi}](\tau_1) +\big(\ep+R^{-1}\big)\M_{r\geq R}[\pmb{\psi}](\tau_1,\tau_2) + R^2 \M_{ R/2, R}[\pmb{\psi}](\tau_1,\tau_2)\nn\\
&&+D_2\sup_{\tau\in[\tau_1, \tau_2]}\int_{\Si(\tau)}r^{-2}|\pmb\psi|^2+D_2\int_{\MM(\tau_1, \tau_2)}r^{-3}|\pmb\psi|^2+D_2\int_{\II_+(\tau_1, \tau_2)}r^{-2}|\pmb\psi|^2\nn\\
&&+ \bigg(\int_{\MM_{ r\geq R/2}(\tau_1,\tau_2)} r^{1+\de}|\N|^2\bigg)\Big(\M_{\de, r\geq R/2}[\pmb\psi](\tau_1,\tau_2)\Big)^{\frac{1}{2}}.
\eea
Next, applying Proposition \ref{prop-app:stadard-comp-Psi} with  $(X=X_{\de}, w=w_{\de})$\footnote{This is an adaption of the choice in \eqref{def:Xandw:improvedMorawetz:extendedRW:Kerrpert}  with the cutoff $\chi_{R_1}$ used therein replaced by $\chi_{R}$ defined in \eqref{def:cutoffcuntionbetweenRandR/2}.}, which is given by 
\bea
\lab{def:Xdeandwde:improvedMorawetz:generaltensorialwave:Kerrpert}
X_{\de}=2\mu \chi_{R}(1-m^{\de}r^{-\de}){\pr}_r^{\text{BL}}, \quad w_{\de}=4\mu \chi_{R} r^{-1} (1-m^{\de}r^{-\de})
\eea
with $\chi_R$ as given in \eqref{def:cutoffcuntionbetweenRandR/2}, and in view of the proof of \cite[Lemma 3.10]{MaSz24} which states\footnote{Compared to the proof of \cite[Lemma 3.10]{MaSz24}, the quantity $\mathcal{J}$ defined here contains one more term $-\frac{1}{2}X(V)|\pmb\psi|^2$, with $V=D_0|q|^{-2}$, which satisfies $-\frac{1}{2}X(V)|\pmb\psi|^2>0$ in $r\geq R$ for $R$ large enough and hence has the good sign.}
\beaa
\int_{\MM(\tau_1, \tau_2)}\mathcal{J}&:=&\int_{\MM(\tau_1, \tau_2)}\left(\frac 1 2 \QQ[\pmb\psi]  \c\piX -\frac{1}{2}X(V)|\pmb\psi|^2+\frac 12  w \LL[\pmb\psi] -\frac 1 4|\pmb\psi|^2   \square_\g  w\right) \nn\\
&\gtrsim& \de\M_{\de,r\geq R}[\pmb{\psi}](\tau_1,\tau_2)-O(\ep)\EM_{r\geq R}[\pmb{\psi}](\tau_1,\tau_2) - O(R^2) \M_{ R/2, R}[\pmb{\psi}](\tau_1,\tau_2),  
\eeaa
we deduce the following Morawetz  estimate near infinity, for $R\geq 20m$ large enough,
\bea
\lab{eq:Moraestimate:case1D1=0:0order:proof:step1}
\de\M_{\de,r\geq R}[\pmb{\psi}](\tau_1,\tau_2)
&\les& \sup_{\tau\in(\tau_1, \tau_2)}\int_{\Si_{R/2,R}(\tau)}r^{-2}\chi_R|\dk^{\leq 1}\pmb\psi|^2+ \EF_{r\geq R}[\pmb{\psi}](\tau_1,\tau_2)\nn\\ 
&&+\ep \EM_{r\geq R}[\pmb{\psi}](\tau_1,\tau_2) 
+ R^2 \M_{ R/2, R}[\pmb{\psi}](\tau_1,\tau_2)\nn\\
&&
+\bigg|\int_{\MM_{r\geq R/2}(\tau_1,\tau_2)} \Re\Big(\ov{\big(-D_2r^{-2}\pmb{\psi}+\F\big)}(\nab_{X_{\de}} \pmb{\psi} + w_{\de}\pmb{\psi})\Big)\bigg|\nn\\
&&
+\bigg|\int_{\MM_{r\geq R/2}(\tau_1,\tau_2)} \frac{k}{2}{}^{(X_{\de})}A_\nu\Im\Big(\pmb\psi\c\ov{\Ddot^{\nu}\pmb\psi}\Big)\bigg|.
\eea
In view of \eqref{eq:relationsbetweennullframeandcoordinatesframe2:moreprecise}, we have
\beaa
\pr_r^{\text{BL}}&=&O(1)\pr_{r} + O(1) \pr_{\tt} + O(r^{-2})\pr_{x^a}\nn\\
&=& O(1)e_4 + O(1) e_3 + O(r^{-1})  e_a,
\eeaa
hence, by Corollary \ref{cor:asymtpticbehavioroftheRdottermintensorialenergyidentity:larger}, 
we infer that in $r\geq R$,
\beaa
 \frac{k}{2}{}^{(X_{\de})}A_\nu \Im\Big(\pmb\psi\c\ov{\Ddot^{\nu}\pmb\psi}\Big)
= \big(O(r^{-3})+r^{-1}\Ga_b\big)\Im\Big(\pmb\psi\c\ov{\pr\pmb\psi }\Big)
\eeaa
and thus
\beaa
&&
\bigg|\int_{\MM_{r\geq R/2}(\tau_1,\tau_2)}  \frac{k}{2}{}^{(X_{\de})}A_\nu \Im\Big(\pmb\psi\c\ov{\Ddot^{\nu}\pmb\psi}\Big)\bigg|\nn\\
&\les&R^2 \M_{ R/2, R}[\pmb{\psi}](\tau_1,\tau_2)
+(R^{-1+\de}+\ep)\M_{\de, r\geq R}[\pmb\psi](\tau_1,\tau_2).
\eeaa
Substituting this back into \eqref{eq:Moraestimate:case1D1=0:0order:proof:step1}, 
we deduce
\beaa
\de\M_{\de,r\geq R}[\pmb{\psi}](\tau_1,\tau_2)
&\les& \sup_{\tau\in(\tau_1, \tau_2)}\int_{\Si_{R/2,R}(\tau)}r^{-2}\chi_R|\dk^{\leq 1}\pmb\psi|^2+ \EF_{r\geq R}[\pmb{\psi}](\tau_1,\tau_2)\\
&&+\ep \EM_{r\geq R}[\pmb{\psi}](\tau_1,\tau_2) + R^2 \M_{ R/2, R}[\pmb{\psi}](\tau_1,\tau_2)\\
&& +(R^{-1+\de}+\ep)\M_{\de, r\geq R}[\pmb\psi](\tau_1,\tau_2)\nn\\
&&+\bigg|\int_{\MM_{r\geq R/2}(\tau_1,\tau_2)} \Re\Big(\ov{\big(-D_2r^{-2}\pmb{\psi}+\F\big)}(\nab_{X_{\de}} \pmb{\psi} + w_{\de}\pmb{\psi})\Big)\bigg|,
\eeaa
and hence, using \eqref{eq:weightedEMnearinf:inhomotensorialwave:simpifiedform:estF} with $\de_0=\de$ and \eqref{eq:weightedEMnearinf:inhomotensorialwave:simpifiedform:estF1}, for $R\geq 20m$ large enough and $\ep$ small enough,
\bea
\lab{eq:Moraestimate:case1D1=0:0order:proof}
\de\M_{\de,r\geq R}[\pmb{\psi}](\tau_1,\tau_2)
&\les& \sup_{\tau\in(\tau_1, \tau_2)}\int_{\Si_{R/2,R}(\tau)}r^{-2}\chi_R|\dk^{\leq 1}\pmb\psi|^2+ \EF_{r\geq R}[\pmb{\psi}](\tau_1,\tau_2)\nn\\
&&+ R^2 \M_{R/2, R}[\pmb{\psi}](\tau_1,\tau_2)+D_2\sup_{\tau\in[\tau_1, \tau_2]}\int_{\Si(\tau)}r^{-2}|\pmb\psi|^2\nn\\
&&+D_2\int_{\MM(\tau_1, \tau_2)}r^{-3}|\pmb\psi|^2+D_2\int_{\II_+(\tau_1, \tau_2)}r^{-2}|\pmb\psi|^2\nn\\
 && + \bigg(\int_{\MM_{ r\geq R/2}(\tau_1,\tau_2)} r^{1+\de}|\N|^2\bigg)\Big(\M_{\de, r\geq R/2}[\pmb\psi](\tau_1,\tau_2)\Big)^{\frac{1}{2}}.
\eea
Combining the above two estimates \eqref{eq:EFestimate:case1D1=0:0order:proof} and \eqref{eq:Moraestimate:case1D1=0:0order:proof} then yields the desired estimate \eqref{eq:weightedEMFnearinf:D1=0:0order}. This concludes the proof of Lemma \ref{lem:zeroorderEMF:inhomogeneoustensorialwave}.
\end{proof}

Next, we show EMF estimates near infinity for weighted derivatives of solutions to a general class of tensorial wave equations. To this end, let us recall from Section \ref{sect:commutatorwithDalembertian} the following commutation formula for $\pmb{\psi} \in \sk_2$, see \eqref{commutator:rDDd2withrescaledwave},
\bea
\lab{commutator:rDDd2withrescaledwave:copy}
&&r\DDd_2(r^2\squared_2\pmb{\psi}) - r^2\squared_1(r\DDd_2\pmb{\psi}) \nn\\
&=& 3r\DDd_2\pmb{\psi} +O(m)\nab_3\dk^{\leq 1}\pmb{\psi}+O(mr^{-1})\dk^{\leq 2}\pmb{\psi}
+ \dk^{\leq 2} (r^2 \Ga_g \c \pmb{\psi})+r\Ga_b\c r^2\squared_2 \pmb{\psi} ,
\eea
 the following commutation formula for $\pmb{\psi} \in \sk_1$, see \eqref{commutator:rDDs2withrescaledwave},
\bea
\lab{commutator:rDDs2withrescaledwave:copy}
&&r\DDs_2 \, (r^2\squared_1\pmb{\psi}) - r^2\squared_2 (r\DDs_2\,\pmb{\psi}) \nn\\
&=&- 3 r\DDs_2\,\pmb{\psi}+O(m)\nab_3\dk^{\leq 1}\pmb{\psi}+O(mr^{-1})\dk^{\leq 2}\pmb{\psi}
+\dk^{\leq 2}(r^2 \Ga_g\c\pmb{\psi})+r\Ga_b\c r^2\squared_1 \pmb{\psi} ,
\eea 
the following commutator for $\pmb{\psi}\in\sk_k$, $k=1,2$, see \eqref{eq:comm-rnab4-squared2},
\bea\label{eq:comm-rnab4-squared2:copy}
[\nab_4 r, r^2\squared_k]\pmb{\psi} &=& -2r\nab_{\pr_{\tt}}\nab_4(r\pmb{\psi}) +O(m)\nab_3\dk^{\leq 1}\pmb{\psi}+O(mr^{-1})\dk^{\leq 2}\pmb{\psi}\nn\\
&&+O(m r^{-1})r^2\squared_k\pmb{\psi} +\dk^{\leq 2}(r^2\Ga_g\c\pmb{\psi}),
\eea
and the following commutator for $\pmb{\psi}\in\sk_k$, $k=1,2$, see \eqref{eq:comm-nabT-squared2},
\bea\label{eq:comm-nabT-squared2:copy}
[\nab_{\pr_{\tt}}, r^2\squared_k]\pmb{\psi} = O(m r^{-2})\dk^{\leq 1}\pmb{\psi} + \dk^{\leq 1}(\dk^{\leq 1}(r^2\Ga_g)\c\dk\psi).
\eea

\begin{lemma}[High-order EMF estimates near infinity for tensorial wave equations]
\lab{lem:EMFestinearinf:generaltensorialwave}
Let $\pmb{\psi}\in\sk_{2}(\mathbb{C})$ and $\pmb{N}\in \sk_{2}(\mathbb{C})$ satisfy
\bea
\lab{eq:weightedEMnearinf:generaltensorialwave}
\squared_2\pmb{\psi} - D_0|q|^{-2} \pmb{\psi} =D_1r^{-1}\nab_{\pr_{\tt}}\pmb{\psi} + O(r^{-2})\nab_3 \pmb{\psi} + O(r^{-3})\dk^{\leq 1}\pmb{\psi} + \dk^{\leq 1}\Ga_g \dk^{\leq 1}\pmb{\psi}+\pmb{N},
\eea
with the constants $D_0> 0$ and $D_1\geq 0$. If $D_1>0$, then we have, for any $\tau_1<\tau_2$, $\reg\leq 14$ and $R\geq 20m$ large enough,
\bea
\lab{eq:weightedEMFnearinf:D1>0}
\EMF_{0,r\geq R}^{(\reg)}[\pmb{\psi}](\tau_1,\tau_2) &\les& \E_{r\geq R/2}^{(\reg)}[\pmb{\psi}](\tau_1) +R^2 \M_{R/2, R}^{(\reg)}[\pmb{\psi}](\tau_1,\tau_2)\nn\\
&&+ \int_{\MM_{r\geq R/2}(\tau_1,\tau_2)} r |\dk^{\leq \reg} \pmb{N}|^2.
\eea
If $D_1=0$, then we have, for any $\tau_1<\tau_2$, $\de\in(0, \frac{1}{3}]$ and $R\geq 20m$ large enough,
\bea
\lab{eq:weightedEMFnearinf:D1=0}
\EMF_{\de,r\geq R}^{(\reg)}[\pmb{\psi}](\tau_1,\tau_2) &\les& \E_{r\geq R/2}^{(\reg)}[\pmb{\psi}](\tau_1) + R^2 \M_{R/2, R}^{(\reg)}[\pmb{\psi}](\tau_1,\tau_2)\nn\\
&&+ \int_{\MM_{r\geq R/2}(\tau_1,\tau_2)} r^{1+\de} |\dk^{\leq \reg} \pmb{N}|^2.
\eea
\end{lemma}

\begin{proof}
To begin with, we shall derive wave equations for weighted derivatives of $\pmb{\psi}$ using the commutation relations \eqref{commutator:rDDd2withrescaledwave:copy} \eqref{commutator:rDDs2withrescaledwave:copy}  \eqref{eq:comm-rnab4-squared2:copy} \eqref{eq:comm-nabT-squared2:copy}. For this purpose, let us multiply both sides of \eqref{eq:weightedEMnearinf:generaltensorialwave} by $r^2$ to infer
\bea
\lab{eq:weightedEMnearinf:generaltensorialwave:rescaled}
r^2\squared_2\pmb{\psi} - D_0 \pmb{\psi} =D_1r \nab_{\pr_{\tt}}\pmb{\psi} + O(1)\nab_3 \pmb{\psi} + O(r^{-1})\dk^{\leq 1}\pmb{\psi} +r^2 \dk^{\leq 1}\Ga_g \dk^{\leq 1}\pmb{\psi}+r^2\pmb{N},
\eea
and recall from Section \ref{sect:HorizontalHodgeopes} the following set of high-order weighted covariant derivatives, for any $\reg\in\mathbb{N}$, see  \eqref{def:highorderweightedderivatives:containhorizontal},
\bea
\lab{def:highorderweightedderivatives:containhorizontal:copy}
\widecheck{\dk}^{\reg}:= {\left\{\dkb^{\reg_1}(\nab_{\pr_{\tt}})^{\reg_2}(\nab_{4}r)^{\reg_3} ,  \,\,\, \reg_1+\reg_2+\reg_3=\reg\right\}},
\eea
where $\dkb^{\reg_1}$ are weighted horizontal Hodge operators given by, see \eqref{eq:definitiondkbpowerj},
\begin{equation}\lab{eq:definitiondkbpowerj:copy}
\begin{split}
\dkb^{\reg_1}:=&(r\DDs_2\,\, r\DDd_2)^{\frac{{\reg_1}}{2}}, \quad \text{if ${\reg_1}$ is even,}\\
\dkb^{\reg_1}:=&r\DDd_2(r\DDs_2\,\, r\DDd_2)^{\frac{{\reg_1}-1}{2}}, \quad \text{if ${\reg_1}$ is odd.}
\end{split}
\end{equation}

Depending on the parity of $\reg_1$, we separate into two cases:
\begin{itemize}
\item {\bf Case 1:} $\reg_1$ is even.
\item {\bf Case 2:} $\reg_1$ is  odd. 
\end{itemize}
We claim that, in Case 1, we have the following
 general form for the wave equation satisfied by $\widecheck\dk^{\reg}\pmb{\psi}\in\sk_2(\mathbb{C})$
\bea
\lab{eq:commutedequation:generalwave:DDdequalDDs}
r^2\squared_2\widecheck\dk^{\reg}\pmb{\psi} - D_0 \widecheck\dk^{\reg}\pmb{\psi} &=& 
\sum_{\reg'\leq\reg, \reg_3'\leq\reg_3}(D_1+2\reg_3')r \nab_{\pr_{\tt}}\widecheck\dk^{\reg'}\pmb{\psi}  + O(1)\nab_3 \widecheck\dk^{\leq \reg}\pmb{\psi} + O(r^{-1})\widecheck\dk^{\leq \reg+1}\pmb{\psi} \nn\\
&&+r^2\widecheck\dk^{\leq \reg+1}\Ga_g \widecheck\dk^{\leq \reg+1}\pmb{\psi} +\sum_{\reg'\leq\reg, \reg_3'\leq\reg_3}(D_1+2\reg_3')r^2\widecheck\dk^{\leq \reg'+1}\Ga_b \nab_{\pr_\tau}\widecheck\dk^{\reg'}\pmb{\psi}\nn\\
&& +O(r^2)\widecheck\dk^{\leq \reg}\pmb{N},  \quad \text{if $\reg_1$ is even},
 \eea
 and in Case 2, we have the following
 general form for the wave equation satisfied by $\widecheck\dk^{\reg}\pmb{\psi}\in\sk_1(\mathbb{C})$
\bea
\lab{eq:commutedequation:generalwave:DDdequalDDsplus1}
r^2\squared_1\widecheck\dk^{\reg}\pmb{\psi} - D_0 \widecheck\dk^{\reg}\pmb{\psi} &=& -3\widecheck\dk^{\reg}\pmb{\psi} +\sum_{\reg'\leq\reg, \reg_3'\leq\reg_3}(D_1+2\reg_3')r \nab_{\pr_{\tt}}\widecheck\dk^{\reg'}\pmb{\psi} + O(1)\nab_3 \widecheck\dk^{\leq \reg}\pmb{\psi}\nn\\
&& + O(r^{-1})\widecheck\dk^{\leq \reg+1}\pmb{\psi} +r^2\widecheck\dk^{\leq \reg+1}\Ga_g \widecheck\dk^{\leq \reg+1}\pmb{\psi}\nn\\
&& +\sum_{\reg'\leq\reg, \reg_3'\leq\reg_3}(D_1+2\reg_3')r^2\widecheck\dk^{\leq \reg'+1}\Ga_b \nab_{\pr_\tau}\widecheck\dk^{\reg'}\pmb{\psi}\nn\\
&& +O(r^2)\widecheck\dk^{\leq \reg}\pmb{N}, \quad \text{if $\reg_1$ is odd} .
\eea

We prove the above equations \eqref{eq:commutedequation:generalwave:DDdequalDDs} and \eqref{eq:commutedequation:generalwave:DDdequalDDsplus1} by induction on the value of $\reg$. First, commuting $\nab_{\pr_{\tt}}$ with the above two equations, and using the commutation relation \eqref{eq:comm-nabT-squared2:copy}, the above equations \eqref{eq:commutedequation:generalwave:DDdequalDDs} and \eqref{eq:commutedequation:generalwave:DDdequalDDsplus1} hold for $\nab_{\pr_{\tt}}\widecheck\dk^{\reg}\pmb{\psi}$, i.e., for $(\reg_1, \reg_2+1, \reg_3)$. Secondly, 
commuting $\nab_4 r$ with the above two equations, and in view of the commutator relation \eqref{eq:comm-rnab4-squared2:copy}, using also the commutators in Corollaries \ref{cor:corollaryofLemmacomm} and \ref{cor:commutatorweightedderivativesrnabandnab4rwithnabprtau}, the above equations \eqref{eq:commutedequation:generalwave:DDdequalDDs} and \eqref{eq:commutedequation:generalwave:DDdequalDDsplus1} hold for $\nab_{4} r \widecheck\dk^{\reg}\pmb{\psi}$, i.e., for $(\reg_1, \reg_2, \reg_3+1)$. Thirdly, we commute with the weighted horizontal Hodge derivatives, which corresponds to $(\reg_1+1, \reg_2, \reg_3)$, and separate the proof based on the above two cases\footnote{Note that the before to last term in \eqref{eq:commutedequation:generalwave:DDdequalDDs} and \eqref{eq:commutedequation:generalwave:DDdequalDDsplus1}, involving $r^2\widecheck\dk^{\leq \reg+1}\Ga_b \nab_{\pr_\tau}\widecheck\dk^{\leq \reg}\pmb{\psi}$, is generated by the commutation of $\dkb$ with $(D_1+2\reg_3)r \nab_{\pr_{\tt}}\widecheck\dk^{\leq \reg}\pmb{\psi}$ in view of the first commutator estimate in Corollary \ref{cor:commutatorweightedderivativesrnabandnab4rwithnabprtau}.}:
\begin{itemize}
\item In Case 1, we commute equation \eqref{eq:commutedequation:generalwave:DDdequalDDs} with $r\DDd_2$, and in view of the commutation relation \eqref{commutator:rDDd2withrescaledwave:copy}, using also the commutators in Corollaries \ref{cor:corollaryofLemmacomm} and \ref{cor:commutatorweightedderivativesrnabandnab4rwithnabprtau}, the equation for $r\DDd_2\widecheck\dk^{\reg}\pmb{\psi}$ is in the form of \eqref{eq:commutedequation:generalwave:DDdequalDDsplus1}. 

\item In Case 2, we commute equation \eqref{eq:commutedequation:generalwave:DDdequalDDsplus1} with $r\DDs_2$ and use the commutation relation \eqref{commutator:rDDs2withrescaledwave:copy}, as well as the commutators in Corollaries \ref{cor:corollaryofLemmacomm} and \ref{cor:commutatorweightedderivativesrnabandnab4rwithnabprtau}, to find that the equation for $r\DDs_2\widecheck\dk^{\reg}\pmb{\psi}$ is in the form of \eqref{eq:commutedequation:generalwave:DDdequalDDs}. 
\end{itemize}
These together prove the equations \eqref{eq:commutedequation:generalwave:DDdequalDDs} and \eqref{eq:commutedequation:generalwave:DDdequalDDsplus1} for $\widecheck{\dk}^{\reg+1}$ and hence, using also the fact that $\reg=0$ holds true in view of \eqref{eq:weightedEMnearinf:generaltensorialwave:rescaled}, conclude the proof of equations \eqref{eq:commutedequation:generalwave:DDdequalDDs} and \eqref{eq:commutedequation:generalwave:DDdequalDDsplus1} for general $\reg\in \mathbb{N}$ by iteration.

We now rely on \eqref{eq:commutedequation:generalwave:DDdequalDDs} and \eqref{eq:commutedequation:generalwave:DDdequalDDsplus1} to derive EMF estimates for $\widecheck{\dk}^{\reg}\pmb{\psi}$. We start first with the case $D_1+2\reg_3>0$. Noticing that:
\begin{itemize}
\item \eqref{eq:commutedequation:generalwave:DDdequalDDs} satisfies \eqref{eq:weightedEMnearinf:inhomotensorialwave} with $D_1$ replaced by $D_1+2\reg_3>0$, $D_2=0$, and $r^2\N$ replaced by  
\bea\lab{eq:mbNreplacedbyhigherorderanalogaftercommutatorestimates:aux}
\nn&&\sum_{\reg'\leq\reg, \reg_3'\leq\reg_3}(D_1+2\reg_3')r \nab_{\pr_{\tt}}\widecheck\dk^{\reg'}\pmb{\psi}  + O(1)\nab_3 \widecheck\dk^{<\reg}\pmb{\psi} + O(r^{-1})\widecheck\dk^{\leq \reg}\pmb{\psi} +r^2\widecheck\dk^{\leq \reg}\Ga_g \widecheck\dk^{\leq \reg}\pmb{\psi}\\
&& +\sum_{\reg'<\reg, \reg_3'\leq\reg_3}(D_1+2\reg_3')r^2\widecheck\dk^{\leq \reg'+1}\Ga_b \nab_{\pr_\tau}\widecheck\dk^{\reg'}\pmb{\psi} +O(r^2)\widecheck\dk^{\leq\reg}\pmb{N},  
\eea
\item and \eqref{eq:commutedequation:generalwave:DDdequalDDsplus1} satisfies \eqref{eq:weightedEMnearinf:inhomotensorialwave} with $D_1$ replaced by $D_1+2\reg_3>0$, $D_2=3$, and $r^2\N$ replaced by  \eqref{eq:mbNreplacedbyhigherorderanalogaftercommutatorestimates:aux},
\end{itemize}
we may apply \eqref{eq:weightedEMFnearinf:D1>0:0order} which yields for $\reg=0$ (where $D_2=0$)
\bea\lab{eq:aux:higherorderenergyMorawetznearinifinityforreg=0:CaseD1plus2reg3pos}
\EMF_{0,r\geq R}[\pmb{\psi}](\tau_1,\tau_2) \les \E_{r\geq R/2}[\pmb{\psi}](\tau_1) + R^2 \M_{R/2, R}[\pmb{\psi}](\tau_1,\tau_2)  + \int_{\MM_{r\geq R/2}(\tau_1,\tau_2)} r |\N|^2,
\eea
and for $1\leq\reg\leq 14$,  making use of the estimate \eqref{eq:equivalencerelationonweightednorms},
\bea\lab{eq:aux:higherorderenergyMorawetznearinifinityforregpos:CaseD1plus2reg3pos}
\nn\EMF_{0,r\geq R}[\dk^{\reg}\pmb{\psi}](\tau_1,\tau_2) &\les& \E^{(\reg)}_{r\geq R/2}[\pmb{\psi}](\tau_1) + R^2 \M^{(\reg)}_{R/2, R}[\pmb{\psi}](\tau_1,\tau_2) +\EM^{(\reg-1)}_{r\geq R}[\pmb{\psi}](\tau_1,\tau_2)\\
\nn&& +\sum_{\reg'<\reg, \reg_3'\leq\reg_3}(D_1+2\reg_3')\M_{0,r\geq R}[\dk^{\reg'}\pmb{\psi}](\tau_1,\tau_2) \\
&& + \int_{\MM_{r\geq R/2}(\tau_1,\tau_2)} r|\dk^{\leq\reg}\N|^2.
\eea

Next, we consider the case $D_1+2\reg_3=0$ which implies $D_1+2\reg_3'=0$ for $\reg_3'\leq\reg_3$ and in particular $D_1=0$. Noticing that:
\begin{itemize}
\item \eqref{eq:commutedequation:generalwave:DDdequalDDs} satisfies \eqref{eq:weightedEMnearinf:inhomotensorialwave} with $D_1$ replaced by $D_1=0$, $D_2=0$, and $r^2\N$ replaced by  
\bea\lab{eq:mbNreplacedbyhigherorderanalogaftercommutatorestimates:aux1}
O(1)\nab_3 \widecheck\dk^{<\reg}\pmb{\psi} + O(r^{-1})\widecheck\dk^{\leq \reg}\pmb{\psi} +r^2\widecheck\dk^{\leq \reg}\Ga_g \widecheck\dk^{\leq \reg}\pmb{\psi}+O(r^2)\widecheck\dk^{\leq\reg}\pmb{N},  
\eea
\item and \eqref{eq:commutedequation:generalwave:DDdequalDDsplus1} satisfies \eqref{eq:weightedEMnearinf:inhomotensorialwave} with $D_1$ replaced by $D_1=0$, $D_2=3$, and $r^2\N$ replaced by  \eqref{eq:mbNreplacedbyhigherorderanalogaftercommutatorestimates:aux1},
\end{itemize}
we may apply \eqref{eq:weightedEMFnearinf:D1=0:0order} which yields for $\reg=0$ (where $D_2=0$)
\bea\lab{eq:aux:higherorderenergyMorawetznearinifinityforreg=0:CaseD1=0}
\EMF_{\de,r\geq R}[\pmb{\psi}](\tau_1,\tau_2) \les \E_{r\geq R/2}[\pmb{\psi}](\tau_1) + R^2 \M_{R/2, R}[\pmb{\psi}](\tau_1,\tau_2)  + \int_{\MM_{r\geq R/2}(\tau_1,\tau_2)} r^{1+\de}|\N|^2,
\eea
and for $1\leq\reg\leq 14$, making use of the estimate \eqref{eq:equivalencerelationonweightednorms},
\bea\lab{eq:aux:higherorderenergyMorawetznearinifinityforregpos:CaseD1=0}
\nn\EMF_{\de,r\geq R}[\dk^{\reg}\pmb{\psi}](\tau_1,\tau_2) &\les& \E^{(\reg)}_{r\geq R/2}[\pmb{\psi}](\tau_1) + R^2 \M^{(\reg)}_{R/2, R}[\pmb{\psi}](\tau_1,\tau_2) +\EM^{(\reg-1)}_{r\geq R}[\pmb{\psi}](\tau_1,\tau_2)\\
&& + \int_{\MM_{r\geq R/2}(\tau_1,\tau_2)} r^{1+\de}|\dk^{\leq\reg}\N|^2,
\eea
where $\de\in(0, \frac{1}{3}]$.

We are now ready to conclude. In the case $D_1>0$, we always have $D_1+\reg_3>0$ so that \eqref{eq:aux:higherorderenergyMorawetznearinifinityforreg=0:CaseD1plus2reg3pos} holds for $\reg=0$ and \eqref{eq:aux:higherorderenergyMorawetznearinifinityforregpos:CaseD1plus2reg3pos} 
 holds for $1\leq\reg\leq 14$. This immediately implies \eqref{eq:weightedEMFnearinf:D1>0} by iteration. Also, in the case $D_1=0$, \eqref{eq:aux:higherorderenergyMorawetznearinifinityforreg=0:CaseD1=0} holds for $\reg=0$, and then:
\begin{itemize}
\item for all $1\leq\reg\leq 14$ with $\reg_3=0$, we apply \eqref{eq:aux:higherorderenergyMorawetznearinifinityforregpos:CaseD1=0},
\item for all $1\leq\reg\leq 14$ with $\reg_3>0$, we apply \eqref{eq:aux:higherorderenergyMorawetznearinifinityforregpos:CaseD1plus2reg3pos}.
\end{itemize}
We infer by iteration for all $\reg\leq 14$, for all $\de\in(0, \frac{1}{3}]$,
\beaa
\nn &&\EMF_{\de,r\geq R}[\dk^{\reg}\pmb{\psi}](\tau_1,\tau_2)+(1-\de_{0k_3})\sum_{\reg'\leq\reg, 0<\reg_3'\leq\reg_3}\M_{0,r\geq R}[\dk^{\reg'}\pmb{\psi}](\tau_1,\tau_2)\\
&\les& \E^{(\reg)}_{r\geq R/2}[\pmb{\psi}](\tau_1) + R^2 \M^{(\reg)}_{R/2, R}[\pmb{\psi}](\tau_1,\tau_2) + \int_{\MM_{r\geq R/2}(\tau_1,\tau_2)} r^{1+\de}|\dk^{\leq\reg}\N|^2,
\eeaa
which yields in particular \eqref{eq:weightedEMFnearinf:D1=0}. This concludes the proof of Lemma \ref{lem:EMFestinearinf:generaltensorialwave}.
\end{proof}

In the rest of this section, we provide the proof of Proposition \ref{prop:EnergyMorawetznearinfinitytensorialTeuk} based on Lemma \ref{lem:EMFestinearinf:generaltensorialwave}.

\begin{proof}[Proof of Proposition \ref{prop:EnergyMorawetznearinfinitytensorialTeuk}]
Recall the wave equations \eqref{eq:TensorialTeuSysandlinearterms:rescaleRHScontaine2:general:Kerrperturbation} for $\{\phis{p}\}_{p=0,1,2}$, and rewrite these wave equations in the form of the tensorial wave equation \eqref{eq:weightedEMnearinf:generaltensorialwave} for a triplet of $\sk_2(\mathbb{C})$ tensors $(\phis{0}, \phis{1}, \phis{2})$ as follows 
\bsub
\lab{eq:TensorialTeuSys:rescaleRHS:pastep=012s-2}
\bea
\lab{eq:TensorialTeuSys:rescaleRHS:pastep=0s-2}
\squared_2 \phis{0}- 2|q|^{-2}\phis{0} &=& O(r^{-2})\nab_3\phis{0}+ O(r^{-3})\dk^{\leq 1}\phis{0} +\N^{(0)}_{s}, \\
\lab{eq:TensorialTeuSys:rescaleRHS:pastep=1s-2}
\squared_2 \phis{1}- 4|q|^{-2}\phis{1} &=&O(r^{-2})\nab_3\phis{1}+ O(r^{-3})\dk^{\leq 1}\phis{1} +\N^{(1)}_{s}, \\
\lab{eq:TensorialTeuSys:rescaleRHS:pastep=2s-2}
\squared_2 \phis{2}- 4|q|^{-2}\phis{2} &=& O(r^{-2})\nab_3\phis{2}+ O(r^{-3})\dk^{\leq 1}\phis{2} +\N^{(2)}_{s}, 
\eea
\esub
where we have, using also \eqref{eq:formofprphi+aprtt:horizontal}, 
\bsub
\lab{formulaofN-2012:nearinf}
\begin{align}
\N^{(0)}_{s}:={}& O(r^{-3})\phis{1}+ \N_{W,s}^{(0)},\\
\N^{(1)}_{s}:={}& O(r^{-3})\phis{2}+ O(r^{-2})(r\nab)^{\leq 1}\phis{0} +O(r^{-2})\nab_{\pr_{\tau}} \phis{0} +\Ga_g\c \dk\phis{0}+\N_{W,s}^{(1)},\\
\N^{(2)}_{s}:={}&O(r^{-2})(r\nab)^{\leq 1}\phis{1} +O(r^{-2})\nab_{\pr_{\tau}} \phis{1} +\Ga_g\c \dk\phis{1} + O(r^{-2})\phis{0}
+\N_{W,s}^{(2)}.
\end{align}
\esub
Now, applying  the estimate \eqref{eq:weightedEMFnearinf:D1=0} to the above system of equations \eqref{eq:TensorialTeuSys:rescaleRHS:pastep=012s-2}, with $\{\N^{(p)}_{s}\}_{p=0,1,2}$  given as in \eqref{formulaofN-2012:nearinf}, we deduce 
\bsub\lab{eq:EMFnearinf:pm2:012:middle}
\bea
\lab{eq:EMFnearinf:pm20:middle}
&&\EMF_{\de,r\geq R}^{(\reg)}[\phis{0}](\tau_1,\tau_2) \nn\\
&\les& \E_{r\geq R/2}^{(\reg)}[\phis{0}](\tau_1) + R^2 \M_{R/2, R}^{(\reg)}[\phis{0}](\tau_1,\tau_2)
+
R^{-2+\de}\M^{(\reg)}_{\de, r\geq R/2}[\phis{1}](\tau_1,\tau_2)\nn\\
&&
+ \int_{\MM_{r\geq R/2}(\tau_1,\tau_2)} r^{1+\de} |\dk^{\leq \reg} \N^{(0)}_{W,s}|^2,
\eea
\bea
\lab{eq:EMFnearinf:pm21:middle}
&&\EMF_{\de,r\geq R}^{(\reg)}[\phis{1}](\tau_1,\tau_2) \nn\\
&\les& \E_{r\geq R/2}^{(\reg)}[\phis{1}](\tau_1) + R^2 \M_{R/2, R}^{(\reg)}[\phis{1}](\tau_1,\tau_2)
+ \int_{\MM_{r\geq R/2}(\tau_1,\tau_2)} r^{1+\de} |\dk^{\leq \reg} \N^{(1)}_{s}|^2\nn\\
&\les& \E_{r\geq R/2}^{(\reg)}[\phis{1}](\tau_1) + R^2 \M_{R/2, R}^{(\reg)}[\phis{1}](\tau_1,\tau_2)\nn\\
&&
+R^{-1+\de}\EM^{(\reg)}_{\de, r\geq R/2}[\phis{0}](\tau_1,\tau_2)
+R^{-2 +\de}\M^{(\reg)}_{\de, r\geq R/2}[\phis{2}](\tau_1,\tau_2)
\nn\\
&&
+ \int_{\MM_{r\geq R/2}(\tau_1,\tau_2)} r^{1+\de} |\dk^{\leq \reg} \N^{(1)}_{W,s}|^2+\Ao_{r\geq R/2}[r^{\frac{\de}{2}} (r\nab)^{\leq 1} \dk^{\leq \reg}\phis{0}](\tau_1,\tau_2),
\eea
\bea
\lab{eq:EMFnearinf:pm22:middle}
&&\EMF_{\de,r\geq R}^{(\reg)}[\phis{2}](\tau_1,\tau_2) \nn\\
&\les& \E_{r\geq R/2}^{(\reg)}[\phis{2}](\tau_1) + R^2 \M_{R/2, R}^{(\reg)}[\phis{2}](\tau_1,\tau_2)
+ \int_{\MM_{r\geq R/2}(\tau_1,\tau_2)} r^{1+\de} |\dk^{\leq \reg} \N^{(2)}_{s}|^2\nn\\
&\les& \E_{r\geq R/2}^{(\reg)}[\phis{2}](\tau_1) + R^2 \M_{R/2, R}^{(\reg)}[\phis{2}](\tau_1,\tau_2)
\nn\\
&&+R^{-1+\de}\EM^{(\reg)}_{\de, r\geq R/2}[\phis{1}](\tau_1,\tau_2)+ \int_{\MM_{r\geq R/2}(\tau_1,\tau_2)} r^{1+\de} |\dk^{\leq \reg} \N^{(2)}_{W,s}|^2\nn\\
&&+\Ao_{r\geq R/2}[r^{\frac{\de}{2}} (r\nab)^{\leq 1} \dk^{\leq \reg}\phis{1}](\tau_1,\tau_2)+\Ao_{r\geq R/2}[r^{\frac{\de}{2}}  \dk^{\leq \reg}\phis{0}](\tau_1,\tau_2).
\eea
\esub
Adding \eqref{eq:EMFnearinf:pm20:middle}, \eqref{eq:EMFnearinf:pm21:middle} and \eqref{eq:EMFnearinf:pm22:middle}, and taking $R\geq 20m$ large enough so that the terms on the RHS whose coefficients  are given by $R^{-1+\de}$ and $R^{-2+\de}$ are absorbed by the LHS, we obtain 
\begin{align}
\lab{eq:EMFnearinf:pm2:012:middle:conditionalestimate}
& \sum_{p=0,1,2}\EMF_{\de,r\geq R}^{(\reg)}[\phis{p}](\tau_1,\tau_2)\nn\\
\les{}& \sum_{p=0,1,2}\bigg(\E_{r\geq R/2}^{(\reg)}[\phis{p}](\tau_1) 
+ R^2 \M_{R/2, R}^{(\reg)}[\phis{p}](\tau_1,\tau_2)+\int_{\MM_{r\geq R/2}(\tau_1,\tau_2)} r^{1+\de} |\dk^{\leq \reg} \N^{(p)}_{W,s}|^2\bigg)\nn\\
&+ \sum_{p=0,1}\Ao_{r\geq R}[r^{\frac{\de}{2}} (r\nab)^{\leq 1} \dk^{\leq \reg}\phis{p}](\tau_1,\tau_2).
\end{align}

It remains to estimate the last term in the above inequality. To this end, in view of the estimates \eqref{eq:improvedestiforangularderiofkerrpert:prop:highorder:largeR}, and using
\beaa
\A_{r\geq R/2}[r^{-\frac{3\de}{2}}\dk^{\leq \reg} \pmb\phi_{s}](\tau_1,\tau_2)\les R^{-2\de}\sum_{p=0}^2{\EMF}^{(\reg)}_{\de, r\geq R/2}[\phis{p}](\tau_1,\tau_2), \quad \forall \, \reg\leq 14, \,\, 0<\de\leq \frac{1}{3},
\eeaa
we deduce, for $\reg\leq 14$ and $0<\de\leq \frac{1}{3}$, 
\beaa
&&\Ao_{r\geq R}[{r^{\frac{\de}{2}}(r\nab)^{\leq 1}}\dk^{\leq \reg}\phis{0}](\tau_1,\tau_2)\nn\\
&\les &\Big({\EMF}^{(\reg)}_{\de, r\geq R/2}[\phis{0}](\tau_1,\tau_2)\Big)^{\frac{1}{2}}\bigg(R^{-2\de}\sum_{p=0,1,2}{\EMF}^{(\reg)}_{\de, r\geq R/2}[\phis{p}](\tau_1,\tau_2)\bigg)^{\frac{1}{2}} \nn\\
&&+R^{-2\de}\sum_{p=0,1,2}{\EMF}^{(\reg)}_{\de, r\geq R/2}[\phis{p}](\tau_1,\tau_2)\nn\\
&&
+\int_{\MM_{r\geq R/2}(\tau_1,\tau_2)}\Big(r^{-1+\de}|\dk^{\leq \reg} \N_{W,s}^{(0)}|+r^{-2+\de}|\dk^{\leq \reg} \nab_{X_s}\N^{(0)}_{T,s}|\Big)|\dk^{\leq \reg} \phis{0}|,\\
&&\Ao_{r\geq R}[r^{\frac{\de}{2}}(r\nab)^{\leq 1}\dk^{\leq \reg} \phis{1}](\tau_1,\tau_2)\nn\\
&\les &\bigg(\sum_{p=0,1}{\EMF}_{\de, r\geq R/2}^{(\reg)}[\phis{p}](\tau_1,\tau_2)\bigg)^{\frac{1}{2}}\bigg(R^{-2\de}\sum_{p=0,1,2}{\EMF}^{(\reg)}_{\de, r\geq R/2}[\phis{p}](\tau_1,\tau_2)\bigg)^{\frac{1}{2}} \nn\\
&&+R^{-2\de}\sum_{p=0,1,2}{\EMF}^{(\reg)}_{\de, r\geq R/2}[\phis{p}](\tau_1,\tau_2)\nn\\
&&
+\int_{\MM_{r\geq R/2}(\tau_1,\tau_2)}\Big(r^{-1+\de}|\dk^{\leq \reg} \N_{W,s}^{(1)}|+r^{-2+\de}|\dk^{\leq \reg} \nab_{X_s}\N^{(1)}_{T,s}|\Big)|\dk^{\leq \reg} \phis{1}|\nn\\
&&+\Ao_{r\geq R/2}[r^{\frac{\de}{2}}(r\nab)^{\leq 1}\dk^{\leq \reg} \phis{0}](\tau_1,\tau_2),
\eeaa
which yields
\bea
&&\sum_{p=0,1}\Ao_{r\geq R}[r^{\frac{\de}{2}} (r\nab)^{\leq 1} \dk^{\leq \reg}\phis{p}](\tau_1,\tau_2)\nn\\
&\les&\bigg(\sum_{p=0,1}{\EMF}_{\de, r\geq R/2}^{(\reg)}[\phis{p}](\tau_1,\tau_2)\bigg)^{\frac{1}{2}}\bigg(R^{-2\de}\sum_{p=0,1,2}{\EMF}^{(\reg)}_{\de, r\geq R/2}[\phis{p}](\tau_1,\tau_2)\bigg)^{\frac{1}{2}} \nn\\
&&
+\sum_{p=0,1}\int_{\MM_{r\geq R/2}(\tau_1,\tau_2)}\Big(r^{-1+\de}|\dk^{\leq \reg} \N_{W,s}^{(p)}|+r^{-2+\de}|\dk^{\leq \reg} \nab_{X_s}\N^{(p)}_{T,s}|\Big)|\dk^{\leq \reg} \phis{p}|\nn\\
&&+R^{-2\de}\sum_{p=0,1,2}{\EMF}^{(\reg)}_{\de, r\geq R/2}[\phis{p}](\tau_1,\tau_2)\nn\\
&\les&\bigg(\sum_{p=0,1}{\EMF}_{\de, r\geq R/2}^{(\reg)}[\phis{p}](\tau_1,\tau_2)\bigg)^{\frac{1}{2}}\bigg(R^{-2\de}\sum_{p=0,1,2}{\EMF}^{(\reg)}_{\de, r\geq R/2}[\phis{p}](\tau_1,\tau_2)\bigg)^{\frac{1}{2}}\nn\\
&&
+\sum_{p=0,1}\bigg(\int_{\MM_{r\geq R/2}(\tau_1,\tau_2)}\Big(r^{1+\de}|\dk^{\leq \reg} \N_{W,s}^{(p)}|^2+r^{-1+\de}|\dk^{\leq \reg} \nab_{X_s}\N^{(p)}_{T,s}|^2\Big)\bigg)^{\frac{1}{2}}\nn\\
&& \times\Big(\Ao_{r\geq R/2}[r^{\frac{\de}{2}}\dk^{\leq \reg} \phis{p}](\tau_1,\tau_2)\Big)^{\frac{1}{2}}+R^{-2\de}\sum_{p=0,1,2}{\EMF}^{(\reg)}_{\de, r\geq R/2}[\phis{p}](\tau_1,\tau_2)
\eea
and hence
\bea
&&\sum_{p=0,1}\Ao_{r\geq R}[r^{\frac{\de}{2}} (r\nab)^{\leq 1} \dk^{\leq \reg}\phis{p}](\tau_1,\tau_2)\nn\\
&\les&\bigg(\sum_{p=0,1}{\EMF}_{\de, r\geq R/2}^{(\reg)}[\phis{p}](\tau_1,\tau_2)\bigg)^{\frac{1}{2}}\bigg(R^{-2\de}\sum_{p=0,1,2}{\EMF}^{(\reg)}_{\de, r\geq R/2}[\phis{p}](\tau_1,\tau_2)\bigg)^{\frac{1}{2}}\nn\\
&&
+\sum_{p=0,1}\int_{\MM_{r\geq R/2}(\tau_1,\tau_2)}\Big(r^{1+\de}|\dk^{\leq \reg} \N_{W,s}^{(p)}|^2+r^{-1+\de}|\dk^{\leq \reg+1} \N^{(p)}_{T,s}|^2\Big)\nn\\
&&+R^{-2\de}\sum_{p=0,1,2}{\EMF}^{(\reg)}_{\de, r\geq R/2}[\phis{p}](\tau_1,\tau_2)+R^{\de}\sum_{p=0,1}\M^{(\reg)}_{R/2, R}[ \phis{p}](\tau_1,\tau_2).
\eea 
Plugging this estimate back into \eqref{eq:EMFnearinf:pm2:012:middle:conditionalestimate} to control the last term on the RHS of \eqref{eq:EMFnearinf:pm2:012:middle:conditionalestimate}, we infer, for $\reg\leq 14$, 
\begin{align*}
& \sum_{p=0,1,2}\EMF_{\de,r\geq R}^{(\reg)}[\phis{p}](\tau_1,\tau_2)+\sum_{p=0,1}\Ao_{r\geq R}[r^{\frac{\de}{2}} (r\nab)^{\leq 1} \dk^{\leq \reg}\phis{p}](\tau_1,\tau_2)\nn\\
\les{}& \sum_{p=0,1,2}\bigg(\E_{r\geq R/2}^{(\reg)}[\phis{p}](\tau_1) 
+ R^2 \M_{R/2, R}^{(\reg)}[\phis{p}](\tau_1,\tau_2)+\int_{\MM_{r\geq R/2}(\tau_1,\tau_2)} r^{1+\de} |\dk^{\leq \reg} \N^{(p)}_{W,s}|^2\bigg)\nn\\
&+\sum_{p=0,1}\int_{\MM_{r\geq R/2}(\tau_1,\tau_2)}r^{-1+\de}|\dk^{\leq \reg+1} \N^{(p)}_{T,s}|^2+R^{-2\de}\sum_{p=0,1,2}{\EMF}^{(\reg)}_{\de, r\geq R/2}[\phis{p}](\tau_1,\tau_2).
\end{align*}
Finally, for $R\geq 20m$ large enough, we may absorb the last term on the RHS which then yields the desired estimate \eqref{eq:EMnearinfinity:highorderweightedderivatives:Teu:pm2}. This concludes the proof of Proposition \ref{prop:EnergyMorawetznearinfinitytensorialTeuk}.
\end{proof}


\section{Weak Morawetz estimates for Teukolsky in Kerr}
\lab{sec:MorawetzestimatesforTeukolskys=pm2inKerrfromMillet}


The goal of this section is to prove Theorem \ref{cor:weakMorawetzforTeukolskyfromMillet:bis}. Throughout the section, we work on Kerr and $(e_3, e_4, e_1, e_2)$ denotes the null frame of Kerr defined in \eqref{def:e3e4inKerr} and \eqref{def:e1e2inKerr}. For $\pmb\phi_s\in\sk_2(\mathbb{C})$, $s=\pm 2$, consider the following inhomogeneous tensorial Teukolsky wave equations in a subextremal Kerr spacetime, for $\tau_0\geq 1$,
\bea
\lab{eq:inhomoTeu:tensorial:Kerr}
\T_s\pmb\phi_s= \N_s,
\eea
where $\T_s$ is the tensorial wave operator on the LHS of equation \eqref{eq:TeukolskyequationforAandAbintensorialforminKerr}, i.e.,
\bea\lab{eq:formofTeukolskyoperatorintensorialform:sec10}
\nn\T_s &:=& \squared_2 -\frac{4ia\cos\th}{|q|^2}\nab_{\pr_t} -  \frac{s}{|q|^2} +\frac{2s}{|q|^2}(r-m)\nab_3\pmb  - \frac{4sr}{|q|^2}\nab_{\pr_t}\\
&&+\frac{4a\cos\th}{|q|^6}\Big(a\cos\th\big(|q|^2+6mr\big) - is\big((r-m)|q|^2+4mr^2\big)\Big).
\eea


\subsection{Basic estimates for the tensorial Teukolsky equation in Kerr}


We first derive weighted local energy estimates for solutions to \eqref{eq:inhomoTeu:tensorial:Kerr}.
\begin{lemma}[Weighted local energy estimates for \eqref{eq:inhomoTeu:tensorial:Kerr}]
\lab{lemma:localenergyestimatewithstrictlynegativeweightsinKerr}
Let $\pmb\phi_s\in\sk_2(\mathbb{C})$, $s=\pm 2$, be a solution to \eqref{eq:inhomoTeu:tensorial:Kerr}. Then, for any $\reg\in\mathbb{N}$, $q>0$ and $p<0$, there exists a constant $C(p,\reg)\geq 0$ large enough such that we have the following weighted local energy estimate for $\pmb\phi_s$ 
\bea\lab{eq:localenergyestimatewithstrictlynegativeweightsinKerr}
\int_{\Si(\tau_0+q)}r^{p-2}|\dk^{\leq \reg+1}\pmb\phi_s|^2\leq e^{C(p,\reg)q}\left(\int_{\Si(\tau_0)}r^{p-2}|\dk^{\leq \reg+1}\pmb\phi_s|^2+\int_{\MM(\tau_0, \tau_0+q)}r^{p+1}|\dk^{\leq\reg}\N_s|^2\right). 
\eea
\end{lemma}

\begin{proof}
First, notice from \eqref{eq:formofTeukolskyoperatorintensorialform:sec10} that 
\beaa
\T_s\pmb\phi_s &=& \squared_2\pmb\phi_s+O(r^{-2})\dk^{\leq 1}\pmb\phi_s+O(r^{-1})\nab_{e_3-2\pr_t}\pmb\phi_s\\
&=& \squared_2\pmb\phi_s+O(r^{-2})\dk^{\leq 1}\pmb\phi_s+O(r^{-1})\nab_4\pmb\phi_s\\
&=& \squared_2\pmb\phi_s+O(r^{-2})\dk^{\leq 1}\pmb\phi_s
\eeaa
which together with \eqref{eq:inhomoTeu:tensorial:Kerr} implies
\bea\lab{eq:veryusefullconsequenceofTeukolskywhichistensorialwaveuptoOrminus2term}
\squared_2\pmb\phi_s = O(r^{-2})\dk^{\leq 1}\pmb\phi_s+\N_s.
\eea
Setting, for $p<0$, $\pmb\phi_{s,p}:=r^{\frac{p}{2}}\pmb\phi_s$, we infer
\beaa
\squared_2\pmb\phi_{s,p} &=& r^{\frac{p}{2}}\left(\squared_2\pmb\phi_s +2\gam^{\a\b}\frac{\pr_\a(r^{\frac{p}{2}})}{r^{\frac{p}{2}}}\nab_{\pr_\b}\pmb\phi_s+\frac{\square_{\gam}(r^{\frac{p}{2}})}{r^{\frac{p}{2}}}\pmb\phi_s\right)\\
&=& r^{\frac{p}{2}}\left(O_p(r^{-2})\dk^{\leq 1}\pmb\phi_s+\N_s +\frac{p}{r}\gam^{r\b}\nab_{\pr_\b}\pmb\phi_s\right)\\
&=& r^{\frac{p}{2}}\left(O_p(r^{-2})\dk^{\leq 1}\pmb\phi_s+\N_s -\frac{p}{r}\nab_{\pr_\tau}\pmb\phi_s\right)
\eeaa
where we used \eqref{eq:assymptiticpropmetricKerrintaurxacoord:1} in the last identity, and hence
\beaa
\squared_2\pmb\phi_{s,p} +\frac{p}{r}\nab_{\pr_\tau}\pmb\phi_{s,p} &=& O_p(r^{-2})\dk^{\leq 1}\pmb\phi_{s,p}+\N_{s,p}, \qquad \N_{s,p}:=r^{\frac{p}{2}}\N_s. 
\eeaa
Next, we scalarize this tensorial equation using the regular triplet $\Om_i$, $j=1,2,3$, in Kerr introduced in Definition \ref{def:regulartripletinKerrOmii=123}. Introducing the notations 
\beaa
\phi_{s,p,ij}:=\pmb\phi_{s,p}(\Om_i, \Om_j), \qquad N_{s,p,ij}:=\N_{s,p}(\Om_i, \Om_j),
\eeaa
we obtain, using Lemmas \ref{lemma:formoffirstordertermsinscalarazationtensorialwaveeq} and \ref{lemma:computationoftheMialphajinKerr}, the following schematic coupled system of scalar wave equations for $\phi_{s,p,ij}$
\beaa
\square_{\gam}(\phi_{s,p,ij})+\frac{p}{r}\pr_\tau(\phi_{s,p,ij})=O_p(r^{-2})\dk^{\leq 1}\phi_{s,p,kl}+N_{s,p,ij}.
\eeaa
Next, we commute by $(\pr_\tau, r\pr_r, \pr_{x^a})^{\leq\reg}$ and obtain in view of Lemma \ref{lem:commutatorwithwave:firstorderderis} restricted to the particular case of Kerr 
\beaa
\square_{\gam}(\dk^{\leq\reg}\phi_{s,p,ij})+\frac{p}{r}\pr_\tau(\dk^{\leq\reg}\phi_{s,p,ij})=O_{p,\reg}(r^{-2})\dk^{\leq \reg+1}\phi_{s,p,kl}+\dk^{\leq\reg}N_{s,p,ij}.
\eeaa
Now, we apply Proposition \ref{prop-app:stadard-comp-Psi} to $\dk^{\leq\reg}\phi_{s,p,ij}\in\sk_0(\mathbb{C})$ with $V=0$, 
\beaa
\N=-\frac{p}{r}\pr_\tau(\dk^{\leq\reg}\phi_{s,p,ij})+O_{p,\reg}(r^{-2})\dk^{\leq \reg+1}\phi_{s,p,kl}+\dk^{\leq\reg}N_{s,p,ij},
\eeaa
and we choose $w=0$, and a vector field $X$ that is globally uniformly timelike in $\MM$ and equals $\pr_{\tau}$ for $r\geq 3m$. By integrating over $\MM(\tau_0, \tau)$, for $\tau\in[\tau_0, \tau_0+q]$, and using the fact that $p<0$, we infer
\beaa
&&\EF^{(\reg)}[\pmb\phi_{s,p}](\tau_0, \tau)+|p|\sum_{i,j}\int_{\MM(\tau_0, \tau)}r^{-1}|\pr_\tau(\dk^{\leq\reg}\phi_{s,p,ij})|^2\\ 
&\les_{p,\reg}& \E^{(\reg)}[\pmb\phi_{s,p}](\tau_0) +\int_{\MM(\tau_0, \tau)}r^{-2}|\dk^{\leq\reg+1}\pmb\phi_{s,p}|^2+\sum_{i,j}\int_{\MM(\tau_0, \tau)}|\pr_\tau(\dk^{\leq\reg}\phi_{s,p,ij})||\dk^{\leq\reg}N_{s,p,ij}|
\eeaa
and hence, as $|p|>0$, we infer
\beaa
\EF^{(\reg)}[\pmb\phi_{s,p}](\tau_0, \tau) \les_{p,\reg} \E^{(\reg)}[\pmb\phi_{s,p}](\tau_0) +\int_{\tau_0}^{\tau}\E^{(\reg)}[\pmb\phi_{s,p}](\tau)d\tau+\int_{\MM(\tau_0, \tau)}r|\dk^{\leq\reg}\N_{s,p}|^2.
\eeaa
Then, using Gr\"onwall, Lemma \ref{lemma:equivalentofmodulussquarepmbpsiandsumijpisijalsowithdkreg} restricted to Kerr, and the fact that 
\beaa
\phi_{s,p,ij}=r^{\frac{p}{2}}\pmb\phi_s(\Om_i, \Om_j), \qquad N_{s,p,ij}=r^{\frac{p}{2}}\N_s(\Om_i, \Om_j),
\eeaa
we infer the existence of a constant $C(p,\reg)\geq 0$ large enough such that we have 
\beaa
\int_{\Si(\tau_0+q)}r^{p-2}|\dk^{\leq \reg+1}\pmb\phi_s|^2\leq e^{C(p,\reg)q}\left(\int_{\Si(\tau_0)}r^{p-2}|\dk^{\leq \reg+1}\pmb\phi_s|^2+\int_{\MM(\tau_0, \tau_0+q)}r^{p+1}|\dk^{\leq\reg}\N_s|^2\right)
\eeaa
as stated in \eqref{eq:localenergyestimatewithstrictlynegativeweightsinKerr}. This concludes the proof of Lemma \ref{lemma:localenergyestimatewithstrictlynegativeweightsinKerr}.
\end{proof}

In the next lemma, we produce a solution to \eqref{eq:inhomoTeu:tensorial:Kerr} with trivial initial data at $\tau=\tau_0$. 
\begin{lemma}\lab{lemma:inhomoTeu:tensorial:Kerr:withcutoffintime}
Let $\pmb\phi_s\in\sk_2(\mathbb{C})$, $s=\pm 2$, be a solution to \eqref{eq:inhomoTeu:tensorial:Kerr}, let $\chi_{\tau_0}=\chi_{\tau_0}(\tau)$ be a cut-off function such that $\chi_{\tau_0}(\tau)=0$ for $\tau\leq\tau_0$ and $\chi_{\tau_0}(\tau)=1$ for $\tau\geq \tau_0+1$, and let 
\bea\lab{eq:relationpmbpsistopmbphisontaugeqtau0withcutoffchitau0}
\pmb\psi_s:=\chi_{\tau_0}(\tau)\pmb\phi_s, \qquad s=\pm 2, \qquad \pmb\psi_s\in\sk_2(\mathbb{C}).
\eea
Then, $\pmb\psi_s$ has trivial initial data at $\tau=\tau_0$ and satisfies 
\bea
\lab{eq:inhomoTeu:tensorial:Kerr:withcutoffintime}
\T_s\pmb\psi_s= \widetilde{\N}_s, 
\eea
where 
\bea\lab{eq:inhomoTeu:tensorial:Kerr:withcutoffintime:RHS}
\widetilde{\N}_s := \chi_{\tau_0}(\tau)\N_s+O(r^{-1})\chi_{\tau_0}'(\tau)\nab_{\pr_r}(r\pmb\phi_s)+O(r^{-2})\Big(\chi_{\tau_0}''(\tau), \chi_{\tau_0}'(\tau)\Big)\dk^{\leq 1}\pmb\phi_s.
\eea
\end{lemma}

\begin{proof}
We have
\beaa
\widetilde{\N}_s &=&  \T_s(\chi_{\tau_0}(\tau)\pmb\phi_s)=\chi_{\tau_0}(\tau)\T_s\pmb\phi_s+[\T_s, \chi_{\tau_0}(\tau)]\pmb\phi_s\\
&=& \chi_{\tau_0}(\tau)\N_s+[\squared_2, \chi_{\tau_0}(\tau)]\pmb\phi_s+O(r^{-2})\Big(\chi_{\tau_0}''(\tau), \chi_{\tau_0}'(\tau)\Big)\dk^{\leq 1}\pmb\phi_s.
\eeaa
Next, we compute 
\beaa
[\squared_2, \chi_{\tau_0}(\tau)]\pmb\phi_s &=& 2\gam^{\a\b}\pr_\a(\chi_{\tau_0})\nab_{\pr_\b}\pmb\phi +\square_{\gam}(\chi_{\tau_0})\pmb\phi_s\\
&=& 2\gam^{\tau\b}\chi_{\tau_0}'(\tau)\nab_{\pr_\b}\pmb\phi +\left(\gam^{\tau\tau}\chi_{\tau_1}''(\tau)+\frac{1}{\sqrt{|\gam|}}\pr_\a\Big(\sqrt{|\gam|}\gam^{\a\tau}\Big)\chi_{\tau_1}'(\tau)\right)\pmb\phi_s.
\eeaa
Now, recall that we have in view of  \eqref{eq:assymptiticpropmetricKerrintaurxacoord:1} 
\beaa
\gam^{\tau\tau}=O(m^2r^{-2}), \qquad \gam^{\tau r}=-1+O(m^2r^{-2}), \qquad \gam^{\tau a}=O(mr^{-2}),
\eeaa
and from \eqref{eq:assymptiticpropmetricKerrintaurxacoord:volumeform} 
\beaa
\frac{1}{\sqrt{|\gam|}}\pr_r\Big(\sqrt{|\gam|}\Big)=\frac{2}{r}(1+O(m^2r^{-2})), \quad \frac{1}{\sqrt{|\gam|}}\pr_{x^a}\Big(\sqrt{|\gam|}\Big)=O(1).
\eeaa
Hence, we infer
\beaa
[\squared_2, \chi_{\tau_0}(\tau)]\pmb\phi_s &=& -2\chi_{\tau_0}'(\tau)\left(\nab_{\pr_r}\pmb\phi_s+\frac{1}{r}\pmb\phi_s\right)+O(r^{-2})\Big(\chi_{\tau_0}''(\tau), \chi_{\tau_0}'(\tau)\Big)\dk^{\leq 1}\pmb\phi_s
\eeaa
which yields 
\beaa
\widetilde{\N}_s &=& \chi_{\tau_0}(\tau)\N_s+[\squared_2, \chi_{\tau_0}(\tau)]\pmb\phi_s+O(r^{-2})\Big(\chi_{\tau_0}''(\tau), \chi_{\tau_0}'(\tau)\Big)\dk^{\leq 1}\pmb\phi_s\\
&=& \chi_{\tau_0}(\tau)\N_s+O(r^{-1})\chi_{\tau_0}'(\tau)\nab_{\pr_r}(r\pmb\phi_s)+O(r^{-2})\Big(\chi_{\tau_0}''(\tau), \chi_{\tau_0}'(\tau)\Big)\dk^{\leq 1}\pmb\phi_s
\eeaa
as stated. This concludes the proof of Lemma \ref{lemma:inhomoTeu:tensorial:Kerr:withcutoffintime}.
\end{proof}


\subsection{Weak Morawetz estimate for Teukolsky in Kerr using \cite{Millet}}
\lab{sec:weakMorawetzestimatesforTeukolskyinKerrusingMillet}


Let $\pmb\psi_s\in\sk_2(\mathbb{C})$, $s=\pm 2$, be the solution to the tensorial Teukolsky equation in Kerr exhibited in Lemma \ref{lemma:inhomoTeu:tensorial:Kerr:withcutoffintime}. As in \eqref{eq:linkteukolskyscalarsandtensors}, we associate to $\pmb\psi_s$ the following complex-valued scalars $\psi_{s,\text{NP}}$, $s=\pm 2$, defined by 
\bea\lab{eq:linkteukolskyscalarsandtensors:bisinsection10}
\psi_{+2,\text{NP}}:=\pmb\psi_{+2}(e_1, e_1), \qquad \psi_{-2,\text{NP}}:=\ov{\pmb\psi_{-2}(e_1, e_1)}.
\eea
In view of Lemma \ref{lemma:Teukolskyintensorialform}, \eqref{eq:inhomoTeu:tensorial:Kerr:withcutoffintime} is equivalent to the following complex-valued Teukolsky equations in NP formalism in a subextremal Kerr spacetime, for $s=\pm 2$ and $\tau_0\geq 1$,
\begin{align}\label{eq:inhomoTME:NP:Kerr}
T_s(\psi_{s,\text{NP}}) =f_s,
\end{align}
where the Teukolsky operator $T_s$ is given by the operator on the LHS of \eqref{eq:TME}, i.e.,
\bea
T_s:=|q|^2\square_{\gam}  -2ias \cos{\theta} \partial_t+\frac{2is \cos{\theta}}{\sin^2{\theta} }  \partial_{\phi} - (s^2 \cot^2{\theta} +s) +2s\big((r-m)e_3-2r\partial_t\big),
\eea
and where $f_s$, $s=\pm 2$, is given by
\bea
\lab{eq:varphif:relationfstoNs}
f_{+2}:=|q|^2\widetilde{\N}_{+2}(e_1,e_1), \qquad f_{-2}:=|q|^2\ov{\widetilde{\N}_{-2}(e_1,e_1)}.
\eea

In this section, we derive a weak Morawetz estimate for the solution $\psi_{s,\text{NP}}$, $s=\pm 2$, to the inhomogeneous Teukolsky wave equation \eqref{eq:inhomoTME:NP:Kerr} by relying on the results of \cite{Millet}. To this end, we recall some notations and results from \cite{Millet}. First, for a scalar function $\psi$, we introduce the following definition of weighted Sobolev spaces on $\Si(\tau)$, see Section 3.1 in \cite{Millet},
\bea\lab{eq:bcalculussobolevspacesiestandardweightedSobolevspaces}
\|\psi\|^2_{\ov{H}^{\tilde{r}, l}_{(b)}}:=\int_{\Sigma{(\tau)}}r^{2l}|(r\pr_r, r\nab)^{\leq \tilde{r}}\psi|^2,
\eea 
where $\tilde{r}\in\mathbb{N}$ and $l\in\mathbb{R}$, and  
\bea\lab{eq:bcalculussobolevspacesiestandardweightedSobolevspaces:semiclassical}
\|\psi\|^2_{\ov{H}^{\tilde{r}, l}_{(b),h}}:=\int_{\Sigma{(\tau)}}r^{2l}|(hr\pr_r, hr\nab)^{\leq \tilde{r}}\psi|^2,
\eea 
where $h>0$ is a constant. Note that we have the following comparison between the two norms
\bea\lab{eq:comparisionbetweenweightedandsemiclassicalnorms}
\|u\|_{\ov{H}^{\tilde{r},l}_{(b),h}}\leq (1+h^{\tilde{r}})\|u\|_{\ov{H}^{\tilde{r},l}_{(b)}}, \qquad \|u\|_{\ov{H}^{\tilde{r},l}_{(b)}}\leq (1+h^{-\tilde{r}})\|u\|_{\ov{H}^{\tilde{r},l}_{(b),h}}.
\eea 
Also, for $\sigma\in\mathbb{C}$ and a scalar function $\psi$, we define, the Fourier-Laplace transform w.r.t. $\tau$ by
\beaa
\widehat{\psi}(\sigma,\c)=\int_{\mathbb{R}}e^{i\sigma\tau}\psi(\tau, \c)d\tau,
\eeaa 
where $(\tau, r, x^a)$ denote the normalized coordinates in Kerr. In particular, taking the Fourier-Laplace transform of \eqref{eq:inhomoTME:NP:Kerr}, we obtain 
\begin{align}\label{eq:inhomoTME:NP:Kerr:Fourier-Laplace}
\widehat{T}_s(\sigma)\widehat{\psi}_{s,\text{NP}}(\sigma) =\widehat{f}_s(\sigma),
\end{align}
where $\widehat{T}_s(\sigma)$ is the second order elliptic operator in $(r,x^a)$ obtained by writing $T_s$ in the normalized coordinates $(\tau, r, x^a)$ of Kerr and then replacing each $\pr_\tau$ derivative by $-i\sigma$. We then introduce the resolvent operator $R_s(\sigma)$ defined for $\Im(\sigma)\geq 0$ by 
\beaa
R_s(\sigma):=\widehat{T_s}(\sigma)^{-1},
\eeaa
which allows to rewrite \eqref{eq:inhomoTME:NP:Kerr:Fourier-Laplace}, for $\Im(\sigma)\geq 0$, as 
\begin{align}\label{eq:inhomoTME:NP:Kerr:Fourier-Laplace:resolvent}
\widehat{\psi}_{s,\text{NP}}(\sigma) =R_s(\sigma)\widehat{f}_s(\sigma).
\end{align}
Finally, we recall Proposition 7.2 in \cite{Millet} on the properties of the resolvent operator $R_s(\sigma)$.

\begin{proposition}[Proposition 7.2 in \cite{Millet}]
For every $\eta\in [0,1]$, the resolvent operator $R_s(\sigma)$ is a bounded operator from $\ov{H}^{\tilde{r}, l}_{(b)}$ to $\ov{H}^{\tilde{r}, l+1-\eta}_{(b)}$ for $\Im(\sigma)\geq 0$, $\sigma\neq 0$, $|\sigma|\leq c$ and 
\bea\lab{eq:rangeparameterstilderandlforboundednessresolvant:1}
-\frac{3}{2}-s-|s|<l+1-\eta<-\frac{1}{2},\qquad \tilde{r}+l+1-\eta>-\frac{1}{2}-2s,\qquad \tilde{r}>\frac{1}{2}+s.
\eea 
Moreover, in this case, we have the following bound (uniform in $|\sigma|\leq c$ for $c$ small enough):
\bea
\|R_s(\sigma)\|_{\LL(\ov{H}^{\tilde{r}, l}_{(b)}, \ov{H}^{\tilde{r}, l+1-\eta}_{(b)})}\leq C|\sigma|^{\eta-1}. 
\eea
It is also a bounded operator from $\ov{H}^{\tilde{r}, l}_{(b)}$ to $\ov{H}^{\tilde{r}, l+1}_{(b)}$ for $\sigma\neq 0$ and 
\bea\lab{eq:rangeparameterstilderandlforboundednessresolvant:2}
l+1<-\frac{1}{2},\qquad \tilde{r}+l+1>-\frac{1}{2}-2s-4m\Im(\sigma),\qquad \tilde{r}>\frac{1}{2}+s-\frac{r_+^2+a^2}{r_+-m}\Im(\sigma),
\eea 
and in this case we have the bound (uniform for $\sigma$ in a strip $\{0\leq\Im(\sigma)\leq A, \,|\sigma|>\frac{1}{A}\}$): 
\bea
\|R_s(\sigma)\|_{\LL(\ov{H}^{\tilde{r}, l}_{(b),|\sigma|^{-1}}, \ov{H}^{\tilde{r}, l+1}_{(b),|\sigma|^{-1}})}\leq C.
\eea
\end{proposition}

We are now ready to prove the following proposition.

\begin{proposition}
\lab{prop:weakMorawetzforTeukolskyfromMillet}
Let $0<\de<1$. Assume that the complex-valued scalars $\psi_{s,\text{NP}}$, $s=\pm 2$, satisfy the inhomogeneous Teukolsky wave equation \eqref{eq:inhomoTME:NP:Kerr} for $\tau\geq \tau_0$, and that $\psi_{s,\text{NP}}$ and the RHS $f_s$ can be smoothly extended to $\tau\leq\tau_0$ by 0. Then, we have 
\bea
\lab{eq:weakmora:rprrandangular:+2}
\int_{\MM(\tau\geq\tau_0)}r^{-11+\de}|(r\pr_r, r\nab)^{\leq 3}{\psi_{+2,\text{NP}}}|^2 &\les&{\int_{\MM(\tau\geq\tau_0)}r^{-11+\de}|(r\pr_r, r\nab)^{\leq 3}\pr^{\leq 3}_{\tau}f_{+2}|^2}
\eea
and 
\bea
\lab{eq:weakmora:rprrandangular:-2}
\int_{\MM{(\tau\geq\tau_0)}}r^{-3+\de}|(r\pr_r, r\nab)^{\leq 5}{\psi_{-2,\text{NP}}}|^2&\les&{\int_{\MM(\tau\geq\tau_0)}r^{-3+\de}|(r\pr_r, r\nab)^{\leq 5}\pr^{\leq 5}_{\tau}f_{-2}|^2}.
\eea
\end{proposition}

\begin{proof}
In view of Lemma \ref{lemma:localenergyestimatewithstrictlynegativeweightsinKerr} applied to \eqref{eq:inhomoTME:NP:Kerr}, using also the fact that $\psi_s$, $s=\pm 2$, have trivial initial data at $\tau=\tau_0$, there exist large enough constants $C_{+2}=C(\de-9,3)>0$ and $C_{-2}=C(\de-1,5)>0$ such that 
\beaa
\bsplit
\|\psi_{+2}(\tau, \c)\|^2_{\ov{H}^{3, \de-11}_{(b)}} &\les e^{C_{+2}(\tau-\tau_0)}\int_{\MM(\tau\geq\tau_0)}r^{-12+\de}|(r\pr_r, r\nab)^{\leq 3}f_{+2}|^2,\\
\|\psi_{-2}(\tau, \c)\|^2_{\ov{H}^{5, \de-3}_{(b)}} &\les e^{C_{-2}(\tau-\tau_0)}\int_{\MM(\tau\geq\tau_0)}r^{-4+\de}|(r\pr_r, r\nab)^{\leq 5}f_{-2}|^2.
\end{split}
\eeaa
In particular, using the Fourier inversion formula, we infer
\beaa
\psi_{s,\text{NP}}=\frac{1}{2\pi}\int_{\Im(\sigma)=2C_s}e^{-i\sigma\tau}\widehat{\psi}_{s,\text{NP}}(\sigma)d\sigma,
\eeaa
which together with \eqref{eq:inhomoTME:NP:Kerr:Fourier-Laplace:resolvent} implies
\beaa
\psi_{s,\text{NP}}=\frac{1}{2\pi}\int_{\Im(\sigma)=2C_s}e^{-i\sigma\tau}R_s(\sigma)\widehat{f_s}(\sigma)d\sigma.
\eeaa

We first prove \eqref{eq:weakmora:rprrandangular:+2} \eqref{eq:weakmora:rprrandangular:-2}  in the particular case where $f_s\in C^\infty_c(\MM)$, i.e., $f_s$ is smooth and compactly supported in $\MM$. In this case, the contour argument, which is based on the holomorphic properties of $\sigma\to R_s(\sigma)\widehat{f_s}(\sigma)$ in $\Im(\sigma)>0$ and outlined at the beginning of Section 8 in \cite{Millet},  applies. This yields
\beaa
\psi_{s,\text{NP}} &=& \frac{1}{2\pi}\int_{\Im(\sigma)=2C_s}e^{-i\sigma\tau}R_s(\sigma)\widehat{f_s}(\sigma)d\sigma\\ 
&=& \frac{1}{2\pi}\int_{\Im(\sigma)=0}e^{-i\sigma\tau}R_s(\sigma)\widehat{f_s}(\sigma)d\sigma,
\eeaa
and hence
\beaa
\psi_{s,\text{NP}} = \FF_\tau^{-1}(R{_s}(\sigma)\widehat{f_s}(\sigma)),
\eeaa
where $\FF_\tau^{-1}$ denotes the inverse Fourier transform w.r.t. $\tau$. We infer
\beaa
\|{\psi_{s,\text{NP}}}\|_{L^2_\tau\ov{H}^{\tilde{r}, l}_{(b)}} &=& \|\FF_\tau^{-1}(R{_s}(\sigma)\widehat{f_s}(\sigma))\|_{L^2_\tau\ov{H}^{\tilde{r}, l}_{(b)}},
\eeaa
which together with Plancherel's lemma  implies, for any $\tilde{r}\in\mathbb{N}$ and $l\in\mathbb{R}$,
\bea\lab{eq:expressingpsiswithinversefouriertransform:afterPlancherel} 
\|{\psi_{s,\text{NP}}}\|_{L^2_\tau\ov{H}^{\tilde{r}, l}_{(b)}} =\|R{_s}(\sigma)\widehat{f_s}(\sigma)\|_{L^2_\sigma\ov{H}^{\tilde{r}, l}_{(b)}}.
\eea

We now separate the cases $s=\pm 2$ and start with the case $s=+2$ with the choice $\tilde{r}=3$ which satisfies the third condition in both \eqref{eq:rangeparameterstilderandlforboundednessresolvant:1} and \eqref{eq:rangeparameterstilderandlforboundednessresolvant:2} in $\Im(\sigma)\geq 0$. Also, we choose $\eta=1$. Then, we choose $l=-\frac{11}{2}+\frac{\de}{2}$ so that $(\tilde{r}=3, l=-\frac{11}{2}+\frac{\de}{2}, \eta=1)$  satisfies all conditions in \eqref{eq:rangeparameterstilderandlforboundednessresolvant:1}   which yields  
\beaa
\|R_{+2}(\sigma)\|_{\LL\left(\ov{H}^{3, -\frac{11}{2}+\frac{\de}{2}}_{(b)}, \ov{H}^{3, -\frac{11}{2}+\frac{\de}{2}}_{(b)}\right)}\leq C, \qquad \Im(\sigma)\geq 0, \qquad |\sigma|\leq c.
\eeaa
Also, we choose $l=-\frac{13}{2}+\frac{\de}{2}$ so that $(\tilde{r}=3, l=-\frac{13}{2}+\frac{\de}{2})$  satisfies all conditions in \eqref{eq:rangeparameterstilderandlforboundednessresolvant:2}   which yields  
\beaa
\|R_{+2}(\sigma)\|_{\LL\left(\ov{H}^{3, -\frac{13}{2}+\frac{\de}{2}}_{(b){,|\sigma|^{-1}}}, \ov{H}^{3, -\frac{11}{2}+\frac{\de}{2}}_{(b){,|\sigma|^{-1}}}\right)}\leq C, \qquad 0\leq\Im(\sigma)\leq A, \qquad |\sigma|>\frac{1}{A}.
\eeaa
Applying the first estimate for $|\sigma|\leq c$ and otherwise the second estimate with $A:=c^{-1}$, we infer, relying also on \eqref{eq:comparisionbetweenweightedandsemiclassicalnorms} with $h=|\sigma|^{-1}$,
\beaa
\|R_{+2}(\sigma)u\|_{\ov{H}^{3, -\frac{11}{2}+\frac{\de}{2}}_{(b)}}&\leq& C\|u\|_{\ov{H}^{3, -\frac{11}{2}+\frac{\de}{2}}_{(b)}}+C(1+c^{-3})(1+|\sigma|^3)\|u\|_{\ov{H}^{3, -\frac{13}{2}+\frac{\de}{2}}_{(b)}}\\
&\les&(1+|\sigma|^3)\|u\|_{\ov{H}^{3, -\frac{11}{2}+\frac{\de}{2}}_{(b)}},
\eeaa
which together with \eqref{eq:expressingpsiswithinversefouriertransform:afterPlancherel} implies
\beaa
\|{\psi_{+2,\text{NP}}}\|_{L^2_\tau\ov{H}^{3,-\frac{11}{2}+\frac{\de}{2}}_{(b)}} &=& \|R{_{+2}}(\sigma)\widehat{f_{+2}}(\sigma)\|_{L^2_\sigma\ov{H}^{3, -\frac{11}{2}+\frac{\de}{2}}_{(b)}}\\
&\les& {\|{\widehat{f_{+2}}(\sigma)}\|_{L^2_\sigma\ov{H}^{3, -\frac{11}{2}+\frac{\de}{2}}_{(b)}}+ \|{\sigma^3\widehat{f_{+2}}(\sigma)}\|_{L^2_\sigma\ov{H}^{3, -\frac{11}{2}+\frac{\de}{2}}_{(b)}}}\\
&\les& {\|{\widehat{f_{+2}}(\sigma)}\|_{L^2_\sigma\ov{H}^{3, -\frac{11}{2}+\frac{\de}{2}}_{(b)}}+ \|{\widehat{\pr_\tau^3f_{+2}}(\sigma)}\|_{L^2_\sigma\ov{H}^{3, -\frac{11}{2}+\frac{\de}{2}}_{(b)}}},
\eeaa
and hence, using Plancherel, we deduce 
\beaa
\|{\psi_{+2,\text{NP}}}\|_{L^2_\tau\ov{H}^{3,-\frac{11}{2}+\frac{\de}{2}}_{(b)}} \les\|{\pr_\tau^{\leq 3}}f_{+2}\|_{L^2_\tau\ov{H}^{3,-\frac{11}{2}+\frac{\de}{2}}_{(b)}},
\eeaa
or 
\beaa
\nn\int_{\MM(\tau\geq\tau_0)}r^{-11+\de}|(r\pr_r, r\nab)^{\leq 3}{\psi_{+2,\text{NP}}}|^2 &\les&\int_{\MM(\tau\geq\tau_0)}r^{-11+\de}|(r\pr_r, r\nab)^{\leq 3}\pr_\tau^{\leq 3}f_{+2}|^2,
\eeaa
which is the  desired estimate \eqref{eq:weakmora:rprrandangular:+2}.

Next, we consider the case $s=-2$ with the choice $\tilde{r}=5$ which satisfies the third condition in both \eqref{eq:rangeparameterstilderandlforboundednessresolvant:1} and \eqref{eq:rangeparameterstilderandlforboundednessresolvant:2} in $\Im(\sigma)\geq 0$. Also, we choose $\eta=1$. Then, we choose $l=-\frac{3}{2}+\frac{\de}{2}$ so that $(\tilde{r}=5, l=-\frac{3}{2}+\frac{\de}{2}, \eta=1)$  satisfies all conditions in \eqref{eq:rangeparameterstilderandlforboundednessresolvant:1}   which yields  
\beaa
\|R_{-2}(\sigma)\|_{\LL\left(\ov{H}^{5, -\frac{3}{2}+\frac{\de}{2}}_{(b)}, \ov{H}^{5, -\frac{3}{2}+\frac{\de}{2}}_{(b)}\right)}\leq C, \qquad \Im(\sigma)\geq 0, \qquad |\sigma|\leq c.
\eeaa
Also, we choose $l=-\frac{5}{2}+\frac{\de}{2}$ so that $(\tilde{r}=5, l=-\frac{5}{2}+\frac{\de}{2})$  satisfies all conditions in \eqref{eq:rangeparameterstilderandlforboundednessresolvant:2}   which yields  
\beaa
\|R_{-2}(\sigma)\|_{\LL\left(\ov{H}^{5, -\frac{5}{2}+\frac{\de}{2}}_{(b){,|\sigma|^{-1}}}, \ov{H}^{5, -\frac{3}{2}+\frac{\de}{2}}_{(b){,|\sigma|^{-1}}}\right)}\leq C, \qquad 0\leq\Im(\sigma)\leq A, \qquad |\sigma|>\frac{1}{A}.
\eeaa
Applying the first estimate for $|\sigma|\leq c$ and otherwise the second estimate with $A:=c^{-1}$, we infer, relying also on \eqref{eq:comparisionbetweenweightedandsemiclassicalnorms} with $h=|\sigma|^{-1}$,
\beaa
\|R_{-2}(\sigma)u\|_{\ov{H}^{5, -\frac{3}{2}+\frac{\de}{2}}_{(b)}}\leq C\|u\|_{\ov{H}^{5, -\frac{3}{2}+\frac{\de}{2}}_{(b)}}+C(1+c^{-5})(1+|\sigma|^5)\|u\|_{\ov{H}^{5, -\frac{5}{2}+\frac{\de}{2}}_{(b)}}\les {(1+|\sigma|^5)}\|u\|_{\ov{H}^{5, -\frac{3}{2}+\frac{\de}{2}}_{(b)}},
\eeaa
which together with \eqref{eq:expressingpsiswithinversefouriertransform:afterPlancherel} implies
\beaa
\|{\psi_{-2,\text{NP}}}\|_{L^2_\tau\ov{H}^{5,-\frac{3}{2}+\frac{\de}{2}}_{(b)}} &=&\|R{_{-2}}(\sigma)\widehat{f_{-2}}(\sigma))\|_{L^2_\sigma\ov{H}^{5,-\frac{3}{2}+\frac{\de}{2}}_{(b)}}\\
&\les&{\|{\widehat{f_{-2}}(\sigma)}\|_{L^2_\sigma\ov{H}^{5, -\frac{3}{2}+\frac{\de}{2}}_{(b)}}+ \|{\sigma^5\widehat{f_{-2}}(\sigma)}\|_{L^2_\sigma\ov{H}^{5, -\frac{3}{2}+\frac{\de}{2}}_{(b)}}}\\
&\les& {\|{\widehat{f_{-2}}(\sigma)}\|_{L^2_\sigma\ov{H}^{5, -\frac{3}{2}+\frac{\de}{2}}_{(b)}}+ \|{\widehat{\pr_\tau^5f_{-2}}(\sigma)}\|_{L^2_\sigma\ov{H}^{5, -\frac{3}{2}+\frac{\de}{2}}_{(b)}}},
\eeaa
and hence, using Plancherel, we deduce 
\beaa
\|{\psi_{-2,\text{NP}}}\|_{L^2_\tau\ov{H}^{5,-\frac{3}{2}+\frac{\de}{2}}_{(b)}} \les\|{\pr_{\tau}^{\leq 5}}f_{-2}\|_{L^2_\tau\ov{H}^{5,-\frac{3}{2}+\frac{\de}{2}}_{(b)}},
\eeaa
or 
\beaa
\int_{\MM}r^{-3+\de}|(r\pr_r, r\nab)^{\leq 5}{\psi_{-2,\text{NP}}}|^2&\les&{\int_{\MM(\tau\geq\tau_0)}r^{-3+\de}|(r\pr_r, r\nab)^{\leq 5}\pr_\tau^{\leq 5}f_{-2}|^2},
\eeaa
which is the desired estimate \eqref{eq:weakmora:rprrandangular:-2}. {We have thus obtained \eqref{eq:weakmora:rprrandangular:+2} \eqref{eq:weakmora:rprrandangular:-2} in the particular case where $f_s\in C^\infty_c(\MM)$ and the general case follows immediately by density.} This concludes the proof of Proposition \ref{prop:weakMorawetzforTeukolskyfromMillet}.
\end{proof}

Given that $\pr_\tau$ commutes with $T_s$, the following is an immediate corollary of Proposition \ref{prop:weakMorawetzforTeukolskyfromMillet}.
\begin{corollary}
\lab{cor:weakMorawetzforTeukolskyfromMillet}
Let $0<\de<1$. Assume that the complex-valued scalars $\psi_{s,\text{NP}}$, $s=\pm 2$, satisfy the inhomogeneous Teukolsky wave equation \eqref{eq:inhomoTME:NP:Kerr} for $\tau\geq \tau_0$, and that $\psi_{s,\text{NP}}$ and the RHS $f_s$ can be smoothly extended to $\tau\leq\tau_0$ by 0. Then, we have 
\beaa
\int_{\MM(\tau\geq\tau_0)}r^{-11+\de}|(r\pr_r, r\nab)^{\leq 3}\pr_\tau^{\leq 3}\psi_{+2,\text{NP}}|^2 &\les& \int_{\MM(\tau\geq\tau_0)}r^{-11+\de}|(r\pr_r, r\nab)^{\leq 3}\pr^{\leq 6}_{\tau}f_{+2}|^2,\\
\int_{\MM(\tau\geq\tau_0)}r^{-3+\de}|(r\pr_r, r\nab)^{\leq 5}\pr_\tau^{\leq 3}\psi_{-2,\text{NP}}|^2&\les& \int_{\MM(\tau\geq\tau_0)}r^{-3+\de}|(r\pr_r, r\nab)^{\leq 5}\pr^{\leq 8}_{\tau}f_{-2}|^2.
\eeaa
\end{corollary}


\subsection{Proof of Theorem \ref{cor:weakMorawetzforTeukolskyfromMillet:bis}}


Let $\pmb\phi_s\in\sk_2(\mathbb{C})$, $s=\pm 2$, be a solution to \eqref{eq:inhomoTeu:tensorial:Kerr}, and let $\pmb\psi_s\in\sk_2(\mathbb{C})$ be the corresponding tensor exhibited in Lemma \ref{lemma:inhomoTeu:tensorial:Kerr:withcutoffintime}, so that $\pmb\psi_s$ satisfies in particular \eqref{eq:inhomoTeu:tensorial:Kerr:withcutoffintime}. Then, let $\psi_{s,\text{NP}}$, $s=\pm 2$, be the complex-valued scalars associated to $\pmb\psi_s$ as in \eqref{eq:linkteukolskyscalarsandtensors:bisinsection10}. Also, let $\widetilde{\N}_s\in\sk_2(\mathbb{C})$ be given by \eqref{eq:inhomoTeu:tensorial:Kerr:withcutoffintime:RHS} and let $f_s$ be the complex-valued scalars associated to $\widetilde{\N}_s$ as in \eqref{eq:varphif:relationfstoNs}. Then, $\psi_{s,\text{NP}}$, $s=\pm 2$, satisfy the inhomogeneous Teukolsky wave equation \eqref{eq:inhomoTME:NP:Kerr} for $\tau\geq \tau_0$, and $\psi_{s,\text{NP}}$ and $f_s$ can be smoothly extended to $\tau\leq\tau_0$ by 0. We may thus apply Corollary \ref{cor:weakMorawetzforTeukolskyfromMillet} which yields 
\beaa
\int_{\MM(\tau\geq\tau_0)}r^{-11+\de}|(r\pr_r, r\nab)^{\leq 3}\pr_\tau^{\leq 3}\psi_{+2,\text{NP}}|^2 &\les& \int_{\MM(\tau\geq\tau_0)}r^{-11+\de}|(r\pr_r, r\nab)^{\leq 3}\pr^{\leq 6}_{\tau}f_{+2}|^2,\\
\int_{\MM(\tau\geq\tau_0)}r^{-3+\de}|(r\pr_r, r\nab)^{\leq 5}\pr_\tau^{\leq 3}\psi_{-2,\text{NP}}|^2&\les& \int_{\MM(\tau\geq\tau_0)}r^{-3+\de}|(r\pr_r, r\nab)^{\leq 5}\pr^{\leq 8}_{\tau}f_{-2}|^2,
\eeaa
and in particular 
\beaa
\int_{\MM(\tau\geq\tau_0)}r^{-11+\de}|(r\pr_r, r\nab, \pr_\tau)^{\leq 3}\psi_{+2,\text{NP}}|^2 &\les& \int_{\MM(\tau\geq\tau_0)}r^{-11+\de}|(r\pr_r, r\nab, \pr_\tau)^{\leq 8}\pr_\tau^{\leq 1}f_{+2}|^2,\\
\int_{\MM(\tau\geq\tau_0)}r^{-3+\de}|(r\pr_r, r\nab, \pr_\tau)^{\leq 3}\psi_{-2,\text{NP}}|^2&\les& \int_{\MM(\tau\geq\tau_0)}r^{-3+\de}|(r\pr_r, r\nab, \pr_\tau)^{\leq 12}\pr_\tau^{\leq 1}f_{-2}|^2.
\eeaa

Now, in view of \eqref{eq:linkteukolskyscalarsandtensors:bisinsection10},   \eqref{eq:varphif:relationfstoNs}, and \eqref{eq:Kerr.La_abc}, we have, for any $\reg\in\mathbb{N}$ and $s=\pm 2$,
\beaa
\bsplit
 |\dk^{\leq\reg}\pmb\psi_s|&\les_{\reg} |(r\pr_r, r\nab, \pr_\tau)^{\leq\reg}\psi_{s,\text{NP}}|\les_{\reg} |\dk^{\leq\reg}\pmb\psi_s|,  \\
 r^2|\dk^{\leq\reg}\nab_{\pr_\tau}^{\leq 1}\widetilde{\N}_s|&\les_{\reg} |(r\pr_r, r\nab, \pr_\tau)^{\leq\reg}\pr_\tau^{\leq 1}f_s |\les_{\reg} r^2|\dk^{\leq\reg}\nab_{\pr_\tau}^{\leq 1}\widetilde{\N}_s|,
 \end{split}
\eeaa
and plugging in the above estimates, we infer
\beaa
\int_{\MM(\tau\geq\tau_0)}r^{-11+\de}|\dk^{\leq 3}\pmb\psi_{+2}|^2 &\les& \int_{\MM(\tau\geq\tau_0)}r^{-7+\de}|\dk^{\leq 8}\nab_{\pr_\tau}^{\leq 1}\widetilde{\N}_{+2}|^2,\\
\int_{\MM(\tau\geq\tau_0)}r^{-3+\de}|\dk^{\leq 3}\pmb\psi_{-2}|^2&\les& \int_{\MM(\tau\geq\tau_0)}r^{1+\de}|\dk^{\leq 12}\nab_{\pr_\tau}^{\leq 1}\widetilde{\N}_{-2}|^2.
\eeaa
Together with \eqref{eq:relationpmbpsistopmbphisontaugeqtau0withcutoffchitau0} and \eqref{eq:inhomoTeu:tensorial:Kerr:withcutoffintime:RHS}, this yields 
\beaa
\int_{\MM(\tau\geq\tau_0+1)}r^{-11+\de}|\dk^{\leq 3}\pmb\phi_{+2}|^2 &\les& \int_{\MM(\tau\geq\tau_0)}r^{-7+\de}|\dk^{\leq 9}\N_{+2}|^2+\int_{\MM(\tau_0, \tau_0+1)}r^{-11+\de}|\dk^{\leq 10}\pmb\phi_{+2}|^2\nn\\
&&+\int_{\MM(\tau_0, \tau_0+1)}r^{-9+\de}|\dk^{\leq 8}\nab_{\pr_\tau}^{\leq 1}\nab_{\pr_r}(r\pmb\phi_{+2})|^2,\\
\int_{\MM(\tau\geq\tau_0+1)}r^{-3+\de}|\dk^{\leq 3}\pmb\phi_{-2}|^2&\les& \int_{\MM(\tau\geq\tau_0)}r^{1+\de}|\dk^{\leq 13}\N_{-2}|^2+\int_{\MM(\tau_0, \tau_0+1)}r^{-3+\de}|\dk^{\leq 14}\pmb\phi_{-2}|^2\\
&&+\int_{\MM(\tau_0, \tau_0+1)}r^{-1+\de}|\dk^{\leq 12}\nab_{\pr_\tau}^{\leq 1}\nab_{\pr_r}(r\pmb\phi_{-2})|^2.
\eeaa
In view of the local energy estimates provided by Lemma \ref{lemma:localenergyestimatewithstrictlynegativeweightsinKerr}, applied on $\tau\in(\tau_0, \tau_0+1)$ to $\pmb\phi_s$ solution to the inhomogeneous tensorial Teukolsky equation \eqref{eq:inhomoTeu:tensorial:Kerr}, we deduce
\bea\lab{eq:beforetolastestimateforgettingthelowfrequencyestimatesfromMillet:plus2}
\int_{\MM(\tau\geq\tau_0)}r^{-11+\de}|\dk^{\leq 3}\pmb\phi_{+2}|^2 &\les& \int_{\MM(\tau\geq\tau_0)}r^{-7+\de}|\dk^{\leq 9}\N_{+2}|^2+\int_{\Sigma(\tau_0)}r^{-11+\de}|\dk^{\leq 10}\pmb\phi_{+2}|^2\nn\\
&&+\int_{\MM(\tau_0, \tau_0+1)}r^{-9+\de}|\dk^{\leq 8}\nab_{\pr_\tau}^{\leq 1}\nab_{\pr_r}(r\pmb\phi_{+2})|^2
\eea
and
\bea\lab{eq:beforetolastestimateforgettingthelowfrequencyestimatesfromMillet:minus2}
\nn\int_{\MM(\tau\geq\tau_0)}r^{-3+\de}|\dk^{\leq 3}\pmb\phi_{-2}|^2&\les& \int_{\MM(\tau\geq\tau_0)}r^{1+\de}|\dk^{\leq 13}\N_{-2}|^2+\int_{\Si(\tau_0)}r^{-3+\de}|\dk^{\leq 14}\pmb\phi_{-2}|^2\\
&&+\int_{\MM(\tau_0, \tau_0+1)}r^{-1+\de}|\dk^{\leq 12}\nab_{\pr_\tau}^{\leq 1}\nab_{\pr_r}(r\pmb\phi_{-2})|^2.
\eea

We begin with proving \eqref{eq:weakMorawetzforTeukolskyfromMillet:bis:minus2}. This requires to estimate the last term on the RHS \eqref{eq:beforetolastestimateforgettingthelowfrequencyestimatesfromMillet:minus2}. We first notice that, for $\tau\in[\tau_0, \tau_0+1]$, we have
\beaa
\bsplit
&\int_{\Si(\tau)}r^{-1+\de}|\dk^{\leq 12}\nab_{\pr_r}(r\pmb\phi_{-2})|^2\\ 
=& \int_{\Si(\tau_0)}r^{-1+\de}|\dk^{\leq 12}\nab_{\pr_r}(r\pmb\phi_{-2})|^2+2\int_{\MM(\tau_0, \tau)}r^{-1+\de}\Re\left(\ov{\dk^{\leq 12}\nab_{\pr_r}(r\pmb\phi_{-2})}\c\dk^{\leq 12}\nab_{\pr_\tau}\nab_{\pr_r}(r\pmb\phi_{-2})\right)
\end{split}
\eeaa
which yields, after using Gr\"onwall and integrating in $\tau$ for $\tau\in[\tau_0, \tau_0+1]$, 
\bea
&&\int_{\MM(\tau_0, \tau_0+1)}r^{-1+\de}|\dk^{\leq 12}\nab_{\pr_\tau}^{\leq 1}\nab_{\pr_r}(r\pmb\phi_{-2})|^2\nn\\ 
&\les& \int_{\Si(\tau_0)}r^{-1+\de}|\dk^{\leq 12}\nab_{\pr_r}(r\pmb\phi_{-2})|^2+\int_{\MM(\tau_0, \tau_0+1)}r^{-1+\de}|\dk^{\leq 12}\nab_{\pr_\tau}\nab_{\pr_r}(r\pmb\phi_{-2})|^2.
\eea
Plugging in \eqref{eq:beforetolastestimateforgettingthelowfrequencyestimatesfromMillet:minus2}, we obtain
\bea\lab{eq:beforetolastestimateforgettingthelowfrequencyestimatesfromMillet:minus2:1}
\nn\int_{\MM(\tau\geq\tau_0)}r^{-3+\de}|\dk^{\leq 3}\pmb\phi_{-2}|^2&\les& \int_{\MM(\tau\geq\tau_0)}r^{1+\de}|\dk^{\leq 13}\N_{-2}|^2+\int_{\Si(\tau_0)}r^{-3+\de}|\dk^{\leq 14}\pmb\phi_{-2}|^2\\
\nn&&+\int_{\Si(\tau_0)}r^{-1+\de}|\dk^{\leq 12}\nab_{\pr_r}(r\pmb\phi_{-2})|^2\\
&&+\int_{\MM(\tau_0, \tau_0+1)}r^{-1+\de}|\dk^{\leq 12}\nab_{\pr_\tau}\nab_{\pr_r}(r\pmb\phi_{-2})|^2.
\eea

We need to estimate the last term on the RHS of \eqref{eq:beforetolastestimateforgettingthelowfrequencyestimatesfromMillet:minus2:1}. Now, using \eqref{eq:simplifiedversionsquaredusefulinlargerregion:bis} in the particular case of Kerr, we have
\beaa
\squared_2\pmb\phi_{-2} = -r^{-1}\nab_3\nab_4(r\pmb\phi_{-2})+O(r^{-2})\dk^{\leq 2}\pmb\phi_{-2}.
\eeaa
Since we have, in view of \eqref{eq:relationsbetweennullframeandcoordinatesframe:1} restricted to Kerr,
\beaa
e_4=\pr_r+O(r^{-2})\dk, \qquad e_3=2\pr_\tau+O(r^{-1})\dk, 
\eeaa
we infer 
\beaa
\squared_2\pmb\phi_{-2} = -2r^{-1}\nab_{\pr_\tau}\nab_{\pr_r}(r\pmb\phi_{-2})+O(r^{-2})\dk^{\leq 2}\pmb\phi_{-2},
\eeaa
which together with \eqref{eq:veryusefullconsequenceofTeukolskywhichistensorialwaveuptoOrminus2term} implies
\beaa
\nab_{\pr_\tau}\nab_{\pr_r}(r\pmb\phi_{-2}) &=& -\frac{r}{2}\N_{-2} +O(r^{-1})\dk^{\leq 2}\pmb\phi_{-2}.
\eeaa
Plugging in \eqref{eq:beforetolastestimateforgettingthelowfrequencyestimatesfromMillet:minus2:1}, we deduce 
\beaa
\nn\int_{\MM(\tau\geq\tau_0)}r^{-3+\de}|\dk^{\leq 3}\pmb\phi_{-2}|^2&\les& \int_{\MM(\tau\geq\tau_0)}r^{1+\de}|\dk^{\leq 13}\N_{-2}|^2+\int_{\Si(\tau_0)}r^{-3+\de}|\dk^{\leq 14}\pmb\phi_{-2}|^2\\
\nn&&+\int_{\Si(\tau_0)}r^{-1+\de}|\dk^{\leq 12}\nab_{\pr_r}(r\pmb\phi_{-2})|^2+\int_{\MM(\tau_0, \tau_0+1)}r^{-3+\de}|\dk^{\leq 14}\pmb\phi_{-2}|^2.
\eeaa
In view of the local energy estimates provided by Lemma \ref{lemma:localenergyestimatewithstrictlynegativeweightsinKerr}, applied on $\tau\in(\tau_0, \tau_0+1)$ to $\pmb\phi_{-2}$, we deduce
\beaa
\nn\int_{\MM(\tau\geq\tau_0)}r^{-3+\de}|\dk^{\leq 3}\pmb\phi_{-2}|^2&\les& \int_{\MM(\tau\geq\tau_0)}r^{1+\de}|\dk^{\leq 13}\N_{-2}|^2+\int_{\Si(\tau_0)}r^{-3+\de}|\dk^{\leq 14}\pmb\phi_{-2}|^2\\
\nn&&+\int_{\Si(\tau_0)}r^{-1+\de}|\dk^{\leq 12}\nab_{\pr_r}(r\pmb\phi_{-2})|^2\\
&\les& \E^{(13)}[\pmb\phi_{-2}](\tau_0)+\E^{(11)}[r^{\frac{1+\de}{2}}\nab_{\pr_r}(r\pmb\phi_{-2})](\tau_0)\\
&&+\int_{\MM(\tau\geq\tau_0)}r^{1+\de}|\dk^{\leq 13}\N_{-2}|^2,
\eeaa
as stated in \eqref{eq:weakMorawetzforTeukolskyfromMillet:bis:minus2}.

It remains to show the estimate \eqref{eq:weakMorawetzforTeukolskyfromMillet:plus2:nolossinr}. To this end, we first obtain from \eqref{eq:beforetolastestimateforgettingthelowfrequencyestimatesfromMillet:plus2} for $s=+2$, in the same manner as proving \eqref{eq:weakMorawetzforTeukolskyfromMillet:bis:minus2} from \eqref{eq:beforetolastestimateforgettingthelowfrequencyestimatesfromMillet:minus2} for $s=-2$, that
\bea\lab{eq:weakMorawetzforTeukolskyfromMillet:bis:plus2}
\int_{\MM(\tau\geq\tau_0)}r^{-11+\de}|\dk^{\leq 3}\pmb\phi_{+2}|^2 &\les& \E^{(9)}[r^{-4}\pmb\phi_{+2}](\tau_0) +\E^{(7)}[r^{-\frac{7-\de}{2}}\nab_{\pr_r}(r\pmb\phi_{+2})](\tau_0)\nn\\
&&+\int_{\MM(\tau\geq\tau_0)}r^{-7+\de}|\dk^{\leq 9}\N_{+2}|^2.
\eea
Next, for $R\geq 10m$, let $\chi_R=\chi_R(r)$ denote a smooth nonnegative cut-off function satisfying
\bea
\lab{properties:chiRcutoff}
\chi_R=1 \,\,\, \text{for} \,\,\,  r\leq R, \qquad \chi_R=0\,\,\,   \text{for} \,\,\,   r\geq 2R, \qquad \pr_r^k \chi_R=O(R^{-k}) \,\,\, \text{for} \,\,\, k\in\mathbb{N}.
\eea
Also, let $\widetilde{\chi}_{\tau_1}=\widetilde{\chi}_{\tau_1}(\tau)$ denote a smooth nonnegative cut-off function  satisfying $\widetilde\chi_{\tau_1}(\tau)=1$ for $\tau\leq \tau_1-1$ and $\widetilde\chi_{\tau_1}=0$ for $\tau\geq \tau_1$.
 Then, denoting by ${\bf{T}}_{+2}$ the Teukolsky wave operator on the LHS of \eqref{eq:TeukolskyequationforAandAbintensorialforminKerr:inhomogenouscase}, 
we consider the following inhomogeneous Teukolsky equation 
\bea
\lab{eq:tensorialTeuinKerr:cutoffinR:+2}
{\bf{T}}_{+2}\phihp=\Nhp,
\eea
where $\phihp$ has the same initial data as $\pmb\phi_{+2}$ on $\Si(\tau_0)$ and where $\Nhp$ is defined by 
\bea
\lab{eq:formulaofNhp:weakMoraKerr:pf}
\Nhp &:=& \widetilde{\chi}_{\tau_1}(\tau) \big(\chi_R\N_{+2} + [{\bf{T}}_{+2}, \chi_R]\pmb\phi_{+2}\big)\nn\\
&=&\widetilde{\chi}_{\tau_1}(\tau) \Big(\chi_R\N_{+2} +O(1)\pr_r(\chi_R)\nab_3\pmb\phi_{+2}+O(r^{-1})\pr_r (\chi_R)\dk^{\leq 1}\pmb\phi_{+2}\nn\\
&&+O(1)\pr_r^2 (\chi_R)\pmb\phi_{+2}\Big).
\eea
 Thus, by causality, we have
\bea
\lab{relation:phihpwithchiRphi+2:pf}
\phihp=\chi_R \pmb\phi_{+2}, \qquad \text{on}\quad \MM(\tau_0,\tau_1-1).
\eea

Next, applying the estimate \eqref{eq:weakMorawetzforTeukolskyfromMillet:bis:plus2}, with the substitution $\de\to \frac{\de}{2}$, to \eqref{eq:tensorialTeuinKerr:cutoffinR:+2}, we deduce
\beaa
\int_{\MM(\tau\geq\tau_0)}r^{-11+\frac{\de}{2}}|\dk^{\leq 3}\phihp|^2 &\les& \E^{(9)}[r^{-4}\phihp](\tau_0) +\E^{(7)}[r^{-\frac{7}{2}+\frac{\de}{4}}\nab_{\pr_r}(r\phihp)](\tau_0)\nn\\
&&+\int_{\MM(\tau\geq\tau_0)}r^{-7+\frac{\de}{2}}|\dk^{\leq 9}\Nhp|^2
\eeaa
which yields, in view of \eqref{properties:chiRcutoff}, \eqref{eq:formulaofNhp:weakMoraKerr:pf} and \eqref{relation:phihpwithchiRphi+2:pf}
\bea
\lab{eq:weakmorawetz:TeuKerr:compactregion:+2}
\int_{\MM_{r\leq R}(\tau_0,\tau_1-1)}r^{-11+\frac{\de}{2}}|\dk^{\leq 3}\pmb\phi_{+2}|^2 
&\les& R^{1+\frac{\de}{2}}\E_{r\leq 2R}^{(9)}[r^{-4}\pmb\phi_{+2}](\tau_0) +\int_{\MM_{r\leq 2R}(\tau_0,\tau_1)}r^{-7+\frac{\de}{2}}|\dk^{\leq 9}\N_{+2}|^2 \nn\\
&&+\int_{\MM_{R, 2R}(\tau_0,\tau_1)}r^{-3+\frac{\de}{2}}|\dk^{\leq 9}\nab_3(r^{-3}\pmb\phi_{+2})|^2\nn\\
&&+\int_{\MM_{R, 2R}(\tau_0,\tau_1)}r^{-11+\frac{\de}{2}}|\dk^{\leq 10}\pmb\phi_{+2}|^2.
\eea

Next, we control the last term on the RHS of  \eqref{eq:weakmorawetz:TeuKerr:compactregion:+2}. Restricting equation \eqref{eq:Laplacianofphiprewrittenbyderivativeofphip+1:pm2:weightedderi} to the case $s=+2$ and $p=0$ and to the Kerr background, which yields that all the coefficients dependent on $\Ga_b$ and $\Ga_g$ vanish and
\begin{equation}
\begin{aligned}
\lab{choiceoftensorsandsource:weakMora:Kerr:pf}
\phiplus{0}={}&|q|^{-4}\pmb\phi_{+2}, & \phiplus{1}={}&\frac{r^2}{|q|^2}\bigg({\frac{q}{\bar{q}}}r \nab_{3} r{\frac{\bar{q}}{q}}\bigg)\bigg(\frac{r^2}{|q|^2}\bigg)^{-2}\pmb\phi_{+2}^{(0)}, \\
 \N_{T,+2}^{(0)}={}&0, & \N_{W,+2}^{(0)}={}&|q|^{-4}\N_{+2},
\end{aligned}
\end{equation} 
we obtain
\bea
\lab{eq:Laplacianofphiprewrittenbyderivativeofphip+1:pm2:weightedderi:inKerr}
&&\qs \De_2\dk^{\reg}\phiplus{0}- 2\dk^{\reg}\phiplus{0}\nn\\
&=& O(r^{-1})(r\nab_4)^{\leq 1}\dk^{\leq \reg}\phiplus{1}
+O(1)\big(\nab_3\dk^{\leq \reg}\phiplus{0}, r^{-1}\dk^{\leq \reg+1}\phiplus{0}\big)
+\de_{\reg\geq 1}O(1)(r\nab)^{\leq 1}\dk^{\leq \reg-1}\phiplus{0}
\nn\\
&&+ O(r^2)\dk^{\leq \reg}\L_{+2}^{(0)}[\pmb\phi_{+2}] + O(r^2)\dk^{\leq \reg}\N_{W,+2}^{(0)}.
\eea 
Multiplying on both sides of \eqref{eq:Laplacianofphiprewrittenbyderivativeofphip+1:pm2:weightedderi:inKerr} by $r^{-3+\de}\ov{\dk^{\reg}\phiplus{0}}$, taking the real part, integrating over $\MM(\tau_1,\tau_2)$, and in view of the expression of $\L_{+2}^{(0)}[\pmb\phi_{+2}]$ in \eqref{eq:tensor:Lsn:onlye_2present:general:Kerrperturbation}, we deduce, for $\reg\leq 10$ and $0<\de\leq \frac{1}{3}$, 
\begin{align*}
\Ao[{r^{\frac{\de}{2}}(r\nab)^{\leq 1}}\dk^{\leq \reg}\phiplus{0}](\tau_0,\tau_1)
\les{} &
\int_{\MM(\tau_0,\tau_1)}\Big(r^{-1+\de}|\dk^{\leq \reg} \N_{W,{+2}}^{(0)}|+r^{-4+\de}|\dk^{\leq \reg+1}\phiplus{1}|\Big)|\dk^{\leq \reg} \phiplus{0}| \nn\\
&
+\widehat{\EMF}^{(\reg)}_{\de}[\phiplus{0}](\tau_0,\tau_1) +\de_{\reg\geq 1}\Ao[{r^{\frac{\de}{2}}(r\nab)^{\leq 1}}\dk^{\leq \reg-1}\phiplus{0}](\tau_0,\tau_1).
\end{align*}
Summing over $\reg\leq 10$, and in view of \eqref{choiceoftensorsandsource:weakMora:Kerr:pf}, this yields
\begin{align*}
\int_{\MM(\tau_0,\tau_1)} r^{-11+\de} |(r\nab)^{\leq 1}\dk^{\leq 10}\pmb\phi_{+2}|^2
\les{}&\widehat{\EMF}^{(10)}_{\de}[r^{-4}\pmb\phi_{+2}](\tau_0,\tau_1)
+\int_{\MM(\tau_0,\tau_1)}r^{-7+\de}|\dk^{\leq 10} \N_{{+2}}|^2\nn\\
&+ \int_{\MM(\tau_0,\tau_1)} r^{-3+\de}|\dk^{\leq 11} \nab_{3}(r^{-3}\pmb\phi_{+2})|^2 .
\end{align*}
In view of the fact that
\beaa
\widehat{\EMF}_{\de}^{(\reg)}[\pmb\psi](\tau_0,\tau_1)\les {\EMF}_{\de}^{(\reg+1)}[\pmb\psi](\tau_0,\tau_1), \qquad \forall \, \reg\in\mathbb{N},
\eeaa
we infer
\bea\lab{eq:weakmorawetz:TeuKerr:tocontrolcompactregion:+2}
\int_{\MM(\tau_0,\tau_1)} r^{-11+\de}|\dk^{\leq 10}\pmb\phi_{+2}|^2
&\les&{\EMF}^{(11)}_{\de}[r^{-4}\pmb\phi_{+2}](\tau_0,\tau_1)
+\int_{\MM(\tau_0,\tau_1)}r^{-7+\de}|\dk^{\leq 10} \N_{{+2}}|^2\nn\\
&&+ \int_{\MM(\tau_0,\tau_1)} r^{-3+\de}|\dk^{\leq 11} \nab_{3}(r^{-3}\pmb\phi_{+2})|^2 .
\eea

Finally, combining the estimates \eqref{eq:weakmorawetz:TeuKerr:compactregion:+2} and  \eqref{eq:weakmorawetz:TeuKerr:tocontrolcompactregion:+2} yields the desired estimate \eqref{eq:weakMorawetzforTeukolskyfromMillet:plus2:nolossinr}. This concludes the proof of Theorem  \ref{cor:weakMorawetzforTeukolskyfromMillet:bis}.


\bigskip
\footnotesize

Siyuan Ma, \par\nopagebreak
\textsc{State Key Laboratory of Mathematical Sciences, Academy of Mathematics and Systems Science, Chinese Academy of Sciences, Beijing 100190, China}\par\nopagebreak
\textit{E-mail address:} \href{mailto:siyuan.ma@amss.ac.cn}{siyuan.ma@amss.ac.cn}

\bigskip

J\'{e}r\'{e}mie Szeftel, \par\nopagebreak
\textsc{CNRS \& Laboratoire Jacques-Louis Lions, Sorbonne Universit\'{e},
4 place Jussieu 75005 Paris, France}\par\nopagebreak
\textit{E-mail address}: \href{mailto:jeremie.szeftel@sorbonne-universite.fr}{jeremie.szeftel@sorbonne-universite.fr}

\end{document}